\documentclass{amsbook}
\pdfoutput=1
\usepackage{amssymb}
\usepackage{mathrsfs}
\usepackage[all]{xy}
\UseComputerModernTips
\CompileMatrices
\usepackage{graphicx}
\usepackage{hyperref}
\newtheorem{theorem}{Theorem}[section]
\newtheorem{lemma}[theorem]{Lemma}
\newtheorem{proposition}[theorem]{Proposition}
\newtheorem{corollary}[theorem]{Corollary}
\theoremstyle{definition}
\newtheorem{definition}[theorem]{Definition}
\newtheorem{example}[theorem]{Example}
\newtheorem{examples}[theorem]{Examples}
\newtheorem{construction}[theorem]{Construction}
\theoremstyle{remark}
\newtheorem{remark}[theorem]{Remark}
\newtheorem{remarks}[theorem]{Remarks}
\newtheorem{note}[theorem]{Note}
\numberwithin{section}{chapter}
\usepackage{ifthen}
\usepackage{latexsym}

\renewcommand{\phi}{\varphi}

\newcommand{\bndry}{\partial}
\DeclareMathSymbol{\boxprod}{\mathbin}{AMSa}{"03} 
\newcommand{\mixprod}{\triangleleft}
\newcommand{\exboxprod}{\boxtimes}
\newcommand{\exsmsh}{\barwedge}

\newcommand{\convto}{\Rightarrow}
\newcommand{\dimpred}{\succcurlyeq}
\newcommand{\dirsum}{\oplus}
\newcommand{\Dirsum}{\bigoplus}
\newcommand{\disjunion}{\sqcup}

\newcommand{\from}{\leftarrow}

\newcommand{\hmtpc}{\simeq}
\newcommand{\homeo}{\approx}
\newcommand{\id}{\mathrm{id}}
\newcommand{\includesin}{\hookrightarrow}
\newcommand{\intersect}{\cap}
\newcommand{\iso}{\cong}
\newcommand{\laxto}{\Rightarrow}

\newcommand{\Loop}{\Omega}

\newcommand{\Mackey}[1]{\overline{#1}\vphantom{#1}}
\newcommand{\MackeyOp}[1]{{\underline {#1}}}

\newcommand{\onto}{\twoheadrightarrow}

\newcommand{\slant}{\setminus}
\newcommand{\smsh}{\wedge}
\newcommand{\Smsh}{\bigwedge}

\newcommand{\Susp}{\Sigma}
\newcommand{\susp}{\Sigma}
\newcommand{\tensor}{\otimes}
\newcommand{\union}{\cup}
\newcommand{\Union}{\bigcup}
\newcommand{\Wedge}{\bigvee}

\newcommand{\ParaQuot}[1]{/\!_{#1}\,}
\newcommand{\PQ}[1]{\ParaQuot{#1}}

\newcommand{\A}{{\mathscr A}}
\newcommand{\B}{{\mathscr B}}
\newcommand{\C}{{\mathscr C}}
\newcommand{\D}{{\mathscr D}}

\newcommand{\F}{{\mathscr F}}
\renewcommand{\H}{{\mathscr H}}
\newcommand{\K}{{\mathscr K}}

\newcommand{\R}{{\mathscr R}}
\newcommand{\RO}{{\mathscr R}O}
\newcommand{\Real}{\mathbb{R}}

\newcommand{\U}{{\mathscr U}}
\newcommand{\V}{{\mathscr V}}
\newcommand{\Vrep}{\mathbb{V}}
\newcommand{\W}{{\mathscr W}}

\newcommand{\sZ}{{\mathscr Z}}
\newcommand{\Z}{\mathbb{Z}}

\newcommand{\Lie}{\mathscr{L}}

\newcommand{\tE}{\tilde E}

\newcommand{\CH}{cH}
\newcommand{\tCH}{c\tilde H}

\newcommand{\Ab}{\text{\textit{Ab}}}
\newcommand{\Op}[1]{{#1}^\mathrm{op}}
\newcommand{\orb}[1]{{\mathscr{O}_{#1}}}

\newcommand{\stab}[1]{{\widehat{#1}}}
\newcommand{\sorb}[1]{{\widehat{\mathscr{O}}_{#1}}}

\newcommand{\ortho}[2]{{\overline{\mathscr V}}_{#1}%
                        \ifthenelse{\equal{#2}{}}{}{({#2})}}
\newcommand{\orthosmall}[2]{\mathscr V_{#1}%
                        \ifthenelse{\equal{#2}{}}{}{({#2})}}
\newcommand{\PL}[2]{{\overline{\mathscr {PL}}}_{#1}%
                        \ifthenelse{\equal{#2}{}}{}{({#2})}}
\newcommand{\PLsmall}[2]{\mathscr {PL}_{#1}%
                        \ifthenelse{\equal{#2}{}}{}{({#2})}}
\newcommand{\Top}[2]{\overline{\mathscr T}\!op_{#1}%
                        \ifthenelse{\equal{#2}{}}{}{({#2})}}
\newcommand{\Topsmall}[2]{{{\mathscr T}op}_{#1}%
                        \ifthenelse{\equal{#2}{}}{}{({#2})}}
\newcommand{\sphere}[2]{\overline{\mathscr F}_{#1}%
                        \ifthenelse{\equal{#2}{}}{}{({#2})}}
\newcommand{\spheresmall}[2]{{\mathscr F}_{#1}%
                        \ifthenelse{\equal{#2}{}}{}{({#2})}}

\newcommand{\ParaU}[1]{\K/{#1}}
\newcommand{\Para}[1]{\K_{#1}}
\newcommand{\Lax}[1]{\K^\lambda_{#1}}
\newcommand{\LaxU}[1]{\K^\lambda/{#1}}

\newcommand{\Spec}[2]{{\mathscr S}#1_{#2}}
\newcommand{\PreSpec}[2]{{\mathscr P}#1_{#2}}
\newcommand{\LaxPreSpec}[2]{{\mathscr P}^{\lambda}#1_{#2}}

\DeclareMathOperator*{\colim}{colim}
\DeclareMathOperator{\Ho}{Ho}
\DeclareMathOperator{\Hom}{Hom}

\DeclareMathOperator{\Inf}{Inf}

\DeclareMathOperator{\Tor}{Tor}
\DeclareMathOperator{\Ext}{Ext}

\begin{document}
\frontmatter
\title{Equivariant ordinary homology and cohomology}

\author{Steven R. Costenoble}
\address{Department of Mathematics\\103 Hofstra University\\
   Hempstead, NY 11549}
\email{Steven.R.Costenoble@hofstra.edu}

\author{Stefan Waner}
\email{Stefan.Waner@hofstra.edu}


\date{\today}
\maketitle
\setcounter{tocdepth}{2}
\tableofcontents

\mainmatter
\chapter*{Introduction}
\label{chap:intro}
%

Poincar\'e duality lies at the heart of the homological study of manifolds.
In the presence of the action of a group, it is well-known that Poincar\'e
duality fails to hold in ordinary, integer-graded (Bredon-Illman) equivariant homology.
(See \cite{Br:cohomology} and \cite{Ill:homology} for the construction
of these homology theories;
see \cite{DR:surgery} and \cite{Lu:transfgroups} for the failure
of Poincar\'e duality.)
It is our purpose here to extend the Bredon theory to one in which
a form of
Poincar\'e duality holds for all smooth compact $G$-manifolds, where
$G$ is any compact Lie group.
We did this for finite groups in \cite{CW:duality}, where we showed
that, for finite group actions and Mackey functor
coefficient systems, Bredon-Illman theory extends
to a theory graded on $RO(\Pi_G X)$, the group of
virtual representations of
the equivariant fundamental groupoid $\Pi_G X$ of a $G$-space $X$,
in which Poincar\'e duality holds for all smooth compact $G$-manifolds. 
That theory is an
extension of a well-known intermediate theory: $RO(G)$-graded equivariant ordinary
(co)homology 
(see \cite{LMM:roghomology}, \cite{Le:Hurewicz},
\cite{Le:emspaces} and \cite[Chapters~X, XI and~XIII]{May:alaska}).
In \cite{CW:spivaknormal} and \cite{CW:simpleduality}
we applied the theory to obtain $\pi$-$\pi$ theorems for
equivariant Poincar\'e duality spaces and equivariant simple Poincar\'e duality spaces.
In other words,
our equivariant ordinary theory is sensitive enough to 
yield a surgery theory that
parallels the non-equivariant theory \cite{Wal:surgery} modulo the usual transversality
issues (\cite{Was:difftopology}, \cite{Pe:pseudo}).

The $RO(G)$-graded extension of Bredon-Illman homology dates back to
1981 and is interesting in its own right: 
It is necessary to extend to $RO(G)$ grading to get
theories represented by equivariant spectra 
(\cite{LMS:eqhomotopy}, \cite{MM:orthogonal}),
and so they provide a minimal substrate in which we can take advantage
of equivariant
Spanier-Whitehead duality. 
However, the documentation of $RO(G)$-graded ordinary homology is
sparse and incomplete.
In particular, its behaviour with regard to products, restriction to fixed sets,
and restriction to subgroups is nowhere discussed in detail.
(An outline of the product structure is given in \cite[XIII.5]{May:alaska}.)

When $G$ is finite, the $RO(G)$-graded theory does yield Poincar\'e duality for smooth
$G$-manifolds that are modeled on a single representation of $G$ 
(see \cite{Pu:rogbordism},
\cite{Ko:rogbordism}, and \cite{Wan:ROGbordism}).
However, as noted in \cite{May:alaska},
ordinary homology and cohomology need not be dual when $G$ is an infinite
compact Lie group---if $DX$ is the $G$-spectrum dual to $X$, $H^G_V(X)$ and
$H_G^{-V}(DX)$ need not be isomorphic. Thus, equivariant
Poincar\'e duality between ordinary homology and ordinary cohomology cannot be
expected to hold
when $G$ is not finite, even for $G$-manifolds modeled on a single representation.
Although ordinary homology and cohomology are not dual, each
possesses a dual theory. These theories, which we refer to as
dual ordinary cohomology and homology, are mentioned in
\cite{May:alaska} but not described in any detail.

In short, the state of the art to date is this:
There are $RO(G)$-graded ordinary homology and cohomology theories
defined for compact Lie groups $G$.
Poincar\'e duality holds with these theories only when $G$ is finite and
the manifold is modeled on a single representation of $G$.
When $G$ is finite, there is an extension of the $RO(G)$-graded theory
to an $RO(\Pi_G X)$-graded theory, in which Poincar\'e duality holds
for all smooth compact $G$-manifolds.
This raises the following questions: How does Poincar\'e duality work
in the $RO(G)$-graded theories when $G$ is infinite?
How do the $RO(\Pi_G X)$-graded theories extend to the case when $G$ is infinite,
and how does Poincar\'e duality work in this context?
The goal of the present work is to answer these questions while also
giving a more complete and coherent account of what is already known.

In Chapter~1 we give
a reasonably complete account of $RO(G)$-graded
ordinary homology and cohomology as well as their associated dual theories,
both geometrically (from the cellular viewpoint) and homotopically (from the
represented viewpoint). 
In fact, we define a variety of cellular theories of which the ordinary
and dual theories appear as extreme cases.
We then discuss various change-of-groups results as well as 
products, accounts that have not appeared
heretofore in the literature.
Finally, we show that
Poincar\'e duality for manifolds modeled on a single representation
holds for general compact Lie groups in that, 
if $M$ is compact, equivariantly orientable (in a sense to be made precise),
and modeled locally on the $G$-representation $V$, one has
duality isomorphisms like
$H^W_G(M,\bndry M) \iso \H_{V-W}^G(M)$,
where $\H^G_*$ is the dual theory.

Our development of the $RO(\Pi_G X)$-graded theories is done in the context
of parametrized spaces and spectra. In Chapter~2 we review facts we need about
parametrized homotopy theory, relying heavily on work by
May and Sigurdsson
\cite{MaySig:parametrized}. However, we emphasize a point of view
not taken by May and Sigurdsson: that, rather than concentrating on maps
that strictly commute with projection to the basespace, we should take
seriously maps that commute only up to specified homotopies in the
basespace. These maps, which we call {\em lax maps}, we discuss in detail
in Chapter~2.

In Chapter~2 we also review basic facts about the equivariant fundamental
groupoid $\Pi_G X$ and its representations, as developed in detail in
\cite{CMW:orientation}.
In the case of $RO(G)$-graded ordinary theories it is well-known that
the coefficient systems we
use must be Mackey functors, i.e., functors on the stable orbit category.
In the $RO(\Pi_G X)$-graded theory the fundamental groupoid replaces
the orbit category, so we must also have a stable version of $\Pi_G X$,
which is developed in Chapter~2.

Chapter~3 then discusses the $RO(\Pi_G X)$-graded theories,
culminating in the Thom isomorphism and Poincar\'e duality theorems.
This chapter parallels Chapter~1 as much as possible and in many cases relies
on results from that chapter.
A word about orientability: The $RO(\Pi_G X)$-graded theories are, 
by their nature, twisted, so that the Thom isomorphism and
Poincar\'e duality theorems require no orientability assumptions.
Indeed, restricting to orientable fibrations or manifolds in this
context gives us no advantage---the fundamental groupoid is
woven so deeply into the theory of equivariant orientations
\cite{CMW:orientation} that there is no significant simplification
in considering the orientable case.

\subsection*{History}
This project has been in development for long enough that it is worth
saying a few words about its history and relation to other works.
We have long believed that equivariant Thom isomorphism and
Poincar\'e duality theorems require homology and cohomology theories
that take into account the varying local representations given
by the fibers of a general equivariant vector bundle.
This belief led to the work with Peter May on equivariant orientations
that became \cite{CMW:orientation},
and early versions of the orientation theory, as well as versions
of $RO(\Pi_G X)$-graded ordinary homology and cohomology,
appeared in \cite{CW:duality} and \cite{CW:thomiso}.
(See \cite{May:thomiso} for a different approach to Thom isomorphisms
suggested by Peter.)

While we were working on \cite{CMW:orientation} we had the present
work in mind as well.
Early on, we did not give parametrized homotopy theory a prominent role,
but as we worked out the details it became clearer that we should,
particularly when discussing how the theories are represented.
In an earlier version of this manuscript
\cite{CW:homologypreprint} we used parametrized spectra 
as the representing objects,
but our naive belief was that the theory of parametrized spectra
was an easy generalization of the nonparametrized theory.
At about the same time we posted that preprint, Peter had started looking
at parametrized homotopy theory for other reasons, and he warned us
that he had already found serious pitfalls and that the parametrized
theory was definitely {\em not} an easy generalization of the
nonparametrized case. Peter was joined in his efforts by
Johann Sigurdsson and their work became \cite{MaySig:parametrized}.
We are pleased that they incorporated some of the ideas from our
preprint, particularly our notion of homological duality,
and we are indebted to them for providing a firm foundation on which
we could build the present work.
Once \cite{MaySig:parametrized} appeared it became clear that
we needed to substantially rewrite our earlier manuscript. The
present version uses parametrized homotopy theory much more extensively,
and we have attempted to do so with the care that Peter has enjoined
us to take.
In several places \cite{MaySig:parametrized} gives better proofs
of results from our earlier preprint, so we have replaced our proofs
with references to theirs.

One thing that took us a long time to fully appreciate is the difficulty
of trying to construct a good theory of CW parametrized spectra,
which Peter and Johann warned about in the final chapter of \cite{MaySig:parametrized}.
After the failures of several attempts to construct such a theory, we have
retreated in this work to discussing only CW parametrized spaces, not spectra.
Thus, the homology and cohomology theories we discuss throughout this work
are restricted to be defined only on spaces.

\vskip 20pt
\noindent\today

\chapter{$RO(G)$-graded Ordinary Homology and Cohomology}
\label{chap:rog}
%
\section*{Introduction}

A construction of $RO(G)$-graded ordinary homology and cohomology was announced
in \cite{LMM:roghomology} and given in \cite{GM:eqhomotopy}.
In this chapter we give another construction, 
versions of which first appeared
in \cite{Wan:gcwv} for finite groups,
in \cite{Le:emspaces} for compact Lie groups,
and in outline form in \cite[Chapter~X]{May:alaska}.
We also give more details about
the behavior of these theories. 
Some of this material is well-known, but a good deal
is new, particularly the material on 
dimension functions, on change of groups, and on products. 

In Chapter~\ref{chap:ordone} we shall generalize this construction to give a theory
graded on a larger group and defined on spaces parametrized by a fixed basespace. 
In a very precise sense, 
the theory discussed in this chapter is the local
case of the more general one. We shall discuss various constructions in detail in this chapter,
where the setting is somewhat simpler,
to allow us to concentrate later
on those things that need to be changed in the parametrized context.

For Bredon's integer-graded equivariant cellular homology, the
appropriate notion of a ``$G$-complex'' is a $G$-CW complex in the
sense of \cite{Br:cohomology} and \cite{Ill:homology}.
This is a $G$-space built from cells
of the form $G/H\times D\sp n$ where $H\leq G$ and $G$ acts
trivially on the unit $n$-disc $D\sp n$. The cell complexes we use
to construct the $RO(G)$-graded extension of Bredon homology 
have cells of the form $G\times_H D(W)$ where $D(W)$ is the unit
disc in a (possibly nontrivial) representation $W$ of $H$.
There are several motivations for doing this that also provide suggestions
for what representations we should use.

One simple thought is that generalizing from integer grading to $RO(G)$ grading
should involve replacing integer discs $D(\Real^n)$ with discs of representations.
Evidence that this should work comes from thinking about what happens when we suspend
a $G$-CW complex by a non-trivial $G$-representation $V$---we turn cells of the
form $G/H\times D^n$ into cells of the form $G/H\times D(V+n)$.

Particularly compelling to us, though, is to consider equivariant Poincar\'e duality
from a geometric point of view when $G$ is infinite. Consider a compact $G$-manifold $M$ and
suppose, for the purposes of this chapter, that it is $V$-dimensional.
That is, there is a representation $V$ of $G$ such that, at each point $x\in M$,
the tangent plane at $x$ is isomorphic to $V$ as a representation of the stabilizer $G_x$
of $x$. Illman showed that we can triangulate $M$, making it a $G$-CW complex
with cells of the form $G/H\times D^n$ as above.
We can form a ``dual cell structure'' by taking 
as the top-dimensional cells the orbits of the closed stars of the original vertices in
the first barycentric subdivision, while the lower dimensional cells are
intersections of these. These cells have interesting actions.
For example, if $x$ is a vertex in the original triangulation, then
the dual cell centered at $x$ will have the form $G\times_{G_x} D(W)$, where
$W \dirsum \Lie(G/G_x) \iso V$ as representations of $G_x$, with
$\Lie(G/G_x)$ denoting the tangent representation of the orbit $G/G_x$ at
the identity coset. We can write this as $W\iso V - \Lie(G/G_x)$.
In general, the dual cells will have the form $G\times_{G_x} D(V-\Lie(G/G_x)-k)$
where $x$ is the center of a $k$-cell in the original triangulation.
Notice that the {\em total} dimension of such a cell, including the dimension of
the orbit, is $V-k$, which is appropriate for a cell dual to a $k$-dimensional
simplex in a $V$-dimensional manifold. 
This differs from the simple cases above where, for example,
a cell of the form $G/H\times D^n$ is considered to be $n$-dimensional in
integer-graded Bredon homology. The difference is in what we consider the orbit
to be contributing to the dimension of the cell. In the case of Bredon homology,
the orbit contributes nothing to the dimension, but in the case of our dual cells,
it contributes its whole geometric dimension.
In order to consider both possibilities and treat them in a uniform way,
we introduce the notion of a {\em dimension function} for $G$, which assigns
to each orbit $G/H$ an $H$-representation $\delta(G/H)$ that tells us how much
the orbit is to contribute to the dimension of a cell. 
For each such $\delta$ there will be corresponding $RO(G)$-graded cellular homology and
cohomology theories.
In the case of Bredon homology
we let $\delta(G/H) = 0$, but in the case of the ``dual homology'' suggested
above we let $\delta(G/H) = \Lie(G/H)$. It will also be useful, particularly when
discussing products, to allows cases between these two extremes,
but the reader is invited to think about only these two cases on first reading
of most of this chapter.

Throughout we understand all spaces to be $k$-spaces.
It is more common to restrict to compactly generated spaces, that is,
weak Hausdorff $k$-spaces, but the extra generality is needed
in subsequent chapters when we discuss parametrized spaces
\cite{MaySig:parametrized}.
Our ambient group $G$ is a compact Lie group and subgroups are
understood to be closed.
We write $\K$ for the category of $k$-spaces, $G\K$ for the category
of $G$-spaces and $G$-maps, $\K_*$ for the category of based
$k$-spaces, and $G\K_*$ for the category of based $G$-spaces
(with $G$-fixed basepoints) and basepoint-preserving $G$-maps.
If $X$ is an unbased $G$-space we write $X_+$ for $X$ with a disjoint
basepoint adjoined. If $Y$ is a based $G$-space and $V$ is a representation
of $G$, we write $\susp_G^V Y = \susp^V Y$ for $Y\smsh S^V$, the smash
product with the one-point compactification of $V$.

\section{Dimension functions}

In the introduction to this chapter we discussed the different possible contributions of the
dimension of an orbit to the dimension of a cell.
The least an orbit can contribute is nothing and the most is its
geometric dimension.
We set up a general setting with the following definitions.

\begin{definition}
If $G$ is a compact Lie group, let $\Lie(G)$ denote its tangent space at the identity
(i.e., its Lie algebra, but we won't use the Lie algebra structure here).
If $H$ is a subgroup of $G$, let $\Lie(G/H)$ denote the tangent space to $G/H$
at the identity coset $eH$.
\end{definition}

We organize all the spaces $\Lie(H/K)$, $K\leq H\leq G$, as follows.
We assume chosen an invariant metric on $G$ so that $\Lie(G)$ is an inner product
space. The (left) conjugation action of $G$ on itself induces an orthogonal action of
$G$ on $\Lie(G)$ (well-known as the {\em adjoint representation} of $G$).
We consider $\Lie(H)\subset \Lie(G)$ for  $H\leq G$, with
$\Lie(H)$ identified with the subspace of vectors in $\Lie(G)$ tangent to $H$.
$\Lie(H)$ is an $N_G H$-subspace of $\Lie(G)$, where $N_G H$ denotes the
normalizer of $H$ in $G$. If $K\leq H\leq G$, we then identify
\[
 \Lie(H/K) \iso \Lie(H)-\Lie(K),
\]
the orthogonal complement of $\Lie(K)$ in $\Lie(H)$.
Writing $N_G(H/K)$ for $N_G H\intersect N_G K$, $\Lie(H/K)$ is an $N_G(H/K)$-subspace
of $\Lie(G)$.
For any $g\in G$, we have $g\cdot\Lie(H/K) = \Lie(H^g/K^g)$ as subsets of $\Lie(G)$,
where $H^g = gHg^{-1}$.

Note that $N_G(H/K)$ acts on $H/K$ by conjugation. The induced action on $\Lie(H/K)$
is the same as the action specified above. 
Restricting to $K$, the conjugation action of $K$ on $H/K$ 
is the same as its action by left multiplication,
and the induced action on $\Lie(H/K)$ is the usual action of $K$ on the fiber
of the tangent bundle $H\times_K \Lie(H/K)$ of $H/K$ over the coset $eK$.

Note also that, if $K\leq H\leq G$, then we have $\Lie(G/H)\subset \Lie(G/K)$ and
\[
 \Lie(G/K) - \Lie(G/H) = \Lie(H/K).
\]
More generally, if $L\leq K\leq H\leq G$, we have 
$\Lie(H/K) \subset \Lie(H/L)$ and
\[
 \Lie(H/L) - \Lie(H/K) = \Lie(K/L).
\]

In the definition of dimension function below, we use a general collection of subgroups
of $G$. The most interesting case is the collection of all subgroups, but
the more general case will be needed at times.

\begin{definition}\label{def:dimensionfunction}
A {\em dimension function $\delta$} for a group $G$
consists of a subset $\F(\delta)$ of subgroups of $G$, 
closed under conjugation (but not necessarily under taking subgroups), and
an assignment to each $H\in\F$ of a sub-$H$-representation $\delta(G/H)\subset \Lie(G/H)$ such that
\begin{enumerate}
\item for each $K\leq H$ such that $K$, $H\in\F$,
\[
 \delta(G/H)\subset\delta(G/K)
\]
and
\[
 \Lie(G/H)-\delta(G/H) \subset \Lie(G/K)-\delta(G/K)
\]
and

\item for each $H\in\F$ and $g\in G$,
$\delta(G/H^g) = g\cdot\delta(G/H)$ (as subspaces of $\Lie(G)$).
\end{enumerate}
\end{definition}

Note that, if we let $\epsilon(G/H) = \Lie(G/H)-\delta(G/H)$, then
$\epsilon$ will also be a dimension function. We call this the
{\em dual} function to $\delta$, and usually write
$\Lie-\delta$ for the dual to $\delta$.

\begin{remarks}\label{rem:dimensionDef}
In the following remarks, let $K\leq H$ with $K$, $H\in\F(\delta)$.

\begin{enumerate}
\item
We could replace the requirement that
$(\Lie-\delta)(G/H)\subset (\Lie-\delta)(G/K)$
with the equivalent requirement that
\[
 \delta(G/K) - \delta(G/H) \subset \Lie(G/K) - \Lie(G/H).
\]
We leave the equivalence of the two conditions as an exercise for the reader.

\item
Rather than think of $\Lie(G/H)\subset\Lie(G/K)$, we could consider the projection
$\Lie(G/K)\onto \Lie(G/H)$ induced by the projection $G/K\to G/H$.
The first conditions in the definition are then equivalent to the conditions that
$\delta(G/K)$ maps onto $\delta(G/H)$ under the projection and that
$(\Lie-\delta)(G/K)$ maps onto $(\Lie-\delta)(G/H)$.

\item
Carrying the preceding idea further, we can consider a dimension function as
an assignment to each orbit $G/H$ of a subbundle $G\times_H\delta(G/H)$
of the tangent bundle $G\times_H\Lie(G/H)$. To each map of orbits
$G/K\to G/H$ we then have a specified bundle epimorphism
\[
 G\times_K\delta(G/K) \onto G\times_H\delta(G/H).
\]
We find it more convenient, however, to think in terms of inclusions as
in the definition above.

\end{enumerate}
\end{remarks}

\begin{definition}
Let $\F$ be a collection of subgroups of $G$ closed under conjugation.
\begin{enumerate}
\item Let $\bar\F$ denote the closure of $\F$ under taking subgroups, i.e.,
\[
 \bar\F = \{ H \leq G \mid H\leq K\text{ for some $K\in\F$} \}.
\]

\item
Let
$\orb{\F}$ denote the category of orbits $G/H$, $H\in\F$, and all $G$-maps between them.

\item
Let $\orb{G}$ denote the category of all orbits of $G$.
\end{enumerate}
\end{definition}

There are various conditions we can require of a dimension function, as listed
in the following definition.

\begin{definition}
Let $\delta$ be a dimension function for $G$.
\begin{enumerate}
\item If $\F(\delta)$ is closed under taking subgroups
(i.e., $\F(\delta) = \overline{\F(\delta)}$), so is a {\em family} of subgroups,
we say that $\delta$ is {\em familial}.

\item If $\F(\delta)$ is the collection of all subgroups of $G$, we say that
$\delta$ is {\em complete}.

\end{enumerate}
\end{definition}

If $K\leq H\leq G$ with $K$ and $H$ in $\F(\delta)$, 
we let
\[
 \delta(H/K) = \delta(G/K) - \delta(G/H).
\]
Then $\delta(H/K)$ is a representation of $K$ and
\[
 \delta(G/K) = \delta(G/H)\dirsum\delta(H/K).
\]
Further, as noted in Remarks~\ref{rem:dimensionDef}, we have
\[
 \delta(H/K) \subset \Lie(H/K).
\]

The two most obvious examples of dimension functions are
the complete functions $\delta = 0$, 
which makes each orbit 0-dimensional,
and $\delta = \Lie$, which assigns to each orbit its geometric dimension $\Lie(G/H)$.
Note that 0 and $\Lie$ are each other's duals.
We conjecture that these are, in fact, the only complete dimension functions.
It's not hard to show that this conjecture is true for tori.
The other interesting examples we have in mind are the product dimension functions
discussed in Section~\ref{subsec:intermediate}, which are not familial,
but which come up while discussing products in cohomology.

\begin{remarks}\label{rem:dimensionproperties}
\begin{enumerate}\item[]
\item
If $K\leq J\leq H \leq G$ with $K$, $J$, $H \in \F(\delta)$, then
\[
 \delta(H/K) = \delta(H/J) \dirsum \delta(J/K).
\]
We sometimes refer to this as the {\em additivity} of $\delta$.

\item
Let $L \leq K\leq J\leq H$ with $K$, $J$, $H\in\F(\delta)$.
From the fact that $\delta(H/J)$ is a subrepresentation of $\delta(H/K)$
and also a subrepresentation of $\Lie(H/J)$, we get the following inequalities
that will be used later:
\begin{align*}
 |\delta(H/J)^L| &\leq |\delta(H/K)^L| \\
 |\delta(H/J)^L| &\leq |\Lie(H/J)^L| \\
\end{align*}

\end{enumerate}
\end{remarks}

Finally, the following relation will be useful.

\begin{definition}\label{def:dimensionpo}
Let $\delta$ and $\epsilon$ be dimension functions for $G$.
We write $\delta\dimpred\epsilon$ if $\F(\delta)\subset \overline{\F(\epsilon)}$ and, 
for each $K\in\F(\delta)$ and $H\in\F(\epsilon)$
such that $K\leq H$, we have that $\epsilon(G/H)$ is $K$-isomorphic to a
subrepresentation of $\delta(G/K)$. (Note that we do not require it to be a subspace.)
\end{definition}

This relation is a partial order on familial dimension functions, but it is not transitive in general.
Also note that, if $\epsilon$ is familial, the last condition of the definition can
be replaced with the requirement that, for each $K\in\F(\delta)$, we have that
$\epsilon(G/K)$ is isomorphic to a subrepresentation of $\delta(G/K)$.

\section{Virtual representations}

For our purposes later when discussing the parametrized case, we also want to use the most
general grading, which is on virtual representations.
We want to view virtual representations $V\ominus W$ as first-class objects in their own right.
The following definition is a restriction of one given in \cite[\S 19]{CMW:orientation}.
Let $\U$ denote a complete $G$-universe, i.e., the sum of countably many copies of each
irreducible representation of $G$.

\begin{definition}\label{def:virtualGReps}
The category of {\em virtual representations of $G$} has as its objects
the pairs $(V,W)$ of finite-dimensional representations of $G$; we think of and often write
a pair $(V,W)$ as a formal difference $V\ominus W$.
A {\em (virtual) map} $V_1\ominus W_1 \to V_2\ominus W_2$ is the equivalence class
of a pair of orthogonal $G$-isomorphisms
\begin{align*}
 f\colon {}& V_1 \dirsum Z_1 \to V_2 \dirsum Z_2 \\
 g\colon {}& W_1 \dirsum Z_1 \to W_2 \dirsum Z_2
\end{align*}
where $Z_1$ and $Z_2$ are finite-dimensional $G$-subspaces of $\U$. The equivalence relation
on such pairs $(f,g)$ is generated by two basic relations, the first being $G$-homotopy
through orthogonal maps.
The second relation is as follows. Let $T_1$ be a finite-dimensional $G$-subspace of $\U$ orthogonal to $Z_1$
and $T_2$ a finite-dimensional $G$-subspace of $\U$ orthogonal to $Z_2$, and let $k\colon T_1\to T_2$
be an orthogonal $G$-isomorphism. Then $(f,g)$ is equivalent to the ``suspension'' $(f\dirsum k, g\dirsum k)$ where
\begin{align*}
 f\dirsum k\colon {}& V_1 \dirsum (Z_1 + T_1) \to V_2 \dirsum (Z_2 + T_2) \\
 g\dirsum k\colon {}& W_1 \dirsum (Z_1 + T_1) \to W_2 \dirsum (Z_2 + T_2).
\end{align*}
Composition of virtual maps is defined by suspending until the pairs can be composed as
pairs of orthogonal $G$-isomorphisms.
\end{definition}

That this gives a well-defined category follows from
\cite[19.2]{CMW:orientation}.
We say that two virtual representations are
{\em stably equivalent} if they
are equivalent in this category.
If $\alpha = V\ominus W$, we write $|\alpha| = |V| - |W|$ for the integer dimension of $\alpha$.

\begin{remark}
We would get an isomorphic category if, instead of using a pair of maps $(f,g)$ as above,
we used a single orthogonal isomorphism $F\colon V_1\dirsum W_2\dirsum Z \to V_2\dirsum W_1\dirsum Z$
and the equivalence relation generated by $G$-homotopy and suspension by identity maps.
\end{remark}

If $V$ is a representation of $G$, we consider $V$ as the virtual representation $V\ominus 0$
and we write $-V$ for $0\ominus V$.
If $n$ is an integer, we write $V\ominus W +n$ for $(V\dirsum\Real^n, W)$ if $n\geq 0$
and $(V,W\dirsum \Real^{|n|})$ if $n<0$.
Note that there will be times where we will have to be careful to specify exactly
what copy of $\Real^n$ we are using.
In general, we write
$(V_1\ominus W_1) + (V_2\ominus W_2) = (V_1\dirsum V_2) \ominus (W_1\dirsum W_2)$
and $-(V \ominus W) = W \ominus V$.

\section{General cell complexes}

If $V$ is an orthogonal representation of $G$, we write $D(V)$ for its unit disc and $S(V)$ for its unit sphere.

\begin{definition}
Let $\alpha=V\ominus W$ be a virtual representation of $G$ and let
$\delta$ be a dimension function for $G$ with underlying collection of subgroups $\F(\delta)$.
\begin{enumerate}
\item
A subgroup $H$ of $G$ is {\em $\delta$-$\alpha$-admissible,} or simply {\em admissible},
if $H\in\F(\delta)$ and $\alpha - \delta(G/H) + n$ is 
stably equivalent to
an actual $H$-representation for some integer $n$.

\item
A {\em $\delta$-$\alpha$-cell} is a pair of $G$-spaces of the form
$(G\times_H D(Z), G\times_H S(Z))$
where $H\in\F(\delta)$ and $Z$ is an actual representation of $H$ such that
$Z$ is stably equivalent to
$\alpha-\delta(G/H)+n$ for some integer $n$.
Note that $H$ must, by definition, be admissible for such a cell to exist.
We say that the {\em dimension} of such a cell is $\alpha+n$.
\end{enumerate}
\end{definition}

As a more concise notation, we shall write
\[
 G\times_H \bar D(Z) = (G\times_H D(Z), G\times_H S(Z)).
\]
Note that the dimension of $G\times_H \bar D(Z)$ is $Z + \delta(G/H)$.
So, we think of $\delta(G/H)$ as the contribution of $G/H$ to the
total dimension of  $G\times_H \bar D(Z)$.

\begin{definition}
Let $\alpha$ be a virtual representation of $G$ and let
$\delta$ be a dimension function for $G$.
\begin{enumerate}
\item
 A {\em $\delta$-$G$-CW($\alpha$) complex} is a $G$-space $X$ together with a
decomposition $X = \colim_n X\sp n$ such that
 \begin{enumerate}
 \item 
$X^0$ is a disjoint union of $G$-orbits $G/H$ such that
$H\in\F(\delta)$, $\delta(G/H) = 0$, and $H$ acts trivially on $\alpha$.
 \item 
$X^n$ is obtained from $X^{n-1}$ by
attaching $\delta$-$\alpha$-cells of dimension 
$\alpha-|\alpha|+n$ along their boundary spheres.
Notice that the boundary sphere may be empty and the disc simply an orbit, in which case attaching the
cell means taking the disjoint union with the orbit.
 \end{enumerate} 
For notational convenience, and to remind ourselves of the role of $\alpha$,
we shall also write $X^{\alpha+n}$ for $X^{|\alpha|+n}$.

\item
A {\em relative} $\delta$-$G$-CW($\alpha$) complex is a
pair $(X,A)$ where $X = \colim_n X^n$, $X^0$ is the disjoint
union of $A$ with orbits as in (a) above, and cells are attached as in
(b). 

\item
A {\em based} $\delta$-$G$-CW($\alpha$) complex is a relative $\delta$-$G$-CW($\alpha$) complex $(X,*)$
where $*$ denotes a $G$-fixed basepoint.

\item
If we allow cells of any dimension to be attached at each stage,
we get the weaker notions of absolute, relative, or based {\em $\delta$-$\alpha$-cell complex}
(or $\delta$-$G$-$\alpha$-cell complex if we need to specify the group).

\item
If $X$ is a $\delta$-$G$-CW($\alpha$) or $\delta$-$\alpha$-cell complex with cells only of dimension
less than or equal to $\alpha+n$, we say that
$X$ is {\em $(\alpha+n)$-dimensional}.
\end{enumerate}
 \end{definition}

A simple observation, but useful:
If $(X,A)$ is a relative cell complex, then $(X/A,*)$ is a based cell complex.

\begin{examples}\label{ex:vcomplexes}
\begin{enumerate}\item[]
\item\label{item:vbundle} 
Suppose that $B$ is a $G$-CW complex in the sense of
\cite{Br:cohomology} and \cite{Ill:homology}, which is to say,
a $0$-$G$-CW($0$) complex by our definition; we shall call
such complexes {\em ordinary $G$-CW complexes} or simply $G$-CW complexes
(see also Section~\ref{sec:particularcomplexes}).
Let $p\colon E\to B$ be a $V$-bundle, that is,
a $G$-vector bundle such that, for each $b\in B$, the fiber $E_b$ over $b$
is isomorphic to $V$ as a representation of the stabilizer $G_b$.
Then the Thom space $T(p)$ is a based $0$-$G$-CW$(V)$ complex
(more simply, a based $G$-CW$(V)$ complex) with the
compactification point as basepoint and with a cell of dimension $V+n$
corresponding to each cell of dimension $n$ in $B$.
This correspondence underlies the Thom isomorphism we shall discuss
in Section~\ref{sec:vThomIso}.

\item\label{item:vmanifold}
Here is an example we mentioned in the introduction to this section.
As in \cite{Pu:rogbordism}, \cite{Ko:rogbordism}, and \cite{Wan:ROGbordism},
we say that a smooth $G$-manifold $M$ is a {\em $V$-manifold} if its tangent
bundle is a $V$-bundle.
In \cite{Ill:triangulation}, Illman showed that
any smooth $G$-manifold has a $G$-triangulation, so fix a triangulation of $M$.
We can them form the ``dual cell complex'' in the usual way, taking as the top-dimensional
cells the orbits of the closed stars of the original vertices in
the first barycentric subdivision of the triangulation, while the lower dimensional cells are
intersections of these.
Corresponding to each simplex in the original triangulation will be one cell
in the dual complex, the intersection of the stars of its vertices,
which intersects the simplex normally at its center. If $x$ is the center point of an $n$-simplex,
then the corresponding dual cell has the form $G\times_{G_x} D(W)$ for some $W$, and
the tangent plane at $x$ decomposes as
\[
 V \iso TM_x \iso \Real^n \dirsum \Lie(G/G_x) \dirsum W.
\]
Thus, the dual cell has the form $G\times_{G_x} D(V-\Lie(G/G_x)-n)$, and the dual
cells give $M$ the structure of an $\Lie$-$G$-CW$(V)$ complex.
The correspondence between the simplices (cells) of the triangulation and the
cells in the dual structure is the geometry underlying Poincar\'e duality,
which we shall discuss in Section~\ref{sec:vThomIso}.

\item
Suppose that $X$ is a $\delta$-$G$-CW$(\alpha)$ complex and $Y$ is
an $\epsilon$-$K$-CW$(\beta)$ complex, for two groups $G$ and $K$.
Then, as we shall discuss in more detail in the following section and in
Section~\ref{sec:vproducts}, the product cell structure on $X\times Y$
makes it a $(\delta\times\epsilon)$-$(G\times K)$-CW$(\alpha+\beta)$ complex
for an appropriate definition of $\delta\times\epsilon$.
(Note that, in general, $\delta\times\epsilon$ will be incomplete
even if $\delta$ and $\epsilon$ are complete.)
Even if we want to restrict our attention to the dimension functions $0$ and $\Lie$,
we still need to consider combinations like $0\times\Lie$, which are intermediate
between the extremes. It is for this reason that we work with general dimension
functions throughout.

\end{enumerate}
\end{examples}

\begin{definition}
Let $\delta$ and $\epsilon$ be dimension functions for $G$.
If $(X,A)$ is a relative $\delta$-$G$-CW($\alpha$) complex and
$(Y,B)$ is a relative $\epsilon$-$G$-CW($\alpha$) complex, then we say that
a map $f\colon (X,A)\to (Y,B)$ is {\em cellular}
if $f(X^{\alpha+n})\subset Y^{\alpha+n}$ for each $n$.
As special cases, we have the notions of cellular maps of absolute or based
CW complexes.
We write $G\W^{\delta,\alpha}$ for the category of $\delta$-$G$-CW$(\alpha)$ complexes
and cellular maps;
we write $G\W_*^{\delta,\alpha}$ for the category of based $\delta$-$G$-CW$(\alpha)$ complexes
and based cellular maps.
\end{definition}

\begin{definition}\label{def:dimensionRestriction}
Let $\delta$ be a dimension function for $G$ and let $H\in\F(\delta)$.
Define $\delta|H$, the {\em restriction of $\delta$ to $H$}, to be the
dimension function with
\[
 \F(\delta|H) = \{ K\in \F(\delta) \mid K\leq H \}
\]
and
\[
 (\delta|H)(H/K) = \delta(H/K) = \delta(G/K) - \delta(G/H)
\]
for $K\in\F(\delta|H)$.
\end{definition}

It is easy to check that $\delta|H$ is a dimension function for $H$.
Where the meaning is clear, we will usually write $\delta$ again for $\delta|H$.

The following observation is key to change-of-group homomorphisms and also simplifies
some arguments.

\begin{proposition}\label{prop:genInduction}
Let $\delta$ be a dimension function for $G$ and let $H\in\F(\delta)$.
If $X$ is a $\delta$-$H$-$(\alpha-\delta(G/H))$-cell complex, then
$G\times_H X$ is a $\delta$-$G$-$\alpha$-cell complex with corresponding cells.
Appying this construction to CW complexes and cellular maps, we get functors
\[
 G\times_H - \colon H\W^{\delta,\alpha-\delta(G/H)} \to G\W^{\delta,\alpha}
\]
and
\[
 G_+\smsh_H - \colon H\W_*^{\delta,\alpha-\delta(G/H)} \to G\W_*^{\delta,\alpha}.
\]
\end{proposition}

\begin{proof}
If $X$ is a $\delta$-$(\alpha-\delta(G/H))$-cell complex,
a typical cell has the form
\[
 H\times_K\bar D(\alpha-\delta(G/H)-\delta(H/K)+n)
\]
for $K\in\F(\delta|H)$. Applying $G\times_H -$ gives a corresponding cell of the form
\[
 G\times_H H\times_K \bar D(\alpha-\delta(G/H)-\delta(H/K) + n) 
  \iso G\times_K \bar D(\alpha-\delta(G/K) + n),
\]
using that $\delta(G/K) \iso \delta(G/H)\dirsum\delta(H/K)$.
Hence, $G\times_H X$ is a $\delta$-$\alpha$-cell complex.
If $X$ is CW, then so is $G\times_H X$ and, if $f\colon X\to Y$ is cellular,
then so is $G\times_H f$.
\end{proof}

\begin{definition}
If $n$ is an integer, a $G$-map $f\colon X\to Y$ is a {\em $\delta$-$(\alpha+n)$-equivalence}
if, for each admissable $H$ and each actual representation $V$ of $H$ stably equivalent to
$\alpha-\delta(G/H)+i$, with $i\leq n$,
every diagram of the following form is homotopic to one in which there exists a lift:
\[
 \xymatrix{
  G\times_H S(V) \ar[r] \ar[d] & X \ar[d]^f \\
  G\times_H D(V) \ar[r] & Y
 }
\]
We say that $f$ is a {\em $\delta$-weak$_{\alpha}$ equivalence}
if it a $\delta$-$(\alpha+n)$-equivalence for all $n$.
\end{definition}

Shortly we shall characterize $\delta$-weak$_{\alpha}$ equivalences
for a familial dimension function $\delta$ in terms
of their behavior on fixed points, but if not all $H$ are $\delta$-$\alpha$-admissible,
$\delta$-weak$_{\alpha}$ equivalences are not in general weak $G$-equivalences.

We have the following variant of the ``homotopy extension and
lifting property'' of \cite{May:foundations}.

\begin{lemma}[H.E.L.P.]\label{lem:HELP}
 Let $r\colon Y \to Z$ be a $\delta$-$(\alpha+n)$-equivalence.
Let $(X,A)$ be a relative 
$\delta$-$\alpha$-cell complex
of dimension $\alpha+n$.
If the following diagram commutes without the
dashed arrows, then there exist maps $\tilde g$ and $\tilde h$
making the diagram commute.
 \[
 \xymatrix{
  A \ar[rr]^-{i_0} \ar[dd] & & A\times I \ar'[d][dd] \ar[dl]_(.6)h & & 
  A \ar[ll]_-{i_1} \ar[dl]_(.6)g \ar[dd]  \\
   & Z & & Y \ar[ll]_(.3)r \\
  X \ar[ur]^f \ar[rr]_-{i_0} & & X\times I \ar@{-->}[ul]^(.6){\tilde h} & &
   X \ar[ll]^-{i_1} \ar@{-->}[ul]^(.6){\tilde g} \\
 }
 \]
The result remains true when $n = \infty$.
 \end{lemma}

\begin{figure}
 \centerline{\includegraphics[scale=0.8]{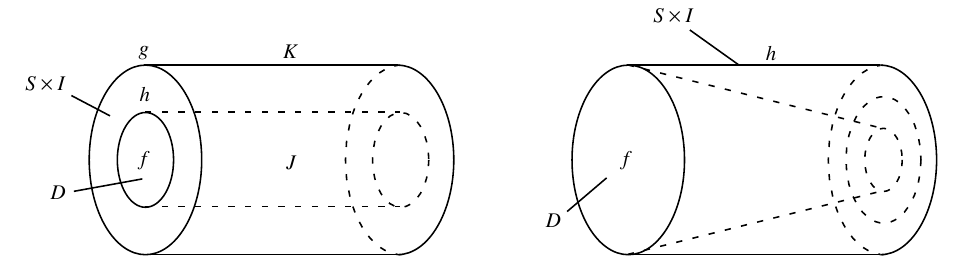}}
 \caption{}\label{fig:HELP}
\end{figure}

\begin{proof}
 The proof is by induction over the cells of $X$
not in $A$, so it suffices to consider the case
$(X, A) = G\times_H \bar D(V)$
where $V$ is stably equivalent to
$\alpha-\delta(G/H)+i$ with $i\leq n$.
As in the picture on the left of Figure~\ref{fig:HELP}
(where $D = D(V)$ and $S=S(V)$),
identify $D(V)\union (S(V)\times I)$ 
with $D(V)$ and use
$f$ and $h$ to define a map $G\times_H D(V)\to Z$ that lifts to $g$ on its boundary:
 \[
 \xymatrix{
  G\times_H S(V) \ar[r]^-g \ar[d] & Y \ar[d]^{r} \\
  G\times_H D(V) \ar[r]_-{f\union h} & Z. \\
 }
 \]
By assumption, this diagram is homotopic to one in which we can
find a lift $k\colon G\times_H D(V)\to Y$. Let $K$ be the homotopy of $g$ and $J$ be the homotopy of $f\union h$, as in Figure~\ref{fig:HELP}.
Now distort the cylinder in that figure to become the picture on the right.
On the far right end of this new cylinder, $k$ and $K$ combine to give a map 
$\tilde g\colon G\times_H D(V)\to Y$, while the whole cylinder describes a map
$\tilde h\colon G\times_H D(V)\times I\to Z$. It is easy to see that these maps make the diagram in the statement of the lemma commute.
 \end{proof}

\begin{theorem}[Whitehead]\label{thm:VWhitehead}
 \begin{enumerate} \item[] 

\item If $f\colon Y\to Z$ is a $\delta$-$(\alpha+n)$-equivalence and $X$ is an
$(\alpha+n-1)$-dimensional $\delta$-$\alpha$-cell complex,
then  
\[
 f_*\colon \pi G\K(X,Y) \to \pi G\K(X,Z)
\]
is an isomorphism, where $\pi G\K(-,-)$ denotes
the set of homotopy classes of $G$-maps.  
It is an epimorphism if $X$ is
$(\alpha+n)$-dimensional.

\item If $f\colon Y\to Z$ is a $\delta$-weak$_{\alpha}$ equivalence and $X$ is a
$\delta$-$\alpha$-cell complex,
then 
\[
 f_*\colon \pi G\K(X,Y) \to \pi G\K(X,Z)
\]
is an
isomorphism.  Therefore, any $\delta$-weak$_{\alpha}$ equivalence of 
$\delta$-$\alpha$-cell complexes
is a $G$-homotopy equivalence.

 \end{enumerate}
 \end{theorem}

\begin{proof}
 For the first part, apply the  H.E.L.P. lemma to the pair $(X,\emptyset)$ to get
surjectivity and to the pair $(X\times I, X\times\bndry I)$ to get
injectivity. The second part follows in the same way from the last statement of
the H.E.L.P.\ lemma.
 \end{proof}
 
\begin{theorem}[Relative Whitehead]\label{thm:relVWhitehead}
Let $A$ be a $G$-space and consider the category $A/G\K$ of $G$-spaces under $A$.
 \begin{enumerate}

\item If $f\colon Y\to Z$ is a $\delta$-$(\alpha+n)$-equivalence of $G$-spaces under $A$ and $(X,A)$ is an
$(\alpha+n-1)$-dimensional 
relative $\delta$-$\alpha$-cell complex,
then  
\[
 f_*\colon \pi A/G\K(X,Y) \to \pi A/G\K(X,Z)
\]
is an isomorphism, where $\pi A/G\K(-,-)$ denotes
the set of homotopy classes of $G$-maps under $A$.  
It is an epimorphism if $(X,A)$ is
$(\alpha+n)$-dimensional.

\item If $f\colon Y\to Z$ is a $\delta$-weak$_\alpha$ equivalence of $G$-spaces under $A$ and $(X,A)$ is a 
relative $\delta$-$\alpha$-cell complex,
then 
\[
 f_*\colon \pi A/G\K(X,Y) \to \pi A/G\K(X,Z)
\]
is an
isomorphism.  Therefore, any $\delta$-weak$_{\alpha}$ equivalence of 
relative $\delta$-$\alpha$-cell complexes $(X,A)\to (Y,A)$
is a $G$-homotopy equivalence rel $A$.

 \end{enumerate}
 \end{theorem}

\begin{proof}
 For the first part, apply the  H.E.L.P. lemma to the pair $(X,A)$ to get
surjectivity and to the pair $(X\times I, X\times\bndry I \union A\times I)$ to get
injectivity. The second part follows in the same way from the last statement of
the H.E.L.P.\ lemma.
 \end{proof}

Note that, when we take $A = *$ in the relative Whitehead theorem,
we get the special case of based $G$-spaces and based cell complexes.
One reason for stating many of the results of this section in relative form
is to include both the based and unbased cases.

Thus far, the results are essentially formal. To get more substantive results
we need to know more about what maps are $\delta$-$(\alpha+n)$-equivalences.
From this point on we need to assume that $\delta$ is familial when
discussing $\delta$-equivalences.
We start with the following important example,
which uses the relation defined in Definition~\ref{def:dimensionpo}.
Note that it applies to the most important case $\epsilon = \delta$ because $\delta\dimpred\delta$;
we will need the more general statements in later sections.

\begin{proposition}\label{prop:SkelEquivalence}
Let $\delta$ be a familial dimension function and let
$\epsilon$ be another dimension function (not necessarily familial)
with $\delta\dimpred\epsilon$.
If $(X,A)$ is a relative $\epsilon$-$G$-CW($\alpha$) complex then the inclusion $X^{\alpha+n}\to X$ is
a $\delta$-$(\alpha+n)$-equivalence.
\end{proposition}

\begin{proof}
By the usual induction, this reduces to showing that, if $Y$ is obtained from $B$ by attaching a cell
of the form $G\times_K D(W)$, where $W$ is stably equivalent
to $\alpha-\epsilon(G/K)+k$,
and $V$ is a representation of $H$ stably equivalent to
$\alpha-\delta(G/H)+i$ where $i<k$,
then any $G$-map of pairs $G\times_H \bar D(V)\to (Y,B)$
is homotopic rel boundary to a map into $B$.
(In the application to the present proposition, $k>n$ and $i\leq n$.)
For this we take an ordinary $G$-triangulation of $G\times_H D(V)$
and show by induction on the cells that we can homotope the map
to miss the orbit $G/K\times 0$ in the attached cell.
We use the fact that the fixed-set dimensions of $V$ must
equal those of $\alpha-\delta(G/H)+i$ and
the fixed set dimensions of $W$ equal those of
$\alpha-\epsilon(G/K)+k$.

The only simplices that might hit the orbit $G/K\times 0$ have the form
$G/J\times\Delta^j$ where $J\leq H$ and $J$ is subconjugate to $K$;
by replacing $K$ with a conjugate we may assume $J\leq K$ also.
Because such a simplex is embedded in $G\times_H D(V)$, we have
\begin{align*}
  j &\leq |[\alpha - \delta(G/H) - \Lie(H/J)]^J| + i \\
    &\leq |[\alpha - \delta(G/H) - \delta(H/J)]^J| + i \\
    &= |[\alpha - \delta(G/J)]^J| + i \\
    &\leq |[\alpha - \epsilon(G/K)]^J| + i \\
    &< |[\alpha - \epsilon(G/K)]^J| + k.
\end{align*}
Here we use that $|\delta(H/J)^J| \leq |\Lie(H/J)^J|$ as in
Remarks~\ref{rem:dimensionproperties} and that
$|\epsilon(G/K)^J| \leq |\delta(G/J)^J|$ because $\delta\dimpred\epsilon$.
Finally, note that $|[\alpha - \epsilon(G/K)]^J| + k$ is the codimension of
$(G/K\times 0)^J$ in $Y^J$,
so the desired homotopy exists for dimensional reasons.
\end{proof}

\begin{theorem}[Cellular Approximation of Maps]\label{thm:dualMapApprox}
Let $\delta$ be a familial dimension function and let
$\epsilon$ be any dimension function such that $\delta\dimpred\epsilon$.
Let $(X,A)$ be a relative $\delta$-$G$-CW($\alpha$) complex and let
$(Y,B)$ be a relative $\epsilon$-$G$-CW($\alpha$) complex.
Let $f\colon (X,A)\to (Y,B)$ be a $G$-map and suppose given
a subcomplex $(Z,A)\subset (X,A)$ and a $G$-homotopy $h$ of $f|Z$ to a cellular map. 
Then $h$ can be extended to a $G$-homotopy of $f$ to a cellular map.
\end{theorem}

\begin{proof}
 This follows by induction on skeleta, using the H.E.L.P. lemma and 
Proposition~\ref{prop:SkelEquivalence}
applied to the inclusion $Y^{\alpha+n}\to Y$.
\end{proof}

We digress now to discuss how $\delta$-$(\alpha+n)$-equivalence is related to conditions
on fixed sets. We first need the following approximation result.

\begin{lemma}\label{lem:vdiskComplex}
Let $\delta$ be a familial dimension function.
If $V$ is an actual representation of $H$ stably equivalent
to $\alpha-\delta(G/H)+n$, then
the pair $G\times_H \bar D(V)$ is $G$-homotopy equivalent to 
a pair of $\delta$-$G$-CW($\alpha$) complexes
$(X,A)$ where $X$ has dimension $\alpha+n$ and $A$ has dimension $\alpha+n-1$.
\end{lemma}

\begin{proof}
We first note that, by Proposition~\ref{prop:genInduction},
it suffices to show that the $H$-pair
$\bar D(V)$ is $H$-homotopy equivalent to a pair of 
$\delta$-$H$-CW($V$) complexes of dimensions $V$ and $V-1$, respectively. 

We'll prove a slightly more general result (necessary for an inductive
argument):
If $W$ is a sub-$H$-representation of $V$,
then the pair $\bar D(V-W)\times D(W)$ is $H$-homotopy
equivalent to a pair of $\delta$-$H$-CW($V$) complexes of dimensions $V$
and $V-1$.
We shall prove this for every subgroup $H$ of $G$ and every $V$ and $W$, by induction on $H$.

The beginning of the induction is the case where $H = e$. In this
case, 
\[
 \bar D(V-W)\times D(W)
  = (D(\Real^n),S(\Real^{n-m})\times D(\Real^m))
\]
for some $n$ and $m$, and $D(\Real^n)$
can be given a nonequivariant triangulation of
dimension $n$ with $S(\Real^{n-m})\times D(\Real^m)$ a subcomplex of dimension $n-1$.

Now assume the result for all proper subgroups of $H$.
If $V^H \neq 0$, write
$n = |(V-W)^H|$ and $m = |W^H|$, and then we can write
\[
\bar D(V-W)\times D(W)
 \iso \bar D(V-W-n)\times D(W-m)\times \bar D(\Real^n)\times D(\Real^m).
\]
Using a non-equivariant triangulation of $\bar D(\Real^n)\times D(\Real^m)$, we see
that it suffices to find the desired cell structure on
$\bar D(V-W-n)\times D(W-m)$.

This reduces us to the case where $V$ contains no $H$-trivial summands.
It suffices now to show that 
$(S(V),S(V-W)\times D(W))$ 
is equivalent to a pair of $\delta$-$G$-CW($V$)
complexes of dimension $V-1$,
since we then get $D(V)$ by attaching one
cell of the form $D(V)$, i.e., the
interior.
Give $(S(V),S(V-W)\times D(W))$ 
an ordinary $H$-triangulation and then take the dual cell structure.
A typical dual cell will have the form
$H\times_J \bar D(V-\Lie(H/J)-j)$ where $J\subset H$ and $j>0$.
Such a cell is $H$-homotopy equivalent to the pair
\begin{multline*}
 H\times_J \bar D(V-\Lie(H/J)-j)\times D(\Lie(H/J)-\delta(H/J)) = \\
  H\times_J \bar D(V-\delta(H/J)-(\Lie(H/J)-\delta(H/J))-j) \times D(\Lie(H/J)-\delta(H/J)).
\end{multline*}
(We use here that $\delta(H/J)\subset\Lie(H/J)$.)
Note that $J$ must be a proper subgroup of $H$ because
$V$ has no $H$-trivial summand, yet its sphere, which contains the cell
$H\times_J \bar D(V-\Lie(H/J)-j)$, has a $J$-fixed point.
With $\Lie(H/J)-\delta(H/J)$ playing the role of $W$ and
$V-\delta(H/J)-j$ playing the role of $V$, we can now apply the inductive hypothesis
to say that the $J$-pair
\[
 \bar D(V-\delta(H/J)-(\Lie(H/J)-\delta(H/J))-j) \times D(\Lie(H/J)-\delta(H/J))
\]
is $J$-homotopy equivalent to a pair of $\delta$-$J$-CW($V-\delta(H/J)$) complexes of
dimensions $V-\delta(H/J)-j$ and $V-\delta(H/J)-j-1$.
Applying $H\times_J -$ and using Proposition~\ref{prop:genInduction} again,
we get the structure we are seeking on $H\times_J \bar D(V-\Lie(H/J)-j)\times D(\Lie(H/J)-\delta(H/J))$.

Using induction on the dual cell structure
of $(S(V),S(V-W)\times D(W))$, 
we replace each dual cell with an
$H$-homotopy equivalent pair of complexes as above, using a cellular approximation of the attaching map.
This constructs a pair of $\delta$-$H$-CW($V$) complexes of dimension $V-1$,
$H$-homotopy equivalent to
$(S(V),S(V-W)\times D(W))$, as required.
 \end{proof}

\begin{corollary}\label{cor:ordinaryCellAsGeneral}
If $\delta$ is familial,
$H$ is a $\delta$-$\alpha$-admissible subgroup of $G$, and $n$ is any nonnegative integer,
then the pair $G\times_H \bar D^n$ is $G$-homotopy equivalent to a pair
of $\delta$-$G$-CW($\alpha$) complexes of dimensions $\alpha - |[\alpha-\delta(G/H)]^H|+n$
and $\alpha - |[\alpha-\delta(G/H)]^H|+n-1$.
\end{corollary}

\begin{proof}
Because $H$ is admissible, $\alpha-\delta(G/H)+k$ is stably 
equivalent to an actual representation for some $k$; in fact,
the smallest such $k$ is $-|[\alpha-\delta(G/H)]^H|$,
so let $V$ be an actual representation of $H$ stably equivalent to
$\alpha-\delta(G/H)-|[\alpha-\delta(G/H)]^H|$.
By the preceding lemma (applied to $H$ rather than $G$),
$D(V)$ is $H$-homotopy equivalent to
a $\delta$-$H$-CW($\alpha-\delta(G/H)$) complex of dimension
$\alpha-\delta(G/H) - |[\alpha-\delta(G/H)]^H|$.
It follows that
\[
 G\times_H \bar D^n \hmtpc G\times_H (\bar D^n\times D(V))
\]
is $G$-homotopy equivalent to a pair of $\delta$-$G$-CW($\alpha$) complexes
of dimensions $\alpha - |[\alpha-\delta(G/H)]^H| + n$
and $\alpha - |[\alpha-\delta(G/H)]^H| + n-1$ as claimed.
\end{proof}

\begin{lemma}\label{lem:generalCellAsOrdinary}
If $\delta$ is familial and $\alpha-\delta(G/H)+n$ is stably equivalent to
an actual representation $V$ of $H$, then
the pair $G\times_H \bar D(V)$ has an ordinary $G$-triangulation
in which the simplices have the form $G\times_J \Delta^j$ with $J$ $\delta$-$\alpha$-admissible
and $j\leq |[\alpha-\delta(G/J)]^J| + n$ [$j\leq |[\alpha-\delta(G/J)]^J| + n - 1$
for a cell in $G\times_H S(V)$].
\end{lemma}

\begin{proof}
Take any $G$-triangulation of the pair $G\times_H \bar D(V)$ and
let $G\times_J \Delta^j$ be any one of its simplices; we may assume that $J\leq H$.
Because $H/J\times\Delta^J$ embeds in $D(V)$ we must have
\[
 \delta(H/J)\subset \Lie(H/J)\subset V
\]
(i.e., $\delta(H/J)$ is isomorphic to a $J$-subspace of $V$), hence
\[
 \delta(G/J) \iso \delta(G/H) + \delta(H/J) \subset V+\delta(G/H),
\]
and $V+\delta(G/H)$ is stably equivalent to $\alpha+n$,
showing that $J$ is admissible.
We also must have
\begin{align*}
 j &\leq |[\alpha - \delta(G/H) - \Lie(H/J) + n]^J| \\
   &\leq |[\alpha - \delta(G/H) - \delta(H/J)]^J| + n \\
   &= |[\alpha-\delta(G/J)]^J| + n
\end{align*}
in general, with the upper limit one lower if the simplex is embedded in the boundary sphere.
\end{proof}

\begin{theorem}\label{thm:weakCharacterization}
Let $\delta$ be a familial dimension function.
A map $f\colon X\to Y$ is a $\delta$-$(\alpha+n)$-equivalence if and only if $f^H$ is
a nonequivariant $(|[\alpha-\delta(G/H)]^H|+n)$-equivalence for each $\delta$-$\alpha$-admissible subgroup $H$.
A map $f$ is a $\delta$-weak$_{\alpha}$ equivalence if and only if
$f^H$ is a weak equivalence for every $\delta$-$\alpha$-admissible subgroup $H$.
\end{theorem}

\begin{proof}
The last statement of the theorem follows directly from the first.

Suppose that $f$ is a $\delta$-$(\alpha+n)$-equivalence, let $H$ be an admissible subgroup,
and let $k = |[\alpha-\delta(G/H)]^H|+n$. If $k < 0$ then $f$ is a $k$-equivalence vacuously,
so suppose $k\geq 0$. Then, by Corollary~\ref{cor:ordinaryCellAsGeneral},
for $0\leq i \leq k$, $G\times_H\bar D^i$ is equivalent to a relative
$\delta$-$G$-CW($\alpha$) complex of dimension
$\alpha - |[\alpha-\delta(G/H)]^H|+i = \alpha + n - (k - i)$.
By the relative Whitehead theorem, it follows that
any map from $\bar D^i$ to $f^H$ is homotopic to one in which we can find a lift,
hence $f^H$ is a $k$-equivalence.

Conversely, suppose that $f$ has the property that $f^H$ is
a nonequivariant $(|[\alpha-\delta(G/H)]^H|+n)$-equivalence for each admissible subgroup $H$.
Consider a pair $G\times_H \bar D(V)$ where $V$ is an actual representation
of $H$ stably equivalent to $\alpha-\delta(G/H)+i$,
where $i\leq n$.
By Lemma~\ref{lem:generalCellAsOrdinary}, this pair has a $G$-triangulation
in which a typical simplex $G/J\times\Delta^j$ has $J$ admissible and
$j\leq |[\alpha-\delta(G/J)]^J| + i$. It follows from induction on this cell structure and
our assumption on $f$, that any map from
$G\times_H \bar D(V)$ to $f$ is homotopic to one in which we can find a lift, hence
$f$ is a $\delta$-$(\alpha+n)$-equivalence.
\end{proof}

One reason we need this theorem is to derive the following corollary, which shows that
$\delta$-weak$_{\alpha}$ equivalence deserves the name.

\begin{corollary}[Two-out-of-three Property]
Let $\delta$ be familial and consider two $G$-maps
$f\colon X\to Y$ and $g\colon Y\to Z$.
If two of $f$, $g$, and $gf$ are $\delta$-weak$_{\alpha}$ equivalences, then so is the third.
\qed
\end{corollary}

Recall from \cite[II.2]{LMS:eqhomotopy} that, 
if $\F$ is a collection of subgroups of $G$, we say that a $G$-map $f\colon X\to Y$
is an {\em $\F$-equivalence} if $f^H\colon X^H\to Y^H$ is a (non-equivariant) weak
equivalence for all $H\in\F$.
Given $\delta$ and $\alpha$, let $\F(\delta,\alpha)$ denote the collection of
$\delta$-$\alpha$-admissible subgroups of $G$. Then a restatement of the last part
of Theorem~\ref{thm:weakCharacterization} is that
$f$ is a $\delta$-weak$_\alpha$ equivalence if and only if it is
an $\F(\delta,\alpha)$-equivalence.
Note for later use that, for a fixed $\alpha$, there exists an actual representation $V$
such that $\F(\delta,\alpha+V) = \F(\delta)$.

Recall also from \cite{LMS:eqhomotopy} that, if $\F$ is a family, there is a universal
$\F$-space $E\F$, characterized by its fixed sets: $(E\F)^H = \emptyset$ if $H\notin\F$
and $(E\F)^H$ is contractible if $H\in\F$.

\begin{corollary}
Let $\delta$ be a familial dimension function for $G$.
\begin{enumerate}
\item
If $X$ is any $G$-space then the projection $E\F(\delta)\times X\to X$ is
a $\delta$-weak$_\alpha$ equivalence.

\item
Suppose that $\F(\delta,\alpha)$ is a family (e.g., in the case $\F(\delta,\alpha)=\F(\delta)$).
If $f\colon X\to Y$ is a $\delta$-weak$_\alpha$ equivalence, then
$1\times f\colon E\F(\delta,\alpha)\times X\to E\F(\delta,\alpha)\times Y$ is
a weak $G$-equivalence.
\qed
\end{enumerate}
\end{corollary}

Returning from our digression into weak equivalences, we now
discuss approximation of spaces, beginning with the following general result.

\begin{theorem}\label{thm:generalCellularApprox}
Let $\delta$ be a familial dimension function,
let $(A,P)$ be a relative $\delta$-$G$-CW($\alpha$) complex, 
let $X$ be a $G$-space, and let $f\colon A\to X$ be a map.
Then there exists a relative $\delta$-$G$-CW($\alpha$) complex $(Y,P)$, containing $(A,P)$ as a subcomplex,
and a $\delta$-weak$_{\alpha}$ equivalence $g\colon Y\to X$ extending $f$.
\end{theorem}

\begin{proof}
We use a variant of the usual technique of killing homotopy groups.
We start by letting
\[
 Y^{\alpha-|\alpha|} = A^{\alpha-|\alpha|}\disjunion \coprod G/K,
\]
where the coproduct runs over all subgroups $K$ such that $\delta(G/K) = 0$
and $K$ acts trivially on $\alpha$, and all maps $G/K\to X$. The map
$g\colon Y^{\alpha-|\alpha|}\to X$ is the one induced by those maps of orbits.
$(Y^{\alpha-|\alpha|},P)$ is then a relative 
$\delta$-$G$-CW($\alpha$) complex of dimension $\alpha-|\alpha|$,
containing $(A^{\alpha-|\alpha|},P)$, and $Y^{\alpha-|\alpha|}\to X$ is a
$\delta$-$(\alpha-|\alpha|)$-equivalence.

Inductively, suppose that we have constructed $(Y^{\alpha+n-1},P)$, a relative
$\delta$-$G$-CW($\alpha$) complex of dimension 
$\alpha+n-1$ containing $(A^{\alpha+n-1},P)$ as a subcomplex,
and a $\delta$-$(\alpha+n-1)$-equivalence 
$g\colon Y^{\alpha+n-1}\to X$ extending $f$ on $A^{\alpha+n-1}$.
Let
\[
 Y^{\alpha+n} = Y^{\alpha+n-1}\union A^{\alpha+n}\union \coprod G\times_K D(V),
\]
where the coproduct runs over all subgroups $K$ of $G$
and isomorphism classes of representations $V$ stably equivalent to
$\alpha - \delta(G/K) +n$, and all diagrams of the form
\[
 \xymatrix{
  G\times_K S(V) \ar[r] \ar[d] & Y^{\alpha+n-1} \ar[d] \\
  G\times_K D(V) \ar[r] & X.
 }
\]
The union that defines $Y^{\alpha+n}$ is along $A^{\alpha+n-1}\to Y^{\alpha+n-1}$ and the maps
$G\times_K S(V)\to Y^{\alpha+n-1}$ displayed above.
By construction, $(Y^{\alpha+n},P)$ is a relative $\delta$-$G$-CW($\alpha$) complex of dimension
$\alpha+n$ containing $(A^{\alpha+n},P)$ as a subcomplex.
We let $g\colon Y^{\alpha+n}\to X$ be the induced map and we claim that $g$
is a $\delta$-$(\alpha+n)$-equivalence.
To see this, consider any diagram of the following form, with
$V$ stably equivalent to $\alpha-\delta(G/K)+i$
and $i\leq n$:
\[
 \xymatrix{
  G\times_K S(V) \ar[r]^-\zeta \ar[d] & Y^{\alpha+n} \ar[d] \\
  G\times_K D(V) \ar[r]_-\xi & X
 }
\]
By Lemma~\ref{lem:vdiskComplex}, the
sphere $G\times_K S(V)$ is $G$-homotopy equivalent to
a $\delta$-$G$-CW($\alpha$) complex of dimension $\alpha+i-1$.
Now $\alpha+i-1 < \alpha+n$, so, by cellular approximation of maps,
the diagram above is homotopic to one in which $\zeta$ maps the sphere into $Y^{\alpha+n-1}$.
We can then find a lift of $\xi$ up to homotopy using
the inductive hypothesis if $i < n$ or the construction of $Y^{\alpha+n}$ if $i = n$.

Finally, $Y = \colim_n Y^{\alpha+n}$ satisfies the claim of the theorem.
\end{proof}

\begin{theorem}[Approximation by $\delta$-$G$-CW($\alpha$) Complexes]
Let $\delta$ be a familial dimension function and
let X be a $G$-space. Then there exists a $\delta$-$G$-CW($\alpha$)
complex $\Gamma X$ and a $\delta$-weak$_{\alpha}$ equivalence $g\colon \Gamma X\to X$.
If $f\colon X\to Y$ is a $G$-map and $g\colon \Gamma Y\to Y$ is an approximation
of $Y$ by a $\delta$-$G$-CW($\alpha$) complex, then
there exists a $G$-map $\Gamma f\colon \Gamma X\to \Gamma Y$, unique up to $G$-homotopy,
such that the following diagram commutes up to $G$-homotopy:
 \[
 \xymatrix{
  {\Gamma X} \ar[r]^{\Gamma f} \ar[d]_g & {\Gamma Y} \ar[d]^g \\
  X \ar[r]_{f} \ar[r] & Y
  }
 \]
 \end{theorem}

\begin{proof}
The existence of $g\colon \Gamma X\to X$ is the special case
of Theorem~\ref{thm:generalCellularApprox} in which we take $A = P = \emptyset$.
The existence and uniqueness of $\Gamma f$ follows from Whitehead's theorem,
which tells us that
$\pi G\K(\Gamma X, \Gamma Y) \iso \pi G\K(\Gamma X, Y)$.
 \end{proof}

\begin{theorem}[Approximation of Based Spaces]\label{thm:basedApproximation}
Let $\delta$ be a familial dimension function and
let X be a based $G$-space. Then there exists a based $\delta$-$G$-CW($\alpha$)
complex $\Gamma X$ and a based map $g\colon \Gamma X\to X$ that is a $\delta$-weak$_{\alpha}$ equivalence.
Further, $\Gamma$ is functorial up to based homotopy.
 \end{theorem}

\begin{proof}
The existence of $g\colon \Gamma X\to X$ is the special case
of Theorem~\ref{thm:generalCellularApprox} in which we take $A = P = *$.

Given $f\colon X\to Y$,
the existence and uniqueness of $\Gamma f\colon \Gamma X\to \Gamma Y$ follows from the relative
version of Whitehead's theorem,
which tells us that
$\pi G\K_*(\Gamma X, \Gamma Y) \iso \pi G\K_*(\Gamma X, Y)$.
 \end{proof}

\begin{theorem}[Approximation of Pairs]
Let $\delta$ be a familial dimension function and
let $(X,A)$ be a pair of $G$-spaces.
Then there exists a pair of $\delta$-$G$-CW($\alpha$) complexes $(\Gamma X, \Gamma A)$
and a pair of $\delta$-weak$_{\alpha}$ equivalences $g\colon (\Gamma X,\Gamma A)\to (X,A)$.
Further, $\Gamma$ is functorial on maps of pairs up to homotopy.
\end{theorem}

\begin{proof}
Take any approximation $g\colon \Gamma A\to A$, then
apply Theorem~\ref{thm:generalCellularApprox} to $\Gamma A\to X$ 
(taking $P = \emptyset$ in that theorem) to get $\Gamma X$ with $\Gamma A$ as a subcomplex.

Given $f\colon (X,A)\to (Y,B)$, we first construct
$\Gamma f\colon \Gamma A\to \Gamma B$ using the Whitehead theorem
and then extend to $\Gamma X\to \Gamma Y$
using the relative Whitehead theorem (considering the category of spaces under $\Gamma A$).
\end{proof}

Finally, we want to show that we can approximate excisive triads by
$\delta$-$G$-CW($\alpha$) triads. 
We need the following pasting result first.

\begin{proposition}\label{prop:genExcisiveTriads}
Let $\delta$ be a familial dimension function.
If $f\colon (X;A,B) \to (X';A',B')$ is a map of excisive triads such that
$f\colon A\intersect B\to A'\intersect B'$, $f\colon A\to A'$, and
$f\colon B\to B'$ are all $\delta$-weak$_{\alpha}$ equivalences,
then $f\colon X\to X'$ is a $\delta$-weak$_{\alpha}$ equivalence.
\end{proposition}

\begin{proof}
This follows from Theorem~\ref{thm:weakCharacterization}
and the corresponding nonequivariant result (see \cite[\S10.7]{May:concise})
applied to $f^H$ for each admissible $H$.
\end{proof}

\begin{theorem}[Approximation of Triads]\label{thm:dualTriadApprox}
Let $\delta$ be a familial dimension function and
let $(X;A,B)$ be an excisive triad.
Then there exists a $\delta$-$G$-CW($\alpha$) triad $(\Gamma X; \Gamma A, \Gamma B)$
and a map of triads 
\[
 g\colon (\Gamma X; \Gamma A, \Gamma B)\to (X;A,B)
\]
such that each of the maps $\Gamma A \intersect \Gamma B\to A\intersect B$,
$\Gamma A\to A$, $\Gamma B\to B$, and $\Gamma X\to X$ is
a $\delta$-weak$_{\alpha}$ equivalence.
$\Gamma$ is functorial on maps of excisive triads up to homotopy.
\end{theorem}

\begin{proof}
Let $C = A\intersect B$. Take a $\delta$-$G$-CW($\alpha$) approximation
$g\colon \Gamma C\to C$. Using Theorem~\ref{thm:generalCellularApprox}, extend $g$ to
approximations $g\colon (\Gamma A, \Gamma C)\to (A,C)$ and
$g\colon (\Gamma B, \Gamma C)\to (B,C)$.
Let $\Gamma X = \Gamma A \union_{\Gamma C} \Gamma B$.
All the statements of the theorem are clear except that the map
$g\colon \Gamma X\to X$ is a $\delta$-weak$_{\alpha}$ equivalence.
This follows from Proposition~\ref{prop:genExcisiveTriads}
using the argument in \cite[\S 10.7]{May:concise}:
We can not apply Proposition~\ref{prop:genExcisiveTriads} directly
to $g$ because $(\Gamma X;\Gamma A,\Gamma B)$ is not excisive. However, if
we replace $\Gamma X$ with the double mapping cylinder of the inclusions of
$\Gamma C$ in $\Gamma A$ and $\Gamma B$, we get an equivalent excisive triad
to which we can apply that result.
\end{proof}

In general, $\delta$-weak$_{\alpha}$ equivalence is weaker than weak $G$-equivalence
even if $\delta$ is complete.
This happens when the set of admissible subgroups
is smaller than the set of all subgroups.
Examples where $\alpha$ is not an actual representation are easy to construct.
The following gives an example with $\alpha$ being an actual representation.

\begin{example}\label{ex:notWeaklyEquivalent}
Let $G = SO(3)$ and let $H = SO(2) < G$.
We can identify $G/H$ with the two-sphere $S^2$, with the action
of $H$ being rotation around the $z$-axis, so
$\Lie(G/H)$ is a two-dimensional representation on which $SO(2)$
acts by rotation in the standard way.
Let $X$ be an $H$-space such that
\[
 X^K \hmtpc
  \begin{cases}
   S^0 & \text{if $K = H$}\\
   * & \text{if $K < H$.}
  \end{cases}
\]
Such a space exists by the construction of \cite{El:fixedpoints}.
Now consider the projection $f\colon G\times_H X\to G/H$.
Clearly, $f$ is not a weak $G$-equivalence, because $f^H$ is not a weak equivalence.
We claim that $f$ is an $\Lie$-weak$_{0}$ equivalence (i.e., we take $\alpha=0$).
The only fixed sets we need to check are by the subgroups of $H$.
If $K < H$, then $f^K$ is a weak equivalence. On the other hand,
$H$ is not $\Lie$-$0$-admissible because $\Lie(G/H)\not\subset n$ for any integer $n$,
so $f^H$ is not relevant. Thus, $f$ is an $\Lie$-weak$_{0}$ equivalence.
\end{example}

It's useful to note that the various notions of equivalence do coincide
if $\alpha$ is large enough.

\begin{theorem}\label{ref:thmStableGenWeakEquiv}
Let $\delta$ be a complete dimension function.
If $\alpha$ is stably equivalent to an actual representation of $G$ containing a copy of $\Lie(G)$, then a $G$-map $f\colon X\to Y$
is a $\delta$-weak$_{\alpha}$ equivalence if and only if it is a weak $G$-equivalence.
\end{theorem}

\begin{proof}
If $V$ contains a copy of $\Lie(G)$, then it contains a copy of $\Lie(G/H)$ for all $H$,
hence a copy of $\delta(G/H)$ for all $H$. Thus, every subgroup $H$
is admissible and the result follows from Theorem~\ref{thm:weakCharacterization}.
\end{proof}

\begin{corollary}\label{cor:equivalenceStability}
Let $\delta$ be a complete dimension function,
let $V$ be a representation of $G$ containing $\Lie(G)$,
let $f\colon X\to Y$ be a $\delta$-weak$_{V}$ equivalence of 
well-based $G$-spaces, and
let $W$ be a representation of $G$. Then
\[
 \susp^W f\colon \susp^W X \to \susp^W Y
\]
is a $\delta$-weak$_{V+W}$ equivalence.
\end{corollary}

\begin{proof}
Under the hypothesis that $V$ contains a copy of $\Lie(G)$, the preceding theorem
shows that both $\delta$-weak$_{V}$ equivalence and $\delta$-weak$_{V+W}$
equivalence coincide with weak $G$-equivalence.
That the suspension of a weak $G$-equivalence of well-based $G$-spaces is again
a weak $G$-equivalence follows on considering the fixed-point maps $f^H$ and
applying the corresponding nonequivariant result.
\end{proof}

The following example shows that this does not hold in general, if $V$ is too small.

\begin{example}\label{ex:dualInstability}
Recall Example~\ref{ex:notWeaklyEquivalent}.
As there, let $G = SO(3)$, $H = SO(2)$, and $V = 0$.
Let $f\colon G\times_H X \to G/H$ be the map constructed
there that is an $\Lie$-weak$_{0}$ equivalence but not a weak $G$-equivalence.
Let $W$ be $\Real^3$ with the usual action of $SO(3)$.
We claim that
\[
 \susp^W f_+ \colon \susp^W (G\times_H X)_+ \to \susp^W G/H_+
\]
is not an $\Lie$-weak$_{W}$ equivalence.
Now $\Lie(G/H)\subset W$ so $H$ is $\Lie$-$W$-admissible. 
Consider the fixed set map
\[
 (\susp^W f_+)^H\colon \susp^W (G\times_H X)_+ \to \susp^W G/H_+.
\]
$(G/H)^H$ consists of two points, while $(G\times_H X)^H$ is homotopy equivalent
to four points, mapping two-to-one to $(G/H)^H$.
It follows that $(\susp^W f_+)^H$ is not a weak equivalence,
hence that $\susp^W f_+$ is not an $\Lie$-weak$_{W}$ equivalence.
\end{example}

\section{Particular cell complexes}\label{sec:particularcomplexes}

We now briefly discuss some of the interesting specializations
of the general cell complexes we've defined.

\subsection{$G$-CW complexes}

Consider the case $\delta=0$ (on all subgroups) and $\alpha=0$: A $0$-$G$-CW($0$) complex
is just a classical $G$-CW complex as in \cite{Br:cohomology}, \cite{Ill:triangulation}, 
and \cite{Wan:thesis1}. A typical cell has the form
$G\times_H D^n = G/H\times D^n$.

A map $f$ is a $0$-$n$-equivalence
if and only if it is an equivariant $n$-equivalence, meaning
that each $f^H$ is a nonequivariant $n$-equivalence.
Similarly, $f$ is a $0$-weak$_{0}$ equivalence
if and only if it is a weak $G$-equivalence.

\subsection{$G$-CW($\alpha$) complexes}

Continuing with $\delta = 0$, we make the following definition.

\begin{definition}
Let $\alpha$ be a virtual representation.
An {\em (ordinary) $G$-CW($\alpha$) complex} is a
$0$-$G$-CW($\alpha$) complex.
An {\em (ordinary) $\alpha$-cell complex} is a
$0$-$\alpha$-cell complex.
\end{definition}

In particular, we have the case when $\alpha = V$ is an actual
representation, where
we recover the
notion of $G$-CW($V$) complex used in \cite{Wan:gcwv} and \cite[Ch. X]{May:alaska}.

Note that a subgroup $H$ is $0$-$\alpha$-admissible if and only if
$\alpha+n$ is stably equivalent to an
actual $H$-representation for some $n$.

\begin{definition}
A map is
an {\em (ordinary) $(\alpha+n)$-equivalence} if it is
a $0$-$(\alpha+n)$-equivalence.
A map is an 
{\em (ordinary) weak$_\alpha$ equivalence} if it is
a $0$-weak$_\alpha$ equivalence.
\end{definition}

From our characterization
of weak equivalence we get the following specialization.

\begin{corollary}
If $\alpha$ is a virtual representation,
then a map $f$ is a weak$_{\alpha}$ equivalence
if and only if $f^H$ is a weak equivalence for all
subgroups $H$ such that $\alpha+n$ is 
stably equivalent to an actual representation of $H$
for some $n$.
In particular, if $V$ is an actual representation,
$f$ is a weak$_V$ equivalence
if and only if it is a weak $G$-equivalence.
\qed
\end{corollary}

\begin{definition}
If $\alpha$ is a virtual representation, we say
that a map $f$ is an
{\em $|\alpha^*|$-equivalence} if, for each $H$
such that $\alpha$ is stably equivalent to an actual representation,
$f^H$ is an $|\alpha^H|$-equivalence.
\end{definition}

In particular, if $V$ is an actual representation,
we say that $f$ is a $|V^*|$-equivalence if
$f^H$ is a $|V^H|$-equivalence for each $H$.

Our characterization of connectivity gives us
the following, versions of which for $\alpha=V$
appear in \cite{Wan:gcwv} and \cite{Le:Hurewicz}.

\begin{corollary}
If $\alpha$ is a virtual representation,
then a map is an $\alpha$-equi\-va\-lence
if and only if it is an $|\alpha^*|$-equivalence.
\end{corollary}

\subsection{Dual $G$-CW($\alpha$) complexes}

Recall the following example. Let $M$ be a closed $V$-manifold,
i.e., a manifold in which, for each $x\in M$, the
tangent representation $T_xM$ at $x$ is isomorphic
to $V$ as a representation of $G_x$.
Then $M$ can be triangulated \cite{Ill:triangulation} and
the dual cell complex is an $\Lie$-$G$-CW($V$) complex
(where we take $\Lie$ to be defined on all subgroups).
With this in mind we make the following definition.

\begin{definition}
Let $\alpha$ be a virtual representation.
A {\em dual $G$-CW($\alpha$) complex} is an
$\Lie$-$G$-CW($\alpha$) complex.
A {\em dual $\alpha$-cell complex} is an
$\Lie$-$\alpha$-cell complex.
\end{definition}

Similarly, we have the notions of a
{\em dual $\alpha$-equivalence} and a
{\em dual weak$_\alpha$ equivalence}.
We can characterize these in terms of fixed sets,
but the general characterizations do not simplify
significantly in this context.

\subsection{Product complexes}\label{subsec:intermediate}

Ordinary complexes and dual complexes are the extreme cases
$\delta=0$ or $\delta=\Lie$ and are the most interesting.
However, other useful examples arise when we form products.

Consider a $\delta$-$H$-CW($\alpha$) complex $X$ and
an $\epsilon$-$K$-CW($\beta$) complex $Y$. Their product $X\times Y$
is a $G = H\times K$ space and has a cell structure in which
the cells have the form
\begin{multline*}
 G\times_{J\times L}(D(V-\delta(H/J)+m)\times D(W-\epsilon(K/L)+n)) \\
  \iso G\times_{J\times L}D(V+W-\delta(H/J)-\epsilon(K/L)+m+n)
\end{multline*}
with $V$ stably equivalent to $\alpha$ and $W$ stably equivalent to $\beta$.
We can therefore naturally view $X\times Y$ as a 
$(\delta\times\epsilon)$-$(H\times K)$-CW($\alpha+\beta$) complex
where $\delta\times\epsilon$ is the dimension function defined as follows.

\begin{definition}\label{def:productDimFcn}
If $\delta$ is a dimension function for $H$ with underlying collection
of subgroups $\F(\delta)$ and
$\epsilon$ is a dimension function for $K$ with underlying collection
of subgroups $\F(\epsilon)$,
let $\delta\times\epsilon$ be the dimension function for $H\times K$
defined on $\F(\delta)\times\F(\epsilon) = 
\{ J\times L \mid J\in\F(\delta) \text{ and } L\in\F(\epsilon) \}$ by
\[
 (\delta\times\epsilon)((H\times K)/(J\times L))
  = \delta(H/J)\dirsum\epsilon(K/L).
\]
\end{definition}

We leave to the reader to show that this is, indeed, a dimension function.

The dimension function $\delta\times\epsilon$ has the defect that it is almost
never complete or even familial, even if $\delta$ and $\epsilon$ are.
In practice, we'll often want to use a dimension function for $H\times K$ defined
on a family of subgroups.
For example, suppose that $X$ is an $H$-space 
with a $\delta$-$H$-CW($\alpha$) approximation
$\Gamma^\delta_\alpha X\to X$
and that $Y$ is a $K$-space with an $\epsilon$-$K$-CW($\beta$) approximation
$\Gamma^\epsilon_\beta Y \to Y$.
Then $\Gamma^\delta_\alpha X \times \Gamma^\epsilon_\beta Y$ is a
$(\delta\times\epsilon)$-$(H\times K)$-CW($\alpha+\beta$) complex, but
what sort of approximation is
$\Gamma^\delta_\alpha X \times \Gamma^\epsilon_\beta Y \to X\times Y$?
In the case that all subgroups in $\F(\delta)$ and $\F(\epsilon)$
are admissible (which we can achieve in the
based case after suitable suspension),
$\Gamma^\delta_\alpha X\to X$ and $\Gamma^\epsilon_\beta Y \to Y$
are weak $\F(\delta)$- and $\F(\epsilon)$-equivalences, respectively, so their product is
a weak $\overline{\F(\delta)\times\F(\epsilon)}$-equivalence. Thus, we get the following.

\begin{proposition}\label{prop:productApproximation}
Let $\delta$ be a familial dimension function for $H$ and
let $\epsilon$ be a familial dimension function for $K$.
Suppose that $\alpha$ is a virtual representation of $H$ such that
every subgroup in $\F(\delta)$ is $\delta$-$\alpha$-admissible and suppose that
$\beta$ is a virtual representation of $K$ such that every subgroup in $\F(\epsilon)$
is $\epsilon$-$\beta$-admissible.
If $X$ is an $H$-space with a $\delta$-$H$-CW$(\alpha)$ approximation
$\Gamma^\delta_\alpha X\to X$
and $Y$ is a $K$-space with an $\epsilon$-$K$-CW$(\beta)$ approximation
$\Gamma^\epsilon_\beta Y\to Y$,
then $\Gamma^\delta_\alpha X \times \Gamma^\epsilon_\beta Y$ is
a $(\delta\times\epsilon)$-$(H\times K)$-CW$(\alpha+\beta)$ complex and
$\Gamma^\delta_\alpha X \times \Gamma^\epsilon_\beta Y \to X\times Y$
is a weak $\overline{\F(\delta)\times\F(\epsilon)}$-equivalence.
\qed
\end{proposition}

\begin{corollary}\label{cor:productMapApproximation}
Let $X$ be an $H$-space, $Y$ a $K$-space, $Z$ an $(H\times K)$-space,
and $f\colon Z\to X\times Y$ an $(H\times K)$-map.
Let $\delta$ be a familial dimension function for $H$,
$\epsilon$ a familial dimension function for $K$,
and $\zeta$ a familial dimension function for $H\times K$ with
$\zeta\dimpred \delta\times\epsilon$.
Let $\alpha$ be a virtual representation of $H$ such that
every subgroup in $\F(\delta)$ is $\delta$-$\alpha$-admissible and let
$\beta$ be a virtual representation of $K$ such that every subgroup in $\F(\epsilon)$
is $\epsilon$-$\beta$-admissible.
Let $\Gamma^\delta_\alpha X\to X$ be a $\delta$-$H$-CW$(\alpha)$ approximation,
$\Gamma^\epsilon_\beta Y\to Y$ an $\epsilon$-$K$-CW$(\beta)$ approximation,
and $\Gamma^\zeta_{\alpha+\beta} Z\to Z$ a $\zeta$-$(H\times K)$-CW$(\alpha+\beta)$
approximation.
Then we may find a cellular map $\Gamma f$, unique up to cellular $(H\times K)$-homotopy, making
the following diagram homotopy commute:
\[
 \xymatrix@C+1em{
  \Gamma^\zeta_{\alpha+\beta} Z \ar@{-->}[r]^-{\Gamma f} \ar[d]
   & \Gamma^\delta_\alpha X \times \Gamma^\epsilon_\beta Y \ar[d] \\
  Z \ar[r] & X\times Y
 }
\]
\end{corollary}

\begin{proof}
The preceding proposition and the Whitehead theorem show that
there exists a unique homotopy class of maps making the diagram commute.
Cellular approximation of maps
(in the form of Theorem~\ref{thm:dualMapApprox})
shows that this homotopy class contains a cellular map and that
any two cellular maps in the class are cellularly homotopic.
\end{proof}

Of course, there are relative and based versions of these statements as well.
Corollary~\ref{cor:productMapApproximation} follows from the universal case
$Z = X\times Y$, but it's useful to think of the result in the more general form.

So, for particular $\delta$ and $\epsilon$ we look for a $\zeta\dimpred\delta\times\epsilon$.
For the approximation in the preceding corollary to be nontrivial,
we should look for as small
a $\zeta$ as possible.
The most interesting case is that of $G\times G$. Here, we're particularly interested
in the diagonal $G$-map $X\to X\times X$, which we view in adjunct form as
the $(G\times G)$-map
$(G\times G)\times_\Delta X \to X\times X$ where $\Delta\leq G\times G$ is the diagonal copy of $G$.

\begin{definition}\label{def:diagonalFamily}
Let $\F_\Delta$ denote the family of subgroups of $G\times G$ given by
\[
 \F_\Delta = \{ H\leq G\times G \mid \text{$H$ is conjugate to a subgroup of $\Delta$} \}.
\]
If $\delta$ is a complete dimension function for $G$, let $\delta_\Delta$ be the 
familial dimension function for $G\times G$, with underlying family $\F_\Delta$, given by
\[
 \delta_\Delta((G\times G)/H) = g^{-1}\cdot\delta(\Delta/H^g)
\]
where $H^g\leq \Delta$.
\end{definition}

Then $\delta_\Delta \dimpred \delta\times 0$ and also $\delta_\Delta\dimpred 0\times\delta$.
To see this, let $H\leq G$ and consider the diagonal copy of $H$, which we'll call
$\Delta_H\leq\Delta$ (the case of a conjugate of $\Delta_H$ is similar).
The smallest subgroup in $\F(\delta\times 0)$ containing $\Delta_H$ is $H\times H$.
By definition,
\[
 (\delta\times 0)((G\times G)/(H\times H)) \iso \delta(G/H)
\]
and
\[
 \delta_\Delta((G\times G)/\Delta_H) \iso \delta(G/H)
\]
as well. Hence, $\delta_\Delta \dimpred \delta\times 0$ as claimed,
and similarly $\delta_\Delta \dimpred 0\times\delta$.

Note also that $\delta_\Delta|\Delta = \delta$ as a dimension function on $G$.

Now, let $X$ be a $G$-space and consider a $0$-$G$-CW($\alpha$) approximation
$\Gamma^0_\alpha X$, a $\delta$-$G$-CW($\beta$) approximation
$\Gamma^\delta_\beta X$, and a $\delta$-$G$-CW($\alpha+\beta$) approximation
$\Gamma^\delta_{\alpha+\beta} X$.
Then, using Proposition~\ref{prop:genInduction} and the fact that
$(\delta_\Delta)((G\times G)/\Delta) = 0$, we have that
$(G\times G)\times_\Delta \Gamma^\delta_{\alpha+\beta} X$ is a
$\delta_\Delta$-$(G\times G)$-CW($\alpha+\beta$) complex.
Assuming all subgroups of $H$ and $K$ admissible,
Corollary~\ref{cor:productMapApproximation}
then gives us a cellular map, unique up to cellular homotopy, making the following diagram
commute and thus approximating the diagonal:
\[
 \xymatrix{
  (G\times G)\times_\Delta \Gamma^\delta_{\alpha+\beta} X \ar@{-->}[r] \ar[d]
   & \Gamma^0_\alpha X \times \Gamma^{\delta}_\beta X \ar[d] \\
  (G\times G)\times_\Delta X \ar[r] & X\times X
 }
\]

Here is a similar case that we'll need.
Suppose again that $\delta$ is a complete dimension function for $G$, with dual
$\Lie-\delta$.
Let $\Lie_\Delta$ be the result of applying Definition~\ref{def:diagonalFamily} to $\Lie$.
Then $\Lie_\Delta|\Delta = \Lie$ as a dimension function for $G$ and
$\Lie_\Delta \dimpred \delta\times(\Lie-\Delta)$. To see the latter relationship,
suppose that $H\leq G$ and $\Delta_H\leq \Delta$ is the corresponding subgroup of the diagonal.
Again, the smallest subgroup of $\F(\delta\times(\Lie-\delta))$ containing $\Delta_H$
is $H\times H$ and we have
\[
 (\delta\times(\Lie-\delta))((G\times G)/(H\times H))
  = \delta(G/H)\dirsum (\Lie-\delta)(G/H)
  \iso \Lie(G/H)
\]
and
\[
 \Lie_\Delta((G\times G)/\Delta_H) \iso \Lie(G/H).
\]
Let $X$ be a $G$-space and consider a $\delta$-$G$-CW($\alpha$) approximation
$\Gamma^\delta_\alpha X$, an $(\Lie-\delta)$-$G$-CW($\beta$) approximation
$\Gamma^{\Lie-\delta}_\beta X$, and an $\Lie$-$G$-CW($\alpha+\beta$) approximation
$\Gamma^\Lie_{\alpha+\beta} X$.
Then, using Proposition~\ref{prop:genInduction} and the fact that
$\Lie_\Delta((G\times G)/\Delta) = 0$, we have that
$(G\times G)\times_\Delta \Gamma^\Lie_{\alpha+\beta} X$ is an
$\Lie_\Delta$-$(G\times G)$-CW($\alpha+\beta$) complex.
Assuming all subgroups of $H$ and $K$ admissibile,
Corollary~\ref{cor:productMapApproximation}
then gives us a cellular map, unique up to cellular homotopy, making the following diagram
commute, so again approximating the diagonal:
\[
 \xymatrix{
  (G\times G)\times_\Delta \Gamma^\Lie_{\alpha+\beta} X \ar@{-->}[r] \ar[d]
   & \Gamma^\delta_\alpha X \times \Gamma^{\Lie-\delta}_\beta X \ar[d] \\
  (G\times G)\times_\Delta X \ar[r] & X\times X
 }
\]

\section{$G$-CW prespectra}\label{sec:CWprespectra}

In a subsequent section we will define cellular homology and cohomology of
$G$-CW complexes. We will then want to extend the definition to arbitrary $G$-spaces.
Our first thought is, given a based $G$-space $X$, to take a
$\delta$-$G$-CW($\alpha$) approximation $\Gamma_{\alpha} X\to X$ and then take
the homology of $\Gamma_\alpha X$ to be the homology of $X$.
However, the resulting theory may not have a suspension isomorphism.
The problem, as illustrated by Example~\ref{ex:dualInstability}, is that
suspension need not preserve $\delta$-weak$_\alpha$ equivalence, so that
$\susp^W\Gamma_\alpha X\to \susp^W X$ need not be a $\delta$-$G$-CW($\alpha+W$) approximation.
One way out of this problem is given by Corollary~\ref{cor:equivalenceStability}:
for a fixed $\delta$ and $\alpha$, we can choose a $V$ such that
$\delta$-weak$_{\alpha+V}$ equivalence is preserved by further suspension, and calculate
the homology  of $X$ as an appropriate shift of the homology of $\Gamma_{\alpha+V}\susp^V X$.

However, this approach does not generalize to the parametrized case we'll be interested
in later, so we take a different tack.
If $\Gamma_{\alpha+W}\susp^W X\to \susp^W X$ is a $\delta$-$G$-CW($\alpha+W$) approximation,
then the Whitehead theorem tells us that there is a map 
$\susp^W\Gamma_\alpha X\to \Gamma_{\alpha+W}\susp^W X$ over $\susp^W X$, unique up to homotopy.
In fact, using Theorem~\ref{thm:generalCellularApprox}, we can arrange that this map
is the inclusion of a subcomplex.
With the knowledge that this process eventually stabilizes, this discussion suggests taking a colimit
to define the chains of $X$.
It also suggests using a special, if old-fashioned, kind of prespectrum to formalize the process.

Recall from \cite{LMS:eqhomotopy} or \cite{MaySig:parametrized} that, if $\U$ is a $G$-universe, an {\em indexing sequence} $\V = \{V_i\}$ in $\U$ is
an expanding sequence
$V_1\subset V_2\subset\cdots$ of finite-dimensional subrepresentations of $\U$ such
that $\Union_i V_i = \U$.
We let $G\PreSpec{\V}{}$ denote the category of $G$-prespectra indexed on $\V$, i.e.,
prespectra $D$ defined by $G$-spaces $D(V_i)$ and structure maps
$\susp^{V_{i}-V_{i-1}} D(V_{i-1}) \to D(V_i)$.

\begin{definition}\label{def:CWprespectra}
Let $\V$ be an indexing sequence in a universe $\U$.
\begin{enumerate}
\item
A map $f\colon D\to E$ in $G\PreSpec{\V}{}$ is a {\em $\delta$-weak$_\alpha$ equivalence}
if, for each $i$, $f_i\colon D(V_i)\to E(V_i)$ is a $\delta$-weak$_{\alpha+V_i}$ equivalence
of $G$-spaces.

\item
A $G$-prespectrum $D$ in $G\PreSpec{\V}{}$ is a {\em $\delta$-$G$-CW$(\alpha)$ prespectrum}
if, for each $i$, $D(V_i)$ is a based $\delta$-$G$-CW$(\alpha+V_i)$ complex and each structure map
$\susp^{V_{i}-V_{i-1}} D(V_{i-1}) \to D(V_i)$
is the inclusion of a subcomplex.

\item
A map $D\to E$ of $\delta$-$G$-CW$(\alpha)$ prespectra is {\em the inclusion of a subcomplex}
if, for each $i$, the map $D(V_i)\to E(V_i)$ is the inclusion of a subcomplex.
We also say simply that $D$ is a subcomplex of $E$.

\item
A map $D\to E$ of $\delta$-$G$-CW$(\alpha)$ prespectra is {\em cellular} if, for each $i$,
the map $D(V_i)\to E(V_i)$ is cellular.

\item
If $D$ is a $G$-prespectrum in $G\PreSpec{\V}{}$, a {\em $\delta$-$G$-CW$(\alpha)$ approximation}
of $D$ is a $\delta$-$G$-CW$(\alpha)$ prespectrum $\Gamma^\delta_\alpha D$ and a
$\delta$-weak$_\alpha$ equivalence $\Gamma^\delta_\alpha D \to D$.
\end{enumerate}
\end{definition}

This notion of CW prespectrum harks back to the old idea of
CW spectra as in \cite{Ad:stablehomotopy}, \cite{Pup:stablehomotopy}, and \cite{Vo:Boardman}.
The notion of $\delta$-weak$_\alpha$ equivalence is similar to the level equivalences
discussed in \cite{MM:orthogonal} or \cite{MaySig:parametrized}.
We will not need or discuss its relationship to stable equivalence except for the following
brief comments:
When $\delta$ is complete, a $\delta$-weak$_\alpha$ equivalence $D\to E$
implies a weak $G$-equivalence $D(V_i)\to E(V_i)$ for sufficiently large $i$,
so is a stable equivalence,
but this simple argument is not available in the parametrized case, where the relationship
appears to be more complicated.
Further, we are not attempting to model the stable category and there are good reasons
that a CW model of that category is out of reach in the parametrized case, as discussed
in \cite[Chap.\ 24]{MaySig:parametrized}.

Our results on $G$-CW spaces give quick proofs of the following results.

\begin{lemma}[H.E.L.P.]
Let $r\colon E \to F$ be a $\delta$-weak$_\alpha$ equivalence of $G$-pre\-spec\-tra and
let $D$ be a $\delta$-$G$-CW$(\alpha)$ prespectrum with subcomplex $C$. 
If the following diagram commutes without the
dashed arrows, then there exist maps $\tilde g$ and $\tilde h$
making the diagram commute.
 \[
 \xymatrix{
  C \ar[rr]^-{i_0} \ar[dd] & & C\smsh I_+ \ar'[d][dd] \ar[dl]_(.6)h & & 
  C \ar[ll]_-{i_1} \ar[dl]_(.6)g \ar[dd]  \\
   & F & & E \ar[ll]_(.3)r \\
  D \ar[ur]^f \ar[rr]_-{i_0} & & D\smsh I_+ \ar@{-->}[ul]^(.6){\tilde h} & &
   D \ar[ll]^-{i_1} \ar@{-->}[ul]^(.6){\tilde g} \\
 }
 \]
\end{lemma}

\begin{proof}
We construct $\tilde g$ and $\tilde h$ inductively on the indexing space $V_i$.
For $i=1$ we simply quote the space-level H.E.L.P. lemma to find
$\tilde g_1$ and $\tilde h_1$.

For the inductive step, we assume that we've constructed $\tilde g_{i-1}$ and
$\tilde h_{i-1}$. We then apply the space-level H.E.L.P. lemma
with (using the notation of Lemma~\ref{lem:HELP}) $Y = E_i$, $Z = F_i$,
$X = D_i$, and $A = \susp X_{i-1}\union C_i$.
\end{proof}

\begin{proposition}[Whitehead]
Suppose that $D$ is a $\delta$-$G$-CW$(\alpha)$ prespectrum and that
$f\colon E\to F$ is a $\delta$-weak$_\alpha$ equivalence. Then
\[
 f_*\colon \pi G\PreSpec{\V}{}(D,E) \to \pi G\PreSpec{\V}{}(D,F)
\]
is an isomorphism, where $\pi G\PreSpec{\V}{}(-,-)$ denotes the group of
homotopy classes of $G$-maps.
Therefore, any $\delta$-weak$_\alpha$ equivalence of $\delta$-$G$-CW$(\alpha)$
prespectra is a $G$-homotopy equivalence.
\end{proposition}

\begin{proof}
We get surjectivity by applying the H.E.L.P. lemma to $D$ and its subcomplex $*$.
We get injectivity by applying it to $D\smsh I_+$ and its subcomplex
$D\smsh \bndry I_+$.
\end{proof}

\begin{proposition}[Cellular Approximation of Maps]
Suppose that $f\colon D\to E$ is a map of $\delta$-$G$-CW$(\alpha)$ prespectra,
$C$ is a subcomplex of $D$, and $h$ is a $G$-homotopy of $f|C$ to a cellular map.
Then $h$ can be extended to a $G$-homotopy of $f$ to a cellular map.
\end{proposition}

\begin{proof}
This follows by induction on the indexing space $V_i$. The first case to consider is
$f_1\colon D(V_1)\to E(V_1)$, and we know from the space-level result that 
we can extend $h_1$ to a $G$-homotopy $k_1$ from $f_1$ to a cellular map $g_1$.
For the inductive step, assume we have a homotopy $k_{i-1}$,
extending $h_{i-1}$, from $f_{i-1}$ to a cellular
map $g_{i-1}$. Then $\susp^{V_i-V_{i-1}}k_{i-1} \union h_i$ is a homotopy on the subcomplex
$\susp^{V_i-V_{i-1}}D(V_{i-1})\union C_i$ of $D(V_i)$. By the space-level result again, we can
extend to a homotopy $k_i$ on $D(V_i)$ from $f_i$ to a cellular map $g_i$.
\end{proof}

\begin{proposition}[Cellular Approximation of Prespectra]\label{prop:approxPrespectra}
If $D$ is a $G$-prespectrum in $G\PreSpec{\V}{}$, then there exists a
$\delta$-$G$-CW$(\alpha)$ approximation $\Gamma D\to D$.
If $f\colon D\to E$ is a map of $G$-prespectra and
$\Gamma E\to E$ is an approximation of $E$,
then there exists a cellular map
$\Gamma f\colon \Gamma D \to \Gamma E$,
unique up to cellular homotopy, making the following diagram homotopy commute:
\[
 \xymatrix{
  \Gamma D \ar[r]^-{\Gamma f} \ar[d]
   & \Gamma E \ar[d] \\
  D \ar[r]_-f & E
 }
\]
\end{proposition}

\begin{proof}
We construct $\Gamma D$ recursively on the indexing space $V_i$.
For $i=1$, we take $(\Gamma D)(V_1) = \Gamma(D(V_1))$ to be any
$\delta$-$G$-CW$(\alpha+V_1)$ approximation of $D(V_1)$.
Suppose that we have constructed $\Gamma D(V_{i-1}) \to D(V_{i-1})$, a
$\delta$-$G$-CW$(\alpha+V_{i-1})$ approximation.
Then $\susp^{V_i-V_{i-1}}\Gamma D(V_{i-1})$ is a $\delta$-$G$-CW$(\alpha+V_i)$ complex
and, by Theorem~\ref{thm:generalCellularApprox}, we can find a
$\delta$-$G$-CW$(\alpha+V_i)$ approximation
$\Gamma D(V_i) \to D(V_i)$ making the following diagram commute, in which the
map $\sigma$ at the top is the inclusion of a subcomplex:
\[
 \xymatrix{
  \susp^{V_i-V_{i-1}}\Gamma D(V_{i-1}) \ar[r]^-\sigma \ar[d]
   & \Gamma D(V_i) \ar[d] \\
  \susp^{V_i-V_{i-1}} D(V_{i-1}) \ar[r]_-\sigma & D(V_i)
 }
\]
Thus, $\Gamma D$ is a $\delta$-$G$-CW$(\alpha)$ prespectrum and the map
$\Gamma D\to D$ so constructed is a $\delta$-weak$_\alpha$ equivalence.

The existence and uniqueness of $\Gamma f$ follow from the Whitehead theorem
and cellular approximation of maps and homotopies.
\end{proof}

As usual, these results imply that we can invert the $\delta$-weak$_\alpha$ equivalences
of prespectra and that the result is equivalent to the ordinary homotopy category
of $\delta$-$G$-CW$(\alpha)$ prespectra.
We emphasize that this is not the stable category and our intention is not
to model the stable category---this notion of equivalence is much stricter
than stable equivalence. These results will be used solely as a mechanism for extending
the definition of ordinary homology and cohomology from $G$-CW complexes to arbitrary
$G$-spaces.

When discussing products we will need to use somewhat more general maps.

\begin{definition}
Let $D$ and $E$ be $G$-prespectra in $G\PreSpec{\V}{}$.
A {\em semistable map}
\[
 f = (M,\{f_i\})\colon D\to E
\]
consists of a nonnegative integer $M$ and,
for each $i\geq M$, maps $f_i\colon D(V_i)\to E(V_i)$
compatible in the usual way.
Say that two semistable maps
$f = (M,\{f_i\})$ and $g = (N,\{g_i\})$ from $D$ to $E$
are {\em equivalent} if there exists a $P\geq \max(M,N)$ such that
$f_i = g_i$ for all $i\geq P$.
Equivalence classes of semistable maps can be composed in the obvious way.
Let $G\PreSpec{\V}{s}$ denote the category
of $G$-prespectra indexed on $\V$ and equivalence classes of semistable maps.
\end{definition}

\section{Modules over preadditive categories}

Before we can construct the cellular homology and cohomology
theories based on general cell complexes, we need to discuss
some algebra.

Beginning with this section
we  start to make occasional use of equivariant spectra.
While \cite{LMS:eqhomotopy} is still a valuable exposition,
there are now a number of models of equivariant spectra available,
including the orthogonal spectra
of \cite{MM:orthogonal}, that have better formal properties.
For our purposes, it will usually not matter much which we use, because they
all give equivalent stable categories.
We borrow freely from notations used in \cite{LMS:eqhomotopy}.
Write $hG\Spec{}{}$ for the stable category of $G$-spectra, i.e.,
the category obtained by inverting stable equivalences in your favorite model.
Write $[E,F]_G = hG\Spec{}{}(E,F)$ for the group of stable $G$-maps between
$G$-spectra.

The {\em orbit category} $\orb{G}$ is the topological category
whose objects are the orbit spaces $G/H$ and whose morphisms
are the $G$-maps between them.
We give $\orb G(G/H,G/K)$ the compact-open topology; it is
homeomorphic to $(G/K)^H$.
Write $\sorb G$ for
the {\em stable orbit category}, the category of $G$-orbits and 
$G$-homotopy clases of stable
$G$-maps between them, i.e. the full subcategory of $hG\Spec{}{}$
on the objects $\susp^\infty_G G/H_+$ (see \cite{LMM:roghomology} and 
\cite[\S V.9]{LMS:eqhomotopy}). As in \cite{LMM:roghomology}
we take a {\em Mackey functor} to be an additive functor 
$\sorb G \to \Ab$ where $\Ab$ is the category of abelian groups.
(see also \cite{Dr:reps}, \cite{Le:Greenfunctors},
and \cite{GM:eqhomotopy}). 
If we want to specify the group $G$, we shall refer to a {\em $G$-Mackey functor}.
Mackey functors can be either covariant or
contravariant. When $G$ is finite the variance is essentially irrelevant because
$\sorb G$ is self-dual (\cite{Le:Greenfunctors}, \cite{LMS:eqhomotopy}). 
However, when $G$ is infinite, $\sorb G$ is not self-dual and the
variance becomes important. We adopt the convention of
writing a bar above or below to indicate variance: $\Mackey T$ will
denote a contravariant functor, and $\MackeyOp S$ will denote a covariant one.

More generally, we make the following definition.

\begin{definition}
If $\delta$ is a dimension function for $G$, let
$\sorb{G,\delta}$ denote the full subcategory of the stable category
on the objects $G_+\smsh_H S^{-\delta(G/H)}$
(with $H$ in the collection of subgroups associated with $\delta$).
To simplify notation, we write $(G/H,\delta)$ for the object
$G_+\smsh_H S^{-\delta(G/H)}$ in $\sorb{G,\delta}$.
We define a {\em $G$-$\delta$-Mackey functor}, or simply a
{\em $\delta$-Mackey functor}, to be an additive
functor $\sorb{G,\delta} \to \Ab$.
A $\delta$-Mackey functor can be either covariant or contravariant;
we use the convention that $\Mackey T$ denotes a contravariant
$\delta$-Mackey functor and $\MackeyOp S$ denotes a covariant one.
\end{definition}

Note the following.

\begin{proposition}\label{prop:orbitDuality}
If $\delta$ is a dimension function on $G$ and $\Lie-\delta$ is its dual,
then $\sorb{G,\delta} \iso \Op{\sorb{G,\Lie-\delta}}$.
\end{proposition}

\begin{proof}
From \cite[II.6.3 \& \S III.2]{LMS:eqhomotopy} we have that the
stable dual of an orbit $G/H$ is
\[
 D(G/H_+) \iso G_+\smsh_H S^{-\Lie(G/H)}.
\]
From this it follows that
\[
 D(G_+\smsh_H S^{-\delta(G/H)}) \iso G_+\smsh_H S^{-(\Lie(G/H)-\delta(G/H))},
\]
which gives the proposition.
\end{proof}

\begin{corollary}\label{cor:MackeyDuality}
Suppose $\delta$ is a dimension function on $G$.
Then the category of contravariant $\delta$-Mackey functors and
natural transformations between them is isomorphic to the category
of covariant $(\Lie-\delta)$-Mackey functors.
\qed
\end{corollary}

In particular, contravariant ordinary Mackey functors are the same
thing as covariant $\Lie$-Mackey functors, and similarly with the variances
reversed.

The definition of a Mackey functor is a special case of the following more general definition.
(See \cite{Mit:rings} or \cite{Mit:functorCats} for much more about these objects.)

\begin{definition}\label{def:indres}
Let $\A$ be a small preadditive category (that is, its hom sets
are abelian groups and composition is bilinear). 
Define an {\em $\A$-module}
to be an additive functor from $\A$ to $\Ab$, the category of abelian
groups. We adopt the convention of writing $\Mackey T$ for a contravariant $\A$-module
and $\MackeyOp S$ for a covariant $\A$-module. We define
 \[
 \Hom_{\A}(\Mackey T, \Mackey U) = \int\nolimits_{a\in\A}\Hom(\Mackey T(a), \Mackey U(a))
 \]
to be the group of natural transformations from $\Mackey T$ to $\Mackey U$, and
similarly for covariant modules. We define
 \[
 \Mackey T \tensor_{\A} \MackeyOp S 
 = \int\nolimits^{a\in\A} \Mackey T(a)\tensor\MackeyOp S(a).
\]
More explicitly,
\[
 \Mackey T \tensor_{\A} \MackeyOp S
  = \left[ \Dirsum_{a} \Mackey T(a)\tensor \MackeyOp S(a) \right] \Big/ \sim,
 \]
where the sum extends over all objects $a$ in $\A$ and
the equivalence relation is generated by
 \[
 (\alpha^* x)\tensor y \sim x\tensor (\alpha_* y)
 \]
when $\alpha\colon a\to a'$ is a map in $\A$, $x\in \Mackey T(a')$ and
$y\in \MackeyOp S(a)$.

Let $\B$ be another preadditive category (not necessarily small) and let
$F\colon\A\to\B$ be an additive functor. If $\Mackey U$ is a contravariant
$\B$-module and $\MackeyOp S$ is a covariant $\B$-module,
define
 \[
 F^*\Mackey U = \Mackey U\circ F
 \]
and
 \[
 F^*\MackeyOp S = \MackeyOp S\circ F.
 \]
If $\Mackey T$ is a contravariant $\A$-module, define the 
contravariant $\B$-module $F_!\Mackey T$ by
\[
 (F_!\Mackey T)(b) = \Mackey T \tensor_{\A} F^* \B(b,-)
  = \int\nolimits^{a\in\A} \Mackey T(a) \tensor \B(b,F(a)).
\]
Similarly, if $\MackeyOp R$ is a covariant $\A$-module, define the
covariant $\B$-module $F_!\MackeyOp R$ by
\[
 (F_!\MackeyOp R)(b) = F^*\B(-,b) \tensor_{\A} \MackeyOp R
   = \int\nolimits^{a\in\A} \B(F(a),b)\tensor \MackeyOp R(a).
\]
Finally, define
\[
 (F_*\Mackey T)(b) = \Hom_\A (F^*\B(-,b), \Mackey T)
\]
and
\[
 (F_*\MackeyOp R)(b) = \Hom_\A(F^*\B(b,-), \MackeyOp R).
\]
 \end{definition}

Very useful examples are the canonical projective modules:
If $a$ is an object of $\A$, let $\Mackey A_a$ be the $\A$-module defined by
$\Mackey A_a = \A(-,a)$, and let $\MackeyOp A^a$ be defined by
$\MackeyOp A^a = \A(a,-)$.

The following facts are standard and straightforward to prove.
If $\B$ is not small, claims of isomorphisms below like
$\Hom_{\B}(F_!\Mackey T, \Mackey U) \iso \Hom_{\A}(\Mackey T, F^*\Mackey U)$
are claims that the large limit or colimit involved on the left exists and is given by the
group on the right.

\begin{proposition}\label{prop:indresadjunction}
Let $\A$ and $\B$ be preadditive categories with $\A$ small, let
$F\colon \A\to\B$ be an additive functor, and let $\Mackey T$, etc.,
denote $\A$- or $\B$-modules of the appropriate variance below.
\begin{enumerate}
\item
For any object $a$ of $\A$ we have the following isomorphisms, 
natural in $a$ as well as in $\Mackey T$ or $\MackeyOp S$.
\begin{align*}
 \Hom_{\A}(\Mackey A_a, \Mackey T) &\iso \Mackey T(a) \\
 \Mackey T \tensor_{\A} \MackeyOp A^a &\iso \Mackey T(a) \\
 \Hom_{\A}(\MackeyOp A^a, \MackeyOp S) &\iso \MackeyOp S(a) \\
 \Mackey A_a \tensor_{\A} \MackeyOp S &\iso \MackeyOp S(a)
\end{align*}

\item
The functor $F_!$ is left adjoint to $F^*$. That is,
 \[
 \Hom_{\B}(F_!\Mackey T, \Mackey U) \iso \Hom_{\A}(\Mackey T, F^*\Mackey U)
 \]
and similarly for covariant modules.

\item
The functor $F_*$ is right adjoint to $F^*$. That is,
\[
 \Hom_\A(F^*\Mackey U, \Mackey T) \iso \Hom_\B(\Mackey U, F_*\Mackey T)
\]
and similarly for covariant modules.

\item
We have isomorphisms
 \[
 F_!\Mackey T \tensor_{\B} \MackeyOp S \iso \Mackey T\tensor_{\A} F^*\MackeyOp S
 \]
and
 \[
 \Mackey U \tensor_{\B} F_!\MackeyOp R \iso F^*\Mackey U \tensor_{\A} \MackeyOp R.
 \]
\item\label{item:indFree}
We have the following isomorphisms, natural in $a$.
 \[
 F_!\Mackey A_a \iso \Mackey A_{F(a)}
 \]
and
 \[
 F_!\MackeyOp A^a \iso \MackeyOp A^{F(a)}
 \]
\item
$F^*$ is exact, $F_!$ is right exact, and $F_*$ is left exact. \qed
 \end{enumerate}
 \end{proposition}

\begin{remark}\label{rem:zeroOrbit}
Recall that a {\em zero object} in a preadditive category 
$\A$ is an object that is both initial and terminal.
Such an object $z$ is characterized by the fact that $\A(z,z)$ is the zero group,
or by the fact that the identity map on $z$ is the zero element of $\A(z,z)$.
Any additive functor on $\A$ must take any zero object to a zero object in the target.
In particular, any module on $\A$ must take any zero object to the zero group.

At times it will be convenient to augment the category $\sorb{G,\delta}$
with a zero object and there is a natural way to do this:
include as an additional object the trivial spectrum $*$, the
zero object in the stable category.
Because a module must take $*$ to 0, modules on $\sorb{G,\delta}$
so augmented are equivalent to modules on the original $\sorb{G,\delta}$.
For most purposes it will not matter which version of $\sorb{G,\delta}$ we use;
we will point out the places where it is useful to include the zero object.
One important example will come while discussing restriction to fixed sets later in this chapter.
\end{remark}

Returning to the context of $\delta$-Mackey functors, we shall quite
often use the Mackey functors
\[
 \Mackey A_{G/K,\delta} = \sorb{G,\delta}(-,(G/K,\delta))
\]
and 
\[
 \MackeyOp A^{G/K,\delta} = \sorb{G,\delta}((G/K,\delta),-).
\]
We think of these as {\em free} Mackey functors because of the isomorphisms
 \[
  \Hom_{\sorb{G,\delta}}(\Mackey A_{G/K,\delta}, \Mackey T) \iso \Mackey T(G/K,\delta)
 \]
and
 \[
  \Hom_{\sorb{G,\delta}}(\MackeyOp A^{G/K,\delta}, \MackeyOp S) \iso \MackeyOp S(G/K,\delta)
 \]
given in the preceding proposition.
The isomorphism is given in each case by sending a homomorphism $f$ to
$f(1_{G/K,\delta})$.
Thus, we can think of $\{1_{G/K,\delta}\}$ as a basis for either $\Mackey A_{G/K,\delta}$
or $\MackeyOp A^{G/K,\delta}$.
In general, we shall use the term {\em free Mackey functor} for a
direct sum of Mackey functors of the form $\Mackey A_{G/K,\delta}$ or
$\MackeyOp A^{G/K,\delta}$.
If, for example, $\Mackey T = \Dirsum_i \Mackey A_{G/K_i,\delta}$, we think of the elements
$1_{G/K_i,\delta} \in \Mackey T(G/K_i,\delta)$ as forming a basis of $\Mackey T$.
As we shall see, the Mackey functors $\Mackey A_{G/K,\delta}$ and $\MackeyOp A^{G/K,\delta}$ 
play much the same central role in the theory
of equivariant homology and cohomology as does the group $\Z$ nonequivariantly.

We end this section with a general result we'll need later while discussing
fixed sets and restriction to subgroups.

\begin{proposition}\label{prop:commutingInducedMaps}
Let $\A$, $\B$, $\C$, and $\D$ be small additive categories and assume that we have
the following commutative diagram of additive functors:
\[
 \xymatrix{
  \A \ar[r]^\phi \ar[d]_\alpha & \B \ar[d]^\beta \\
  \C \ar[r]_\psi & \D
 }
\]
Then, as functors on contravariant $\C$-modules, there is a natural transformation
$\xi\colon \phi_!\alpha^* \to \beta^*\psi_!$,
which is an isomorphism if and only if the map
\[
 \int\nolimits^{a\in\A} \C(\alpha(a),c)\tensor \B(b,\phi(a)) \to \D(\beta(b),\psi(c))
\]
is an isomorphism for all $b\in \B$ and $c\in \C$.
For covariant $\C$-modules we have the dual condition: $\xi$ is an isomorphism if and only if
\[
 \int\nolimits^{a\in\A} \B(\phi(a),b)\tensor \C(c,\alpha(a)) \to \D(\psi(c),\beta(b))
\]
is an isomorphism for all $b\in\B$ and $c\in\C$.
\end{proposition}

\begin{proof}
The map $\xi$ is the adjoint of the map
\[
 \alpha^* \xrightarrow{\eta} \alpha^*\psi^*\psi_! \iso \phi^*\beta^*\psi_!;
\]
it is also the adjoint of
\[
 \beta_!\phi_!\alpha^* \iso \psi_!\alpha_!\alpha^* \xrightarrow{\epsilon} \psi_!,
\]
as can be seen from an appropriate diagram. Explicitly, if $\Mackey T$ is a $\C$-module,
\[
 \phi_!\alpha^*\Mackey T =
  \int\nolimits^{a\in\A}\Mackey T(\alpha(a))\tensor\B(-,\phi(a))
\]
while
\[
 \beta^*\psi_!\Mackey T =
  \int\nolimits^{c\in \C}\Mackey T(c)\tensor\D(\beta(-),\psi(c))
\]
and the map $\phi_!\alpha^*\Mackey T \to \beta^*\psi_!\Mackey T$ is given by
$t\tensor f \mapsto t\tensor\beta(f)$ for $t\in\Mackey T(\alpha(a))$, 
setting $c = \alpha(a)$ so $\psi(c) = \beta\phi(a)$.

Now, if $\xi$ is a natural isomorphism, we can look at the special case
$\Mackey T = \C(-,c)$. From above, we have
\[
 (\phi_!\alpha^*\C(-,c))(b) =
  \int\nolimits^{a\in\A}\Mackey \C(\alpha(a),c)\tensor\B(b,\phi(a))
\]
for each object $b\in\B$, while
\[
 (\beta^*\psi_!\C(-,c))(b) =
  \int\nolimits^{c'\in \C}\C(c',c)\tensor\D(\beta(b),\psi(c'))
  \iso \D(\beta(b),\psi(c)).
\]
Therefore,
\[
 \int\nolimits^{a\in\A} \C(\alpha(a),c)\tensor\B(b,\phi(a))
  \iso \D(\beta(b),\psi(c)).
\]

For the converse, suppose that this isomorphism holds. Then, for any $\Mackey T$,
the following shows that $\xi$ is an isomorphism:
\begin{align*}
  (\phi_!\alpha^*\Mackey T)(b) 
    &= \int\nolimits^{a\in\A}\Mackey T(\alpha(a))\tensor\B(b,\phi(a)) \\
    &\iso \int\nolimits^{a\in\A}
       \int\nolimits^{c\in\C}\Mackey T(c)\tensor\C(\alpha(a),c)\tensor\B(b,\phi(a)) \\
    &\iso \int\nolimits^{c\in\C} \Mackey T(c) \tensor
       \int\nolimits^{a\in\A}\C(\alpha(a),c)\tensor\B(b,\phi(a)) \\
    &\iso \int\nolimits^{c\in\C} \Mackey T(c) \tensor \D(\beta(b),\psi(c)) \\
    &= (\beta^*\psi_!\Mackey T)(b).
\end{align*}

The proof for covariant modules is similar or follows by duality.
\end{proof}

\begin{remark}
An appealing way of formulating this proof is to write
\[
 \phi_!\alpha^*\Mackey T = \Mackey T\tensor_\A \B
\]
and
\[
 \beta^*\psi_!\Mackey T = \Mackey T\tensor_\C \D.
\]
Using $\Mackey T \iso \Mackey T\tensor_\C \C$, the map $\xi$ can be written as
\[
 \Mackey T\tensor_\A \B \iso (\Mackey T\tensor_\C \C)\tensor_\A \B
  \iso \Mackey T\tensor_\C (\C\tensor_\A \B)
  \to \Mackey T\tensor_\C \D,
\]
which is an isomorphism for all $\Mackey T$ if and only if $\C\tensor_\A \B \to \D$ is an isomorphism.
When you make sense of this notation you get the explicit description given in the proof.
\end{remark}

\section{Cellular homology and cohomology}\label{sec:cellhomology}

We now turn to the construction of cellular equivariant
homology and cohomology graded on $RO(G)$. 

Let $\delta$ be a dimension function for $G$ and let
$\alpha = V\ominus W$ be a virtual representation of $G$.
Let $(X,A)$ be a relative $\delta$-$G$-CW($\alpha$) complex with skeleta $X^{\alpha+n}$. Then
 \[
 X^{\alpha+n}/X^{\alpha+n-1} = \Wedge_i G_+\smsh_{H_i} S^{V_i},
 \]
where the wedge runs over the $(\alpha+n)$-cells of $X$ and
$V_i$ is stably equivalent to $\alpha-\delta(G/H_i)+n$.
The suspension spectrum is thus
\begin{align*}
 \susp_G^\infty X^{\alpha+n}/X^{\alpha+n-1} 
   &\hmtpc \Wedge_i \susp_G^\infty G_+\smsh_{H_i} S^{V} \\
   &\hmtpc \Wedge_i \susp_G^{\alpha+n}(G_+\smsh_{H_i}S^{-\delta(G/H_i)}).
\end{align*}

To avoid sign ambiguities we need to be careful and rather precise about definitions.
We choose a specific sequence of trivial representations
\[
 \Real\subset \Real^2\subset \cdots \subset \Real^n \subset \cdots
\]
disjoint from $V$ and $W$,
with $\Real^{n+1} = \Real^n \dirsum \Real$ specifying how each sits inside the next.
We choose another sequence
\[
 \tilde\Real\subset \tilde\Real^2\subset \cdots \subset \tilde\Real^n \subset \cdots
\]
disjoint from $V$ and $W$,
with $\tilde\Real^{n+1} = \Real\dirsum\tilde\Real^n$.
Below, when we write $S^n$ we mean the one-point compactification of this particular choice of
$\Real^n$, so that
we have the identification $S^{n+1} = S^n\smsh S^1$.
We write $\tilde S^n$ for the one-point compactification of $\tilde \Real^n$, so that
we have $\tilde S^{n+1} = S^1\smsh \tilde S^n$.

\begin{definition}\label{def:vchains}
Let $(X,A)$ be a relative $\delta$-$G$-CW($\alpha$) complex, where $\alpha = V \ominus W$. 
The {\em cellular chain complex of $(X,A)$}, 
$\Mackey {C}^{G,\delta}_{\alpha+*} (X,A) = \Mackey C_{\alpha+*}(X,A)$, is the differential graded contravariant $\delta$-Mackey functor
specified on orbits by
\begin{multline*}
 \Mackey {C}_{\alpha+n} (X,A)(G/H,\delta) \\
  = \begin{cases}
    [ G_+\smsh_H S^{-\delta(G/H)}\smsh S^V\smsh S^n, \susp_G^\infty X\sp {\alpha+n}/X\sp {\alpha+n-1}\smsh S^W ]_G & \text{if $n \geq 0$} \\
    [ G_+\smsh_H S^{-\delta(G/H)}\smsh S^V, \susp_G^\infty X\sp {\alpha+n}/X\sp {\alpha+n-1}\smsh S^W\smsh\tilde S^{-n} ]_G & \text{if $n < 0$,}
    \end{cases}
\end{multline*}
where $[-,-]_G$ denotes the group of stable $G$-maps.
The differential $d$ is given by the $G$-map 
$X\sp {\alpha+n}/X\sp {\alpha+n-1}\to \susp X^{\alpha+n-1}_+\to \susp X\sp{\alpha+n-1}/X\sp {\alpha+n-2}$
as follows:
If $n>0$ we take $d$ to be the composite
\begin{align*}
\Mackey C_{\alpha+n}&(X,A) \\
  &= [ G_+\smsh_H S^{-\delta(G/H)}\smsh S^V\smsh S^n, \susp_G^\infty X\sp {\alpha+n}/X\sp {\alpha+n-1}\smsh S^W ]_G\\
  &\xrightarrow{} [ G_+\smsh_H S^{-\delta(G/H)}\smsh S^V\smsh S^n, \susp_G^\infty X\sp {\alpha+n-1}/X\sp {\alpha+n-2}\smsh S^1\smsh S^W ]_G \\
  &\xrightarrow{} [ G_+\smsh_H S^{-\delta(G/H)}\smsh S^V\smsh S^n, \susp_G^\infty X\sp {\alpha+n-1}/X\sp {\alpha+n-2}\smsh S^W\smsh S^1 ]_G \\
  &\xrightarrow{\susp^{-1}} [ G_+\smsh_H S^{-\delta(G/H)}\smsh S^V\smsh S^{n-1}, \susp_G^\infty X\sp {\alpha+n-1}/X\sp {\alpha+n-2}\smsh S^W ]_G \\
  &= \Mackey C_{\alpha+n-1}(X,A).
\end{align*}
If $n \leq 0$, $d$ is the composite
\begin{align*}
\Mackey C_{\alpha+n}&(X,A) \\
  &= [ G_+\smsh_H S^{-\delta(G/H)}\smsh S^V, \susp_G^\infty X\sp {\alpha+n}/X\sp {\alpha+n-1}\smsh S^W\smsh\tilde S^{-n} ]_G\\
  &\xrightarrow{} [ G_+\smsh_H S^{-\delta(G/H)}\smsh S^V, \susp_G^\infty X\sp {\alpha+n-1}/X\sp {\alpha+n-2}\smsh S^1\smsh S^W\smsh\tilde S^{-n} ]_G \\
  &\xrightarrow{} [ G_+\smsh_H S^{-\delta(G/H)}\smsh S^V, \susp_G^\infty X\sp {\alpha+n-1}/X\sp {\alpha+n-2}\smsh S^W\smsh S^1\smsh\tilde S^{-n} ]_G \\
  &= [ G_+\smsh_H S^{-\delta(G/H)}\smsh S^V, \susp_G^\infty X\sp {\alpha+n-1}/X\sp {\alpha+n-2}\smsh S^W\smsh\tilde S^{-(n-1)} ]_G \\
  &= \Mackey C_{\alpha+n-1}(X,A).
\end{align*}
\end{definition}

For simplicity of notation, we will write
\[
 \Mackey C_{\alpha+n}(X,A)(G/H,\delta) =
 [G_+\smsh_H S^{-\delta(G/H)+\alpha+n}, \susp_G^\infty X^{\alpha+n}/X^{\alpha+n-1}]_G,
\]
but it is important to keep in mind that this is shorthand for the precise definition above.

\begin{remarks}\label{rem:chains}
\begin{enumerate}\item[]
\item
If $f\colon (X,A)\to (Y,B)$ is a cellular map of relative $\delta$-$G$-CW($\alpha$) complexes,
the maps $X^{\alpha+n}\to Y^{\alpha+n}$ induce a map of chain complexes
\[
 f_*\colon \Mackey C_{\alpha+*}(X,A)\to \Mackey C_{\alpha+*}(Y,B).
\]
This makes $\Mackey C_{\alpha+*}(X,A)$ a covariant functor on the category of
relative $\delta$-$G$-CW($\alpha$) complexes and cellular maps.

\item\label{item:alphaFunctoriality}
If $\zeta\colon \alpha\to \alpha'$ is a virtual map,
any relative $\delta$-$G$-CW($\alpha$) structure on $(X,A)$ is also a relative $\delta$-$G$-CW($\alpha'$) structure, and vice versa.
The induced stable map $S^\alpha\to S^{\alpha'}$
(i.e., the pair of homeomorphisms $S^{V\dirsum Z_1} \to S^{V'\dirsum Z_2}$ and
$S^{W\dirsum Z_1} \to S^{W'\dirsum Z_2}$ given by $\zeta$)
then gives a chain isomorphism
\[
 \zeta^*\colon \Mackey C_{\alpha'+*}(X,A)\to \Mackey C_{\alpha+*}(X,A).
\]
This makes $\Mackey C_{\alpha+*}(X,A)$ a contravariant functor on the category of
virtual representations of $G$ equivalent to a specified $\alpha$.
\item\label{item:chainSuspension}
A special case of the preceding remark is the following:
If $Z$ is a representation of $G$, suspension by $Z$ induces a chain isomorphism
\[
 \Mackey C_{V\ominus W + *}(X,A) \iso
 \Mackey C_{(V\dirsum Z)\ominus(Z\dirsum W)+*}(X,A).
\]
A reinterpretation of the latter chain complex in the based case gives us a suspension isomorphism
\[
 \sigma^Z\colon \Mackey C_{V\ominus W + *}(X,*) \iso
 \Mackey C_{(V\dirsum Z)\ominus W +*}(\susp^Z X, *)
\]
or
\[
 \sigma^Z\colon \Mackey C_{\alpha + *}(X,*) \iso
 \Mackey C_{\alpha + Z +*}(\susp^Z X, *).
\]
Note that we use the fact that, if $X$ is a based $\delta$-$G$-CW($\alpha$) complex,
then $\susp^Z X$ is a based $\delta$-$G$-CW($\alpha+Z$) complex with
$(\susp^Z X)^{\alpha+Z+n} = \susp^Z X^{\alpha+n}$.

\item
The functor $\Mackey C_{V\ominus W + n}(X)$ should really be thought of as a functor
of the triple $(V,W,n)$. Our notation, although convenient and suggestive, tends to hide
this fact. On the other hand, if we're careful, we see that we can shift trivial summands
from $V$ or $W$ to $n$ or back. For example, for $n\geq 0$, we can identify
$\Mackey C_{(\alpha+1)+n}(X,*)$ with $\Mackey C_{\alpha+(n+1)}(X,*)$ via the isomorphism
\begin{multline*}
 [G_+\smsh_H S^{-\delta(G/H)}\smsh S^{V+1}\smsh S^n,
    \susp_G^\infty X^{(\alpha+1)+n}/X^{(\alpha+1)+n-1}\smsh S^W]_G \\
 \iso [G_+\smsh_H S^{-\delta(G/H)}\smsh S^{V}\smsh S^{1+n},
    \susp_G^\infty X^{\alpha+n+1}/X^{\alpha+n}\smsh S^W]_G
\end{multline*}
where $S^{1+n} = S^1\smsh S^n$. Taking the smash product in this order makes
the identification a chain map. Similarly, if $n<0$, we identify
\begin{multline*}
 [G_+\smsh_H S^{-\delta(G/H)}\smsh S^{V+1},
    \susp_G^\infty X^{(\alpha+1)+n}/X^{(\alpha+1)+n-1}\smsh S^W\smsh S^{-n}]_G \\
 \iso [G_+\smsh_H S^{-\delta(G/H)}\smsh S^{V},
    \susp_G^\infty X^{\alpha+n+1}/X^{\alpha+n}\smsh S^W\smsh S^{-n-1}]_G
\end{multline*}
where $S^{-n-1}\smsh S^1 = S^{-n}$.

\end{enumerate}
\end{remarks}

$\Mackey C_{\alpha+n}(X,A)$ is a free Mackey functor. Precisely, if
 \[
 \susp_G^\infty X^{\alpha+n}/X^{\alpha+n-1} 
   \iso \Wedge_i \susp_G^{\alpha+n}(G_+\smsh_{H_i}S^{-\delta(G/H_i)}).
 \]
as above, then
 \[
 \Mackey C_{\alpha+n}(X,A) \iso \Dirsum_i \Mackey A_{G/H_i,\delta}.
 \]
A basis is given by the set of inclusions
 \[
  \susp_G^{\alpha+n}(G_+\smsh_{H_j}S^{-\delta(G/H_j)}) \includesin \Wedge_i \susp_G^{\alpha+n}(G_+\smsh_{H_i}S^{-\delta(G/H_i)}).
 \]
It follows that, if $\MackeyOp S$ is a covariant $\delta$-Mackey functor and
$\Mackey T$ is a contravariant one, then
 \[
 \Hom_{\sorb{G,\delta}}(\Mackey{C}_{\alpha+n}(X,A), \Mackey T)
  = \Dirsum_i \Mackey T(G/H_i,\delta)
 \]
and
 \[
 \Mackey{C}_{\alpha+n}(X,A) \tensor_{\sorb {G,\delta}} \MackeyOp S 
  = \Dirsum_i \MackeyOp S(G/H_i,\delta).
 \]
In both
cases the induced differential can be understood as coming from the attaching
maps of the cells.

\begin{definition}\label{def:alphaHomology}
 Let $\Mackey T$ be a contravariant $\delta$-Mackey functor and let  
$\MackeyOp S$ be a covariant $\delta$-Mackey functor.
\begin{enumerate}
\item
Let $(X,A)$ be a relative $\delta$-$G$-CW($\alpha$) complex.
We define the {\em $(\alpha+n)$th cellular homology} of
$(X,A)$, with coefficients in $\MackeyOp S$, to be
 \[
 \CH^{G,\delta}_{\alpha+n}(X,A;\MackeyOp S) =
  H_{\alpha+n}( \Mackey{C}^{G,\delta}_{\alpha+*}(X,A) \tensor_{\sorb{G,\delta}} \MackeyOp S).
 \]
 and we define the {\em $(\alpha+n)$th cellular cohomology} of $(X,A)$, with
coefficients in $\Mackey T$, to be
 \[
 \CH_{G,\delta}^{\alpha+n} (X,A; \Mackey T) = 
  H^{\alpha+n}(\Hom_{\sorb{G,\delta}}( \Mackey{C}^{G,\delta}_{\alpha+*}(X,A), \Mackey T)).
 \]
where we introduce a sign in the differential:
$(da)(x) = (-1)^{n+1}a(dx)$ if 
$a\in \Hom_{\sorb{G,\delta}}( \Mackey{C}^{G,\delta}_{\alpha+n}(X,A), \Mackey T)$.
(The sign is necessary to make evaluation be a chain map.)
Homology is covariant in cellular maps of $(X,A)$
while cohomology is contravariant in cellular maps of $(X,A)$.

\item
If $X$ is a $\delta$-$G$-CW($\alpha$) complex, so $(X,\emptyset)$ is a relative
$\delta$-$G$-CW($\alpha$) complex, we define
\[
 \CH^{G,\delta}_{\alpha+n}(X; \MackeyOp S) = \CH^{G,\delta}_{\alpha+n}(X,\emptyset; \MackeyOp S)
\]
and
\[
 \CH_{G,\delta}^{\alpha+n}(X; \Mackey T) = \CH_{G,\delta}^{\alpha+n}(X,\emptyset; \Mackey T).
\]

\item
If $X$ is a based $\delta$-$G$-CW($\alpha$) complex, so $(X,*)$ is a relative
$\delta$-$G$-CW($\alpha$) complex, we define the {\em reduced} homology and cohomology
of $X$ to be
\[
 \tCH^{G,\delta}_{\alpha+n}(X; \MackeyOp S) = \CH^{G,\delta}_{\alpha+n}(X,*; \MackeyOp S)
\]
and
\[
 \tCH_{G,\delta}^{\alpha+n}(X; \Mackey T) = \CH_{G,\delta}^{\alpha+n}(X,*; \Mackey T).
\]
\end{enumerate}
\end{definition}

Note: We use the notation $\CH$ to distinguish the theories defined on $G$-CW complexes and cellular
maps from the theories we'll discuss shortly defined on all $G$-spaces. Among other differences,
the theories on $G$-CW complexes can be defined for any dimension function $\delta$
while the theories defined on $G$-spaces will require that $\delta$ be familial.

For the moment we concentrate on the reduced theories.
In order to have sufficient letters to use, we no longer reserve $V$ and $W$ for the
representations defining $\alpha$.

\begin{theorem}[Reduced Homology and Cohomology of Complexes]\label{thm:reducedHomologyComplexes}
Let $\delta$ be a dimension function for $G$, let
$\alpha$ be a virtual representation of $G$, and
let $\MackeyOp S$ and $\Mackey T$ be respectively a covariant and a contravariant
$\delta$-Mackey functor. Then the abelian groups
$\tCH^{G,\delta}_\alpha(X;\MackeyOp S)$ and $\tCH_{G,\delta}^\alpha(X;\Mackey T)$ are
respectively covariant and contravariant
functors on the homotopy category of 
based $\delta$-$G$-CW($\alpha$) complexes and cellular maps and homotopies.
They are also respectively contravariant and covariant functors 
of $\alpha$.
These functors satisfy the following properties.
\begin{enumerate}

\item (Exactness)
If $A$ is a based subcomplex of $X$, then the following sequences are exact:
\[
 \tCH^{G,\delta}_\alpha(A; \MackeyOp S) \to \tCH^{G,\delta}_\alpha(X; \MackeyOp S)
  \to \tCH^{G,\delta}_\alpha(X/A; \MackeyOp S)
\]
and
\[
 \tCH_{G,\delta}^\alpha(X/A; \Mackey T) \to \tCH_{G,\delta}^\alpha(X; \Mackey T)
  \to \tCH_{G,\delta}^\alpha(A; \Mackey T).
\]

\item (Additivity)
If $X = \Wedge_i X_i$ is a wedge of based $\delta$-$G$-CW($\alpha$) complexes,
then the inclusions of the wedge summands induce isomorphisms
\[
 \Dirsum_i\tCH^{G,\delta}_\alpha(X_i; \MackeyOp S) \iso \tCH^{G,\delta}_\alpha(X; \MackeyOp S)
\]
and
\[
 \tCH_{G,\delta}^\alpha(X; \Mackey T) \iso \prod_i \tCH_{G,\delta}^\alpha(X_i; \Mackey T).
\]

\item (Suspension)
There are suspension isomorphisms
 \[
  \sigma^V\colon \tCH^{G,\delta}_\alpha(X; \MackeyOp S) 
     \xrightarrow{\iso} \tCH^{G,\delta}_{\alpha + V}(\susp^V X; \MackeyOp S)
 \]
and
 \[
 \sigma^V\colon \tCH_{G,\delta}^\alpha(X; \Mackey T) 
     \xrightarrow{\iso} \tCH_{G,\delta}^{\alpha + V}(\susp^V X; \Mackey T).
 \]
These isomorphisms satisfy $\sigma^0 = \id$, $\sigma^{W}\circ \sigma^V = \sigma^{V\oplus W}$,
and the following naturality condition: If $\zeta\colon V\to V'$ is an isomorphism,
then the following diagrams commute (cf \cite[XIII.1.1]{May:alaska}):
\[
 \xymatrix@C+2em{
  \tCH^{G,\delta}_\alpha(X; \MackeyOp S) \ar[r]^-{\sigma^V} \ar[d]_{\sigma^{V'}}
    & \tCH^{G,\delta}_{\alpha + V}(\susp^{V}X; \MackeyOp S)
          \ar[d]^{\tCH_{\id}(\id\smsh\zeta)} \\
  \tCH^{G,\delta}_{\alpha + V'}(\susp^{V'}X; \MackeyOp S)
       \ar[r]_-{\tCH_{\id+\zeta}(\id)}
    & \tCH^{G,\delta}_{\alpha + V}(\susp^{V'}X; \MackeyOp S)
 }
\]
and
\[
 \xymatrix@C+2em{
  \tCH_{G,\delta}^\alpha(X; \Mackey T) \ar[r]^-{\sigma^V} \ar[d]_{\sigma^{V'}}
    & \tCH_{G,\delta}^{\alpha + V}(\susp^{V}X; \Mackey T)
          \ar[d]^{\tCH^{\id+\zeta}(\id)} \\
  \tCH_{G,\delta}^{\alpha + V'}(\susp^{V'}X; \Mackey T)
       \ar[r]_-{\tCH^{\id}(\id\smsh\zeta)}
    & \tCH_{G,\delta}^{\alpha + V'}(\susp^{V}X; \Mackey T)
 }
\]

\item (Dimension Axiom)
If $H\in \F(\delta)$ and $V$ is a representation of $H$ so large that
$V-\delta(G/H)+n$ is an actual representation, then
there are natural isomorphisms
\[
 \tCH^{G,\delta}_{V+k}(G_+\smsh_H S^{V-\delta(G/H)+n};\MackeyOp S)
  \iso
   \begin{cases}
    \MackeyOp S(G/H,\delta) & \text{if $k = n$} \\
    0 & \text{if $k\neq n$}
   \end{cases}
\]
and
\[
 \tCH_{G,\delta}^{V+k}(G_+\smsh_H S^{V-\delta(G/H)+n};\Mackey T)
  \iso
   \begin{cases}
    \Mackey T(G/H,\delta) & \text{if $k = n$} \\
    0 & \text{if $k\neq n$.}
   \end{cases}
\]

\end{enumerate}
\end{theorem}

\begin{proof}
That $\tCH^{G,\delta}_\alpha(X; \MackeyOp S)$ and $\tCH_{G,\delta}^\alpha(X; \Mackey T)$
are functorial on cellular maps and invariant under cellular homotopies is standard.
We also need to show that they
are functorial in $\alpha$ as stated in the theorem.
Let $\zeta\colon \alpha\to \alpha'$ be a virtual map.
As noted in Remark~\ref{rem:chains}(\ref{item:alphaFunctoriality}),
there is an induced chain isomorphism
$\zeta^*\colon\Mackey C_{\alpha'+*}(X,*) \to \Mackey C_{\alpha+*}(X,*)$.
This chain map induces the maps in homology and cohomology with
the variance claimed.

{\em Exactness:}
If $A$ is a based subcomplex of $X$, then $\Mackey C_{\alpha+*}(A,*)$ is a sub-chain
complex of $\Mackey C_{\alpha+*}(X,*)$ and
$\Mackey C_{\alpha+*}(X/A,*) \iso \Mackey C_{\alpha+*}(X,*)/\Mackey C_{\alpha+*}(A,*)$.
The exactness of the homology and cohomology sequences follows from familiar
homological algebra.

{\em Additivity:}
If $X = \Wedge_i X_i$, then $X^{\alpha+n} = \Wedge_i X_i^{\alpha+n}$, hence
\[
 \Mackey C_{\alpha+*}(X,*) \iso \Dirsum_i\Mackey C_{\alpha+*}(X_i,*).
\]
The additivity of homology and cohomology follows from this.

{\em Suspension:}
The suspension isomorphisms in homology and cohomology are induced by the
chain isomorphism $\sigma^V$ defined in
Remark~\ref{rem:chains}(\ref{item:chainSuspension}).
Clearly, $\sigma^0 = \id$ and $\sigma^W\circ\sigma^V = \sigma^{V\dirsum W}$.

If $\zeta\colon V\to V'$, then the following diagram commutes for any spectra $E$ and $F$
--- around either side it takes a map $f$ to $f\smsh S^\zeta$:
\[
 \xymatrix@C+2em{
  [E,F]_G \ar[r]^-{\sigma^V} \ar[d]_{\sigma^{V'}}
   & [\susp^V E, \susp^V F]_G \ar[d]^{[\id,\id\smsh\zeta]} \\
  [\susp^{V'}E,\susp^{V'}F]_G \ar[r]_-{[\id\smsh\zeta,\id]}
   & [\susp^V E, \susp^{V'} F]_G
 }
\]
From this it follows that the following diagram commutes:
\[
 \xymatrix@C+2em{
  \Mackey C_{\alpha+*}(X,*) \ar[r]^-{\sigma^V} \ar[d]_{\sigma^{V'}}
    & \Mackey C_{\alpha+V+*}(\susp^V X,*) \ar[d]^{\Mackey C_{\id}(\id\smsh\zeta)} \\
  \Mackey C_{\alpha+V'+*}(\susp^{V'} X,*)
    \ar[r]_-{\Mackey C_{\id+\zeta}(\id)}
    & \Mackey C_{\alpha+V+*}(\susp^{V'}X,*)
 }
\]
The naturality diagrams for homology and cohomology now follow.

{\em Dimension Axiom:}
Let $X = G_+\smsh_H S^{V-\delta(G/H)+n}$. Then
$X$ is a based $\delta$-$G$-CW($V$) complex with a single relative cell
of dimension $V+n$. We have then
\[
 \susp_G^\infty X^{V+n}/X^{V+n-1} \iso
  \susp^{V+n} \susp_G^\infty G_+\smsh_H S^{-\delta(G/H)}
\]
while the other filtration quotients are trivial, hence
\[
 \Mackey C_{V+k}(X,*) \iso
  \begin{cases}
   \MackeyOp A^{G/H,\delta} & \text{if $k = n$} \\
   0 & \text{if $k\neq n$.}
  \end{cases}
\]
The statements about homology and cohomology follow.
\end{proof}

We now want to extend the definition of cellular homology and cohomology to
arbitrary $G$-spaces. As mentioned at the beginning of
Section~\ref{sec:CWprespectra}, we cannot do this simply by taking
$\delta$-$G$-CW$(\alpha)$ approximations: If $\Gamma_\alpha X\to X$ is
a $\delta$-$G$-CW$(\alpha)$ approximation of $X$ and $W$ is a reprentation of $G$,
the map $\susp^W\Gamma_\alpha X\to \susp^W X$ need not be a
$\delta$-weak$_{\alpha+W}$ equivalence, so the resulting theory would not have
suspension isomorphisms.
Instead, we partially stabilize and use the
$\delta$-$G$-CW$(\alpha)$ prespectra discussed in that earlier section.

Let $E$ be a $\delta$-$G$-CW$(\alpha)$ prespectrum as in Section~\ref{sec:CWprespectra}.
For each $i$ we can consider the chain complex 
$\Mackey C_{\alpha+V_i+*}^{G,\delta}(E(V_i))$. The structure maps induce maps (inclusions, actually)
\[
 \Mackey C_{\alpha+V_{i-1}+*}^{G,\delta}(E(V_{i-1}))
  \iso \Mackey C_{\alpha+V_i+*}^{G,\delta}(\susp^{V_i - V_{i-1}}E(V_{i-1}))
  \to \Mackey C_{\alpha+V_{i}+*}^{G,\delta}(E(V_{i})).
\]

\begin{definition}
Let $E$ be a $\delta$-$G$-CW$(\alpha)$ prespectrum. We define the 
{\em cellular chain complex of $E$} to be the colimit
\[
 \Mackey C_{\alpha+*}^{G,\delta}(E) = \colim_i \Mackey C_{\alpha+V_i+*}(E(V_i)).
\]
If $F$ is an arbitrary prespectrum and $\delta$ is familial, we define
\[
 \Mackey C_{\alpha+*}^{G,\delta}(F) = \Mackey C_{\alpha+*}^{G,\delta}(\Gamma F)
\]
where $\Gamma F\to F$ is a $\delta$-$G$-CW$(\alpha)$ approximation of $F$.
If $X$ is a based $G$-space and $\susp^\infty X$ denotes its suspension prespectrum, we define
\[
 \Mackey C_{\alpha+*}^{G,\delta}(X) = \Mackey C_{\alpha+*}^{G,\delta}(\susp^\infty X)
  = \Mackey C_{\alpha+*}^{G,\delta}(\Gamma\susp^\infty X).
\]
\end{definition}

Our main interest is in the chains of an arbitrary $G$-space, so we will concentrate on that case,
using more general prespectra only as necessary. 
Because we are not really modeling the stable category by looking at $\delta$-weak$_\alpha$
equivalences,
the extra generality of discussing prespectra does not seem that useful in itself.

Of course, we want to make these chains functorial, so we make the following definitions.

\begin{definition}
Let $f\colon E\to F$ be a cellular map of $\delta$-$G$-CW$(\alpha)$ prespectra.
For each $i$ we have an induced chain map
\[
 (f_i)_*\colon \Mackey C_{\alpha+V_i+*}(E(V_i)) 
 		\to \Mackey C_{\alpha+V_i+*}(F(V_i))
\]
and these maps are compatible under suspension. We define
\[
 f_* = \colim_i (f_i)_*\colon \Mackey C_{\alpha+*}^{G,\delta}(E)
    \to \Mackey C_{\alpha+*}^{G,\delta}(F).
\]
If $f\colon E\to F$ is a map of arbitrary prespectra, we define
$f_*$ to be the map induced by a cellular approximation
$\Gamma f\colon \Gamma E \to \Gamma F$.
If $f\colon X\to Y$ is a map of based $G$-spaces, we let
\[
 f_*\colon \Mackey C_{\alpha+*}^{G,\delta}(X)
   \to \Mackey C_{\alpha+*}^{G,\delta}(Y)
\]
be the map associated with the suspension $\susp^\infty f\colon \susp^\infty X\to \susp^\infty Y$.
\end{definition}

It's easy to see that this makes $\Mackey C_{\alpha+*}^{G,\delta}$ functorial on the
category of $\delta$-$G$-CW$(\alpha)$ prespectra and cellular maps.
Note that we can extend the definition easily to semistable cellular maps.
The following is
also reassuring and is essentially just the observation that, if
$X$ is a based $\delta$-$G$-CW$(\alpha)$ complex, then
$\susp^\infty X$ is a $\delta$-$G$-CW$(\alpha)$ prespectrum.

\begin{proposition}
Suspension $\susp^\infty$ defines a functor from the category of 
based $\delta$-$G$-CW$(\alpha)$ complexes and cellular maps to the category
of $\delta$-$G$-CW$(\alpha)$ prespectra and cellular maps.
When $X$ is a based $\delta$-$G$-CW$(\alpha)$ complex we have a natural isomorphism
\[
 \Mackey C_{\alpha+*}^{G,\delta}(X) \iso \Mackey C_{\alpha+*}^{G,\delta}(\susp^\infty X)
\]
where the chains on the left are those of $X$ as a based $\delta$-$G$-CW$(\alpha)$ complex.
\qed
\end{proposition}

For prespectra or spaces that first need to be approximated, there are choices involved.
We start to unravel to what extent those choices matter.

\begin{proposition}\label{prop:chainsFunctorial}
\begin{enumerate}\item[]
\item
If $E$ is any $G$-prespectrum and $\Gamma_1 E\to E$ and $\Gamma_2 E \to E$ are two
approximations by $\delta$-$G$-CW$(\alpha)$ prespectra, 
then there is a canonical chain isomorphism
$\Mackey C_{\alpha+*}^{G,\delta}(\Gamma_1 E) \iso
 \Mackey C_{\alpha+*}^{G,\delta}(\Gamma_2 E)$.

\item
If $f\colon E\to F$ is any map of $G$-prespectra and
$\Gamma E\to E$ and $\Gamma F\to F$ are chosen approximations, the map
$f_*\colon \Mackey C_{\alpha+*}^{G,\delta}(\Gamma E)\to \Mackey C_{\alpha+*}^{G,\delta}(\Gamma F)$
is well-defined up to chain homotopy.
\end{enumerate}
\end{proposition}

\begin{proof}
For the second part, we note by
the last part of Proposition~\ref{prop:approxPrespectra} that the approximating cellular map
$\Gamma f$ is well-defined up to cellular homotopy. This implies that the chain map it induces
is well-defined up to chain isomorphism.
The first part follows from the second applied to the identity map on $E$.
\end{proof}

Put another way, if we {\em choose} for each prespectrum $E$ an approximation $\Gamma E\to E$,
then we get a functor from the category of prespectra to the category of chain complexes
modulo chain homotopy. Any two collections of choices of approximations lead to canonically
naturally isomorphic functors.
Here, we can interpret ``canonically'' as meaning that, if we consider the category
(groupoid, actually)
whose objects are collections of choices of approximations and in which there is a unique
morphism from any collection to any other, then the natural isomorphisms are functorial
on this groupoid.
  
Now, if $X$ is a based $G$-space,
$\Mackey C_{\alpha+*}^{G,\delta}(X)$ appears to depend also on the choice of universe $\U$ and
indexing sequence $\V$, but we now argue that, up to canonical natural chain homotopy equivalence,
it really does not.
We first show that the chains are independent of the choice of $\V$.
To emphasize the possible dependence on $\V$, we write
$\Mackey C_{\alpha+*}^{\V,\delta}(X)$ for chains defined using $\V$
as the indexing sequence.

\begin{proposition}\label{prop:sequenceIndependence}
Let $\V$ and $\W$ be two indexing sequences in the same universe $\U$.
Then the chain complexes $\Mackey C_{\alpha+*}^{\V,\delta}(X)$ and
$\Mackey C_{\alpha+*}^{\W,\delta}(X)$ are canonically naturally chain homotopy equivalent.
\end{proposition}

\begin{proof}
Let $\susp^\infty_\V X$ denote the suspension prespectrum of $X$ indexed on $\V$ and
let $\susp^\infty_\W X$ denote the one indexed on $\W$.
Choose $\delta$-$G$-CW$(\alpha)$ approximations
$\Gamma \susp^\infty_\V X \to \susp^\infty_\V X$ and
$\Gamma \susp^\infty_\W X \to \susp^\infty_\W X$,
so that, for example, $\Gamma \susp^\infty_\V X(V_i) \to \susp^{V_i} X$
is is an approximation by a $\delta$-$G$-CW$(\alpha+V_i)$ complex.
Because $\W$ is an indexing sequence, we can find a strictly increasing
sequence $\{j(i)\}$ such that $V_i\subset W_{j(i)}$ for each $i$.
Using the relative Whitehead theorem, we choose cellular maps
\[
 \zeta_i\colon \susp^{W_{j(i)}-V_i}\Gamma \susp^\infty_\V X(V_i) \to
   \Gamma \susp^\infty_\W X(W_{j(i)})
\]
over $\susp^{W_{j(i)}}X$ inductively so that the following diagram commutes for each $i>1$:
\[
 \xymatrix{
  \susp^{W_{j(i)}-V_{i-1}}\Gamma \susp^\infty_\V X(V_{i-1})
   \ar[d]_{\susp\sigma} \ar[r]^-{\susp \zeta_{i-1}}
   & \susp^{W_{j(i)}-W_{j(i-1)}}\Gamma \susp^\infty_\W X(W_{j(i-1)}) \ar[d]^\sigma \\
  \susp^{W_{j(i)}-V_i}\Gamma \susp^\infty_\V X(V_{i}) \ar[r]_{\zeta_i}
   & \Gamma \susp^\infty_\W X(W_{j(i)})
 }
\]
This determines a chain map 
$\zeta\colon\Mackey C_{\alpha+*}^{\V,\delta}(X)\to \Mackey C_{\alpha+*}^{\W,\delta}(X)$,
unique up to chain homotopy --- on each term in the colimit it is the composite
\begin{align*}
 \Mackey C_{\alpha+V_i+*}(\Gamma \susp^\infty_\V X(V_i))
  &\xrightarrow{\susp} \Mackey C_{\alpha+W_{j(i)}+*}(\susp^{W_{j(i)}-V_i}\Gamma \susp^\infty_\V X(V_i)) \\
  &\xrightarrow{(\zeta_i)_*} \Mackey C_{\alpha+W_{j(i)}+*}(\Gamma \susp^\infty_\W X(W_{j(i)}))
\end{align*}

When we say that these chain maps are canonical, we mean that the composite of
the map  $\Mackey C_{\alpha+*}^{\V,\delta}(X) \to \Mackey C_{\alpha+*}^{\W,\delta}(X)$
with the map  $\Mackey C_{\alpha+*}^{\W,\delta}(X) \to\Mackey C_{\alpha+*}^{\sZ,\delta}(X)$
agrees, up to chain homotopy, with the map 
$\Mackey C_{\alpha+*}^{\V,\delta}(X) \to \Mackey C_{\alpha+*}^{\sZ,\delta}(X)$.
This is a straightforward, if tedious, check.

In particular, the composite
\[
 \Mackey C_{\alpha+*}^{\V,\delta}(X) \to \Mackey C_{\alpha+*}^{\W,\delta}(X)
  \to \Mackey C_{\alpha+*}^{\V,\delta}(X)
\]
is chain homotopic to the identity, and similarly when $\V$ and $\W$ are reversed,
so the maps we've constructed are chain homotopy equivalences.
\end{proof}

Now suppose that $\U$ and $\U'$ are two complete $G$-universes.
(The following line of argument leads to a result similar in concept to the more general
\cite[II.1.7]{LMS:eqhomotopy}, but the latter result deals with the stable category,
so we give a separate argument here for our case.)
Because they are both
complete, they are isomorphic;
choose a linear isometric isomorphism $f\colon \U\to\U'$. If $\V$ is an indexing sequence in $\U$, let
$V'_i = f(V_i)$ for all $i$ and let $\V' = f(\V) = \{V'_i\}$, which is an indexing sequence in $\U'$.
If $\Gamma \susp^\infty_\V X$ is a CW approximation to $\susp^\infty_\V X$,
we get a CW approximation indexed on $\V'$ by taking the same spaces with the composite maps
\[
 \Gamma \susp^\infty_\V X(V_i) \to \susp^{V_i}X \xrightarrow{f}\susp^{V'_i}X.
\]
We freely use the fact that a $\delta$-$G$-CW($\alpha+V_i$) structure is also a $\delta$-$G$-CW($\alpha+V'_i$) structure.
In what follows, any CW approximation $\Gamma \susp^\infty_{\V'} X$ can be used in place
of the one just given and we do not assume any particular choice. 
In any case, the relative Whitehead theorem allows us to inductively define cellular maps $\bar f$ making the top square in the following
diagram commute and the bottom square commute up to homotopy:
\[
 \xymatrix{
 \susp^{V_i-V_{i-1}}\Gamma \susp^\infty_\V X(V_{i-1}) \ar[r]^{\bar f}\ar[d]
  & \susp^{V'_i-V'_{i-1}}\Gamma \susp^\infty_{\V'} X(V'_{i-1}) \ar[d] \\
 \Gamma \susp^\infty_\V X(V_i) \ar[r]^{\bar f}\ar[d]
  & \Gamma \susp^\infty_{\V'} X(V'_i) \ar[d] \\
 \susp^{V_{i}}X \ar[r]^f & \susp^{V'_{i}}X
 }
\]
(If we simply make the right hand side a copy of the left, the map commutes on the nose.)
The maps $\bar f$ then induce a chain map 
$f_*\colon \Mackey C^{\V,\delta}_{\alpha+*}(X) \to \Mackey C^{\V',\delta}_{\alpha+*}(X)$,
unique up to chain homotopy --- on each term in the colimit it is given by
\[
 \Mackey C_{\alpha+V_i+*}(\Gamma \susp^\infty_\V X(V_i))
  \xrightarrow{\bar f_*}
  \Mackey C_{\alpha+V'_i+*}(\Gamma \susp^\infty_{\V'} X(V'_i))
\]
The similar map $(f^{-1})_*$ is clearly a
chain homotopy inverse.

To what extent does the chain homotopy equivalence
$\Mackey C^{\V,\delta}_{\alpha+*}(X) \to \Mackey C^{\V',\delta}_{\alpha+*}(X)$
depend on the choice of $f$?
Hardly at all, in the following sense.

\begin{proposition}\label{prop:universeIndependence}
Let $f$ and $g$ be two linear isometric isomorphisms $\U\to\U'$.
Let $\V$ be an indexing sequence in $\U$, let $\V' = f(\V)$, and let $\W' = g(\V)$.
Let $\zeta\colon \Mackey C^{\V',\delta}_{\alpha+*}(X) \to \Mackey C^{\W',\delta}_{\alpha+*}(X)$
be the canonical chain homotopy equivalence of Proposition~\ref{prop:sequenceIndependence}. 
Then the following diagram
of chain homotopy equivalences commutes up to chain homotopy:
\[
 \xymatrix@C-3em{
   & \Mackey C^{\V,\delta}_{\alpha+*}(X) \ar[dl]_{f_*} \ar[dr]^{g_*} \\
   \Mackey C^{\V',\delta}_{\alpha+*}(X) \ar[rr]_\zeta
     && \Mackey C^{\W',\delta}_{\alpha+*}(X)
 }
\]
\end{proposition}

\begin{proof}
Choose a strictly increasing sequence $j(i)$ such that $V'_i\subset W'_{j(i)}$
and $j(i)\geq i$.
The map $\zeta f_*$ is represented by the composite
$(\zeta_i)_*\susp\bar f_*$ down the left side of the following commutative
diagram, hence by the composite across the top and down the right side.
\[
 \xymatrix{
  \Mackey C_{\alpha+V_i+*}(\Gamma\susp^\infty_{\V}X(V_i)) \ar[r]^-{\susp} \ar[d]_{\bar f_*}
   & \Mackey C_{\alpha+V_i+(W'_{j(i)}-V'_i)+*}(\susp^{W'_{j(i)}-V'_i}\Gamma\susp^\infty_{\V}X(V_i))
      \ar[d]^{(\overline{f+1})_*} \\
  \Mackey C_{\alpha+V'_i+*}(\Gamma\susp^\infty_{\V'}X(V'_i)) \ar[r]^-{\susp}
   & \Mackey C_{\alpha+W'_{j(i)}+*}(\susp^{W'_{j(i)}-V'_i}\Gamma\susp^\infty_{\V'}X(V'_i))
      \ar[d]^{(\zeta_i)_*} \\
  & \Mackey C_{\alpha+W'_{j(i)}+*}(\Gamma\susp^\infty_{\W'}X(W'_{j(i)}))
 }
\]
In order to compare this with $g_*$, we can assume that $W'_{j(i)}$ is so large that
there is a linear isometry $h\colon W'_{j(i)} - V'_i \to W'_{j(i)} - W'_i$ such that the map
\[
 V_i\dirsum (W'_{j(i)} - V'_i) \xrightarrow{g+h} W'_i\dirsum (W'_{j(i)} - W'_i) = W'_{j(i)}
\]
is homotopic through linear isometries to
\[
 V_i\dirsum (W'_{j(i)} - V'_i) \xrightarrow{f+1} V'_i\dirsum (W'_{j(i)} - V'_i) = W'_{j(i)}.
\]
Now consider the following cube:
\[
 \xymatrix@C-3.5em{
  & \susp^{W'_{j(i)}-V'_i}\Gamma\susp^\infty_{\V}X(V_i)
      \ar[dl]_{\susp\bar f} \ar[rr]^-{\susp^h\bar g} \ar'[d][dd]
   && \susp^{W'_{j(i)}-W'_i}\Gamma\susp^\infty_{\W'}X(W'_i) \ar[dd] \ar[dl]_\sigma \\
  \susp^{W'_{j(i)}-V'_i}\Gamma\susp^\infty_{\V'}X(V'_i) \ar[rr]^(0.65){\zeta_i} \ar[dd]
   && \Gamma\susp^\infty_{\W'}X(W'_{j(i)}) \ar[dd] \\
  & \susp^{W'_{j(i)}-V'_i}\susp^{V_i}X \ar'[r]_-{g+h}[rr] \ar[dl]_{f+1}
   && \susp^{W'_{j(i)}-W'_i}\susp^{W'_i}X \ar@{=}[dl] \\
  \susp^{W'_{j(i)}-V'_i}\susp^{V'_i}X \ar@{=}[rr]
   && \susp^{W'_{j(i)}}X
 }
\]
Because the four vertical squares and the bottom square commute up to homotopy,
the Whitehead theorem implies that the top square also commutes up to homotopy.
Moreover, we can choose a cellular homotopy.
Finally, we now have the following chain homotopy commutative diagram:
\[
 \def\objectstyle{\scriptstyle}
 \def\labelstyle{\scriptstyle}
 \xymatrix@C-1.5em{
  \Mackey C_{\alpha+V_i+*}(\Gamma\susp^\infty_{\V}X(V_i)) \ar[d]_{\susp} \ar[r]^{\bar g_*}
   & \Mackey C_{\alpha+W'_i+*}(\Gamma\susp^\infty_{\W'}X(W'_i)) \ar[d]^{\susp} \\
  \Mackey C_{\alpha+V_i+(W'_{j(i)}-V'_i)+*}(\susp^{W'_{j(i)}-V'_i}\Gamma\susp^\infty_{\V}X(V_i))
    \ar[r]^-{\bar g_*} \ar[dd]_{(\overline{f+1})_*}
   & \Mackey C_{\alpha+W'_i+(W'_{j(i)}-V'_i)+*}(\susp^{W'_{j(i)}-V'_i}\Gamma\susp^\infty_{\W'}X(W'_i))
    \ar[d]^{\susp^h} \\
  & \Mackey C_{\alpha+W'_i+(W'_{j(i)}-W'_i)+*}(\susp^{W'_{j(i)}-W'_i}\Gamma\susp^\infty_{\W'}X(W'_i))
    \ar[d]^{\sigma_*} \\
  \Mackey C_{\alpha+V'_i+(W'_{j(i)}-V'_i)+*}(\susp^{W'_{j(i)}-V'_i}\Gamma\susp^\infty_{\V'}X(V'_i))
    \ar[r]_-{(\zeta_i)_*}
   & \Mackey C_{\alpha+W'_{j(i)}+*}(\Gamma\susp^\infty_{\W'}X(W'_{j(i)}))
 }
\]
The composite down the left side represents $\zeta f_*(u)$ while the composite down the right
represents $g_*(u)$. Here, we use that the composite $\susp^h\circ\susp$ on the right
coincides with the simple suspension by $W'_{j(i)}-W'_i$, which follows by
inspecting the definition of the chain complexes.
Hence, $\zeta f_*$ is cellularly homotopic to $g_*$ as claimed.
\end{proof}

There is one more technical detail to deal with. Because suspension 
is guaranteed to respect weak $G$-equivalences
only for well-based $G$-spaces, we need to recall the following common definitions.

\begin{definition}\label{def:mappingCylinders}
\begin{enumerate}
\item[]
\item
Let $f\colon A\to X$ be a map of unbased $G$-spaces.
The {\em (unreduced) mapping cylinder, $Mf$,} is the pushout in the following diagram:
\[
 \xymatrix{
  A \ar[r]^f \ar[d]_{i_0} & X \ar[d] \\
  A\times I \ar[r] & Mf
 }
\]
The {\em (unreduced) mapping cone}
is $Cf = Mf/(A\times 1)$ with basepoint the image of $A\times 1$.

\item
Let $f\colon A\to X$ be a map of based $G$-spaces.
The {\em reduced mapping cylinder, $\tilde Mf$,} is the pushout in the following diagram:
\[
 \xymatrix{
  A \ar[r]^f \ar[d]_{i_0} & X \ar[d] \\
  A\smsh I_+ \ar[r] & \tilde Mf
 }
\]
The {\em reduced mapping cone, $\tilde Cf$,}
is the pushout in the following diagram, in which $I$ has basepoint 1:
\[
 \xymatrix{
  A \ar[r]^f \ar[d]_{i_0} & X \ar[d] \\
  A\smsh I \ar[r] & \tilde Cf
 }
\]
\end{enumerate}
In the special case in which $p\colon *\to X$ is the inclusion of the basepoint,
we call $Mp$ the {\em whiskering construction} and write $X_w = Mp$ for $X$ with a whisker attached.
\end{definition}

The point of the mapping cylinder is, of course, that the inclusion
$A\to Mf$ induced by $i_1\colon A\to A\times I$ is always a cofibration.
In particular, if we give $X_w$ the basepoint at the end of the whisker, it is
well-based.
Collapsing the whisker gives an unbased $G$-homotopy equivalence $X_w\to X$ that is
a based $G$-homotopy equivalence if $X$ is well-based. 
Hence we can use the whiskering construction as a functorial way of replacing any based $G$-space
with a weakly equivalent well-based space.

Likewise, if $f\colon A\to X$ is a based map, then the inclusion $A\to \tilde Mf$ is a based
cofibration. If $A$ and $X$ are well-based it is in fact an unbased cofibration as well.

With these definitions and results in place, we can now define the homology and cohomology of arbitrary $G$-spaces.

\begin{definition}
Let $\delta$ be a familial dimension function for $G$, let $\alpha$ be a virtual representation of $G$,
let $\Mackey T$ be a contravariant $\delta$-Mackey functor, and let $\MackeyOp S$ be a covariant $\delta$-Mackey functor.
If $X$ is a based $G$-space, let
\[
 \tilde H^{G,\delta}_{\alpha+n}(X;\MackeyOp S)
   = H_{\alpha+n}(\Mackey C^{G,\delta}_{\alpha+*}(X_w) \tensor_{\sorb{G,\delta}} \MackeyOp S)
\]
and
\[
 \tilde H_{G,\delta}^{\alpha+n}(X;\Mackey T)
   = H^{\alpha+n}(\Hom_{\sorb{G,\delta}}(\Mackey C^{G,\delta}_{\alpha+*}(X_w),\Mackey T)).
\]
\end{definition}

\begin{proposition}\label{prop:homologyWellDefined}
$\tilde H^{G,\delta}_{\alpha+n}(X;\MackeyOp S)$ is a well-defined covariant homotopy functor of $X$, while
$\tilde H_{G,\delta}^{\alpha+n}(X;\Mackey T)$ is a well-defined contravariant homotopy functor of $X$.
If $X$ is a $\delta$-$G$-CW$(\alpha)$ complex, these groups are naturally isomorphic to those
given by Definition~\ref{def:alphaHomology}.
\end{proposition}

\begin{proof}
That these groups are well-defined is shown by Propositions~\ref{prop:chainsFunctorial},
\ref{prop:sequenceIndependence}, and \ref{prop:universeIndependence}.
That they are homotopy functors of $X$ follows by the usual argument.

If $X$ is a $\delta$-$G$-CW$(\alpha)$ complex, then $\susp^\infty X$ is a
$\delta$-$G$-CW$(\alpha)$ prespectrum so serves as an approximation of itself.
There is then an isomorphism of chains
$\Mackey C^{G,\delta}_{\alpha+*}(X) \iso \Mackey C^{G,\delta}_{\alpha+*}(\susp^\infty X)$
given by the suspension maps.
\end{proof}

\begin{theorem}[Reduced Homology and Cohomology of Spaces]\label{thm:reducedHomologySpaces}
Let $\delta$ be a familial dimension function for $G$, let $\alpha$ be a virtual representation of $G$,
and let $\MackeyOp S$ and $\Mackey T$ be respectively a covariant and a contravariant
$\delta$-Mackey functor. Then the abelian groups
$\tilde H^{G,\delta}_{\alpha}(X;\MackeyOp S)$ and $\tilde H_{G,\delta}^{\alpha}(X;\Mackey T)$ are
respectively covariant and contravariant
functors on the homotopy category of 
based $G$-spaces.
They are also respectively contravariant and covariant functors 
of $\alpha$.
These functors satisfy the following properties.
\begin{enumerate}

\item (Weak Equivalence)
If $f\colon X\to Y$ is an $\F(\delta)$-equivalence of based $G$-spaces, then
\[
 f_*\colon \tilde H^{G,\delta}_{\alpha}(X;\MackeyOp S) \to \tilde H^{G,\delta}_{\alpha}(Y;\MackeyOp S)
\]
and
\[
 f^*\colon \tilde H_{G,\delta}^{\alpha}(Y;\Mackey T) \to \tilde H_{G,\delta}^{\alpha}(X;\Mackey T)
\]
are isomorphisms.

\item (Exactness)
If $A\to X$ is a cofibration, then the following sequences are exact:
\[
 \tilde H^{G,\delta}_{\alpha}(A;\MackeyOp S) \to \tilde H^{G,\delta}_{\alpha}(X;\MackeyOp S)
  \to \tilde H^{G,\delta}_{\alpha}(X/A;\MackeyOp S)
\]
and
\[
 \tilde H_{G,\delta}^{\alpha}(X/A;\Mackey T) \to \tilde H_{G,\delta}^{\alpha}(X;\Mackey T)
  \to \tilde H_{G,\delta}^{\alpha}(A;\Mackey T).
\]

\item (Additivity)
If $X = \Wedge_i X_i$ is a wedge of well-based $G$-spaces,
then the inclusions of the wedge summands induce isomorphisms
\[
 \Dirsum_i\tilde H^{G,\delta}_{\alpha}(X_i;\MackeyOp S) \iso \tilde H^{G,\delta}_{\alpha}(X;\MackeyOp S)
\]
and
\[
 \tilde H_{G,\delta}^{\alpha}(X;\Mackey T) \iso \prod_i \tilde H_{G,\delta}^{\alpha}(X_i;\Mackey T).
\]

\item (Suspension)
If $X$ is well-based, there are suspension isomorphisms
 \[
  \sigma^V\colon \tilde H^{G,\delta}_{\alpha}(X;\MackeyOp S) 
     \xrightarrow{\iso} \tilde H^{G,\delta}_{\alpha+V}(\susp^V X;\MackeyOp S)
 \]
and
 \[
 \sigma^V\colon \tilde H_{G,\delta}^{\alpha}(X;\Mackey T) 
     \xrightarrow{\iso} \tilde H_{G,\delta}^{\alpha+V}(\susp^V X;\Mackey T)
 \]
These isomorphisms satisfy $\sigma^0 = \id$, $\sigma^{W}\circ \sigma^V = \sigma^{V\oplus W}$,
and the following naturality condition: If $\zeta\colon V\to V'$ is an isomorphism,
then the following diagrams commute (cf \cite[XIII.1.1]{May:alaska}):
\[
 \xymatrix@C+2em{
  \tilde H^{G,\delta}_\alpha(X; \MackeyOp S) \ar[r]^-{\sigma^V} \ar[d]_{\sigma^{V'}}
    & \tilde H^{G,\delta}_{\alpha + V}(\susp^{V}X; \MackeyOp S)
          \ar[d]^{\tilde H_{\id}(\id\smsh\zeta)} \\
  \tilde H^{G,\delta}_{\alpha + V'}(\susp^{V'}X; \MackeyOp S)
       \ar[r]_-{\tilde H_{\id+\zeta}(\id)}
    & \tilde H^{G,\delta}_{\alpha + V}(\susp^{V'}X; \MackeyOp S)
 }
\]
\[
 \xymatrix@C+2em{
  \tilde H_{G,\delta}^\alpha(X; \Mackey T) \ar[r]^-{\sigma^V} \ar[d]_{\sigma^{V'}}
    & \tilde H_{G,\delta}^{\alpha + V}(\susp^{V}X; \Mackey T)
          \ar[d]^{\tilde H^{\id+\zeta}(\id)} \\
  \tilde H_{G,\delta}^{\alpha + V'}(\susp^{V'}X; \Mackey T)
       \ar[r]_-{\tilde H^{\id}(\id\smsh\zeta)}
    & \tilde H_{G,\delta}^{\alpha + V'}(\susp^{V}X; \Mackey T)
 }
\]

\item (Dimension Axiom)
If $H\in\F(\delta)$ and $V$ is a representation of $H$ so large that
$V-\delta(G/H)+n$ is an actual representation, then
there are natural isomorphisms
\[
 \tilde H^{G,\delta}_{V+k}(G_+\smsh_H S^{V-\delta(G/H)+n};\MackeyOp S)
  \iso
   \begin{cases}
    \MackeyOp S(G/H,\delta) & \text{if $k = n$} \\
    0 & \text{if $k\neq n$}
   \end{cases}
\]
and
\[
 \tilde H_{G,\delta}^{V+k}(G_+\smsh_H S^{V-\delta(G/H)+n};\Mackey T)
  \iso
   \begin{cases}
    \Mackey T(G/H,\delta) & \text{if $k = n$} \\
    0 & \text{if $k\neq n$.}
   \end{cases}
\]

\end{enumerate}
\end{theorem}

\begin{proof}
That homology and cohomology are homotopy functors is Proposition~\ref{prop:homologyWellDefined}.
Functoriality in $\alpha$ on the chain level
follows from the same property for spaces
shown in the proof of Theorem~\ref{thm:reducedHomologyComplexes}.

{\em Weak Equivalence:}
If $f\colon X\to Y$ is an $\F(\delta)$-equivalence, then any approximation
$\Gamma f\colon \Gamma\susp^\infty X_w 
\to \Gamma\susp^\infty Y_w$
is a $\delta$-weak$_\alpha$ equivalence. 
(Here, we use that $X_w$ and $Y_w$ are well-based to conclude that each
$\susp^{V_i}X_w\to \susp^{V_i}Y_w$ is an $\F(\delta)$-equivalence,
hence a $\delta$-weak$_{\alpha+V_i}$-equivalence.)
We've shown previously that $\Gamma f$ then has a $G$-homotopy inverse,
hence that $f_*$ and $f^*$ are isomorphisms
in homology and cohomology, respectively.

{\em Exactness}
Let $i\colon A\to X$ be a cofibration.
Using Theorem~\ref{thm:generalCellularApprox} we can choose approximations
$\Gamma i\colon\Gamma\susp^\infty A_w 
\to \Gamma\susp^\infty X_w$ with $\Gamma i$ the inclusion of a subcomplex.
At each level, by Proposition~\ref{prop:genExcisiveTriads} we have a weak equivalence of mapping cones
\[
 C(\Gamma i(V_k)) \to C(\susp^{V_k}A_w \to \susp^{V_k} X_w).
\]
We also have homotopy equivalences
\[
 C(\Gamma i(V_k)) \hmtpc (\Gamma\susp^\infty_\V X(V_i))/(\Gamma\susp^\infty_\V A(V_i))
\]
and
\[
 C(\susp^{V_k}A_w \to \susp^{V_k} X_w) \hmtpc \susp^{V_k}(X/A)_w
\]
We conclude that
\[
 \Gamma\susp^\infty X_w/\Gamma\susp^\infty A_w
  \to \susp^\infty (X/A)_w
\]
is a CW approximation, from which it follows that
\[
 \Mackey C^{G,\delta}_{\alpha+*}((X/A)_w) 
  = \Mackey C^{G,\delta}_{\alpha+*}(X_w)/\Mackey C^{G,\delta}_{\alpha+*}(A_w).
\]
(That is, we get equality when we choose the CW approximation of $\susp^\infty (X/A)_w$ 
to be $\Gamma\susp^\infty X_w/\Gamma\susp^\infty A_w$.)
The exact sequences now follow from familiar homological algebra.

{\em Additivity}
This follows from the fact that the wedge of well-based $G$-spaces preserves
$\delta$-weak$_{\alpha+V_i}$ equivalence, so that the wedge of CW approximations is again a CW approximation,
and the chains of a wedge can be taken to be the direct sum of the chains of the wedge summands.

{\em Suspension}
This is where we need the full force of our approximation by CW prespectra.
Assuming that $X$ is well-based, we may use $X$ in place of $X_w$.
To define $\sigma^V$, we begin with the case of $V\subset \U$.
Let $\V$ be an indexing sequence in $\U$ such that $V\subset V_i$ for all $i$.
Then $\V-V = \{ V_i - V \}$ is an indexing sequence in $\U-V$.
Given a CW approximation $\Gamma\susp^\infty_\V X$, we can reinterpret
each $\Gamma\susp^\infty_\V X(V_i) \to \susp^{V_i}X$ as
\[
 (\Gamma\susp^\infty_{\V-V} \susp^V X)(V_i-V) \to \susp^{V_i-V}\susp^V X.
\]
This then gives $\Gamma\susp^\infty_{\V-V} \susp^V X \to \susp^\infty_{\V-V}\susp^V X$,
a $\delta$-$G$-CW$(\alpha+V)$ approximation of $\susp^V X$ indexed on $\V-V$.
From this we get an isomorphism
\[
 \sigma^V \colon \Mackey C^{\V,\delta}_{\alpha+*}(X) \iso
   \Mackey C^{\V-V,\delta}_{\alpha+V+*}(\susp^V X).
\]
For a general $V$, choose an isomorphism $\xi\colon V\to V_0$ where $V_0\subset\U$.
Define
$\sigma^V\colon \Mackey C^{\V,\delta}_{\alpha+*}(X)
 \to \Mackey C^{\V-V_0,\delta}_{\alpha+V+*}(\susp^V X)$
to be the composite
\[
 \Mackey C^{\V,\delta}_{\alpha+*}(X)
  \xrightarrow{\sigma^{V_0}} \Mackey C^{\V-V_0,\delta}_{\alpha+V_0+*}(\susp^{V_0}X)
  \xrightarrow{\Mackey C_\xi(\id\smsh\xi^{-1})}
    \Mackey C^{\V-V_0,\delta}_{\alpha+V+*}(\susp^{V}X).
\]
This is a chain isomorphism so induces isomorphisms in homology and cohomology that we also
call $\sigma^V$. Note that we use the fact that the chains are essentially independent of
choice of universe and indexing sequence to justify using chains based on $\U-V_0$ and $\V-V_0$.

Transitivity, $\sigma^W\circ\sigma^V = \sigma^{V\oplus W}$, can be checked directly from
this definition.

To show naturality, suppose $\zeta\colon V\to V'$, choose $\xi'\colon V'\to V_0$
with $V_0\subset\U$ as above, and let $\xi = \xi'\circ\zeta\colon V\to V_0$. Then naturality
follows from the (chain homotopy) commutativity of the following diagram:
\[
 \xymatrix@C+2em{
  \Mackey C^{\V,\delta}_{\alpha+*}(X) \ar[d]_{\sigma^{V_0}} \\
  \Mackey C^{\V-V_0,\delta}_{\alpha+V_0+*}(\susp^{V_0}X)
     \ar[d]_{\Mackey C_{\xi'}(\id\smsh{\xi'}^{-1})}
     \ar[r]^{\Mackey C_{\xi}(\id\smsh{\xi}^{-1})}
   & \Mackey C^{\V-V_0,\delta}_{\alpha+V+*}(\susp^{V}X)
       \ar[d]^{\Mackey C_{\id}(\id\smsh\zeta)} \\
  \Mackey C^{\V-V_0,\delta}_{\alpha+V'+*}(\susp^{V'}X)
     \ar[r]_{\Mackey C_{\zeta}(\id)}
   & \Mackey C^{\V-V_0,\delta}_{\alpha+V+*}(\susp^{V'}X)
 }
\]

{\em Dimension Axiom}
This was already shown in Theorem~\ref{thm:reducedHomologyComplexes}.
(We apply Proposition~\ref{prop:homologyWellDefined} to the based
complex $G_+\smsh_H S^{V-\delta(G/H)+n}$.)
\end{proof}

\begin{corollary}\label{cor:exactsequences}
\begin{enumerate}
\item[]
\item
If $f\colon A\to X$ is a map of well-based $G$-spaces, then the following sequences are exact:
\[
 \tilde H^{G,\delta}_\alpha(A; \MackeyOp S) \to \tilde H^{G,\delta}_\alpha(X; \MackeyOp S)
  \to \tilde H^{G,\delta}_\alpha(\tilde Cf; \MackeyOp S)
\]
and
\[
 \tilde H_{G,\delta}^\alpha(\tilde Cf; \Mackey T) \to \tilde H_{G,\delta}^\alpha(X; \Mackey T)
  \to \tilde H_{G,\delta}^\alpha(A; \Mackey T).
\]

\item
If $f\colon A\to X$ is any map of $G$-spaces, then the following sequences are exact,
where $f_w\colon A_w\to X_w$ is the induced map:
\[
 \tilde H^{G,\delta}_\alpha(A; \MackeyOp S) \to \tilde H^{G,\delta}_\alpha(X; \MackeyOp S)
  \to \tilde H^{G,\delta}_\alpha(\tilde C(f_w); \MackeyOp S)
\]
and
\[
 \tilde H_{G,\delta}^\alpha(\tilde C(f_w); \Mackey T) \to \tilde H_{G,\delta}^\alpha(X; \Mackey T)
  \to \tilde H_{G,\delta}^\alpha(A; \Mackey T).
\]

\end{enumerate}
\end{corollary}

\begin{proof}
If $f\colon A\to X$ is a map of well-based $G$-spaces, then $f$ factors as
$A\to \tilde Mf \to X$, where the first map is a cofibration and the second is
an equivalence. Because $\tilde Cf = \tilde Mf/A$, the claimed exact sequences
in (1) now follow from the preceding theorem.

If $f\colon A\to X$ is any based map, we apply (1) to $f_w$ to get
the exact sequences claimed in (2).
\end{proof}

Finally, we define unreduced homology and cohomology of pairs.

\begin{definition}
If $(X,A)$ is a pair of $G$-spaces, write $i\colon A\to X$ for the inclusion and let
\[
 H^{G,\delta}_\alpha(X,A;\MackeyOp S) = \tilde H^{G,\delta}_\alpha(Ci;\MackeyOp S)
\]
and
\[
 H_{G,\delta}^\alpha(X,A;\Mackey T) = \tilde H_{G,\delta}^\alpha(Ci;\Mackey T).
\]
In particular, we write 
\[
 H^{G,\delta}_\alpha(X;\MackeyOp S) = H^{G,\delta}_\alpha(X,\emptyset;\MackeyOp S)
   = \tilde H^{G,\delta}_\alpha(X_+;\MackeyOp S)
\]
and
\[
 H_{G,\delta}^\alpha(X;\Mackey T) = H_{G,\delta}^\alpha(X,\emptyset;\Mackey T) 
   = \tilde H_{G,\delta}^\alpha(X_+;\Mackey T).
\]
\end{definition}

\begin{theorem}[Unreduced Homology and Cohomology]\label{thm:unreducedHomology}
Let $\delta$ be a familial dimension function for $G$, let $\alpha$ be a virtual representation of $G$,
and let $\MackeyOp S$ and $\Mackey T$ be respectively a covariant and a contravariant
$\delta$-Mackey functor. Then the abelian groups
$H^{G,\delta}_{\alpha}(X,A;\MackeyOp S)$ and $H_{G,\delta}^{\alpha}(X,A;\Mackey T)$ are
respectively covariant and contravariant
functors on the homotopy category of 
pairs of $G$-spaces.
They are also respectively contravariant and covariant functors 
of $\alpha$.
These functors satisfy the following properties.
\begin{enumerate}

\item (Weak Equivalence)
If $f\colon (X,A)\to (Y,B)$ is an $\F(\delta)$-equivalence of pairs of $G$-spaces, then
\[
 f_*\colon H^{G,\delta}_{\alpha}(X,A;\MackeyOp S) \to H^{G,\delta}_{\alpha}(Y,B;\MackeyOp S)
\]
and
\[
 f^*\colon H_{G,\delta}^{\alpha}(Y,B;\Mackey T) \to H_{G,\delta}^{\alpha}(X,A;\Mackey T)
\]
are isomorphisms.

\item (Exactness)
If $(X,A)$ is a pair of $G$-spaces, then there are natural homomorphisms
\[
 \bndry\colon H^{G,\delta}_{\alpha+n}(X,A;\MackeyOp S) \to H^{G,\delta}_{\alpha+n-1}(A;\MackeyOp S)
\]
and
\[
 d\colon H_{G,\delta}^{\alpha+n}(A;\Mackey T) \to H_{G,\delta}^{\alpha+n+1}(X,A;\Mackey T)
\]
and long exact sequences
\[
 \cdots \to H^{G,\delta}_{\alpha+n}(A; \MackeyOp S) \to H^{G,\delta}_{\alpha+n}(X; \MackeyOp S)
  \to H^{G,\delta}_{\alpha+n}(X,A; \MackeyOp S) \to H^{G,\delta}_{\alpha+n-1}(A;\MackeyOp S) \to \cdots
\]
and
\[
 \cdots \to H_{G,\delta}^{\alpha+n-1}(A;\Mackey T) \to H_{G,\delta}^{\alpha+n}(X,A; \Mackey T) 
  \to H_{G,\delta}^{\alpha+n}(X; \Mackey T)  \to H_{G,\delta}^{\alpha+n}(A; \Mackey T) \to \cdots.
\]

\item (Excision)
If $(X;A,B)$ is an excisive triad, i.e., $X$ is the union of the interiors of $A$ and $B$,
then the inclusion $(A,A\intersect B)\to (X,B)$ induces isomorphisms
\[
 H^{G,\delta}_{\alpha}(A,A\intersect B;\MackeyOp S) \iso H^{G,\delta}_{\alpha}(X,B;\MackeyOp S)
\]
and
\[
 H_{G,\delta}^{\alpha}(X,B;\Mackey T) \iso H_{G,\delta}^{\alpha}(A,A\intersect B;\Mackey T).
\]

\item (Additivity)
If $(X,A) = \coprod_k (X_k,A_k)$ is a disjoint union of pairs of $G$-spaces,
then the inclusions $(X_k,A_k)\to (X,A)$ induce isomorphisms
\[
 \Dirsum_k H^{G,\delta}_{\alpha}(X_k,A_k; \MackeyOp S) \iso H^{G,\delta}_{\alpha}(X,A; \MackeyOp S)
\]
and
\[
 H_{G,\delta}^{\alpha}(X,A; \Mackey T) \iso \prod_k H_{G,\delta}^{\alpha}(X_k,A_k; \Mackey T).
\]

\item (Suspension)
If $A\to X$ is a cofibration, there are suspension isomorphisms
 \[
  \sigma^V\colon H^{G,\delta}_{\alpha}(X,A; \MackeyOp S) 
     \xrightarrow{\iso} H^{G,\delta}_{\alpha+V}((X,A)\times(D(V),S(V)); \MackeyOp S)
 \]
and
 \[
 \sigma^V\colon H_{G,\delta}^{\alpha}(X,A; \Mackey T) 
     \xrightarrow{\iso} H_{G,\delta}^{\alpha+V}((X,A)\times(D(V),S(V)); \Mackey T),
 \]
where $(X,A)\times (D(V),S(V)) = (X\times D(V), X\times S(V)\union A\times D(V))$.
These isomorphisms satisfy $\sigma^0 = \id$ and $\sigma^{W}\circ \sigma^V = \sigma^{V\oplus W}$,
under the identification $\bar D(V) \times \bar D(W) \homeo \bar D(V\oplus W)$
(here we use the notation $\bar D(V) = (D(V),S(V))$).
They also satisfy the following naturality condition: If $\zeta\colon V\to V'$ is an isomorphism,
then the following diagrams commute:
\[
 \xymatrix@C+2em{
  H^{G,\delta}_{\alpha}(X,A; \MackeyOp S) \ar[r]^-{\sigma^V} \ar[d]_{\sigma^{V'}}
    & H^{G,\delta}_{\alpha+V}((X,A)\times \bar D(V); \MackeyOp S)
          \ar[d]^{H^G_{\id}(\id\times\zeta)} \\
  H^{G,\delta}_{\alpha+V'}((X,A)\times \bar D(V'); \MackeyOp S)
       \ar[r]_-{H^G_{\id\dirsum\zeta}(\id)}
    & H^{G,\delta}_{\alpha+V}((X,A)\times \bar D(V'); \MackeyOp S)
 }
\]
\[
 \xymatrix@C+2em{
  H_{G,\delta}^{\alpha}(X,A; \Mackey T) \ar[r]^-{\sigma^V} \ar[d]_{\sigma^{V'}}
    & H_{G,\delta}^{\alpha+V}((X,A)\times \bar D(V); \Mackey T)
          \ar[d]^{H_G^{\id\dirsum\zeta}(\id)} \\
  H_{G,\delta}^{\alpha+V'}((X,A)\times \bar D(V'); \Mackey T)
       \ar[r]_-{H_G^{\id}(\id\times\zeta)}
    & H_{G,\delta}^{\alpha+V'}((X,A)\times \bar D(V); \Mackey T)
 }
\]

\item (Dimension Axiom)
If $H\in\F(\delta)$ and $V$ is a representation of $H$ so large that
$V-\delta(G/H)+n$ is an actual representation, then there are natural isomorphisms
\[
 H^{G,\delta}_{V+k}(G\times_H \bar D(V-\delta(G/H)+n);\MackeyOp S)
  \iso
   \begin{cases}
    \MackeyOp S(G/H) & \text{if $k = n$} \\
    0 & \text{if $k\neq n$}
   \end{cases}
\]
and
\[
 H_{G,\delta}^{V+k}(G\times_H \bar D(V-\delta(G/H)+n);\Mackey T)
  \iso
   \begin{cases}
    \Mackey T(G/H) & \text{if $k = n$} \\
    0 & \text{if $k\neq n$.}
   \end{cases}
\]

\end{enumerate}
\end{theorem}

\begin{proof}
Let $i\colon A\to X$ be the inclusion.
The claims in the first paragraph translate into claims about the reduced homology
and cohomology of $Ci$ that we have already proven.

{\em Weak Equivalence}
If $i\colon A\to X$ and $j\colon B\to Y$ are the inclusions, then $Ci\to Cj$ is
an $\F(\delta)$-equivalence by an argument using Proposition~\ref{prop:genExcisiveTriads}.
The result follows because reduced homology and cohomology take $\F(\delta)$-equivalences to isomorphisms.

{\em Exactness}
If $i\colon A\to X$ is the inclusion, we consider the cofiber sequence
\[
 A_+ \to X_+ \to Ci \to \susp A_+ \to \susp X_+ \to \cdots
\]
where we use that $\tilde C(i_+) \homeo Ci$.
As in the nonequivariant case (cf \cite[\S8.4]{May:concise}), each pair of maps is, up to sign,
$G$-homotopy equivalent to a map followed by the inclusion of its target in its mapping cone.
Thus, on applying reduced homology and cohomology, we get exact sequences,
which we identify as the long exact sequences claimed, using the suspension isomorphisms.
But, note that we must use the {\em internal} suspension discussed after
Remark~\ref{rem:chains}.
The maps $\bndry$ and $d$ are then the natural homomorphisms induced by
$Ci\to \susp A_+$.

{\em Excision}
Let $i\colon A\intersect B\to A$ and $j\colon B\to X$ be the inclusions.
Consider the following map of triads: $(Ci;Ci,C(A\intersect B))\to (Cj;Ci,CB)$.
Each of these triads is excisive, and the three maps
$Ci\to Ci$, $C(A\intersect B)\to CB$, and
$Ci\intersect C(A\intersect B) = C(A\intersect B) \to Ci\intersect CB = C(A\intersect B)$
are all weak equivalences. Therefore, by Proposition~\ref{prop:genExcisiveTriads},
the inclusion $Ci\to Cj$ is a weak equivalence.
The excision isomorphisms follow.

{\em Additivity}
Let $i\colon A\to X$ and $i_k\colon A_k \to X_k$ be the inclusions.
Then $Ci = \Wedge_k Ci_k$, so additivity follows from the reduced case.

{\em Suspension}
Under the assumption that $i\colon A\to X$ is a cofibration,
$j\colon X\times S(V)\union A\times D(V)\to X\times D(V)$ is also a cofibration, and
$Ci \hmtpc X/A$ and $Cj \hmtpc \susp^V X/A$. The result now follows from the reduced
suspension isomorphism.

{\em Dimension Axiom}
If $i\colon G\times_H S(V-\delta(G/H)+n)\to G\times_H D(V-\delta(G/H)+n)$, 
then $Ci \homeo G_+\smsh_H S^{V-\delta(G/H)+n}$,
so the result follows from the reduced case.
\end{proof}

\begin{definition}\label{def:genMackeyvalued}
Recall from \cite{Wi:duality} or \cite[XIII.1]{May:alaska}
that a {\em (generalized) reduced $RO(G)$-graded homology or cohomology theory}
is a collection of functors satisfying all of Theorem~\ref{thm:reducedHomologySpaces}
except that it is required to take weak $G$-equivalences
(not necessarily $\F(\delta)$-equivalences) to isomorphisms and
it need not obey the dimension axiom.
Given a reduced $RO(G)$-graded homology theory $\tilde h^G_*(-)$
and a dimension function $\delta$ for $G$, 
we can regard $\tilde h^G_*(-)$ as $G$-$\delta$-Mackey functor-valued by defining
the covariant functor
$\MackeyOp h^{G,\delta}_\alpha(X)$, for a well-based $G$-space $X$, by
\[
 \MackeyOp h^{G,\delta}_\alpha(X)(G/H,\delta) = \tilde h^G_{\alpha}(X\smsh G_+\smsh S^{-\delta(G/H)})
   = \tilde h^G_{\alpha+V}(X\smsh G_+\smsh S^{V-\delta(G/H)}),
\]
where $V$ is a representation of $G$ so large that $\delta(G/H)\subset V$.
In particular, we write
$\MackeyOp h^{G,\delta}_* = \MackeyOp h^{G,\delta}_*(S^0)$ for the 
{\em coefficient system} of $\tilde h^G_*(-)$.
Similarly, if $\tilde h_G^*(-)$ is a reduced $RO(G)$-graded cohomology theory,
we write $\Mackey h_{G,\delta}^\alpha(X)$ for the contravariant $\delta$-Mackey functor defined by
\[
 \Mackey h_{G,\delta}^\alpha(X)(G/H,\delta) = \tilde h_G^\alpha(X\smsh G_+\smsh S^{-\delta(G/H)})
   = \tilde h_G^{\alpha+V}(X\smsh G_+\smsh S^{V-\delta(G/H)})
\]
and we write $\Mackey h_{G,\delta}^* = \Mackey h_{G,\delta}^*(S^0)$.
\end{definition}

In these terms, the dimension axioms take the form of the following statements
about the coefficient systems of cellular homology and cohomology:
\[
 \MackeyOp H^{G,\delta}_{n}(S^0;\MackeyOp S)
  \iso \begin{cases}
          \MackeyOp S &\text{if $n = 0$} \\
          0 &\text{if $n \neq 0$}
       \end{cases}
\]
and
\[
 \Mackey H_{G,\delta}^{n}(S^0;\Mackey T)
  \iso \begin{cases}
          \Mackey T &\text{if $n = 0$} \\
          0 &\text{if $n \neq 0$.}
       \end{cases}
\]

The Atiyah-Hirzebruch spectral sequence generalizes as follows.

\begin{theorem}[Atiyah-Hirzebruch Spectral Sequence]\label{thm:AtiyahHirzebruch}
Suppose that $\tilde h^G_*(-)$ is a reduced $RO(G)$-graded homology theory,
let $\delta$ be a familial dimension function for $G$, and
let $\alpha$ and $\beta$ be virtual representations of $G$.
Assume that $\tilde h^G_*(-)$ takes $\F(\delta)$-equivalences to isomorphisms.
Then there is a strongly convergent spectral sequence
\[
 E^2_{p,q} = \tilde H^{G,\delta}_{\alpha+p}(X;\MackeyOp h^{G,\delta}_{\beta+q}) \convto
   \tilde h^G_{\alpha+\beta+p+q}(X).
\]
Similarly, if $\tilde h_G^*(-)$ is a reduced $RO(G)$-graded cohomology theory
taking $\F(\delta)$-equivalences to isomorphisms,
there is a conditionally convergent spectral sequence
\[
 E_2^{p,q} = \tilde H_{G,\delta}^{\alpha+p}(X;\Mackey h_{G,\delta}^{\beta+q}) \convto
   \tilde h_G^{\alpha+\beta+p+q}(X).
\]
\end{theorem}

\begin{proof}
These spectral sequences arise in the usual way.
(See, for example, \cite[\S12]{Bo:spectralsequences}.)
Recall that, for sufficiently large $V$, $\delta$-weak$_{\alpha+V}$ equivalence
coincides with $\F(\delta)$-equivalence;
by taking suspensions,
it suffices to consider the case where $\delta$-weak$_\alpha$ equivalence
coincides with $\F(\delta)$-equivalence.
By taking $\delta$-$G$-CW($\alpha$) approximation,
using that $\tilde h^G_*(-)$ takes $\F(\delta)$-equivalences to isomorphisms,
it suffices to consider the case
where $X$ is a based $\delta$-$G$-CW($\alpha$) complex.
The skeletal filtration $\{X^{\alpha+*}\}$ leads to an exact couple on applying
either $\tilde h^G_{\alpha+\beta+*}(-)$ or $\tilde h_G^{\alpha+\beta+*}(-)$.
We identify the $E^1$ and $E_1$ terms
using the natural isomorphisms
\[
 E^1_{p,q} = \tilde h^G_{\alpha+\beta+p+q}(X^{\alpha+p}/X^{\alpha+p-1})
  \iso \Mackey C^{G,\delta}_{\alpha+p}(X,*) \tensor_{\sorb{G,\delta}} \MackeyOp h^{G,\delta}_{\beta+q}
\]
and
\[
 E_1^{p,q} = \tilde h_G^{\alpha+\beta+p+q}(X^{\alpha+p}/X^{\alpha+p-1})
  \iso \Hom_{\sorb{G,\delta}}(\Mackey C^{G,\delta}_{\alpha+p}(X,*), \Mackey h_{G,\delta}^{\beta+q}),
\]
and the spectral sequences follow.

It will be useful later to show explicitly how the natural isomorphisms displayed above arise.
In the case of homology, recall that
\begin{align*}
 \Mackey C^{G,\delta}_{\alpha+p}&(X,*) \tensor_{\sorb{G,\delta}} \MackeyOp h^{G,\delta}_{\beta+q} \\
  &\iso \int^{(G/H,\delta)} [G_+\smsh_H S^{-\delta(G/H)+\alpha+p}, \susp_G^\infty X^{\alpha+p}/X^{\alpha+p-1}]_G \\
       &\hphantom{\int^{(G/H,\delta)}}\qquad \tensor \tilde h^G_{\beta+q}(G_+\smsh_H S^{-\delta(G/H)}) \\
  &\iso \int^{(G/H,\delta)} [G_+\smsh_H S^{-\delta(G/H)+\alpha+p}, \susp_G^\infty X^{\alpha+p}/X^{\alpha+p-1}]_G \\
       &\hphantom{\int^{(G/H,\delta)}}\qquad \tensor \tilde h^G_{\alpha+\beta+p+q}(G_+\smsh_H S^{-\delta(G/H)+\alpha+p}).
\end{align*}
Evaluation defines a natural map
\begin{multline*}
 [G_+\smsh_H S^{-\delta(G/H)+\alpha+p}, \susp_G^\infty X^{\alpha+p}/X^{\alpha+p-1}]_G \\
        {}\tensor \tilde h^G_{\alpha+\beta+p+q}(G_+\smsh_H S^{-\delta(G/H)+\alpha+p}) 
  \to \tilde h^G_{\alpha+\beta+p+q}(X^{\alpha+p}/X^{\alpha+p-1}).
\end{multline*}
Note that, because $G/H$ is compact, any stable map
from $G_+\smsh_H S^{-\delta(G/H)+\alpha+p}$
into $\susp_G^\infty X^{\alpha+p}/X^{\alpha+p-1}$ is represented
by a map of $G$-spaces
\[
 G_+\smsh_H S^{V-\delta(G/H)+\alpha+p} \to \susp_G^V X^{\alpha+p}/X^{\alpha+p-1}
\]
for some $V$, which we need because we are assuming only that $\tilde h^G_*$
is defined on spaces.
The evaluation maps are compatible as $G/H$ varies, so define a map
\[
 \Mackey C^{G,\delta}_{\alpha+p}(X,*) \tensor_{\sorb{G,\delta}} \MackeyOp h^{G,\delta}_{\beta+q}
  \to \tilde h^G_{\alpha+\beta+p+q}(X^{\alpha+p}/X^{\alpha+p-1})
\]
which is an isomorphism because $X^{\alpha+p}/X^{\alpha+p-1}$ is a wedge of spheres.

In cohomology, we start with the evaluation map
\begin{align*}
[G_+\smsh_H S^{-\delta(G/H)+\alpha+p},{}& \susp_G^\infty X^{\alpha+p}/X^{\alpha+p-1}]_G
  \tensor \tilde h_G^{\alpha+\beta+p+q}(X^{\alpha+p}/X^{\alpha+p-1}) \\
   &\to \tilde h_G^{\alpha+\beta+p+q}(G_+\smsh_H S^{-\delta(G/H)+\alpha+p}) \\
   &\iso \tilde h_G^{\beta+q}(G_+\smsh_H S^{-\delta(G/H)}) \\
   &= \Mackey h_{G,\delta}^{\beta+q}(G/H,\delta).
\end{align*}
In adjoint form this gives a homomorphism
\[
 \tilde h_G^{\alpha+\beta+p+q}(X^{\alpha+p}/X^{\alpha+p-1})
 \to \Hom(\Mackey C^{G,\delta}_{\alpha+p}(X,*)(G/H,\delta), \Mackey h_{G,\delta}^{\beta+q}(G/H,\delta))
\]
for each $G/H$. These homomorphisms are compatible as $G/H$ varies, giving
a map
\[
 \tilde h_G^{\alpha+\beta+p+q}(X^{\alpha+p}/X^{\alpha+p-1})
  \to \Hom_{\sorb{G,\delta}}(\Mackey C^{G,\delta}_{\alpha+p}(X,*), \Mackey h_{G,\delta}^{\beta+q})
\]
which is again an isomorphism because $X^{\alpha+p}/X^{\alpha+p-1}$ is a wedge of spheres.
\end{proof}

\begin{corollary}[Uniqueness of Cellular $RO(G)$-Graded Homology and Cohomology]
\label{cor:ordinaryUniqueness}
Let $\delta$ be a familial dimension function for $G$.
Let $\tilde h^G_*(-)$ be a reduced $RO(G)$-graded homology theory that
takes $\F(\delta)$-equivalences to isomorphisms and
obeys a $\delta$-dimension axiom, in the sense that
\[
 \MackeyOp h^{G,\delta}_n = 0 \quad\text{for integers $n\neq 0$.}
\]
Then there is a natural isomorphism 
\[
 \tilde h^G_*(-) \iso \tilde H^{G,\delta}_*(-;\MackeyOp h^{G,\delta}_0)
\]
of reduced $RO(G)$-graded homology theories.
Similarly, if $\tilde h_G^*(-)$ is a reduced $RO(G)$-graded cohomology theory
that takes $\F(\delta)$-equivalences to isomorphisms and
satisfies
\[
 \Mackey h_{G,\delta}^n = 0 \quad\text{for integers $n\neq 0$,}
\]
then there is a natural isomorphism
\[
 \tilde h_G^*(-) \iso \tilde H_{G,\delta}^*(-;\Mackey h_{G,\delta}^0)
\]
of reduced $RO(G)$-graded cohomology theories.
\end{corollary}

\begin{proof}
There is an Atiyah-Hirezebruch spectral sequence
\[
 E^2_{p,q} = \tilde H^{G,\delta}_{\alpha+p}(X;\MackeyOp h^{G,\delta}_q) \convto \tilde h^G_{\alpha+p+q}(X).
\]
If $\tilde h^G_*(-)$ obeys a $\delta$-dimension axiom, then the $E^2$ term collapses to the line $q=0$,
giving the isomorphism. 
(Note that $\MackeyOp h^{G,\delta}_*$ involves only those subgroups in $\F(\delta)$.)
The argument for cohomology is the same.
\end{proof}

We also have a universal coefficients spectral sequence. The universal coefficients
in this case are the canonical projectives $\MackeyOp A^{G/K,\delta}$. Write
$\Mackey H^{G,\delta}_\alpha(X)$ for the contravariant Mackey functor given by
\[
 \Mackey H^{G,\delta}_\alpha(X)(G/K,\delta) = \tilde H^{G,\delta}_\alpha(X;\MackeyOp A^{G/K,\delta}).
\]
With this notation, we have
\[
 \Mackey C^{G,\delta}_{\alpha+n}(X,*) \iso \Mackey H^{G,\delta}_{\alpha+n}(X^{\alpha+n}/X^{\alpha+n-1})
\]
for any based $\delta$-$G$-CW($\alpha$) complex $X$.
In the following, $\Tor^{\sorb{G,\delta}}_*$ and $\Ext_{\sorb{G,\delta}}^*$ are the derived functors
of $\tensor_{\sorb{G,\delta}}$ and $\Hom_{\sorb{G,\delta}}$, respectively.

\begin{theorem}[Universal Coefficients Spectral Sequence]\label{thm:ucc}
If $\MackeyOp S$ is a covariant $\delta$-Mackey functor and $\Mackey T$ is a contravariant
$\delta$-Mackey functor, there are spectral sequences
\[
 E^2_{p,q} = \Tor^{\sorb{G,\delta}}_p(\Mackey H^{G,\delta}_{\alpha+q}(X),\MackeyOp S)
  \convto \tilde H^{G,\delta}_{\alpha+p+q}(X;\MackeyOp S)
\]
and
\[
 E_2^{p,q} = \Ext_{\sorb{G,\delta}}^p(\Mackey H^{G,\delta}_{\alpha+q}(X),\Mackey T)
  \convto \tilde H_{G,\delta}^{\alpha+p+q}(X;\Mackey T)
\]
\end{theorem}

\begin{proof}
These spectral sequences are constructed in the usual way by taking
resolutions of $\MackeyOp S$ and $\Mackey T$ and then taking the
spectral sequences of the bigraded complexes obtained from
the chain complex of (an approximation to) $X$ and the resolutions.
\end{proof}

\section{Ordinary and dual homology and cohomology}

We now briefly give names to and discuss the specializations of cellular homology and cohomology to
particular choices of the dimension function $\delta$.

\begin{definition}
If we set $\delta=0$ in cellular homology and cohomology, we call the resulting
$RO(G)$-graded theories
{\em ordinary homology and cohomology}. We use the notations 
\[
 \tilde H^G_*(X;\MackeyOp S) = \tilde H^{G,0}_*(X;\MackeyOp S)
\]
and
\[
 \tilde H_G^*(X;\Mackey T) = \tilde H_{G,0}^*(X;\Mackey T)
\]
and similarly for the unreduced theories, where $\MackeyOp S$ is
a covariant ($\sorb G$-)Mackey functor and $\Mackey T$ is a contravariant Mackey functor.
\end{definition}

The ordinary theories are the ones constructed using ordinary $G$-CW($\alpha$) complexes,
in which the orbits do not contribute to dimension.

The ordinary theories are characterized by the following dimension axioms:
\[
 H^G_{k}(G/H;\MackeyOp S)
  \iso
   \begin{cases}
    \MackeyOp S(G/H) & \text{if $k = 0$} \\
    0 & \text{if $k\neq 0$}
   \end{cases}
\]
and
\[
 H_G^{k}(G/H;\Mackey T)
  \iso
   \begin{cases}
    \Mackey T(G/H) & \text{if $k = 0$} \\
    0 & \text{if $k\neq 0$}
   \end{cases}
\]
for $k$ an integer. Thus, the ordinary theories constructed here coincide with the theories
of the same name in \cite{LMM:roghomology} and \cite{May:alaska}.

Turning to the case $\delta=\Lie$, recall from Corollary~\ref{cor:MackeyDuality} that
a covariant $\Lie$-Mackey functor is the same thing as a contravariant Mackey functor
and vice versa.

\begin{definition}
If we set $\delta=\Lie$ in cellular homology and cohomology, we call the resulting
$RO(G)$-graded theories {\em dual homology and cohomology}. We use the notations
\[
 \tilde \H^G_*(X;\Mackey T) = \tilde H^{G,\Lie}_*(X;\Mackey T)
\]
and
\[
 \tilde \H_G^*(X;\MackeyOp S) = \tilde H_{G,\Lie}^*(X;\MackeyOp S)
\]
and similarly for the unreduced theories, where $\MackeyOp S$ is
a covariant (0-)Mackey functor
(hence a contravariant $\Lie$-Mackey functor)
and $\Mackey T$ is a contravariant Mackey functor.
\end{definition}

The dual theories are the ones constructed using dual $G$-CW($\alpha$) complexes,
in which an orbit $G/H$ contributes its full geometric dimension, $\Lie(G/H)$.
Thus, the dimension of a cell $G\times_H D(V)$ is $V + \Lie(G/H)$ in the context of
the dual theories.

The reason for the name will become clear in the following section but is suggested by
the following characterization. Recall that the stable dual of the orbit $G/H$ is the spectrum
\[
 D(G/H_+) = G_+\smsh_H S^{-\Lie(G/H)}.
\]
In the following, we write as usual
\[
 \tilde h^G_\alpha(D(G/H_+)) = \tilde h^G_\alpha(G_+\smsh_H S^{-\Lie(G/H)})
  = \tilde h^G_{\alpha+V}(G_+\smsh_H S^{V-\Lie(G/H)})
\]
where $V$ is any representation large enough that $\Lie(G/H)\subset V$.
The dual theories are then characterized by the following dimension axioms:
\[
 \tilde\H^G_{k}(D(G/H_+);\Mackey T)
  \iso
   \begin{cases}
    \Mackey T(G/H) & \text{if $k = 0$} \\
    0 & \text{if $k\neq 0$}
   \end{cases}
\]
and
\[
 \tilde\H_G^{k}(D(G/H_+);\MackeyOp S)
  \iso
   \begin{cases}
    \MackeyOp S(G/H) & \text{if $k = 0$} \\
    0 & \text{if $k\neq 0$}
   \end{cases}
\]
for $k$ an integer. (Notice that the variances match correctly due to the use of duality.)

\section{The representing spectra}\label{sec:representing}

In this section we discuss the equivariant spectra representing cellular
homology and cohomology of spaces.
We can use any of several available models of equivariant spectra,
for example the orthogonal spectra of \cite{MM:orthogonal} or the older
model expounded in \cite{LMS:eqhomotopy}; all models give equivalent stable categories.
(Note that we do not attempt to connect our use of spectra here with
our use of $G$-CW prespectra in defining the chains of a $G$-space.)
We write $G\Spec{}{}$ for one of these categories of $G$-spectra.
Recall that a $G$-spectrum $E$ represents an $RO(G)$-graded cohomology theory
$\tilde h_G^*(-)$ defined on based $G$-spaces if there is a natural isomorphism of reduced theories
\[
 \tilde h_G^\alpha(X) \iso [\susp_G^\infty X, \susp_G^\alpha E]_G.
\]
Similarly, $E$ represents an $RO(G)$-grade
homology theory $\tilde h^G_*(-)$ if there is a natural isomorphism
\[
 \tilde h^G_\alpha(X) \iso [S^\alpha, E\smsh X]_G.
\]
(Here and elsewhere, $S^\alpha$ or $S^V$ denotes either a space or a spectrum, as appropriate
from the context.)
Given any $RO(G)$-graded cohomology theory $\tilde h_G^*(-)$, there exists
a $G$-spectrum representing $\tilde h_G^*(-)$, which is unique up to non-unique stable equivalence,
and similarly for homology (\cite[XIII.3]{May:alaska}).

In particular, for each familial dimension function $\delta$ and each
contravariant $\delta$-Mackey functor $\Mackey T$, there is a
$G$-spectrum $H_\delta \Mackey T$ representing 
cellular cohomology, $\tilde H_{G,\delta}^*(-;\Mackey T)$.
Similarly, if $\MackeyOp S$ is a covariant $\delta$-Mackey functor,
there is a $G$-spectrum $H^\delta\MackeyOp S$ representing cellular
homology, $\tilde H^{G,\delta}_*(-;\MackeyOp S)$.
(This establishes only that $H_\delta\Mackey T$ is unique up to
non-unique stable equivalence. In fact, the stable equivalence is unique as well.
This was shown in \cite{May:alaska} for $H_0\Mackey T$ and we will return to this point
in \S\ref{sec:OrdinaryRemarks}.)
Explicit constructions of the spectrum $H\Mackey T = H_0\Mackey T$ have
previously been announced or given in several places, including
\cite{LMM:roghomology} (announcement),
\cite{CW:loopspaces} (for finite $G$ only), \cite{Le:emspaces},
\cite{GM:eqhomotopy}, and \cite{May:alaska}.
Below we shall give constructions of $H_\delta\Mackey T$ and $H^\delta\MackeyOp S$,
but first we discuss how these spectra are characterized.

By Corollary~\ref{cor:ordinaryUniqueness} and the uniqueness of representing spectra,
$H_\delta\Mackey T$ is characterized up to stable equivalence by its stable homotopy, in the following
sense.
For any $G$-spectrum $E$, we write $\Mackey\pi^{G,\delta}_n(E)$ for the $\delta$-Mackey functor defined by
\[
 \Mackey\pi^{G,\delta}_n(E)(G/H,\delta) = [G_+\smsh_H S^{-\delta(G/H)+n}, E]_G.
\]
With this notation, $H_\delta\Mackey T$ is characterized by
 \[
 \Mackey\pi^{G,\delta}_n(H_\delta\Mackey T) \iso
 \begin{cases}  \Mackey T & \text{if $n = 0$} \\
                0 & \text{if $n\neq 0$} \end{cases}
 \]
for integers $n$, together with the fact that the cohomology theory it represents takes
$\F(\delta)$-equivalences to isomorphisms. The latter property can be expressed as the
fact that the map $H_\delta\Mackey T \to F(E\F(\delta)_+,H_\delta\Mackey T)$
is a stable equivalence, where $F(-,-)$ denotes the mapping spectrum.
Note the strength of this characterization: Given any two $G$-spectra with these properties,
there exists a stable equivalence between them.
We shall call such a spectrum a {\em $\delta$-Eilenberg-Mac\,Lane
spectrum} of type $\Mackey T$.

There is a similar characterization of $H^\delta\MackeyOp S$ but with some crucial changes.
By Corollary~\ref{cor:ordinaryUniqueness}, the relevant homotopy groups are
\begin{align*}
 [S^n, H^\delta\MackeyOp S \smsh G_+\smsh_H S^{-\delta(G/H)}]_G
  &\iso [D(G_+\smsh_H S^{-\delta(G/H)})\smsh S^n, H^\delta\MackeyOp S]_G \\
  &\iso [G_+\smsh_H S^{-(\Lie(G/H) - \delta(G/H)) + n}, H^\delta\MackeyOp S]_G \\
  &\iso \Mackey \pi^{G,\Lie-\delta}_n(H^\delta\MackeyOp S)(G/H),
\end{align*}
where we use that the stable dual of an orbit $G/H_+$ is $G_+\smsh_H S^{-\Lie(G/H)}$.
Thus, we require that
\[
 \Mackey\pi^{G,\Lie-\delta}_n(H^\delta\MackeyOp S) \iso
 \begin{cases}  \MackeyOp S & \text{if $n = 0$} \\
                0 & \text{if $n\neq 0$} \end{cases}
\]
(using that $\MackeyOp S$ may be considered a contravariant $(\Lie-\delta)$-Mackey functor).
The other requirement, corresponding to the fact that the represented homology
theory must take $\F(\delta)$-equivalences to isomorphisms, becomes the fact that
$H^\delta\MackeyOp S\smsh E\F(\delta)_+ \to H^\delta\MackeyOp S$ is a stable equivalence.

It has been noted before, in \cite{May:alaska} for example, that, for infinite groups $G$,
ordinary homology and cohomology are not dual, i.e., they are not represented
by the same spectrum. 
The characterizations above establish the following result, which shows
how the various cellular homology and cohomology theories are related.

\begin{theorem}\label{thm:dualCohomologyRep}
Let $\delta$ be a familial dimension function for $G$ with dual $\Lie-\delta$.
Let $\Mackey T$ be a contravariant $\delta$-Mackey functor, which we may
also consider as a covariant $(\Lie-\delta)$-Mackey functor.
Then there is a $G$-equivalence
\[
 H_\delta\Mackey T\smsh E\F(\delta)_+ \hmtpc H^{\Lie-\delta}\Mackey T
\]
and a $G$-equivalence
\[
 H_\delta\Mackey T \hmtpc F(E\F(\delta)_+,H^{\Lie-\delta}\Mackey T).
\]
\qed\end{theorem}

Thus, when $\delta$ is complete (so $E\F(\delta)_+ \hmtpc S^0$),
we have $H_\delta\Mackey T \hmtpc H^{\Lie-\delta}\Mackey T$, so
$\tilde H_{G,\delta}^*(-;\Mackey T)$ and $\tilde H^{G,\Lie-\delta}_*(-;\Mackey T)$
are dual theories, in the sense that they are represented by the same spectrum.
In particular, we have the following.

\begin{corollary}
The following are pairs of dual theories.
\begin{enumerate}
\item
$\tilde H_G^*(-;\Mackey T)$ and $\tilde\H^G_*(-;\Mackey T)$ are both represented by 
$H\Mackey T = H_0\Mackey T \hmtpc H^\Lie \Mackey T$.

\item
$\tilde H^G_*(-;\MackeyOp S)$ and $\tilde\H_G^*(-;\MackeyOp S)$ are both represented by
$H\MackeyOp S = H_\Lie\MackeyOp S \hmtpc H^0\MackeyOp S$.
\qed
\end{enumerate}
\end{corollary}

\begin{note}
If $\delta$ is not complete, then, in general, $H_\delta\Mackey T$ and $H^{\Lie-\delta}\Mackey T$
need not be equivalent.
Consider a simple example:
Let $G = \Z/2$ and let $\delta$ be the dimension
function that is 0 (as it must be) on $\F(\delta) = \{e\}$.
Let $\Mackey T$ be the $\delta$-Mackey functor that is constant at $\Z/2$.
Then, in integer dimensions we have
\begin{align*}
 H_{G,\delta}^*(G/G;\Mackey T) \iso H_{G,\delta}^*(EG;\Mackey T) &\iso H^*(B\Z/2;\Z/2) \\
\intertext{and}
 H_{G,\delta}^k(G/e;\Mackey T) &\iso
  \begin{cases}
   \Z/2 & \text{if $k=0$} \\
   0 & \text{otherwise.}
  \end{cases}
\end{align*}
A similar calculation shows that
\begin{align*}
 H^{G,\delta}_*(G/G;\Mackey T) &\iso H_*(B\Z/2;\Z/2) \\
\intertext{and}
 H^{G,\delta}_k(G/e;\Mackey T) &\iso
  \begin{cases}
   \Z/2 & \text{if $k=0$} \\
   0 & \text{otherwise,}
  \end{cases}
\end{align*}
where we can consider $\Mackey T$ as both a contravariant and a covariant
$\delta$-Mackey functor because $G$ is finite.
Looking at the values at $G/G$ these two theories are clearly not dual.
(Note that they are, in fact, Borel homology and cohomology with $\Z/2$ coefficients.)
This is related to the observation in \cite[V.8]{LMS:eqhomotopy} that
there are two distinct ways to restrict cohomology theories to (pairs of) families.
In the present case, $\Mackey T$ extends to a $G$-Mackey functor with
an associated Eilenberg-Mac\,Lane spectrum $H\Mackey T$
(for example, $\Mackey T = \Mackey A_{G/G}\tensor\Z/2$).
The spectrum representing $H^{G,\delta}_*(-;\Mackey T)$ is then
$H\Mackey T\smsh E\Z/2_+$ while the spectrum representing
$H_{G,\delta}^*(-;\Mackey T)$ is the function spectrum $F(E\Z/2_+, H\Mackey T)$.
\end{note}

Finally, we give an elementary, functorial, explicit construction of
$H_\delta\Mackey T$.
This particular structure will be of use later.

\begin{construction}\label{con:EilenbergMacLane}
Let $\delta$ be a familial dimension function for $G$ and let
$\Mackey T$ be a contravariant $\delta$-Mackey functor.
Define
\[
 F_\delta\Mackey T = \Wedge_{\Mackey T(G/H,\delta)} G_+\smsh_H S^{-\delta(G/H)},
\]
where the wedge runs over all objects $(G/H,\delta)$ in $\sorb{G,\delta}$
and all elements in $\Mackey T(G/H,\delta)$.
Then $F_\delta$ is a functor taking $\delta$-$G$-Mackey functors to $G$-spectra 
and we have a natural epimorphism
\[
 \epsilon\colon \Mackey\pi_0^{G,\delta} F_\delta\Mackey T \to \Mackey T.
\]
Let
\[
 K = \{ \kappa\colon G_+\smsh_{H_\kappa} S^{-\delta(G/H_\kappa)}\to F_\delta\Mackey T
       \mid \epsilon(\kappa) = 0 \},
\]
and let
\[
 R_\delta\Mackey T = \Wedge_{\kappa\in K} G_+\smsh_{H_\kappa} S^{-\delta(G/H_\kappa)}.
\]
Then $R_\delta$ is also a functor,
there is a natural transformation $R_\delta\Mackey T\to F_\delta\Mackey T$
(given by the maps $\kappa$), and we have an
exact sequence
\[
 \Mackey\pi_0^{G,\delta}R_\delta\Mackey T \to \Mackey\pi_0^{G,\delta} F_\delta\Mackey T \to \Mackey T \to 0.
\]
It follows that, if we let $C_\delta\Mackey T$ be the cofiber of $R_\delta\Mackey T\to F_\delta\Mackey T$,
then $\Mackey\pi_0^{G,\delta}C_\delta\Mackey T \iso \Mackey T$
and $\Mackey\pi_n^{G,\delta}C_\delta\Mackey T = 0$ for $n<0$.
The vanishing of the homotopy groups for $n<0$ depends on the calculation
\[
 [G_+\smsh_H S^{-\delta(G/H)+n}, G_+\smsh_K S^{-\delta(G/K)}]_G = 0
\]
for $n<0$. This calculation follows from our being able to write this stable mapping group as
a colimit of maps of spaces under suspension (using compactness) 
and then using the same argument as was used in
the proof of Proposition~\ref{prop:SkelEquivalence}.

We can then functorially kill all the homotopy $\Mackey\pi_n^{G,\delta}C_\delta\Mackey T$ for $n>0$, obtaining
a functor $P_\delta$ with
\[
 \Mackey\pi_n^{G,\delta}P_\delta\Mackey T =
   \begin{cases}
      \Mackey T &\text{if $n=0$} \\
      0 &\text{if $n\neq 0$.}
   \end{cases}
\]
Finally, we let
\[
 H_\delta\Mackey T = F(E\F(\delta)_+,P_\delta\Mackey T)
\]
and
\[
 H^\delta\MackeyOp S = P_{\Lie-\delta}\MackeyOp S \smsh E\F(\delta)_+
\]
if $\MackeyOp S$ is a covariant $\delta$-Mackey functor.
Then $H_\delta$ and $H^\delta$ are functors and produce spectra that
satisfy the characterizations of the Eilenberg-Mac\,Lane spectra representing
cohomology and homology, respectively.
\qed
\end{construction}

\section{Change of groups}\label{sec:vsubgroups}

We now discuss various change-of-groups maps.
We first discuss restriction to subgroups, where the underlying result is
the Wirthm\"uller isomorphism.

\subsection{Subgroups}

Given a $G$-homology theory $\tilde h^G_*(-)$ and a subgroup $K$ of $G$,
$\tilde h^G_*(G_+\smsh_K -)$ defines a homology theory on $K$-spaces,
and similarly for cohomology.
When we apply this construction to our cellular theories,
the Wirthm\"uller isomorphisms identify the resulting theories
as again cellular.
This can be shown by calculating 
$\tilde H_G^*(G_+\smsh_K K/L;\Mackey T)$ and using the dimension axiom,
but it is useful to see that the isomorphism starts at the chain level.
Later we shall look at how it is represented on the spectrum level.

Continue to let $K$ be a subgroup of $G$.
Recall from \cite[II.4]{LMS:eqhomotopy} or \cite[V.2]{MM:orthogonal} that
there is a forgetful functor $E\mapsto E|K$ taking $G$-spectra to $K$-spectra, which has
a left adjoint denoted $G_+\smsh_K -$ 
(we use the notation from \cite{MM:orthogonal}).
This adjunction passes to the stable category to give the adjunction
\[
 [D, E|K]_K \iso [G_+\smsh_K D, E]_G
\]
for $K$-spectra $D$ and $G$-spectra $E$.
When the meaning is clear from context, we shall simply write $E$ for $E|K$,
$E$ regarded as a $K$-spectrum.
We also have that $G_+\smsh_K\susp_K^\infty X = \susp_G^\infty (G_+\smsh_K X)$,
relating the spectrum- and space-level functors.

On the algebraic side,
suppose that $\delta$ is a dimension function for $G$ and that $K\in\F(\delta)$.
Write $\delta$ again for $\delta|K$ as defined in Definition~\ref{def:dimensionRestriction}.

\begin{definition}\label{def:algRestrictToSubgroup}
Suppose $K\in\F(\delta)$.
\begin{enumerate}
\item
Let $i_K^G\colon \sorb{K,\delta}\to \sorb{G,\delta}$ be
the functor $G_+\smsh_K\susp^{-\delta(G/K)}(-)$.
Note that
\[
 G_+\smsh_K\susp^{-\delta(G/K)}(K_+\smsh_L S^{-\delta(K/L)})
  \iso G_+\smsh_L S^{-\delta(G/L)},
\]
using that $\delta(G/K) \dirsum \delta(K/L) = \delta(G/L)$.

\item
If $\Mackey T$ is a contravariant $G$-$\delta$-Mackey functor, let
\[
 \Mackey T|K = (i_K^G)^*\Mackey T = \Mackey T\circ i_K^G,
\]
a contravariant $K$-$\delta$-Mackey functor. If $\MackeyOp S$ is a covariant $G$-$\delta$-Mackey functor,
let $\MackeyOp S|K = (i_K^G)^*\MackeyOp S$ similarly.

\item
If $\Mackey C$ is a contravariant $K$-$\delta$-Mackey functor, let
\[
 G\times_K\Mackey C = (i_K^G)_!\Mackey C
\]
as in Definition~\ref{def:indres}, a contravariant $G$-$\delta$-Mackey functor.
If $\MackeyOp D$ is a covariant $K$-$\delta$-Mackey functor, define
$G\times_K\MackeyOp D = (i_K^G)_!\MackeyOp D$ similarly.
\end{enumerate}
\end{definition}

It follows from Proposition~\ref{prop:indresadjunction} that,
if $\Mackey C$ is a contravariant $K$-$\delta$-Mackey functor
and $\Mackey T$ is a contravariant $G$-$\delta$-Mackey functor, then
 \[
 \Hom_{\sorb{G,\delta}}(G\times_K \Mackey C, \Mackey T)
 \iso \Hom_{\sorb{K,\delta}}(\Mackey C, \Mackey T|K).
 \]
If $\MackeyOp S$ is a covariant $G$-$\delta$-Mackey functor and
$\MackeyOp D$ is a covariant $K$-$\delta$-Mackey functor,
we also have the isomorphism
 \[
 (G\times_K \Mackey C)\tensor_{\sorb{G,\delta}} \MackeyOp S
 \iso \Mackey C\tensor_{\sorb{K,\delta}} (\MackeyOp S|K).
 \]
Proposition~\ref{prop:indresadjunction} gives us the calculation
 \[
 G\times_K \Mackey A_{K/J,\delta} \iso \Mackey A_{G/J,\delta}.
 \]
The restrictions of the functors $\Mackey A_{G/J,\delta}$ are complicated in general,
but we do have the following simple special cases:
\[
 \Mackey A_{G/G,0}|K \iso \Mackey A_{K/K,0}
\]
and
\[
 \MackeyOp A^{G/G,\Lie}|K \iso \MackeyOp A^{K/K,\Lie}.
\]

The following is the main calculation that leads to the Wirthm\"uller isomorphisms.

\begin{proposition}\label{prop:wirthmulleronchains}
Let $\delta$ be a dimension function for $G$,
let $K\in\F(\delta)$,
and let $\alpha\in RO(G)$.
For notational simplicity, write $\alpha$ again for $\alpha|K$.
Let $X$ be a based $\delta$-$K$-CW$(\alpha-\delta(G/K))$ complex and give
$G_+\smsh_K X$ the $\delta$-$G$-CW$(\alpha)$ structure from Proposition~\ref{prop:genInduction}.
Then
\[
 G\times_K\Mackey C^{K,\delta}_{\alpha-\delta(G/K)+*}(X,*) 
   \iso \Mackey C^{G,\delta}_{\alpha+*}(G_+\smsh_K X,*).
\]
This isomorphism respects suspension in the sense that, if $W$ is a representation of $G$,
then the following diagram commutes:
\[
 \xymatrix{
  G\times_K\Mackey C^{K,\delta}_{\alpha-\delta(G/K)+*}(X,*) \ar[r]^-\iso \ar[d]_{\sigma^W}
   & \Mackey C^{G,\delta}_{\alpha+*}(G_+\smsh_K X,*) \ar[d]^{\sigma^W} \\
  G\times_K\Mackey C^{K,\delta}_{\alpha-\delta(G/K)+W+*}(\susp^W X,*) \ar[r]^-\iso
   & \Mackey C^{G,\delta}_{\alpha+W+*}(G_+\smsh_K \susp^W X,*)
 }
\]  
\end{proposition}

\begin{proof}
Geometrically, it's clear that we have a one-to-one correspondence between the cells
in the two cell complexes.
To see that the algebra tracks the geometry, we define a map 
 \[
  G\times_K\Mackey C^{K,\delta}_{\alpha-\delta(G/K)+*}(X,*) \to \Mackey C^{G,\delta}_{\alpha+*}(G_+\smsh_K X,*)
 \]
as follows. 
By Definition~\ref{def:indres}, 
\begin{multline*}
 (G\times_K\Mackey C^{K,\delta}_{\alpha-\delta(G/K)+*}(X,*))(G/J,\delta) \\
    = \int^{(K/L,\delta)} \Mackey C^{K,\delta}_{\alpha-\delta(G/K)+*}(X,*)(K/L)\tensor \sorb{G,\delta}(G/J,G\times_K K/L).
\end{multline*}
We have the extension and evaluation map
\begin{multline*}
 \Mackey C^{K,\delta}_{\alpha-\delta(G/K)+n}(X,*)(K/L)\tensor \sorb{G,\delta}(G/J,G\times_K K/L) \\
 \to \Mackey C^{G,\delta}_{\alpha+n}(G_+\smsh_K X,*)(G/J)
\end{multline*}
defined as follows.
Let $c\in \Mackey C^{K,\delta}_{\alpha-\delta(G/K)+n}(X,*)(K/L)$, so $c$ is a stable map
\begin{multline*}
 c\colon K_+\smsh_L S^{\alpha-\delta(G/K)-\delta(K/L)+n}
 = K_+\smsh_L S^{\alpha-\delta(G/L)+n} \\
 \to \susp_K^\infty X^{\alpha-\delta(G/K)+n}/X^{\alpha-\delta(G/K)+n-1},
\end{multline*}
and let $d\in \sorb{G,\delta}(G/J,G\times_K K/L)$, so $d$ is a stable map
\[
 d\colon G_+\smsh_J S^{-\delta(G/J)} \to 
 G_+\smsh_K K_+\smsh_L S^{-\delta(G/L)}.
\]
Then the extension and evaluation map takes $c\tensor d$
to the stable map
\begin{multline*}
 (G_+\smsh_K c)\circ d \colon
  G_+\smsh_J S^{\alpha-\delta(G/J)+n} \\
  \to
  \susp_G^\infty (G_+\smsh_K X^{\alpha-\delta(G/K)+n})/(G_+\smsh_K X^{\alpha-\delta(G/K)+n-1}).
\end{multline*}
As mentioned above, we can take $G_+\smsh_K X^{\alpha-\delta(G/K)+n}$ as the
$(\alpha+n)$-skeleton of a $G$-$\delta$-CW($\alpha$) structure on $G_+\smsh_K X$, i.e.,
\[
 G_+\smsh_K X^{\alpha-\delta(G/K)+n} = (G_+\smsh_K X)^{\alpha+n}.
\]
Thus, $(G_+\smsh_K c)\circ d$ defines an element of
$\Mackey C^{G,\delta}_{\alpha+n}(G_+\smsh_K X,*)(G/J)$ as claimed.
The extension and evaluation maps are compatible as $K/L$ varies, defining a map out of the coend.
The isomorphism $G\times_K \Mackey A_{K/J,\delta} \iso \Mackey A_{G/J,\delta}$
then implies that the map so defined is an isomorphism
\[
 G\times_K\Mackey C^{K,\delta}_{\alpha-\delta(G/K)+n}(X,*) \iso \Mackey C^{G,\delta}_{\alpha+n}(G_+\smsh_K X,*).
\]

That $G_+\smsh_K\susp^W X \iso \susp^W G_+\smsh_K X$
implies that the isomorphism respects suspension.
\end{proof}

\begin{theorem}[Wirthm\"uller Isomorphisms]\label{thm:wirthmuller}
Let $\delta$ be a dimension function for $G$,
let $K\in\F(\delta)$,
and let $\alpha\in RO(G)$.
For notational simplicity, write $\alpha$ again for $\alpha|K$.
Then, for $X$ in $K\W_*^{\delta,\alpha-\delta(G/K)}$ there are natural isomorphisms
 \[
 \tCH^{G,\delta}_{\alpha}(G_+\smsh_K X;\MackeyOp S)
 \iso  \tCH^{K,\delta}_{\alpha - \delta(G/K)}(X;\MackeyOp S|K)
 \]
and
 \[
 \tCH_{G,\delta}^{\alpha}(G_+\smsh_K X;\Mackey T)
 \iso  \tCH_{K,\delta}^{\alpha-\delta(G/K)}(X;\Mackey T|K).
 \]
These isomorphisms respect suspension in the sense that, if $W$ is a representation of $G$,
then the following diagram commutes:
\[
 \xymatrix{
  \tCH^{G,\delta}_{\alpha}(G_+\smsh_K X;\MackeyOp S) \ar[r]^-{\iso} \ar[d]_{\sigma^W}
    & \tCH^{K,\delta}_{\alpha - \delta(G/K)}(X;\MackeyOp S|K) \ar[d]^{\sigma^W} \\
  \tCH^{G,\delta}_{\alpha+W}(G_+\smsh_K \susp^W X;\MackeyOp S) \ar[r]^-{\iso}
    & \tCH^{K,\delta}_{\alpha + W - \delta(G/K)}(\susp^W X;\MackeyOp S|K)
 }
\]
and similarly for cohomology.
If $\delta$ is familial, for $X$ in $K\K_*$ there are natural isomorphisms
 \[
 \tilde H^{G,\delta}_{\alpha}(G_+\smsh_K X;\MackeyOp S)
 \iso  \tilde H^{K,\delta}_{\alpha - \delta(G/K)}(X;\MackeyOp S|K)
 \]
and
 \[
 \tilde H_{G,\delta}^{\alpha}(G_+\smsh_K X;\Mackey T)
 \iso  \tilde H_{K,\delta}^{\alpha-\delta(G/K)}(X;\Mackey T|K).
 \]
These isomorphisms respect suspension in the sense above.
 \end{theorem}

\begin{proof}
Consider first the case where 
$X$ is a $K$-$\delta$-CW($\alpha-\delta(G/K)$) complex.
Using the isomorphism of the preceding proposition and
isomorphisms we noted earlier coming from Proposition~\ref{prop:indresadjunction}, 
we see that there are chain isomorphisms
 \[
 \Mackey C^{G,\delta}_{\alpha+*}(G_+\smsh_K X,*)\tensor_{\sorb{G,\delta}} \MackeyOp S
 \iso \Mackey C^{K,\delta}_{\alpha-\delta(G/K)+*}(X,*)\tensor_{\sorb{K,\delta}} (\MackeyOp S|K)
 \]
and
 \[
 \Hom_{\sorb{G,\delta}}(\Mackey C^{G,\delta}_{\alpha+*}(G_+\smsh_K X,*), \Mackey T)
 \iso \Hom_{\sorb{K,\delta}}(\Mackey C^{K,\delta}_{\alpha-\delta(G/K)+*}(X,*), \Mackey T|K).
 \]
These isomorphisms induce the Wirthm\"uller isomorphisms.
That the Wirthm\"uller isomorphisms respect suspension follows from the similar statement
in the preceding proposition.

Now consider a general based $K$-space $X$, which we may assume is well-based.
Let $\U$ be a complete $G$-universe, which will also be complete when considered
as a $K$-universe. Let $\V$ be an indexing sequence in $\U$.
Let $\Gamma\susp_K^\infty X\to \susp_K^\infty X$ be a
$\delta$-$K$-CW$(\alpha-\delta(G/K))$ approximation indexed on $\V$.
We can then construct a $\delta$-$G$-CW$(\alpha)$ approximation
$\Gamma\susp_G^\infty X\to \susp_G^\infty X$ indexed on $\V$ such that there are inclusions
\[
 G_+\smsh_K (\Gamma\susp_K^\infty X)(V_i) \to (\Gamma\susp_G^\infty (G_+\smsh_K X))(V_i)
\]
over $\susp^{V_i}(G_+\smsh X)$ and compatible with the structure maps.
For sufficiently large $i$, $(\Gamma\susp_K^\infty X)(V_i) \to \susp^{V_i} X$
is an $\F(\delta)$-weak equivalence of $K$-spaces. It follows that
$G_+\smsh_K (\Gamma\susp_K^\infty X)(V_i) \to \susp^{V_i} (G_+\smsh_K X)$
is an $\F(\delta)$-weak equivalence of $G$-spaces for such $i$,
hence a $\delta$-weak$_{\alpha+V_i}$ equivalence.
Using the preceding proposition, we then get a chain homotopy equivalence
\[
 G\times_K \Mackey C_{\alpha-\delta(G/K)+*}^{\V,\delta}(X) \hmtpc
  \Mackey C_{\alpha+*}^{\V,\delta}(G_+\smsh_K X).
\]
The theorem now follows as it did in the CW case.
 \end{proof}

Recall from Definition~\ref{def:genMackeyvalued} that we can consider any
$RO(G)$-graded theory as $G$-$\delta$-Mackey functor-valued.
Applying that definition to cellular homology and cohomology we get, for $K\in\F(\delta)$,
\begin{align*}
 \MackeyOp H^{G,\delta}_\alpha(X;\MackeyOp S)(G/K,\delta)
  &= \tilde H^{G,\delta}_\alpha(X\smsh G_+\smsh_K S^{-\delta(G/K)};\Mackey S) \\
  &= \tilde H^{G,\delta}_{\alpha+V}(X\smsh G_+\smsh_K S^{V-\delta(G/K)};\Mackey S)
\end{align*}
for $V$ so large that $\delta(G/K)\subset V$, and
\[
 \Mackey H_{G,\delta}^\alpha(X;\Mackey T)(G/K,\delta)
  = \tilde H_{G,\delta}^\alpha(X\smsh G_+\smsh_K S^{-\delta(G/K)};\Mackey T)
\]
similarly.
We can use the Wirthm\"uller isomorphisms to interpret the components of these functors.

\begin{corollary}\label{cor:MackeyStructure}
If $X$ is a $G$-space, $\delta$ is a familial dimension function for $G$,
$\MackeyOp S$ is a covariant $G$-$\delta$-Mackey functor, and
$\Mackey T$ is a contravariant $G$-$\delta$-Mackey functor,
then, for $K\in\F(\delta)$, we have
\[
 \MackeyOp H^{G,\delta}_\alpha(X;\MackeyOp S)(G/K,\delta) \iso
  \tilde H^{K,\delta}_\alpha(X;\MackeyOp S|K)
\]
and
\[
 \Mackey H_{G,\delta}^\alpha(X;\Mackey T)(G/K,\delta) \iso
  \tilde H_{K,\delta}^\alpha(X;\Mackey T|K).
\]
\end{corollary}

\begin{proof}
This is simply a restatement of the Wirthm\"uller isomorphisms. Explicitly,
\begin{align*}
 \MackeyOp H^{G,\delta}_\alpha(X;\MackeyOp S)(G/K,\delta)
  &= \tilde H^{G,\delta}_\alpha(X\smsh G_+\smsh_K S^{-\delta(G/K)};\Mackey S) \\
  &= \tilde H^{G,\delta}_{\alpha+V}(X\smsh G_+\smsh_K S^{V-\delta(G/K)};\Mackey S) \\
  &\iso \tilde H^{K,\delta}_{\alpha+V-\delta(G/K)}(X\smsh S^{V-\delta(G/K)}; \Mackey S|K) \\
  &\iso \tilde H^{K,\delta}_\alpha(X; \Mackey S|K)
\end{align*}
and similarly for cohomology.
\end{proof}

Note that this is {\em not} the structure used in Theorem~\ref{thm:ucc}, 
the universal coefficients theorem.

Another application of the Wirthm\"uller isomorphism,
or inspection of the preceding corollary, gives the following.

\begin{corollary}
Let $\delta$ be a familial dimension function for $G$ and let $K\in\F(\delta)$.
If $X$ is a $G$-space,
$\MackeyOp S$ is a covariant $G$-$\delta$-Mackey functor, and
$\Mackey T$ is a contravariant $G$-$\delta$-Mackey functor,
then we have
\[
 \MackeyOp H^{G,\delta}_\alpha(X;\MackeyOp S)|K \iso
  \MackeyOp H^{K,\delta}_\alpha(X;\MackeyOp S|K)
\]
and
\[
 \Mackey H_{G,\delta}^\alpha(X;\Mackey T)|K \iso
  \Mackey H_{K,\delta}^\alpha(X;\Mackey T|K).
\]
\qed
\end{corollary}

We now look at how the Wirthm\"uller isomorphisms are represented. The main result
we need is the following.

\begin{proposition}\label{prop:emsubgroups}
Let $\delta$ be a familial dimension function for $G$,
let $H_\delta\Mackey T$ be a $G$-$\delta$-Eilenberg-Mac\,Lane spectrum of type $\Mackey T$, and
let $K\in\F(\delta)$. 
Then
\[ 
 (H_\delta\Mackey T)|K \hmtpc \susp_K^{-\delta(G/K)}H_\delta(\Mackey T|K).
\]
Similarly, for a covariant $\MackeyOp S$, we have
\[
 (H^\delta\MackeyOp S)|K \hmtpc \susp_K^{-(\Lie(G/K)-\delta(G/K))} H^\delta(\MackeyOp S|K).
\]
\end{proposition}

\begin{proof}
We verify that $\susp_K^{\delta(G/K)}(H_\delta\Mackey T)|K$ has the homotopy
that characterizes $H_\delta(\Mackey T|K)$.
If $L$ is a subgroup of $K$ we have
\begin{align*}
 \Mackey\pi^{K,\delta}_n(\susp_K^{\delta(G/K)}(H_\delta\Mackey T)|K)(K/L,\delta)
   &=[K_+\smsh_L S^{-\delta(K/L)+n}, \susp_K^{\delta(G/K)}(H_\delta\Mackey T)|K]_K \\
   &\iso[K_+\smsh_L S^{-\delta(K/L)-\delta(G/K)+n}, (H_\delta\Mackey T)|K]_K \\
   &\iso[K_+\smsh_L S^{-\delta(G/L)+n}, (H_\delta\Mackey T)|K]_K \\
   &\iso[G_+\smsh_L S^{-\delta(G/L)+n}, H_\delta\Mackey T]_G \\
   &\iso \Mackey\pi_n^{G,\delta}(H_\delta\Mackey T)(G/L,\delta) \\
   &\iso \begin{cases} \Mackey T(G/L,\delta) & \text{if $n = 0$} \\
                       0 & \text{if $n \neq 0$.}
         \end{cases}
 \end{align*}
Further, $\F(\delta|K)$ contains all the subgroups of $K$ because $\delta$ is familial,
so the homotopy above does characterize $H_\delta(\Mackey T|K)$.

The argument for $H^\delta\MackeyOp S$ is similar.
\end{proof}

The cohomology Wirthm\"uller isomorphism is then the adjunction
\begin{align*}
 [G_+\smsh_K \susp_K^\infty X,\susp_G^\alpha H_\delta\Mackey T]_G 
  &\iso [\susp_K^\infty X,\susp_K^\alpha(H_\delta\Mackey T)|K]_K \\
  &\iso [\susp_K^\infty X,\susp_K^{\alpha-\delta(G/K)}H_\delta(\Mackey T|K)]_K,
\end{align*}
using the equivalence $(H_\delta\Mackey T)|K \iso \susp_K^{-\delta(G/K)}H_\delta(\Mackey T|K)$ shown above. 
For the homology isomorphism, we have
\begin{align*}
 [S^\alpha, H^\delta\MackeyOp S \smsh (G_+\smsh_K X)]_G
 &\iso [S^\alpha, G_+\smsh_K ((H^\delta\MackeyOp S)|K \smsh X)]_G \\
 &\iso [S^\alpha, \susp_K^{\Lie(G/K)}(H^\delta\MackeyOp S)|K \smsh X]_K \\
 &\iso [S^\alpha, \susp_K^{\Lie(G/K)} \susp_K^{-(\Lie(G/K)-\delta(G/K))}H^\delta(\MackeyOp S|K)\smsh X]_K \\
 &\iso [S^\alpha, \susp_K^{\delta(G/K)}H^\delta(\MackeyOp S|K)\smsh X]_K \\
 &\iso [S^{\alpha-\delta(G/K)}, H^\delta(\MackeyOp S|K)\smsh X]_K.
\end{align*}
Here, the second isomorphism is shown in \cite[II.6.5]{LMS:eqhomotopy} and the third is
again from Proposition~\ref{prop:emsubgroups}.
That these isomorphisms agree with the ones constructed on the chain level above
follows by comparing their behavior on 
the filtration quotients $X^{\alpha-\delta(G/K)+n}/X^{\alpha-\delta(G/K)+n-1}$.

We now distinguish a particular map from $G$ homology to $K$-homology.
The definition of a ``restriction to subgroups'' map is clear on the represented level.
Let $\delta$ be a familial dimension function for $G$
and let $\Mackey T$ be a contravariant $\delta$-Mackey functor.
Then restriction of cohomology from $G$ to a subgroup $K\in\F(\delta)$ should be the following map:
\begin{align*}
 \tilde H_{G,\delta}^\alpha(X;\Mackey T)
  &\iso [\susp_G^\infty X, \susp_G^\alpha H_\delta\Mackey T]_G \\
  &\to [\susp_K^\infty X, \susp_K^\alpha (H_\delta\Mackey T)|K ]_K \\
  &\iso [\susp_K^\infty X, \susp_K^\alpha\susp_K^{-\delta(G/K)}H_\delta(\Mackey T|K)]_K \\
  &\iso \tilde H_{K,\delta}^{\alpha-\delta(G/K)}(X;\Mackey T|K).
\end{align*}
The arrow above is restriction of $G$-stable maps to $K$-stable maps and
we use Proposition~\ref{prop:emsubgroups} to identify the restriction of the
Eilenberg-Mac\,Lane spectrum.
For homology, we have the following map:
\begin{align*}
 \tilde H^{G,\delta}_\alpha(X;\MackeyOp S)
  &\iso [S^\alpha, H^\delta\MackeyOp S \smsh X]_G \\
  &\to [S^\alpha, (H^\delta\MackeyOp S)|K \smsh X]_K \\
  &\iso [S^\alpha, \susp_K^{-(\Lie(G/K)-\delta(G/K))}H^\delta(\MackeyOp S|K) \smsh X]_K \\
  &\iso [S^{\alpha+\Lie(G/K)-\delta(G/K)}, H^\delta(\MackeyOp S|K) \smsh X]_K \\
  &\iso \tilde H^{K,\delta}_{\alpha+\Lie(G/K)-\delta(G/K)}(X;\MackeyOp S|K).
\end{align*}
We write these maps as $a \mapsto a|K$ for $a$ an element of homology or cohomology.

We can also describe restriction to subgroups entirely in terms of space-level maps.
The key to doing so is noticing that the restriction map $[E,F]_G \to [E,F]_K$ factors
as
\[
 [E,F]_G \to [G/K_+\smsh E, F]_G \iso [E,F]_K,
\]
where the first map is induced by the projection $G/K\to *$;
using duality we see that it also factors as
\[
 [E,F]_G \to [E, D(G/K_+)\smsh F]_G \iso [E,F]_K
\]
where the first map is induced by the stable map $S\to D(G/K_+)$ dual to the projection
and the isomorphism is the one shown in \cite[II.6.5]{LMS:eqhomotopy}.
(Note that, if $G\not\in\F(\delta)$ or $\delta(G/K)\neq 0$, the projection
$G/K\to *$ will not be a map in $\sorb{G,\delta}$, so the restriction map
we're now considering may not
be one of the maps in the Mackey functor structure on $\Mackey H_{G,\delta}^*$.)
With the first factorization in mind, we can see that the restriction in cohomology is
given by
\begin{align*}
 \tilde H_{G,\delta}^\alpha(X;\Mackey T)
  &\to \tilde H_{G,\delta}^\alpha(G/K_+\smsh X;\Mackey T) \\
  &\iso \tilde H_{K,\delta}^{\alpha-\delta(G/K)}(X;\Mackey T|K)
\end{align*}
where the first map is induced by the projection and the second map is the
Wirth\-m\"uller isomorphism.
Note that, if $X$ is a based $\delta$-$G$-CW$(\alpha)$ complex, there is, in general, no
canonical cell structure on $G/K_+\smsh X$, so we do not attempt to describe
restriction on the chain level.
(Put another way, if $X$ is a based $\delta$-$G$-CW$(\alpha)$ complex,
there is no canonical $\delta$-$K$-CW$(\alpha-\delta(G/K))$ structure on $X$.)

To describe the map in homology, we need a space-level map representing
the stable map $S\to D(G/K_+)$ dual to the projection $G/K \to *$.
Here's the classical construction: Let $V$ be a representation large enough that
there is an embedding $G/K\subset V$. A tubular neighborhood has the
form $G\times_K D(V-\Lie(G/K))$, so there is a collapse map
\[
 c\colon S^V \to G_+\smsh_K S^{V-\Lie(G/K)}.
\]
This represents $S\to D(G/K_+)$.
The restriction in homology can now be described as follows:
\begin{align*}
 \tilde H^{G,\delta}_\alpha(X;\MackeyOp S)
  &\iso \tilde H^{G,\delta}_{\alpha+V}(\susp_G^V X; \MackeyOp S) \\
  &\xrightarrow{c_*} \tilde H^{G,\delta}_{\alpha+V}(G_+\smsh_K \susp^{V-\Lie(G/K)}X;\MackeyOp S) \\
  &\iso \tilde H^{K,\delta}_{\alpha+V-\delta(G/K)}(\susp^{V-\Lie(G/K)}X;\MackeyOp S|K) \\
  &\iso \tilde H^{K,\delta}_{\alpha+\Lie(G/K)-\delta(G/K)}(X;\MackeyOp S|K).
\end{align*}
This again uses the Wirthm\"uller isomorphism.

\begin{remark}
The specializations of these maps to the ordinary and dual theories give us the
following maps, each of which we denote by $a\mapsto a|K$:
\begin{align*}
 \tilde H^G_\alpha(X;\MackeyOp S) &{}\to \tilde H^K_{\alpha+\Lie(G/K)}(X;\MackeyOp S|K) 
  & \tilde\H_G^\alpha(X;\MackeyOp S) &{}\to \tilde\H_K^{\alpha-\Lie(G/K)}(X;\MackeyOp S|K) \\
 \tilde H_G^\alpha(X;\Mackey T) &{}\to \tilde H_K^{\alpha}(X;\Mackey T|K)
  & \tilde\H^G_\alpha(X;\Mackey T) &{}\to \tilde\H^K_{\alpha}(X;\Mackey T|K) \\
\end{align*}
\end{remark}

\subsection{Quotient groups}\label{subsec:quotientgroups}

Let $N$ be a normal subgroup of $G$ and let $\epsilon\colon G\to G/N$ denote the
quotient map. The results of this section are essentially an elaboration of
two observations: First, that, if $Y$ is a $G/N$-CW complex, then it can be considered
a $G$-CW complex via $\epsilon$. Second, that, if $X$ is a $G$-CW complex,
then $X^N$ has a natural structure as a $G/N$-CW complex.

We use repeatedly the fact that
\[
 (G/H)^N =
  \begin{cases}
    G/H & \text{if $N\leq H$} \\
    \emptyset & \text{if $N\not\leq H$,}
  \end{cases}
\]
which can be seen by observing that $G$, hence $G/N$, acts transitively on $(G/H)^N$,
hence $(G/H)^N$ must be either empty or an orbit of $G/N$.
In the case that $G/H$ is fixed by $N$, $N$ will act trivially on $\Lie(G/H)$, hence
on $\delta(G/H)$ for any dimension function $\delta$ for $G$.

This calculation allows us to make the following definition.

\begin{definition}
Let $\delta$ be a dimension function for $G$ and let $N$ be a normal subgroup of $G$.
We write $\delta|G/N$ for the dimension function on $G/N$ defined by
\[
 \F(\delta|G/N) = \{ H/N \mid N\leq H\text{ and }H\in\F(\delta) \}
\]
and
\[
 (\delta|G/N)((G/N)/(H/N)) = \delta(G/H).
\]
We will usually write $\delta$ again for $\delta|G/N$, for simplicity of notation.
\end{definition}

The following two results record the observations we mentioned at the beginning of the section.
If $Y$ is a $G/N$-space, we write $\epsilon^* Y$ for $Y$ considered as a $G$-space
via $\epsilon$. (We use this notation only when we want to be very explicit; when the meaning
is understood we will write $Y$ for $\epsilon^* Y$.)

\begin{proposition}\label{prop:inducedCWstructure}
Let $\delta$ be a dimension function for $G$ and let $N$ be a normal subgroup of $G$.
Then, for $\alpha$ a virtual representation of $G/N$, $\epsilon^*$ defines functors
\[
 \epsilon^*\colon (G/N)\W^{\delta,\alpha} \to G\W^{\delta,\alpha}
\]
and
\[
 \epsilon^*\colon (G/N)\W_*^{\delta,\alpha} \to G\W_*^{\delta,\alpha}.
\]
\end{proposition}

\begin{proof}
This follows from the fact that
a $G/N$-$\delta$-$\alpha$-cell of the form
\[
 G/N\times_{H/N}\bar D(\alpha-\delta((G/N)/(H/N))+n)
\]
can be considered a
$G$-$\delta$-$\alpha$-cell of the form
\[
 G\times_{H}\bar D(\alpha-\delta(G/H)+n).
\]
If $f$ is a cellular $G/N$-map, it remains cellular when considered as a $G$-map.
\end{proof}

\begin{proposition}\label{prop:CWfixedsets}
Let $\delta$ be a dimension function for $G$ and let $N$ be a normal subgroup of $G$.
Then, for $\alpha$ a virtual representation of $G$, $(-)^N$ defines functors
\[
 (-)^N\colon G\W^{\delta,\alpha} \to (G/N)\W^{\delta,\alpha^N}
\]
and
\[
 (-)^N\colon G\W_*^{\delta,\alpha} \to (G/N)\W_*^{\delta,\alpha^N}.
\]
\end{proposition}

\begin{proof}
Let $X$ be a based $\delta$-$G$-CW($\alpha$) complex.
The fixed set of a cell in $X$ of the form
$G\times_H \bar D(\alpha-\delta(G/H)+n)$ with $N\leq H$ is a cell in $X^N$ of the form
$G/N\times_{H/N} \bar D(\alpha^N-\delta(G/H)+n)$.
Any such cell in $X$ with $N\not\leq H$ will not be fixed by $N$ hence not appear in $X^N$.
Thus, $X^N$ has a canonical structure as a $\delta$-$G/N$-CW$(\alpha^N)$ complex.
If a map $f$ is cellular, then so is $f^N$.
\end{proof}

To describe the algebra involved on chains, we describe two related functors,
\[
 \theta\colon \sorb{G/N,\delta}\to\sorb{G,\delta}
\]
and
\[
 \Phi^N\colon \sorb{G,\delta} \to \sorb{G/N,\delta}.
\]
Recall from \cite{LMS:eqhomotopy} or \cite{MM:orthogonal}
that there is a functor $\epsilon^\sharp$ taking $G/N$-spectra to $G$-spectra,
obtained by letting $G$ act via $\epsilon$ and then extending to a complete $G$-universe.
This functor has a right adjoint $(-)^N$, called the {\em categorical} fixed point functor,
and there results a stable adjunction
\[
  [\epsilon^\sharp D, E]_{G} \iso [D, E^N]_{G/N}.
\]
If $Y$ is a based $G/N$-space, then 
\[
 \epsilon^\sharp\susp_{G/N}^\infty Y = \susp_{G}^\infty \epsilon^*Y.
\]
However, the similar statement for $N$-fixed
points does not hold: $(\susp_G^\infty X)^N \not\iso \susp_{G/N}^\infty X^N$ in general
in the stable category.
We also have $(E\smsh D)^N \not\iso E^N\smsh D^N$ in general, but we do have the following
useful fact, not appearing elsewhere to our knowledge.

\begin{proposition}
If $E$ is a $G$-spectrum and $A$ is a $G/N$-spectrum, then there is an equivalence
\[
 E^N\smsh A \hmtpc (E\smsh \epsilon^\sharp A)^N 
\]
in the stable category.
\end{proposition}

\begin{proof}
We begin by noticing that there is a stable equivalence
\[
 \epsilon^\sharp(A\smsh B) \hmtpc \epsilon^\sharp A\smsh \epsilon^\sharp B
\]
for $G/N$-spectra $A$ and $B$,
which follows from the definitions and the fact that both sides preserve
cofibrations and acyclic cofibrations. By adjunction, we also get an equivalence
\[
 F(A,E^N) \hmtpc F(\epsilon^\sharp A, E)^N.
\]

There is a stable map $\epsilon^\sharp F(X,B) \to F(\epsilon^\sharp X, \epsilon^\sharp B)$,
adjoint to the composite
\[
 \epsilon^\sharp F(X,B)\smsh \epsilon^\sharp X 
 \hmtpc \epsilon^\sharp( F(X,B)\smsh X) \to \epsilon^\sharp B,
\]
that need not be an equivalence in general. However, if $X$ is finite (i.e., dualizable) and $B = S$,
then $\epsilon^\sharp F(X,S) \hmtpc F(\epsilon^\sharp X, S)$
(note that $\epsilon^\sharp S = S$) because the dual of $\epsilon^\sharp X$
is equivalent to $\epsilon^\sharp DX$.

We can now define $E^N\smsh A \to (E\smsh \epsilon^\sharp A)^N$ as the composite
\begin{align*}
 E^N\smsh A
  &\to [\epsilon^\sharp(E^N\smsh A)]^N \\
  &\hmtpc [\epsilon^\sharp (E^N) \smsh \epsilon^\sharp A]^N \\
  &\to (E\smsh \epsilon^\sharp A)^N,
\end{align*}
where the first map is the unit of the $\epsilon^\sharp$-$(-)^N$ adjunction and
the last map is induced by the counit.
Writing down an appropriate diagram shows that, when $A = F(X,S)$ with $X$ (hence $A$) finite,
this composite agrees with the composite
\begin{align*}
 E^N\smsh F(X,S) 
  &\hmtpc F(X,E^N) \\
  &\hmtpc F(\epsilon^\sharp X,E)^N \\
  &\hmtpc [E \smsh F(\epsilon^\sharp X, S)]^N \\
  &\hmtpc [E \smsh \epsilon^\sharp F(X,S)]^N,
\end{align*}
where the first and third maps are equivalences by duality.
Therefore, $E^N\smsh A \to (E\smsh \epsilon^\sharp A)^N$ is an equivalence when $A$
is finite, so, in particular, when $A$ is an orbit. 
Because both sides preserve wedges and cofibration sequences, it follows
that the map is an equivalence for all cell complexes, thence for all spectra.
(Put another way, the stable homotopy groups of $E^N\smsh A$ and $(E\smsh\epsilon^\sharp A)^N$ define
$G/N$-homology theories in $A$ that agree on orbits, hence are isomorphic.)
\end{proof}

Better than the categorial fixed point functor for many purposes is the so-called
{\em geometric} fixed-point construction discussed in
\cite[II.9]{LMS:eqhomotopy}, \cite[XVI.3]{May:alaska}, and \cite[V.4]{MM:orthogonal}.
The following is the quickest definition.

\begin{definition}\label{def:geoFixedPoints}
If $N$ is a normal subgroup of $G$ and $E$ is a $G$-spectrum, we define
the {\em geometric fixed point spectrum} of $E$ to be the $G/N$-spectrum given by
\[
 \Phi^N(E) = (E \smsh \tilde E\F[N])^N,
\]
where $\tilde E\F[N]$ is the based $G$-space characterized by
\[
 \tilde E\F[N]^L \hmtpc
 \begin{cases} S^0 & \text{if $N \leq L$} \\
               * & \text{if $N \not\leq L$.}
 \end{cases}
\]
\end{definition}

There is another useful definition of $\Phi^N$---we recall the definition for the $G$-prespectra
of \cite{LMS:eqhomotopy}, the details being somewhat different for the orthogonal spectra of
\cite{MM:orthogonal}.
Consider then a $G$-prespectrum $D$ indexed on a sequence $\V = \{V_i\}$ in a
complete $G$-universe $\U$.
Then $\V^N = \{V_i^N\}$ is an indexing sequence in the complete $G/N$-universe $\U^N$.
We define
\[
 (\Phi^N D)(V_i^N) = D(V_i)^N
\]
and define the structure maps to be
\[
 \susp^{V_{i+1}^N - V_i^N}D(V_i)^N = [\susp^{V_{i+1}-V_i}D(V_i)]^N
  \to D(V_i)^N,
\]
the restriction to the $N$-fixed points of the original structure map.
The equivalence of this to Definition~\ref{def:geoFixedPoints}
is shown in \cite[II.9.8]{LMS:eqhomotopy}.

The main facts we need to know
about $\Phi^N$ are that we have the following natural stable equivalences for based
$G$-spaces $X$ and $G$-spectra $E$ and $E'$:
\begin{align*}
 \Phi^N(\susp_G^\infty X) &\hmtpc \susp_{G/N}^\infty X^N \\
 \Phi^N(E \smsh X) &\hmtpc \Phi^N(E) \smsh X^N \\
 \Phi^N(E\smsh E') &\hmtpc \Phi^N(E)\smsh \Phi^N(E')
\end{align*}
In particular, $\Phi^N$ respects suspension in the sense that
$\Phi^N(\susp^W E) \hmtpc \susp^{W^N}\Phi^N(E)$.
(See \cite[II.9]{LMS:eqhomotopy} and \cite[V.4]{MM:orthogonal}.)
Not shown in those references, but following directly from the alternate definition of $\Phi^N$,
we have the following equivalence if $N\leq L\leq G$:
\[
 \Phi^N(G_+\smsh_L E) \hmtpc (G/N)_+\smsh_{L/N} \Phi^N E.
\]
We also have the following relationship between $\Phi^N$ and $\epsilon^\sharp$, which
is \cite[II.9.10]{LMS:eqhomotopy}:
\[
 \Phi^N\epsilon^\sharp D \hmtpc D,
\]
which can also be written as
\[
 (\epsilon^\sharp D \smsh \tE\F[N])^N \hmtpc D.
\]

We have the following behavior of these functors on spheres:
\[
 \epsilon^\sharp (G/N_+\smsh_{H/N} S^{-\delta((G/N)/(G/H))})
  = G_+\smsh_H S^{-\delta(G/H)}
\]
for $N\leq H\leq G$, and
\[
 \Phi^N(G_+\smsh_H S^{-\delta(G/H)}) \hmtpc
  \begin{cases}
    (G/N)_+\smsh_{H/N} S^{-\delta(G/H)} & \text{if $N\leq H$} \\
    * & \text{if $N\not\leq H$.}
  \end{cases}
\]
These allow us to make the following definitions.

\begin{definition}
Let $N$ be a normal subgroup of $G$. We write
\[
 \theta\colon \sorb{G/N,\delta}\to\sorb{G,\delta}
\]
for the restriction of $\epsilon^\sharp$.
As in Remark~\ref{rem:zeroOrbit}, we can consider $\sorb{G,\delta}$
and $\sorb{G/N,\delta}$ as augmented with
a zero object $*$ given by the trivial spectrum.
We write 
\[
 \Phi^N\colon \sorb{G,\delta}\to \sorb{G/N,\delta}
\]
for the restriction
of $\Phi^N$ to the augmented $\delta$-orbit category.
In this context we write $((G/H)^N,\delta) \in \sorb{G/N}$ 
for $\Phi^N(G/H,\delta)$.
Note that $((G/H)^N,\delta)$ is $*$ if $N\not\leq H$
and is $((G/N)/(H/N),\delta)$ if $N\leq H$.
\end{definition}

We use these functors to define operations on Mackey functors.
We have the operations $\theta^*$ on $G$-$\delta$-Mackey functors and
$\theta_!$ on $G/N$-$\delta$-Mackey functors. We give special names
to the functors associated with $\Phi^N$:

\begin{definition}
If $\MackeyOp S$ is a covariant and $\Mackey T$ a contravariant
$G$-$\delta$-Mackey functor, let
\begin{align*}
 \MackeyOp S^N &= \Phi^N_!\MackeyOp S \qquad\text{and}\\
 \Mackey T^N &= \Phi^N_!\Mackey T.
\end{align*}
We call this the {\em $N$-fixed point functor.}
If $\MackeyOp U$ is a covariant and $\Mackey V$ is a contravariant
$G/N$-$\delta$-Mackey functor, let
\begin{align*}
 \Inf_{G/N}^G \MackeyOp U &= (\Phi^N)^*\MackeyOp U \qquad\text{and}\\
 \Inf_{G/N}^G \Mackey V &= (\Phi^N)^*\Mackey V.
\end{align*}
We call this the {\em inflation functor.}
 \end{definition}

The fixed point and inflation constructions are, of course, adjoint:
\[
 \Hom_{\sorb{G/N,\delta}}(\MackeyOp S^N, \MackeyOp U)
  \iso \Hom_{\sorb{G,\delta}}(\MackeyOp S, \Inf_{G/N}^G\MackeyOp U)
\]
and similarly for contravariant functors.

Th\'evenaz and Webb gave essentially this
definition of the inflation functors for finite groups in \cite[\S5]{TW:mackey}.
They defined the left adjoint, which they wrote as $\MackeyOp S^+$ rather than $\MackeyOp S^N$,
in a different way than we did above but, given the uniqueness of adjuncts, their $\MackeyOp S^+$ must
agree with $\Phi^N_!\MackeyOp S$.

Note that we have $\Phi^N\theta = 1$, which gives
\[
 \theta^*\Inf_{G/N}^G\MackeyOp U \iso \MackeyOp U
\]
and
\[
 (\theta_!\MackeyOp S)^N \iso \MackeyOp S,
\]
and similarly for contravariant Mackey functors.

We next record the effects of $\theta_!$ and $(-)^N$ on the canonical projectives.
To make explicit the varying underlying group, we write
\[
 \MackeyOp A^{G/H,\delta}_G = \sorb{G,\delta}((G/H,\delta),-)
\]
and
\[
 \Mackey A_{G/H,\delta}^G = \sorb{G}(-,(G/H,\delta)).
\]

\begin{proposition}\label{prop:burnsidefixedsets}
If $N$ is normal in $G$, then
 \[
 (\MackeyOp A_G^{G/H,\delta})^N \iso 
  	\begin{cases} \MackeyOp A_{G/N}^{G/H,\delta} & \text{if $N\leq H$} \\
                  0 & \text{if $N\not\leq H$}
    \end{cases}
 \]
and
 \[
 (\Mackey A^G_{G/H,\delta})^N \iso 
  	\begin{cases} \Mackey A^{G/N}_{G/H,\delta} & \text{if $N\leq H$} \\
                  0 & \text{if $N\not\leq H$}
    \end{cases}
 \]
for any subgroup $H$ of $G$.
For $N\leq L\leq G$ we have
\[
 \theta_! \MackeyOp A_{G/N}^{G/L,\delta} \iso \MackeyOp A_G^{G/L,\delta}
\]
and
\[
 \theta_! \Mackey A^{G/N}_{G/L,\delta} \iso \Mackey A^G_{G/L,\delta}.
\]
\end{proposition}

\begin{proof}
These are special cases of Proposition~\ref{prop:indresadjunction}(\ref{item:indFree}).
 \end{proof}

Our main calculations are then the following two results.

\begin{proposition}\label{prop:inducedchains}
Let $N$ be a normal subgroup of $G$, let $\delta$ be a dimension function for $G$,
and let $\alpha$ be a virtual representation of $G/N$.
Let $Y$ be a based $\delta$-$G/N$-CW($\alpha$) complex and give
$\epsilon^* Y$ the corresponding $\delta$-$G$-CW($\alpha$) structure as in
Proposition~\ref{prop:inducedCWstructure}. Then we have a natural chain isomorphism
\[
 \theta_!\Mackey C^{G/N,\delta}_{\alpha+*}(Y,*) \iso
  \Mackey C^{G,\delta}_{\alpha+*}(Y,*),
\]
where we write $Y$ for $\epsilon^* Y$ on the right.
This isomorphism respects suspension in the sense that, if $W$ is a representation
of $G/N$, then the following diagram commutes:
\[
 \xymatrix{
   \theta_!\Mackey C^{G/N,\delta}_{\alpha+*}(Y,*) \ar[r]^\iso \ar[d]_{\sigma^W}
    & \Mackey C^{G,\delta}_{\alpha+*}(Y,*) \ar[d]^{\sigma^W} \\
   \theta_!\Mackey C^{G/N,\delta}_{\alpha+W+*}(\susp^W Y,*) \ar[r]_\iso
    & \Mackey C^{G,\delta}_{\alpha+W+*}(\susp^W Y,*).
 }
\]
\end{proposition}

\begin{proof}
The isomorphism is the chain map adjoint to the map
\[
 \Mackey C^{G/N,\delta}_{\alpha+*}(Y,*) \to \theta^*\Mackey C^{G,\delta}_{\alpha+*}(Y,*)
\]
given by
\begin{multline*}
 \epsilon^\sharp\colon
  [G/N_+\smsh_{H/N} S^{\alpha-\delta((G/N)/(H/N))+n}, 
     \susp_{G/N}^\infty Y^{\alpha+n}/Y^{\alpha+n-1}]_{G/N} \\
  \to [G_+\smsh_H S^{\alpha-\delta(G/H)+n}, \susp_G^\infty Y^{\alpha+n}/Y^{\alpha+n-1}]_G
\end{multline*}
when $N\leq H\leq G$. That the adjoint is an isomorphism follows from 
Proposition~\ref{prop:burnsidefixedsets}.
That the isomorphism respects suspension follows from the fact that $\epsilon^\sharp$ does.
\end{proof}

\begin{proposition}\label{prop:fixedSetChainIso}
Let $N$ be a normal subgroup of $G$, let $\delta$ be a dimension function for $G$,
and let $\alpha$ be a virtual representation of $G$.
Let $X$ be a based $\delta$-$G$-CW$(\alpha)$ complex
and give $X^N$ the $\delta$-$G/N$-CW$(\alpha^N)$
structure from Proposition~\ref{prop:CWfixedsets}.
With this structure, we have a natural chain isomorphism
\[
 \Mackey C^{G,\delta}_{\alpha+*}(X,*)^N \iso \Mackey C^{G/N,\delta}_{\alpha^N+*}(X^N,*).
\]
This isomorphism respects suspension in the sense that, if $W$ is a representation $G$,
then the following diagram commutes:
\[
 \xymatrix{
  \Mackey C^{G,\delta}_{\alpha+*}(X,*)^N \ar[r]^-\iso \ar[d]_{(\sigma^W)^N}
   & \Mackey C^{G/N,\delta}_{\alpha^N+*}(X^N,*) \ar[d]^{\sigma^{W^N}} \\
  \Mackey C^{G,\delta}_{\alpha+W+*}(\susp^W X,*)^N \ar[r]_-\iso
   & \Mackey C^{G/N,\delta}_{\alpha^N+W^N+*}(\susp^{W^N}X^N,*)
 }
\]  
\end{proposition}

\begin{proof}
To see the chain isomorphism, note that,
for each subgroup $L$ of $G$ containing $N$, we have the map
\begin{multline*}
 \Phi^N\colon [G_+\smsh_L S^{\alpha-\delta(G/L)+n}, \susp_G^\infty X^{\alpha+n}/X^{\alpha+n-1}]_G \\
  \to [ G/N_+\smsh_{L/N} S^{\alpha^N-\delta(G/L)+n}, \susp_{G/N}^\infty (X^{\alpha+n}/X^{\alpha+n-1})^N]_{G/N};
\end{multline*}
if $L$ does not contain $N$ then the result of applying $\Phi^N$ is 0.
This defines (the adjoint of) the map 
$\Mackey C^{G,\delta}_{\alpha+n}(X,*)^N \to \Mackey C^{G/N,\delta}_{\alpha^N+n}(X^N,*)$,
which we see is an isomorphism using Proposition~\ref{prop:burnsidefixedsets}.
That this isomorphism respects suspension follows from the fact that $\Phi^N$ does.
\end{proof}

Before stating the consequences for homology and cohomology we insert the following
definition and lemma.

\begin{definition}\label{def:Nclosed}
Let $N$ be a normal subgroup of $G$ and
let $\delta$ be a dimension function for $G$.
We say that $\delta$ is {\em $N$-closed} if, whenever $K\in\F(\delta)$,
we also have $KN\in\F(\delta)$.
\end{definition}

\begin{lemma}
Let $N$ be a normal subgroup of $G$, let $\delta$ be an
$N$-closed familial dimension function for $G$,
and let $\alpha$ be a virtual representation of $G/N$.
If $f\colon X\to Y$ is a $(\delta|G/N)$-weak$_\alpha$ equivalence of $G/N$-spaces,
then it is a $\delta$-weak$_\alpha$ equivalence of $G$-spaces.

Further, if all subgroups in $\F(\delta|G/N)$ are admissible and $X$ and $Y$ are
well-based, then
$\susp_G^W f$ is a $\delta$-weak$_\alpha$ equivalence of $G$-spaces for any
representation $W$ of $G$.
\end{lemma}

\begin{proof}
Suppose $K$ is a $\delta$-$\alpha$-admissible subgroup of $G$. Because $X$ and $Y$ are fixed
by $N$ we have $f^K = f^{KN}$. To say that $K$ is admissible is to say that
$\alpha - \delta(G/K) + n$ is equivalent to an actual representation for some $n$,
which implies that $\alpha - \delta(G/KN) + n \iso (\alpha-\delta(G/K)+n) + \delta(KN/K)$
is also equivalent to an actual representation. Hence, $KN$ is also admissible.
By Theorem~\ref{thm:weakCharacterization} we have that $f^{KN}$ is a weak equivalence,
hence $f^K$ is, too. Applying Theorem~\ref{thm:weakCharacterization} again, we
conclude that $f$ is a $\delta$-weak$_\alpha$ equivalence of $G$-spaces.

If all subgroups in $\F(\delta|G/N)$ are $\delta$-$\alpha$-admissible, then the argument above
shows that all subgroups in $\F(\delta)$ are $\delta$-$\alpha$-admissible as well.
Thus, in this case, $f^K$ is a weak equivalence for all $K\in\F(\delta)$, so
$(\susp_G^W f)^K$ is a weak equivalence for all $K\in\F(\delta)$, hence
$\susp_G^W f$ is again a $\delta$-weak$_\alpha$ equivalence.
\end{proof}

\begin{theorem}\label{thm:inducedHomology}
Let $N$ be a normal subgroup of $G$,
let $\delta$ be a dimension function for $G$,
let $\alpha$ be a virtual representation of $G/N$,
let $\MackeyOp S$ be a covariant $G$-$\delta$-Mackey functor,
and let $\Mackey T$ be a contravariant $G$-$\delta$-Mackey functor.
Then, for $Y$ in $(G/N)\W_*^{\delta,\alpha}$ we have natural isomorphisms
\begin{align*}
 \tCH^{G,\delta}_{\alpha}(Y;\MackeyOp S)
   &\iso \tCH^{G/N,\delta}_{\alpha}(Y; \theta^* \MackeyOp S) \qquad\text{and} \\
 \tCH_{G,\delta}^{\alpha}(Y;\Mackey T)
   &\iso \tCH_{G/N,\delta}^{\alpha}(Y; \theta^*\Mackey T).
\end{align*}
These isomorphisms respect suspension in the sense that, if
$W$ is a representation of $G/N$,
then the following diagram commutes:
\[
 \xymatrix{
  \tCH^{G,\delta}_{\alpha}(Y;\MackeyOp S) \ar[r]^-\iso \ar[d]_{\sigma^W}
   & \tCH^{G/N,\delta}_{\alpha}(Y; \theta^*\MackeyOp S) \ar[d]^{\sigma^{W}} \\
  \tCH^{G,\delta}_{\alpha+W}(\susp^W Y;\MackeyOp S) \ar[r]_-\iso
   & \tCH^{G/N,\delta}_{\alpha+W}(\susp^{W}Y; \theta^*\MackeyOp S)
 }
\]
and similarly for cohomology.
If $\delta$ is $N$-closed and familial,
then, for $Y$ in $(G/N)\K_*$ we have natural isomorphisms
\begin{align*}
 \tilde H^{G,\delta}_{\alpha}(Y;\MackeyOp S)
   &\iso \tilde H^{G/N,\delta}_{\alpha}(Y; \theta^*\MackeyOp S) \qquad\text{and} \\
 \tilde H_{G,\delta}^{\alpha}(Y;\Mackey T)
   &\iso \tilde H_{G/N,\delta}^{\alpha}(Y; \theta^*\Mackey T).
\end{align*}
These isomorphisms respect suspension in the sense that, if
$Y$ is well-based and $W$ is a representation of $G/N$,
then the following diagram commutes:
\[
 \xymatrix{
  \tilde H^{G,\delta}_{\alpha}(Y;\MackeyOp S) \ar[r]^-\iso \ar[d]_{\sigma^W}
   & \tilde H^{G/N,\delta}_{\alpha}(Y; \theta^*\MackeyOp S) \ar[d]^{\sigma^{W}} \\
  \tilde H^{G,\delta}_{\alpha+W}(\susp^W Y;\MackeyOp S) \ar[r]_-\iso
   & \tilde H^{G/N,\delta}_{\alpha+W}(\susp^{W}Y; \theta^*\MackeyOp S)
 }
\]
and similarly for cohomology.
\end{theorem}

\begin{proof}
For $Y$ in $(G/N)\W_*^{\delta,\alpha}$, the theorem follows from
Propositions~\ref{prop:inducedchains} and~\ref{prop:indresadjunction}, which give
\begin{align*}
 \Mackey C^{G,\delta}_{\alpha+*}(Y) \tensor_{\sorb{G,\delta}}\MackeyOp S
  &\iso \theta_!\Mackey C^{G/N,\delta}_{\alpha+*}(Y)\tensor_{\sorb{G,\delta}}\MackeyOp S \\
  &\iso \Mackey C^{G/N,\delta}_{\alpha+*}(Y)\tensor_{\sorb{G/N,\delta}}\theta^*\MackeyOp S
\end{align*}
and
\begin{align*}
 \Hom_{\sorb{G,\delta}}(\Mackey C^{G,\delta}_{\alpha+*}(Y), \Mackey T)
  &\iso \Hom_{\sorb{G,\delta}}(\theta_!\Mackey C^{G/N,\delta}_{\alpha+*}(Y), \Mackey T) \\
  &\iso \Hom_{\sorb{G/N,\delta}}(\Mackey C^{G/N,\delta}_{\alpha+*}(Y), \theta^*\Mackey T).
\end{align*}
That the isomorphisms respect suspension follows from the similar statement in
Proposition~\ref{prop:inducedchains}.

Now consider a general based $G/N$-space $Y$, which we may assume is well-based.
Let $\U$ be a complete $G$-universe so that $\U^N$ is a complete $G/N$-universe.
Let $\V$ be an indexing sequence in $\U$ so that $\V^N$ is an indexing sequence in $\U^N$.
Let $\Gamma\susp_{G/N}^\infty Y \to \susp_{G/N}^\infty Y$ be a
$\delta$-$G/N$-CW($\alpha$) approximation indexed on $\V^N$.
We can then construct a $\delta$-$G$-CW($\alpha$) approximation
$\Gamma\susp_G^\infty Y\to \susp_G^\infty Y$ indexed on $\V$ such that
there are inclusions
\[
 \susp_G^{V_i - V_i^N}(\Gamma\susp_{G/N}^\infty Y)(V_i^N) \to
   (\Gamma\susp_G^\infty Y)(V_i)
\]
over $\susp_G^{V_i}Y$ and compatible with the structure maps.
For sufficiently large $i$, all subgroups in
$\F(\delta|G/N)$ are $\delta$-$(\alpha+V_i^N)$-admissible. It follows from
the preceding lemma that
\[
 \susp_G^{V_i - V_i^N}(\Gamma\susp_{G/N}^\infty Y)(V_i^N) \to \susp_G^{V_i} Y
\]
is a $\delta$-weak$_{\alpha+V_i}$ equivalence.
Using Proposition~\ref{prop:inducedchains}, we have a chain homotopy equivalence
\[
 \theta_!\Mackey C^{\V^N,\delta}_{\alpha+*}(Y) \hmtpc \Mackey C^{\V,\delta}_{\alpha+*}(Y).
\]
The theorem now follows as it did in the CW case.
\end{proof}

This theorem leads to a statement about Mackey functor-valued homology and cohomology.

\begin{corollary}
Let $N$ be a normal subgroup of $G$,
let $\delta$ be an $N$-closed familial dimension function for $G$,
let $\alpha$ be a virtual representation of $G/N$,
let $\MackeyOp S$ be a covariant $G$-$\delta$-Mackey functor, and
let $\Mackey T$ be a contravariant $G$-$\delta$-Mackey functor.
Then, for any based $G/N$-space $Y$, we have
\[
 \theta^* \MackeyOp H^{G,\delta}_\alpha(Y;\MackeyOp S)
   \iso \MackeyOp H^{G/N,\delta}_\alpha(Y;\theta^*\MackeyOp S)
\]
and
\[
 \theta^* \Mackey H_{G,\delta}^\alpha(Y;\Mackey T)
   \iso \Mackey H_{G/N,\delta}^\alpha(Y;\theta^*\Mackey T).
\]
\end{corollary}

\begin{proof}
If $N\leq L\leq G$, we have the following commutative diagram:
\[
 \xymatrix{
  \sorb{L/N,\delta} \ar[r]^-\theta \ar[d]_{i_{L/N}^{G/N}}
     & \sorb{L,\delta} \ar[d]^{i_L^G} \\
  \sorb{G/N,\delta} \ar[r]_-\theta 
     & \sorb{G,\delta}
 }
\]
From this we get that $\theta^*(i_L^G)^* \iso (i_{L/N}^{G/N})^*\theta^*$, so
$\theta^*(\MackeyOp S|L) \iso (\theta^*\MackeyOp S)|L$ and the following diagram commutes,
in which we write the spectrum $S^{-\delta(G/L)}$ but implicitly take suspensions
by a $G/N$-representation containing $\delta(G/L)$
to work entirely with spaces:
\[
 \xymatrix{
  \tilde H_\alpha^{G,\delta}(Y\smsh G_+\smsh_L S^{-\delta(G/L)};\MackeyOp S) \ar[r]^-\iso \ar[d]_\iso
    & \tilde H_\alpha^{L,\delta}(Y; \MackeyOp S|L) \ar[d]^\iso \\
  \tilde H_\alpha^{G/N,\delta}(Y\smsh G_+\smsh_L S^{-\delta(G/L)};\theta^*\MackeyOp S) \ar[r]_-\iso
    & \tilde H_\alpha^{L/N,\delta}(Y;(\theta^*\MackeyOp S)|L).
 }
\]
This shows that
$\theta^* \MackeyOp H^{G,\delta}_\alpha(Y;\MackeyOp S)
   \iso \MackeyOp H^{G/N,\delta}_\alpha(Y;\theta^*\MackeyOp S)$
and the argument for cohomology is the same.
\end{proof}

\begin{theorem}\label{thm:fixedSetIso}
Let $\delta$ be a dimension function for $G$,
let $N$ be a normal subgroup of $G$, let $\MackeyOp S$ be a covariant $G/N$-$\delta$-Mackey functor,
and let $\Mackey T$ be a contravariant $G/N$-$\delta$-Mackey functor.
Then, for $X$ in $G\W_*^{\delta,\alpha}$ we have natural isomorphisms
\begin{align*}
 \tCH^{G,\delta}_{\alpha}(X;\Inf_{G/N}^G \MackeyOp S)
   &\iso \tCH^{G/N,\delta}_{\alpha^N}(X^N; \MackeyOp S) \qquad\text{and} \\
 \tCH_{G,\delta}^{\alpha}(X;\Inf_{G/N}^G \Mackey T)
   &\iso \tCH_{G/N,\delta}^{\alpha^N}(X^N; \Mackey T).
\end{align*}
These isomorphisms respect suspension in the sense that, if
$W$ is a representation of $G$,
then the following diagram commutes:
\[
 \xymatrix{
  \tCH^{G,\delta}_{\alpha}(X;\Inf_{G/N}^G \MackeyOp S) \ar[r]^-\iso \ar[d]_{\sigma^W}
   & \tCH^{G/N,\delta}_{\alpha^N}(X^N; \MackeyOp S) \ar[d]^{\sigma^{W^N}} \\
  \tCH^{G,\delta}_{\alpha+W}(\susp^W X;\Inf_{G/N}^G \MackeyOp S) \ar[r]_-\iso
   & \tCH^{G/N,\delta}_{\alpha^N+W^N}(\susp^{W^N}X^N; \MackeyOp S)
 }
\]
and similarly for cohomology.
If $\delta$ is familial, then, for $X$ in $G\K_*$ we have natural isomorphisms
\begin{align*}
 \tilde H^{G,\delta}_{\alpha}(X;\Inf_{G/N}^G \MackeyOp S)
   &\iso \tilde H^{G/N,\delta}_{\alpha^N}(X^N; \MackeyOp S) \qquad\text{and} \\
 \tilde H_{G,\delta}^{\alpha}(X;\Inf_{G/N}^G \Mackey T)
   &\iso \tilde H_{G/N,\delta}^{\alpha^N}(X^N; \Mackey T).
\end{align*}
These isomorphisms respect suspension in the sense that, if
$X$ is well-based and $W$ is a representation of $G$,
then the following diagram commutes:
\[
 \xymatrix{
  \tilde H^{G,\delta}_{\alpha}(X;\Inf_{G/N}^G \MackeyOp S) \ar[r]^-\iso \ar[d]_{\sigma^W}
   & \tilde H^{G/N,\delta}_{\alpha^N}(X^N; \MackeyOp S) \ar[d]^{\sigma^{W^N}} \\
  \tilde H^{G,\delta}_{\alpha+W}(\susp^W X;\Inf_{G/N}^G \MackeyOp S) \ar[r]_-\iso
   & \tilde H^{G/N,\delta}_{\alpha^N+W^N}(\susp^{W^N}X^N; \MackeyOp S)
 }
\]
and similarly for cohomology.
\end{theorem}

\begin{proof}
For $X$ in $G\W_*^{\delta,\alpha}$,
the theorem follows from the following
isomorphisms from Proposition~\ref{prop:indresadjunction}:
\begin{align*}
 \Mackey C^{G,\delta}_{\alpha+*}(X)\tensor_{\sorb{G,\delta}} \Inf_{G/N}^G \MackeyOp S
  &\iso \Mackey C^{G,\delta}_{\alpha+*}(X)^N\tensor_{\sorb{G/N,\delta}} \MackeyOp S \\
  &\iso \Mackey C^{G/N,\delta}_{\alpha^N+*}(X^N)\tensor_{\sorb{G/N,\delta}} \MackeyOp S
\end{align*}
and
\begin{align*}
 \Hom_{\sorb{G,\delta}}(\Mackey C^{G,\delta}_{\alpha+*}(X), \Inf_{G/N}^G \Mackey T)
  &\iso \Hom_{\sorb{G/N,\delta}}(\Mackey C^{G,\delta}_{\alpha+*}(X)^N, \Mackey T) \\
  &\iso \Hom_{\sorb{G/N,\delta}}(\Mackey C^{G/N,\delta}_{\alpha^N+*}(X^N), \Mackey T).
\end{align*}
That the isomorphisms respect suspension follows from the similar statement
in Proposition~\ref{prop:fixedSetChainIso}.

Now consider a general based $G$-space $X$, which we may assume is well-based.
Let $\V$ be an indexing sequence in a complete $G$-universe and let
$\Gamma\susp_G^\infty X\to \susp_G^\infty X$ be a $\delta$-$G$-CW$(\alpha)$
approximation.
Recall that $\V^N$ is an indexing sequence in the complete $G/N$-universe $\U^N$.
We can therefore form the $\delta$-$G/N$-CW$(\alpha^N)$ prespectrum
$\Phi^N(\Gamma\susp_G^\infty X)$ indexed on $\V^N$ given by
\[
 \Phi^N(\Gamma\susp_G^\infty X)(V_i^N) = [(\Gamma\susp_G^\infty X)(V_i)]^N,
\]
which comes with a map $\Phi^N(\Gamma\susp_G^\infty X) \to \susp_{G/N}^\infty X^N$.
This map may not be a $\delta$-$G$-CW$(\alpha)$ approximation because
taking $N$-fixed points does not necessarily preserve $\delta$-weak$_\alpha$ equivalence.
However, for fixed $\alpha$, $(-)^N$ does preserve $\delta$-weak$_{\alpha+V}$
equivalence for sufficiently large $V$. Therefore, if
$\Gamma\susp_{G/N}^\infty X^N \to \susp_{G/N}^\infty X^N$ 
is a $\delta$-$G/N$-CW$(\alpha^N)$ approximation, the induced map
$\Phi^N(\Gamma\susp_G^\infty X) \to \Gamma\susp_{G/N}^\infty X^N$ is
a $\delta$-weak$_{\alpha+V_i^N}$ equivalence at level $i$ for sufficiently large $i$.
Using Proposition~\ref{prop:fixedSetChainIso}, we see that the map induces a chain homotopy equivalence
\[
 \Mackey C^{\V,\delta}_{\alpha+*}(X)^N \hmtpc \Mackey C^{\V^N,\delta}_{\alpha^N+*}(X^N).
\]
The theorem now follows as it did in the CW case.
\end{proof}

Again, we have a corresponding statement about Mackey functor-valued homology and cohomology.

\begin{corollary}
Let $N$ be a normal subgroup of $G$, let $\delta$ be a familial dimension function for $G$,
let $\alpha$ be a virtual representation of $G$,
let $\MackeyOp S$ be a covariant $G/N$-$\delta$-Mackey functor, and
let $\Mackey T$ be a contravariant $G/N$-$\delta$-Mackey functor.
Then, for $X$ any based $G$-space, we have
\[
 \Inf_{G/N}^G\MackeyOp H^{G/N,\delta}_{\alpha^N}(X^N; \MackeyOp S)
  \iso \MackeyOp H^{G,\delta}_\alpha(X; \Inf_{G/N}^G\MackeyOp S)
\]
and
\[
 \Inf_{G/N}^G\Mackey H_{G/N,\delta}^{\alpha^N}(X^N; \Mackey T)
  \iso \Mackey H_{G,\delta}^\alpha(X; \Inf_{G/N}^G\Mackey T).
\]
\end{corollary}

\begin{proof}
If $N\leq L\leq G$, we have the following commutative diagram:
\[
 \xymatrix{
  \sorb{L,\delta} \ar[d]_{i_L^G} \ar[r]^-{\Phi^N}
    & \sorb{L/N,\delta}  \ar[d]^{i_{L/N}^{G/N}} \\
  \sorb{G,\delta} \ar[r]_-{\Phi^N}
    & \sorb{G/N,\delta} 
 }
\]
From this we get that $(\Phi^N)^*(i_{L/N}^{G/N})^* \iso (i_{L}^{G})^*(\Phi^N)^*$, so
\[
 \Inf_{L/N}^L(\MackeyOp S|L/N) \iso (\Inf_{G/N}^G\MackeyOp S)|L
\]
and the following diagram commutes:
\[
 \xymatrix{
   \tilde H_{\alpha^N}^{G/N,\delta}(X^N\smsh G_+\smsh_L S^{-\delta(G/L)}; \MackeyOp S)
     \ar[r]^-\iso \ar[d]_\iso
     & \tilde H_{\alpha^N}^{L/N,\delta}(X^B; \MackeyOp S|L/N) \ar[d]^\iso \\
   \tilde H_\alpha^{G,\delta}(X\smsh G_+\smsh_L S^{-\delta(G/L)}; \Inf_{G/N}^G\MackeyOp S)
     \ar[r]_-\iso
     & \tilde H_\alpha^{L,\delta}(X; (\Inf_{G/N}^G\MackeyOp S)|L ).
 }
\]
This shows that
$ \Inf_{G/N}^G\MackeyOp H^{G/N,\delta}_{\alpha^N}(X^N; \MackeyOp S)
  \iso \MackeyOp H^{G,\delta}_\alpha(X; \Inf_{G/N}^G\MackeyOp S)$
and the argument for cohomology is the same.
\end{proof}

\begin{remarks}\label{rem:inductionFixed}
\begin{enumerate}\item[]
\item
The isomorphisms in Theorems~\ref{thm:inducedHomology}
and~\ref{thm:fixedSetIso} are compatible, in the sense that the composite of
the following isomorphisms is the identity
if $X$ is a $G/N$ space, $\alpha$ is a representation of $G/N$,
and $\Mackey T$ is a $G/N$-$\delta$-Mackey functor:
\begin{align*}
 \tilde H_\alpha^{G/N,\delta}(X;\Mackey T)
  &= \tilde H_{\alpha^N}^{G/N,\delta}(X^N;\Mackey T) \\
  &\iso \tilde H_\alpha^{G,\delta}(X;\Inf_{G/N}^G\Mackey T) \\
  &\iso \tilde H_\alpha^{G/N,\delta}(X;\theta^*\Inf_{G/N}^G\Mackey T) \\
  &\iso \tilde H_\alpha^{G/N,\delta}(X;\Mackey T).
\end{align*}
That the composite is the identity follows from considering the chain-level
isomorphisms. On that level we are simply composing the $(\Phi^N)_!$-$(\Phi^N)^*$ adjunction
and the $\theta_!$-$\theta^*$ adjunction, which gives the
$(\Phi^N\theta)_!$-$(\Phi^N\theta)^*$ adjunction, which is the identity.
The similar statement for homology is also true.

\item
The two isomorphisms combine to give a third isomorphism, if
$X$ is a $G$-space and $\Mackey T$ is a $G/N$-$\delta$-Mackey functor:
\begin{align*}
 \tilde H_\alpha^G(X;\Inf_{G/N}^G\Mackey T)
  &\iso \tilde H_{\alpha^N}^{G/N,\delta}(X^N; \Mackey T) \\
  &\iso \tilde H_{\alpha^N}^{G/N,\delta}(X^N; \theta^*\Inf_{G/N}^G\Mackey T) \\
  &\iso \tilde H_{\alpha^N}^{G,\delta}(X^N; \Inf_{G/N}^G\Mackey T).
\end{align*}
We get a similar isomorphism in homology.

\end{enumerate}
\end{remarks}

Using Theorems~\ref{thm:inducedHomology}
and~\ref{thm:fixedSetIso} we can define induction from $G/N$ to $G$
and restriction to fixed sets.

\begin{definition}\label{def:induction}
Let $\delta$ be an $N$-closed familial dimension function for $G$,
let $\alpha$ be a virtual representation of $G/N$,
let $\MackeyOp S$ be a covariant $G/N$-$\delta$-Mackey functor, and
let $\Mackey T$ be a contravariant $G/N$-$\delta$-Mackey functor.
If $Y$ is a based $G/N$-space,
we define {\em induction from $G/N$ to $G$} to be the composites
\[
 \epsilon^*\colon \tilde H^{G/N,\delta}_\alpha(Y;\MackeyOp S)
  \to \tilde H^{G/N,\delta}_\alpha(Y;\theta^*\theta_!\MackeyOp S)
  \iso \tilde H^{G,\delta}_\alpha(Y;\theta_!\MackeyOp S)
\]
and
\[
 \epsilon^*\colon \tilde H_{G/N,\delta}^\alpha(Y;\Mackey T)
  \to \tilde H_{G/N,\delta}^\alpha(Y;\theta^*\theta_!\Mackey T)
  \iso \tilde H_{G,\delta}^\alpha(Y;\theta_!\Mackey T).
\]
The first map in each case is induced by the unit of the $\theta_!$-$\theta^*$ adjunction.
\end{definition}

\begin{definition}\label{def:restrictToFixed}
Let $\delta$ be a familial dimension function for $G$,
let $\alpha$ be a virtual representation of $G$,
let $\MackeyOp S$ be a covariant $G$-$\delta$-Mackey functor,
and let $\Mackey T$ be a contravariant $G$-$\delta$-Mackey functor.
If $X$ is a based $G$-space, define
{\em restriction to the $N$-fixed set} to be the composites
\[
 (-)^N\colon \tilde H^{G,\delta}_{\alpha}(X;\MackeyOp S) \to
 \tilde H^{G,\delta}_{\alpha}(X;\Inf_{G/N}^G\MackeyOp S^N) \xrightarrow{\iso}
 \tilde H^{G/N,\delta}_{\alpha^N}(X^N;\MackeyOp S^N)
\]
and
\[
 (-)^N\colon \tilde H_{G,\delta}^{\alpha}(X;\Mackey T) \to
 \tilde H_{G,\delta}^{\alpha}(X;\Inf_{G/N}^G\Mackey T^N) \xrightarrow{\iso}
 \tilde H_{G/N,\delta}^{\alpha^N}(X^N;\Mackey T^N).
\]
The first map in each case is induced by the unit of the $(-)^N$-$\Inf_{G/N}^G$ adjunction.
\end{definition}

The reader should check that if, for example, $\MackeyOp S = \Inf_{G/N}^G \MackeyOp U$,
then the composite
\[
 \tilde H^{G,\delta}_\alpha(X;\Inf_{G/N}^G\MackeyOp U)
  \to \tilde H^{G/N,\delta}_{\alpha^N}(X^N;(\Inf_{G/N}^G\MackeyOp U)^N)
  \to \tilde H^{G/N,\delta}_{\alpha^N}(X^N;\MackeyOp U)
\]
agrees with the isomorphism of Theorem~\ref{thm:fixedSetIso}.

It follows from Theorem~\ref{thm:inducedHomology} that induction respects suspension, 
in the sense that, if $X$ is well-based and $W$ is a representation of $G/N$, then the following
diagram commutes, as does the similar one for cohomology:
\[
 \xymatrix{
  {\tilde H^{G/N,\delta}_{\alpha}(X;\MackeyOp S)} \ar[d]_-{\sigma^W} \ar[r]^{\epsilon^*}
   & {\tilde H^{G,\delta}_{\alpha}(X;\theta_!\MackeyOp S)} \ar[d]^{\sigma^{W}} \\
  {\tilde H^{G/N,\delta}_{\alpha+W}(\susp^W X;\MackeyOp S)} \ar[r]^-{\epsilon^*}
   & {\tilde H^{G,\delta}_{\alpha + W}(\susp^{W} X;\theta_!\MackeyOp S)}
  }
\]

It follows from Theorem~\ref{thm:fixedSetIso} that restriction to fixed sets respects suspension, 
in the sense that, if $X$ is well-based, the following
diagram commutes, as does the similar one for cohomology:
\[
 \xymatrix{
  {\tilde H^{G,\delta}_{\alpha}(X;\MackeyOp S)} \ar[d]_-{\sigma^W} \ar[r]^{(-)^N}
   & {\tilde H^{G/N,\delta}_{\alpha^N}(X^N;\MackeyOp S^N)} \ar[d]^{\sigma^{W^N}} \\
  {\tilde H^{G,\delta}_{\alpha+W}(\susp^W X;\MackeyOp S)} \ar[r]^-{(-)^N}
   & {\tilde H^{G/N,\delta}_{\alpha^N + W^N}(\susp^{W^N} X^N;\MackeyOp S^N)}
  }
\]

We also have the following relationship between induction and restriction to fixed sets.

\begin{proposition}
Let $\delta$ be an $N$-closed familial dimension function for $G$,
let $\alpha$ be a virtual representation of $G/N$,
let $\MackeyOp S$ be a covariant $G/N$-$\delta$-Mackey functor, and
let $\Mackey T$ be a contravariant $G/N$-$\delta$-Mackey functor.
Then, if $Y$ is a $G/N$-space, each of the following composites is the identity:
\[
 \tilde H_{\alpha}^{G/N,\delta}(Y;\MackeyOp S)
  \xrightarrow{\epsilon^*} \tilde H_{\alpha}^{G,\delta}(Y;\theta_!\MackeyOp S)
  \xrightarrow{(-)^N} \tilde H_{\alpha}^{G/N,\delta}(Y;(\theta_!\MackeyOp S)^N)
  = \tilde H_{\alpha}^{G/N,\delta}(Y;\MackeyOp S)
\]
and
\[
 \tilde H^{\alpha}_{G/N,\delta}(Y;\Mackey T)
  \xrightarrow{\epsilon^*} \tilde H^{\alpha}_{G,\delta}(Y;\theta_!\Mackey T)
  \xrightarrow{(-)^N} \tilde H^{\alpha}_{G/N,\delta}(Y;(\theta_!\Mackey T)^N)
  = \tilde H^{\alpha}_{G/N,\delta}(Y;\Mackey T)
\]
\end{proposition}

\begin{proof}
As in Remarks~\ref{rem:inductionFixed}, this comes down to
the fact that the $(\Phi^N\theta)_!$-$(\Phi^N\theta)^*$ adjunction
is the identity.
\end{proof}

Now we look at how the isomorphisms of Theorems~\ref{thm:inducedHomology}
and~\ref{thm:fixedSetIso} are represented on the spectrum level.

\begin{proposition}\label{prop:inductEMspectra}
Let $N$ be a normal subgroup of $G$, let $\delta$ be an $N$-closed familial dimension
function for $G$, let $\MackeyOp S$ be a covariant $G$-$\delta$-Mackey functor
and let $\Mackey T$ be a contravariant $G$-$\delta$-Mackey functor.
Then
\[
 (H^\delta\MackeyOp S)^N \hmtpc H^\delta(\theta^*\MackeyOp S)
\]
and
\[
 (H_\delta\Mackey T)^N \hmtpc H_\delta(\theta^*\Mackey T).
\]
\end{proposition}

\begin{proof}
Consider $(H^\delta\MackeyOp S)^N$ first. For $K/N\in\F(\delta|G/N)$ we have
\begin{align*}
 \Mackey\pi_n^{G/N,\Lie-\delta}(G/K,\Lie-\delta)(&(H^\delta\MackeyOp S)^N) \\
  &= [G_+\smsh_K S^{-(\Lie(G/K)-\delta(G/K))+n}, (H^\delta\MackeyOp S)^N ]_{G/N} \\
  &\iso [G_+\smsh_K S^{-(\Lie(G/K)-\delta(G/K))+n}, H^\delta\MackeyOp S ]_{G} \\
  &\iso \begin{cases}
  			\MackeyOp S(G/K,\delta) &\text{if $n=0$} \\
			0 &\text{if $n\neq 0$}
		\end{cases}
\end{align*}
Thus, we have
\[
 \Mackey\pi_n^{G/N,\Lie-\delta}((H^\delta\MackeyOp S)^N)
  \iso \begin{cases}
  			\theta^*\MackeyOp S &\text{if $n=0$} \\
			0 &\text{if $n\neq 0$.}
	   \end{cases}
\]
We also have the equivalences
\begin{align*}
 (H^\delta\MackeyOp S)^N\smsh E\F(\delta|G/N)_+
  &\hmtpc (H^\delta\MackeyOp S \smsh \epsilon^* E\F(\delta|G/N)_+)^N \\
  &\hmtpc (H^\delta\MackeyOp S \smsh E\F(\delta)_+ \smsh \epsilon^* E\F(\delta|G/N)_+)^N \\
  &\hmtpc (H^\delta\MackeyOp S \smsh E\F(\delta)_+)^N \\
  &\hmtpc (H^\delta\MackeyOp S)^N.
\end{align*}
The $G$-homotopy equivalence 
$E\F(\delta)_+ \smsh \epsilon^* E\F(\delta|G/N)_+ \hmtpc E\F(\delta)_+$ follows because,
if $K\in\F(\delta)$, then
\[
 (\epsilon^* E\F(\delta|G/N)_+)^K = (\epsilon^* E\F(\delta|G/N)_+)^{KN} \hmtpc S^0,
\]
using the assumption that $\delta$ is $N$-closed.
So, we have verified the conditions that characterize the Eilenberg-Mac\,Lane spectrum and conclude that
\[
 (H^\delta\MackeyOp S)^N \hmtpc H^\delta(\theta^*\MackeyOp S).
\]

Now consider $(H_\delta\Mackey T)^N$. The computation of the homotopy groups is the same and
we get that
\[
 \Mackey\pi_n^{G/N,\delta}((H_\delta\Mackey T)^N)
  \iso \begin{cases}
  			\theta^*\Mackey T &\text{if $n=0$} \\
			0 &\text{if $n\neq 0$}.
	   \end{cases}
\]
On the other hand, we have
\begin{align*}
 F(E\F(\delta|G/N)_+, (H_\delta\Mackey T)^N)
  &\hmtpc F(\epsilon^*E\F(\delta|G/N)_+, H_\delta\Mackey T)^N \\
  &\hmtpc F(\epsilon^*E\F(\delta|G/N)_+\smsh E\F(\delta)_+, H_\delta\Mackey T)^N \\
  &\hmtpc F(E\F(\delta)_+, H_\delta\Mackey T)^N \\
  &\hmtpc (H_\delta\Mackey T)^N.
\end{align*}
From the characterization of the Eilenberg-Mac\,Lane spectrum we now conclude that
\[
 (H_\delta\Mackey T)^N \hmtpc H_\delta(\theta^*\Mackey T).
\]
\end{proof}

The isomorphisms of Theorem~\ref{thm:inducedHomology} are represented as follows,
for a $G/N$-space $Y$ and a virtual representation $\alpha$ of $G/N$:
\begin{align*}
 [S^\alpha, H^\delta\MackeyOp S \smsh Y]_G
  &\iso [S^\alpha, (H^\delta\MackeyOp S \smsh Y)^N]_{G/N} \\
  &\iso [S^\alpha, (H^\delta\MackeyOp S)^N \smsh Y]_{G/N} \\
  &\iso [S^\alpha, H^\delta(\theta^*\MackeyOp S)\smsh Y]_{G/N}
\end{align*}
for homology, and
\begin{align*}
 [\susp^\infty Y, \susp^\alpha H_\delta\Mackey T]_G
  &\iso [\susp^\infty Y, (\susp^\alpha H_\delta\Mackey T)^N]_{G/N} \\
  &\iso [\susp^\infty Y, \susp^\alpha (H_\delta\Mackey T)^N]_{G/N} \\
  &\iso [\susp^\infty Y, \susp^\alpha H_\delta(\theta^*\Mackey T)]_{G/N}
\end{align*}
for cohomology.

\begin{proposition}\label{prop:fixedEMSpectra}
Let $\delta$ be a familial dimension function for $G$,
let $N$ be a normal subgroup of $G$, let $\Mackey T$ be a contravariant
$G/N$-$\delta$-Mackey functor, and let $\MackeyOp S$ be a covariant
$G/N$-$\delta$-Mackey functor. Then
\[
 \epsilon^\sharp H^\delta\MackeyOp S \smsh \tE\F[N] \hmtpc
   H^\delta\Inf_{G/N}^G\MackeyOp S
\]
and
\[
 \epsilon^\sharp H_\delta\Mackey T\smsh \tE\F[N] \hmtpc
   H_\delta\Inf_{G/N}^G\Mackey T.
\]
Therefore,
\[
 (H^\delta\Inf_{G/N}^G\MackeyOp S)^N \hmtpc \Phi^N H^\delta\Inf_{G/N}^G\MackeyOp S 
   \hmtpc H^\delta\MackeyOp S
\]
and
\[
 (H_\delta\Inf_{G/N}^G\Mackey T)^N \hmtpc \Phi^N H_\delta\Inf_{G/N}^G\Mackey T \hmtpc H_\delta\Mackey T.
\]
\end{proposition}

\begin{proof}
First consider $\epsilon^\sharp H^\delta\MackeyOp S \smsh \tE\F[N]$.
We know that 
\[
 H^\delta\MackeyOp S \hmtpc H^\delta\MackeyOp S\smsh E\F(\delta|G/N)_+,
\]
so 
\[
 \epsilon^\sharp H^\delta\MackeyOp S \hmtpc \epsilon^\sharp H^\delta\MackeyOp S
   \smsh \epsilon^* E\F(\delta|G/N)_+.
\]
Comparison of fixed points shows that
$\epsilon^* E\F(\delta|G/N)_+\smsh\tE\F[N] \hmtpc \tE\F[N] \smsh E\F(\delta)_+$, so
\[
 \epsilon^\sharp H^\delta\MackeyOp S \smsh \tE\F[N] \hmtpc
  \epsilon^\sharp H^\delta\MackeyOp S \smsh \tE\F[N] \smsh E\F(\delta)_+.
\]
If $K$ does not contain $N$, then $\tE\F[N]$ is contractible as a $K$-space, hence
\[
 \Mackey\pi_n^{G,\Lie-\delta}(\epsilon^\sharp H^\delta\MackeyOp S \smsh \tE\F[N])(G/K,\Lie-\delta)
  = 0
\]
for $K\in\F(\delta)$ not containing $N$. On the other hand, if $N\leq K$ and $K\in\F(\delta)$, we have
\begin{align*}
 \Mackey\pi_n^{G,\Lie-\delta}(&\epsilon^\sharp H^\delta\MackeyOp S \smsh \tE\F[N])(G/K,\Lie-\delta) \\
  &= [G_+\smsh_K S^{-(\Lie(G/K)-\delta(G/K))+n},
    \epsilon^\sharp H^\delta\MackeyOp S \smsh \tE\F[N] ]_G \\
  &\iso [G_+\smsh_K S^{-(\Lie(G/K)-\delta(G/K))+n},
    (\epsilon^\sharp H^\delta\MackeyOp S \smsh \tE\F[N])^N ]_{G/N} \\
  &\iso [G_+\smsh_K S^{-(\Lie(G/K)-\delta(G/K))+n},
    H^\delta\MackeyOp S ]_{G/N} \\
  &\iso \begin{cases}
  			\MackeyOp S(G/K,\delta) &\text{if $n=0$} \\
			0 &\text{if $n\neq 0$.}
		\end{cases}
\end{align*}
Therefore,

\[
 \Mackey\pi_n^{G,\Lie-\delta}(\epsilon^\sharp H^\delta\MackeyOp S \smsh \tE\F[N])
 \iso \begin{cases}
  			\Inf_{G/N}^G\MackeyOp S &\text{if $n=0$} \\
			0 &\text{if $n\neq 0$,}
	  \end{cases}
\]
and we have verified the conditions that characterize the Eilenberg-Mac\,Lane spectrum, so
\[
 \epsilon^\sharp H^\delta\MackeyOp S \smsh \tE\F[N] \hmtpc
   H^\delta\Inf_{G/N}^G\MackeyOp S.
\]

Turning to $\epsilon^\sharp H_\delta\Mackey T\smsh \tE\F[N]$, we calculate exactly as above
that
\[
 \Mackey\pi_n^{G,\delta}(\epsilon^\sharp H_\delta\Mackey T\smsh \tE\F[N]) 
 \iso \begin{cases}
  			\Inf_{G/N}^G\Mackey T &\text{if $n=0$} \\
			0 &\text{if $n\neq 0$,}
	  \end{cases}
\]
Further, for an arbitrary $G$-spectrum $E$ we have the following isomorphisms,
which follow from ones we've shown previously together with
\cite[II.9.2]{LMS:eqhomotopy} and \cite[II.9.6]{LMS:eqhomotopy}:
\begin{align*}
 [E\smsh E\F(\delta)_+,{} &\epsilon^\sharp H_\delta\Mackey T\smsh \tE\F[N] ]_G \\
  &\iso [E\smsh E\F(\delta)_+\smsh\tE\F[N], \epsilon^\sharp H_\delta\Mackey T\smsh \tE\F[N] ]_G \\
  &\iso [E\smsh \epsilon^*E\F(\delta|G/N)_+\smsh\tE\F[N], \epsilon^\sharp H_\delta\Mackey T\smsh \tE\F[N] ]_G \\
  &\iso [\Phi^N(E\smsh \epsilon^*E\F(\delta|G/N)_+), \Phi^N\epsilon^\sharp H_\delta\Mackey T]_{G/N} \\
  &\iso [\Phi^N E \smsh E\F(\delta|G/N)_+, H_\delta\Mackey T]_{G/N} \\
  &\iso [\Phi^N E, H_\delta\Mackey T]_{G/N} \\
  &\iso [E, \epsilon^\sharp H_\delta\Mackey T\smsh \tE\F[N] ]_G
\end{align*}
This shows that
\[
 \epsilon^\sharp H_\delta\Mackey T\smsh \tE\F[N] \hmtpc
  F(E\F(\delta)_+, \epsilon^\sharp H_\delta\Mackey T\smsh \tE\F[N]),
\]
hence we've verified the conditions that characterize the Eilenberg-Mac\,Lane spectrum, so
\[
 \epsilon^\sharp H_\delta\Mackey T\smsh \tE\F[N] \hmtpc
   H_\delta\Inf_{G/N}^G\Mackey T.
\]
The last statement of the proposition follows on taking the $N$-fixed points of
the equivalences already shown.
\end{proof}

The isomorphisms of Theorem~\ref{thm:fixedSetIso} are then represented as
\begin{align*}
 [S^\alpha, H^{\delta}\Inf_{G/N}^G\MackeyOp S \smsh X]_G
  &\iso [\Phi^N(S^{\alpha}), \Phi^N(H^{\delta}\Inf_{G/N}^G\MackeyOp S\smsh X)]_{G/N} \\
  &\iso [S^{\alpha^N}, \Phi^N H^{\delta}\Inf_{G/N}^G\MackeyOp S \smsh X^N]_{G/N} \\
  &\iso [S^{\alpha^N}, H^{\delta}\MackeyOp S \smsh X^N]_{G/N}
\end{align*}
and
\begin{align*}
 [\susp_G^\infty X, \susp^\alpha H_\delta\Inf_{G/N}^G\Mackey T]_G
  &\iso [\Phi^N(\susp_G^\infty X), \Phi^N(\susp^\alpha H_\delta\Inf_{G/N}^G\Mackey T)]_{G/N} \\
  &\iso [\susp_{G/N}^\infty X^N, \susp^{\alpha^N} H_\delta\Mackey T]_{G/N}.
\end{align*}
The first isomorphism in each case follows because $H_\delta\Inf_{G/N}^G\Mackey T$, for example,
is concentrated over $N$ in the language of \cite[\S II.9]{LMS:eqhomotopy}.

Finally, a word about taking $K$-fixed sets if $K$ is not normal in $G$.
The natural way to do this is first to restrict from $G$ to $NK$
(assuming $NK\in \F(\delta)$) and
then take $K$ fixed points, defining $a^K = (a|NK)^K$. 
In cohomology, for example, this gives
\[
 (-)^K\colon \tilde H_{G,\delta}^\alpha(X;\Mackey T)
    \to \tilde H_{WK,\delta}^{\alpha^K-\delta(G/NK)^K}(X^K;(\Mackey T|NK)^K),
\]
where $WK = NK/K$.
However, because $(G/K)^K = NK/K$, we have that $\Lie(G/NK)^K = 0$ and
so also $\delta(G/NK)^K = 0$. Therefore, the possible shift in dimensions goes away
and restriction to fixed sets defines maps
\begin{align*}
 (-)^K\colon \tilde H^{G,\delta}_\alpha(X;\MackeyOp S)
    &\to \tilde H^{WK,\delta}_{\alpha^K}(X^K;\MackeyOp S^K)\qquad\text{and}\\
 (-)^K\colon \tilde H_{G,\delta}^\alpha(X;\Mackey T)
    &\to \tilde H_{WK,\delta}^{\alpha^K}(X^K;\Mackey T^K),
\end{align*}
where we write $\MackeyOp S^K = (\MackeyOp S|NK)^K$ and similarly for
$\Mackey T^K$.

\subsection{Subgroups of quotient groups}

We now look at how induction and restriction to fixed sets interact
with restriction to subgroups. Consider restriction to fixed sets, for example.
Suppose that $N$ is a normal subgroup of $G$, $N\leq L \leq G$, and
$\delta$ is familial with $L\in\F(\delta)$.
We would like to say that $(a|L)^N = a^N|(L/N)$, but, at first glance,
these elements appear to live in homology groups with different coefficient systems.
In cohomology, for example, if $a\in \tilde H_{G,\delta}^\alpha(X;\Mackey T)$,
then
\[
 (a|L)^N \in \tilde H_{L/N,\delta}^{\alpha^N-\delta(G/L)}(X^N; (\Mackey T|L)^N)
\]
while
\[
 a^N|(L/N) \in \tilde H_{L/N,\delta}^{\alpha^N-\delta(G/L)}(X^N; \Mackey T^N|(L/N)).
\]
We shall use Proposition~\ref{prop:commutingInducedMaps} to show that
$(\Mackey T|L)^N \iso \Mackey T^N|(L/N)$ and that the elements in question coincide
under this identification.

The proofs for both induction and restriction to fixed sets depend on a more in-depth understanding of the
stable orbit category.
The following result will be a special case of Theorem~\ref{thm:StableMapsOrbitsOverB}
so we defer the proof until then;
the case $\delta=0$ is a special case of
\cite[V.9.4]{LMS:eqhomotopy}.

\begin{proposition}\label{prop:stableOrbitCategory}
Let $\delta$ be a familial dimension function for $G$.
Then $\sorb{G,\delta}((G/H,\delta), (G/K,\delta))$ 
is the free abelian group generated by the equivalence
classes of diagrams of orbits of the form
\[
 G/H \xleftarrow{p} G/J \xrightarrow{q} G/K,
\]
where $J\leq H$, $p$ is the projection, $q$ is a (space-level) $G$-map with, say,
$q(eJ) = gK$, and $\delta(N_HJ/J) \iso \Lie(N_HJ/J)$ and $\delta(N_{K^g}J/J) = 0$. 
Two such diagrams are equivalent if
there is a diagram of the following form in which
$G/J\to G/J'$ is a $G$-homeomorphism, the left triangle
commutes, and the right triangle commutes up to $G$-homotopy:
\[
 \xymatrix@R-5ex{
  & G/J \ar[dl]_{p} \ar[dr]^{q} \ar[dd] \\
  G/H && G/K \\
  & G/J' \ar[ul]^{p'} \ar[ur]_{q'}
 }
\]
\qed
\end{proposition}

Explicitly, a pair of maps $(p,q)$ as above corresponds to the stable map constructed as follows.
Take $V$ so large that $H/J$ embeds in $V-\delta(G/H)$ as an $H$-space. We then have the following
composite, in which the first map is the collapse map and the last is induced by $q$:
\begin{align*}
 G_+\smsh_H S^{V-\delta(G/H)}
  &\to G_+\smsh_J S^{V-\delta(G/H)-\Lie(H/J)} \\
  &\includesin G_+\smsh_J S^{V-\delta(G/H)-\delta(H/J)} \\
  &= G_+\smsh_J S^{V-\delta(G/J)} \\
  &\includesin G_+\smsh_J S^{V-\delta(G/K^g)} \\
  &\to G_+\smsh_K S^{V-\delta(G/K)}.
\end{align*}
Note that the first inclusion would be null-homotopic if
\[
 [\Lie(H/J)-\delta(H/J)]^J = \Lie(N_HJ/J)-\delta(N_HJ/J) \neq 0
\]
and the second would be
null-homotopic if
\[
 \delta(K^g/J)^J = \delta(N_{K^g}J/J) \neq 0.
\]
Here, we use that $\delta(H/J)^J \iso \delta(N_H J/J)$. This is clear for $\delta=\Lie$
because $(H/J)^J = N_H J/J$. In general, it follows from the fact that
$\delta(H/J)$ is isomorphic to a $J$-subspace of $\Lie(H/J)$ and
$\delta(N_H J/J)$ is isomorphic to a subspace of $\Lie(N_H J/J)$, which imply
that
\[
 \delta(H/J)^J \iso \delta(H/N_H J)^J \dirsum \delta(N_H J/J)^J = 0 \dirsum \delta(N_H J/J).
\]
The stable map corresponding to $(p,q)$ is then the suspension of the space-level
composite map above.

Because the corresponding stable map is 0 if $G/J$ does not meet the dimension requirements,
we shall allow ourselves to write $[G/H\from G/J\to G/K]$ for
a stable map $(G/H,\delta)\to (G/K,\delta)$, with any $J$, with the understanding that the element
may be 0.

Now, let's turn to the interaction of induction and restriction to subgroups.
If $N\leq L\leq G$ and $L\in\F(\delta)$, we have the following commutative diagram (which we've used before):
\[
 \xymatrix{
  \sorb{L/N,\delta} \ar[r]^-\theta \ar[d]_{i_{L/N}^{G/N}}
     & \sorb{L,\delta} \ar[d]^{i_L^G} \\
  \sorb{G/N,\delta} \ar[r]_-\theta 
     & \sorb{G,\delta}
 }
\]
By Proposition~\ref{prop:commutingInducedMaps}, we have a natural map
$\xi\colon \theta_! i^*\to i^*\theta_!$, where we write $i$ for either
$i_L^G$ or $i_{L/N}^{G/N}$, as appropriate.

\begin{lemma}\label{lem:inductFixedIso}
The natural transformation $\xi\colon \theta_! i^*\to i^*\theta_!$
is an isomorphism.
\end{lemma}

\begin{proof}
We give the argument for contravariant Mackey functors; the proof for
covariant functors is similar or we can appeal to duality (replacing $\delta$
with $\Lie-\delta$).

By Proposition~\ref{prop:commutingInducedMaps}, the result will follow if we show that
\[
 \xi\colon \int\nolimits^{L/J\in\sorb{L/N,\delta}}
  \sorb{G/N,\delta}(G/J, G/K) \tensor \sorb{L,\delta}(L/H,L/J)
   \to \sorb{G,\delta}(G/H, G/K),
\]
given by $\xi(f\tensor g) = \epsilon^\sharp f \circ (G_+\smsh_L g)$,
is an isomorphism for all $L/H$ in $\sorb{L,\delta}$ and $G/K$ in $\sorb{G/N,\delta}$.

Define a map $\zeta$ inverse to $\xi$ as follows:
Let $[G/H \xleftarrow{p} G/M \xrightarrow{q} G/K]$ be a generator of $\sorb{G,\delta}(G/H, G/K)$,
so $M\leq H\leq L$ and $p$ is the projection. Then $p = G\times_L p'$ where $p'\colon L/M\to L/H$.
By assumption, $N\leq K$, so $q$ factors
as $G/M \to G/MN \to G/K$. We let
\[
 \zeta[G/H\from G/M\to G/K] = [G/MN\to G/K] \tensor [L/H \xleftarrow{p'} L/M \to L/MN].
\]
(Notice that $MN\leq L$ because $L$ contains both $M$ and $N$.)
Clearly, $\xi\circ\zeta$ is the identity. On the other hand, a typical element in the coend
is a sum of elements of the form
\[
 [G/J \xleftarrow{p} G/M \to G/K] \tensor g,
\]
where $N\leq M\leq J$ and $p = G\times_L p'$. We then have
\begin{align*}
 [G/J \xleftarrow{p} G/M &{}\to G/K] \tensor g \\
  &\sim [G/M\to G/K] \tensor [L/J \xleftarrow{p'} L/M]\circ g \\
  &= \sum_i [G/M\to G/K] \tensor [L/H \from L/P_i \to L/M] \\
  &= \sum_i [G/M\to G/K] \tensor [L/H \from L/P_i \to L/P_iN \to L/M] \\
  &\sim \sum_i [G/P_iN\to G/M \to G/K] \tensor [L/H \from L/P_i \to L/P_iN]
\end{align*}
which is in the image of $\zeta$. So, $\zeta$ is an epimorphism, hence an isomorphism
and the inverse of $\xi$.
\end{proof}

\begin{proposition}\label{prop:inductionfixedsets}
Let $N$ be a normal subgroup of $G$,
let $\delta$ be an $N$-closed familial dimension function for $G$,
let $\alpha$ be a virtual representation of $G/N$,
let $\MackeyOp S$ be a covariant $G/N$-$\delta$-Mackey functor, and
let $\Mackey T$ be a contravariant $G/N$-$\delta$-Mackey functor.
Let $N\leq L\leq G$ with $L\in\F(\delta)$.
If $Y$ is a based $G/N$-space and $y\in\tilde H_\alpha^{G/N,\delta}(Y;\MackeyOp S)$, then
\[
 \epsilon^*(y|L/N) = (\epsilon^* y)|L
  \in \tilde H_\alpha^{L,\delta}(Y; \theta_!(\MackeyOp S|L/N))
  \iso \tilde H_\alpha^{L,\delta}(Y; (\theta_!\MackeyOp S)|L).
\]
Similarly, if $y\in \tilde H^\alpha_{G/N,\delta}(Y;\Mackey T)$, then
\[
 \epsilon^*(y|L/N) = (\epsilon^* y)|L
  \in \tilde H^\alpha_{L,\delta}(Y; \theta_!(\Mackey T|L/N))
  \iso \tilde H^\alpha_{L,\delta}(Y; (\theta_!\Mackey T)|L).
\]
\end{proposition}

\begin{proof}
We concentrate on cohomology; the proof for homology is similar.
We need to show that the following diagram commutes:
\[
 \xymatrix{
   \tilde H^\alpha_{G/N,\delta}(Y;\Mackey T) \ar[r]^-{\epsilon^*} \ar[d]_{-|L/N}
     & \tilde H^\alpha_{G,\delta}(Y;\theta_!\Mackey T) \ar[d]^{-|L} \\
   \tilde H^\alpha_{L/N,\delta}(Y;\Mackey T|L/N) \ar[dr]_{\epsilon^*}
     & \tilde H^\alpha_{L,\delta}(Y;(\theta_!\Mackey T)|L) \\
   & \tilde H^\alpha_{L,\delta}(Y;\theta_!(\Mackey T|L/N)) \ar[u]_{\xi_*}^\iso
 }
\]
We can expand the top two rows of this diagram to get the following commutative diagram:
\[
 \xymatrix@C-1em{
  \tilde H^\alpha_{G/N,\delta}(Y;\Mackey T) \ar[r] \ar[d]
    & \tilde H^\alpha_{G/N,\delta}(Y;\theta^*\theta_!\Mackey T) \ar[r]^-\iso \ar[d]
    & \tilde H^\alpha_{G,\delta}(Y;\theta_!\Mackey T) \ar[d] \\
  \tilde H^\alpha_{G/N,\delta}(G_+\smsh_L Y;\Mackey T) \ar[r] \ar[d]_\iso
    & \tilde H^\alpha_{G/N,\delta}(G_+\smsh_L Y;\theta^*\theta_!\Mackey T) \ar[r]^-\iso \ar[d]_\iso
    & \tilde H^\alpha_{G,\delta}(G_+\smsh_L Y;\theta_!\Mackey T) \ar[d]^\iso \\
  \tilde H_*^{L/N,\delta}(Y; i^*\Mackey T) \ar[r]
    & \tilde H^\alpha_{L/N,\delta}(Y;i^*\theta^*\theta_!\Mackey T) \ar[r]_-\iso
    & \tilde H^\alpha_{L,\delta}(Y; i^*\theta_!\Mackey T)
 }
\]
Here, we write $i^*$ instead of $-|L$ or $-|L/N$ for clarity. 
The two leftmost squares commute by naturality while the
two rightmost squares commute because of the
commutation $i^*\theta^* \iso \theta^* i^*$.
It remains to show, then that the following diagram commutes:
\[
 \xymatrix{
  \tilde H_*^{L/N,\delta}(Y; i^*\Mackey T) \ar[r] \ar[dr]
    & \tilde H^\alpha_{L/N,\delta}(Y;i^*\theta^*\theta_!\Mackey T) \ar[r]^-\iso
    & \tilde H^\alpha_{L,\delta}(Y; i^*\theta_!\Mackey T) \\
  & \tilde H^\alpha_{L/N,\delta}(Y; \theta^*\theta_!i^*\Mackey T) \ar[r]_-\iso
    & \tilde H^\alpha_{L,\delta}(Y; \theta_! i^*\Mackey T) \ar[u]_{\xi_*}
 }
\]
That this commutes follows from the following diagram of natural transformations,
which commutes by the definition of $\xi$ and properties of adjunctions:
\[
 \xymatrix{
  i^* \ar[r]^{i^*\eta} \ar[d]_{\eta i^*}
    & i^* \theta^*\theta_! \ar[d]^\iso \\
  \theta^*\theta_! i^* \ar[r]_{\theta^* \xi} & \theta^* i^*\theta_!
 }
\]
\end{proof}

We now consider the interaction of restriction to fixed sets and restriction to subgroups.
Again, let $N\leq L\leq G$ with $L\in\F(\delta)$.
We have the following commutative diagram:
\[
 \xymatrix{
  \sorb{L,\delta} \ar[d]_{i_L^G} \ar[r]^{\Phi^N} & \sorb{L/N,\delta} \ar[d]^{i_{L/N}^{G/N}} \\
  \sorb{G,\delta} \ar[r]_{\Phi^N} & \sorb{G/N,\delta}
 }
\]
By Proposition~\ref{prop:commutingInducedMaps}, we have a natural map
\[
 \xi\colon \Phi^N_!i^* \to i^*\Phi^N_!.
\]
From Proposition~\ref{prop:stableOrbitCategory} we get 
the following explicit description of $\Phi^N$.
For simplicity of notation we write $G/H$ for the object $(G/H,\delta)$ in
$\sorb{G,\delta}$.

\begin{proposition}\label{prop:stableFixedPoints}
Consider the map
\[
 \Phi^N\colon \sorb{G,\delta}(G/H,G/K) \to \sorb{G/N,\delta}((G/H)^N,(G/K)^N).
\]
If either $N\not\leq H$ or $N\not\leq K$, the target is the trivial group.
If $N\leq H$ and $N\leq K$, then,
on generators, $\Phi^N$ is given by
\begin{align*}
 \Phi^N[G/H\from G/J\to G/K]
  &= [(G/H)^N \from (G/J)^N \to (G/K)^N] \\
  &= \begin{cases}
       [G/H\from G/J\to G/K] & \text{if $N\leq J$} \\
       0 & \text{otherwise.}
     \end{cases}
\end{align*}
Therefore, it is a split epimorphism with kernel generated by those
diagrams $[G/H\from G/J\to G/K]$ with $N\not\leq J$.
\end{proposition}

\begin{proof}
That $\Phi^N$ takes the suspension of a map of spaces to the suspension of
its $N$-fixed set map is the naturality of the equivalence
$\Phi^N(\susp_G^\infty X) \hmtpc \susp_{G/N}^\infty X^N$.
Applying this to the space-level description above of the stable map
corresponding to a diagram $[G/H\from G/J\to G/K]$, we get the
description of $\Phi^N$ on generators given in the statement of the proposition.

The only other comment necessary is that, if $N\leq J$, the conditions and the equivalence
relations on diagrams $G/H\from G/J\to G/K$ are the same whether we consider them
as diagrams of $G/N$-orbits or $G$-orbits. Therefore, the generators of
$\sorb{G/N,\delta}((G/H)^N,(G/K)^N)$ are in one-to-one correspondence with those
generators of $\sorb{G,\delta}(G/H,G/K)$ having the form
$G/H\from G/J\to G/K$ with $N\leq J$.
\end{proof}

Now we can prove that $\xi\colon \Phi^N_!i^* \to i^*\Phi^N_!$ is an isomorphism.

\begin{lemma}\label{lem:fixedAndSubgroupEquality}
If $N$ is a normal subgroup of $G$, $N\leq L\leq G$, then
$\xi\colon \Phi^N_!i^* \to i^*\Phi^N_!$
is an isomorphism on both covariant and contravariant 
$G/N$-$\delta$-Mackey functors.
\end{lemma}

\begin{proof}
We give the argument for contravariant Mackey functors; the proof for covariant functors is similar
or we can appeal to duality (replacing $\delta$ with $\Lie-\delta$).

By Proposition~\ref{prop:commutingInducedMaps}, the result will follow if we show that
\begin{multline*}
 \xi\colon \int\nolimits^{L/J\in\sorb{L,\delta}} \sorb{G,\delta}(G/J,G/K)\tensor \sorb{L/N,\delta}(L/H,(L/J)^N) \\
   \to \sorb{G/N,\delta}(G/H,(G/K)^N),
\end{multline*}
given by $\xi(f\tensor g) = \Phi^Nf\circ(G_+\smsh_L g)$, is an isomorphism for all
$L/H$ in $\sorb{L/N,\delta}$ and $G/K$ in $\sorb{G,\delta}$.

If $N\not\leq K$, then $(G/K)^N = *$ and the target of $\xi$ is 0.
For the source,
consider a typical generator $f = [G/J \xleftarrow{p} G/M\to G/K]$ of $\sorb{G,\delta}(G/J,G/K)$.
Because $p$ is a projection and $M\leq J\leq L$, we can write $p = G\times_L p'$
where $p'\colon L/M\to L/J$ is a projection.
Because $M$ is subconjugate to $K$, $N\not\leq M$, hence $(L/M)^N = 0$.  
Therefore, for any $g$,
\begin{align*}
 f\tensor g 
   &= [G/J \xleftarrow{p} G/M \to G/K] \tensor g \\
   &\sim [G/M \to G/K] \tensor [L/J \xleftarrow{p'} L/M]^N\circ g \\
   &= [G/M \to G/K] \tensor [(L/J)^N \from *] \circ g \\
   &= 0
\end{align*}
in the coend. Therefore, the coend is 0 and $\xi$ is an isomorphism in this case.

So, assume that $N\leq K$ so $(G/K)^N = G/K$.
Define a map $\zeta$ inverse to $\xi$ as follows: If $[G/H\xleftarrow{p} G/J \to G/K]$
is a generator of $\sorb{G/N,\delta}(G/H,G/K)$ with $N\leq J\leq H \leq L$, 
write $p = G\times_L p'$ where $p'\colon L/J\to L/H$, and let
\[
 \zeta[G/H\xleftarrow{p} G/J \to G/K] = [G/J \to G/K] \tensor [L/H \xleftarrow{p'} L/J].
\]
Clearly, $\xi\circ\zeta$ is the identity.
On the other hand, a typical element in the coend is a sum of elements of the form
\[
 [G/J\xleftarrow{p}G/M\to G/K]\tensor g,
\]
where we may assume $N\leq J$ because otherwise
the element would live in a 0 group. 
Writing $p = G\times_L p'$, we then have
\begin{align*}
 [G/J\xleftarrow{p}G/M\to G/K]\tensor g
  &\sim [G/M\to G/K]\tensor [L/J \xleftarrow{p'} L/M]\circ g \\
  &= \sum_i [G/M\to G/K]\tensor [L/H \from L/P_i \to L/M] \\
  &\sim \sum_i [G/P_i \to G/M \to G/K] \tensor [L/H\from L/P_i]
\end{align*}
which is in the image of $\zeta$. So, $\zeta$ is an epimorphism, hence an isomorphism
and the inverse of $\xi$.
\end{proof}

\begin{proposition}\label{prop:fixedandsubgroup}
Let $N$ be a normal subgroup of $G$, let $\delta$ be a familial dimension function for $G$,
let $\alpha$ be a virtual representation of $G$,
let $\MackeyOp S$ be a covariant $G$-$\delta$-Mackey functor, and
let $\Mackey T$ be a contravariant $G$-$\delta$-Mackey functor.
Let $N\leq L\leq G$ with $L\in\F(\delta)$.
If $X$ is a based $G$-space and $x\in \tilde H_\alpha^{G,\delta}(X;\MackeyOp S)$, then
\[
 (x|L)^N = x^N|L/N \in 
  \tilde H_{\alpha^N-\delta(G/L)}^{L/N,\delta}(X^N; (\MackeyOp S|L)^N) \iso
  \tilde H_{\alpha^N-\delta(G/L)}^{L/N,\delta}(X^N; \MackeyOp S^N|L/N).
\]
Similarly, if $x\in \tilde H^\alpha_{G,\delta}(X;\Mackey T)$, then
\[
 (x|L)^N = x^N|L/N \in 
  \tilde H^{\alpha^N-\delta(G/L)}_{L/N,\delta}(X^N; (\Mackey T|L)^N) \iso
  \tilde H^{\alpha^N-\delta(G/L)}_{L/N,\delta}(X^N; \Mackey T^N|L/N).
\]
\end{proposition}

\begin{proof}
The proof is the same as the proof of Proposition~\ref{prop:inductionfixedsets},
but using Lemma~\ref{lem:fixedAndSubgroupEquality}.
\end{proof}

\section{Products}\label{sec:vproducts}

We now turn to various pairings, including cup products, evaluation maps, and cap products.
There are quite a few variations of such pairings available.

\subsection{Cup products}

The product of two orbits of $G$ is not generally an orbit of $G$; 
if $G$ is infinite the product
is not generally even a disjoint union of orbits. Thus, if $X$ and $Y$ are
$G$-CW complexes, their product does not have a canonical $G$-CW
structure. However, it does have a canonical $(G\times G)$-CW
structure. Following where the geometry leads us, we therefore
start with an {\em external} cup product, pairing a based $G$-space $X$
and a based $K$-space $Y$ to get a based $(G\times K)$-space $X\smsh Y$. In
the case of two $G$-spaces, we can then internalize by restricting to the
diagonal $G\leq G\times G$ using the restriction to subgroups
discussed in detail in the preceding section.

We begin by defining the appropriate tensor product of Mackey
functors (see also \cite{Dr:reps} and \cite{Le:Greenfunctors}).

\begin{definition}\label{def:mackeyproducts}
Let $G$ and $K$ be two compact Lie groups.
Let $\delta$ be a dimension function
for $G$, let $\epsilon$ be a dimension function for $K$,
and let $\delta\times\epsilon$ denote their product
as in Definition~\ref{def:productDimFcn}.
Let $\zeta$ be a dimension function for $G\times K$ with
$\zeta\dimpred \delta\times\epsilon$,
as in Definition~\ref{def:dimensionpo}.
\begin{enumerate}
\item
Let $\sorb{G,\delta}\tensor\sorb{K,\epsilon}$ denote the preadditive category whose objects
are pairs of objects $(a,b)$, as in the product category, with
\[
 (\sorb{G,\delta}\tensor\sorb{K,\epsilon})((a,b),(c,d)) 
  = \sorb{G,\delta}(a,c)\tensor \sorb{K,\epsilon}(b,d).
\]
Let 
\[
 p\colon \sorb{G,\delta}\tensor\sorb{K,\epsilon} \to h(G\times K)\Spec{}{}
\]
be the restriction of the external smash product, i.e., the additive functor taking 
\begin{align*}
 ((G/H,\delta), (K/L,\epsilon)) 
  &\mapsto  (G_+\smsh_H S^{-\delta(G/H)})\smsh (K_+\smsh_L S^{-\epsilon(K/L)}) \\
  &\iso (G\times K)_+\smsh_{H\times L} S^{-(\delta\times\epsilon)((G\times K)/(H\times L))}.
\end{align*}
Similarly, pairs of maps are taken to their smash products.

Let $i\colon \sorb{G\times K,\zeta}\to h(G\times K)\Spec{}{}$ denote the inclusion of the subcategory.

\item
If $\Mackey T$ is a contravariant $\delta$-$G$-Mackey functor and $\Mackey U$ is a
contravariant $\epsilon$-$K$-Mackey functor, let $\Mackey T\tensor \Mackey U$ denote the 
$(\sorb{G,\delta}\tensor\sorb{K,\epsilon})$-module defined by
 \[
  (\Mackey T\tensor \Mackey U)(a,b) = \Mackey T(a)\tensor \Mackey U(b).
 \]
Define $\MackeyOp S\tensor\MackeyOp V$ similarly for covariant Mackey functors.

\item
Let $\Mackey T \exboxprod \Mackey U= \Mackey T \exboxprod_\zeta \Mackey U$ 
(the {\em external box product}) be
the contravariant $\zeta$-$(G\times K)$-Mackey functor defined by
 \[
 \Mackey T \exboxprod \Mackey U= \Mackey T \exboxprod_\zeta \Mackey U 
   = i^*p_!(\Mackey T\tensor \Mackey U),
 \]
where $i^*$ and $p_!$ are defined in Definition~\ref{def:indres}.
Define the external box product of covariant Mackey functors similarly.

\item
Let $\Delta\leq G\times G$ denote the diagonal subgroup.
If $K=G$, so $\Mackey T$ and $\Mackey U$ are both $G$-Mackey functors,
and $\Delta\in\F(\zeta)$, define the
{\em (internal) box product} to be
 \[
  \Mackey T \boxprod \Mackey U = \Mackey T \boxprod_\zeta \Mackey U
   = (\Mackey T\exboxprod \Mackey U) | \Delta.
 \]
Define the box product of covariant Mackey functors similarly.
\end{enumerate}
 \end{definition}

\begin{remark}
In order to define the internal box product we required $\Delta\in\F(\zeta)$,
hence $\Delta\in\overline{\F(\delta)\times\F(\epsilon)}$. However, the smallest
product subgroup containing $\Delta$ is all of $G\times G$, so the
internal box product can be defined only if $\delta$ and $\epsilon$ are both complete.
Hence, this is the case of most interest to us.
We introduced more general dimension functions largely so that we
could properly handle $\delta\times\epsilon$.
\end{remark}

\begin{proposition}\label{prop:burnsideboxproduct}
Let $\delta$ be a dimension function for $G$, 
let $\epsilon$ be a dimension function for $K$,
and let $\zeta \dimpred \delta\times\epsilon$. Then
\[
 \Mackey A_{G/H,\delta} \exboxprod_\zeta \Mackey A_{K/L,\epsilon} 
   \iso h(G\times K)\Spec{}{}(-,(G\times K)_+\smsh_{H\times L} S^{-(\delta\times\epsilon)((G\times K)/(H\times L))})
\]
and
\[
 \MackeyOp A^{G/H,\delta} \exboxprod_\zeta \MackeyOp A^{K/L,\epsilon} 
   \iso h(G\times K)\Spec{}{}((G\times K)_+\smsh_{H\times L} S^{-(\delta\times\epsilon)((G\times K)/(H\times L))},-).
\]
\end{proposition}

\begin{proof}
It's clear from the definitions that
\[
 \Mackey A_{G/H,\delta} \tensor \Mackey A_{K/L,\epsilon}
  \iso (\sorb{G,\delta}\tensor\sorb{K,\epsilon})(-, ((G/H,\delta),(K/L,\epsilon))).
\]
Applying $p_!$ and using Proposition~\ref{prop:indresadjunction} gives the result.
The proof for covariant functors is similar.
\end{proof}

The following proposition shows that $\Mackey A_{G/G} = \Mackey A_{G/G,0}$ 
acts as a unit for the internal box product.
If $\delta$ is a complete dimension function for $G$, 
recall from Definition~\ref{def:diagonalFamily}
that we have the familial dimension function $\delta_\Delta$ on $G\times G$ with the properties that
$\delta_\Delta \dimpred 0\times\delta$,
$\delta_\Delta(G\times G/\Delta) = 0$, and $\delta_\Delta|\Delta = \delta$.

\begin{proposition}\label{prop:burnsideIdentity}
Let $\delta$ be a complete dimension function for $G$ and let $\Mackey T$ be
a contravariant $\delta$-$G$-Mackey functor. Then
\[
 \Mackey A_{G/G} \boxprod_{\delta_\Delta} \Mackey T \iso \Mackey T.
\]
\end{proposition}

\begin{proof}
By writing 
$\Mackey T \iso \int\nolimits^{(G/H,\delta)} \Mackey T(G/H,\delta) \tensor \Mackey A_{G/H,\delta}$,
we see that it suffices to prove the result for $\Mackey T = \Mackey A_{G/H,\delta}$.
From the preceding proposition we know that 
\[
 \Mackey A_{G/G}\exboxprod_{\delta_\Delta}\Mackey A_{G/H,\delta} 
   \iso h(G\times G)\Spec{}{}(-,(G\times G)_+\smsh_{G\times H}S^{-\delta(G/H)}),
\]
with $G\times H$ acting on $\delta(G/H)$ via the projection $G\times H\to H$.
Let $i\colon\Delta\to G\times G$ be the inclusion of the diagonal. 
Then, using that the restriction of the spectrum
$(G\times G)_+\smsh_{G\times H} S^{-\delta(G/H)}$ 
to the diagonal is
$G_+\smsh_H S^{-\delta(G/H)}$, we have
\begin{align*}
(&\Mackey A_{G/G}\boxprod_{\delta_\Delta}\Mackey A_{G/H,\delta})(G/K,\delta) \\
 &= [(h(G\times G)\Spec{}{}(-,(G\times G)_+\smsh_{G\times H}S^{-\delta(G/H)}))|\Delta](G/K,\delta) \\
 &= h(G\times G)\Spec{}{}((G\times G)_+\smsh_\Delta (G_+\smsh_K S^{-\delta(G/K)}),
              (G\times G)_+\smsh_{G\times H} S^{-\delta(G/H)}) \\
 &\iso \sorb{G,\delta}(G_+\smsh_K S^{-\delta(G/K)}, G_+\smsh_H S^{-\delta(G/H)}) \\
 &= \Mackey A_{G/H,\delta}(G/K,\delta).
\end{align*}
Thus $\Mackey A_{G/G}\boxprod_{\delta_\Delta}\Mackey A_{G/H,\delta} \iso \Mackey A_{G/H,\delta}$, which implies
the result in general.
\end{proof}

The following result identifies the chain complex of a product of CW complexes.

\begin{proposition}\label{prop:productChains}
Let $\delta$ be a dimension function for $G$ and
let $\epsilon$ be a dimension function for $K$.
If $X$ is a based $\delta$-$G$-CW($\alpha$) complex and
$Y$ is a based $\epsilon$-$K$-CW($\beta$) complex, then
$X\smsh Y$ is a based $(\delta\times\epsilon)$-$(G\times K)$-CW($\alpha+\beta$) complex with
\[
 \Mackey C^{G,\delta}_{\alpha+*}(X,*) \exboxprod_{\delta\times\epsilon} \Mackey C^{K,\epsilon}_{\beta+*}(Y,*)
  \iso \Mackey C^{G\times K,\delta\times\epsilon}_{\alpha+\beta + *} (X\smsh Y,*).
\]
Moreover, this isomorphism respects suspension in each of $X$ and $Y$.
\end{proposition}

\begin{proof}
We give $X\smsh Y$ the filtration
\[
 (X\smsh Y)^{\alpha+\beta+n} = \Union_{i + j = n} X^{\alpha+i}\smsh Y^{\alpha+j}.
\]
With this filtration, $X\smsh Y$ is a based 
$(\delta\times\epsilon)$-$(G\times K)$-CW($\alpha+\beta$) complex, with
\[
 (X\smsh Y)^{\alpha+\beta+n}/(X\smsh Y)^{\alpha+\beta+n-1}
  = \Wedge_{i+j=n} X^{\alpha+i}/X^{\alpha+i-1}\smsh Y^{\beta+j}/Y^{\beta+j-1}.
\]
We now define
\[
 \kappa\colon 
 \Mackey C^{G,\delta}_{\alpha+*}(X,*) \exboxprod_{\delta\times\epsilon} \Mackey C^{K,\epsilon}_{\beta+*}(Y,*)
  \to \Mackey C^{G\times K,\delta\times\epsilon}_{\alpha+\beta + *} (X\smsh Y,*)
\]
by $\kappa(x\tensor y) = (-1)^{pq} x\smsh y$, where
\[
 x\in \Mackey C^{G,\delta}_{\alpha+p}(X,*)(G/H)
  = [G_+\smsh_H S^{-\delta(G/H)+\alpha+p}, \susp_G^\infty X^{\alpha+p}/X^{\alpha+p-1}]_G,
\]
\[
 y\in \Mackey C^{K,\epsilon}_{\beta+q}(Y,*)(K/L)
  = [K_+\smsh_L S^{-\epsilon(K/L)+\beta+q}, \susp_K^\infty Y^{\beta+q}/Y^{\beta+q-1}]_K,
\]
and $x\smsh y\in \Mackey C^{G\times K,\delta\times\epsilon}_{\alpha+\beta + *} (X\smsh Y,*)(G\times K/H\times L)$ is
their smash product.
With the sign convention
$d(x\tensor y) = dx\tensor y + (-1)^{p}x\tensor dy$, we see that $\kappa$ is an
isomorphism of chain complexes, using Proposition~\ref{prop:burnsideboxproduct}
and the calculation $d(x\smsh y) = (-1)^q dx\smsh y + x\smsh dy$ as in
\cite[\S13.4]{May:concise}.
That the isomorphism respects suspension is straightforward.
\end{proof}

Because of the limitations of $\delta\times\epsilon$, this result by itself is not
as useful as we would like, in particular for discussing general spaces or
internal products. Therefore, we need the following results.
First, we note that level-wise smash product induces a map
\[
 \smsh\colon G\PreSpec{\V}{}\tensor K\PreSpec{\W}{}
  \to (G\times K)\PreSpec{(\V\dirsum\W)}{}
\]
where $\V\dirsum\W = \{ V_i\dirsum W_i \}$.
This pairing extends to semistable maps as well.

\begin{proposition}
Let $\delta$ be a familial dimension function for $G$ and let
$\epsilon$ be a familial dimension function for $K$.
Let $\Gamma^\delta\susp_G^\infty X\to \susp_G^\infty X$ and
$\Gamma^\epsilon\susp_K^\infty Y\to \susp_K^\infty Y$
be, respectively, $\delta$-$G$-CW($\alpha$) and
$\epsilon$-$K$-CW($\beta$) approximations.
Then, there exists a nonnegative integer $N$ such that
\[
 (\Gamma^\delta\susp_G^\infty X\smsh \Gamma^\epsilon\susp_K^\infty Y)(V_i\dirsum W_i)
  \to \susp_G^{V_i} X \smsh \susp_K^{W_i} Y
\]
is an $\overline{\F(\delta)\times\F(\epsilon)}$-equivalence for all $i\geq N$.
\end{proposition}

\begin{proof}
This is an immediate corollary of Proposition~\ref{prop:productApproximation}
and the fact that, for large enough $i$,
$(\Gamma^\delta\susp_G^\infty X)(V_i) \to \susp_G^{V_i} X$ is
an $\F(\delta)$-approximation, and similarly for $Y$.
\end{proof}

\begin{corollary}
Let $\delta$ be a familial dimension function for $G$, let
$\epsilon$ be a familial dimension function for $K$,
and let $\zeta$ be a familial dimension function for $G\times K$ with $\zeta\dimpred\delta\times\epsilon$.
Let $\Gamma^\delta\susp_G^\infty X\to \susp_G^\infty X$ and
$\Gamma^\epsilon\susp_K^\infty Y\to \susp_K^\infty Y$
be, respectively, $\delta$-$G$-CW($\alpha$) and
$\epsilon$-$K$-CW($\beta$) approximations.
and let $\Gamma^\zeta\susp_{G\times K}^\infty(X\smsh Y) \to \susp_{G\times K}^\infty(X\smsh Y)$
be a $\zeta$-$(G\times K)$-CW($\alpha+\beta$) approximation.
Then there exists a semistable cellular map
\[
 \mu\colon \Gamma^\zeta\susp_{G\times K}^\infty(X\smsh Y) \to
  \Gamma^\delta\susp_G^\infty X \smsh \Gamma^\epsilon\susp_K^\infty Y
\]
over $\susp_{G\times K}^\infty(X\smsh Y)$, unique up to semistable cellular homotopy.
\end{corollary}

\begin{proof}
This follows from the preceding proposition and
Corollary~\ref{cor:productMapApproximation}.
\end{proof}

\begin{definition}\label{def:chainProduct}
Let $\delta$ be a familial dimension function for $G$, let
$\epsilon$ be a familial dimension function for $K$,
and let $\zeta$ be a familial dimension function for $G\times K$ with $\zeta\dimpred\delta\times\epsilon$.
If $X$ is a based $G$-space and $Y$ is a based $K$-space, let
\[
 \mu_*\colon \Mackey C^{G\times K,\zeta}_{\alpha+\beta+*}(X\smsh Y)
  \to \Mackey C^{G,\delta}_{\alpha+*}(X)\exboxprod_{\zeta} \Mackey C^{K,\epsilon}_{\beta+*}(Y)
\]
be the chain map induced by the map $\mu$ from the preceding corollary
(using Proposition~\ref{prop:productChains} to identify the chain complex on the right).
It is well-defined up to chain homotopy.
\end{definition}

We can now define pairings in cohomology.
Let $\delta$, $\epsilon$, and $\zeta$ be familial with $\zeta\dimpred \delta\times\epsilon$, as above,
let $X$ be a based $G$-space, let $Y$ be a based $K$-space,
let $\Mackey T$ be a $\delta$-$G$-Mackey functor, and
let $\Mackey U$ be an $\epsilon$-$K$-Mackey functor.
The external box product $\exboxprod = \exboxprod_\zeta$ and the
map $\mu_*$ induce a natural chain map
\begin{align*}
 \Hom_{\sorb{G,\delta}}(\Mackey C^{G,\delta}_{\alpha+*}(X), \Mackey T) \tensor {}&
       \Hom_{\sorb{K,\epsilon}}(\Mackey C^{K,\epsilon}_{\beta+*}(Y), \Mackey U) \\
 &\to
   \Hom_{\sorb {G\times K,\zeta}}(
    \Mackey C^{G,\delta}_{\alpha+*}(X)\exboxprod \Mackey C^{K,\epsilon}_{\beta+*}(Y), 
    \Mackey T \exboxprod \Mackey U) \\
 &\to \Hom_{\sorb {G\times K,\zeta}}(
    \Mackey C^{G\times K,\zeta}_{\alpha+\beta+*}(X\smsh Y), 
    \Mackey T \exboxprod \Mackey U).
\end{align*}
This induces the (external) cup product
 \[
 -\cup - \colon 
 \tilde H^{\alpha}_{G,\delta}(X;\Mackey T) \tensor
 \tilde H^{\beta}_{K,\epsilon}(Y;\Mackey U) 
  \to \tilde H^{\alpha+\beta}_{G\times K,\zeta}(X\smsh Y;\Mackey T\exboxprod \Mackey U).
 \]
When $G = K$ and $\F(\zeta)$ contains the diagonal subgroup $\Delta\leq G\times G$, we can follow the external cup product with the restriction to $\Delta$. This gives the internal cup product
\[
 -\cup - \colon 
 \tilde H^{\alpha}_{G,\delta}(X;\Mackey T) \tensor
 \tilde H^{\beta}_{G,\epsilon}(Y;\Mackey U) 
 \to \tilde H^{\alpha+\beta-\zeta(G\times G/\Delta)}_{G,\zeta|\Delta}(X\smsh Y;\Mackey T\boxprod \Mackey U).
\]
Of course, when $X = Y$ we can apply restriction along the diagonal $X\to X\smsh X$ to
get
\[
 -\cup - \colon 
 \tilde H^{\alpha}_{G,\delta}(X;\Mackey T) \tensor
 \tilde H^{\beta}_{G,\epsilon}(X;\Mackey U) 
 \to \tilde H^{\alpha+\beta-\zeta(G\times G/\Delta)}_{G,\zeta|\Delta}(X;\Mackey T\boxprod \Mackey U).
\]
A useful special case is $\delta_\Delta \dimpred 0\times \delta$ with $\delta$ complete, which gives us pairings
\[
 -\cup - \colon 
 \tilde H^{\alpha}_{G}(X;\Mackey T) \tensor
 \tilde H^{\beta}_{G,\delta}(Y;\Mackey U) 
 \to \tilde H^{\alpha+\beta}_{G,\delta}(X\smsh Y;\Mackey T\boxprod \Mackey U)
\]
and
\[
 -\cup - \colon 
 \tilde H^{\alpha}_{G}(X;\Mackey T) \tensor
 \tilde H^{\beta}_{G,\delta}(X;\Mackey U) 
 \to \tilde H^{\alpha+\beta}_{G,\delta}(X;\Mackey T\boxprod \Mackey U).
\]

Recall that a contravariant $G$-Mackey functor (i.e., 0-$G$-Mackey functor) $\Mackey T$ is a {\em ring} 
(also called a {\em Green functor}) if there is an associative multiplication 
$\Mackey T\boxprod \Mackey T\to \Mackey T$. 
This then makes $\tilde H_G^*(X;\Mackey T)$ a ring.
For example, consider $\Mackey A_{G/G} = \Mackey A_{G/G,0}$, which is a ring
by Proposition~\ref{prop:burnsideIdentity}.
By that same proposition,
$\Mackey A_{G/G} \boxprod_{\delta_\Delta} \Mackey T \iso \Mackey T$ for any contravariant
$\delta$-$G$-Mackey functor $\Mackey T$. 
This makes $\Mackey T$ a {\em module} over
$\Mackey A_{G/G}$, in the sense that there is an associative pairing
$\Mackey A_{G/G}\boxprod_{\delta_\Delta} \Mackey T\to \Mackey T$ (namely, the isomorphism).
This makes every cellular cohomology theory with complete $\delta$
a module over ordinary cohomology
with coefficients in $\Mackey A_{G/G}$, using the cup product
\[
 -\cup - \colon 
 \tilde H^{\alpha}_{G}(X;\Mackey A_{G/G}) \tensor
 \tilde H^{\beta}_{G,\delta}(X;\Mackey T) 
 \to \tilde H^{\alpha+\beta}_{G,\delta}(X;\Mackey T).
\]

Obviously, there are other variations possible.
The following theorem records the main properties of the external cup product, 
from which similar properties
of the other products follow by naturality.

\begin{theorem}\label{thm:vcupproduct}
Let $\delta$ be a familial dimension function for $G$, 
let $\epsilon$ be a familial dimension function for $K$,
and let $\zeta$ be a familial dimension function for $G\times K$
with $\zeta\dimpred \delta\times\epsilon$.
The external cup product
 \[
 -\cup - \colon 
 \tilde H^{\alpha}_{G,\delta}(X;\Mackey T) \tensor
 \tilde H^{\beta}_{K,\epsilon}(Y;\Mackey U) 
  \to \tilde H^{\alpha+\beta}_{G\times K,\zeta}(X\smsh Y;\Mackey T\exboxprod \Mackey U)
 \]
generalizes the nonequivariant cup product and
satisfies the following.
 \begin{enumerate}
 \item It is natural: $f^*(x)\cup g^*(y) = (f\smsh g)^*(x\cup y)$.

 \item It respects suspension: For any representation $V$ of $G$, the following diagram commutes:
\[
 \xymatrix{
 \tilde H^{\alpha}_{G,\delta}(X;\Mackey T) \tensor \tilde H^{\beta}_{K,\epsilon}(Y;\Mackey U)
  \ar[r]^\cup \ar[d]_{\sigma^V\tensor\id}
  & \tilde H^{\alpha+\beta}_{G\times K,\zeta}(X\smsh Y;\Mackey T\exboxprod \Mackey U)
     \ar[dd]^{\sigma^V} \\
 \tilde H^{\alpha+V}_{G,\delta}(\susp^V X;\Mackey T) \tensor \tilde H^{\beta}_{K,\epsilon}(Y;\Mackey U)
  \ar[d]_\cup \\
 \tilde H^{\alpha+V+\beta}_{G\times K,\zeta}(\susp^V X\smsh Y;\Mackey T\exboxprod \Mackey U)
     \ar[r]_-\iso
  & \tilde H^{\alpha+\beta+V}_{G\times K,\zeta}(\susp^V(X\smsh Y);\Mackey T\exboxprod \Mackey U)
 }
\]
The horizontal isomorphism at the bottom of the diagram comes from the identification
$\alpha+V+\beta \iso \alpha+\beta+V$.
The similar diagram for suspension of $Y$ also commutes.

 \item It is associative: $(x\cup y)\cup z = x\cup (y\cup z)$ when we identify
gradings using the obvious identification $(\alpha+\beta)+\gamma \iso \alpha + (\beta+\gamma)$.

 \item\label{item:commutativity} It is commutative: If $x\in \tilde H^\alpha_{G,\delta}(X;\Mackey T)$ and
$y\in \tilde H^\beta_{K,\epsilon}(Y;\Mackey U)$ then $x\cup y = \iota(y\cup x)$ where $\iota$ is the
evident isomorphism between
$\tilde H^{\alpha+\beta}_{G\times K,\zeta}(X\smsh Y;\Mackey T\exboxprod \Mackey U)$
and
$\tilde H^{\beta+\alpha}_{K\times G,\tilde\zeta}(Y\smsh X;\Mackey U\exboxprod \Mackey T)$;
$\tilde\zeta$ is the dimension function on $K\times G$ induced by $\zeta$
and $\iota$ uses the isomorphism of $\alpha+\beta$ and $\beta+\alpha$ that
switches the direct summands.

 \item It is unital: The map
\[
 \tilde H^0_G(S^0;\Mackey A_{G/G}) \tensor \tilde H^\alpha_{G,\delta}(X;\Mackey T)
   \to \tilde H^\alpha_{G,\delta}(X;\Mackey T)
\]
takes $1\tensor x \mapsto x$, where $1 \in \tilde H^0_G(S^0;\Mackey A_{G/G})\iso A(G)$
is the unit.

\item\label{item:vcupWirthmuller}
It respects the Wirthm\"uller isomorphism:
Suppose that $J\leq G$ and $L\leq K$, that
$J\in\F(\delta)$, $L\in\F(\epsilon)$, and $J\times L\in \F(\zeta)$, and
that $\zeta(G\times K/J\times L) = \delta(G/K)\oplus \epsilon(K/L)$.
Then the following diagram commutes:
\[
 \def\objectstyle{\scriptstyle}
 \def\labelstyle{\scriptstyle}
 \xymatrix{
  \tilde H_{G,\delta}^\alpha(G_+\smsh_J X; \Mackey T)
    \tensor \tilde H_{K,\epsilon}^\beta(K_+\smsh_L Y; \Mackey U) \ar[dd]_{-\cup-} \ar[r]^-\iso
    & \tilde H_{J,\delta}^{\alpha-\delta(G/J)}(X; \Mackey T|J)
          \tensor \tilde H_{L,\epsilon}^{\beta-\epsilon(K/L)}(Y;\Mackey U|L) \ar[d]^{-\cup-} \\
  & \tilde H_{J\times L,\zeta}^{\alpha+\beta-\zeta(G\times K/J\times L)}(X\smsh Y; (\Mackey T|J)\exboxprod(\Mackey U|L)) \ar[d] \\
  \tilde H_{G\times K,\zeta}^{\alpha+\beta}((G\times K)_+\smsh_{J\times L}(X\smsh  Y);
         \Mackey T\exboxprod \Mackey U) \ar[r]_-{\iso}
    & \tilde H_{J\times L,\zeta}^{\alpha+\beta-\zeta(G\times K/J\times L)}(X\smsh Y; (\Mackey T\exboxprod\Mackey U)|(J\times L))
 }
\]

 \item\label{item:vcupsubgroup}
 It respects restriction to subgroups: Under the same assumptions as the previous point,
\[
 (x|J)\cup (y|L) = (x\cup y)|(J\times L).
\]
(But see Remark~\ref{rem:prodCoeffs}.)

 \item\label{item:vcupfixedsets}
 It respects restriction to fixed sets:
\[
 x^J \cup y^L = (x\cup y)^{J\times L}.
\]
(But, again, see Remark~\ref{rem:prodCoeffs}.)

 \end{enumerate}
 \end{theorem}

The proofs are all standard except for the last three points.

\begin{proof}[Proof of Parts (\ref{item:vcupWirthmuller}) and (\ref{item:vcupsubgroup}) of Theorem \ref{thm:vcupproduct}]
(Refer to  Section~\ref{sec:vsubgroups}
for the algebra involved in the construction of the
restriction to subgroups.) 
Note that the assumptions given are necessary for the statement to even make sense;
see Proposition~\ref{prop:mackeycupproduct} for a related result.
From the commutativity of the diagram
 \[
 \xymatrix{
  {\sorb{J,\delta} \tensor \sorb{L,\epsilon}} \ar[r]^-{p} \ar[d]_{i_J^G\tensor i_L^K}
   & h(J\times L)\Spec{}{} \ar[d]^{i_{J\times L}^{G\times K}}
   & {\sorb{J\times L,\zeta}} \ar[l]_-{j} \ar[d]^{i_{J\times L}^{G\times K}} \\
  {\sorb{G,\delta} \tensor \sorb{K,\epsilon}} \ar[r]_-{p}
   & h(G\times K)\Spec{}{}
   & {\sorb{G\times K,\zeta}} \ar[l]^-{j}
 }
 \]
(where $j$ denotes the inclusions) it follows that there is a natural isomorphism
 \[
 (G\times_J \Mackey C)\exboxprod (K\times_L \Mackey D)
 \iso (G\times K)\times_{J\times L} (\Mackey C\exboxprod \Mackey D).
 \]
The diagram and the $p_!$-$p^*$ adjunction lead to a natural homomorphism
 \[
 (\Mackey T|J)\exboxprod (\Mackey U|L) \to (\Mackey T\exboxprod \Mackey U)|(J\times L)
 \]
(not necessarily an isomorphism). Together these give the following
diagram.
 \[
 \def\objectstyle{\scriptstyle}
 \def\labelstyle{\scriptstyle}
 \xymatrix{
  {\Hom_{\sorb{G,\delta}}(G\times_J\Mackey C, \Mackey T)\tensor
    \Hom_{\sorb{K,\epsilon}}(K\times_L\Mackey D, \Mackey U)}
    \ar[r]^-{\iso} \ar[d]
    & {\Hom_{\sorb{J,\delta}}(\Mackey C, \Mackey T|J)\tensor
        \Hom_{\sorb{L,\epsilon}}(\Mackey D, \Mackey U|L)} \ar[d] \\
  {\Hom_{\sorb{G\times K,\zeta}}((G\times_J\Mackey C) \exboxprod (K\times_L\Mackey D),
      \Mackey T \exboxprod \Mackey U)} \ar[d]_{\iso}
    & {\Hom_{\sorb {J\times L,\zeta}}(\Mackey C\exboxprod \Mackey D,
      (\Mackey T|J) \exboxprod (\Mackey U|L))} \ar[d] \\
  {\Hom_{\sorb {G\times K,\zeta}}((G\times K)\times_{J\times L}
       (\Mackey C \exboxprod \Mackey D), \Mackey T \exboxprod \Mackey U)}
    \ar[r]^-{\iso}
    & {\Hom_{\sorb {J\times L,\zeta}}(\Mackey C\exboxprod \Mackey D,
      (\Mackey T \exboxprod \Mackey U)|(J\times L))}
  }
 \]
It is an elementary exercise to show that this diagram
commutes. Letting $\Mackey C = \Mackey C_*^{J,\delta}(X)$ and
$\Mackey D = \Mackey C_*^{L,\epsilon}(Y)$ shows that the cup product respects the
Wirthm\"uller isomorphism on the chain level, so does so also in cohomology.

Since restriction from
$G$ to $J$ is defined as the composite
 \[
 \tilde H_{G,\delta}^\alpha(X;\Mackey T) 
 \to \tilde H_{G,\delta}^\alpha(G_+\smsh_J X;\Mackey T)
 \iso \tilde H_{J,\delta}^{\alpha-\delta(G/J)}(X;\Mackey T|J)
 \]
and the cup product is natural, that the cup product respects restriction to subgroups
follows from the fact that it respects the
Wirthm\"uller isomorphism. 
 \end{proof}

\begin{proof}[Proof of Part (\ref{item:vcupfixedsets}) of Theorem \ref{thm:vcupproduct}]
(Refer to Section~\ref{sec:vsubgroups}
for the algebra involved in the construction of
restriction to fixed points.) 
By part~(\ref{item:vcupsubgroup}) we know that the cup product respects the
restriction to the normalizers $NJ$ and $NL$, so we may as well assume that
$J$ is normal in $G$ and $L$ is normal in $K$.
The naturality of the cup product means that it suffices to show that the cup
product preserves the isomorphism of Theorem~\ref{thm:fixedSetIso}.

From the commutativity of the diagram
\[
 \xymatrix{
  \sorb{G,\delta} \tensor \sorb{K,\epsilon} \ar[r]^-p \ar[d]_{\Phi^J\tensor\Phi^L}
    & h(G\times K)\Spec{}{} \ar[d]^{\Phi^{J\times L}}
    & \sorb{G\times K,\zeta} \ar[l]_-{j} \ar[d]^{\Phi^{J\times L}} \\
  \sorb{G/J,\delta} \tensor \sorb{K/L,\epsilon} \ar[r]_-p
    & h(G\times K/J\times L)\Spec{}{}
    & \sorb{G\times K/J \times L,\zeta} \ar[l]^-{j}
 }
\]
(where $j$ denotes the inclusions) it follows that there is a natural isomorphism
\[
 \Mackey C^J \exboxprod \Mackey D^L \iso (\Mackey C \exboxprod \Mackey D)^{J\times L}
\]
and a natural homomorphism
\[
 \Inf_{G/J}^G\Mackey T \exboxprod \Inf_{K/L}^K\Mackey U
  \to \Inf_{G\times K/J\times L}^{G\times K}(\Mackey T \exboxprod \Mackey U).
\]
The result now follows from the following commutative diagram:
 \[
 \def\objectstyle{\scriptstyle}
 \def\labelstyle{\scriptstyle}
  \xymatrix{
   {\Hom_{\sorb{G,\delta}}(\Mackey C, \Inf_{G/J}^G\Mackey T)\tensor
     \Hom_{\sorb{K,\epsilon}}(\Mackey D, \Inf_{K/L}^K\Mackey U)} \ar[r]^-{\iso} \ar[d]
     & {\Hom_{\sorb {G/J,\delta}}(\Mackey C^J, \Mackey T)\tensor
          \Hom_{\sorb {K/L,\epsilon}}(\Mackey D^L, \Mackey U)} \ar[d] \\
   {\Hom_{\sorb {G\times K,\zeta}}(\Mackey C \exboxprod \Mackey D,
      \Inf_{G/J}^G\Mackey T \exboxprod \Inf_{K/L}^K\Mackey U)} \ar[d]
     & {\Hom_{\sorb {G\times K/J\times L,\zeta}}(\Mackey C^J\exboxprod \Mackey D^L,
         \Mackey T \exboxprod \Mackey U)} \ar[d]^\iso \\
   {\Hom_{\sorb {G\times K,\zeta}}(\Mackey C \exboxprod \Mackey D,
      \Inf_{G\times K/J\times L}^{G\times K}(\Mackey T \exboxprod \Mackey U))} \ar[r]^-{\iso}
     & {\Hom_{\sorb {G\times K/J\times L,\zeta}}((\Mackey C \exboxprod \Mackey D)^{J\times L},
       \Mackey T \exboxprod \Mackey U)}
  }
 \]
 \end{proof}

\begin{remark}\label{rem:prodCoeffs}
Part (\ref{item:vcupsubgroup}) of Theorem~\ref{thm:vcupproduct} is stated
a bit loosely. In fact,
 \[
 (x|J)\cup(y|L) \in \tilde H_{J\times L,\zeta}^*(X\smsh Y; 
			(\Mackey T|J)\exboxprod (\Mackey U|L))
 \]
while
 \[
 (x\cup y)|(J\times L) \in \tilde H_{J\times L,\zeta}^*(X\smsh Y;
  (\Mackey T\exboxprod \Mackey U)|(J\times L)).
 \]
The theorem should say that, when we apply the map induced by the homomorphism
$(\Mackey T|J)\exboxprod (\Mackey U|L) \to 
(\Mackey T\exboxprod \Mackey U)|(J\times L)$,
the element $(x|J)\cup(y|L)$ maps to $(x\cup y)|(J\times L)$.

A similar comment applies to part (\ref{item:vcupfixedsets}) in general,
when we have to first restrict to normalizers.
However, if $J$ is normal in $G$ and $L$ is normal in $K$, we can use the isomorphism
$\Mackey T^J \exboxprod \Mackey U^L \iso (\Mackey T \exboxprod \Mackey U)^{J\times L}$
to identify the two cohomology groups.
\end{remark}

\begin{remark}
Part (\ref{item:commutativity}) of Theorem~\ref{thm:vcupproduct} states that the external
product is commutative, in the sense given there. What about the internal cup product
(in the case $\delta = \epsilon = 0$)
\[
 \tilde H^{\alpha}_{G}(X;\Mackey T) \tensor
 \tilde H^{\beta}_{G}(X;\Mackey T) 
 \to \tilde H^{\alpha+\beta}_{G}(X;\Mackey T)?
\]
Recall that this product is defined when $\Mackey T$ is a ring, meaning that there is
an associative product $\Mackey T\boxprod \Mackey T\to \Mackey T$.
Let us further assume that the product is {\em commutative}, meaning that the following diagram
commutes, where $\gamma$ is the interchange map:
\[
 \xymatrix@C1em{
  \Mackey T \boxprod \Mackey T \ar[rr]^\gamma \ar[rd]
    && \Mackey T \boxprod \Mackey T \ar[ld] \\
  & \Mackey T.
 }
\]
After restricting to diagonals, Theorem~\ref{thm:vcupproduct} tells us that,
if $x\in \tilde H^{\alpha}_{G}(X;\Mackey T)$ and $y \in \tilde H^{\beta}_{G}(X;\Mackey T)$,
then $x\cup y = \iota(y\cup x)$, where $\iota$ is the self map of
$H^{\alpha+\beta}_{G}(X;\Mackey T)$
induced by reversing the summands in $\alpha+\beta$ and then identifying $\beta+\alpha$ with
$\alpha+\beta$ again. To be very precise: Write
\[
 \alpha = \Dirsum_{i} V_i^{m_i} \qquad\text{and}\qquad
 \beta = \Dirsum_{i} V_i^{n_i}
\]
where the $V_i$ are the irreducible representations of $G$ and $m_i, n_i\in \Z$.
Then $\iota$ is given by multiplication by the unit in $A(G)$ determined by the stable self-map
\[
 \Smsh_i S^{(m_i+n_i)V_i} = \Smsh_i S^{m_iV_i}\smsh S^{n_iV_i}
  \xrightarrow{\gamma} \Smsh_i S^{n_iV_i}\smsh S^{m_iV_i}
  = \Smsh_i S^{(m_i+n_i)V_i}.
\]
Call this element $\gamma(\alpha,\beta)\in A(G)$. Explicitly, if $\gamma_i = \gamma(V_i,V_i)$
is the unit given by the interchange map on $S^{V_i}\smsh S^{V_i}$, then
\[
 \gamma(\alpha,\beta) = \prod_i \gamma_i^{m_i n_i}.
\]
We then we have the anti-commutativity rule
\[
 x\cup y = \gamma(\alpha,\beta)y\cup x
\]
for the internal cup product.
\end{remark}

We have the following result for Mackey-valued cohomology.

\begin{proposition}\label{prop:mackeycupproduct}
Let $\delta$ be a familial dimension function for $G$, 
let $\epsilon$ be a familial dimension function for $K$,
and let $\zeta$ be a familial dimension function for $G\times K$
with $\zeta\dimpred \delta\times\epsilon$.
Then the external cup product gives a map of $(G\times K)$-$\zeta$-Mackey functors
\[
 -\cup-\colon \Mackey H_{G,\delta}^\alpha(X;\Mackey T) 
   \exboxprod_\zeta \Mackey H_{K,\epsilon}^\beta(Y;\Mackey U)
   \to \Mackey H_{G\times K,\zeta}^{\alpha+\beta}(X\smsh Y;\Mackey T\exboxprod_\zeta\Mackey U).
\]
When $G=K$ and $\F(\zeta)$ contains the diagonal subgroup $\Delta$, the cup product restricts
to give a map of $G$-$\zeta$-Mackey functors
\[
 -\cup-\colon \Mackey H_{G,\delta}^\alpha(X;\Mackey T) 
   \boxprod_\zeta \Mackey H_{G,\epsilon}^\beta(Y;\Mackey U)
   \to \Mackey H_{G,\zeta}^{\alpha+\beta-\zeta(G\times G/\Delta)}(X\smsh Y;\Mackey T\boxprod_\zeta\Mackey U).
\]
\end{proposition}

\begin{proof}
For $J\in\F(\delta)$, $L\in\F(\epsilon)$, and $M\in\F(\zeta)$, to make the following readable
we introduce the abbreviations
\begin{align*}
 A &= G_+\smsh_J S^{-\delta(G/J)}, \\
 B &= K_+\smsh_L S^{-\delta(K/L)},\\
\intertext{and}
 C &= (G\times K)_+\smsh_M S^{-\zeta(G\times K/M)}.
\end{align*}
We then have the composite
\begin{align*}
 \tilde H_{G,\delta}^\alpha(&X\smsh A;\Mackey T) \tensor
 \tilde H_{K,\epsilon}^\beta(Y\smsh B;\Mackey U)
 \tensor
 (G\times K)\Spec{}{}(C, A\smsh B) \\
 &\xrightarrow{-\cup-}
  \tilde H_{G\times K,\zeta}^{\alpha+\beta}(X\smsh Y\smsh
     A\smsh B; \Mackey T\exboxprod\Mackey U)
      \tensor (G\times K)\Spec{}{}(C, A\smsh B) \\
 &\to \tilde H_{G\times K,\zeta}^{\alpha+\beta}(
       X\smsh Y\smsh C; \Mackey T\exboxprod\Mackey U)
\end{align*}
Naturality of the cup product implies that this passes to the coend 
over $A$ and $B$ to give a map of Mackey functors
\[
 -\cup-\colon \Mackey H_{G,\delta}^\alpha(X;\Mackey T) 
   \exboxprod \Mackey H_{K,\epsilon}^\beta(Y;\Mackey U)
   \to \Mackey H_{G\times K,\zeta}^{\alpha+\beta}(X\smsh Y;\Mackey T\exboxprod\Mackey U).
\]
Restriction to the diagonal subgroup gives the internal version.
\end{proof}

The pairings in the preceding proposition were defined in terms of the pairing
of $G$- and $K$-cohomology.
Note, however, that Theorem~\ref{thm:vcupproduct}(\ref{item:vcupWirthmuller})
says that, if $J\times L\in\F(\zeta)$ and
$\zeta(G\times K/J\times L) = \delta(G/J)\dirsum\epsilon(K/L)$, then the component
\[
 \tilde H_{J,\delta}^\alpha(X;\Mackey T|J) \tensor
 \tilde H_{L,\epsilon}^\beta(Y;\Mackey U|L)
 \to H_{J\times L,\zeta}^{\alpha+\beta}(X\smsh Y;(\Mackey T\exboxprod\Mackey U)|J\times L)
\]
of the Mackey functor pairing,
corresponding to the identity map on $(G\times K)_+\smsh_{J\times L}S^{-\zeta(G\times K/J\times L)}$,
agrees with the cup product pairing $J$- and $L$-cohomology.

Now we look at how the cup product is represented.
Given the $G$-spectrum $H_\delta\Mackey T$ and $K$-spectrum $H_\epsilon\Mackey U$,
we can form the $(G\times K)$-spectrum
$H_\delta\Mackey T \smsh H_\epsilon\Mackey U$.
The external cup product should then be represented by a $(G\times K)$-map
$H_\delta\Mackey T \smsh H_\epsilon\Mackey U \to H_\zeta(\Mackey T\exboxprod \Mackey U)$
that is an isomorphism in $\Mackey\pi_0^{G\times K,\zeta}$. 
An explicit construction
is based on the following calculation.

\begin{proposition}\label{prop:EMsmshhomotopy}
Let $\delta$ be a familial dimension function for $G$, let $\epsilon$ be
a familial dimension function for $K$, and let $\zeta\dimpred \delta\times\epsilon$.
Let $\Mackey T$ be a $\delta$-$G$-Mackey functor and 
let $\Mackey U$ be an $\epsilon$-$K$-Mackey functor.
Then we have
\[
 \Mackey\pi_n^{G\times K,\zeta}(H_\delta\Mackey T \smsh H_\epsilon\Mackey U)
  = \begin{cases}
      0 & \text{if $n < 0$} \\
      \Mackey T \exboxprod_\zeta \Mackey U & \text{if $n = 0$.}
    \end{cases}
\]
\end{proposition}

\begin{proof}
There should be a proof based on a K\"unneth spectral sequence as developed in
\cite{EKMM:Foundations} and \cite{LM:SpectralSequences}. However, the former is the
nonequivariant case and the latter deals only with finite groups.
Developing the spectral sequence in the compact Lie case would take us too far afield
and may not be available in the parametrized case we discuss later, anyway.
So, we give a more elementary argument here.

We can assume that our Eilenberg-Mac\,Lane spectra are constructed
as in Construction~\ref{con:EilenbergMacLane} as, for example
$H_\delta\Mackey T = F(E\F(\delta)_+,P_\delta\Mackey T)$.
Consider the cofibrations
\[
 R_\delta\Mackey T \to F_\delta\Mackey T \to C_\delta\Mackey T
\]
and
\[
 R_\epsilon\Mackey U \to F_\epsilon\Mackey U \to C_\epsilon\Mackey U
\]
used in Construction~\ref{con:EilenbergMacLane}, which,
on applying $\Mackey\pi_0$, give the exact sequences
\[
 \Mackey\pi_0^{G,\delta}R_\delta\Mackey T \to \Mackey\pi_0^{G,\delta}F_\delta\Mackey T
  \to \Mackey T \to 0
\]
and
\[
 \Mackey\pi_0^{K,\epsilon}R_\epsilon\Mackey U \to \Mackey\pi_0^{K,\epsilon}F_\epsilon\Mackey U
  \to \Mackey U \to 0.
\]
If $C$ is the cofiber of
\[
 (R_\delta\Mackey T\smsh F_\epsilon\Mackey U) \vee (F_\delta\Mackey T\smsh R_\epsilon\Mackey U)
  \to F_\delta\Mackey T \smsh F_\epsilon\Mackey U,
\]
we have an exact sequence
\begin{multline*}
 \Mackey\pi_0^{G\times K,\zeta}(R_\delta\Mackey T\smsh F_\epsilon\Mackey U) \dirsum
  \Mackey\pi_0^{G\times K,\zeta}(F_\delta\Mackey T\smsh R_\epsilon\Mackey U) \\
  \to \Mackey\pi_0^{G\times K,\zeta}(F_\delta\Mackey T \smsh F_\epsilon\Mackey U)
  \to \Mackey\pi_0^{G\times K,\zeta}(C) \to 0.
\end{multline*}
Because $F_\delta$, $R_\delta$, $F_\epsilon$, and $R_\epsilon$ are wedges of spheres,
using the calculation
\[
 hG\Spec{}{}(-,(G/J,\delta)\smsh(K/L,\epsilon))
        \iso \Mackey A_{G/J,\delta} \exboxprod_\zeta \Mackey A_{K/L,\delta}
\]
on $\sorb{G\times K,\zeta}$,
the exact sequence above is isomorphic to
\begin{multline*}
 (\Mackey\pi_0^{G,\delta}(R_\delta\Mackey T)\exboxprod \Mackey\pi_0^{K,\epsilon}(F_\epsilon\Mackey U))
  \dirsum
  (\Mackey\pi_0^{G,\delta}(F_\delta\Mackey T)\exboxprod \Mackey\pi_0^{K,\epsilon}(R_\epsilon\Mackey U)) \\
  \to \Mackey\pi_0^{G,\delta}(F_\delta\Mackey T) \exboxprod \Mackey\pi_0^{K,\epsilon}(F_\epsilon\Mackey U)
  \to \Mackey\pi_0^{G\times K,\zeta}(C) \to 0,
\end{multline*}
from which it follows that
$\Mackey\pi_0^{G\times K,\zeta}(C) \iso \Mackey T\exboxprod\Mackey U$.
Further, looking at the connectivities, we have
\[
 \Mackey\pi_0^{G\times K,\zeta}(C_\delta\Mackey T\smsh C_\epsilon\Mackey U)
 \iso \Mackey\pi_0^{G\times K,\zeta}(C) 
  \iso \Mackey T\exboxprod\Mackey U.
\]
Passing to $P$ and then $H$,
we conclude that
\[
 \Mackey\pi_0^{G\times K,\zeta}(H_\delta\Mackey T \smsh H_\epsilon\Mackey U)
  \iso \Mackey T\exboxprod\Mackey U
\]
as claimed. The vanishing of the homotopy groups for negative $n$ follows from
the construction.
\end{proof}

As usual, for example, by killing higher homotopy groups, it follows that there
is a map of $(G\times K)$-spectra
\[
 H_\delta\Mackey T \smsh H_\epsilon\Mackey U \to H_\zeta(\Mackey T\exboxprod_\zeta \Mackey U)
\]
that is an isomorphism in $\Mackey\pi_0^{G\times K,\zeta}$.
That this represents the cup product in cohomology that we constructed on the chain level
follows by considering its effect on the quotients
$X^{\alpha+m}/X^{\alpha+m-1}$ and $Y^{\beta+n}/Y^{\beta+n-1}$
when $X$ and $Y$ are CW complexes.

When $G = K$ we can restrict to the diagonal $\Delta\leq G\times G$. Doing so
gives us a map of $G$-spectra
\[
 H_\delta\Mackey T \smsh H_\epsilon\Mackey U 
   \to \susp^{-\zeta(G\times G/\Delta)}H_\zeta(\Mackey T\boxprod_\zeta \Mackey U)
\]
that is an isomorphism in $\Mackey\pi_0^{G,\zeta}$.

Finally, let's point out the specializations to several interesting choices of $\delta$ and $\epsilon$,
using the cup products internal in $G$.

\begin{remark}
The general cup product gives us the following special cases.
\begin{enumerate}
\item
Taking $\delta = \epsilon = 0$ on $G$ and $\zeta=0$ on $G\times G$, we have the cup product
\[
 \tilde H_G^\alpha(X;\Mackey T)\tensor \tilde H_G^\beta(Y;\Mackey U)
   \to \tilde H_G^{\alpha+\beta}(X\smsh Y; \Mackey T\boxprod \Mackey U).
\]
The Mackey functor-valued version is a pairing
\[
 \Mackey H_G^\alpha(X;\Mackey T)\boxprod \Mackey H_G^\beta(Y;\Mackey U)
   \to \Mackey H_G^{\alpha+\beta}(X\smsh Y;\Mackey T\boxprod\Mackey U).
\]
This product is represented by a $G$-map
\[
 H\Mackey T \smsh H\Mackey U \to H(\Mackey T \boxprod \Mackey U).
\]

\item
Taking $\delta=\epsilon=\Lie$ on $G$ and $\zeta=\Lie$ on $G\times G$, we have the cup product
\[
 \tilde\H_G^\alpha(X;\MackeyOp R)\tensor \tilde\H_G^\beta(Y;\MackeyOp S)
   \to \tilde\H_G^{\alpha+\beta-\Lie(G)}(X\smsh Y; \MackeyOp R\boxprod_\Lie \MackeyOp S).
\]
The Mackey functor-valued version is a pairing
\[
 \MackeyOp \H_G^\alpha(X;\MackeyOp R)\boxprod \MackeyOp\H_G^\beta(Y;\MackeyOp S)
   \to \MackeyOp\H_G^{\alpha+\beta-\Lie(G)}(X\smsh Y; \MackeyOp R\boxprod_\Lie \MackeyOp S).
\]
These Mackey functors can all be thought of as either covariant $G$-Mackey functors or
contravariant $G$-$\Lie$-Mackey functors. Here and in the following cases the notation
indicates their variance as $G$-Mackey functors.
This product is represented by a $G$-map
\[
 H_\Lie\MackeyOp R \smsh H_\Lie\MackeyOp S \to \susp^{-\Lie(G)}H_\Lie(\MackeyOp R \boxprod_\Lie \MackeyOp S).
\]

\item
Taking $\delta=0$, $\epsilon=\Lie$, and $\zeta= \Lie_\Delta$, we have the cup product
\[
 \tilde H_G^\alpha(X;\Mackey T)\tensor \tilde\H_G^\beta(Y;\MackeyOp S)
   \to \tilde\H_G^{\alpha+\beta}(X\smsh Y; \Mackey T \boxprod_{\Lie_\Delta} \MackeyOp S).
\]
In terms of Mackey functor-valued theories we get
\[
 \Mackey H_G^\alpha(X;\Mackey T)\boxprod_{\Lie_\Delta} \MackeyOp\H_G^\beta(Y;\MackeyOp S)
   \to \MackeyOp\H_G^{\alpha+\beta}(X\smsh Y; \Mackey T \boxprod_{\Lie_\Delta} \MackeyOp S).
\]
This product is represented by a $G$-map
\[
 H\Mackey T \smsh H_\Lie\MackeyOp S \to H_\Lie(\Mackey T \boxprod_{\Lie_\Delta} \MackeyOp S).
\]
Note that $\Lie_\Delta$ restricts to $\Lie$ on the diagonal, so $\Mackey T \boxprod_{\Lie_\Delta} \MackeyOp S$ is a contravariant
$\Lie$-Mackey functor.

\item
Taking $\delta = 0$, $\epsilon=\Lie$, and $\zeta= \Lie - \Lie_\Delta$, we have the cup product
\[
 \tilde H_G^\alpha(X;\Mackey T)\tensor \tilde\H_G^\beta(Y;\MackeyOp S)
   \to \tilde H_G^{\alpha+\beta-\Lie(G)}(X\smsh Y; \Mackey T \boxprod_{\Lie - \Lie_\Delta} \MackeyOp S).
\]
In terms of Mackey functor-valued theories we get
\[
 \Mackey H_G^\alpha(X;\Mackey T)\boxprod_{\Lie - \Lie_\Delta} \MackeyOp\H_G^\beta(Y;\MackeyOp S)
   \to \Mackey H_G^{\alpha+\beta-\Lie(G)}(X\smsh Y; \Mackey T \boxprod_{\Lie - \Lie_\Delta} \MackeyOp S).
\]
This product is represented by a $G$-map
\[
 H\Mackey T \smsh H_\Lie\MackeyOp S \to \susp^{-\Lie(G)}H(\Mackey T \boxprod_{\Lie - \Lie_\Delta} \MackeyOp S).
\]
Note that $\Mackey T \boxprod_{\Lie - \Lie_\Delta} \MackeyOp S$ is an
ordinary contravariant Mackey functor because $\Lie - \Lie_\Delta$ restricts to 0 on the diagonal.
\end{enumerate}
\end{remark}

\subsection{Slant products, evaluations, and cap products}\label{sec:capproducts}

We now construct evaluation maps and cap products,
but we start with the slant product that underlies both.
As May points out in \cite{May:additivity}, the literature does not agree
on how the slant product should be formulated.
Unfortunately, if we were to follow May's suggestion of using the earliest definition
that appeared, the resulting cap product would not make homology be
a module over cohomology. Given that, we use what we think is the most
useful formulation for us.

A little bit of algebra first.

\begin{definition}
Let $\delta$ be a dimension function for $G$, let $\epsilon$ be a
dimension function for $K$, and let $\zeta$ be a dimension function for $G\times K$
with $\zeta \dimpred \delta\times\epsilon$.
Suppose that $\Mackey T$ is a contravariant $\epsilon$-$K$-Mackey functor and
$\MackeyOp U$ is a covariant $\zeta$-$(G\times K)$-Mackey functor.
Then we define a covariant $\delta$-$G$-Mackey functor $\Mackey T\mixprod \MackeyOp U$ by
\[
 (\Mackey T\mixprod \MackeyOp U)(G/J,\delta) 
  = \Mackey T \tensor_{\sorb{K,\epsilon}} (p^*i_!\MackeyOp U)((G/J,\delta)\tensor -)
\]
where $p\colon \sorb{G,\delta}\tensor\sorb{K,\epsilon}\to h(G\times K)\Spec{}{}$
and $i\colon\sorb{G\times K,\zeta} \to h(G\times K)\Spec{}{}$
are the functors given in Definition~\ref{def:mackeyproducts}.
\end{definition}

\begin{example}\label{ex:mixprod}
As an example and as a calculation we'll need later, we show that
\[
 \Mackey A_{K/L,\epsilon}\mixprod \MackeyOp A^{G\times K/M,\zeta}
  \iso h(G\times K)\Spec{}{}((G\times K/M,\zeta), -\smsh (K/L,\epsilon)).
\]
For, if $(G/J,\delta)$ is an object in $\sorb{G,\delta}$, we have
\begin{align*}
 (\Mackey A_{K/L,\epsilon}\mixprod \MackeyOp A^{G\times K/M,\zeta})&(G/J,\delta) \\
  &= \Mackey A_{K/L,\epsilon}\tensor_{\sorb{K,\epsilon}} 
           (p^*i_!\MackeyOp A^{G\times K/M,\zeta})((G/J,\delta)\tensor-) \\
  &= \Mackey A_{K/L,\epsilon}\tensor_{\sorb{K,\epsilon}} 
           h(G\times K)\Spec{}{}((G\times K/M,\zeta), (G/J,\delta) \smsh -) \\
  &\iso h(G\times K)\Spec{}{}((G\times K/M,\zeta), (G/J,\delta) \smsh (K/L,\epsilon)).
\end{align*}
\end{example}

For $X$ a based $G$-space and $Y$ a based $K$-space, the external slant product will be a map
\[
 -\slant - \colon 
  \tilde H_{K,\epsilon}^\beta(Y;\Mackey T) \tensor
    \tilde H^{G\times K,\zeta}_{\alpha+\beta}(X\smsh Y; \MackeyOp U)
  \to \tilde H^{G,\delta}_\alpha(X; \Mackey T \mixprod \MackeyOp U).
\]
On the chain level, we take the following map:
\begin{align*}
 \Hom_{\sorb{K,\epsilon}}&(\Mackey C^{K,\epsilon}_{\beta+*}(Y),\Mackey T) \tensor
  (\Mackey C^{G\times K,\zeta}_{\alpha+\beta+*}(X\smsh Y)\tensor_{\sorb{G\times K,\zeta}} \MackeyOp U) \\
  &\to \Hom_{\sorb{K,\epsilon}}(\Mackey C^{K,\epsilon}_{\beta+*}(Y),\Mackey T) \tensor
    ((\Mackey C^{G,\delta}_{\alpha+*}(X)\exboxprod_\zeta \Mackey C^{K,\epsilon}_{\beta+*}(Y))
    \tensor_{\sorb{G\times K,\zeta}} \MackeyOp U) \\
  &\xrightarrow{\nu} (\Mackey C^{G,\delta}_{\alpha+*}(X)\exboxprod_\zeta\Mackey T )
    \tensor_{\sorb{G\times K,\zeta}} \MackeyOp U \\
  &\iso (\Mackey C^{G,\delta}_{\alpha+*}(X)\tensor\Mackey T)
    \tensor_{\sorb{G,\delta}\tensor\sorb{K,\epsilon}} p^*i_!\MackeyOp U \\
  &\iso \Mackey C^{G,\delta}_{\alpha+*}(X) \tensor_{\sorb{G,\delta}}
    (\Mackey T \tensor_{\sorb{K,\epsilon}} p^*i_!\MackeyOp U) \\
  &= \Mackey C^{G,\delta}_{\alpha+*}(X) \tensor_{\sorb{G,\delta}}
    (\Mackey T\mixprod \MackeyOp U).
\end{align*}
The first arrow is induced by the map $\mu_*$ from Definition~\ref{def:chainProduct}.
The map $\nu$ is evaluation, with the sign
\[
 \nu(y\tensor a\tensor b\tensor u) = (-1)^{pq}a\tensor y(b)\tensor u
\]
if $y\in \Hom_{\sorb{K,\epsilon}}(\Mackey C^{K,\epsilon}_{\beta+p}(Y),\Mackey T)$ and
$a\in \Mackey C^{G,\delta}_{\alpha+q}(X)$.
Our various sign conventions imply that the composite above is a chain map.
Taking homology defines our slant product.
The following properties follow easily from the definition,
except for the last, which follows from examining chains and arguing much as
in the proof of Part~(\ref{item:vcupWirthmuller}) of Theorem~\ref{thm:vcupproduct}.

\begin{theorem}\label{thm:slantproduct}
Let $\delta$ be a familial dimension function for $G$, let $\epsilon$ be a familial
dimension function for $K$, and let $\zeta$ be a familial dimension function for $G\times K$
with $\zeta \dimpred \delta\times\epsilon$.
Let $\alpha$ be a virtual representation of $G$ and let $\beta$ be a virtual representation
of $K$.
The slant product
\[
 -\slant - \colon 
  \tilde H_{K,\epsilon}^\beta(Y;\Mackey T) \tensor
    \tilde H^{G\times K,\zeta}_{\alpha+\beta}(X\smsh Y; \MackeyOp U)
  \to \tilde H^{G,\delta}_\alpha(X; \Mackey T \mixprod \MackeyOp U).
\]
has the following properties.
\begin{enumerate}
\item
 It is natural in the following sense: Given a $G$-map
$f\colon X\to X'$, a $K$-map $g\colon Y\to Y'$, and elements
$y'\in \tilde H_{K,\epsilon}^\beta(Y';\Mackey T)$ and
$z\in \tilde H^{G\times K,\zeta}_{\alpha+\beta}(X\smsh Y; \MackeyOp U)$,
we have
\[
 y' \slant (f \smsh g)_*(z) = f_*(g^*(y') \slant z).
\]
Put another way, the slant product is a natural transformation in its adjoint form
\[
 \tilde H^{G\times K,\zeta}_{\alpha+\beta}(X\smsh Y; \MackeyOp U)
  \to \Hom(\tilde H_{K,\epsilon}^\beta(Y;\Mackey T), 
  		\tilde H^{G,\delta}_\alpha(X; \Mackey T \mixprod \MackeyOp U)).
\]

\item
It respects suspension in the sense that
\[
 (\sigma^W y) \slant (\sigma^{V+W} z) = \sigma^V(y \slant z).
\]

\item
It is associative in the following sense. 
Suppose given three groups, $G$, $K$, and $L$,
with respective familial dimension functions $\delta$, $\epsilon$, and $\zeta$.
Suppose that $\eta\dimpred\delta\times\epsilon$, $\theta\dimpred \epsilon\times\zeta$,
and $\kappa\dimpred \eta\times\zeta$, $\kappa\dimpred \delta\times\theta$
are also familial.
Given
$y\in \tilde H_{K,\epsilon}^\beta(Y;\Mackey R)$,
$z\in \tilde H_{L,\zeta}^\gamma(Z;\Mackey S)$, and
$w\in \tilde H^{G\times K\times L,\theta}_{\alpha+\beta+\gamma}(X\smsh Y\smsh Z;\MackeyOp U)$,
we have 
\[
 (y\cup z)\slant w = y\slant(z\slant w).
\]
That is, the following diagram commutes:
\[
 \def\objectstyle{\scriptstyle}
 \def\labelstyle{\scriptstyle}
 \xymatrix{
  \tilde H_{K,\epsilon}^\beta(Y;\Mackey R)\tensor \tilde H_{L,\zeta}^\gamma(Z;\Mackey S)
   \tensor \tilde H^{G\times K\times L,\kappa}_{\alpha+\beta+\gamma}(X\smsh Y\smsh Z;\MackeyOp U)
   \ar[r] \ar[d]
   & \tilde H_{K,\epsilon}^\beta(Y;\Mackey R)\tensor
    \tilde H^{G\times K,\eta}_{\alpha+\beta}(X\smsh Y;\Mackey S\mixprod \MackeyOp U) \ar[d] \\
  \tilde H_{K\times L,\theta}^{\beta+\gamma}(Y\smsh Z;\Mackey R\exboxprod\Mackey S)
   \tensor \tilde H^{G\times K\times L,\kappa}_{\alpha+\beta+\gamma}(X\smsh Y\smsh Z;\MackeyOp U)
   \ar[r]
   & \tilde H^{G,\delta}_\alpha(X;\Mackey R\mixprod(\Mackey S\mixprod\MackeyOp U))
  }
\]
Note that we're using the algebraic identity
$(\Mackey R\exboxprod\Mackey S)\mixprod \MackeyOp U 
\iso \Mackey R\mixprod(\Mackey S\mixprod \MackeyOp U)$.

\item
It respects the Wirthm\"uller isomprhism: Suppose that
$J\leq G$ and $L\leq K$,
that $J\in\F(\delta)$, $K\in\F(\epsilon)$, and $J\times L\in\F(\zeta)$,
and that $\zeta(G\times K/J\times L) = \delta(G/K)\dirsum \epsilon(K/L)$.
Then the following diagram commutes:
\[
 \def\objectstyle{\scriptstyle}
 \def\labelstyle{\scriptstyle}
 \xymatrix{
  \tilde H_{K,\epsilon}^\beta(K_+\smsh_L Y;\Mackey T) \tensor
    \tilde H^{G\times K,\zeta}_{\alpha+\beta}(G_+\smsh_J X\smsh K_+\smsh_L Y; \MackeyOp U)
    \ar[r] \ar[d]_\iso 
   & \tilde H^{G,\delta}_\alpha(G_+\smsh_J X; \Mackey T \mixprod \MackeyOp U) \ar[d]^\iso \\
  \tilde H_{L,\epsilon}^{\beta-\epsilon(K/L)}(Y;\Mackey T|L) \tensor
    \tilde H^{J\times L,\zeta}_{\alpha+\beta-\zeta(G\times K/J\times L)}(X\smsh Y; \MackeyOp U|J\times L)
    \ar[r]
   & \tilde H^{J,\delta}_{\alpha-\delta(G/J)}(X;(\Mackey T\mixprod\MackeyOp U)|J)
   }
\]
In the bottom row we are implicitly using a map
$(\Mackey T|L)\mixprod(\MackeyOp U|J\times L) \to (\Mackey T\mixprod\MackeyOp U)|J$
that need not be an isomorphism.
\qed
 \end{enumerate}
 \end{theorem}

Naturality implies that we have the following version of the slant product
for the Mackey functor-valued theories:
\[
 -\slant - \colon 
  \Mackey H_{K,\epsilon}^\beta(Y;\Mackey T) \mixprod
    \MackeyOp H^{G\times K,\zeta}_{\alpha+\beta}(X\smsh Y; \MackeyOp U)
  \to \MackeyOp H^{G,\delta}_\alpha(X; \Mackey T \mixprod \MackeyOp U).
\]

We're most interested in the internalization of the slant product to the diagonal $\Delta\leq G\times G$.
So, we also call the following a slant product.

\begin{definition}\label{def:internalslant}
Let $\delta$ and $\epsilon$ be familial dimension functions for $G$ and let
$\zeta\dimpred \delta\times\epsilon$ be a familial dimension function for $G\times G$;
assume that $\Delta\in\F(\zeta)$ and
write $\zeta$ again for $\zeta|\Delta$.
If $\Mackey T$ is a contravariant $\epsilon$-$G$-Mackey functor and
$\MackeyOp U$ is a covariant $\zeta$-$G$-Mackey functor, write
\[
 \Mackey T\mixprod \MackeyOp U = \Mackey T \mixprod [(G\times G)\times_\Delta \MackeyOp U],
\]
a covariant $\delta$-$G$-Mackey functor defined using the external version of $\mixprod$ on the right.
The {\em internal slant product}
\[
 -\slant- \colon 
   \tilde H^\beta_{G,\epsilon}(Y;\Mackey T) \tensor
   \tilde H^{G,\zeta}_{\alpha+\beta-\zeta(G\times G/\Delta)}(X\smsh Y;\MackeyOp U)
  \to \tilde H^{G,\delta}_\alpha(X;\Mackey T\mixprod \MackeyOp U)
\]
is then the composite
\begin{align*}
 \tilde H^\beta_{G,\epsilon}(Y;\Mackey T) &\tensor
		\tilde H^{G,\zeta}_{\alpha+\beta-\zeta(G\times G/\Delta)}(X\smsh Y;\MackeyOp U) \\
  &\to \tilde H^\beta_{G,\epsilon}(Y;\Mackey T) \tensor
       \tilde H^{G,\zeta}_{\alpha+\beta-\zeta(G\times G/\Delta)}(X\smsh Y;((G\times G)\times_\Delta\MackeyOp U)|\Delta) \\
  &\iso \tilde H^\beta_{G,\epsilon}(Y;\Mackey T) \tensor
       \tilde H^{G\times G,\zeta}_{\alpha+\beta}((G\times G)_+\smsh_\Delta(X\smsh Y); (G\times G)\times_\Delta\MackeyOp U) \\
  &\to \tilde H^\beta_{G,\epsilon}(Y;\Mackey T) \tensor
       \tilde H^{G\times G,\zeta}_{\alpha+\beta}(X\smsh Y; (G\times G)\times_\Delta\MackeyOp U) \\
  &\to \tilde H^{G,\delta}_\alpha(X;\Mackey T\mixprod \MackeyOp U).
\end{align*}
The first map is the unit $\MackeyOp U\to ((G\times G)\times_\Delta\MackeyOp U)|\Delta$,
the second is the Wirthm\"uller isomorphism, the third is induced by the $(G\times G)$-map
$(G\times G)_+\smsh_\Delta(X\smsh Y) \to X\smsh Y$, and the last map is the external slant product.
\end{definition}

Of course, this gives an internal slant product of Mackey functor-valued theories as well:
\[
 -\slant- \colon 
   \Mackey H^\beta_{G,\epsilon}(Y;\Mackey T) \mixprod
   \MackeyOp H^{G,\zeta}_{\alpha+\beta-\zeta(G\times G/\Delta)}(X\smsh Y;\MackeyOp U)
  \to \MackeyOp H^{G,\delta}_\alpha(X;\Mackey T\mixprod \MackeyOp U)
\]

We can now use the internal slant product to define evaluation and the cap product.

\begin{definition}
Let $\delta$ be a familial dimension function for $G$.
The {\em evaluation map}
\[
 \langle -,- \rangle \colon
  \tilde H_{G,\delta}^\beta(X;\Mackey T) 
     \tensor \tilde H^{G,\delta}_{\alpha+\beta}(X;\MackeyOp U)
  \to \tilde H^{G,0}_\alpha(S^0;\Mackey T\mixprod \MackeyOp U)
\]
is the slant product
\[
 -\slant- \colon
  \tilde H_{G,\delta}^\beta(X;\Mackey T) 
     \tensor \tilde H^{G,\delta_\Delta}_{\alpha+\beta}(S^0\smsh X;\MackeyOp U)
  \to \tilde H^{G,0}_\alpha(S^0;\Mackey T\mixprod \MackeyOp U).
\]
Here we are using $\delta_\Delta\dimpred 0\times\delta$ and the fact that
$\delta_\Delta(G\times G/\Delta) = 0$ and $\delta_\Delta|\Delta = \delta$.
\end{definition}

Note that we can express the naturality of evaluation by saying that the adjoint map
\[
 \tilde H_{G,\delta}^\beta(X;\Mackey T) \to
  \Hom(\tilde H^{G,\delta}_{\alpha+\beta}(X;\MackeyOp U), \tilde H^{G,0}_\alpha(S^0;\Mackey T\mixprod \MackeyOp U))
\]
is natural in $X$. When $\alpha = n\in\Z$, 
$\tilde H^{G,0}_n(S^0;\Mackey T\mixprod\MackeyOp U)$ is nonzero only when $n = 0$ and is then
$(\Mackey T\mixprod\MackeyOp U)(G/G) \iso \Mackey T\tensor_{\sorb G}\MackeyOp U$, giving the evaluation
\[
 \tilde H_{G,\delta}^\beta(X;\Mackey T) \to
  \Hom(\tilde H^{G,\delta}_{\beta}(X;\MackeyOp U), \Mackey T\tensor_{\sorb G}\MackeyOp U).
\]
There are many other interesting variations available, which we leave to the imagination
of the reader.

\begin{definition}\label{def:capproduct}
Let $\delta$, $\epsilon$, and $\zeta$ be as in Definition~\ref{def:internalslant}.
The {\em cap product}
\[
 -\cap- \colon
  \tilde H_{G,\epsilon}^{\beta}(X;\Mackey T)
    \tensor \tilde H^{G,\zeta}_{\alpha+\beta-\zeta(G\times G/\Delta)}(X;\MackeyOp U)
  \to \tilde H^{G,\delta}_{\alpha}(X;\Mackey T\mixprod\MackeyOp U)
\]
is the composite
\begin{align*}
 \tilde H_{G,\epsilon}^{\beta}(X;\Mackey T)
    &\tensor \tilde H^{G,\zeta}_{\alpha+\beta-\zeta(G\times G/\Delta)}(X;\MackeyOp U) \\
 &\to \tilde H_{G,\epsilon}^{\beta}(X;\Mackey T)
    \tensor \tilde H^{G,\zeta}_{\alpha+\beta-\zeta(G\times G/\Delta)}(X\smsh X;\MackeyOp U) \\
 &\to \tilde H^{G,\delta}_{\alpha}(X;\Mackey T\mixprod\MackeyOp U)
\end{align*}
where the first map is induced by the diagonal $X\to X\smsh X$.
\end{definition}

We get interesting special cases by considering particular choices of $\delta$, $\epsilon$, and $\zeta$.
One we'll use later is the case $\delta = \Lie - \epsilon$ and $\zeta = \Lie_\Delta$.
Using the facts that $\Lie_\Delta(G\times G/\Delta) = 0$ and
$\Lie_\Delta|\Delta = \Lie$, the cap product in this case takes the form
\[
 -\cap- \colon
  \tilde H_{G,\epsilon}^{\beta}(X;\Mackey T)
    \tensor \tilde \H^G_{\alpha+\beta}(X;\MackeyOp U)
  \to \tilde H^{G,\Lie-\epsilon}_{\alpha}(X;\Mackey T\mixprod\MackeyOp U).
\]
Specializing further, we can use $\MackeyOp U = \Mackey A_{G/G}$
(considered as a covariant $\Lie$-$G$-Mackey functor) and
the following calculation.

\begin{proposition}\label{prop:mixprodidentity}
If $\Mackey T$ is an $\epsilon$-$G$-Mackey functor, then
$\Mackey T\mixprod \Mackey A_{G/G} \iso \Mackey T$ as covariant
$(\Lie-\epsilon)$-$G$-Mackey functors.
\end{proposition}

\begin{proof}
Write $\delta = \Lie-\epsilon$. We first compute
\begin{align*}
 (p^*i_! &(G\times G)\times_\Delta \Mackey A_{G/G})((G/H,\delta)\tensor (G/K,\epsilon)) \\
  &\iso h(G\times G)\Spec{}{}(\susp^\infty (G\times G)/\Delta_+, 
  							(G_+\smsh_H S^{-\delta(G/H)})\smsh (G_+\smsh_K S^{-\epsilon(G/K)})) \\
  &\iso hG\Spec{}{}(S, (G_+\smsh_H S^{-\delta(G/H)})\smsh (G_+\smsh_K S^{-\epsilon(G/K)})) \\
  &\iso hG\Spec{}{}(G_+\smsh_H S^{-\epsilon(G/H)}, G_+\smsh_K S^{-\epsilon(G/K)})\\
  &= \MackeyOp A^{(G/H,\epsilon)}(G/K,\epsilon),
\end{align*}
where we use duality to get the third isomorphism. Therefore, we have
\begin{align*}
 (\Mackey T\mixprod \Mackey A_{G/G})(G/H,\delta)
  &= \Mackey T\tensor_{\sorb{G,\epsilon}} 
  		(p^*i_!(G\times G)\times_\Delta \Mackey A_{G/G})((G/H,\delta)\tensor -) \\
  &\iso \Mackey T\tensor_{\sorb{G,\epsilon}} \MackeyOp A^{(G/H,\epsilon)} \\
  &\iso \Mackey T(G/H,\epsilon),
\end{align*} 
as claimed.
\end{proof}

Thus, we have a cap product
\[
 -\cap- \colon
  \tilde H_{G,\epsilon}^{\beta}(X;\Mackey T)
    \tensor \tilde \H^G_{\alpha+\beta}(X;\Mackey A_{G/G})
  \to \tilde H^{G,\Lie-\epsilon}_{\alpha}(X;\Mackey T).
\]

Both the evaluation map and the cap product inherit properties
from Theorem~\ref{thm:slantproduct}. In particular, the associativity property
gives us the following:
\[
 (x\cup y)\cap z = x\cap (y\cap z)
\]
and
\[
 \langle x\cup y, z \rangle = \langle x, y\cap z \rangle
\]
when $x$, $y$, and $z$ lie in appropriate groups.

Let us now look at these pairings on the spectrum level.
The main result is the following calculation.

\begin{proposition}\label{prop:EMmixedproduct}
Let $\delta$ be a familial dimension function for $G$,
let $\epsilon$ be a familial dimension function for $K$, and
let $\zeta$ be a familial dimension function for $G\times K$ with $\zeta\dimpred \delta\times\epsilon$.
Let $\Mackey T$ be a contravariant $K$-$\epsilon$-Mackey functor and
let $\MackeyOp U$ be a covariant $(G\times K)$-$\zeta$-Mackey functor.
Then
\[
 \Mackey\pi_n^{G,\Lie-\delta}((H_\epsilon\Mackey T\smsh H^\zeta\MackeyOp U)^K) \iso
    \begin{cases}
      0 &\text{if $n<0$} \\
      \Mackey T\mixprod \MackeyOp U &\text{if $n=0$.}
    \end{cases}
\]
\end{proposition}

\begin{proof}
The proof is, in outline, the same as that of Proposition~\ref{prop:EMsmshhomotopy}.
The calculational input this time is that
\begin{multline*}
 h(G\times K)\Spec{}{}((G/J,\Lie-\delta),(K/L,\epsilon)\smsh (G\times K/M,\Lie-\zeta)) \\
   \iso h(G\times K)\Spec{}{}((G\times K/M,\zeta), (G/J,\delta)\smsh (K/L,\epsilon))
\end{multline*}
so that
\[
 h(G\times K)\Spec{}{}(-,(K/L,\epsilon)\smsh (G\times K/M,\Lie-\zeta))
   \iso \Mackey A_{K/L,\epsilon} \mixprod \MackeyOp A^{G\times K/M,\zeta}
\]
as modules over $\sorb{G,\delta}$, using Example~\ref{ex:mixprod}.
The result then follows by analyzing the structure provided by
Construction~\ref{con:EilenbergMacLane}, as in the proof of Proposition~\ref{prop:EMsmshhomotopy}.
\end{proof}

It follows that there is a map
\[
 (H_\epsilon\Mackey T\smsh H^\zeta\MackeyOp U)^K \to P^\delta(\Mackey T\mixprod\MackeyOp U)
\]
of $G$-spectra that is an isomorphism on $\Mackey\pi_0^{G,\Lie-\delta}$,
where $P^\delta(\Mackey T\mixprod\MackeyOp U)$ is a spectrum with $\Mackey\pi_n^{G,\Lie-\delta}$
homotopy concentrated in dimension 0. 
We then use that $E\F(\zeta)\times p_1^*E\F(\delta) \hmtpc E\F(\zeta)$ because
$\F(\zeta)\subset \F(\delta)\times\A(K)$, where $\A(K)$ is the collection of all subgroups of $K$,
so that $H^\zeta\MackeyOp U\smsh p_1^*E\F(\delta)_+ \hmtpc H^\zeta\MackeyOp U$,
to get a map
\begin{align*}
  (H_\epsilon\Mackey T\smsh H^\zeta\MackeyOp U)^K
  &\hmtpc (H_\epsilon\Mackey T\smsh H^\zeta\MackeyOp U\smsh p_1^* E\F(\delta)_+)^K \\
  &\hmtpc (H_\epsilon\Mackey T\smsh H^\zeta\MackeyOp U)^K\smsh E\F(\delta)_+ \\
  &\to P^\delta(\Mackey T\mixprod\MackeyOp U) \smsh E\F(\delta)_+ \\
  &\hmtpc H^\delta(\Mackey T\mixprod\MackeyOp U).
\end{align*}
The slant product
\[
 -\slant - \colon 
  \tilde H_{K,\epsilon}^\beta(Y;\Mackey T) \tensor
    \tilde H^{G\times K,\zeta}_{\alpha+\beta}(X\smsh Y; \MackeyOp U)
  \to \tilde H^{G,\delta}_\alpha(X; \Mackey T \mixprod \MackeyOp U).
\]
is then represented as follows:
\begin{align*}
 [Y, H_\epsilon\Mackey T\smsh S^\beta]_K
   \tensor{}& [S^{\alpha+\beta}, H^\zeta\MackeyOp U \smsh X\smsh Y]_{G\times K} \\
  &\to [S^{\alpha+\beta}, H^\zeta\MackeyOp U \smsh X\smsh H_\epsilon\Mackey T\smsh S^\beta]_{G\times K} \\
  &\iso [S^\alpha, H_\epsilon\Mackey T\smsh H^\zeta\MackeyOp U\smsh X]_{G\times K} \\
  &\iso [S^\alpha, (H_\epsilon\Mackey T\smsh H^\zeta\MackeyOp U\smsh X)^K]_G \\
  &\iso [S^\alpha, (H_\epsilon\Mackey T\smsh H^\zeta\MackeyOp U)^K \smsh X]_G \\
  &\to [S^\alpha, H^\delta(\Mackey T\mixprod\MackeyOp U)\smsh X]_G.
\end{align*}

\section{The Thom isomorphism and Poincar\'e duality}\label{sec:vThomIso}

Two of the most useful calculational results in nonequivariant homology are
the Thom isomorphism and Poincar\'e duality.
We give versions here for the $RO(G)$-graded theories.
These results apply only under very restrictive conditions, so we present them
here not so much for their own sake but as prelude to the more general
results we shall show for the $RO(\Pi B)$-graded theories we shall discuss later.

\subsection{The Thom isomorphism}

\begin{definition}
Let $p\colon E\to B$ be a $V$-bundle as in 
Example~\ref{ex:vcomplexes}(\ref{item:vbundle}).
Let $D(p)$ and $S(p)$ denote the disk and sphere bundles of $p$, respectively,
and let $T(p) = D(p)/S(p)$, the {\em Thom space} of $p$, be the quotient space.
A {\em Thom class} for $p$ is an element
$t\in \tilde H_G^V(T(p);\Mackey A_{G/G})$ such that, for each
$G$-map $b\colon G/K\to B$,
 \begin{align*}
 b^*(t) &{}\in \tilde H_G^V(T(b^*p);\Mackey A_{G/G}) \\
 &\iso \tilde H_G^V(\susp^V G/K_+;\Mackey A_{G/G}) \\
 &\iso A(K)
 \end{align*}
is a generator of $A(K)$ as an $A(K)$-module, i.e., a unit.
 \end{definition}

\begin{remarks}\label{rem:Thomclasshome}
A Thom class must live in ordinary cohomology $\tilde H^*_{G,\delta}$
with $\delta = 0$.
The isomorphism 
$ \tilde H_G^V(\susp^V G/K_+;\Mackey A_{G/G})
\iso A(K)$
would not hold with any other $\delta$.

The requirement that $p$ be a $V$-bundle is highly restrictive.
One of our motivations for discussing $RO(\Pi B)$-graded cohomology later is to
remove this restriction.
 \end{remarks}

A Thom class for $p$ is related to Thom classes for the fixed-point bundles
$p^K\colon E^K\to B^K$ as follows.

\begin{proposition}\label{prop:vThomclass}
The following are equivalent for a cohomology
class $t\in \tilde H_G^V(T(p);\Mackey A_{G/G})$.
 \begin{enumerate}
 \item
$t$ is a Thom class for $p$.
 \item
For every subgroup $K\leq G$,
$t|K\in \tilde H_K^V(T(p);\Mackey A_{K/K})$
is a Thom class for $p$ as a $K$-bundle.
 \item
For every subgroup $K\leq G$, 
$t^K\in \tilde H_{WK}^{V^K}(T(p^K); \Mackey A_{WK/WK})$
is a Thom class for $p^K$ as a $WK$-bundle.
 \item
For every subgroup $K\leq G$,
$t^K|e \in \tilde H^{|V^K|}(T(p^K);\Z)$
is a Thom class for $p^K$ as a nonequivariant bundle.
 \end{enumerate}
 \end{proposition}

\begin{proof}
To show that $(1) \Rightarrow (2)$,
let $b\colon K/J\to B$ be a $K$-map,
let $\tilde b\colon G/J\to B$ be the extension of $b$ to a $G$-map, and
consider the following diagram. 
\[
 \xymatrix{
  \tilde H_G^V(T(p);\Mackey A_{G/G}) \ar[r] \ar[d]_{\tilde b^*}
   & \tilde H_K^V(T(p);\Mackey A_{K/K}) \ar[d]^{b^*} \\
  \tilde H_G^V(T(\tilde b^*p);\Mackey A_{G/G}) \ar[r]_{\iso}
   & \tilde H_K^V(T(b^*p);\Mackey A_{K/K})
 }
\]
We need to show that the image of $t$, under the composite across the top
and right, is a generator of the copy of $A(J)$ on the bottom right. However,
by definition, the image of $t$ is a generator in the bottom left, and the
bottom map is an example of the Wirthm\"uller isomorphism.

To show that $(2) \Rightarrow (3)$,
it suffices to assume that $K$ is normal
in $G$ (otherwise, begin by restricting to $NK$).
Let $b\colon G/L\to B^K$ where $K\leq L\leq G$, and consider the following
diagram:
\[
 \xymatrix{
  \tilde H_G^V(T(p);\Mackey A_{G/G}) \ar[r] \ar[d]_{b^*}
   & \tilde H_{G/K}^{V^K}(T(p^K);\Mackey A_{(G/K)/(G/K)}) \ar[d]^{b^*} \\
  \tilde H_G^V(T(b^*p);\Mackey A_{G/G}) \ar[r] \ar[d]_{\iso}
   & \tilde H_{G/K}^{V^K}(T(b^*p^K);\Mackey A_{(G/K)/(G/K)}) \ar[d]^{\iso} \\
  \tilde H_G^V(\susp^V G/L_+;\Mackey A_{G/G}) \ar[r]_-{\iso}
   & \tilde H_{G/K}^{V^K}(\susp^{V^K} G/L_+;\Mackey A_{(G/K)/(G/K)})
  }
\]
Again, we need to show that the image of $t$ across the top and down the
right is a generator. By definition, the image down the left is a generator,
and the bottom map is an isomorphism because restriction to fixed sets
respects suspension, as shown in Section~\ref{subsec:quotientgroups}.

That $(3) \Rightarrow (4)$ follows from the implication $(1) \Rightarrow (2)$.
That $(4) \Rightarrow (1)$
follows from the fact that, for every $K$,
an element of $A(K) \iso \colim_W[S^W,S^W]_K$ is a generator if and only
if its restriction to $\colim_W[S^{W^J},S^{W^J}]$ is a unit (that is, $\pm 1\in\Z$)
for each subgroup $J$ of $K$.
(This follows from, for example, \cite[V.1.9]{LMS:eqhomotopy}.)
 \end{proof}

In order to state the Thom isomorphism, we need the correct definition of orientation.
The appropriate notion for $RO(G)$-graded cohomology is the following, which was called
a ``naive'' orientation in \cite[11.2]{CMW:orientation}.

\begin{definition}
Let $p\colon E\to B$ be a $V$-bundle.
If $b\colon G/H\to B$ is a $G$-map, 
let $E_b = p^{-1}(b(eH))$ be the corresponding fiber and
define an orientation of $E_b$ to be
a choice of a homotopy class $\phi(b)$ of
stable spherical $H$-maps $E_b\to V$, i.e., stable homotopy equivalences
$\susp_H^\infty S^{E_b}\to \susp_H^\infty S^V$.
A {\em strict orientation} of $p$ is a compatible collection 
$\{\phi(b)\}$ of orientations of each $E_b$, where compatibility means
that the following two conditions are satisfied:
\begin{enumerate}
\item
Let $b\colon G/H\to B$ and $\alpha\colon G/K\to G/H$, with $\alpha(eK) = gH$.
Then
\[
 \phi(b\circ\alpha) = g\phi(b)g^{-1}.
\]

\item
Let $\omega\colon G/H\times I\to B$ be a homotopy from $\omega_0$ to $\omega_1$ and
let $\bar\omega\colon E_{\omega_0}\to E_{\omega_1}$ be the induced
homotopy class of $H$-linear isometries. 
Then 
\[
 \phi(\omega_0) = \phi(\omega_1)\circ\bar\omega.
\]
\end{enumerate}
We say that $p$ is {\em strictly orientable} if there exists
a strict orientation of $p$.
\end{definition}

\begin{theorem}[Thom Isomorphism]\label{thm:vthomiso}
If $p\colon E\to B$ is a $V$-bundle, then $p$ is strictly orientable
if and only if it has a Thom class $t\in \tilde H_G^V(T(p);\Mackey A_{G/G})$.
Moreover, there is a one-to-one correspondence between strict orientations
of $p$ and Thom classes for $p$.
For any Thom class $t$,
the map
 \[
 t\cup - \colon H_{G,\delta}^\alpha(B;\Mackey U)
 \to \tilde H_{G,\delta}^{\alpha+V}(T(p);\Mackey U)
 \]
is an isomorphism for any familial $\delta$.
\end{theorem}

In the statement, the cup product uses the special case of $p_2^*\delta \dimpred 0\times\delta$
(see the discussion preceding Theorem~\ref{thm:vcupproduct})
and is internalized along the diagonal map $T(p) \to T(p)\smsh B_+$.

\begin{proof}
The theorem is clear when $p$ is trivial;
the correspondence between Thom classes and orientations comes from the correspondence between
units in $A(K)$ and stable homotopy classes of self-equivalences of spheres.
The general case follows, as it does nonequivariantly, by a Mayer-Vietoris
patching argument (see \cite{MS:characteristicclasses} or 
\cite{CW:duality}).
 \end{proof}

\subsection{Poincar\'e duality}

We now outline a proof of Poincar\'e duality for compact oriented $V$-manifolds.
The noncompact case can be handled in the usual way using cohomology with
compact supports, as in \cite{MS:characteristicclasses} or 
\cite{CW:duality}.

\begin{definition}
Let $M$ be a closed $V$-manifold as in Example~\ref{ex:vcomplexes}(\ref{item:vmanifold}). 
A {\em fundamental class} of $M$ is
a class $[M]\in \H_V^G(M;\Mackey A_{G/G})$ such that, for each 
point $m\in M$ with
$G$-orbit $Gm \subset M$ and tangent plane $\tau_m\iso V$,
the image of $[M]$ in
 \begin{align*}
 \H_V^G(M,M-Gm;{}\Mackey A_{G/G})
 &\iso \tilde\H_V^G(G_+\smsh_{G_m} S^{V-\Lie(G/G_m)};
    \Mackey A_{G/G}) \\
 &\iso \tilde\H_V^{G_m}(S^V; \Mackey A_{G_m/G_m}) \\
 &\iso A(G_m)
 \end{align*}
is a generator.
 \end{definition}

Recall that $\H^G_* = H^{G,\Lie}_*$.
The first isomorphism above comes from excision, using a tubular neighborhood of the orbit $Gm$.
The second isomorphism is the Wirth\-m\"uller isomorphism, which would not
be true in this form if we were to use any $\delta$ other than $\Lie$, because there would
then be a shift in dimension.
A fundamental class must live in dual homology.

The fundamental class $[M]$ is related to fundamental classes of the fixed
submanifolds $M^K$ as follows. The proof is similar to that of
Proposition~\ref{prop:vThomclass}.

\begin{proposition}\label{prop:vFundclass}
Let $M$ be a closed $V$-manifold.
The following are equivalent for a homology
class $\mu\in \H_V^G(M;\Mackey A_{G/G})$:
 \begin{enumerate}
 \item
$\mu$ is a fundamental class for $M$ as a $G$-manifold.
 \item
For every subgroup $K\leq G$,
$\mu|K \in \H_V^K(M;\Mackey A_{K/K})$ is a fundamental class for
$M$ as a $K$-manifold.
 \item
For every subgroup $K\leq G$,
$\mu^K \in \H_{V^K}^{WK} (M^K; \Mackey A_{WK/WK})$ is a fundamental class for
$M^K$ as a $WK$-manifold.
 \item
For every subgroup $K\leq G$,
$\mu^K|e \in H_{|V^K|}(M^K;\Z)$ is a fundamental class for $M^K$
as a nonequivariant manifold. \qed
 \end{enumerate}
\end{proposition}

We say that a $V$-manifold is strictly orientable if its tangent bundle
is strictly orientable in the sense of the preceding section.

\begin{theorem}[Poincar\'e Duality]
A closed $V$-dimensional $G$-manifold $M$ is strictly orientable if and only if
it has a fundamental class $[M]\in\H_V^G(M;\Mackey A_{G/G})$.
If $M$ does have a fundamental class $[M]$, then
 \[
 -\cap [M] \colon H_{G,\delta}^\alpha(M;\Mackey T) \to H_{V-\alpha}^{G,\Lie-\delta}(M;\Mackey T)
 \]
is an isomorphism for every familial $\delta$ such that $\F(\delta)$ contains
every isotropy subgroup of $M$.
 \end{theorem}

Here, we are using the special case of the general cap product given
after Proposition~\ref{prop:mixprodidentity}.
Notice the following particular cases:
When $\delta = 0$ we get the isomorphism
 \[
 -\cap [M] \colon H_{G}^\alpha(M;\Mackey T) \to \H_{V-\alpha}^{G}(M;\Mackey T)
 \]
and when $\delta = \Lie$ we get the isomorphism
 \[
 -\cap [M] \colon \H_{G}^\alpha(M;\Mackey T) \to H_{V-\alpha}^{G}(M;\Mackey T).
 \]

\begin{proof}
We can adapt the standard geometric proof, as given in 
\cite{MS:characteristicclasses} or \cite{CW:duality}.
Recall that the argument starts with the local case of a 
tubular neighborhood of an orbit $G/G_m$ in $M$.
This neighborhood will have the form $G\times_{G_m}\bar D(V-\Lie(G/G_m))$.
Letting $\mu \in \H^G_V(G\times_{G_m}\bar D(V-\Lie(G/G_m));\Mackey A_{G/G})$ be
the restriction of $[M]$, we then have the following diagram
(recall that we write $\bar D$ as shorthand for the pair $(D,S)$):
\[
 \def\objectstyle{\scriptstyle}
 \def\labelstyle{\scriptstyle}
 \xymatrix{
  H_{G,\delta}^\alpha(G\times_{G_m}\bar D(V-\Lie(G/G_m));\Mackey T)
    \ar[r]^-{-\cap\mu} \ar[d]_\iso
   & H^{G,\Lie-\delta}_{V-\alpha}(G\times_{G_m} D(V-\Lie(G/G_m);\Mackey T) \ar[d]^\iso \\
  H_{G_m,\delta}^{\alpha-\delta(G/G_m)}(\bar D(V-\Lie(G/G_m));\Mackey T|G_m)
    \ar[r]^-{-\cap\mu|G_m} \ar[d]_\iso
   & H^{G_m,\Lie-\delta}_{V-\alpha-\Lie(G/G_m)+\delta(G/G_m)}(D(V-\Lie(G/G_m)); \Mackey T|G_m)
    \ar[d]^\iso \\
  H_{G_m,\delta}^{\alpha-V+\Lie(G/G_m)-\delta(G/G_m)}(*;\Mackey T|G_m)
    \ar[r]^-{-\cap u}_\iso
   & H^{G_m,\Lie-\delta}_{V-\alpha-\Lie(G/G_m)+\delta(G/G_m)}(*; \Mackey T|G_m)
 }
\]
The top square commutes because the cap product respects the Wirthm\"uller isomorphism;
the bottom square commutes because it respects suspension.
The element $u\in \H^{G_m}_0(*;\Mackey A_{G_m/G_m}) \iso A(G_m)$ is the image
of $\mu$ and is, by assumption, a unit, making the map along the bottom an isomorphism
between two copies of $\Mackey T(G/G_m)$.
It follows that the other two horizontal maps are isomorphisms as well.
Note, however, that this argument requires that $G_m\in\F(\delta)$; we cannot
expect to get a local isomorphism otherwise.

The argument then proceeds via a Mayer-Vietoris patching argument
as in the references above.
 \end{proof}

As with the Thom isomorphism theorem, we will lift both the orientability
assumption and the assumption that $M$ is $V$-dimensional in
Section~\ref{sec:genThom}.

We can also prove Poincar\'e duality from the Thom isomorphism by using
the representing spectra. For this, we use the fact that, if $M$ is embedded
in a $G$-representation $W$ with normal bundle $\nu$, then the stable dual of $M_+$
is $\susp^{-W}T\nu$ \cite{LMS:eqhomotopy}. 
Thus, we have the following chain of isomorphisms, where the first
is the Thom isomorphism.
\begin{align*}
 H_{G,\delta}^\alpha(M;\Mackey T)
 &\iso \tilde H_{G,\delta}^{\alpha+W-V}(T\nu;\Mackey T) \\
 &\iso [T\nu, \susp^{\alpha+W-V} H_\delta\Mackey T]_G \\
 &\iso [\susp^W D(M_+), \susp^{\alpha+W-V} H_\delta\Mackey T]_G \\
 &\iso [S^W, \susp^{\alpha+W-V}H_\delta\Mackey T\smsh M_+]_G \\
 &\iso H^{G,\Lie-\delta}_{V-\alpha}(M;\Mackey T).
\end{align*}
For the last isomorphism, we use the fact that, because $\F(\delta)$ contains
every isotropy subgroup of $M$, we have $M\hmtpc E\F(\delta) \times M$,
hence
\[
 H_\delta\Mackey T \smsh M_+ \hmtpc H_\delta\Mackey T\smsh E\F(\delta)_+\smsh M_+
  \hmtpc H^{\Lie-\delta}\Mackey T\smsh M_+.
\]
As usual, we can take as a fundamental class 
$[M] \in \H^G_{V}(M;\Mackey A_{G/G})$ the image of a Thom class
in $\tilde H_G^{W-V}(T\nu;\Mackey A_{G/G})$. The maps displayed above
can then be shown to be given by capping with this fundamental class.

If $M$ is a compact orientable $V$-manifold with boundary, then we can prove relative,
or Lefschetz, duality. We state the results.

\begin{definition}
Let $M$ be a compact $V$-manifold with boundary. 
A {\em fundamental class} of $M$ is
a class $[M,\bndry M]\in \H_V^G(M,\bndry M;\Mackey A_{G/G})$ such that, for each 
point $m\in M-\bndry M$ with
$G$-orbit $Gm \subset M$ and tangent plane $\tau_m\iso V$,
the image of $[M,\bndry M]$ in
 \begin{align*}
 \H_V^G(M,M-Gm;\Mackey A_{G/G})
 &\iso \tilde\H_V^G(G_+\smsh_{G_m}S^{V-\Lie(G/G_m)};
    \Mackey A_{G/G}) \\
 &\iso \tilde\H_V^{G_m}(S^{V}; \Mackey A_{G_m/G_m}) \\
 &\iso A(G_m)
 \end{align*}
is a generator.
 \end{definition}

There is an obvious relative version of Proposition~\ref{prop:vFundclass}.
Finally, we have the following relative duality.

\begin{theorem}[Lefschetz Duality]
A compact $V$-dimensional $G$-manifold $M$ is strictly orientable if and only if
it has a fundamental class $[M,\bndry M]$.
If $M$ does have a fundamental class, then the following are
isomorphisms for every familial dimension function $\delta$
such that $\F(\delta)$ contains every isotropy subgroup of $M$:
 \[
 -\cap [M,\bndry M] \colon H_{G,\delta}^\alpha(M;\Mackey T) \to
   H_{V-\alpha}^{G,\Lie-\delta}(M,\bndry M;\Mackey T)
 \]
and
\[
 -\cap [M,\bndry M] \colon H_{G,\delta}^\alpha(M,\bndry M;\Mackey T) \to
   H_{V-\alpha}^{G,\Lie-\delta}(M;\Mackey T).
\]
\qed
 \end{theorem}

Note that these versions of the cap product are obtained by internalizing the
slant product along the diagonal map $M/\bndry M \to M/\bndry M\smsh M_+$
in the first case and
$M/\bndry M \to M_+\smsh M/\bndry M$ in the second.

Again, we have the following special cases when we take $\delta=0$ or $\delta=\Lie$:
\begin{align*}
 -\cap [M,\bndry M] &\colon H_{G}^\alpha(M;\Mackey T) \to
   \H_{V-\alpha}^{G}(M,\bndry M;\Mackey T),
\\
 -\cap [M,\bndry M] &\colon H_{G}^\alpha(M,\bndry M;\Mackey T) \to
   \H_{V-\alpha}^{G}(M;\Mackey T),
\\
 -\cap [M,\bndry M] &\colon \H_{G}^\alpha(M;\Mackey T) \to
   H_{V-\alpha}^{G}(M,\bndry M;\Mackey T), \qquad\text{and}
\\
 -\cap [M,\bndry M] &\colon \H_{G}^\alpha(M,\bndry M;\Mackey T) \to
   H_{V-\alpha}^{G}(M;\Mackey T).
\end{align*}

\section{Miscellaneous remarks}\label{sec:OrdinaryRemarks}

\subsection{Ordinary homology of $G$-spectra}

In the nonparametrized context, there is no significant difficulty to
developing a theory of $G$-CW spectra and resulting cellular homology and cohomology
theories. The reason we didn't do this is that there are signficant problems
in the parametrized case, as discussed in \cite[Ch.~24]{MaySig:parametrized}.
We sketch here how the nonparametrized case would go.

Given a dimension function $\delta$ for $G$ and a virtual representation
$\alpha = V \ominus W$, we define a {\em $\delta$-$\alpha$-cell} to be a spectrum
of the form
\[
 G_+\smsh_H S^{\alpha-\delta(G/H)+q}
  = G_+\smsh_H \susp_H^{V+q} \susp^\infty_{W+\delta(G/H)} S^0,
\]
where $\susp^\infty_{W+\delta(G/H)}$ is the shift desuspension functor
of \cite[I.4.1]{LMS:eqhomotopy} or \cite[II.4.7]{MM:orthogonal} and
$G_+\smsh_H-$ is the functor defined in 
\cite[II.4.1]{LMS:eqhomotopy} or \cite[V.2.3]{MM:orthogonal}.
If $E$ is a $G$-spectrum, we say that the cofiber of a map
$G_+\smsh_H S^{\alpha-\delta(G/H)+q-1} \to E$ is the result of
{\em attaching a cell of dimension $\alpha+q$} to $E$.

We can then generalize the definitions of \cite[I.5]{LMS:eqhomotopy}
A {\em $\delta$-$\alpha$-cell spectrum} is a $G$-spectrum $E$ together
with a sequence of subspectra $E_n$, $n\geq 0$, such that $E$ is the union of the $E_n$,
with $E_0 = *$ and
$E_{n+1}$ obtained from $E_{n}$ as the cofiber of a map $J_n\to E_n$
where $J_n$ is a wedge of $\delta$-$\alpha$-cells.
A {\em $\delta$-$G$-CW$(\alpha)$ spectrum} is a $G$-cell spectrum $E$
such that each attaching map
\[
 G_+\smsh_H S^{\alpha-\delta(G/H)+q-1} \to E_n
\]
factors through a cell subspectrum containing only cells of dimension $< \alpha+q$.
We define the {\em $(\alpha+n)$-skeleton} $E^{\alpha+n}$ to be the union of the cells
of dimension at most $\alpha+n$.

If $E$ is any $G$-spectrum, let
\[
 \pi_{\alpha+n}^{H,\delta}(E) = [G_+\smsh_H S^{\alpha-\delta(G/H)+n}, E]_G.
\]
We say that a map $f\colon E\to F$ is
a {\em $\delta$-weak$_\alpha$ equivalence} if, for each subgroup $H\leq G$ and integer $n$,
$\pi_{\alpha+n}^{H,\delta}(E) \to \pi_{\alpha+n}^{H,\delta}(F)$
is an isomorphism.
Using Corollary~\ref{cor:ordinaryCellAsGeneral} and
Lemma~\ref{lem:generalCellAsOrdinary} we can show that $f\colon E\to F$ is a
$\delta$-weak$_\alpha$ equivalence if and only if it is
an $\F(\delta)$-equivalence.

With these definitions we get a HELP lemma, approximation of
spectra by $\delta$-$G$-CW$(\alpha)$ spectra
(up to $\delta$-weak$_\alpha$ equivalence), and cellular approximation
of maps between $\delta$-$G$-CW$(\alpha)$ spectra.

To define the chains of an arbitrary $G$-spectrum $E$ we first take a
$\delta$-$G$-CW$(\alpha)$ approximation $\Gamma E\to E$ and let
\[
 \Mackey C_{\alpha+n}^{G,\delta}(E)(G/H,\delta) =
  \pi_{\alpha+n}^{H,\delta}((\Gamma E)^{\alpha+n}/(\Gamma E)^{\alpha+n-1}).
\]
This gives a chain complex, independent of the choice of approximation $\Gamma E$
up to chain homotopy equivalence. We then define homology and cohomology groups of $E$
using this chain complex.

Homology and cohomology theories defined on spectra are representable.
In particular, for each contravariant $\delta$-Mackey functor $\Mackey T$ we get a spectrum
$H_\delta\Mackey T$ with the same characterization as in Section~\ref{sec:representing}. 
Of importance to the theories defined on spaces, we can now conclude
that $H_\delta\Mackey T$ is unique {\em up to unique stable equivalence}.
On the other hand, if $\MackeyOp S$ is a covariant $\delta$-Mackey functor,
then Theorem~\ref{thm:dualCohomologyRep} tells us that
$H^\delta\MackeyOp S \hmtpc H_{\Lie-\delta}\MackeyOp S\smsh E\F(\delta)_+$
and we get a similar strong uniqueness result for this spectrum.

\subsection{Model categories}

Although we've deliberately kept the discussion as elementary as possible,
there were various points at which
we could have phrased things using the language of model categories.
Indeed, for a given dimension function $\delta$ and virtual representation $\alpha$,
we could define a compactly generated model category structure on $G$-spaces in which
the weak equivalences are the $\delta$-weak$_\alpha$ equivalences and the generating
cofibrations are the maps
$G\times_H S(Z) \to G\times_H D(Z)$ where $G\times_H\bar D(Z)$ is a $\delta$-$\alpha$-cell.
We leave the details to the interested reader.

Although this is an interesting point of view, it's not clear that it would add much to
the exposition. Also, CW objects are at the heart of the definition of cellular homology
and cohomology, but are rarely mentioned in discussions of model categories.
It would be interesting to see if a general theory of CW objects in model categories
would be of any use in other contexts.

\chapter{Parametrized Homotopy Theory and Fundamental Groupoids}
\label{chap:homotopytheory}

\section*{Introduction}

In this chapter we discuss some needed background material on
parametrized spaces and spectra. Much of this is summarized from
\cite{MaySig:parametrized}. We add some material on {\em lax maps,}
which requires additional discussion, beyond that provided by May and Sigurdsson, of the Moore path
fibration functor.

We also summarize material on the fundamental groupoid.
The equivariant fundamental groupoid was considered by
tom~Dieck in \cite{tD:transfgroups}.
It was used in \cite{CMW:orientation}
to give a theory of equivariant orientations,
where it provides the necessary machinery to keep track of
the local representations appearing in a $G$-vector bundle
or a smooth $G$-manifold.
In particular, a collection of local representations assembles into
what we call a representation of the fundamental groupoid.
Such representations can be thought of as the natural dimensions
of $G$-vector bundles or $G$-manifolds, and provide
the grading for our extension of ordinary homology and cohomology.

The stable orbit category played a prominent role
in Chapter~\ref{chap:rog}. In our extended theory a similar role
is played by the {\em stable fundamental groupoid}, which we introduce and discuss in
Section~5.

In Section~6 we discuss parametrized homology and cohomology theories in general,
in anticipation of discussing the particular case of cellular homology in the following chapter.
In Section~7 we discuss an alternative to fiberwise duality that is more appropriate
for parametrized homology theories.


\section{Parametrized spaces and lax maps}

In \cite{MaySig:parametrized} May and Sigurdsson give a careful exposition of the homotopy theory
of parametrized spaces and spectra. We recall the parts we need and refer the reader to \cite{MaySig:parametrized} for further details.

\begin{definition}
Let $B$ be a compactly generated $G$-space.
\begin{enumerate}
\item
Let $G\ParaU B$ be the category of {\em $G$-spaces over $B$}:
Its objects are pairs $(X,p)$ where $X$ is a $G$-$k$-space
and $p\colon X\to B$ is a $G$-map.
A map $(X,p)\to (Y,q)$ is a $G$-map
$f\colon X\to Y$ such that $q\circ f = p$, i.e., a $G$-map over $B$.
\item
Let $G\Para B$ be the category of {\em ex-$G$-spaces over $B$}:
Its objects are triples $(X,p,\sigma)$ where $(X,p)$ is a $G$-space
over $B$ and $\sigma$ is a section of $p$,
i.e., $p\circ\sigma$ is the identity.
A map $(X,p,\sigma)\to (Y,q,\tau)$ is a section-preserving $G$-map
$f\colon X\to Y$ over $B$, i.e., a $G$-map over and under $B$.
\end{enumerate}
When the meaning is clear, we shall write just $X$ for $(X,p)$ or $(X,p,\sigma)$. If $(X,p)$ is a space over $B$, we write
$(X,p)_+$ for the ex-$G$-space obtained by adjoining a disjoint section.
We shall also write $X_+$ for $(X,p)_+$, a notation that
May and Sigurdsson denigrate but should not cause confusion in context.
\end{definition}

As May and Sigurdsson point out, in order to make these categories closed symmetric monoidal we need $B$ to be weak Hausdorff but cannot require it of the spaces over $B$. Hence we assume that $B$ is compactly generated (i.e., it is a $k$-space and weakly Hausdorff) but assume of the spaces over it only that they are $k$-spaces.

On either category we can define a model category structure in which a map
$f$ is a weak equivalence, fibration, or cofibration if it is one as
a map of $G$-spaces, ignoring $B$. Here, as usual, a weak equivalence is
a weak equivalence on all components of all fixed sets, a fibration is
an equivariant Serre fibration, and a cofibration is a retract of
a relative $G$-CW complex. This is what May and Sigurdsson call the
$q$-model structure and they call these classes of maps $q$-equivalences, $q$-fibrations, and $q$-cofibrations.
However, they point out serious technical difficulties with this model structure, stemming from the fact that these cofibrations do not interact well with fiberwise homotopy. To rectify this, they define another model structure, the $qf$-model structure, with the same weak equivalences, but with a more restrictive class of cofibrations, the $qf$-cofibrations, and a corresponding notion of $qf$-fibration.
In any case, we are entitled to invert these weak equivalences, getting the homotopy categories
$\Ho G\ParaU B$ and $\Ho G\Para B$.
If $X$ and $Y$ are ex-$G$-spaces over $B$ we write
$[X,Y]_{G,B}$ for the set of maps in $\Ho G\Para B$.

In this context, inverting the weak equivalences introduces many more maps than are evident if one pays attention only to the fiber-preserving maps over $B$. For example, let $(X\times I,h)$ and $(Y,q)$ be spaces over $B$,
where $h\colon X\times I\to B$ does not necessarily factor through $X$.
Let $i_0, i_1\colon X\to X\times I$ be the inclusions of the two endpoints and let $p_0 = hi_0$ and $p_1 = hi_1$. Then the maps
$i_0\colon (X,p_0)\to (X\times I,h)$ and $i_1\colon (X,p_1)\to (X\times I,h)$ are both $q$-equivalences.
(That we really want these to be weak equivalences is implicit
already in Dold's axiom {\bf CYL} in \cite{Dold:parahomology}.)
Hence, if $f\colon (X,p_0)\to (Y,q)$ is a map over $B$, we have in the homotopy category the map
\[
 [f]\circ [i_0]^{-1}\circ [i_1]\in [(X,p_1), (Y,q)]_{G,B},
\]
which might not be represented by any fiberwise map $(X,p_1)\to (Y,q)$ over $B$.

There are a couple of useful ways of viewing this example.
Conceptually, the pair $(f,h)$ specifies a map $f\colon X\to Y$
and a homotopy $h\colon qf\to p_1$. This is an example of what we
shall call a {\em lax map} from $(X,p_1)$ to $(Y,q)$.
Technically, it helps to view $(f,h)$ as specifying a map over $B$
from $(X,p_1)$ to the Hurewicz fibration associate with $q$.
It will be most convenient to use the following variant of the associated fibration, given in \cite[\S8.3]{MaySig:parametrized}. Recall that the space of
{\em Moore paths} in $B$ is
\[
 \Lambda B = \{ (\lambda, l) \in B^{[0,\infty]}\times [0,\infty)
                \mid \lambda(r) = \lambda(l) \text{ for $r\geq l$} \}.
\]
We usually write $\lambda$ for $(\lambda,l)$ and $l_\lambda$ for $l$, the
{\em length} of $\lambda$.
Moore paths are useful largely because they have an associative composition:
Suppose $\lambda$ and $\mu$ are Moore paths with $\lambda(\infty) = \mu(0)$. Then we write $\lambda\mu$ for the Moore path with length
$l_{\lambda\mu} = l_\lambda + l_\mu$ given by
\[
 (\lambda\mu)(t) =
  \begin{cases}
    \lambda(t) & \text{if $t\leq l_\lambda$} \\
    \mu(t - l_\lambda) & \text{if $t\geq l_\lambda$.}
  \end{cases}
\]
(We reverse the order of composition here compared to
\cite[\S8.3]{MaySig:parametrized}, as it will be more useful to us in this order.)
If $p\colon X\to B$, the {\em Moore path fibration} $LX = L(X,p)$ is given by
\[
 LX = X\times_B \Lambda B
  = \{ (x,\lambda)\in X\times \Lambda B \mid p(x) = \lambda(0) \}.
\]
We define the projection $Lp\colon LX \to B$ using the endpoint projection: $Lp(x,\lambda) = \lambda(\infty)$. If $X$ is an ex-space with section $\sigma$, we define a section $L\sigma$ of $LX$ by $L\sigma(b) = (\sigma(b),b)$, where we write $b$ for the path of length 0 at the point $b$. There is an inclusion $\iota\colon X\to LX$ given by $\iota(x) = (x,p(x))$
and composition of paths defines a map
$\chi\colon L^2 X \to LX$, $\chi(x,\lambda,\mu) = (x,\lambda\mu)$.
With these maps, $(L,\iota,\chi)$ is a monad.
The unit $\iota$ is a $G$-homotopy equivalence under (but not over) $B$, hence a $q$-equivalence, and the projection $LX\to B$ is a Hurewicz fibration, hence a $q$-fibration and a $qf$-fibration.

For simplicity, restrict now to the ex-space case. (For the unbased case, just ignore sections in the following discussions.)
Because $\iota$ is a $q$-equivalence, it induces an isomorphism
\[
 [X,Y]_{G,B} \iso [X,LY]_{G,B}
\]
for any $(X,p,\sigma)$ and $(Y,q,\tau)$. So, we can represent more maps in $[X,Y]_{G,B}$ by considering maps from $X$ to $LY$ over $B$. A map $X\to LY$ over $B$ is given by a pair $(f,\lambda)$ where $f\colon X\to Y$ and $\lambda\colon X\to \Lambda B$ are $G$-maps under $B$, with $\lambda(x)(0) = qf(x)$ and $\lambda(x)(\infty) = p(x)$ for all $x$. 
We can view $\lambda$ as specifying a {\em Moore homotopy} rel $\sigma(B)$ from $q\circ f$ to $p$, by which we mean a $G$-map
\[
 \lambda\colon X\times [0,\infty] \to B
\]
and a $G$-invariant map
\[
 l_\lambda\colon X\to [0,\infty)
\]
for which each $\lambda(x,-)$ is a Moore path of length $l_\lambda(x)$, and such that $\lambda(-,0) = qf$, $\lambda(-,\infty) = p$, and $\lambda(\sigma(b),-)$ is the 0-length path at $b$ for every $b\in B$.

\begin{definition}\label{def:laxmap}
A {\em lax map} $(X,p,\sigma)\laxto (Y,q,\tau)$ is a map $X\to LY$ under and over $B$. It can also be thought of as a pair $(f,\lambda)$ where $f\colon X\to Y$ is a $G$-map under $B$ and $\lambda\colon qf\to p$ is a Moore homotopy rel $\sigma(B)$.
\end{definition}

In this context we shall call maps under and over $B$ {\em strict} maps. We can also think of strict maps as those lax maps in which the homotopies have zero length.
To avoid confusion we will use ``$\laxto$" in the text to indicate lax maps,
reserving ``$\to$" for strict maps.

We have an associative composition of lax maps: If
$(f,\lambda)\colon (X,p,\sigma)\laxto (Y,q,\tau)$ and
$(g,\mu)\colon (Y,q,\tau) \laxto (Z,r,\upsilon)$, then
\[
 (g,\mu)\circ (f,\lambda) = (gf, (\mu f)\lambda),
\]
where $\mu f\colon rgf \to qf$ is the evident Moore homotopy.
Thus, we can form the category $G \Lax B$ of ex-$G$-spaces and lax maps
(similarly, we have the category $G\LaxU B$ of $G$-spaces over $B$ and lax maps).
Thinking of lax maps as strict maps $X\to LY$, we have
\[
 G\Lax{B}(X,Y) = G\Para{B}(X,LY) \quad\text{and}\quad
 (G\LaxU{B})(X,Y) = (G\ParaU{B})(X,LY)
\]
and we topologize $G\Lax{B}$ and $G\LaxU{B}$ using these identifications.
The composite of $f\colon X\to LY$ and $g\colon Y\to LZ$ is then the composite
\[
 X \xrightarrow{f} LY \xrightarrow{Lg} L^2 Z \xrightarrow{\chi} LZ.
\]
Thus, the category of lax maps is the Kleisli category of the monad $L$
(see \cite{Kle:KleisliCat}, \cite[\S VI.5]{Mac:categories}, and \cite[\S5]{Str:MonadsI}).

One of the points of the Kleisli category is that the
inclusion $G\Para B\to G\Lax B$ has a right adjoint
$\bar L\colon G\Lax B\to G\Para B$. On objects it is given by $L$,
and, if $f\colon X\to LY$ represents a lax map, $\bar Lf$ is the composite
\[
 LX \xrightarrow{Lf} L^2Y \xrightarrow{\chi} LY.
\]
Both the inclusion $G\Para B\to G\Lax B$ and
$\bar L\colon G\Lax B\to G\Para B$ are continuous.

There is a notion of homotopy in $G\Lax B$, defined using lax maps out of cylinders $X\smsh_B I_+$ or paths in $G\Lax{B}(X,Y)$.
We've already noted that $\iota\colon X\to LX$ is a weak equivalence;
it's not hard to see that, considered as a lax map, it is a lax homotopy equivalence,
i.e., a homotopy equivalence in $G\Lax{B}$. Its inverse is the map
$\kappa\colon LX\laxto X$ represented by the identity $1\colon LX\to LX$.
We record the following useful facts.

\begin{lemma}\label{lem:laximpliesweak}
\begin{enumerate}\item[]
\item
If $f\colon X\to Y$ is a strict map that is also a lax homotopy
equivalence, then it is a weak equivalence (i.e., a $q$-equivalence).
\item
If $f\colon X\laxto Y$ is a lax homotopy equivalence,
then its representing map $\hat f\colon X\to LY$ is also a 
lax homotopy equivalence, hence a weak equivalence.
\end{enumerate}
\end{lemma}

\begin{proof}
To show the first claim, note that $\bar L f = Lf\colon LX\to LY$
is a fiberwise homotopy equivalence, because $\bar L$, being continuous,
takes lax homotopy equivalences to fiberwise homotopy equivalences.
It then follows, from the following diagram and the fact that $\iota$ is
a weak equivalence, that $f$ is a weak equivalence:
\[
 \xymatrix{
  X \ar[d]_\iota \ar[r]^f & Y \ar[d]^\iota \\
  LX \ar[r]_{Lf} & LY
 }
\]

For the second claim, an easy calculation shows that 
$\hat f = \bar L f\circ \iota$. Again, $\bar L$ takes the lax homotopy
equivalence $f$ to a fiberwise homotopy equivalence, hence a
lax homotopy equivalence. Combined with the fact that $\iota$ is
a lax homotopy equivalence, we get that $\hat f$ is a lax homotopy equivalence.
\end{proof}

If $Y$ is any ex-$G$-space, then $LY$ is $q$-fibrant.
If $X$ is $q$-cofibrant, we then have
\begin{align*}
 [X,Y]_{G,B}
  &\iso [X,LY]_{G,B} \\
  &\iso \pi G\Para B(X,LY) \\
  &\iso \pi G\Lax B(X,Y),
\end{align*}
where $\pi$ denotes homotopy classes of maps.
In other words, if $X$ is $q$-cofibrant and $Y$ is any space over $B$, then
$[X,Y]_{G,B}$ is precisely the set of lax homotopy classes of lax
maps from $X$ to $Y$.

We should point out that the category $G\Lax B$ is wanting in several respects. It does not, in general, have limits. In particular, it has no terminal object.
It also lacks colimits in general, although any diagram of strict maps has a colimit in $G\Lax B$, namely the image of the colimit in $G\Para B$ (because $G\Para B\to G\Lax B$ is a left adjoint).
It is, therefore, not a candidate to be a model category. We shall use it primarily as a place to represent maps in $\Ho G\Para B$ that cannot easily be represented in $G\Para B$.

We end this section with some results about lax cofibrations.

\begin{definition}
A lax map $i\colon A\laxto X$ is a {\em lax cofibration} if it satisfies the lax homotopy extension property, meaning that we can always fill in the dashed arrow in the following commutative diagram of lax maps:
\[
 \xymatrix{
  A \ar[r]^-{i_0} \ar@{=>}[d]_i & A\smsh_B I_+ \ar@{=>}[d]\ar@/^/@{=>}[ddr] \\
  X \ar[r] \ar@/_/@{=>}[drr] & X\smsh_B I_+ \ar@{==>}[dr] \\
   && Y
 }
\]
\end{definition}

We shall show that the lax cofibrations are precisely the strict maps that are cofibrations under $B$.

\begin{lemma}
If $(i,\lambda)\colon A\laxto X$ is a lax cofibration, then
$\lambda$ is a constant homotopy, so $i$ is a map over $B$.
\end{lemma}

\begin{proof}
Suppose that $(i,\lambda)\colon (A,p,\sigma)\laxto (X,q,\tau)$ is a lax cofibration. Construct $\tilde Mi$ as the following variant of the mapping cylinder: As a $G$-space, let
$\tilde Mi = X\union_i (A\smsh_B I_+)$, the pushout along the inclusion $i_0\colon A \to A\smsh_B I_+$. Give $\tilde Mi$ the evident section, but define $r\colon \tilde Mi\to B$ as follows:
\begin{align*}
 r(x) &= q(x)
 			\qquad\text{if $x\in X$} \\
 r(a,s) &= \lambda(a, l_\lambda(a)\cdot s)
 			\qquad\text{if $(a,s)\in A\smsh_B I_+$}.
\end{align*}
Define a lax map $(\tilde\imath,\tilde\lambda)\colon A\smsh_B I_+\laxto \tilde Mi$ by letting $\tilde\imath$ be the usual map and defining $\tilde\lambda$ as follows:
\begin{align*}
 \tilde\lambda(a,s,t) &= \lambda(a, l_\lambda(a)\cdot s + t) \\
 l_{\tilde\lambda}(a,s) &= l_\lambda(a)\cdot(1 - s).
\end{align*}
Note that the restriction of $(\tilde\imath,\tilde\lambda)$ to $A\smsh 1$ is a strict map over $B$.
Now, from the assumption that $A\laxto X$ is a lax cofibration, we can fill in a lax map $f$ in the following diagram:
\[
 \xymatrix{
  A \ar[r] \ar@{=>}[d]_{(i,\lambda)}
    & A\smsh_B I_+ \ar@{=>}[d]\ar@/^/@{=>}[ddr]^{(\tilde\imath,\tilde\lambda)} \\
  X \ar[r] \ar@/_/[drr] & X\smsh_B I_+ \ar@{==>}[dr]^f \\
   && \tilde Mi
 }
\]
Consider the rightmost triangle restricted to $A\smsh 1$. From the definition of composition and the fact that the restriction of $\tilde\imath$ is strict, it follows that the map $A\smsh 1\laxto X\smsh 1$ must be strict, from which it follows that the original map $A\laxto X$ must be strict.
\end{proof}

\begin{theorem}\label{thm:laxcofib}
A lax map $i\colon A\laxto X$ is a lax cofibration if and only $i$ is a strict map and a cofibration under $B$.
\end{theorem}

\begin{proof}
Suppose that $i\colon A\laxto X$ is a lax cofibration.
The lemma above shows that $i$ is a strict map. We now form $Mi=X\union_A (A\smsh_B I_+)$, the usual mapping cylinder of $i$ in $G\Para B$, which is a subspace of $X\smsh_B I_+$. From the definition of lax cofibration we get a lax retraction $r\colon X\smsh_B I_+\laxto Mi$. Forgetting the maps to $B$, this retraction under $B$ shows that $i$ is a cofibration under $B$.

Conversely, let $i\colon A\to X$ be a strict map and a cofibration under $B$.
Let $p\colon X\to B$ be the projection.
Let $q\colon Mi\to B$ be the mapping cylinder of $i$ and let
$r\colon X\smsh_B I_+ \to Mi$ be a retraction under $B$.
We want to make $r$ into a lax retraction.
Define $h\colon (X\smsh_B I_+)\times [0,\infty]\to B$ by
\[
 h(x,s,t) =
  \begin{cases}
   qr(x,s-t) & \text{if $t \leq s$} \\
   q(x) & \text{if $t\geq s$}
  \end{cases}
\]
The map $h$ is clearly continuous, but is not quite what we need.
Let $u\colon X\to I$ be a $G$-invariant map such that $u^{-1}(0) = A$,
e.g., $u(x) = \sup\{t - p_2 r(x,t) \mid t\in I\}$, where $p_2$ is projection to $I$.
Define $\lambda\colon (X\smsh_B I_+)\times [0,\infty]\to B$ by
\[
 \lambda(x,s,t) =
  \begin{cases}
   h(x,s,t/u(x)) & \text{if $u(x) > 0$} \\
   q(x) & \text{if $u(x) = 0$}
  \end{cases}
\]
and let
\[
 l_\lambda(x,s) = su(x).
\]
Assuming that $\lambda$ is continuous, the pair $(r,\lambda)$ defines a lax retraction, from which it follows that $i$ is a lax cofibration. So, we need to check that $\lambda$ is continuous.

Continuity at points $(x,s,t)$ such that $u(x) \neq 0$ is clear, so we need to check continuity at points $(a,s,t)$ with $a\in A$. For such a point we have $\lambda(a,s,t) = p(a)$. Let $V$ be any neighborhood of $p(a)$ in $B$. Because $I\times [0,\infty]$ is compact, we can find a neighborhood $U$ of $a$ such that $U\times I\times [0,\infty] \subset h^{-1}(V)$. It follows that
$U\times I\times [0,\infty] \subset \lambda^{-1}(V)$, showing that $\lambda$ is continuous.
\end{proof}

We are now entitled to the following result, in which $C_B$ denotes the mapping cone over $B$ and $/_B$ indicates the quotient over $B$.

\begin{proposition}
Suppose that $i\colon A\to X$ is a strict map that is a cofibration under $B$. Then the natural map $c\colon C_B i\to X/_B A$ is a lax homotopy equivalence.
\end{proposition}

\begin{proof}
Because $i$ is a lax cofibration,
the obvious map
$k\colon X\union_A (A\smsh_B I_+) \to C_B i$ over $B$ extends to a
lax map $(h,\lambda)\colon X\smsh_B I_+\laxto C_B i$.
Note that $\lambda(a,t)$ is a 0-length path for every $a\in A$ and $t\in I$.
If $(h_1,\lambda_1)$ is the restriction to
$X\smsh 1$, note that $h_1$ takes $A$ to the section, hence induces a lax map
$\bar h_1\colon X/_B A \laxto C_B i$, which we claim is a lax homotopy inverse to $c$.

The composite $c\circ h\colon X\smsh_B I_+\laxto X/_B A$ factors through
$k\colon X/_B A\smsh_B I_+\laxto X/_B A$, which clearly gives a lax homotopy from the identity to $c\circ \bar h_1$. (This uses the fact that $\lambda$ is constant on $A$.)

On the other hand, it is easy to extend $h$ to a lax homotopy
$h'\colon C_B i \smsh_B I_+\to C_B i$ from the identity to $\bar h_1\circ c$.
\end{proof}

We also have the following results on the homotopy invariance
of pushouts and colimits.

\begin{proposition}
Suppose given the following diagram in which the horizontal maps
are strict, the vertical maps are lax homotopy equivalences, and
$i$ and $i'$ are lax cofibrations.
\[
 \xymatrix{
  X \ar@{=>}[d] & A \ar[l]_i \ar[r] \ar@{=>}[d] & Y \ar@{=>}[d] \\
  X' & A' \ar[l]^{i'} \ar[r] & Y'
 }
\]
Then the induced map of pushouts
$Y\union_A X \laxto Y'\union_{A'} X'$ is a lax homotopy equivalence.
\end{proposition}

\begin{proof}
The proof of the classical result, as in
\cite[7.5.7]{Br:topologyBook}, depends only on formal properties of
the homotopy extension property and homotopies, and can
be adapted to the context of $G\Lax B$.
The precursor results to this result also appear in the more recent
\cite[\S6.5]{May:concise}; Peter tells us that the result itself
will appear in the second edition.
\end{proof}

\begin{proposition}
Suppose given the following diagram in which the horizontal maps
are lax cofibrations and the vertical maps are lax homotopy equivalences.
\[
 \xymatrix{
  X_0 \ar[r] \ar@{=>}[d] & X_1 \ar[r] \ar@{=>}[d] & X_2 \ar[r] \ar@{=>}[d] & \cdots \\
  Y_0 \ar[r] & Y_1 \ar[r] & Y_2 \ar[r] & \cdots
 }
\]
Then the induced map $\colim_n X_n \laxto \colim_n Y_n$ is a
lax homotopy equivalence.
\end{proposition}

\begin{proof}
This follows from
\cite[7.4.1]{Br:topologyBook}, specifically the special case
considered in Step~1 of its proof.
(It would also follow from the last proposition
of \cite[\S6.5]{May:concise} using details
about the forms of the homotopies used in the proof.)
Again, the proof generalizes to our context.
\end{proof}


\section{The fundamental groupoid}\label{sec:fundgrpd}

The following definition was given in \cite{CMW:orientation}.

\begin{definition}
Let $B$ be a $G$-space. The {\em equivariant fundamental groupoid}
$\Pi_G B$ of $B$ is the category whose objects are the $G$-maps
$x\colon G/H\to B$ and whose morphisms $x\to y$, $y\colon G/K\to B$,
are the pairs $(\omega,\alpha)$, where
$\alpha\colon G/H\to G/K$ is a $G$-map and $\omega$ is an equivalence
class of paths $x\to y\circ \alpha$ in $B^H$.
As usual, two paths are equivalent if they are homotopic rel endpoints.
Composition is induced by composition of maps of orbits and the usual
composition of path classes.
\end{definition}

Recall that $\orb G$ is the category of $G$-orbits and $G$-maps.
We have a functor
$\pi\colon \Pi_G B\to \orb G$ defined by
$\pi(x\colon G/H\to B) = G/H$ and $\pi(\omega,\alpha) = \alpha$.
We topologize the mapping sets in $\Pi_G B$ as in
\cite[3.1]{CMW:orientation} so, in particular, $\pi$ is continuous.
(The details of the topology are not important here.)
For each subgroup $H$, the subcategory
$\pi^{-1}(G/H)$ of objects mapping to $G/H$ and morphisms mapping to
the identity is a copy of $\Pi(B^H)$, the nonequivariant fundamental
groupoid of $B^H$ (with discrete topology).
$\Pi_G B$ itself is not a groupoid in the usual
sense, but a ``cat\'egories fibr\'ees en groupoides'' \cite{Gr:groupefond}
or a ``bundle of groupoids'' \cite[5.1]{CMW:orientation} over $\orb G$
(which we usually shorten to {\em groupoid over $\orb G$}).

The {\em homotopy fundamental groupoid} of $B$ is the category
$h\Pi_G B$ obtained by replacing each morphism space with its
set of path components. Equivalently, we identify any two morphisms
that are homotopic in the following sense:
Maps $(\omega,\alpha)$, $(\xi, \beta)\colon x\to y$, where
$x\colon G/H\to B$ and $y\colon G/K\to B$, are homotopic
if there is a homotopy
$j\colon G/H\times I\to G/K$ from $\alpha$ to $\beta$
and a homotopy $k\colon G/H\times I\times I\to B$ from
a representative of $\omega$ to a representative of $\xi$
such that
$k(a,0,t) = x(a)$ and $k(a,1,t) = yj(a,t)$ for all
$a\in G/H$ and $t\in I$.
Thus, $h\Pi_G B$ coincides with tom~Dieck's
``discrete fundamental groupoid'' \cite[10.9]{tD:transfgroups}.
When $G$ is finite, $h\Pi_G B = \Pi_G B$.

The main reason we mention the homotopy fundamental groupoid now is to make the following connection: $h\Pi_G B$ is isomorphic to the full subcategory of
$\Ho G\Para B$ determined by the objects $x\colon G/H\to B$.
This follows from the description of $\Ho G\Para B$ in terms of lax maps given in the preceding section.

\begin{definition}\label{def:ortho}
Let $\ortho G{}$ be the category whose objects are the
orthogonal $G$-vector bundles over orbits of
$G$ and whose morphisms are the equivalence classes of orthogonal
$G$-vector bundle maps between them.
Here, two maps are equivalent if they are $G$-bundle homotopic,
with the homotopy inducing the constant homotopy on base spaces.
Let $\pi\colon \ortho G{} \to \orb G$ be the functor taking the
bundle $p\colon E\to G/H$ to
its base space $G/H$, and taking a bundle
map to the underlying map of base spaces. Let $\ortho Gn$ be the full
subcategory of $\ortho G{}$ consisting of the $n$-dimensional bundles.
 \end{definition}

We topologize the mapping sets in $\ortho G{}$ as in
\cite[4.1]{CMW:orientation}. The map
\[
 \pi\colon \ortho G{}(p,q)\to \orb G(\pi(p),\pi(q))
\]
is then a bundle with discrete fibers.

The categories $\ortho G{}$ and $\ortho Gn$ are not small, but have small skeleta
which we call
$\orthosmall G{}$ and $\orthosmall Gn$,
obtained by choosing one representative in each equivalence class.
An explicit choice is given in \cite[2.2]{CMW:orientation}.
The following is \cite[2.7]{CMW:orientation}.

\begin{proposition} \label{prop:bundlereps}
A $G$-vector bundle $p\colon E\to B$ determines by pullbacks a map
\[
 p^*\colon \Pi_G B \to \ortho G{}
\]
 over $\orb G$. A $G$-vector bundle map
$(\tilde f,f)\colon p\to q$, with $\tilde f\colon E\to E'$ the map of total
spaces and $f\colon B\to B'$ the map of base spaces, determines a natural
isomorphism $\tilde f_*\colon p^*\to q^*f_*$ over the identity.
If $(\tilde h, h)\colon (\tilde f_0,f_0)\to (\tilde f_1,f_1)$ is a homotopy of
$G$-vector bundle maps, then
$(\tilde f_1)_* = (p')^*h_* \circ (\tilde f_0)_*$.
\qed
\end{proposition}

The map $p^*$ is an example of a {\em representation} of
the fundamental groupoid $\Pi_G B$:

\begin{definition}
An {\em $n$-dimensional orthogonal representation} of $\Pi_G B$ is a
continuous functor $\gamma\colon \Pi_G B\to \orthosmall Gn$
over $\orb G$.
If $\gamma$ is an $n$-dimensional orthogonal representation of
$\Pi_G B$ and
$\gamma'$ is an $n$-dimensional orthogonal representation of
$\Pi_G B'$, then a {\em map} from $\gamma$ to $\gamma'$ is a pair
$(f,\eta)$ where $f\colon B\to B'$ is a $G$-map and
$\eta\colon \gamma\to \gamma'\circ f_*$ is a natural isomorphism
over the identity of $\orb G$.
\end{definition}

For example, if $V$ is an orthogonal $H$-representation, then
we write $G\times_H \Vrep$ for the orthogonal representation of
$\Pi_G(G/H)$ associated to the bundle
$G\times_H V\to G/H$.
If $B$ is a $G$-space and $V$ is $G$-representation, we also write
$\Vrep$ for the pullback to $\Pi_G B$ of the representation
$\Vrep$ of $\Pi_G(*) = \orb G$.

In addition to $\orthosmall G{}$, there are several other small categories
that are useful to use as the targets of representations of $\Pi_G B$.
They fit into the following commutative diagram
(see \cite[\S\S 18 \& 19]{CMW:orientation} for details).
\[
 \xymatrix{
  {\orthosmall Gn} \ar[r] \ar[d] & {\PLsmall Gn} \ar[r] \ar[d] &
      {\Topsmall Gn} \ar[r] \ar[d] & {\spheresmall Gn} \ar[d] \\
  {s\orthosmall Gn} \ar[r] \ar[d] & {s\PLsmall Gn} \ar[r] \ar[d] &
      {s\Topsmall Gn} \ar[r] \ar[d] & {s\spheresmall Gn} \ar[d] \\
  {v\orthosmall Gn} \ar[r] & {v\PLsmall Gn} \ar[r]  &
      {v\Topsmall Gn} \ar[r]  & {v\spheresmall Gn} 
 }
\]
To form the top row,
we take the categories of $G$-vector bundles over orbits, with morphisms
being, respectively, the fiber-homotopy classes of $G$-vector bundle maps, PL $G$-bundle
maps, topological $G$-bundle maps, and spherical $G$-fibration maps.
In the case of $\spheresmall Gn$, what we mean is that the objects are bundles, but a map
from $G\times_H V$ to $G\times_K W$ is a based
fiber-homotopy class of homotopy equivalences
$G\times_H S^V\to G\times_K S^W$.
The middle row is the stabilization of the top row and
the bottom row gives virtual bundles.
Maps into these categories give us, respectively orthogonal, PL, topological, and
spherical representations, and their stable and virtual analogues.

In particular, spherical representations arise naturally from spherical fibrations.
Here, by a spherical fibration we mean the following
(see also \cite[\S 23]{CMW:orientation}).

\begin{definition}
A {\em spherical $G$-fibration} is a sectioned $G$-fibration $p\colon E\to B$ such that
each fiber $p^{-1}(b)$ is based $G_b$-homotopy equivalent to $S^V$ for some
orthogonal $G_b$-representation $V$.
\end{definition}

The various categories over $\orb G$ that we've discussed, including
$\Pi_G B$ and the twelve categories in the diagram above, are all
groupoids over $\orb G$.
Associated to a groupoid $\R$ over $\orb G$, there is a classifying
$G$-space $B\R$, defined in \cite[\S 20]{CMW:orientation}.
The following classification result is proved for finite groups as
\cite[24.1]{CMW:orientation} and extended to compact Lie groups
in \cite{Co:monoids}.

\begin{theorem}\label{thm:classification}
The $G$-space $B\R$ classifies representations of $\Pi_G B$ in $\R$.
That is, for $G$-CW complexes $B$, $[B, B\R]_G$ is in natural bijective
correspondence with the set of isomorphism classes of representations
$R\colon \Pi_G B\to \R$.
\end{theorem}


\section{Parametrized spectra}\label{sec:paraspec}

May and Sigurdsson give a very careful treatment of equivariant parametrized spectra
in \cite{MaySig:parametrized}.
They concentrate on orthogonal $G$-spectra over $B$, which give a model category
with good formal properties.
They also discuss $G$-prespectra over $B$, which
arise more naturally as representing objects for cohomology theories.
As in the first chapter,
it is therefore more convenient for us to use the latter, keeping in mind that the two kinds
of spectra give equivalent stable categories.
As explained in \cite{MaySig:parametrized}, we need to restrict our parametrizing
spaces to be compactly generated and have the homotopy types of $G$-CW complexes so that all of the
functors and adjunctions we need pass to homotopy categories.
We make this restriction from this point on.

Let $\RO(G)$ be the collection of all finite-dimensional representations of $G$. (May and Sigurdsson consider collections closed under direct sum but not necessarily containing all representations. For our purposes we do not need the greater generality.) The following definition is
\cite[11.2.16 \& 12.3.6]{MaySig:parametrized}.

\begin{definition}
\begin{enumerate}
\item[]
\item
A {\em $G$-prespectrum $E$ over $B$} consists of a collection of
ex-$G$-spaces $E(V)\to B$, one
for each $V\in\RO(G)$, together with structure
maps given by $G$-maps
\[
 \sigma\colon \susp_B^W E(V) \to E(V\oplus W)
\]
over $B$.
The $\sigma$ are required to be unital and transitive in the usual way.
\item
An {\em $\Omega$-$G$-prespectrum} is a $G$-prespectrum in which each
ex-$G$-space $E(V)$ is $qf$-fibrant and the adjoint structure maps
\[
 \tilde\sigma \colon E(V) \to \Loop_B^W E(V\oplus W)
\]
are $G$-weak equivalences.
\end{enumerate}
We let $G\PreSpec{} B$ denote the category of $G$-prespectra over $B$ and $G$-maps between them, i.e., levelwise $G$-maps over $B$ that respect the structure maps.
\end{definition}

May and Sigurdsson
define a stable model structure on $G\PreSpec{} B$
\cite[\S 12.3]{MaySig:parametrized}.
The stable equivalences
({\em $s$-equivalences} for short) are the fiberwise stable equivalences,
where the fibers are those of levelwise fibrant approximations.
The fibrant objects are exactly the $\Omega$-$G$-prespectra.

We write $[E,F]_{G,B}$ for the group of stable $G$-maps
between two $G$-prespectra over $B$,
that is, $\Ho G\PreSpec{} B(E,F)$.
By \cite[12.4.5]{MaySig:parametrized}, $[E,F]_{G,B}$ is stable in the sense that, if
$V$ is any representation of $G$ then there is an isomorphism
\[
 [E,F]_{G,B} \iso [\susp_B^V E, \susp_B^V F]_{G,B}.
\]
It is useful to note that these groups are stable in a stronger sense: they are
stable under suspension by any spherical fibration over $B$.
The following is a direct consequence of \cite[15.1.5]{MaySig:parametrized}.

\begin{proposition}\label{prop:stability}
Let $\xi$ be a spherical $G$-fibration over
a $G$-CW complex $B$ and let
$\susp^\xi$ denote the fiberwise smash product with the total space of $\xi$.
Then, if $E$ and $F$ are any two $G$-prespectra over $B$ we have a natural isomorphism
\[
 [E,F]_{G,B} \iso [\susp^\xi E, \susp^\xi F]_{G,B}.
\]
\qed
\end{proposition}

The following functors and adjunctions are discussed in detail in
\cite{MaySig:parametrized}.
Much of that book is aimed at showing that these and other relationships descend to homotopy categories.
We have the functor $\susp_B^\infty\colon G\Para B \to G\PreSpec{} B$ taking an ex-$G$-space $X$ over $B$ to the prespectrum with
$(\susp_B^\infty X)(V) = \Susp_B^V X$.
If $X$ and $Y$ are ex-$G$-spaces over $B$, we write
\[
 \{X,Y\}_{G,B} = [\Susp_B^\infty X, \Susp_B^\infty Y]_{G,B}
\]
for the group of stable maps from $X$ to $Y$.
More generally, if $W$ is a representation of $G$, we have the
{\em shift desuspension} functor
$\susp^\infty_W\colon G\Para B \to G\PreSpec{} B$ defined by
\[
 \susp^\infty_W X(V) =
  \begin{cases}
   \susp_B^{V-W}X & \text{when $W\subset V$} \\
   B & \text{when $W\not\subset V$.}
  \end{cases}
\]
This functor is left adjoint to the evaluation functor $\Omega^\infty_W$ given by evaluation at $W$.
The adjunction descends to the homotopy categories, by
\cite[12.6.2]{MaySig:parametrized}. We insert here the standard warning:  Functors used on homotopy categories are the derived functors, obtained by first taking a cofibrant or fibrant approximation as appropriate. In particular, on the homotopy category, $\Omega^\infty_W$ does not return the $W$th space of a prespectrum but the $W$th space of a stably equivalent $\Omega$-$G$-prespectrum.

Suppose that $\alpha\colon A\to B$ is a $G$-map. Then
there are functors
\[
 \alpha^*\colon G\Para {B} \to G\Para A
\]
and
\[
 \alpha^*\colon G\PreSpec{}{B} \to G\PreSpec{} A
\]
given by taking pullbacks.
These functors have left adjoints
\[
 \alpha_!\colon G\Para A \to G\Para {B}
\]
and
\[
 \alpha_!\colon G\PreSpec{} A\to G\PreSpec{} {B}
\]
given by composition with $\alpha$ and identification of base points.
Precisely, on $G\Para A$, $\alpha_!(X,p,\sigma)$ is given by taking the pushout in the
top square below: 
\[
 \xymatrix{
  A \ar[d]_{\sigma} \ar[r]^{\alpha} & B \ar@{-->}[d] \\
  X \ar[d]_p \ar@{-->}[r] & \alpha_! X \ar@{-->}[d] \\
  A \ar[r]_{\alpha} & B
 }
\]  
$(\alpha_!,\alpha^*)$ is a Quillen adjoint pair.
Moreover, $\alpha^*$ and $\alpha_!$ both commute with suspension, in the sense that
\[
 \susp_A^\infty \alpha^* X \iso \alpha^*\susp_{B}^\infty X
\]
and
\[
 \susp_{B}^\infty \alpha_! X \iso \alpha_!\susp_A^\infty X.
\]

Another useful property is the natural homeomorphism
 \[
  X\smsh_{B} \alpha_!Y \homeo \alpha_!(\alpha^*X\smsh_A Y)
 \]
for ex-$G$-spaces $X$ over $B$ and $Y$ over $A$.
(The smash product used is the usual
fiberwise smash product.)
This homeomorphism descends to an isomorphism in the homotopy category.
Similarly, if $E$ is a prespectrum over $B$, then
 \[
  E\smsh_{B} \alpha_!Y \iso \alpha_!(\alpha^*E\smsh_A Y).
 \]
We shall apply this most often to the case of
$\rho\colon B\to *$, the projection to a point.
(We shall use $\rho$ generically for any projection to a point.)
$\rho$ induces functors
$\rho^*\colon G\PreSpec{}{} \to G\PreSpec{} B$ and
$\rho_!\colon G\PreSpec{} B\to G\PreSpec{}{}$. Notice that
 \[
 \rho^*\susp^\infty Y = \susp^\infty_B \rho^* Y = \susp^\infty_B (B\times Y)
 \]
for a based $G$-space $Y$, and
 \[
 \rho_!\susp^\infty_B X = \susp^\infty \rho_! X = \susp^\infty (X/\sigma(B))
 \]
for an ex-$G$-space $(X,p,\sigma)$ over $B$.


\section{Lax maps of prespectra}

We would like to define a category of lax maps of prespectra,
similar to the category $G\Lax B$. May and Sigurdsson
extended $L$ to a functor on prespectra in
\cite[\S13.3]{MaySig:parametrized}, but there is a mistake in their
definition. Moreover, they were not trying to make 
their $L$ a monad, so, for example, took no care to allow for an
associative composition.
We begin by giving a definition that does give us a monad,
so that we can again define the category of lax maps to be the
Kleisli category of $L$.

The first difficulty with extending $L$ to a functor on prespectra is that, if $K$ is an ex-space, there is no canonical map $LK\smsh_B S^V\to L(K\smsh_B S^V)$. So, we need to define an ad hoc map, which we call $\beta^V$.
(May and Sigurdsson write $\beta_V$ for their similar map.) In defining $\beta^V$, we shall think of $S^V = D(V)/S(V)$, so that a $v\in S^V$ will have $\|v\|\leq 1$. Here we differ from May and Sigurdsson; the difference is unimportant, but we find this choice easier to work with here.
We first choose any homeomorphism $\phi\colon [0,\infty]\to [1/2,1]$ with
$\phi(0) = 1$ and $\phi(\infty) = 1/2$.
Now we define $\beta^V$ by the following formula, in which we write
$l$ for $l_\lambda$ to simplify the notation:
\[
 \beta^V((x,\lambda)\smsh v) =
  \begin{cases}
    ( x\smsh v, \lambda )
    & \text{if $\|v\| \leq 1/2$}
  \\[\bigskipamount]
   \left( x\smsh \phi(\phi^{-1}(\|v\|) - l)\frac{v}{\|v\|}, \lambda \right)
    & \text{if $1/2 < \|v\| \leq \phi(l)$}
  \\[\bigskipamount]
   \left(
     \sigma\lambda(l-\phi^{-1}(\|v\|)),
      \lambda|^{l}_{l-\phi^{-1}(\|v\|)} \right)
    & \text{if $\|v\| \geq \phi(l)$}
  \end{cases}
\]
Here, as in \cite{MaySig:parametrized}, we write $\lambda|_a^b$ for the Moore path of length $b-a$ given by 
\[
 (\lambda|_a^b)(t) =
  \begin{cases}
   \lambda(a+t) & \text{if $t \leq b-a$} \\
   \lambda(b) & \text{if $t \geq b-a$}.
  \end{cases}
\]

Now, to extend $L$ to prespectra we need to deal with the fact that 
$\beta^V\circ\beta^W \neq \beta^{V\dirsum W}$, so we do as in 
\cite[\S13.3]{MaySig:parametrized} and 
consider prespectra indexed on a fixed countable cofinal sequence
$\W = \{0 = V_0 \subset V_1\subset\cdots\}$.
Write $W_i = V_i-V_{i-1}$
(note that our indexing of the $W_i$ is slightly different from that
of May and Sigurdsson).
As explained in \cite[13.3.5]{MaySig:parametrized},
the category of prespectra indexed on $\W$ gives us the same
stable homotopy category as usual.
If $X$ is a $G$-prespectrum, we define an indexed prespectrum $LX$
to have spaces $(LX)(V_i) = L(X(V_i))$ and
structure maps
$\sigma_i\colon LX(V_i)\smsh_B S^{W_{i+1}}\to LX(V_{i+1})$ given by the composites
\[
 \xymatrix{
 LX(V_i)\smsh_B S^{W_{i+1}}
  \ar[r]^{\beta^{W_{i+1}}} &
 L(X(V_i)\smsh_B S^{W_{i+1}})
  \ar[r] &
 LX(V_{i+1}).
 }
\]
In general, define the structure map
$\Susp^{V_j-V_i}LX(V_i) \to LX(V_j)$ to be the composite
$\sigma_{j-1}\circ\sigma_{j-2}\circ\cdots\circ\sigma_i$.

\begin{proposition}\label{prop:LSpecMonad}
The levelwise inclusion $\iota\colon X\to LX$ and levelwise
composition $\chi\colon L^2 X\to LX$ make $(L,\iota,\chi)$
a monad on the category of $\W$-indexed $G$-prespectra over $B$.
\end{proposition}

\begin{proof}
The main difficulty is showing that $\iota$ and $\chi$ are actually
maps of prespectra, and the definition of $\beta^V$ was crafted to
make them so. The fact that $\iota$ is a map of prespectra follows from
the fact that the following diagram commutes for any ex-$G$-space $K$:
\[
 \xymatrix@C-3em{
  & K\smsh_B S^V \ar[dl]_{\iota\smsh 1} \ar[dr]^{\iota} \\
  LK\smsh_B S^V \ar[rr]_-{\beta^V} && L(K\smsh_B S^V)
 }
\]
The fact that $\chi$ is a map of prespectra comes down to the commutativity
of the following diagram:
\[
 \xymatrix{
  L^2K \smsh_B S^V \ar[r]^{(\beta^V)^2} \ar[d]_{\chi\smsh 1}
    & L^2(K\smsh_B S^V) \ar[d]^{\chi} \\
  LK \smsh_B S^V \ar[r]_{\beta^V} & L(K\smsh_B S^V)
 }
\]
The verification that these diagrams commute is by
straightforward, if tedious, calculation using the definition of $\beta^V$.

That $(L,\iota,\chi)$ is a monad now follows directly
from the analogous fact for spaces.
\end{proof}

By analogy with spaces, we now make the following definition.

\begin{definition}
The category of {\em lax maps of $G$-prespectra over $B$}, denoted
$G\LaxPreSpec{}B$, is the Kleisli category of the monad $L$
on the category of $\W$-indexed $G$-prespectra over $B$.
\end{definition}

Explicitly, then, a lax map $X\laxto Y$ is a strict map
$X\to LY$. The composite of the maps $f\colon X\laxto Y$ and 
$g\colon Y\laxto Z$ is the following composite of strict maps:
\[
 X \xrightarrow{f} LY \xrightarrow{Lg} L^2 Z \xrightarrow{\chi} LZ.
\]
A lax map of prespectra can be thought of as a collection of lax maps
of ex-spaces at each level, compatible under suspension, using $\beta^{W_i}$ to interpret what compatibility means. Composition of lax maps is given levelwise
by composition of lax maps of ex-spaces.

Since $LY$ is a levelwise fibrant approximation but not, in general, a $qf$-fibrant approximation of $Y$, we don't expect to represent all stable maps by lax maps. 
However, lax maps give us useful examples of stable maps that are not represented by strict maps.

There are a number of cases where we will want to know that a functor
on strict maps extends to one on lax maps.
The general context is this: Suppose that $\C$ and $\D$
are categories,
$(S,\eta,\mu)$ is a monad on $\C$, and $(T,\theta,\nu)$ is a monad
on $\D$. Write $\C_S$ for the Kleisli category of $S$ and $\D_T$ for
that of $T$, so $\C_S(X,Y) = \C(X,SY)$ and $\D_T(X,Y) = \D(X,TY)$.
We ask what we need to extend a functor $F\colon\C\to\D$ to a functor $\C_S\to \D_T$.
The answer is the following.

\begin{definition}
A {\em map of monads} $(F,\psi)\colon (\C,S)\to (\D,T)$
consists of a functor $F\colon \C\to \D$ and a natural transformation
$\psi\colon FS\to TF$ such that the following diagrams commute:
\[
 \xymatrix@C-1.5em{
  & F \ar[dl]_{F\eta} \ar[dr]^{\theta F} \\
  FS \ar[rr]_-{\psi} && TF
 }
\]
\[
 \xymatrix{
  FS^2 \ar[r]^{\psi^2} \ar[d]_{F\mu}
   & T^2 F \ar[d]^{\nu F} \\
  FS \ar[r]_\psi & TF
 }
\]
(This is what Street calls a {\em monad opfunctor} in
\cite{Str:MonadsI}. His monad functors have the natural transformation
going in the other direction, and induce maps of Eilenberg-Moore categories.)
\end{definition}

The general theory of \cite{Str:MonadsI} shows that a
map of monads $(F,\psi)\colon (\C,S)\to (\D,T)$ gives
an extension of $F$ to $\bar F\colon \C_S\to \D_T$.
Explicitly, $\bar F(X) = F(X)$ on objects and, if
$f\colon X\to SY$ represents a map in
$\C_S(X,Y)$, then $\bar F(f)$ is the composite
\[
 FX \xrightarrow{Ff} FSY \xrightarrow{\psi} TFY.
\]
It is straightforward to check that this defines a functor
extending $F$.
In fact, though we shall not need it, extensions of $F$
are in one-to-one correspondence with natural transformations
$\psi\colon FS\to TF$ such that $(F,\psi)$ is a map of monads.

We have already seen one example of a map of monads, namely
\[
 (\susp^V_B, \beta^V)\colon (G\Para B,L) \to (G\Para B,L).
\]
The diagrams shown in the proof of Proposition~\ref{prop:LSpecMonad}
can be rewritten as the diagrams showing that $(\susp^V_B,\beta^V)$
is a map of monads:
\[
 \xymatrix@C-1.5em{
  & \susp^V_B \ar[dl]_{\susp^V_B\iota} \ar[dr]^{\iota\susp^V_B} \\
  \susp^V_B L \ar[rr]_-{\beta^V} && L\susp^V_B
 }
\]
\[
 \xymatrix{
  \susp^V_B L^2 \ar[r]^{(\beta^V)^2} \ar[d]_{\susp^V_B \chi}
   & L^2 \susp^V_B \ar[d]^{\chi \susp^V_B} \\
  \susp^V_B L \ar[r]_{\beta^V} & L\susp^V_B
 }
\]
Thus, we get the following.

\begin{proposition}
The suspension functor $\susp^V_B$ extends to a continuous functor we shall
again call $\susp^V_B\colon G\Lax B\to G\Lax B$.
\qed
\end{proposition}

Explicitly,
the suspension of a map $f\colon X\laxto Y$ is given by the composite
\[
 \susp^V_B X \xrightarrow{\susp^V_B f} \susp^V_B LY
   \xrightarrow{\beta^V} L\susp^V_B Y.
\]
Note, however, that $\susp^W_B\susp^V_B \neq \susp^{V\dirsum W}_B$
in general, because $\beta^W\circ\beta^V \neq \beta^{V\dirsum W}$.

We insert here a useful fact about lax cofibrations.

\begin{corollary}\label{cor:suspcofibr}
Suspension preserves lax cofibrations. I.e., if $A\to X$ is
a lax cofibration of ex-$G$-spaces, then so is $\susp^V_B A\to \susp^V_B X$.
\end{corollary}

\begin{proof}
Let $i\colon A\to X$ be a lax cofibration.
Because $\susp^V_B$ is a left adjoint (on $G\Para B$),
it preserves colimits, hence $\susp^V_B Mi \iso M(\susp^V_B i)$,
where $Mi$ denotes the usual (fiberwise) mapping cylinder.
Because $i$ is a lax cofibration, $Mi$ is a retract of $X\smsh_B I_+$
in $G\Lax B$,
hence $M(\susp^V_B i) \iso \susp^V_B Mi$ is a retract of $\susp^V_B X\smsh I_+$
in $G\Lax B$. Therefore $\susp^V_B i$ is a lax cofibration.
\end{proof}

We would also like to extend the functors $\susp^\infty_Z$
to lax maps.

\begin{remark}\label{rem:indexedSuspension}
Note that, because we are using indexed prespectra, the functor
$\susp^\infty_Z$ is now defined by
\[
 \susp^\infty_Z X(V_i) =
  \begin{cases}
   \susp_B^{V_i-Z}X & \text{when $Z\subset V_i$} \\
   B & \text{when $Z\not\subset V_i$.}
  \end{cases}
\]
Let $\bar Z$ represent the smallest $V_i$ such that $Z\subset V_i$.
Then $\susp^\infty_Z K = \susp^\infty_{\bar Z}\susp_B^{\bar Z-Z}K$.
From this it is clear that, in the context of indexed prespectra,
$\susp^\infty_Z$ still has a right adjoint, but that right adjoint
is $\Loop_B^{\bar Z-Z} \Loop^\infty_{\bar Z}$.
\end{remark}

\begin{definition}
The natural transformation
$\beta_Z^\infty \colon \susp^\infty_Z L \to L \susp^\infty_Z$
is defined as follows. Let $k$ be the least integer such that
$Z\subset V_k$. For $i\geq k$ we let
\[
 \beta_Z^\infty
 = \beta^{W_i}\circ \beta^{W_{i-1}}\circ\dots\circ\beta^{W_{k+1}}
   \circ \beta^{V_k-Z}
 \colon \susp_B^{V_i-Z}L \to L\susp_B^{V_i-Z},
\]
using the decomposition
\[
 V_i - Z = (V_k - Z)\dirsum W_{k+1}\dirsum \cdots \dirsum W_i.
\]
For $i<k$, $\beta_Z^\infty$ is the trivial map between trivial ex-spaces.
\end{definition}

\begin{proposition}
The natural transformation $\beta_Z^\infty$ defines a map of monads
$(\susp^\infty_Z, \beta_Z^\infty)\colon
 (G\Para B, L) \to (G\PreSpec{}B, L)$.
Hence, $\susp^\infty_Z$ extends to a continuous functor
$\susp^\infty_Z\colon G\Lax B \to G\LaxPreSpec{}B$.
\end{proposition}

\begin{proof}
We first need to verify that $\beta_Z^\infty$ is a map of prespectra.
For this we need the following to commute:
\[
 \xymatrix{
  \susp_B^{W_{i+1}}\susp_B^{V_i-Z} LK
    \ar[r]^{\susp\beta_Z^\infty}
    \ar@{=}[d]
   & \susp_B^{W_{i+1}} L\susp_B^{V_i-Z} K \ar[d]^{\beta^{W_{i+1}}} \\
  \susp_B^{V_{i+1}-Z} LK
    \ar[r]_{\beta_Z^\infty}
   & L\susp_B^{V_{i+1}-Z} K
 }
\]
We defined $\beta_Z^\infty$ precisely to make this diagram commute.

Now we need to know that $(\susp^\infty_Z,\beta_Z^\infty)$ is a map
of monads. That $\beta_Z^\infty$ preserves units and multiplication
follows by iterating the corresponding facts about the $\beta^V$s.
\end{proof}

\begin{corollary}\label{cor:suspPreservesHomotopy}
If $f\colon X\laxto Y$ is a lax homotopy equivalence of ex-$G$-spaces,
then $\susp^\infty_Z f\colon \susp^\infty_Z X\laxto \susp^\infty_Z Y$
is a a lax homotopy equivalence of $G$-prespectra over $B$.
\qed
\end{corollary}

The proof of the following is almost identical to that
of Corollary~\ref{cor:suspcofibr}.

\begin{corollary}\label{cor:suspcofibrspec}
If $A\to X$ is
a lax cofibration of ex-$G$-spaces, then 
$\susp^\infty_Z A\to \susp^\infty_Z X$ is a lax cofibration
of $G$-prespectra over $B$.
\qed
\end{corollary}

We will use this result, together with the fact that pushouts
of lax cofibrations are lax cofibrations, to construct
lax cofibrations of $G$-prespectra over $B$.
The following homotopy invariance results are proved the same way
as for the corresponding results for ex-$G$-spaces.

\begin{proposition}\label{prop:hominvpushouts}
Suppose given the following diagram
of $G$-prespectra, in which the horizontal maps
are strict, the vertical maps are lax homotopy equivalences, and
$i$ and $i'$ are lax cofibrations.
\[
 \xymatrix{
  X \ar@{=>}[d] & A \ar[l]_i \ar[r] \ar@{=>}[d] & Y \ar@{=>}[d] \\
  X' & A' \ar[l]^{i'} \ar[r] & Y'
 }
\]
Then the induced map of pushouts
$Y\union_A X \laxto Y'\union_{A'} X'$ is a lax homotopy equivalence.
\qed
\end{proposition}

\begin{proposition}\label{prop:hominvcolimits}
Suppose given the following diagram of $G$-prespectra over $B$,
in which the horizontal maps are strict and
are lax cofibrations, and the vertical maps are lax homotopy equivalences.
\[
 \xymatrix{
  X_0 \ar[r] \ar@{=>}[d] & X_1 \ar[r] \ar@{=>}[d] & X_2 \ar[r] \ar@{=>}[d] & \cdots \\
  Y_0 \ar[r] & Y_1 \ar[r] & Y_2 \ar[r] & \cdots
 }
\]
Then the induced map $\colim_n X_n \laxto \colim_n Y_n$ is a
lax homotopy equivalence.
\qed
\end{proposition}

We mentioned earlier, and it's easy to show, that
$\iota\colon X\to LX$ is a lax homotopy equivalence if $X$ is
an ex-$G$-space. The same is true if $X$ is a $G$-prespectrum over $B$,
but the proof is not as straightforward.
(Note that it's easy to show that $\iota$ is a stable equivalence.)

\begin{proposition}
If $X$ is a $G$-prespectrum over $B$, then $\iota\colon X\to LX$ is
a lax homotopy equivalence.
Its inverse is the lax map $\kappa\colon LX\laxto X$ represented
by the identity map $1\colon LX\to LX$.
\end{proposition}

\begin{proof}
The composite $\kappa\iota\colon X\laxto X$ is the map
represented by $\iota\colon X\to LX$, i.e., the identity map
in $G\LaxPreSpec{}B$.

The composite $\iota\kappa\colon LX\laxto LX$ is represented
by the map $L\iota\colon LX\to L^2X$. This is not the identity map,
but we will show that it is (fiberwise) homotopic to $\iota L$,
which does represent the identity.
To this end, for any ex-$G$-space $K$ we define a homotopy
$H\colon \iota L\to L\iota$,
\[
 H\colon LK\smsh_B [0,\infty]_+ \to L^2K,
\]
natural in $K$.
Here, we use $[0,\infty]$ rather than the homeomorphic $[0,1]$
to make the formulas defining $H$ simpler:
\[
 H(x,\lambda,t) =
  \begin{cases}
    (x, \lambda|_0^{l-t}, \lambda|_{l-t}^l )
       & \text{if $t\leq l$}
    \\[\bigskipamount]
    (x, \lambda|_0^0, \lambda|_0^l)
       & \text{if $t\geq l$.}
  \end{cases}
\]
In order for $H$ to define a map of prespectra $LX\to L^2X$,
we need the following diagram to commute
(where we have in mind $K = X(V_i)$ and $V = W_{i+1}$):
\[
 \xymatrix{
  {\susp^V_B LK\smsh_B I_+} \ar[r]^-{\susp^V H} \ar[d]_{\beta^V}
   & {\susp^V_B L^2 K} \ar[d]^{(\beta^V)^2} \\
  L\susp^V_B K\smsh_B I_+ \ar[r]_-H
   & L^2 \susp^V_B K
 }
\]
To check that this diagram commutes is now a computation based on
the formulas for $\beta^V$ and $H$, which we will leave to the reader.
\end{proof}

We can now state the prespectrum version of Lemma~\ref{lem:laximpliesweak};
the proof is exactly the same.

\begin{lemma}\label{lem:speclaximpliesweak}
\begin{enumerate}\item[]
\item
If $f\colon X\to Y$ is a strict map that is also a lax homotopy
equivalence, then it is a stable equivalence.
\item
If $f\colon X\laxto Y$ is a lax homotopy equivalence,
then its representing map $\hat f\colon X\to LY$ is a 
lax homotopy equivalence, hence a stable equivalence.
\qed
\end{enumerate}
\end{lemma}

\begin{corollary}
$\beta^\infty_Z\colon \susp^\infty_Z LK\to L\susp^\infty_Z K$
is a lax homotopy equivalence.
\end{corollary}

\begin{proof}
Write $\kappa\colon LK\laxto K$ again for the lax map of
ex-$G$-spaces represented by the identity $1\colon LK\to LK$.
We know that $\kappa$ is a lax homotopy equivalence (inverse
to $\iota\colon K\laxto LK$), hence
$\susp^\infty_Z\kappa\colon \susp^\infty_Z LK\laxto \susp^\infty_Z K$ is
a lax homotopy equivalence. By the preceding lemma,
$\widehat{\susp^\infty_Z\kappa}$ is also a lax homotopy equivalence,
but it is easy to see from the definitions that
$\widehat{\susp^\infty_Z\kappa} = \beta^\infty_Z$.
\end{proof}

Finally, some results about computing the stable homotopy groups used to define stable equivalence of $G$-prespectra.
Recall from \cite[12.3.4]{MaySig:parametrized} that the homotopy groups of a $G$-prespectrum $X$ are the stable homotopy groups of the fibers of a levelwise fibrant approximation of $X$, for which we can use $LX$. So, for each $b\in B$ and $H\subset G_b$, we have
\begin{align*}
 \pi_q^H(X_b)
   &\iso \colim_i \pi_q^H(\Omega^{V_i} LX(V_i)_b) \\
   &\iso \colim_i \pi H\Para *(S^{q+V_i}, LX(V_i)_b) \\
   &\iso \colim_i \pi H\Para B (b_! S^{q+V_i}, LX(V_i)) \\
   &\iso \colim_i \pi H\Lax B(b_! S^{q+V_i}, X(V_i)).
\end{align*}
Thus, we can represent an element of $\pi_q^H(X_b)$ as a lax $H$-map $(f,\lambda)\colon b_! S^{q+V_i}\to X(V_i)$ for some $i$. The equivalence relation is then the one generated by lax homotopy and suspension of lax maps.

We also have the following result that simplifies computing the stable homotopy groups of an ex-$G$-space:
If $K$ is an ex-$G$-space, $b\in B$, and $H\subset G_b$, then
\begin{align*}
 \pi_q^H((\susp^\infty_B K)_b)
  &\iso \pi_q^H((L\susp^\infty_B K)_b) \\
  &\iso \pi_q^H((\susp^\infty_B LK)_b) \\
  &\iso \pi_q^H(\susp^\infty_H LK_b) \\
  &\iso \colim_V \pi H\Para * (S^{q+V}, \susp^V LK_b)
\end{align*}
is just the stable homotopy group of the (non-parametrized) $H$-space $LK_b$.
This also follows from \cite[13.7.4]{MaySig:parametrized}.


\section{The stable fundamental groupoid}\label{sec:stabgrpd}

As mentioned earlier, we need a stable version of the fundamental
groupoid. First, we need some notation.

\begin{definition}\label{def:sphereSpec}
Let $x\colon G/H\to B$ be a $G$-map and let $V$ be a representation of $H$.
We have the ex-$G$-space $G\times_H S^V$ over $G/H$ and we let
\[
 G_+\smsh_H S^{V,x} = x_!(G\times_H S^V)
\]
over $B$.
If $\alpha$ is a virtual representation of $H$, then,
by \cite[11.5.4 \& 13.7.9]{MaySig:parametrized}, corresponding to the
$H$-spectrum $S^\alpha$ there is a $G$-spectrum over $G/H$ we write as
$G\times_H S^\alpha$. We then let
\[
 G_+\smsh_H S^{\alpha,x} = x_!(G\times_H S^\alpha)
\]
over $B$.
\end{definition}

In one important case, we shall be much more specific:

\begin{definition}
Let $\delta$ be a dimension function for $G$. For each orbit $G/H$, choose
an embedding $G/H\subset V$ in a $G$-representation $V$;
note that any two embeddings are isotopic if $V$ is sufficiently large.
Then 
\[
 G/H\times V \iso G\times_H((V-\Lie(G/H))\dirsum\Lie(G/H)),
\]
identifying $\Lie(G/H)$ with the vectors in $V$ tangent to $G/H$ at $eH$.
Using this decomposition allows us to consider $G\times_H \delta(G/H)$
as a subbundle of $G/H\times V$ and to write
$G\times_H(V-\delta(G/H))$ as its orthogonal complement.
We then write
\[
 G\times_H S^{-\delta(G/H)} = G\times_H\susp_V^\infty S^{V-\delta(G/H)},
\]
a $G$-spectrum over $G/H$,
and, if $x\colon G/H\to B$ as in the preceding definition, we let
\[
 G_+\smsh_H S^{-\delta(G/H),x} = x_! (G\times_H S^{-\delta(G/H)}).
\]
\end{definition}

\begin{definition}
\begin{enumerate}\item[]
\item
The {\em stable fundamental groupoid} of a $G$-space $B$ is the category
$\stab\Pi_G B$ whose objects are the $G$-maps $x\colon G/H\to B$ (the same
objects as the unstable fundamental groupoid) and whose maps
from $(G/H,x)$ to $(G/K,y)$ are the stable $G$-maps
$\{(G/H,x)_+,(G/K,y)_+\}_{B,G}$.
\item
More generally, suppose that $\delta$ is an additive dimension function for $G$.
The {\em $\delta$-stable fundamental groupoid} of $B$ is the category
$\stab\Pi_{G,\delta}B = \stab\Pi_\delta B$ with the same objects as above, but whose maps from
$x$ to $y$ are the stable $G$-maps
\[
 [ G_+\smsh_H S^{-\delta(G/H),x}, G_+\smsh_K S^{-\delta(G/K),y} ]_{B,G}.
\]
\item
Even more generally, if $\gamma$ is a virtual orthogonal representation of
$\Pi B$, the {\em $\delta$-$\gamma$-stable fundamental groupoid} of $B$
is the category $\stab\Pi_{\delta,\gamma} B$ with the same objects as above but
whose maps from $x$ to $y$ are the stable $G$-maps
\[
 [ G_+\smsh_H S^{\gamma_0(x)-\delta(G/H),x}, G_+\smsh_K S^{\gamma_0(y)-\delta(G/K),y} ]_{B,G}.
\]
Here,
$G_+\smsh_H S^{\gamma_0(x)-\delta(G/H),x}
 = x_!(G\times_H \susp^{\gamma_0(x)}\susp_V^\infty S^{V-\delta(G/H)})$
for a chosen embedding of $G/H$ in a sufficiently large $V$.
\end{enumerate}
\end{definition}

The purpose of this section is to prove the following result.

\begin{theorem}\label{thm:StableMapsOrbitsOverB}
For any $\delta$ and $\gamma$, $\stab\Pi_\delta B$ is naturally isomorphic to
$\stab\Pi_{\delta,\gamma} B$.
Moreover, if $x\colon G/H\to B$ and $y\colon G/K\to B$ are orbits in $B$, then
$\stab\Pi_\delta B(x,y)$ is isomorphic to a free abelian group on equivalence classes of
diagrams of maps of orbits over $B$ of the form
$[x \from z \laxto y]$ where $z\colon G/L\to B$. We may assume that $L \leq H$ and
$L\leq K^g$ for some $g$, and then we have the conditions that
$\Lie(N_H L/L) - \delta(N_H L/L) = 0$ and
$\delta(N_{K^g}L/L) = 0$.
The diagram $x \from z \laxto y$ is equivalent to $x \from z' \laxto y$, where
$z'\colon G/L'\to B$ if there is a $G$-homeomorphism $z\to z'$ such that
the left triangle below commutes and the right one lax homotopy commutes:
\[
 \xymatrix@R-4ex
 {
  & z \ar[dl] \ar@{=>}[dr] \ar[dd] \\
  x && y \\
  & z' \ar[ul] \ar@{=>}[ur]
 }
\]
\end{theorem}

\begin{remark}
We conjecture that the preceding theorem is true as stated if we allow $\gamma$ to be a
virtual spherical representation rather than a virtual orthogonal representation.
If true, that would allow us to grade ordinary homology and cohomology on
virtual spherical representations throughout this work.
\end{remark}

The calculation of $\stab\Pi_\delta B(x,y)$ generalizes the calculation of stable maps
between orbits in \cite[V.9.4]{LMS:eqhomotopy}. The case where $B$ is a point and
$\delta=0$ gives that result, because the condition in the theorem amounts
to saying that $N_H L/L$ is finite, as in \cite[V.9.4]{LMS:eqhomotopy}.
At the other extreme, when $\delta = \Lie$, we get a calculation of the stable maps between the duals of orbits,
which, as it should, simply has the condition reversed to say that
$N_{K^g}L/L$ is finite.

A crucial ingredient in the proof is a general position result.
This result is in the same vein as our results on ``pseudotransversality'' in
\cite{CW:vectorfields}.
In the following we shall mention several times that $V$ should be ``sufficiently large''
that any two embeddings $G/L\subset V$ are isotopic and that any two such isotopies
are themselves homotopic through isotopies. It suffices that $V$ satisfy the {\em gap hypothesis}
that, if $L\leq J$ are two isotropy groups of $V$, then either $V^J = V^L$ or
\[
 |V^J|\leq |V^L| - |G| -3.
\]
This is true, for example, if each irreducible representation contained
in $V$ appears at least $|G|+3$ times.

Consider, then, the following context.
Let
\[
 f\colon M = G\times_H N \to G_+\smsh_K S^{V-\delta(G/K)}
\]
be a $G$-map where $N$ is a compact $H$-manifold with boundary whose tangent bundle
is trivial and equal to $V-\delta(G/H)+n$ for some integer $n$.
(As earlier, we assume chosen an embedding of $G/H$ in $V$ so that
$\delta(G/H)\subset \Lie(G/H)\subset V$.)
Suppose $x = [e,x']\in G\times_H N$ is a point mapping to $eH$ in the projection to $G/H$.
Let $L$ be the isotropy subgroup of $x$, so $L\leq H$.
Consider the tangent plane to $x'$ in $N$, which is
$V - \delta(G/H) + n$ by assumption. The tangent vectors to the $H$-orbit of $x'$ then give us an
inclusion $\Lie(H/L) \subset V - \delta(G/H) + n$ and the vectors perpendicular
to $\Lie(H/L)$, i.e., in $V - \delta(G/H) + n - \Lie(H/L)$, are the vectors
normal to the orbit of $x'$.
We can then write the space of tangent vectors to $x$ in $M$  normal to the orbit of $x$ as
\[
	\tau_x - \Lie(G/L) = V - \delta(G/H) + n - \Lie(H/L).
\]
Note that we have the inclusion
\[
 \delta(G/L) = \delta(G/H) + \delta(H/L) \subset \delta(G/H) + \Lie(H/L),
\]
which leads to the inclusion
\[
 V - \delta(G/H)-\Lie(H/L) \subset V - \delta(G/L).
\]
This copy of $V-\delta(G/L)$ can be thought of as coming from an embedding
of $G/L$ in $V$ close to the given embedding of $G/H$ in $V$;
if $V$ is sufficiently large,
then there is a canonical identification (up to homotopy through orthogonal $L$-maps)
of this copy of $V-\delta(G/L)$ with the one
we get from any other embedding of $G/L$ in $V$.

Now suppose that $f(x) \in G/K\times 0$, say, $f(x) = gK\times 0$.
Then we have $L\leq K^g = gKg^{-1}$.
We can write the space of tangent vectors at $f(x)$ normal to $G/K\times 0$ as
\[
 \tau_{f(x)} - \Lie(G/K^g) = V - \delta(G/K^g).
\]
Note that we have an injection
\[
 i\colon V - \delta(G/H)-\Lie(H/L) \subset V - \delta(G/L) \includesin V - \delta(G/K^g).
\]
The last map we should define using an isotropy from an embedding of $G/L$ near $G/H$ to
one near $G/K^g$; this is well-defined up to homotopy through orthogonal maps if
$V$ is sufficiently large.

At a point $x$ with $f(x)\in G/K\times 0$, we know that $f$ must take the whole orbit of $x$
to $G/K\times 0$. To specify what we mean by $f$ being in general position, we are left with
specifying the behaviour of the derivative of $f$ on the space of tangent vectors at $x$
normal to the orbit of $x$. We decompose this space into two pieces, the $L$-fixed vectors and the
vectors perpendicular to those:
\begin{align*}
 \tau_x - \Lie(G/L) &= [\tau_x - \Lie(G/L)]^L + [\tau_x - \Lie(G/L)]_L \\
  &= [V - \delta(G/H) + n - \Lie(H/L)]^L + [V - \delta(G/H) - \Lie(H/L)]_L
\end{align*}

\begin{definition}
A map $f\colon M = G\times_H N\to G_+\smsh_K S^{V-\delta(G/K)}$, with $N$ as above,
is said to be in {\em general position} if it is smooth in a neighborhood of
$f^{-1}(G/K\times 0)$ and, for each point $x\in M$ over $eH$ with $f(x)\in G/K\times 0$,
we have, writing $L$ for the isotropy subgroup of $x$ and 
$Df\colon \tau_x\to \tau_{f(x)}$ for the induced map of tangent planes,
\begin{enumerate}
\item
the following map induced by $Df^L$ is an epimorphism:
\[
 [\tau_x - \Lie(G/L)]^L \to \tau_{f(x)}^L \twoheadrightarrow [\tau_{f(x)} - \Lie(G/K^g)]^L,
\]
where the second map is the orthogonal projection, and
\item
on the piece $[V - \delta(G/H)-\Lie(H/L)]_L$, $Df$ is the injection $i_L$ with $i$ as above.
\end{enumerate}
\end{definition}

We also have the following definition of general position of a vector field
from \cite{CW:vectorfields}, where it was called canonical pseudotransversality.

\begin{definition}
Let $M$ be a smooth $G$-manifold and let $\tau\colon TM\to M$ be its tangent bundle. 
We say that a smooth section $s\colon M\to TM$ 
is in {\em general position} if, for each $x\in M$ with $s(x) = x$
(i.e., equal to $0\in \tau_x$), the composite
\[
 \tau_x \xrightarrow{Ds} \tau_x\dirsum \tau_x \xrightarrow{p_2} \tau_x
\]
is the identity on $[\tau_x - \Lie(G/G_x)]_{G_x}$ and takes
$[\tau_x - \Lie(G/G_x)]^{G_x}$ isomorphically onto itself.
\end{definition}

\begin{proposition}[General Position]\label{prop:generalPosition}
Let $M = G\times_H N$ where $N$ is a compact $H$-manifold with boundary whose tangent bundle is trivial
and equal to $V-\delta(G/H)+n$ for some integer $n$. Let $f\colon M\to G_+\smsh_K S^{V-\delta(G/K)}$
be a $G$-map that is in general position on a neighborhood of $\bndry M$.
Then, if $V$ is sufficiently large in the sense above,
$f$ is $G$-homotopic rel $\bndry M$ to a map in general position.
\end{proposition}

\begin{proof}
The argument is similar to those given by Wasserman \cite{Was:difftopology} and by
Hauschild \cite{Haus:Transversality}, \cite{Haus:Splitting}.
The proof is by induction on the isotropy groups appearing in $M$.
Let $L$ be a maximal isotropy group appearing in $M$; we may assume that $L\leq H$.
Using the argument in \cite[I.4]{Haus:Transversality}, for example, we can
homotope $f$ rel boundary to a $G$-map transverse on the $L$-fixed sets.
We can then $G$-homotope $f$ rel boundary and $L$-fixed sets
so that, at each $L$-fixed point mapping to $G/K\times 0$, in the directions
normal to the fixed sets and the orbit of the point, $f$ agrees with $i_L$ as above.

Then $f$ is in general position in a neighborhood of $GM^L$ with no points on the boundary
of the neighborhood mapping to $G/K\times 0$. We can then excise this neighborhood and,
by induction, put $f$ in regular position on the remainder of $M$.
\end{proof}

\begin{proof}[Proof of Theorem~\ref{thm:StableMapsOrbitsOverB}]
We first define an auxiliary category $\D_\delta$.
Its objects are the same as those of $\stab\Pi_\delta B$, that is,
maps $x\colon G/H\to B$. 
The set of maps from $x$ to $y\colon G/K\to B$ is the 
Grothendieck group on the set of equivalence
classes of diagrams $x\from M\laxto y$ where $M$ is a closed $G$-manifold
thought of as a manifold over $B$ via $M\to G/H \xrightarrow{x} B$.
$M$ has to satisfy a certain condition to be explained, and the equivalence
relation must also be explained.

Name the maps from $M$ to the orbits $\alpha\colon M\to G/H$ and
$\beta\colon M\to G/K$. Let $N$ be the fiber of $\alpha$, so that
$M = G\times_H N$. We can then decompose the tangent bundle of $M$ as
\[
 \tau_M \iso \alpha^*(G\times_H \Lie(G/H)) \dirsum (G\times_H \tau_N).
\]
In other words, we identify $\alpha^*(G\times_H \Lie(G/H))$ with the subbundle
of tangent vectors normal to the fiber of $\alpha$.
Similarly, we can identify $\beta^*(G\times_K\Lie(G/K))$ with the subbundle
of tangent vectors normal to the fiber of $\beta$,
and within this bundle we have its subbundle $\beta^*(G\times_K\delta(G/K))$.
The condition we impose on $M$ is then that
\[
 \beta^*(G\times_K\delta(G/K)) 
   \subset \alpha^*(G\times_H \delta(G/H)) \dirsum (G\times_H \tau_N).
\]
Another way of putting this is that $\alpha$ carries
$\beta^*(G\times_K\delta(G/K))$ into $G\times_H \delta(G/H)$.
This is automatically true if $M$ is an orbit of $G$, by the additivity of $\delta$.

Now for the equivalence relation. Given the diagram
$x\from M\laxto y$, we define a stable map as follows.
Choose a $G$-representation $V$ so large that we can embed $G/H$, so
that $V$ has a distinguished $H$-subspace isomorphic to $\Lie(G/H)$.
We assume $V$ is also large enough that we can embed $N$ in $V-\delta(G/H)$,
which we extend to an embedding of $M$ in $G\times_H(V-\delta(G/H))$.
The normal bundle $\nu$ to this embedding will then be
\[
 \nu \iso V - \alpha^*(G\times_H \delta(G/H)) - G\times_H \tau_N.
\]
By the condition imposed on $M$, this bundle is a subbundle of
$V - \beta^*(G\times_K\delta(G/K))$. Hence, we can write down the following
map:
\[
 G\times_H S^{V-\delta(G/H)}
  \xrightarrow{c} S^\nu
  \to S^{V - \beta^*(G\times_K\delta(G/K))}
  \xrightarrow{\beta} G\times_K S^{V-\delta(G/K)},
\]
where $S^\nu$ denots the sphere bundle over $M$ obtained by taking the one-point
compactification of each fiber of $\nu$, and similarly for
$S^{V - \beta^*(G\times_K\delta(G/K))}$.
Now, using the homotopy involved in the lax map $M\laxto y$, we get
a lax stable map
\[
 \chi(x \from M\laxto y)\colon G_+\smsh_H S^{-\delta(G/H),x}
  \laxto G_+\smsh_K S^{-\delta(G/K),y},
\]
an element of $\stab\Pi_\delta B(x,y)$, which we call the
{\em Euler characteristic} of the diagram.
We consider two diagrams to be equivalent if they have the same
Euler characteristic.
(Note that, if $G$ is trivial and $B$ is a point, then
$\stab\Pi_\delta B(x,y) \iso \Z$ and $\chi(M)$ is
the classical Euler characteristic of $M$.)

We define the sum of two such diagrams $[x\from M\laxto y]$ and $[x\from N\laxto y]$
to be $[x\from M\disjunion N \laxto y]$ and then,
as mentioned above, $\D_\delta(x,y)$ is
defined to be the Grothendieck group of the resulting monoid.

Composition is defined by taking pullbacks.
Precisely, if given two diagrams $[x\from M \laxto y]$ and
$[y\from N \laxto z]$, we let $P$ be the pullback in the following diagram:
\[
 \xymatrix@R-1em@C-1em{
  &&P \ar[dl] \ar@{=>}[dr] \\
  & M \ar[dl] \ar@{=>}[dr] && N \ar[dl] \ar@{=>}[dr] \\
  x \ar[drr] && y \ar[d] && z \ar[dll] \\
  && B
  }
\]
The homotopy involved in the lax map $P\laxto N$ comes from the homotopy
underlying $M\laxto y$.
We then define the composite
\[
 [y\from N\laxto z]\circ [x\from M\laxto y] = [x \from P \laxto z].
\]
It's a formal check that $x \from P \laxto z$ satisfies the condition required of it.
It follows from the construction that
\[
 \chi(x \from P\laxto z) = \chi(y\from N\laxto z) \circ \chi(x\from M\laxto y),
\]
from which it follows that composition respects the equivalence relation,
hence we have a well-defined category $\D_\delta$.

In the process of defining $\D_\delta$, note that we've also defined a functor
$\chi\colon \D_\delta \to \stab\Pi_\delta B$, which,
by definition, is a monomorphism on each mapping group.

Now suppose that
$f\colon G_+\smsh_H S^{V-\delta(G/H),x} \laxto G_+\smsh_K S^{V-\delta(G/K),y}$
represents an element of $\stab\Pi_\delta B(x,y)$.
By Proposition~\ref{prop:generalPosition}, we may assume that $f$ is in general position
(away from the basepoint).
The inverse image of $G/K\times 0$ is then a disjoint union of orbits of $G$:
If $m\in S^{V-\delta(G/H)}$ is a point such that $f(m)\in G/K\times 0$, with $Gm \iso G/L$, 
then the tangent plane at $m$ normal to the orbit of $m$ is isomorphic to
$V - \delta(G/H) - \Lie(H/L)$.
If $f(m) = gK\times 0$, then we can write the codimension of $G/K\times 0$ in the target
of $f$ as $V - \delta(G/K^g)$.
As noted while defining general position, we have
an inclusion
\[
 V - \delta(G/H) - \Lie(H/L) \subset V-\delta(G/L) \subset V-\delta(G/K^g).
\]
On $L$-fixed sets, this means that the dimension
$|[V - \delta(G/H) - \Lie(H/L)]^L|$ is no more than $|[V-\delta(G/K^g)]^L|$.
In order to satisfy the first condition of being in general position, these dimensions
must therefore be the same and $Df^L$ must be an isomorphism.
Together with the second condition, this implies that the only tangent directions
at $m$ that can be taken into $G/K\times 0$ are those along the orbit of $m$,
hence the orbit of $m$ is isolated in $f^{-1}(G/K\times 0)$.
Thus, $f^{-1}(G/K\times 0)$ is a disjoint union of orbits as claimed.

The inclusions above also tell us exactly which orbits may appear in the inverse image,
namely those for which the inclusion are equalities on fixed sets:
We need to have both
\[
 [\Lie(H/L) - \delta(H/L)]^L = 0
\]
(using $\delta(G/L) = \delta(G/H) + \delta(H/L)$) and
\[
 \delta(K^g/L)^L = 0
\]
(using $\delta(G/L) = \delta(G/K^g) + \delta(K^g/L)$). We can rewrite these as
\[
 \Lie(N_H L/L) - \delta(N_H L/L) = 0
\]
and
\[
 \delta(N_{K^g} L/L) = 0
\]
as in the statement of the theorem.

Each orbit in the inverse image determines a diagram $x \from G/L \laxto y$;
we write $z_i = [x \from G/L \laxto y]$ for the class of this diagram in $\D_\delta(x,y)$.
Also, for each orbit, the map $Df^L$ above is an isomorphism, so determines
a sign $\epsilon_i = \pm 1$, depending on whether $Df^L$ preserves orientation or not,
using the canonical identification of
$|[V - \delta(G/H) - \Lie(H/L)]^L|$ with $|[V-\delta(G/K^g)]^L|$ for $V$ sufficiently large.
It then follows that
\[
 [f] = \chi\left(\sum_i \epsilon_i z_i\right)
\]
where the sum runs over the orbits in the inverse image of $G/K\times 0$.
This shows that $\chi$ is an epimorphism on mapping groups, hence is an isomorphism
of categories.

Moreover, the argument above shows that $\D_\delta(x,y) \iso \stab\Pi_\delta B(x,y)$ is generated by the
classes $[x\from G/L \laxto y]$ for orbits $G/L$ satisfying the conditions in the statement
of the theorem.
To show that these classes freely generate, suppose that we have a sum
$\sum_i n_i z_i$ such that $\chi(\sum_i n_i z_i) = 0$, where the $z_i$ are all distinct
generators. Let $f$ be the map constructed above to represent
$\chi(\sum_i n_i z_i)$, so that $f$ is in general position and the inverse image
of $G/K\times 0$ consists of $n_i$ copies of each orbit underlying $z_i$.
We are assuming that $f$ is null homotopic and we may assume we have a null homotopy $h$
in general position.
For each orbit $G/L$ that appears in a $z_i$, we know that
\[
 |(V-\delta(G/H)-\Lie(H/L))^L| = |(V-\delta(G/K^g))^L|,
\]
hence that
\[
 |(V-\delta(G/H)-\Lie(H/L)+1)^L| = |(V-\delta(G/K^g))^L| + 1.
\]
So, each of these orbits appearing in $f^{-1}(G/K\times 0)$ will be one end
of a manifold of the form $G/L\times I \subset h^{-1}(G/K\times 0)$.
Because $h$ is a null homotopy, the other end of this manifold must be on the same
face of $G_+\smsh_H S^{V-\delta(G/)}\smsh I_+$, so will have opposite orientation.
The only way this can happen without contradiction is to have each $n_i = 0$.
Thus, $\D_\delta(x,y)$ is free as claimed in the statement of the theorem.

We now define a functor $\psi\colon\D_\delta\to \stab\Pi_{\delta,\gamma}B$.
On objects it is the same as $\chi$; on morphisms we define what it does on a generator:
Given $x\from G/L \laxto y$, let $\alpha\colon G/L\to G/H$ and $\beta\colon G/L\to G/K$
be the corresponding maps of orbits. 
Writing $z = x\circ\alpha$, we define $\psi(x\from G/L \laxto y)$ to be
the lax stable composite
\begin{align*}
 G_+\smsh_H S^{\gamma(x)-\delta(G/H),x}
  & \xrightarrow{c} G_+\smsh_L S^{\gamma(x)-\delta(G/H)-\Lie(H/L),z} \\
  & \includesin G_+\smsh_L S^{\gamma(x)-\delta(G/L),z} \\
  & \xrightarrow{\gamma} G_+\smsh_L S^{\gamma(y)-\delta(G/L),z} \\
  & \includesin G_+\smsh_L S^{\gamma(y) - \delta(G/K),z} \\
  & \laxto G_+\smsh_K S^{\gamma(y)-\delta(G/K),y}.
\end{align*}
The collapse map $c$ is the same as used in the definition of $\chi$.
The map labeled $\gamma$ is the one obtained by applying $\gamma$ to the
lax map $G/L\laxto y$.
The inclusions come from the fact that $\delta$ is a dimension function,
while the last map comes from the projection $G/L\to G/K$.

We next have to show that $\psi$ is a functor, i.e., that it respects composition.
Let $x\colon G/H\to B$, $y\colon G/K\to B$, and $z\colon G/J\to B$.
Consider two generators
$x\from G/L \laxto y$ and $y\from G/I \laxto z$.
Let $P$ be the pullback in the following diagram:
\[
 \xymatrix@R-1em@C-1em{
  &&P \ar[dl] \ar@{=>}[dr] \\
  & G/L \ar[dl] \ar@{=>}[dr] && G/I \ar[dl] \ar@{=>}[dr] \\
  x \ar[drr] && y \ar[d] && z \ar[dll] \\
  && B
  }
\]
In order to understand what $\psi$ does to the composite, we need to
describe the composite $x\from P\laxto z$ in terms of the generators.
In \cite{CW:vectorfields}, we showed that an equivariant manifold like $P$
has a vector field in general position. In particular, its set of zeros
is a disjoint union of orbits in $P$ and to each of the orbits we can assign a sign
$\pm 1$ determined by the behaviour of the vector field in a neighborhood of the orbit.
This determines a signed sum of generators.
If we apply $\chi$ to $x\from P\laxto z$, we can use the vector field on $P$ to homotope
the resulting map to the image under $\chi$ of that sum of generators.
Therefore, $x\from P\laxto z$ equals this sum of generators in $\D_\delta$.
Now, we can use the vector field in the same way to homotope
$\psi(y\from G/I\laxto z)\circ \psi(x\from G/L\laxto y)$ to the image under $\psi$
of the same sum of generators (choosing a representative isomorphism in $\gamma$
for each component of $P/G$), showing that
\[
 \psi(y\from G/I\laxto z)\circ \psi(x\from G/L\laxto y) = 
 \psi(x\from P\laxto z),
\]
i.e., that $\psi$ respects composition.

Suppose that $f\colon G_+\smsh_H S^{V+\gamma_0(x)-\delta(G/H),x}
 \laxto G_+\smsh_K S^{V+\gamma_0(y)-\delta(G/K),y}$
represents an element of $\stab\Pi_{\delta,\gamma} B$.
Let $U\subset G_+\smsh_K S^{V+\gamma_0(y)-\delta(G/K),y}$ be a $G$-invariant open neighborhood
of $G_+\smsh_K 0$ not including the base section. Then $f^{-1}(U)$ is an open subspace of 
$G_+\smsh_H S^{V+\gamma_0(x)-\delta(G/H),x}$ not including the base section.
Let $G\times_J A$ be the orbit of a component of $f^{-1}(U)$.
Each point in $A$ determines, via $f$, a homotopy in $B^J$, and any two points
in $A$, being connected by a path, determine the same homotopy, up to path
homotopy rel endpoints. Therefore, $A$ determines a class of paths, hence, via $\gamma$,
a homotopy class of isomorphisms of $\gamma_0(x)$ with $\gamma_0(y)$.
Pick an isomorphism in that class.
Using that isomorphism to identify $V+\gamma_0(x)$ with $V+\gamma_0(y)$,
we can put $f$ in general position inside $G\times_J A$. Doing this for each
component of $f^{-1}(U)$, we put $f$ in general position.
It follows that $\psi$ is an epimorphism. A similar argument applied to homotopies
shows that $\psi$ is a monomorphism. Hence, $\psi$ is an isomorphism of categories.

We have therefore shown that have the following isomorphisms:
\[
 \stab\Pi_\delta B \iso \D_\delta \iso \stab\Pi_{\delta,\gamma} B,
\]
which proves the theorem.
\end{proof}


\section{Parametrized homology and cohomology theories}

The cellular theories we shall define in the following chapter can be thought
of as consisting of a series of related ``$RO(G)$-graded parametrized homology and cohomology
theories.'' We define exactly what we mean by that, concentrating on reduced theories.

In \cite{Dold:parahomology} Dold gave axioms
for nonequivariant homology and cohomology theories defined on
spaces parametrized by $B$.
Several authors, including Wirthm\"uller \cite{Wi:duality}, have given
axioms for $RO(G)$-graded theories defined on $G$-spaces;
we follow those given by May in \cite[XIII.1]{May:alaska}.
We combine the two approaches to give axioms for
equivariant theories defined on parametrized spaces;
this treatment is essentially the same as that of
May and Sigurdsson \cite[\S21.1]{MaySig:parametrized}.

Let $\RO(G)$ be the category of finite-dimensional
linear representations of $G$ and $G$-linear isometric
isomorphisms between them.
Say that two maps $V\to W$ are homotopic if the associated $G$-maps
$S^V\to S^W$ are stably homotopic. Let $h\RO(G)$ denote the resulting
homotopy category.
(Comparing with \cite{May:alaska}, we are assuming that we are using a complete
$G$-universe $\U$. Variants for smaller universes are possible, but we are interested
in duality, which at present requires a complete universe.)

\begin{definition}\label{def:parametrizedtheory}
An {\em $RO(G)$-graded parametrized (reduced) homology theory} consists of a 
functor
\[
 \tilde h^G_*\colon \Op{h\RO(G)}\times \Ho G\Para B \to \Ab,
\]
written $(V,X) \mapsto \tilde h^G_V(X)$ on objects. 
Similarly, an {\em $RO(G)$-graded parametrized (reduced)
cohomology theory} consists of a functor
\[
 \tilde h_G^*\colon h\RO(G)\times \Op{\Ho G\Para B} \to \Ab.
\]
These functors are required to satisfy the following
axioms.
 \begin{enumerate}
\item For each $V$, the functor $\tilde h^G_V$ is exact on cofiber sequences
and sends wedges to direct sums. The functor $\tilde h_G^V$ is exact on cofiber sequences
and sends wedges to products.
\item For each $W$ there are natural isomorphisms
\[
 \sigma^W\colon \tilde h^G_V(X)\to \tilde h^G_{V\dirsum W}(\susp_B^W X)
\]
and
\[
 \sigma^W\colon \tilde h_G^V(X)\to \tilde h_G^{V\dirsum W}(\susp_B^W X).
\]
$\sigma^0 = \id$ and
$\sigma^Z\circ \sigma^W = \sigma^{W\dirsum Z}$
for every pair of representations $W$ and $Z$.
\item The $\sigma$ are natural in $W$ in the following sense. If
$\alpha\colon W\to W'$ is a map in $h\RO(G)$, then the following diagrams commute.
\[
 \xymatrix{
 {\tilde h_V^G(X)} \ar[r]^-{\sigma^{W'}} \ar[d]_{\sigma^{W}} &
   {\tilde h_{V\dirsum W'}^G(\susp_B^{W'} X)} \ar[d]^{\tilde h_{\id\dirsum\alpha}^G(\id)} \\
 {\tilde h_{V\dirsum W}^G(\susp_B^W X)} \ar[r]_{(\susp_B^\alpha\id)_*} &
   {\tilde h_{V\dirsum W}^G(\susp_B^{W'} X)} \\
 }
\]
\[
 \xymatrix{
 {\tilde h^V_G(X)} \ar[r]^-{\sigma^{W}} \ar[d]_{\sigma^{W'}} &
   {\tilde h^{V\dirsum W}_G(\susp_B^{W} X)} \ar[d]^{\tilde h^{\id\dirsum\alpha}_G(\id)} \\
 {\tilde h^{V\dirsum W'}_G(\susp_B^{W'} X)} \ar[r]_{(\susp_B^\alpha\id)^*} &
   {\tilde h^{V\dirsum W'}_G(\susp_B^{W} X)} \\
 }
\]
 \end{enumerate}
\end{definition}
\bigskip

As in \cite{May:alaska} we extend the grading to ``formal differences'' $V\ominus W$ by setting
\[
 \tilde h^G_{V\ominus W}(X) = \tilde h^G_V(\susp_B^W X)
\]
and similarly in cohomology. Rigorously, we extend the grading from
$\Op{h\RO(G)}$ to $\Op{h\RO(G)}\times h\RO(G)$. Now, $RO(G)$ can be obtained from
$\RO(G)$ by considering $V\ominus W$ to be equivalent to $V'\ominus W'$ when there
is a $G$-linear isometric isomorphism
\[
 \alpha\colon V\dirsum W'\to V'\dirsum W.
\]
Given such an isomorphism we get the following diagram.
\[
 \xymatrix{
  {\tilde h^G_{V'}(\susp_B^{W'} X)} \ar[r]^-{\sigma^W} \ar[d] &
   {\tilde h^G_{V'\dirsum W}(\susp_B^{W'\dirsum W} X)} \ar[d]^{\tilde h^G_\alpha(\susp^\tau\id)} \\
  {\tilde h^G_V(\susp_B^W X)} \ar[r]_-{\sigma^{W'}} &
   {\tilde h^G_{V\dirsum W'}(\susp_B^{W\dirsum W'} X)}
 }
\]
Here, $\tau\colon W\dirsum W' \to W'\dirsum W$ is the transposition map.
Thus,
\[
 \tilde h^G_{V\ominus W}(X) \iso \tilde h^G_{V'\ominus W'}(X).
\]
(There is a similar isomorphism in cohomology.)
It is in this sense
that we can think of grading on $RO(G)$, but note that the isomorphism depends on the
choice of $\alpha$, so is not canonical.

Dold pointed out that
$\tilde h^G_*(-)$ and
$\tilde h_G^*(-)$ restrict to functors on $\Pi B$. 
Actually, the suspension isomorphisms imply that
homology and cohomology theories extend to stable maps, and so we 
get restrictions to functors on the stable fundamental groupoid $\stab\Pi B$.
More generally, we can make the following definitions.

\begin{definition}
If $\tilde h^G_*(-)$ is an $RO(G)$-graded parametrized homology theory on
ex-$G$-spaces over $B$ and
$H\leq G$ is a subgroup, define an
$RO(H)$-graded parametrized homology theory $\tilde h^H_*(-)$ on ex-$H$-spaces over $B$ by setting
\[
 \tilde h^H_V(X) = \tilde h^G_{W}(G_+\smsh_H \susp^{W-V}X)
\]
where $W$ is any representation of $G$ containing $V$ as a sub-$H$-representation.
Similarly,
if $\tilde h_G^*(-)$ is an $RO(G)$-graded parametrized cohomology theory, define an
$RO(H)$-graded parametrized cohomology theory $\tilde h_H^*(-)$ by setting
\[
 \tilde h_H^V(X) = \tilde h_G^{W}(G_+\smsh_H \susp^{W-V}X).
\]
\end{definition}

\begin{definition}\label{def:homologymackey}
If $\tilde h^G_*(-)$ is an $RO(G)$-graded parametrized homology theory,
$\gamma$ is a virtual representation of $\Pi B$,
and $\delta$ is a dimension function for $G$,
we define a covariant $\stab\Pi_\delta B$-module $\MackeyOp h^{G,\delta}_\gamma$ by
\begin{align*}
 \MackeyOp h^{G,\delta}_\gamma(b\colon G/H\to B) 
  &= \tilde h^H_{\gamma_0(b)}(S^{-\delta(G/H),b}) \\
  &= \tilde h^H_{\gamma_0(b)+V}(S^{V-\delta(G/H),b}),
\end{align*}
where $V$ is a representation of $H$ so large that $\delta(G/H)\subset V$.
This is a functor on $\stab\Pi_{\delta,-\gamma} B$, so we make it a functor
on $\stab\Pi_\delta B$ using Theorem~\ref{thm:StableMapsOrbitsOverB}.
We call these modules the {\em coefficient systems} of $\tilde h^G_*(-)$.
Similarly, if $\tilde h_G^*(-)$ is an $RO(G)$-graded parametrized cohomology theory,
we write $\Mackey h_{G,\delta}^\gamma$ for the contravariant 
$\stab\Pi_\delta B$-module defined by
\begin{align*}
 \Mackey h_{G,\delta}^\gamma(b\colon G/H\to B) 
  &= \tilde h_H^{\gamma_0(b)}(S^{-\delta(G/H),b}) \\
  &= \tilde h_H^{\gamma_0(b)+V}(S^{V-\delta(G/H),b}).
\end{align*}
\end{definition}

As usual, if $X$ is an unbased space over $B$, or $(X,A)$ is a pair of
based or unbased $G$-spaces over $B$, we define
 \[
  h^G_*(X) = \tilde h^G_*(X_+)
 \]
and
 \[
  h^G_*(X,A) = \tilde h^G_*(C_B i)
 \]
where $i\colon A\to X$ is the inclusion and $C_B i$ is its cofiber over $B$. 
We use a similar convention for cohomology.

We are interested in how these theories are represented. 
It is well-known how to represent parametrized cohomology
theories using spectra over $B$ 
(see, for example, \cite{Cl:paraspectra} and \cite{CP:paraspectra}).
The following is part of Theorem~21.2.3 in \cite{MaySig:parametrized},
and follows from Brown's representability theorem in the usual way.

\begin{theorem}\label{thm:cohomrep}
If $\tilde h_G^*(-)$ is an $RO(G)$-graded cohomology theory on $\Ho G\Para B$, then it is
represented by an $\Omega$-$G$-prespectrum $E$ over $B$, in the sense
that, for each $V$, there is a natural isomorphism
 \[
 \tilde h_G^V(X) \iso [\Susp_B^\infty X,\Susp_B^V E]_{G,B}.
 \]
The prespectrum $E$ is determined by $h$ up to non-unique equivalence.
\qed
\end{theorem}

We now consider homology theories.
Curiously, previously to \cite{MaySig:parametrized} and this paper,
there appears to have been no discussion in the literature
of representing homology theories on parametrized
spaces. The next proposition gives a way of defining a homology theory using
a parametrized prespectrum. Before that we should say a few words about
smash products. (We summarize the more detailed discussion \cite{MaySig:parametrized}.)

The fiberwise smash product behaves
poorly from the point of view of homotopy theory; it does not pass directly
to homotopy categories. However, the {\em external} fiberwise smash product does
behave well. If $D$ is a $G$-prespectrum over $A$ and $E$ is a $G$-prespectrum
over $B$ we can form the external smash product $D\exsmsh E$ over $A\times B$.
(Since we are using prespectra and not orthogonal spectra, the smash product
we use is either ``external'' with respect to indexing spaces as well, or
is a ``handicrafted'' smash product \cite[13.7.1]{MaySig:parametrized}.)
The internal smash product of spectra $D$ and $E$ in $\Ho G\PreSpec{} B$
is then defined by
$D\smsh_B E = \Delta^*(D\exsmsh E)$,
where $\Delta\colon B\to B\times B$ is the diagonal. 
As usual, this notation is convenient but misleading.
In the homotopy category, $\Delta^*$ and $\exsmsh$ must be understood as
derived functors of their point-set level precursors. In particular,
$\Delta^*$ is defined by first applying a fibrant approximation functor and
then applying the point-set level $\Delta^*$.
Fortunately, May and Sigurdsson show that the familiar compatibility
relations among the smash product and the change of basespace functors
remain true in the homotopy category.

\begin{proposition}\label{prop:representedHomology}
If $E$ is a $G$-prespectrum over $B$ then the groups
 \[
 \tilde E^G_V(X) = [S^V, \rho_!(E\smsh_B X)]_G
 \]
define an $RO(G)$-graded homology theory on $\Ho G\Para B$.
Here, $\rho\colon B\to *$ denotes projection to a point.
\qed
\end{proposition}

\begin{proof}
(See also \cite[21.2.22]{MaySig:parametrized}.)
The main problem is to see that $\tilde E^G_V(X)$ is exact on cofiber
sequences and takes wedges to direct sums. Write
$\rho_!(E\smsh_B X) = \rho_!\Delta^*(E\exsmsh X)$.
As in \cite{MaySig:parametrized},
each of $\rho_!$, $\Delta^*$, and $E\exsmsh -$ 
is exact on $\Ho G\PreSpec{} B$ and also takes wedges to wedges.
The remainder of the proof is as in \cite[XIII.2.6]{May:alaska}.
\end{proof}

For a discussion of the representability of general
$RO(G)$-graded parametrized homology theories, see \cite[21.5]{MaySig:parametrized}.
Unfortunately, the story is not completely straightforward.
Fortunately, this will not be a problem when we want to discuss the spectra
representing cellular homology.

At first glance, the construction of Proposition~\ref{prop:representedHomology}
is not an obvious way to define a homology
theory using a spectrum. In particular, the use of
$\rho_!$ on the right might not be expected,
given that $\rho_!$ is a left adjoint. One
justification is that it works: We do get a homology theory and,
with some assumptions on the base space, every homology theory
is representable in this way.
Here is, perhaps, a better justification.
The homology and cohomology theories determined by a parametrized
spectrum $E$ should, locally, reflect the nonparametrized theories
given by the fibers of $E$. To make this precise, we first record
a change of basespace result.

\begin{proposition}\label{prop:hombasechange}
Let $\alpha\colon A\to B$ be a $G$-map, let $X$ be an ex-$G$-space
over $A$, and let $E$ be a $G$-prespectrum over $B$. Then
\[
 \tilde E^G_V(\alpha_! X) \iso (\widetilde{\alpha^*E})^G_V(X)
\]
and
\[
 \tilde E_G^V(\alpha_! X) \iso (\widetilde{\alpha^*E})_G^V(X).
\]
\end{proposition}

\begin{proof}
We begin with cohomology, which is a little more obvious:
\begin{align*}
 \tilde E_G^V(\alpha_! X)
 &= [\alpha_! X, \susp^V_B E]_{G,B} \\
 &\iso [X, \alpha^*(\susp^V_B E)]_{G,A} \\
 &\iso [X, \susp^V_A \alpha^*E]_{G,A} \\
 &= (\widetilde{\alpha^*E})_G^V(X).
\end{align*}
For homology we have the following, using one of the
stable equivalences listed in Section~\ref{sec:paraspec}:
\begin{align*}
 \tilde E^G_V(\alpha_! X)
 &= [S^V, \rho_!(E\smsh_B \alpha_! X)]_G \\
 &\iso [S^V, \rho_!\alpha_!(\alpha^*E\smsh_A X)]_G \\
 &= [S^V, \rho_!(\alpha^*E\smsh_A X)]_G \\
 &= (\widetilde{\alpha^*E})^G_V(X).
\end{align*}
\end{proof}

We can now say that these theories behave correctly
over individual points of $B$:

\begin{corollary}
Let $E$ be a $G$-prespectrum over $B$, let $b\colon G/H\to B$,
and let $X$ be a based $H$-space. Then
\[
 \tilde E^G_V(b_!(G \times_H X)) \iso (\tilde E_b)^H_V(X)
\]
and
\[
 \tilde E_G^V(b_!(G \times_H X)) \iso (\tilde E_b)_H^V(X)
\]
where $b^* E = G \times_H E_b$.
\qed
\end{corollary}

Thus, the homology and cohomology theories defined by $E$
restrict over each point $b$ of $B$ to the nonparametrized
theories defined by the fiber $E_b$.

Specializing further, we can examine the coefficient systems
of the homology and cohomology theories determined by $E$.

\begin{corollary}\label{cor:coefficientsfibers}
Let $\gamma$ be a virtual representation of $\Pi B$ and
let $\delta$ be a dimension function for $G$.
If $E$ is a $G$-prespectrum over $B$, then the coefficient systems
of its associated homology and cohomology theories are given by
\[
 \MackeyOp E^{G,\delta}_\gamma(b\colon G/H\to B)
  \iso \pi^H_{\gamma_0(b)-(\Lie-\delta)(G/H)}(E_b)
\]
and
\[
 \Mackey E_{G,\delta}^\gamma(b) \iso \pi^H_{-\gamma_0(b)-\delta(b)}(E_b).
\]
\end{corollary}

\begin{proof}
The calculation for cohomology follows directly from the
definitions, including Definition~\ref{def:homologymackey},
and the preceding corollary.
The calculation for homology uses the
isomorphism of \cite[II.6.5]{LMS:eqhomotopy}.
\end{proof}

We can reframe this result as follows.
For a fixed $b\colon G/H\to B$, define a covariant
$H$-$\delta$-Mackey functor $\MackeyOp E^{G,\delta}_\gamma|b$ by
\[
 (\MackeyOp E^{G,\delta}_\gamma|b)(H/K) =
  \MackeyOp E^{G,\delta}_\gamma(G/K\to G/H\to B).
\]
Then we have
\[
 \MackeyOp E^{G,\delta}_\gamma|b \iso \Mackey\pi^{H,\Lie-\delta}_{\gamma_0(b)}(E_b).
\]
Defining $\Mackey E_{G,\delta}^\gamma|b$ similarly, we get
\[
 \Mackey E_{G,\delta}^\gamma|b \iso \Mackey\pi^{H,\delta}_{-\gamma_0(b)}(E_b).
\]

Finally, we give one more useful application
of Proposition~\ref{prop:hombasechange}:
The groups $\tilde E^G_*(X,p,\sigma)$ and $\tilde E_G^*(X,p,\sigma)$
depend only on the spectrum $p^*E$ parametrized by $X$. To see
this, note first that 
$(X,p,\sigma) = p_!(X',q,\tau)$ where $X' = X\union_{\sigma(B)} X$, 
$q\colon X'\to X$ is the evident map, and $\tau\colon X\to X'$ takes $X$ onto
the second copy of $X$ in $X'$. From this it follows that
\[
 \tilde E^G_V(X) = \tilde E^G_V(p_! X')
  \iso (\widetilde{p^*E})^G_V(X')
\]
and similarly for cohomology.
In particular, if $p^*E$ is trivial, in the sense that it is
equivalent to $\rho^* D = X\times D$ for a nonparametrized
$G$-prespectrum $D$, then (assuming that $X$
is well-sectioned)
 \[
 \tilde E^G_V(X) \iso \tilde D^G_V(\rho_! X') \iso \tilde D^G_V(X/\sigma(B))
 \]
is just the usual $D$-homology of $X/\sigma(B)$, and
 \[
 \tilde E_G^V(X) \iso \tilde D_G^V(X/\sigma(B))
 \]
is just the usual $D$-cohomology of $X/\sigma(B)$.


\section{Duality}

Let $E$ be a $G$-prespectrum parametrized by $B$. 
We would like to relate the homology and cohomology theories determined by
$E$ by finding, for each suitably ``finite'' $G$-prespectrum $X$ over $B$,
a dual prespectrum $\bar D_B X$ over $B$ such that
$E_G^*(X) \iso E^G_{-*}(\bar D_B X)$.
We call this kind of duality {\em homological duality}
to distinguish it from the better-known fiberwise duality
(see, for example, \cite{Cl:paraspectra} or \cite{CP:paraspectra} as well as
\cite{MaySig:parametrized}).
What we call homological duality, May and Sigurdsson have seen fit
to call Costenoble-Waner duality, for which we thank them.

Homological duality does not fit the usual formalism of duality in closed symmetric
monoidal categories---fiberwise duality does. However, we can set up a very similar
formalism for our purposes. 
May and Sigurdsson \cite{MaySig:parametrized}
have developed a more general notion of duality, in what they call closed symmetric bicategories,
that includes both fiberwise duality and homological duality as special cases.
Even though we recommend their treatment, we thought it worth presenting the
one we gave before \cite{MaySig:parametrized} appeared.

It will be convenient to generalize as follows.
Let $A$ and $B$ be basespaces, let
$W$ be a $G$-prespectrum over $A$, and let $E$ be a $G$-prespectrum over $A\times B$.
Writing in terms of mapping sets, we would like a natural isomorphism
of the form
\[
 [W\exsmsh X, E]_{G,A\times B} \iso [W, (1\times\rho)_!(1\times\Delta)^*(E\exsmsh \bar D_B X)]_{G,A}.
\]
Here, $\exsmsh$ denotes the external smash product while,
as is our custom, $\rho$ denotes projection to a point, so
\[
 1\times\Delta\colon A\times B \to A\times B\times B
\]
and
\[
 1\times\rho \colon A\times B \to A.
\]
We are most interested in the case $A = *$, but the extra generality comes at no
cost and gives a considerably more satisfying result.

Recall from
\cite{MaySig:parametrized} that there is an external function spectrum functor
$\bar F$: If $Y$ is a $G$-prespectrum over $B$ and $Z$ is a $G$-prespectrum over
$A\times B$, then $\bar F(Y,Z)$ is a $G$-prespectrum over $A$.
$\bar F(Y,-)$ is right adjoint to $-\exsmsh Y$.

\begin{definition}
If $X$ is a $G$-prespectrum over $B$, let $\bar D_B X = \bar F(X,\Delta_! S_B)$
where $S_B$ is the unit spectrum over $B$ (i.e., $\susp^\infty_B B_+$) and
$\Delta\colon B\to B\times B$ is the diagonal.
\end{definition}

We shall justify this definition over the course of this section.
It's interesting to note that the fiberwise dual of $X$ can be written
as $\bar F(X,\Delta_* S_B)$ where $\Delta_*$ is the right adjoint of $\Delta^*$.
We should reiterate our standard warning at this point: 
All of our constructions are
done in the stable category. The point-set definition of
$\bar F(X,\Delta_! S_B)$ gives a misleading picture of $\bar D_B X$,
so we shall rely heavily on the formal properties of the functors we use,
as given in \cite{MaySig:parametrized}.

Evaluation, the counit of the $\exsmsh$-$\bar F$ adjunction, gives a map
\[
 \epsilon\colon \bar D_B X\exsmsh X \to \Delta_! S_B.
\]

The following innocuous-looking lemma is a key technical point in this circle of ideas.

\begin{lemma}
$(1\times\rho\times 1)_!(1\times\Delta\times 1)^*(E\exsmsh\Delta_! S_B) \iso E$
in the stable category, for every $G$-prespectrum $E$ over $A\times B$.
\end{lemma}

\begin{proof}
Here and elsewhere we use the commutation relation $j^*f_! \iso g_!i^*$ that holds
when the following is a pullback diagram and at least one of $f$ or $j$ is a fibration:
\[
 \xymatrix{
  C \ar[r]^g \ar[d]_i & D \ar[d]^j \\
  A \ar[r]_f & B
 }
\]
(This is \cite[13.7.7]{MaySig:parametrized}.) The proof of the present lemma would be
considerably easier if the fibration condition were not required, but
the result if false in general without it.
The following chain of isomorphisms proves the lemma.
\begin{align*}
 (1\times\rho\times 1)_!&(1\times\Delta\times 1)^*(E\exsmsh\Delta_! S_B) \\
  &\iso (1\times\rho\times 1)_!\Delta_{A\times B\times B}^*
     [(1\times 1\times\rho)^*E\exsmsh (\rho\times 1\times 1)^* \Delta_! S_B] \\
  &= (1\times\rho\times 1)_![(1\times 1\times\rho)^*E\smsh_{A\times B\times B} 
     (\rho\times 1\times 1)^*\Delta_! S_B] \\
  &\iso (1\times\rho\times 1)_![(1\times 1\times\rho)^* E\smsh_{A\times B\times B}
     (1\times\Delta_!)(\rho\times 1)^* S_B] \\
  &\iso (1\times\rho\times 1)_!(1\times\Delta)_!
     [(1\times\Delta)^*(1\times 1\times\rho)^*E\smsh_{A\times B} (\rho\times 1)^*S_B] \\
  &\iso E\smsh_{A\times B} S_{A\times B} \\
  &\iso E
\end{align*}
Here, $\Delta$ generally denotes the diagonal map $B\to B\times B$, but
$\Delta_{A\times B\times B}$ denotes the diagonal map on $A\times B\times B$.
\end{proof}

\begin{definition}
Let $X$ and $Y$ be $G$-prespectra over $B$, let $W$ be a $G$-prespectrum over $A$,
and let $E$ be a $G$-prespectrum over $A\times B$.
\begin{enumerate}
\item
If $\epsilon\colon Y\exsmsh X\to \Delta_! S_B$ is a $G$-map over $B\times B$, let
\[
 \epsilon_\sharp\colon [W,(1\times\rho)_!(1\times\Delta)^*(E\exsmsh Y)]_{G,A}
  \to [W\exsmsh X, E]_{G,A\times B}
\]
be defined by letting $\epsilon_\sharp(f)$ be the composite
\begin{align*}
 W\exsmsh X
  &\xrightarrow{f\smsh 1} (1\times\rho)_!(1\times\Delta)^*(E\exsmsh Y)\exsmsh X \\
  &\iso (1\times\rho\times 1)_!(1\times\Delta\times 1)^*(E\exsmsh Y\exsmsh X) \\
  &\xrightarrow{1\smsh\epsilon} 
    (1\times\rho\times 1)_!(1\times\Delta\times 1)^*(E\exsmsh \Delta_! S_B) \\
  &\iso E.
\end{align*}
\item
If $\eta\colon S\to \rho_!\Delta^*(X\exsmsh Y) = \rho_!(X\smsh_B Y)$ is a $G$-map, let
\[
 \eta_\sharp\colon [W\exsmsh X, E]_{G,A\times B} \to 
   [W, (1\times\rho)_!(1\times\Delta)^*(E\exsmsh Y)]_{G,A}
\]
be defined by letting $\eta_\sharp(f)$ be the composite
\begin{align*}
 W \iso W\exsmsh S
  &\xrightarrow{1\smsh\eta} W\exsmsh\rho_!\Delta^*(X\exsmsh Y) \\
  &\iso (1\times\rho)_!(1\times\Delta)^*(W\exsmsh X\exsmsh Y) \\
  &\xrightarrow{f\smsh 1} (1\times\rho)_!(1\times\Delta)^*(E\exsmsh Y).
\end{align*}
\end{enumerate}
\end{definition}

\begin{definition}
A $G$-prespectrum $X$ over $B$ is {\em homologically finite}
if
\[
 \epsilon_\sharp\colon [W, (1\times\rho)_!(1\times\Delta)^*(E\exsmsh \bar D_B X)]_{G,A}
  \to [W\exsmsh X, E]_{G,A\times B}
\]
is an isomorphism for all $W$ and $E$,
where $\epsilon\colon \bar D_B X\exsmsh X\to \Delta_! S_B$ is the evaluation map.
\end{definition}

We have the following characterizations of the dual,
similar to \cite[III.1.6]{LMS:eqhomotopy}.

\begin{theorem}\label{thm:dualcharacterization}
Let $X$ and $Y$ be $G$-prespectra over $B$. The following are equivalent.
\begin{enumerate}
\item\label{item:XIsFinite}
$X$ is homologically finite and $Y\iso \bar D_B X$
in the stable category.

\item\label{item:epsilon}
There exists a map $\epsilon\colon Y\exsmsh X\to \Delta_! S_B$
such that
\[
 \epsilon_\sharp\colon [W, (1\times\rho)_!(1\times\Delta)^*(E\exsmsh \bar Y)]_{G,A}
  \to [W\exsmsh X, E]_{G,A\times B}
\]
is an isomorphism for all $W$ and $E$.

\item\label{item:eta}
There exists a map $\eta\colon S\to \rho_!\Delta^*(X\exsmsh Y)$ such that
\[
 \eta_\sharp\colon [W\exsmsh X, E]_{G,A\times B} \to 
   [W, (1\times\rho)_!(1\times\Delta)^*(E\exsmsh Y)]_{G,A}
\]
is an isomorphism for all $W$ and $E$.

\item\label{item:both}
There exist maps $\epsilon\colon Y\exsmsh X\to \Delta_! S_B$ and
$\eta\colon S\to \rho_!\Delta^*(X\exsmsh Y)$ such that
the following two composites are each the identity in the stable category:
\begin{align*}
 X
 &\xrightarrow{\eta\smsh 1} \rho_!\Delta^*(X\exsmsh Y)\exsmsh X \\
 &\iso (\rho\times 1)_!(\Delta\times 1)^*(X\exsmsh Y\exsmsh X) \\
 &\xrightarrow{1\smsh\epsilon} (\rho\times 1)_!(\Delta\times 1)^*(X\exsmsh\Delta_! S_B) \\
 &\iso X
\end{align*}
and
\begin{align*}
 Y
 &\xrightarrow{1\smsh\eta} Y \exsmsh \rho_!\Delta^*(X\exsmsh Y) \\
 &\iso (1\times\rho)_!(1\times\Delta)^*(Y\exsmsh X\exsmsh Y) \\
 &\xrightarrow{\epsilon\smsh 1} (1\times\rho)_!(1\times\Delta)^*(\Delta_! S_B\smsh Y) \\
 &\iso Y
\end{align*}

\end{enumerate}

Moreover, when these conditions are satisfied, the following are also true:
The maps $\epsilon_\sharp$ and $\eta_\sharp$ are inverse isomorphisms;
the adjoint of $\epsilon\colon Y\exsmsh X\to \Delta_! S_B$ is an
equivalence $\tilde\epsilon\colon Y\to \bar D_B X$;
and the maps $\epsilon\gamma$ and $\gamma\eta$ satisfy the same conditions
with the roles of $X$ and $Y$ reversed, where $\gamma$ is the twist map.
\end{theorem}

\begin{proof}
(\ref{item:XIsFinite}) implies (\ref{item:epsilon}) trivially.

Suppose that (\ref{item:both}) is true.
The following diagram shows that $\epsilon_\sharp\eta_\sharp(f) = f$
for any $f\colon W\exsmsh X\to E$:
\[
 \def\objectstyle{\scriptstyle}
 \def\labelstyle{\scriptstyle}
 \xymatrix{
  W\exsmsh X \ar[d]_{1\smsh\eta\smsh 1} \\
  W\exsmsh \rho_!\Delta^*(X\exsmsh Y)\exsmsh X \ar[d]_\iso \\
  (1\times\rho\times 1)_!(1\times\Delta\times 1)^*(W\exsmsh X\exsmsh Y\exsmsh X)
    \ar[r]^{f\smsh 1} \ar[d]_{1\smsh 1\smsh\epsilon}
  & (1\times\rho\times 1)_!(1\times\Delta\times 1)^*(E\exsmsh Y\exsmsh X)
    \ar[d]^{1\smsh\epsilon} \\
  (1\times\rho\times 1)_!(1\times\Delta\times 1)^*(W\exsmsh X\exsmsh\Delta_! S_B)
    \ar[r]^{f\smsh 1} \ar[d]_\iso
  & (1\times\rho\times 1)_!(1\times\Delta\times 1)^*(E\exsmsh\Delta_! S_B) \ar[d]^\iso \\
  W\exsmsh X \ar[r]^{f} & E
 }
\]
On the other hand, the following diagram shows that
$\eta_\sharp\epsilon_\sharp(g) = g$ for any
$g\colon W \to (1\times\rho)_!(1\times\Delta^*)(E\exsmsh Y)$:
\[
 \xymatrix{
  W \ar[r]^-{g} \ar[d]_{1\smsh\eta}
   & (1\times\rho)_!(1\times\Delta)^*(E\exsmsh Y) \ar[d]^{1\smsh\eta} \\
  W\exsmsh\rho_!\Delta^*(X\exsmsh Y) \ar[r]^-{g\smsh 1} \ar[d]_\iso
   & (1\times\rho)_!(1\times\Delta)^*(E\exsmsh Y)\exsmsh\rho_!\Delta^*(X\exsmsh Y) \ar[d]^\iso \\
  (1\times\rho)_!(1\times\Delta)^*(W\exsmsh X\exsmsh Y) \ar[r]^-{g\smsh 1}
   & (1\times\rho\times\rho)_!(1\times\Delta\times\Delta)^*(E\exsmsh Y\exsmsh X\exsmsh Y)
     \ar[d]^{1\smsh\epsilon \smsh 1} \\
   & (1\times\rho\times\rho)_!(1\times\Delta\times\Delta)^*(E\exsmsh\Delta_!S_B\exsmsh Y)
     \ar[d]^\iso \\
   & (1\times\rho)_!(1\times\Delta)^*(E\exsmsh Y)
  }
\]
Here we use the natural isomorphism 
\[
 (1\times\rho)_!(1\times\Delta)^*(1\times\rho\times 1\times 1)_!(1\times\Delta\times 1\times 1)^*
\iso (1\times\rho\times\rho)_!(1\times\Delta\times\Delta)^*, 
\]
which follows from the
isomorphism
$(1\times\Delta)^*(1\times\rho\times 1\times 1)_! 
\iso (1\times\rho\times 1)_!(1\times 1\times\Delta)^*$.
We also use the isomorphism
\[
 (1\times\rho)_!(1\times\Delta)^*(1\times 1\times 1\times\rho)_!(1\times 1\times 1\times\Delta)^*
\iso (1\times\rho\times\rho)_!(1\times\Delta\times\Delta)^*, 
\]
which is true for a similar reason.
Therefore, $\epsilon_\sharp$ and $\eta_\sharp$ are inverse isomorphisms, so
(\ref{item:both}) implies both (\ref{item:epsilon}) and (\ref{item:eta}).

Assume now that (\ref{item:epsilon}) is true.
Taking $W = S$ and $E = X$, let 
\[
 \eta\colon S\to\rho_!\Delta^*(X\exsmsh Y) 
\]
be the
map such that $\epsilon_\sharp(\eta) = 1$, the identity map of $X$.
The first composite in (\ref{item:both}) is just $\epsilon_\sharp(\eta)$,
so is the identity by construction.
The first of the diagrams above shows that $\epsilon_\sharp\circ\eta_\sharp = 1$,
so $\eta_\sharp$ is the inverse isomorphism to $\epsilon_\sharp$.
Now take $W = Y$ and $E = \Delta_! S_B$.
Using the identification $(1\times\rho)_!(1\times\Delta)^*(\Delta_!S_B\exsmsh Y) \iso Y$,
we have $\epsilon_\sharp(1) = \epsilon$, so $\eta_\sharp(\epsilon) = 1$.
Written out, this says that the second composite in (\ref{item:both}) is also
the identity. Thus, (\ref{item:epsilon}) implies (\ref{item:both}).

Now suppose that (\ref{item:eta}) is true.
Take $W = Y$ and $E = \Delta_! S_B$ and
let 
\[
 \epsilon\colon Y\exsmsh X\to \Delta_!S_B 
\]
be the map such that
$\eta_\sharp(\epsilon) = 1$.
The remainder of the proof that (\ref{item:eta}) implies (\ref{item:both})
is now similar to the preceding argument.

Assuming (\ref{item:both}), hence (\ref{item:epsilon}), let
$\tilde\epsilon\colon Y\to \bar D_B X$ be the adjoint of
$\epsilon\colon Y\exsmsh X\to \Delta_! S_B$.
A reasonably straightforward diagram chase shows that $\tilde\epsilon$
is an equivalence with inverse the composite
\begin{align*}
 \bar D_B X
  &\xrightarrow{1\smsh\eta} \bar D_B X\exsmsh \rho_!\Delta^*(X\exsmsh Y) \\
  &\iso (1\times\rho)_!(1\times\Delta)^*(\bar D_B X\exsmsh X\exsmsh Y) \\
  &\xrightarrow{\epsilon\smsh 1} (1\times\rho)_!(1\times\Delta)^*(\Delta_! S_B\exsmsh Y) \\
  &\iso Y.
\end{align*}
(\ref{item:XIsFinite}) now follows. 

The claim in the theorem about $\epsilon\gamma$ and $\gamma\eta$ follows from
(\ref{item:both}) by symmetry. The other claims in the last part of the theorem
have been shown along the way.
\end{proof}

\begin{corollary}
If $X$ is homologically finite then so is $\bar D_B X$, and
$\bar D_B\bar D_B X \iso X$.
\qed
\end{corollary}

We have the following useful compatibility result.

\begin{proposition}\label{prop:pushdual}
If $X$ is a homologically finite $G$-prespectrum over $B$ and
$\beta\colon B\to B'$ is a $G$-map, then $\beta_! X$ is homologically
finite and $\beta_!\bar D_B X \iso \bar D_{B'}\beta_! X$.
\end{proposition}

\begin{proof}
Let $E$ be a $G$-prespectrum over $A\times B'$.
We first show that, for any $G$-prespectrum $Y$ over $B$, we have an isomorphism
\[
 (1\times\rho)_!(1\times\Delta)^*(E\exsmsh\beta_! Y)
  \iso (1\times\rho)_!(1\times\Delta)^*[(1\times\beta)^*E\exsmsh Y].
\]
Once more, we need to detour carefully around commutation relations we wish
were true but need not be in general. The following chain of isomorphisms works:
\begin{align*}
 (1\times\rho)_!(1\times\Delta)^*(E\exsmsh \beta_! Y)
  &\iso (1\times\rho)_![E\smsh_{A\times B'}(\rho\times 1)^*\beta_! Y] \\
  &\iso (1\times\rho)_![E\smsh_{A\times B'}(1\times\beta)_!(\rho\times 1)^*Y] \\
  &\iso (1\times\rho)_!(1\times\beta)_![(1\times\beta)^*E\smsh_{A\times B}(\rho\times 1)^*Y] \\
  &\iso (1\times\rho)_!(1\times\Delta)^*[(1\times\beta)^*E\exsmsh Y]
\end{align*}
Let $W$ be a $G$-prespectrum over $A$. We now have the following chain of isomorphisms:
\begin{align*}
 [W,(1\times\rho)_!(1\times\Delta)^*&(E\exsmsh\beta_! \bar D_B X)]_{G,A} \\
 &\iso [W,(1\times\rho)_!(1\times\Delta)^*[(1\times\beta)^*E\exsmsh \bar D_B X]]_{G,A} \\
 &\iso [W\exsmsh X, (1\times\beta)^*E]_{G,A\times B} \\
 &\iso [W\exsmsh \beta_! X, E]_{G,A\times B'}
\end{align*}
Tracing through the adjunctions, we can see that this is $\epsilon'_\sharp$ where
$\epsilon'$ is the composite
\begin{align*}
 \beta_!Y\exsmsh\beta_!X &\iso (\beta\times\beta)_!(Y\exsmsh X) \\
  &\xrightarrow{\epsilon} (\beta\times\beta)_!\Delta_! S_B \\
  &\iso \Delta_!\beta_!S_B \\
  &\to \Delta_!S_{B'},
\end{align*}
where the last map is induced by the adjoint of the isomorphism
$S_B \xrightarrow{\iso} \beta^*S_{B'}$.
The proposition now follows from
Theorem~\ref{thm:dualcharacterization}.
\end{proof}

In particular, if $X$ is homologically finite over $B$,
then $\rho_!X$ is a finite nonparametrized spectrum and
$D\rho_! X \iso \rho_!\bar D_B X$.
Thus, the homological dual of $X$ is essentially the ordinary dual of the
underlying nonparametrized spectrum $\rho_! X$, except that it is constructed
so as to be parametrized by $B$.
For example, we have the following sequence of results, giving examples
of duals of many spaces.
(Our original proofs were modeled on
\cite[\S III.4]{LMS:eqhomotopy} and were rather involved.
Peter informs us that he no longer trusts the argument of that section
and we're not confident that our arguments were correct, either.
Given that May and Sigurdsson have since provided better proofs,
we shall defer to them for the proofs of the following results.)

We first introduce the following notion for spaces.

\begin{definition}
Let $X$ and $Y$ be ex-$G$-spaces over $B$.
We say that $X$ and $Y$ are {\em $V$-dual} if there exist
weak maps
\[
 \epsilon\colon Y\exsmsh X\to \Delta_! S_B \smsh S^V
\]
and
\[
 \eta\colon S^V\to \rho_!(X\smsh_B Y)
\]
such that the following diagrams stably commute:
\[
 \xymatrix{
  S^V\smsh X \ar[r]^-{\eta\smsh 1} \ar[d]_\gamma
   & (\rho\times 1)_![(1\times\rho)^*X\smsh_{B\times B}(Y\exsmsh X)]
      \ar[d]^{1\smsh\epsilon} \\
  X\smsh S^V 
   & (\rho\times 1)_![(1\times\rho)^*X\smsh_{B\times B}\Delta_! S_B \smsh S^V]
     \ar[l]^-\iso
  }
\]
and
\[
 \xymatrix{
  Y\smsh S^V \ar[r]^-{1\smsh\eta} \ar[d]_\gamma
   & (1\times\rho)_![(Y\exsmsh X)\smsh_{B\times B}(\rho\times 1)^*Y]
     \ar[dd]^{\epsilon\smsh 1} \\
  S^V\smsh Y \ar[d]_{\sigma\smsh 1} \\
  S^V\smsh Y 
   & (1\times\rho)_![(\Delta_!S_B\smsh S^V)\smsh_{B\times B}(\rho\times 1)^*Y]
     \ar[l]^-\iso
 }
\]
Here, $\gamma$ is transposition and $\sigma(v) = -v$.
\end{definition}

\begin{proposition}\cite[18.3.2]{MaySig:parametrized}
If $X$ and $Y$ are $V$-dual, then $\susp^\infty_B X$ is homologically finite
and $\bar D_B X \iso \susp^{-V}_B\susp^\infty_B Y$.
\qed
\end{proposition}

If $(X,A)\to B$ is any pair of (unbased) spaces over $B$, let $C_B(X,A)$ denote the
fiberwise mapping cone of $A\to X$ over $B$, with section given by the conepoints.
In particular, $C_B(X,\emptyset) = X_+$.

Let $X$ be a compact $G$-ENR and let $p\colon X\to B$. Take an embedding of $X$ as
an equivariant neighborhood retract in a $G$-representation $V$.
Let $r\colon N\to X$ be a $G$-retraction of an open neighborhood of $X$ in $V$.
Via $r$, we think of $N$ as a space over $X$, hence over $B$.

\begin{theorem}\label{thm:GENRduality}\cite[18.5.2]{MaySig:parametrized}
Let $(X,A)$ be a compact $G$-ENR pair over $B$ and let $V$ and $N$ be as above.
Then $C_B(X,A)$ is $V$-dual to $C_B(N-A,N-X)$, hence
$\susp^\infty_B C_B(X,A)$ is homologically finite with 
\[
 \bar D_B\susp^\infty_B C_B(X,A) \iso \susp^{-V}\susp^\infty_B C_B(N-A, N-X).
\]
\qed
\end{theorem}

As usual, we can specialize to smooth compact manifolds.
Let $M$ be a smooth compact $G$-manifold and let $p\colon M\to B$ be a $G$-map.
Take a proper embedding of $M$ in
$V'\times[0,\infty)$ for some $G$-representation $V'$,
so $\bndry M$ embeds in $V'\times 0$ and $M$ meets
$V'\times 0$ transversely.
Write $V = V'\oplus\Real$.
Write $\nu$ for the normal bundle of the embedding of $M$ 
in $V'\times[0,\infty)$ and write
$\bndry\nu$ for the normal bundle of the embedding of $\bndry M$ in $V'$,
so $\bndry\nu = \nu|\bndry M$.
Write $S^\nu$ for the sectioned bundle over $M$ formed by taking the one-point
compactification of each fiber of $\nu$, the section being given by the
compactification points. Let $S^\nu_B = p_!S^\nu$.
In particular, $S^\nu_* = \rho_! S^\nu = T\nu$, the Thom space of $\nu$.
We have the following version of Atiyah duality.

\begin{theorem}\label{thm:manifoldduality}\cite[18.7.2]{MaySig:parametrized}
Let $M$ be a compact $G$-manifold and let $p\colon M\to B$ be a $G$-map.
Embed $M$ in a $G$-representation $V$ as above.
Then $\susp_B^\infty M_+$ is homologically finite, with dual
$\susp^{-V}S^\nu_B\PQ{B}S^{\bndry\nu}_B$,
and $\susp_B^\infty M\PQ{B}\bndry M$ is homologically finite with dual
$\susp^{-V}S^\nu_B$.
The cofibration sequence
\[
 \bndry M_+ \to M_+ \to M\PQ{B}\bndry M \to \susp_B \bndry M_+
\]
is $V$-dual to the cofibration sequence
\[
 \susp_B S^{\bndry\nu}_B \from S^\nu_B\PQ{B}S^{\bndry\nu}_B
  \from S^\nu_B \from S^{\bndry\nu}_B.
\]
\qed
\end{theorem}

Note: In \cite{MaySig:parametrized}, Atiyah duality is shown first
and used to prove the more general duality theorem for $G$-ENRs.

In particular, if $b\colon G/H\to B$ and $W$ is a representation of $H$,
for example, $W = \delta(G/H)$ for a dimension function $\delta$, then we have
\[
 \bar D_B (G_+\smsh_H S^{-W,b}) \iso G_+\smsh_H S^{-(\Lie(G/H)-W),b}.
\]

One last observation:

\begin{proposition}
If $X$ and $Y$ are homologically finite spectra over $B$, then
\[
 [X, Y]_{G,B} \iso [\bar D_B Y, \bar D_B X]_{G,B}.
\]
\end{proposition}

\begin{proof}
We have
\begin{align*}
 [X, Y]_{G,B} &\iso [S, \rho_!(Y\smsh_B \bar D_B X)]_G \\
   &\iso [S, \rho_!(\bar D_B X\smsh_B Y)]_G \\
   &\iso [\bar D_B Y, D_B X]_{G,B}.
\end{align*}
\end{proof}

This allows us to identify the opposites of the stable orbit categories.

\begin{corollary}\label{cor:paramMackeyDuality}
If $\delta$ is a dimension function for $G$ then
\[
 \Op{(\stab\Pi_\delta B)} \iso \stab\Pi_{\Lie-\delta} B.
\]
\qed
\end{corollary}

\chapter{$RO(\Pi B)$-graded Ordinary Homology and Cohomology}
\label{chap:ordone}
%
\section*{Introduction}

As we've already mentioned, $RO(G)$-graded ordinary cohomology
is not adequate to give Poincar\'e duality for $G$-manifolds
except in very restrictive settings.
To get Poincar\'e duality we need to extend to a theory
indexed on representations of $\Pi X$.
(For simplicity of notation we shall now write
$\Pi X$ for $\Pi_G X$.)
That is, the homology and cohomology of $X$ should be graded
on representations of $\Pi X$.
A construction of the $RO(\Pi X)$-graded theory for finite $G$
was given in \cite{CW:duality}, and the theory was used in
\cite{CW:spivaknormal} and \cite{CW:simpleduality}
to obtain $\pi$-$\pi$ theorems for
equivariant Poincar\'e duality spaces and equivariant simple Poincar\'e duality spaces.
In this chapter we give a construction for all compact Lie groups.

The construction uses CW spaces in which the cells are
modeled on disks of varying representations, as specified by a 
virtual representation of $\Pi X$.
It will be convenient to have the
representation (and the coefficient system)
carried by a parametrizing space. Thus, we work in the category
of spaces parametrized by a fixed basespace $B$ and
grade our theories on $RO(\Pi B)$.

When deciding what to use as the grading group, there are several possibilities.
We've chosen to work with the most general possible, the group virtual representations of
$\Pi B$. Another possibility would be to work with the Grothendieck group constructed
from the orthogonal representations, which is isomorphic when $B$ is compact.
However, for noncompact $B$, there can be virtual representations that cannot be expressed
as differences of actual representations.

As much as possible, this chapter parallels Chapter~1. In many cases we will
refer to results from that chapter or point out how proofs given there can be
generalized to the present context, rather than repeating them in full.

\section{$G$-CW($\gamma$) complexes}

Fix a compactly generated $G$-space $B$ of the homotopy type of a $G$-CW complex,
a virtual representation $\gamma$ of $\Pi B$, and
a dimension function $\delta$ for $G$.
In this section we describe a theory of CW complexes in $G\ParaU B$,
the cells being locally modeled on $\gamma$, with cells having dimensions given by $\delta$.

If $\gamma$ is a virtual representation of $\Pi B$ and $b\colon G/H\to B$ is an orbit over $B$,
write $\gamma_0(b)$ for the virtual representation of $H$ given by restricting $\gamma(b)$
to $eH$.

\begin{definition}
Let $\gamma$ be a virtual representation of $\Pi B$ and
let $\delta$ be a dimension function for $G$.
\begin{enumerate}
\item
An orbit $b\colon G/H\to B$ is {\em $\delta$-$\gamma$-admissible}, or simply {\em admissible},
if $H\in\F(\delta)$ and $\gamma_0(b)-\delta(G/H)+n$ is stably equivalent to an actual $H$-representation
for some integer $n$.

\item
A {\em $\delta$-$\gamma$-cell} is a pair of objects 
$c = (G\times_H \bar D(Z),p)$ in $G\ParaU B$ where
$Z$ is an actual representation of $H$ such that
$Z$ is stably equivalent to
\[
 \gamma_0(p|G\times_H 0)-\delta(G/H)+n
\]
for some integer $n$.
The {\em dimension} of $c$ is $\gamma+n$.
The {\em boundary} of $c$ is $\bndry c = (G\times_H S(V),p)$.
\end{enumerate}
\end{definition}

Note that $p\colon G\times_H D(Z)\to B$ can be any $G$-map in the definition above.
By definition, $p|G\times_H 0$ must be admissible.

Because $G\ParaU B$ has pushouts, we can now speak of attaching a cell $c$ to a space $X$ over $B$ along an attaching map $\bndry c\to X$ over $B$.

\begin{definition}
\begin{enumerate}\item[]
\item
A {\em $\delta$-$G$-CW($\gamma$) complex} is a $G$-space $(X,p)$ over $B$
together with a decomposition
\[
 (X,p) = \colim_n (X^n,p^n)
\]
in $G\ParaU B$, where 
\begin{enumerate}
\item
$X^0$ is a union of $(\gamma-|\gamma|)$-dimensional 
$\delta$-$\gamma$-cells, which is to say, a union of orbits $(G/H,p)$ over $B$
such that $\delta(G/H) = 0$ and $H$ acts trivially on $\gamma_0(G/H,p)$, and
\item
each
$(X^n,p^n)$, $n>0$, is obtained from $(X^{n-1},p^{n-1})$ by attaching
$(\gamma-|\gamma|+n)$-dimensional $\delta$-$\gamma$-cells.
\end{enumerate}
For notational convenience and to remind ourselves of the role of $\gamma$, we
shall also write $X^{\gamma+n}$ for $X^{|\gamma|+n}$.

\item
A {\em relative} $\delta$-$G$-CW$(\gamma)$ complex is a pair $(X,A,p)$ over $B$
where $(X,p) = \colim_n (X^n,p^n)$, $X^0$ is the disjoint union of
$A$ with $(\gamma-|\gamma|)$-dimensional $\gamma$-cells as in (a) above,
and cells are attached as in (b).

\item
An {\em ex-$\delta$-$G$-CW($\gamma$) complex} is an ex-$G$-space $(X,p,\sigma)$
with a relative $\delta$-$G$-CW($\gamma$) structure on $(X,\sigma(B),p)$.

\item
If we allow cells of any dimension to be attached at each stage,
we get the weaker notions of absolute, relative, or ex- {\em $\delta$-$\gamma$-cell complex}.

\item
If $X$ is a $\delta$-$G$-CW($\gamma$) complex or $\delta$-$\gamma$-cell complex with cells
only of dimension less than or equal to $\gamma+n$, we say that
$X$ is {\em $(\gamma+n)$-dimensional}.
\end{enumerate}
\end{definition}

As in the nonparametrized case, we give names to the complexes we get for particular
cases of $\delta$:
An {\em (ordinary) $G$-CW($\gamma$) complex} is a $0$-$G$-CW($\gamma$) complex, i.e.,
one with $\delta=0$.
A {\em dual $G$-CW($\gamma$) complex} is an $\Lie$-$G$-CW($\gamma$) complex.

\begin{examples}\label{ex:gammacomplexes}
 \begin{enumerate}\item[]

\item If $\alpha$ is a virtual representation of $G$ and
$X$ is a (nonparametrized) $\delta$-$G$-CW($\alpha$) complex, as in Chapter~1,
then any map $X\to B$ gives $X$ the structure of a $\delta$-$G$-CW($\bar\alpha$)
complex over $B$, where $\bar\alpha$ is the constant representation at $\alpha$.

\item If $p\colon E \to B$ is a $G$-vector bundle, let $\rho$ be
the associated representation of $\Pi B$ as in \ref{prop:bundlereps}.
Let $T_B(p)$ be the fiberwise one-point compactification of the bundle,
with basepoint section $\sigma\colon B\to T_B(p)$, so $(T_B(p),p,\sigma)$ is an
ex-$G$-space.
If $B$ is an ordinary $G$-CW complex,
then $T_B(p)$ has an evident structure as an
ex-$G$-CW($\rho$) complex (taking $\delta=0$), with cells in one-to-one correspondence
with the cells of $B$. 
More generally, if $B$ has a $\delta$-$G$-CW($\gamma$) structure over itself,
then $T_B(p)$ is an ex-$\delta$-$G$-CW($\gamma+\rho$) complex over $B$.
This is the geometry underlying the Thom isomorphism.
We use $T_B(p)$ as our
model for the Thom space of $p$ in this context because, unlike the 
actual Thom space (which is $T(p) = T_B(p)/\sigma(B)$), it
comes with a natural parametrization, hence a sensible way to consider it as
a $\delta$-$G$-CW($\gamma+\rho$) complex.

\item If $M$ is any smooth $G$-manifold, let $\tau$ be the tangent
representation of $\Pi M$.
Then the dual of a smooth $G$-triangulation
(as in Example~\ref{ex:vcomplexes}) gives an explicit dual $G$-CW($\tau$) structure on
$M$, considering $M$ as a space over itself.
This is the geometry underlying Poincar\'e duality.

 \end{enumerate}
 \end{examples}

\begin{definition}
If $n$ is an integer, a lax map $f\colon X\laxto Y$ is a {\em $\delta$-$(\gamma+n)$-equivalence}
if, for every $\delta$-$\gamma$-cell $c$ of dimension $\gamma+i$ with $i\leq n$,
every diagram of the following form is lax homotopic to one in which there exists a lift:
\[
 \xymatrix{
  {\bndry c} \ar@{=>}[r] \ar[d] & X \ar@{=>}[d]^f \\
  c \ar@{=>}[r] & Y.
 }
\]
We say that $f$ is a {\em $\delta$-weak$_{\gamma}$ equivalence}
if it a $\delta$-$(\gamma+n)$-equivalence for all $n$.
\end{definition}

Note that, because $\bndry c\to c$ is a lax cofibration
by Theorem~\ref{thm:laxcofib}, we can say that
$f$ is a $\delta$-$(\gamma+n)$-equivalence if every diagram as above is
homotopic {\em rel $\bndry c$} to one in which we can find a lift.

We insert the following reassuring result.

\begin{theorem}\label{thm:gammaequivalence}
Let $f\colon X\laxto Y$ be a lax map.
Then $f$ is a $\delta$-weak$_\gamma$ equivalence if and only if,
for each $b\in B$, $f_b\colon LX_b\to LY_b$ is a
$\delta$-weak$_{\gamma_0(b)}$ equivalence.
\end{theorem}

\begin{proof}
A $\delta$-$\gamma$-cell $c = (G\times_H \bar D(Z),p)$ with $p(eH\times 0) = b$ is lax homotopy equivalent
to the cell 
\[
 c' = (G\times_H\bar D(Z),p')
\]
of the same dimension with 
$p'(G\times_H\bar D(Z)) = Gb$.
It follows that $f$ is a $\delta$-weak$_\gamma$ equivalence if and only if
each map of fibers $f_b\colon LX_b\to LY_b$ is a
$\delta$-weak$_{\gamma_0(b)}$ equivalence.
\end{proof}

We have the following variant of the ``homotopy extension and
lifting property''.

\begin{lemma}[H.E.L.P.]\label{lem:HELPParaSpace}
 Let $r\colon Y \laxto Z$ be a $\delta$-$(\gamma+n)$-equivalence.
Let $(X,A)$ be a relative
$\delta$-$\gamma$-cell complex
of dimension
$\gamma+n$. If the following diagram commutes in $G\LaxU B$ without the
dashed arrows, then there exist lax maps $\tilde g$ and $\tilde h$
making the diagram commute.
 \[
 \xymatrix{
  A \ar[rr]^-{i_0} \ar[dd] & & A\times I \ar'[d][dd] \ar@{=>}[dl]_(.6)h & &
  A \ar[ll]_-{i_1} \ar@{=>}[dl]_(.6)g \ar[dd]  \\
   & Z & & Y \ar@{=>}[ll]_(.3)r \\
  X \ar@{=>}[ur]^f \ar@{-)}[rr]_-{i_0} & & X\times I \ar@{==>}[ul]^(.6){\tilde h} & &
   X \ar[ll]^-{i_1} \ar@{==>}[ul]^(.6){\tilde g} \\
 }
 \]
The result remains true when $n = \infty$.
 \end{lemma}

\begin{proof}
We can use the same proof as that given for
Lemma~\ref{lem:HELP}, understanding all maps to be in the lax category.
\end{proof}

To state the Whitehead theorem, we need the following notations.
Write $\pi(G\LaxU B)(X,Y)$ for the set of lax homotopy classes
of maps $X\laxto Y$.
If $A$ is a space over $B$, write $A/G\LaxU B$ for the category of
spaces under $A$ in $G\LaxU B$, meaning spaces $Y$ over $B$ equipped
with lax maps $A\laxto Y$.
We then write $\pi(A/G\LaxU B)(X,Y)$ for the set of lax homotopy classes of
maps under $A$ and over $B$.

\begin{theorem}[Whitehead]\label{thm:ParaWhitehead}
 \begin{enumerate}\item[]

\item If $f\colon Y\laxto Z$ is a $\delta$-$(\gamma+n)$-equivalence and $X$ is a
$(\gamma+n-1)$-dimensional
$\delta$-$\gamma$-cell complex,
then
\[
 f_*\colon \pi(G\LaxU B)(X,Y) \to \pi(G\LaxU B)(X,Z)
\]
is an isomorphism.
It is an epimorphism if $X$ is
$(\gamma+n)$-dimensional.

\item If $f\colon Y\laxto Z$ is a $\delta$-weak$_\gamma$ equivalence and $X$ is a
$\delta$-$\gamma$-cell complex,
then
\[
 f_*\colon \pi(G\LaxU B)(X,Y) \to \pi(G\LaxU B)(X,Z)
\]
is an
isomorphism.  In particular any $\delta$-weak$_\gamma$ equivalence of $\delta$-$\gamma$-cell complexes
is a lax homotopy equivalence.

\item If $f\colon Y\laxto Z$ is a $\delta$-$(\gamma+n)$-equivalence of spaces under $A$ and $(X,A)$ is a
$(\gamma+n-1)$-dimensional relative
$\delta$-$\gamma$-cell complex, then
\[
 f_*\colon \pi(A/G\LaxU B)(X,Y) \to \pi(A/G\LaxU B)(X,Z)
\]
is an isomorphism.
It is an epimorphism if $(X,A)$ is
$(\gamma+n)$-dimensional.

\item If $f\colon Y\laxto Z$ is a $\delta$-weak$_\gamma$ equivalence of spaces under $A$ and $(X,A)$ is a
relative $\delta$-$\gamma$-cell complex, then 
\[
 f_*\colon \pi(A/G\LaxU B)(X,Y) \to \pi(A/G\LaxU B)(X,Z)
\]
is an
isomorphism.

\end{enumerate}
\end{theorem}

\begin{proof}
This follows from the H.E.L.P. lemma in exactly the same way as
in Theorems~\ref{thm:VWhitehead} and~\ref{thm:relVWhitehead}.
 \end{proof}
 
Note that, when we take $A = B$ in the relative Whitehead theorem, we get the special
case of ex-$G$-spaces and ex-$\delta$-$\gamma$-cell complexes.

Although $\delta$-weak$_\gamma$ equivalence of parametrized spaces is, in general,
weaker than weak $G$-equivalence, for $\delta$-$\gamma$-cell complexes these notions coincide
and are equivalent to lax homotopy equivalence.
The point is that $\delta$-$\gamma$-cell complexes are limited in what
cells they can use, and this limitation
allows $\delta$-weak$_\gamma$ equivalence to detect weak $G$-equivalence
of $\delta$-$\gamma$-cell complexes.

As in the nonparametrized case, we have the following example of a $\delta$-$(\gamma+n)$-equivalence.

\begin{proposition}\label{prop:gammaSkelEquivalence}
Let $\delta$ be a familial dimension function and let $\epsilon$ be another dimenions function
(not necessarily familial) with $\delta\dimpred\epsilon$.
If $(X,A)$ is a relative $\epsilon$-$G$-CW($\gamma$) complex, then the inclusion
$X^{\gamma+n}\to X$ is a $\delta$-$(\gamma+n)$-equivalence.
\end{proposition}

\begin{proof}
By the usual induction, this reduces to showing that, if $Y$ is obtained from $C$ by attaching a $(\gamma+k)$-dimensional $\epsilon$-$\gamma$-cell
of the form $(G\times_K \bar D(W),q)$,
and $(G\times_H \bar D(V),p)$ is a $(\gamma+i)$-dimensional $\delta$-$\gamma$-cell,
with $i<k$,
then any lax $G$-map of pairs $G\times_H \bar D(V)\to (Y,C)$
is homotopic rel boundary to a map into $B$.
For this we take an ordinary $G$-triangulation of $G\times_H D(V)$
and show by induction on the cells that we can (lax) homotope the map
to miss the orbit $G/K\times 0 = (G/K,b')$ in the attached cell.
We use the fact that the fixed-set dimensions of $V$ must
equal those of $\gamma_0(b)-\delta(G/H)+i$, where $(G/H,b)$ is the center of $G\times_H D(V)$, and
the fixed set dimensions of $W$ equal those of
$\gamma_0(b')-\epsilon(G/K)+k$.

The only simplices that might hit the orbit have the form
$G/J\times\Delta^j$ where $J\leq H$ and $J$ is subconjugate to $K$;
by replacing $K$ with a conjugate we may assume $J\leq K$ also.
If such a simplex meets the orbit, then we get a lax map
$(G/J,b)\laxto (G/K,b')$, where $G/J$ is any orbit in $G/J\times \Delta^J$.
Applying $\gamma$ we get a virtual equivalence of $J$-representations
$\gamma_0(b) \hmtpc \gamma_0(b')$, the import of which is that the $J$-fixed points
of these representations have the same dimension.
Because such a simplex is embedded in $G\times_H D(V)$, we now have
\begin{align*}
  j &\leq |[\gamma_0(b) - \delta(G/H) - \Lie(H/J)]^J| + i \\
    &\leq |[\gamma_0(b) - \delta(G/H) - \delta(H/J)]^J| + i \\
    &= |[\gamma_0(b) - \delta(G/J)]^J| + i \\
    &\leq |[\gamma_0(b') - \epsilon(G/K)]^J| + i \\
    &< |[\gamma_0(b') - \epsilon(G/K)]^J| + k
\end{align*}
as in the proof of Proposition~\ref{prop:SkelEquivalence}.
Finally, note that $|[\gamma_0(b') - \epsilon(G/K)]^J| + k$ is the codimension of
$(G/K\times 0)^J$ in $Y^J$,
so the desired homotopy exists for dimensional reasons.
\end{proof}

\begin{definition}
A lax map $f\colon (X,A)\laxto (Y,C)$ from a relative $\delta$-$G$-CW($\gamma$) complex
to a relative $\epsilon$-$G$-CW($\gamma$) complex is {\em cellular} if
$f(X^{\gamma+n})\subset Y^{\gamma+n}$ for each~$n$.
 \end{definition}

\begin{theorem}[Cellular Approximation of Maps]\label{thm:cellapprox}
Let $\delta$ be a familial dimension function and let $\epsilon$ be any dimension function
such that $\delta\dimpred\epsilon$.
Let $(X,A)$ be a relative $\delta$-$G$-CW($\gamma$) complex
and let $(Y,C)$ be a relative $\epsilon$-$G$-CW($\gamma$) complex.
Let $f\colon (X,A)\laxto (Y,C)$ be a lax map and suppose given
a subcomplex $(Z,A)\subset (X,A)$ and a lax homotopy $h$ of $f|Z$ to a cellular map. 
Then $h$ can be extended to a lax homotopy of $f$ to a cellular map.
 \end{theorem}

\begin{proof}
 This follows by induction on skeleta, using the H.E.L.P. lemma and Proposition~\ref{prop:gammaSkelEquivalence}
applied to the inclusion $Y^{\gamma+n}\to Y$.
\end{proof}

\begin{theorem}\label{thm:generalGammaCellularApprox}
Let $\delta$ be a familial dimension function,
let $(A,P)$ be a relative $\delta$-$G$-CW($\gamma$) complex, 
let $X$ be a $G$-space over $B$,
and let $f\colon A\to X$ be a (strict) map.
Then there exists a relative $\delta$-$G$-CW($\gamma$) complex $(Y,P)$, 
containing $(A,P)$ as a subcomplex,
and a $\delta$-weak$_\gamma$ equivalence $g\colon Y\to X$ extending $f$.
\end{theorem}

\begin{proof}
The proof is similar to that of Theorem~\ref{thm:generalCellularApprox}.
We start by letting
\[
 Y^{\gamma-|\gamma|} = A^{\gamma-|\gamma|}\disjunion \coprod (G/K,b),
\]
where the coproduct runs over all $b\colon G/K\to B$ such that $\delta(G/K) = 0$ and
$K$ acts trivially on $\gamma(b)$, and all maps $G/K\to X$ over $B$. The map
$g\colon Y^{\gamma-|\gamma|}\to X$ is the one induced by those maps of orbits.
$(Y^{\gamma-|\gamma|},P)$ is then a relative $\delta$-$G$-CW($\gamma$) complex 
of dimension $\gamma-|\gamma|$,
containing $(A^{\gamma-|\gamma|},P)$, and $Y^{\gamma-|\gamma|}\to X$ is a
$\delta$-$(\gamma-|\gamma|)$-equivalence.

Inductively, suppose that we have constructed $(Y^{\gamma+n-1},P)$, a relative
$\delta$-$G$-CW($\gamma$) complex of dimension $\gamma+n-1$ containing 
$(A^{\gamma+n-1},P)$ as a subcomplex,
and a $\delta$-$(\gamma+n-1)$-equivalence $g\colon Y^{\gamma+n-1}\to X$ 
extending $f$ on $A^{\gamma+n-1}$.
Let
\[
 Y^{\gamma+n} = Y^{\gamma+n-1}\union A^{\gamma+n}\union \coprod G\times_K D(V),
\]
where the coproduct runs over all $(|\gamma|+n)$-dimensional
$\delta$-$\gamma$-cells $G\times_K D(V)$ and all diagrams of the form
\[
 \xymatrix{
  G\times_K S(V) \ar[r] \ar[d] & Y^{\gamma+n-1} \ar[d] \\
  G\times_K D(V) \ar[r] & X
 }
\]
over $B$.
The union that defines $Y^{\gamma+n}$ is along $A^{\gamma+n-1}\to Y^{\gamma+n-1}$ and the maps
$G\times_K S(V)\to Y^{\gamma+n-1}$ displayed above.
By construction, $(Y^{\gamma+n},P)$ is a relative $\delta$-$G$-CW($\gamma$) complex of dimension
$\gamma+n$ containing $(A^{\gamma+n},P)$ as a subcomplex.
We let $g\colon Y^{\gamma+n}\to X$ be the induced map and we claim that $g$
is a $\delta$-$(\gamma+n)$-equivalence.
To see this, consider any diagram of the following form, with $G\times_K D(V)$
being a $\delta$-$\gamma$-cell of dimension $\leq \gamma+n$:
\[
 \xymatrix{
  G\times_K S(V) \ar@{=>}[r]^-\alpha \ar[d] & Y^{\gamma+n} \ar[d] \\
  G\times_K D(V) \ar@{=>}[r]_-\beta & X
 }
\]
Let $b\colon G/K\to B$ be the map of the center of the cell.
By Lemma~\ref{lem:vdiskComplex}, $S(V)$ is $K$-homotopy equivalent to
a $\delta$-$K$-CW($\gamma_0(b)$) complex of dimension strictly less than $\gamma_0(b)+n$. Taking $G\times_K -$ and composing
with the map of the cell into $B$, we can consider this a $\delta$-$G$-CW($\gamma$) structure.
So, by cellular approximation of maps,
the diagram above is lax homotopic to one in which $\alpha$ maps the sphere into $Y^{\gamma+n-1}$.
We can then find a lift of $\beta$ up to lax homotopy using
the inductive hypothesis if $|V| < |\gamma|+n$ or the construction of $Y^{\gamma+n}$ 
if $|V| = |\gamma|+n$.

Finally, $Y = \colim_n Y^{\gamma+n}$ satisfies the claim of the theorem.
\end{proof}

\begin{theorem}[Approximation by $\delta$-$G$-CW($\gamma$) complexes]
Let $\delta$ be a familial dimension function and
let X be a $G$-space over $B$. Then there exists a $\delta$-$G$-CW($\gamma$)
complex $\Gamma X$ and a $\delta$-weak$_\gamma$ equivalence $g\colon \Gamma X\to X$ over $B$.
If $f\colon X\laxto Y$ is a lax $G$-map and $g\colon \Gamma Y\to Y$ is an approximation
of $Y$ by a $\delta$-$G$-CW($\gamma$) complex, then
there exists a lax $G$-map $\Gamma f\colon \Gamma X\laxto \Gamma Y$, unique up to 
lax $G$-homotopy,
such that the following diagram commutes up to lax $G$-homotopy:
 \[
 \xymatrix{
  {\Gamma X} \ar@{=>}[r]^{\Gamma f} \ar[d]_g & {\Gamma Y} \ar[d]^g \\
  X \ar@{=>}[r]_{f} & Y
  }
 \]
 \end{theorem}

\begin{proof}
The existence of $g\colon \Gamma X\to X$ is the special case
of Theorem~\ref{thm:generalGammaCellularApprox} in which we take $A = P = \emptyset$.
The existence and uniqueness of $\Gamma f$ follows from Whitehead's theorem.
 \end{proof}

\begin{theorem}[Approximation of Ex-Spaces]\label{thm:basedGammaApproximation}
Let $\delta$ be a familial dimension function and
let X be an ex-$G$-space over $B$. Then there exists an ex-$\delta$-$G$-CW($\gamma$)
complex $\Gamma X$ and a $\delta$-weak$_\gamma$ equivalence $g\colon \Gamma X\to X$ over and under $B$.
Further, $\Gamma$ is functorial up to lax homotopy under $B$.
 \end{theorem}

\begin{proof}
The existence of $g\colon \Gamma X\to X$ is the special case
of Theorem~\ref{thm:generalGammaCellularApprox} in which we take $A = P = B$.

Given $f\colon X\laxto Y$,
the existence and uniqueness of $\Gamma f\colon \Gamma X\laxto \Gamma Y$ follows from the relative
part of Whitehead's theorem,
which tells us that
\[
 \pi(B/G\LaxU B)(\Gamma X, \Gamma Y) \iso \pi(B/G\LaxU B)(\Gamma X, Y).
\]
 \end{proof}

\begin{theorem}[Approximation of Pairs]\label{thm:pairGammaApproximation}
Let $\delta$ be a familial dimension function and
Let $(X,A)$ be a pair of $G$-spaces over $B$.
Then there exists a pair of $\delta$-$G$-CW($\gamma$) complexes $(\Gamma X, \Gamma A)$
and a pair of $\delta$-weak$_\gamma$ equivalences $g\colon (\Gamma X,\Gamma A)\to (X,A)$ over $B$.
Further, $\Gamma$ is functorial on lax maps of pairs up to lax homotopy.
\end{theorem}

\begin{proof}
Take any approximation $g\colon \Gamma A\to A$, then
apply Theorem~\ref{thm:generalGammaCellularApprox} to $\Gamma A\to X$ 
(taking $P = \emptyset$ in that theorem) to get $\Gamma X$ with $\Gamma A$ as a subcomplex.

Given $f\colon (X,A)\laxto (Y,C)$, we first construct
$\Gamma f\colon \Gamma A\laxto \Gamma C$ using the Whitehead theorem
and then extend to $\Gamma X\laxto \Gamma Y$
using the relative Whitehead theorem
(considering the category of spaces under $\Gamma A$ and lax maps over $B$).
\end{proof}

We say that a triad $(X; A, C)$ over $B$ is excisive if it is excisive as a triad of $G$-spaces,
ignoring $B$.

\begin{theorem}[Approximation of Triads]\label{thm:triadGammaApproximation}
Let $\delta$ be a familial dimension function and
let $(X;A,C)$ be an excisive triad over $B$.
Then there exists a $\delta$-$G$-CW($\gamma$) triad $(\Gamma X; \Gamma A, \Gamma C)$
and a map of triads 
\[
 g\colon (\Gamma X; \Gamma A, \Gamma C)\to (X;A,C)
\]
over $B$ such that each of the maps $\Gamma A \intersect \Gamma C\to A\intersect C$,
$\Gamma A\to A$, $\Gamma C\to C$, and $\Gamma X\to X$ is
a $\delta$-weak$_\gamma$ equivalence.
$\Gamma$ is functorial on lax maps of excisive triads up to lax homotopy.
\end{theorem}

\begin{proof}
Let $D = A\intersect C$. Take a $\delta$-$G$-CW($\gamma$) approximation
$g\colon \Gamma D\to D$. Using Theorem~\ref{thm:generalGammaCellularApprox}, extend to
approximations
\[
 g\colon (\Gamma A, \Gamma D)\to (A,D)
\]
and
\[
 g\colon (\Gamma C, \Gamma D)\to (C,D).
\]
Let $\Gamma X = \Gamma A \union_{\Gamma D} \Gamma C$.
All the statements of the theorem are clear except that the map
$g\colon \Gamma X\to X$ is a $\delta$-weak$_\gamma$ equivalence.
By Theorem~\ref{thm:gammaequivalence}, it suffices to show that $g$
induces a $\delta$-weak$_{\gamma_0(b)}$ equivalence on each homotopy fiber.
This follows from the observation that, for example,
$L(A\union_D C)_b = LA_b \union_{LD_b} LC_b$ and an application of
Theorem~\ref{thm:gammaequivalence} to say that $L\Gamma A_b\to LA_b$ is
a $\delta$-weak$_{\gamma_0(b)}$ equivalence, and similary for the other approximations.
\end{proof}


\section{$G$-CW prespectra}\label{sec:paramCWprespectra}

As in Chapter~1, we will have a technical problem extending the definition of
cellular homology and cohomology from $G$-CW complexes to arbitrary $G$-spaces,
stemming from Example~\ref{ex:dualInstability}.
In the parametrized case we have no analogue of Corollary~\ref{cor:equivalenceStability},
so the use of the CW prespectra we now describe is essential.

If $\V$ is an indexing sequence in a $G$-universe $\U$, recall that
$G\LaxPreSpec{}{B}$ is the category of $G$-prespectra over $B$ indexed on $\V$
and lax maps.
We write $G\LaxPreSpec{\V}{B}$ if we want to make the indexing sequence explicit.
We have the following analogue of Definition~\ref{def:CWprespectra}.

\begin{definition}\label{def:paramCWprespectra}
Let $\V$ be an indexing sequence in a universe $\U$,
let $\delta$ be a dimension function for $G$, and let
$\gamma$ be a virtual representation of $\Pi B$.
\begin{enumerate}
\item
A lax map $f\colon D\laxto E$ in $G\LaxPreSpec{\V}{B}$ is a {\em $\delta$-weak$_\gamma$ equivalence}
if, for each $i$, $f_i\colon D(V_i)\laxto E(V_i)$ is a $\delta$-weak$_{\gamma+V_i}$ equivalence
of $G$-spaces over $B$.

\item
A $G$-prespectrum $D$ in $G\LaxPreSpec{\V}{B}$ is a {\em $\delta$-$G$-CW$(\gamma)$ prespectrum}
if, for each $i$, $D(V_i)$ is an ex-$\delta$-$G$-CW$(\gamma+V_i)$ complex and if each structure map
$\susp^{V_{i}-V_{i-1}} D(V_{i-1}) \to D(V_i)$
is the inclusion of a subcomplex.

\item
A (strict) map $D\to E$ of $\delta$-$G$-CW$(\gamma)$ prespectra is {\em the inclusion of a subcomplex}
if, for each $i$, the map $D(V_i)\to E(V_i)$ is the inclusion of a subcomplex.
We also say simply that $D$ is a subcomplex of $E$.

\item
A lax map $D\laxto E$ of $\delta$-$G$-CW$(\gamma)$ prespectra is {\em cellular} if, for each $i$,
the map $D(V_i)\to E(V_i)$ is cellular.

\item
If $D$ is a $G$-prespectrum in $G\LaxPreSpec{\V}{B}$, a {\em $\delta$-$G$-CW$(\gamma)$ approximation}
of $D$ is a $\delta$-$G$-CW$(\gamma)$ prespectrum $\Gamma^\delta_\alpha D$ and a strict
$\delta$-weak$_\gamma$ equivalence $\Gamma^\delta_\alpha D \to D$.
\end{enumerate}
\end{definition}

Our results on $G$-CW parametrized spaces give quick proofs of the following results.

\begin{lemma}[H.E.L.P.]
Let $r\colon E \laxto F$ be a $\delta$-weak$_\gamma$ equivalence of $G$-pre\-spectra over $B$ and
let $D$ be a $\delta$-$G$-CW$(\gamma)$ prespectrum over $B$ with subcomplex $C$. 
If the following diagram commutes without the
dashed arrows, then there exist maps $\tilde g$ and $\tilde h$
making the diagram commute.
 \[
 \xymatrix{
  C \ar[rr]^-{i_0} \ar[dd] & & C\smsh I_+ \ar'[d][dd] \ar@{=>}[dl]_(.6)h & & 
  C \ar[ll]_-{i_1} \ar@{=>}[dl]_(.6)g \ar[dd]  \\
   & F & & E \ar@{=>}[ll]_(.3)r \\
  D \ar@{=>}[ur]^f \ar[rr]_-{i_0} & & D\smsh I_+ \ar@{==>}[ul]^(.6){\tilde h} & &
   D \ar[ll]^-{i_1} \ar@{==>}[ul]^(.6){\tilde g} \\
 }
 \]
\end{lemma}

\begin{proof}
We construct $\tilde g$ and $\tilde h$ inductively on the indexing space $V_i$.
For $i=1$ we simply quote the space-level H.E.L.P. lemma to find
$\tilde g_1$ and $\tilde h_1$.

For the inductive step, we assume that we've constructed $\tilde g_{i-1}$ and
$\tilde h_{i-1}$. We then apply the space-level H.E.L.P. lemma
with (using the notation of Lemma~\ref{lem:HELPParaSpace}) $Y = E_i$, $Z = F_i$,
$X = D_i$, and $A = \susp X_{i-1}\union C_i$.
\end{proof}

\begin{proposition}[Whitehead]
Suppose that $D$ is a $\delta$-$G$-CW$(\gamma)$ prespectrum over $B$ and that
$f\colon E\laxto F$ is a $\delta$-weak$_\gamma$ equivalence. Then
\[
 f_*\colon \pi G\LaxPreSpec{\V}{B}(D,E) \to \pi G\LaxPreSpec{\V}{B}(D,F)
\]
is an isomorphism, where $\pi G\LaxPreSpec{\V}{B}(-,-)$ denotes the group of
lax homotopy classes of lax $G$-maps.
Therefore, any $\delta$-weak$_\gamma$ equivalence of $\delta$-$G$-CW$(\gamma)$
prespectra over $B$ is a lax $G$-homotopy equivalence.
\end{proposition}

\begin{proof}
We get surjectivity by applying the H.E.L.P. lemma to $D$ and its subcomplex $*$.
We get injectivity by applying it to $D\smsh I_+$ and its subcomplex
$D\smsh \bndry I_+$.
\end{proof}

\begin{proposition}[Cellular Approximation of Maps]
Suppose that $f\colon D\laxto E$ is a lax map of $\delta$-$G$-CW$(\gamma)$ prespectra over $B$,
$C$ is a subcomplex of $D$, and $h$ is a lax $G$-homotopy of $f|C$ to a cellular map.
Then $h$ can be extended to a lax $G$-homotopy of $f$ to a cellular map.
\end{proposition}

\begin{proof}
This follows by induction on the indexing space $V_i$. The first case to consider is
$f_1\colon D(V_1)\laxto E(V_1)$, and we know from the space-level result that 
we can extend $h_1$ to a lax $G$-homotopy $k_1$ from $f_1$ to a cellular map $g_1$.
For the inductive step, assume we have a lax homotopy $k_{i-1}$,
extending $h_{i-1}$, from $f_{i-1}$ to a cellular
map $g_{i-1}$. Then $\susp^{V_i-V_{i-1}}k_{i-1} \union h_i$ is a lax homotopy on the subcomplex
$\susp^{V_i-V_{i-1}}D(V_{i-1})\union C_i$ of $D(V_i)$. By the space-level result again, we can
extend to a lax homotopy $k_i$ on $D(V_i)$ from $f_i$ to a cellular map $g_i$.
\end{proof}

\begin{proposition}[Cellular Approximation of Prespectra]\label{prop:paramApproxPrespectra}
If $D$ is a $G$-prespectrum in $G\LaxPreSpec{\V}{B}$, then there exists a
$\delta$-$G$-CW$(\gamma)$ approximation $\Gamma D\to D$.
If $f\colon D\laxto E$ is a lax map of $G$-prespectra over $B$ and
$\Gamma E\to E$ is an approximation of $E$,
then there exists a lax cellular map
$\Gamma f\colon \Gamma D \laxto \Gamma E$,
unique up to lax cellular homotopy, making the following diagram lax homotopy commute:
\[
 \xymatrix{
  \Gamma D \ar@{=>}[r]^-{\Gamma f} \ar[d]
   & \Gamma E \ar[d] \\
  D \ar@{=>}[r]_-f & E
 }
\]
\end{proposition}

\begin{proof}
We construct $\Gamma D$ recursively on the indexing space $V_i$.
For $i=1$, we take $(\Gamma D)(V_1) = \Gamma(D(V_1))$ to be any
$\delta$-$G$-CW$(\gamma+V_1)$ approximation of $D(V_1)$.
Suppose that we have constructed $\Gamma D(V_{i-1}) \to D(V_{i-1})$, a
$\delta$-$G$-CW$(\gamma+V_{i-1})$ approximation.
Then $\susp^{V_i-V_{i-1}}\Gamma D(V_{i-1})$ is a $\delta$-$G$-CW$(\gamma+V_i)$ complex
and, by Theorem~\ref{thm:generalGammaCellularApprox}, we can find a
$\delta$-$G$-CW$(\gamma+V_i)$ approximation
$\Gamma D(V_i) \to D(V_i)$ making the following diagram commute, in which the
map $\sigma$ at the top is the inclusion of a subcomplex:
\[
 \xymatrix{
  \susp^{V_i-V_{i-1}}\Gamma D(V_{i-1}) \ar[r]^-\sigma \ar[d]
   & \Gamma D(V_i) \ar[d] \\
  \susp^{V_i-V_{i-1}} D(V_{i-1}) \ar[r]_-\sigma & D(V_i)
 }
\]
Thus, $\Gamma D$ is a $\delta$-$G$-CW$(\gamma)$ prespectrum and the map
$\Gamma D\to D$ so constructed is a $\delta$-weak$_\gamma$ equivalence.

The existence and uniqueness of $\Gamma f$ follow from the Whitehead theorem
and cellular approximation of maps and homotopies.
\end{proof}

As usual, these results imply that we can invert the $\delta$-weak$_\gamma$ equivalences
of prespectra and that the result is equivalent to the ordinary lax homotopy category
of $\delta$-$G$-CW$(\gamma)$ prespectra.
We emphasize, once again, that this is not the parametrized stable category.


\section{Cellular homology and cohomology}\label{sec:ordinarytheories}

We shall now construct the cellular chains of $\delta$-$G$-CW($\gamma$) spaces.
As usual, we work with a fixed basespace $B$, a virtual representation $\gamma$
of $\Pi B$, and a dimension function $\delta$ for $G$.
Referring to Definition~\ref{def:indres}, our coefficient systems will be
$\stab\Pi_{\delta}B$-modules.
(We might call these ``$G$-$\delta$-Mackey functors over $B$,'' but
that's a bit long to repeat as many times as we will need to refer to them.
Recall that we write $\stab\Pi_{G,\delta}B$ when we need to specify the group.)
Such modules can be either contravariant or covariant and, as previously,
we adopt the convention that contravariant modules are written with a bar
on top and covariant ones with a bar underneath.

We have the following generalization of Corollary~\ref{cor:MackeyDuality}.

\begin{corollary}
If $\delta$ is a dimension function for $G$,
then the category of contravariant $\stab\Pi_{\delta}B$-modules and
natural transformations between them is isomorphic to the category
of covariant $\stab\Pi_{\Lie-\delta}B$-modules.
\end{corollary}

\begin{proof}
This is an immediate consequence of Corollary~\ref{cor:paramMackeyDuality}.
\end{proof}

If $(X,A)$ is a pair of $G$-spaces over $B$,
we write $X\ParaQuot B A$ for the ex-$G$-space
over $B$ obtained by taking the fiberwise quotient.
This is the same ex-$G$-space as $X_+\ParaQuot B A_+$.
Recall also the notation of Definition~\ref{def:sphereSpec}.

\begin{lemma}\label{lem:skeletalquotient}
Let $(X,A,p)$ be a relative $\delta$-$G$-CW($\gamma$) complex. Then
 \[
 X^{\gamma+n}\ParaQuot B X^{\gamma+n-1} \hmtpc_\lambda 
   \Wedge_x G_+\smsh_H S^{\gamma_0(px)-\delta(G/H)+n,px}
 \]
as ex-$G$-spaces, where the wedge runs through the centers $x$ of
the $(\gamma+n)$-cells of $X$
and the equivalence is lax homotopy equivalence.
(Here we think of $x\colon G/H\to X$ and $px$ denotes the composite
$px\colon G/H\to B$.
When we write $\gamma_0(px)-\delta(G/H)+n$, we mean an actual representation of $H$
stably equivalent to $\gamma_0(px)-\delta(G/H)+n$.)
 \end{lemma}

\begin{proof}
$X^{\gamma+n}\ParaQuot B X^{\gamma+n-1}$ consists of a copy of $B$,
the image of the
section, with a copy of each $(\gamma+n)$-cell of $X$ adjoined via
the composite of its attaching map and the projection to $B$.
The attaching map for a cell, as a map into $B$, may be deformed to be the
constant map at the image of the center of the cell, and then the
assertion is clear.
 \end{proof}

\begin{definition}\label{def:gammachains}
Let $(X,A)$ be a relative $\delta$-$G$-CW($\gamma$) complex.
The {\em cellular chain complex of $(X,A)$},
$\Mackey C_{\gamma+*}(X,A)$ (which we also
write as $\Mackey C_{\gamma+*}^{G,\delta}(X,A)$ when we need to emphasize the group involved or $\delta$),
is the chain complex of contravariant $\stab\Pi_{\delta}B$-modules defined by
 \begin{align*}
 \Mackey C_{\gamma+n}(X,A)(b)
   &= [ \susp^\infty_B G_+\smsh_H S^{\gamma_0(b)-\delta(G/H)+n,b},
           \susp^\infty_B (X^{\gamma+n}\ParaQuot B X^{\gamma+n-1}) ]_{G,B} \\
   &\iso \pi^H_{\gamma_0(b)-\delta(G/H)+n}(\susp^\infty_B (X^{\gamma+n}\ParaQuot B X^{\gamma+n-1})_b)
 \end{align*}
if $b\colon G/H\to B$.
(We are using shorthand notation---we should really be as careful as
in Definition~\ref{def:vchains}.)
This is obviously a contravariant functor on $\stab\Pi_{\delta,\gamma} B$, and we use
the isomorphism of Theorem~\ref{thm:StableMapsOrbitsOverB} to consider it as a functor on
$\stab\Pi_{\delta} B$. Let $d\colon \Mackey C_{\gamma+n}(X,A)\to \Mackey C_{\gamma+n-1}(X,A)$ be
the natural transformation induced by the composite
 \[
 X^{\gamma+n}\ParaQuot B X^{\gamma+n-1} 
  \to \Susp_B (X^{\gamma+n-1}\ParaQuot B A)
  \to \Susp_B (X^{\gamma+n-1}\ParaQuot B X^{\gamma+n-2}).
 \]
 \end{definition}
 
The analogues of Remarks~\ref{rem:chains} hold here:
$\Mackey C_{\gamma+*}(X,A)$ is covariant in cellular maps $(X,A)\to (Y,B)$
and contravariant in virtual maps $\gamma\to\gamma'$.

If $\MackeyOp S$ is a covariant $\stab\Pi_{\delta}B$-module, and
$\Mackey T$ and $\Mackey U$ are contravariant $\stab\Pi_{\delta}B$-modules,
then Definition~\ref{def:indres} defines for us the groups
 \[
 \Hom_{\stab\Pi_\delta B}(\Mackey T,\Mackey U)
 \]
and
 \[
 \Mackey T \tensor_{\stab\Pi_\delta B} \MackeyOp S.
 \]

Recall that a $\stab\Pi_{\delta}B$-module is said to be
{\em free} if it is a sum of functors of the form $\stab\Pi_\delta B(-,b)$. The
following follows from Lemma~\ref{lem:skeletalquotient}.

\begin{proposition}
$\Mackey C_{\gamma+*}(X)$ is a chain complex of free
$\stab\Pi_{\delta}B$-modules.
If $\Mackey T$ is a contravariant $\stab\Pi_{\delta}B$-module then
 \[
 \Hom_{\stab\Pi_\delta B}(\Mackey C_{\gamma+n}(X),\Mackey T) \iso
     \Dirsum_x \Mackey T(px)
 \]
where the direct sum runs over the centers $x$ of the $(\gamma+n)$-cells of $X$.
Similarly, if $\MackeyOp S$ is a covariant $\stab\Pi_{\delta}B$-module then
 \[
 \Mackey C_{\gamma+n}(X) \tensor_{\stab\Pi_\delta B} \MackeyOp S \iso
     \Dirsum_x \MackeyOp S(px)
 \]
where again $x$ runs over the centers of the $(\gamma+n)$-cells of $X$.
 \qed \end{proposition}

\begin{definition}\label{def:gammaHomology}
 Let $\Mackey T$ be a contravariant $\stab\Pi_{\delta}B$-module and let  
$\MackeyOp S$ be a covariant $\stab\Pi_{\delta}B$-module.
\begin{enumerate}
\item
Let $(X,A)$ be a relative $\delta$-$G$-CW($\gamma$) complex.
We define the {\em $(\gamma+n)$th cellular homology} of
$(X,A)$, with coefficients in $\MackeyOp S$, to be
 \[
 \CH^{G,\delta}_{\gamma+n}(X,A;\MackeyOp S) =
  H_{\gamma+n}( \Mackey{C}^{G,\delta}_{\gamma+*}(X,A) \tensor_{\stab\Pi_\delta B} \MackeyOp S).
 \]
 and we define the {\em $(\gamma+n)$th cellular cohomology} of $(X,A)$, with
coefficients in $\Mackey T$, to be
 \[
 \CH_{G,\delta}^{\gamma+n} (X,A; \Mackey T) = 
  H^{\gamma+n}(\Hom_{\stab\Pi_\delta B}( \Mackey{C}^{G,\delta}_{\gamma+*}(X,A), \Mackey T)),
 \]
where we introduce the same sign as in Definition~\ref{def:alphaHomology}.
Homology is covariant in cellular maps of $(X,A)$
while cohomology is contravariant in cellular maps of $(X,A)$.

\item
If $X$ is a $\delta$-$G$-CW($\gamma$) complex, so $(X,\emptyset)$ is a relative
$\delta$-$G$-CW($\gamma$) complex, we define
\[
 \CH^{G,\delta}_{\gamma+n}(X; \MackeyOp S) = \CH^{G,\delta}_{\gamma+n}(X,\emptyset; \MackeyOp S)
\]
and
\[
 \CH_{G,\delta}^{\gamma+n}(X; \Mackey T) = \CH_{G,\delta}^{\gamma+n}(X,\emptyset; \Mackey T).
\]

\item
If $X$ is an ex-$\delta$-$G$-CW($\gamma$) complex, so $(X,B)$ is a relative
$\delta$-$G$-CW($\gamma$) complex, we define the {\em reduced} homology and cohomology
of $X$ to be
\[
 \tCH^{G,\delta}_{\gamma+n}(X; \MackeyOp S) = \CH^{G,\delta}_{\gamma+n}(X,B; \MackeyOp S)
\]
and
\[
 \tCH_{G,\delta}^{\gamma+n}(X; \Mackey T) = \CH_{G,\delta}^{\gamma+n}(X,B; \Mackey T).
\]
\end{enumerate}
\end{definition}

\begin{theorem}[Reduced Homology and Cohomology of Complexes]\label{thm:reducedHomologyComplexesOverB}
Let $\delta$ be a dimension function for $G$, let
$\gamma$ be a virtual representation of $\Pi_G B$, and
let $\MackeyOp S$ and $\Mackey T$ be respectively a covariant and a contravariant
$\stab\Pi_{\delta}B$-module. Then the abelian groups
$\tCH^{G,\delta}_\gamma(X;\MackeyOp S)$ and $\tCH_{G,\delta}^\gamma(X;\Mackey T)$ are
respectively covariant and contravariant
functors on the homotopy category of 
ex-$\delta$-$G$-CW($\gamma$) complexes and cellular maps and homotopies.
They are also respectively contravariant and covariant functors 
of $\gamma$.
These functors satisfy the following properties.
\begin{enumerate}

\item (Exactness)
If $A$ is a subcomplex of $X$, then the following sequences are exact:
\[
 \tCH^{G,\delta}_\gamma(A; \MackeyOp S) \to \tCH^{G,\delta}_\gamma(X; \MackeyOp S)
  \to \tCH^{G,\delta}_\gamma(X/_B A; \MackeyOp S)
\]
and
\[
 \tCH_{G,\delta}^\gamma(X/_B A; \Mackey T) \to \tCH_{G,\delta}^\gamma(X; \Mackey T)
  \to \tCH_{G,\delta}^\gamma(A; \Mackey T).
\]

\item (Additivity)
If $X = \Wedge_i X_i$ is a fiberwise wedge of ex-$\delta$-$G$-CW($\gamma$) complexes,
then the inclusions of the wedge summands induce isomorphisms
\[
 \Dirsum_i\tCH^{G,\delta}_\gamma(X_i; \MackeyOp S) \iso \tCH^{G,\delta}_\gamma(X; \MackeyOp S)
\]
and
\[
 \tCH_{G,\delta}^\gamma(X; \Mackey T) \iso \prod_i \tCH_{G,\delta}^\gamma(X_i; \Mackey T).
\]

\item (Suspension)
There are suspension isomorphisms
 \[
  \sigma^V\colon \tCH^{G,\delta}_\alpha(X; \MackeyOp S) 
     \xrightarrow{\iso} \tCH^{G,\delta}_{\alpha + V}(\susp_B^V X; \MackeyOp S)
 \]
and
 \[
 \sigma^V\colon \tCH_{G,\delta}^\gamma(X; \Mackey T) 
     \xrightarrow{\iso} \tCH_{G,\delta}^{\gamma + V}(\susp_B^V X; \Mackey T).
 \]
These isomorphisms satisfy $\sigma^0 = \id$, $\sigma^{W}\circ \sigma^V = \sigma^{V\oplus W}$,
and the following naturality condition: If $\zeta\colon V\to V'$ is an isomorphism,
then the following diagrams commute:
\[
 \xymatrix@C+2em{
  \tCH^{G,\delta}_\gamma(X; \MackeyOp S) \ar[r]^-{\sigma^V} \ar[d]_{\sigma^{V'}}
    & \tCH^{G,\delta}_{\gamma + V}(\susp_B^{V}X; \MackeyOp S)
          \ar[d]^{\tCH_{\id}(\id\smsh\zeta)} \\
  \tCH^{G,\delta}_{\gamma + V'}(\susp_B^{V'}X; \MackeyOp S)
       \ar[r]_-{\tCH_{\id+\zeta}(\id)}
    & \tCH^{G,\delta}_{\gamma + V}(\susp_B^{V'}X; \MackeyOp S)
 }
\]
\[
 \xymatrix@C+2em{
  \tCH_{G,\delta}^\gamma(X; \Mackey T) \ar[r]^-{\sigma^V} \ar[d]_{\sigma^{V'}}
    & \tCH_{G,\delta}^{\gamma + V}(\susp_B^{V}X; \Mackey T)
          \ar[d]^{\tCH^{\id+\zeta}(\id)} \\
  \tCH_{G,\delta}^{\gamma + V'}(\susp_B^{V'}X; \Mackey T)
       \ar[r]_-{\tCH^{\id}(\id\smsh\zeta)}
    & \tCH_{G,\delta}^{\gamma + V'}(\susp_B^{V}X; \Mackey T)
 }
\]

\item (Dimension Axiom)
If $H\in\F(\delta)$, $b\colon G/H\to B$, and $V$ is a representation of $H$ so large that
$\gamma_0(b) + V - \delta(G/H) + n$ is an actual representation, then
there are natural isomorphisms
\[
 \tCH^{G,\delta}_{\gamma+V+k}(G_+\smsh_H S^{\gamma_0(b)+V-\delta(G/H)+n};\MackeyOp S)
  \iso
   \begin{cases}
    \MackeyOp S(b) & \text{if $k = n$} \\
    0 & \text{if $k\neq n$}
   \end{cases}
\]
and
\[
 \tCH_{G,\delta}^{\gamma+V+k}(G_+\smsh_H S^{\gamma_0(b)+V-\delta(G/H)+n};\Mackey T)
  \iso
   \begin{cases}
    \Mackey T(b) & \text{if $k = n$} \\
    0 & \text{if $k\neq n$.}
   \end{cases}
\]

\end{enumerate}
\end{theorem}

The proof is the same as Theorem~\ref{thm:reducedHomologyComplexes}
with appropriate changes.

We now extend to arbitrary $G$-ex-spaces over $B$.
As in the nonparametrized case, we do so by approximating by
$G$-CW prespectra.

\begin{definition}
Let $E$ be a $\delta$-$G$-CW$(\gamma)$ prespectrum over $B$. We define the 
{\em cellular chain complex of $E$} to be the colimit
\[
 \Mackey C_{\gamma+*}^{G,\delta}(E) = \colim_i \Mackey C_{\gamma+V_i+*}(E(V_i)).
\]
If $F$ is an arbitrary prespectrum over $B$ and $\delta$ is familial, we define
\[
 \Mackey C_{\gamma+*}^{G,\delta}(F) = \Mackey C_{\gamma+*}^{G,\delta}(\Gamma F)
\]
where $\Gamma F\to F$ is a $\delta$-$G$-CW$(\gamma)$ approximation of $F$.
If $X$ is an ex-$G$-space over $B$ and $\susp_B^\infty X$ denotes its suspension prespectrum, we define
\[
 \Mackey C_{\gamma+*}^{G,\delta}(X) = \Mackey C_{\gamma+*}^{G,\delta}(\susp_B^\infty X)
  = \Mackey C_{\gamma+*}^{G,\delta}(\Gamma\susp_B^\infty X).
\]
\end{definition}

The analogues of Propositions~\ref{prop:chainsFunctorial}
through \ref{prop:universeIndependence} are true, so that
$\Mackey C_{\gamma+*}^{G,\delta}(X)$ is well-defined up to chain homotopy equivalence and functorial up to chain homotopy, independent of the choices of universe, indexing sequence,
and CW approximation.

We also have the analogue of Definition~\ref{def:mappingCylinders}:

\begin{definition}
\begin{enumerate}
\item[]
\item
Let $f\colon A\to X$ be a map of $G$-spaces over $B$
The {\em (unreduced) mapping cylinder, $Mf$,} is the pushout in the following diagram:
\[
 \xymatrix{
  A \ar[r]^f \ar[d]_{i_0} & X \ar[d] \\
  A\times I \ar[r] & Mf
 }
\]
The {\em (unreduced) mapping cone}
is $Cf = Mf/(A\times 1)$ with basepoint the image of $A\times 1$.

\item
Let $f\colon A\to X$ be a map of  ex-$G$-spaces.
The {\em reduced mapping cylinder, $\tilde Mf$,} is the pushout in the following diagram:
\[
 \xymatrix{
  A \ar[r]^f \ar[d]_{i_0} & X \ar[d] \\
  A\smsh I_+ \ar[r] & \tilde Mf
 }
\]
The {\em reduced mapping cone, $\tilde Cf$,}
is the pushout in the following diagram, in which $I$ has basepoint 1:
\[
 \xymatrix{
  A \ar[r]^f \ar[d]_{i_0} & X \ar[d] \\
  A\smsh I \ar[r] & \tilde Cf
 }
\]
\end{enumerate}
In the special case in which $\sigma\colon B\to X$ is the section,
we call $M\sigma$ the {\em whiskering construction}
(perhaps better, the {\em beard construction}) and write $X_w = M\sigma$ for $X$ with a whisker attached.
\end{definition}

\begin{definition}
Let $\delta$ be a familial dimension function for $G$, let $\gamma$ be a virtual representation of $G$,
let $\Mackey T$ be a contravariant $\stab\Pi_{\delta}B$-module, 
and let $\MackeyOp S$ be a covariant $\stab\Pi_{\delta}B$-module.
If $X$ is an ex-$G$-space, let
\[
 \tilde H^{G,\delta}_{\gamma+n}(X;\MackeyOp S)
   = H_{\gamma+n}(\Mackey C^{G,\delta}_{\gamma+*}(X_w) \tensor_{\stab\Pi_\delta B} \MackeyOp S)
\]
and
\[
 \tilde H_{G,\delta}^{\gamma+n}(X;\Mackey T)
   = H^{\gamma+n}(\Hom_{\stab\Pi_\delta B}(\Mackey C^{G,\delta}_{\gamma+*}(X_w),\Mackey T)).
\]
\end{definition}

The following proposition follows just as in the non-parametrized case.

\begin{proposition}
$\tilde H^{G,\delta}_{\gamma+n}(X;\MackeyOp S)$ is a well-defined covariant homotopy functor of $X$, while
$\tilde H_{G,\delta}^{\gamma+n}(X;\Mackey T)$ is a well-defined contravariant homotopy functor of $X$.
If $X$ is a $\delta$-$G$-CW$(\gamma)$ complex, these groups are naturally isomorphic to those
given by Definition~\ref{def:gammaHomology}.
\qed
\end{proposition}

\begin{theorem}[Reduced Homology and Cohomology of Spaces over $B$]
Let $\delta$ be a familial dimension function for $G$, let
$\gamma$ be a virtual representation of $\Pi_G B$, and
let $\MackeyOp S$ and $\Mackey T$ be respectively a covariant and a contravariant
$\stab\Pi_{\delta}B$-module. Then the abelian groups
$\tilde H^{G,\delta}_\gamma(X;\MackeyOp S)$ and $\tilde H_{G,\delta}^\gamma(X;\Mackey T)$ are
respectively covariant and contravariant
functors on the homotopy category of 
ex-$G$-spaces over $B$.
They are also respectively contravariant and covariant functors 
of $\gamma$.
These functors satisfy the following properties.
\begin{enumerate}

\item (Weak Equivalence)
If $f\colon X\to Y$ is a an $\F(\delta)$-equivalence of ex-$G$-spaces, then
\[
 f_*\colon \tilde H^{G,\delta}_\gamma(X;\MackeyOp S) \to \tilde H^{G,\delta}_\gamma(Y;\MackeyOp S)
\]
and
\[
 f^*\colon \tilde H_{G,\delta}^\gamma(Y;\Mackey T) \to \tilde H_{G,\delta}^\gamma(X;\Mackey T)
\]
are isomorphisms.

\item (Exactness)
If $A\to X$ is a cofibration, then the following sequences are exact:
\[
 \tilde H^{G,\delta}_\gamma(A; \MackeyOp S) \to \tilde H^{G,\delta}_\gamma(X; \MackeyOp S)
  \to \tilde H^{G,\delta}_\gamma(X/_B A; \MackeyOp S)
\]
and
\[
 \tilde H_{G,\delta}^\gamma(X/_B A; \Mackey T) \to \tilde H_{G,\delta}^\gamma(X; \Mackey T)
  \to \tilde H_{G,\delta}^\gamma(A; \Mackey T).
\]

\item (Additivity)
If $X = \Wedge_i X_i$ is a fiberwise wedge of well-based ex-$G$-spaces,
then the inclusions of the wedge summands induce isomorphisms
\[
 \Dirsum_i\tilde H^{G,\delta}_\gamma(X_i; \MackeyOp S) \iso \tilde H^{G,\delta}_\gamma(X; \MackeyOp S)
\]
and
\[
 \tilde H_{G,\delta}^\gamma(X; \Mackey T) \iso \prod_i \tilde H_{G,\delta}^\gamma(X_i; \Mackey T).
\]

\item (Suspension)
If $X$ is well-based, there are suspension isomorphisms
 \[
  \sigma^V\colon \tilde H^{G,\delta}_\alpha(X; \MackeyOp S) 
     \xrightarrow{\iso} \tilde H^{G,\delta}_{\alpha + V}(\susp_B^V X; \MackeyOp S)
 \]
and
 \[
 \sigma^V\colon \tilde H_{G,\delta}^\gamma(X; \Mackey T) 
     \xrightarrow{\iso} \tilde H_{G,\delta}^{\gamma + V}(\susp_B^V X; \Mackey T).
 \]
These isomorphisms satisfy $\sigma^0 = \id$, $\sigma^{W}\circ \sigma^V = \sigma^{V\oplus W}$,
and the following naturality condition: If $\zeta\colon V\to V'$ is an isomorphism,
then the following diagrams commute:
\[
 \xymatrix@C+2em{
  \tilde H^{G,\delta}_\gamma(X; \MackeyOp S) \ar[r]^-{\sigma^V} \ar[d]_{\sigma^{V'}}
    & \tilde H^{G,\delta}_{\gamma + V}(\susp_B^{V}X; \MackeyOp S)
          \ar[d]^{\tilde H_{\id}(\id\smsh\zeta)} \\
  \tilde H^{G,\delta}_{\gamma + V'}(\susp_B^{V'}X; \MackeyOp S)
       \ar[r]_-{\tilde H_{\id+\zeta}(\id)}
    & \tilde H^{G,\delta}_{\gamma + V}(\susp_B^{V'}X; \MackeyOp S)
 }
\]
\[
 \xymatrix@C+2em{
  \tilde H_{G,\delta}^\gamma(X; \Mackey T) \ar[r]^-{\sigma^V} \ar[d]_{\sigma^{V'}}
    & \tilde H_{G,\delta}^{\gamma + V}(\susp_B^{V}X; \Mackey T)
          \ar[d]^{\tilde H^{\id+\zeta}(\id)} \\
  \tilde H_{G,\delta}^{\gamma + V'}(\susp_B^{V'}X; \Mackey T)
       \ar[r]_-{\tilde H^{\id}(\id\smsh\zeta)}
    & \tilde H_{G,\delta}^{\gamma + V'}(\susp_B^{V}X; \Mackey T)
 }
\]

\item (Dimension Axiom)
If $H\in\F(\delta)$, $b\colon G/H\to B$, and $V$ is a representation of $H$ so large that
$\gamma_0(b) + V - \delta(G/H) + n$ is an actual representation, then
there are natural isomorphisms
\[
 \tilde H^{G,\delta}_{\gamma+V+k}(G_+\smsh_H S^{\gamma_0(b)+V-\delta(G/H)+n,b};\MackeyOp S)
  \iso
   \begin{cases}
    \MackeyOp S(b) & \text{if $k = n$} \\
    0 & \text{if $k\neq n$}
   \end{cases}
\]
and
\[
 \tilde H_{G,\delta}^{\gamma+V+k}(G_+\smsh_H S^{\gamma_0(b)+V-\delta(G/H)+n,b};\Mackey T)
  \iso
   \begin{cases}
    \Mackey T(b) & \text{if $k = n$} \\
    0 & \text{if $k\neq n$.}
   \end{cases}
\]

\end{enumerate}
\end{theorem}

The proof is the same as that of Theorem~\ref{thm:reducedHomologySpaces} with appropriate changes.

The following analogue of Corollary~\ref{cor:exactsequences} follows.

\begin{corollary}
\begin{enumerate}
\item[]
\item
If $f\colon A\to X$ is a map of well-based ex-$G$-spaces, then the following sequences are exact:
\[
 \tilde H^{G,\delta}_\gamma(A; \MackeyOp S) \to \tilde H^{G,\delta}_\gamma(X; \MackeyOp S)
  \to \tilde H^{G,\delta}_\gamma(\tilde Cf; \MackeyOp S)
\]
and
\[
 \tilde H_{G,\delta}^\gamma(\tilde Cf; \Mackey T) \to \tilde H_{G,\delta}^\gamma(X; \Mackey T)
  \to \tilde H_{G,\delta}^\gamma(A; \Mackey T).
\]

\item
If $f\colon A\to X$ is any map of ex-$G$-spaces, then the following sequences are exact,
where $f_w\colon A_w\to X_w$ is the induced map:
\[
 \tilde H^{G,\delta}_\gamma(A; \MackeyOp S) \to \tilde H^{G,\delta}_\gamma(X; \MackeyOp S)
  \to \tilde H^{G,\delta}_\gamma(\tilde C(f_w); \MackeyOp S)
\]
and
\[
 \tilde H_{G,\delta}^\gamma(\tilde C(f_w); \Mackey T) \to \tilde H_{G,\delta}^\gamma(X; \Mackey T)
  \to \tilde H_{G,\delta}^\gamma(A; \Mackey T).
\]
\qed
\end{enumerate}
\end{corollary}

Finally, we define unreduced homology and cohomology of pairs.

\begin{definition}
If $(X,A)$ is a pair of $G$-spaces over $B$, write $i\colon A\to X$ for the inclusion and let
\[
 H^{G,\delta}_\gamma(X,A;\MackeyOp S) = \tilde H^{G,\delta}_\gamma(Ci;\MackeyOp S)
\]
and
\[
 H_{G,\delta}^\gamma(X,A;\Mackey T) = \tilde H_{G,\delta}^\gamma(Ci;\Mackey T).
\]
In particular, we write 
\[
 H^{G,\delta}_\gamma(X;\MackeyOp S) = H^{G,\delta}_\gamma(X,\emptyset;\MackeyOp S)
   = \tilde H^{G,\delta}_\gamma(X_+;\MackeyOp S)
\]
and
\[
 H_{G,\delta}^\gamma(X;\Mackey T) = H_{G,\delta}^\gamma(X,\emptyset;\Mackey T) 
   = \tilde H_{G,\delta}^\gamma(X_+;\Mackey T).
\]
\end{definition}

\begin{theorem}[Unreduced Homology and Cohomology of Spaces over $B$]
Let $\delta$ be a familial dimension function for $G$, let
$\gamma$ be a virtual representation of $\Pi_G B$, and
let $\MackeyOp S$ and $\Mackey T$ be respectively a covariant and a contravariant
$\stab\Pi_{\delta}B$-module. Then the abelian groups
$H^{G,\delta}_{\gamma}(X,A;\MackeyOp S)$ and $H_{G,\delta}^{\gamma}(X,A;\Mackey T)$ are
respectively covariant and contravariant
functors on the homotopy category of 
pairs of $G$-spaces over $B$.
They are also respectively contravariant and covariant functors 
of $\gamma$.
These functors satisfy the following properties.
\begin{enumerate}

\item (Weak Equivalence)
If $f\colon (X,A)\to (Y,B)$ is an $\F(\delta)$-equivalence of pairs of $G$-spaces over $B$, then
\[
 f_*\colon H^{G,\delta}_{\gamma}(X,A;\MackeyOp S) \to H^{G,\delta}_{\gamma}(Y,B;\MackeyOp S)
\]
and
\[
 f^*\colon H_{G,\delta}^{\gamma}(Y,B;\Mackey T) \to H_{G,\delta}^{\gamma}(X,A;\Mackey T)
\]
are isomorphisms.

\item (Exactness)
If $(X,A)$ is a pair of $G$-spaces over $B$, then there are natural homomorphisms
\[
 \bndry\colon H^{G,\delta}_{\gamma+n}(X,A;\MackeyOp S) \to H^{G,\delta}_{\gamma+n-1}(A;\MackeyOp S)
\]
and
\[
 d\colon H_{G,\delta}^{\gamma+n}(A;\Mackey T) \to H_{G,\delta}^{\gamma+n+1}(X,A;\Mackey T)
\]
and long exact sequences
\[
 \cdots \to H^{G,\delta}_{\gamma+n}(A; \MackeyOp S) \to H^{G,\delta}_{\gamma+n}(X; \MackeyOp S)
  \to H^{G,\delta}_{\gamma+n}(X,A; \MackeyOp S) \to H^{G,\delta}_{\gamma+n-1}(A;\MackeyOp S) \to \cdots
\]
and
\[
 \cdots \to H_{G,\delta}^{\gamma+n-1}(A;\Mackey T) \to H_{G,\delta}^{\gamma+n}(X,A; \Mackey T) 
  \to H_{G,\delta}^{\gamma+n}(X; \Mackey T)  \to H_{G,\delta}^{\gamma+n}(A; \Mackey T) \to \cdots.
\]

\item (Excision)
If $(X;A,B)$ is an excisive triad, i.e., $X$ is the union of the interiors of $A$ and $B$,
then the inclusion $(A,A\intersect B)\to (X,B)$ induces isomorphisms
\[
 H^{G,\delta}_{\gamma}(A,A\intersect B;\MackeyOp S) \iso H^{G,\delta}_{\gamma}(X,B;\MackeyOp S)
\]
and
\[
 H_{G,\delta}^{\gamma}(X,B;\Mackey T) \iso H_{G,\delta}^{\gamma}(A,A\intersect B;\Mackey T).
\]

\item (Additivity)
If $(X,A) = \coprod_k (X_k,A_k)$ is a disjoint union of pairs of $G$-spaces over $B$,
then the inclusions $(X_k,A_k)\to (X,A)$ induce isomorphisms
\[
 \Dirsum_k H^{G,\delta}_{\gamma}(X_k,A_k; \MackeyOp S) \iso H^{G,\delta}_{\gamma}(X,A; \MackeyOp S)
\]
and
\[
 H_{G,\delta}^{\gamma}(X,A; \Mackey T) \iso \prod_k H_{G,\delta}^{\gamma}(X_k,A_k; \Mackey T).
\]

\item (Suspension)
If $A\to X$ is a cofibration, there are suspension isomorphisms
 \[
  \sigma^V\colon H^{G,\delta}_{\gamma}(X,A; \MackeyOp S) 
     \xrightarrow{\iso} H^{G,\delta}_{\gamma+V}((X,A)\times(D(V),S(V)); \MackeyOp S)
 \]
and
 \[
 \sigma^V\colon H_{G,\delta}^{\gamma}(X,A; \Mackey T) 
     \xrightarrow{\iso} H_{G,\delta}^{\gamma+V}((X,A)\times(D(V),S(V)); \Mackey T),
 \]
where $(X,A)\times (D(V),S(V)) = (X\times D(V), X\times S(V)\union A\times D(V))$.
These isomorphisms satisfy $\sigma^0 = \id$ and $\sigma^{W}\circ \sigma^V = \sigma^{V\oplus W}$,
under the identification $\bar D(V) \times \bar D(W) \homeo \bar D(V\oplus W)$
(here we use the notation $\bar D(V) = (D(V),S(V))$).
They also satisfy the following naturality condition: If $\zeta\colon V\to V'$ is an isomorphism,
then the following diagrams commute:
\[
 \xymatrix@C+2em{
  H^{G,\delta}_{\gamma}(X,A; \MackeyOp S) \ar[r]^-{\sigma^V} \ar[d]_{\sigma^{V'}}
    & H^{G,\delta}_{\gamma+V}((X,A)\times \bar D(V); \MackeyOp S)
          \ar[d]^{H^G_{\id}(\id\times\zeta)} \\
  H^{G,\delta}_{\gamma+V'}((X,A)\times \bar D(V'); \MackeyOp S)
       \ar[r]_-{H^G_{\id\dirsum\zeta}(\id)}
    & H^{G,\delta}_{\gamma+V}((X,A)\times \bar D(V'); \MackeyOp S)
 }
\]
\[
 \xymatrix@C+2em{
  H_{G,\delta}^{\gamma}(X,A; \Mackey T) \ar[r]^-{\sigma^V} \ar[d]_{\sigma^{V'}}
    & H_{G,\delta}^{\gamma+V}((X,A)\times \bar D(V); \Mackey T)
          \ar[d]^{H_G^{\id\dirsum\zeta}(\id)} \\
  H_{G,\delta}^{\gamma+V'}((X,A)\times \bar D(V'); \Mackey T)
       \ar[r]_-{H_G^{\id}(\id\times\zeta)}
    & H_{G,\delta}^{\gamma+V'}((X,A)\times \bar D(V); \Mackey T)
 }
\]

\item (Dimension Axiom)
If $H\in\F(\delta)$, $b\colon G/H\to B$, and $V$ is a representation of $H$ so large that
$\gamma_0(b)+V-\delta(G/H)+n$ is an actual representation, write
\[
  \bar D_b(\gamma_0(b)+V-\delta(G/H)+n) = b_!(G\times_H \bar D(\gamma_0(b)+V-\delta(G/H)+n)),
\]
i.e., the pair $G\times_H \bar D(\gamma_0(b)+V-\delta(G/H)+n)$ mapping to $B$ via $b$.
Then there are natural isomorphisms
\[
 H^{G,\delta}_{\gamma+V+k}(\bar D_b(\gamma_0(b)+V-\delta(G/H)+n);\MackeyOp S)
  \iso
   \begin{cases}
    \MackeyOp S(b) & \text{if $k = n$} \\
    0 & \text{if $k\neq n$}
   \end{cases}
\]
and
\[
 H_{G,\delta}^{\gamma+V+k}(\bar D_b(\gamma_0(b)+V-\delta(G/H)+n);\Mackey T)
  \iso
   \begin{cases}
    \Mackey T(b) & \text{if $k = n$} \\
    0 & \text{if $k\neq n$.}
   \end{cases}
\]

\end{enumerate}
\end{theorem}

The proof is the same as Theorem~\ref{thm:unreducedHomology} with appropriate changes.

\begin{definition}\label{def:paramMackeyvalued}
Given a ex-$G$-space $X$ over $B$, define the covariant $\stab\Pi_\delta B$-module
$\MackeyOp H^{G,\delta}_\gamma(X;\MackeyOp S)$ by
\begin{align*}
 \MackeyOp H^{G,\delta}_\gamma(X;\MackeyOp S)(b\colon G/H\to B)
  &= \tilde H^{G,\delta}_\gamma(X\smsh (G_+\smsh_H S^{-\delta(G/H),b});\MackeyOp S) \\
  &= \tilde H^{G,\delta}_{\gamma+V}(X\smsh (G_+\smsh_H S^{V-\delta(G/H),b});\MackeyOp S)
\end{align*}
for a $V$ sufficiently large that $V-\delta(G/H)$ is an actual representation.
Similarly, define the contravariant $\stab\Pi_\delta B$-module
$\Mackey H_{G,\delta}^\gamma(X;\Mackey T)$ by
\begin{align*}
 \Mackey H_{G,\delta}^\gamma(X;\Mackey T)(b)
  &= \tilde H_{G,\delta}^\gamma(X\smsh (G_+\smsh_H S^{-\delta(G/H),b});\Mackey T) \\
  &= \tilde H_{G,\delta}^{\gamma+V}(X\smsh (G_+\smsh_H S^{V-\delta(G/H),b});\Mackey T).
\end{align*}
\end{definition}

In these terms, the dimension axioms take following form:
\[
 \MackeyOp H^{G,\delta}_{n}(B_+;\MackeyOp S)
  \iso \begin{cases}
          \MackeyOp S &\text{if $n = 0$} \\
          0 &\text{if $n \neq 0$}
       \end{cases}
\]
and
\[
 \Mackey H_{G,\delta}^{n}(B_+;\Mackey T)
  \iso \begin{cases}
          \Mackey T &\text{if $n = 0$} \\
          0 &\text{if $n \neq 0$.}
       \end{cases}
\]

The following theorem is proved in the same way as Theorem~\ref{thm:AtiyahHirzebruch}.
Recall that, 
if $\tilde h^G_*(-)$ is a generalized $RO(G)$-graded homology theory on ex-$G$-spaces over $B$
(as in Definition~\ref{def:parametrizedtheory}),
$\gamma$ is a virtual representation of $\Pi B$,
and $\delta$ is a dimension function for $G$, then
Definition~\ref{def:homologymackey} defines the coefficient system
$\MackeyOp h^{G,\delta}_\gamma$.
Similarly, if $\tilde h_G^*(-)$ is a reduced $RO(G)$-graded cohomology theory,
we have the coefficient system $\Mackey h_{G,\delta}^\gamma$.

\begin{theorem}[Atiyah-Hirzebruch Spectral Sequence]
Suppose that $\tilde h^G_*(-)$ is an $RO(G)$-graded homology theory on
ex-$G$-spaces over $B$.
Let $\delta$ be a familial dimension function for $G$,
let $\alpha$ be an element of $RO(G)$, and let
$\gamma$ be a virtual representation of $\Pi B$.
Assume that $\tilde h^G_*(-)$ takes $\F(\delta)$-equivalences to isomorphisms.
Then there is a strongly convergent spectral sequence
 \[
 E^2_{p,q} = \tilde H^{G,\delta}_{\alpha-\gamma+p}(X;\MackeyOp h^{G,\delta}_{\gamma+q})
  \convto \tilde h^G_{\alpha+p+q}(X).
 \]
Similarly, if $h_G^*(-)$ is an $RO(G)$-graded cohomology theory
on ex-$G$-spaces over $B$,
taking $\F(\delta)$-equivalences to isomorphisms, then
there is a conditionally convergent spectral sequence
 \[
 E_2^{p,q} = \tilde H_{G,\delta}^{\alpha-\gamma+p}(X;\Mackey h_{G,\delta}^{\gamma+q})
  \convto h_G^{\alpha+p+q}(X).
 \]
\qed
\end{theorem}

Uniqueness follows in the usual way.

\begin{corollary}[Uniqueness of Ordinary $RO(\Pi B)$-Graded Homology]\label{cor:tCohomUniqueness}
Let $\delta$ be a familial dimension function for $G$.
Let $\tilde h^G_*(-)$ be an $RO(G)$-graded homology theory on ex-$G$-spaces over $B$
that takes $\F(\delta)$-equivalences to isomorphisms
and let $\gamma$ be a virtual representation of $\Pi B$.
Suppose that, for integers $n$,
\[
 \MackeyOp h^{G,\delta}_{\gamma+n} = 0\qquad\text{for $n\neq 0$}
\]
Then, for $\alpha\in RO(G)$, there is a natural isomorphism
\[
 \tilde h^G_\alpha(-) \iso \tilde H^{G,\delta}_{\alpha-\gamma}(-;\MackeyOp h^{G,\delta}_\gamma).
\]
There is a similar statement for cohomology theories.
\qed
 \end{corollary}

Finally, we have a universal coefficients spectral sequence, constructed in the usual way.
Write $\Mackey H^{G,\delta}_\gamma(X)$ for the contravariant $\stab\Pi_{\delta}B$-module given by
 \[
 \Mackey H^{G,\delta}_\gamma(X)(b) =
   \tilde H^{G,\delta}_\gamma(X;\stab\Pi_\delta B(b,-)).
 \]
$\Tor_*^{\stab\Pi_\delta B}$ and $\Ext^*_{\stab\Pi_\delta B}$
below are the derived functors of $\tensor_{\stab\Pi_\delta B}$ and
$\Hom_{\stab\Pi_\delta B}$ respectively.

\begin{theorem}[Universal Coefficients Spectral Sequence]
If $\MackeyOp S$ is a covariant $\stab\Pi_{\delta}B$-module, $\Mackey T$ is a
contravariant $\stab\Pi_{\delta}B$-module, and $\gamma$
is a virtual representation of $\Pi B$,
there are spectral sequences
 \[
 E^2_{p,q} = \Tor^{\stab\Pi_\delta B}_p(\Mackey H^{G,\delta}_{\gamma+q}(X), \MackeyOp S)
 \convto
 \tilde H^{G,\delta}_{\gamma+p+q}(X;\MackeyOp S)
 \]
and
 \[
 E_2^{p,q} = \Ext_{\stab\Pi_\delta B}^p(\Mackey H^{G,\delta}_{\gamma+q}(X), \Mackey T)
 \convto
 \tilde H_{G,\delta}^{\gamma+p+q}(X;\Mackey T).
 \]
\qed
\end{theorem}

\section{The representing spectra}\label{sec:paramrepspectra}

Let $\delta$ be a familial dimension function for $G$, let $\gamma$
be a virtual representation of $\Pi_G B$, and let $\Mackey T$ be a contravariant
$\stab\Pi_{\delta}B$-module.
By Theorem~\ref{thm:cohomrep}, there exists a parametrized prespectrum representing
the $RO(G)$-graded cohomology theory 
$\tilde H_{G,\delta}^{\gamma+*}(-;\Mackey T)$, 
which we will call $H_\delta\Mackey T^\gamma$. 

If $E$ is any $G$-spectrum over $B$, define the 
contravariant $\stab\Pi_\delta B$-module $\Mackey\pi^\delta_\gamma E$ by
\[
 \Mackey\pi^\delta_\gamma E(b\colon G/K\to B) 
  = [G_+\smsh_K S^{\gamma_0(b) - \delta(G/K),b}, E]_{G,B}.
\]
This obviously defines a functor on
$\stab\Pi_{\delta,\gamma} B$, and we make it a functor on
$\stab\Pi_\delta B$ using the isomorphism of Theorem~\ref{thm:StableMapsOrbitsOverB}.
From the dimension axiom for cellular cohomology, we have
\[
 \Mackey\pi^\delta_{\gamma+n} H_\delta\Mackey T^\gamma \iso
 \begin{cases} \Mackey T & n = 0 \\ 0 & n \neq 0 \end{cases}
\]
and, from the uniqueness of ordinary cohomology, this together with the equivalence
$H_\delta\Mackey T^\gamma \hmtpc F(E\F(\delta)_+,H_\delta\Mackey T^\gamma)$ characterizes
$H_\delta\Mackey T^\gamma$ up to nonunique homotopy equivalence.
We call any spectrum satisfying these conditions a
{\em parametrized Eilenberg-Mac\,Lane spectrum.}
Using Corollary~\ref{cor:coefficientsfibers} and the discussion following it,
we can rewrite the calculation above as
\[
 \Mackey\pi^{K,\delta}_{\gamma_0(b)+n}((H_\delta\Mackey T^\gamma)_b)
  \iso \begin{cases} b^*\Mackey T & n = 0 \\ 0 & n \neq 0 \end{cases}
\]
for $b\colon G/K\to B$. Thus, we can identify each fiber as the suspension of
a nonparametrized Eilenberg-Mac\,Lane spectrum:
\[
 (H_\delta\Mackey T^\gamma)_b \hmtpc \susp^{\gamma_0(b)}H_\delta(b^*\Mackey T)
\]
as $K$-spectra.
Again, this characterizes $H_\delta\Mackey T^\gamma$.

What about cellular homology?
It follows from
\cite[18.6.7 \& 18.1.5]{MaySig:parametrized} that
\begin{multline*}
 [S^n, \rho_!(H_\delta\Mackey T^\gamma \smsh_B G_+\smsh_H S^{-\gamma_0(b) - (\Lie-\delta)(G/H),b})]_{G,B} \\
 \iso [G_+\smsh_H S^{\gamma_0(b) - \delta(G/H)+n,b}, H_\delta\Mackey T^\gamma]_{G,B}.
\end{multline*}
It follows, as in Theorem~\ref{thm:dualCohomologyRep}, that
\[
 H^{\Lie-\delta}\Mackey T_{-\gamma} = H_\delta\Mackey T^\gamma \smsh E\F(\delta)_+
\]
represents $\tilde H^{G,\Lie-\delta}_{-\gamma+*}(-;\Mackey T)$, or that
\[
 H^\delta\MackeyOp S_\gamma = H_{\Lie-\delta}\MackeyOp S^{-\gamma}\smsh E\F(\delta)_+
\]
represents $\tilde H^{G,\delta}_{\gamma+*}(-;\MackeyOp S)$.
The characterization of the homotopy groups then looks like
\[
 \Mackey\pi^{\Lie-\delta}_{-\gamma+n} H^\delta\MackeyOp S_\gamma \iso
 \begin{cases} \MackeyOp S & n = 0 \\ 0 & n \neq 0 \end{cases}
\]
The spectrum $H^\delta\MackeyOp S_\gamma$ can also be characterized by the fact that
its fibers are nonparametrized Eilenberg-Mac\,Lane spectra representing ordinary homology theories.

It will be useful to have a simple, explicit construction of these spectra.
We simply have to elaborate Construction~\ref{con:EilenbergMacLane}.

\begin{construction}\label{con:paramEilenbergMacLane}
Let $\delta$ be a familial dimension function for $G$,
let $\gamma$ be a virtual representation of $\Pi_G B$,
and let $\Mackey T$ be a contravariant $\stab\Pi_{G,\delta}B$-module.
Define
\[
 F_\delta\Mackey T^\gamma = \Wedge_{\Mackey T(b)} G_+\smsh_{H_b} S^{\gamma_0(b)-\delta(G/H),b},
\]
where the wedge runs over all objects $b$ in $\stab\Pi_{G,\delta}B$
and all elements in $\Mackey T(b)$.
Then $F_\delta(-)^\gamma$ is a functor taking contravariant $\stab\Pi_{G,\delta}B$-modules to $G$-spectra 
and we have a natural epimorphism
\[
 \epsilon\colon \Mackey\pi_\gamma^{\delta} F_\delta\Mackey T \to \Mackey T.
\]
Let
\[
 K = \{ \kappa\colon G_+\smsh_{H_\kappa} S^{\gamma_0(b_\kappa)-\delta(G/H_\kappa),b_\kappa}
           \to F_\delta\Mackey T
       \mid \epsilon(\kappa) = 0 \},
\]
and let
\[
 R_\delta\Mackey T^\gamma = \Wedge_{\kappa\in K} S^{\gamma_0(b_\kappa)-\delta(G/H_\kappa),b_\kappa}.
\]
Then $R_\delta(-)^\gamma$ is also a functor,
there is a natural transformation $R_\delta\Mackey T^\gamma\to F_\delta\Mackey T^\gamma$
(given by the maps $\kappa$), and we have an
exact sequence
\[
 \Mackey\pi_\gamma^{\delta}R_\delta\Mackey T^\gamma 
   \to \Mackey\pi_\gamma^{\delta} F_\delta\Mackey T^\gamma 
   \to \Mackey T \to 0.
\]
It follows that, if we let $C_\delta\Mackey T^\gamma$ be the cofiber of 
$R_\delta\Mackey T^\gamma\to F_\delta\Mackey T^\gamma$,
then $\Mackey\pi_\gamma^{\delta}C_\delta\Mackey T^\gamma \iso \Mackey T$
and $\Mackey\pi_{\gamma+n}^{\delta}C_\delta\Mackey T^\gamma = 0$ for $n<0$.

We can then functorially kill all the homotopy 
$\Mackey\pi_{\gamma+n}^{\delta}C_\delta\Mackey T^\gamma$ for $n>0$, obtaining
a functor $P_\delta(-)^\gamma$ with
\[
 \Mackey\pi_{\gamma+n}^{\delta}P_\delta\Mackey T^\gamma =
   \begin{cases}
      \Mackey T &\text{if $n=0$} \\
      0 &\text{if $n\neq 0$.}
   \end{cases}
\]
Finally, we let
\[
 H_\delta\Mackey T^\gamma = F(E\F(\delta)_+,P_\delta\Mackey T^\gamma)
\]
and
\[
 H^\delta\MackeyOp S_\gamma = P_{\Lie-\delta}\MackeyOp S^{-\gamma} \smsh E\F(\delta)_+
\]
if $\MackeyOp S$ is a covariant $\delta$-Mackey functor.
Then $H_\delta(-)^\gamma$ and $H^\delta(-)^\gamma$ are functors and produce spectra that
satisfy the characterizations of the Eilenberg-Mac\,Lane spectra representing
cohomology and homology, respectively.
\qed
\end{construction}

\section{Change of base space}

Let $f\colon A\to B$ be a $G$-map.
The main results of this section come from the observation that,
if $X$ is a CW complex over $A$, then $f_! X$ is a CW complex over $B$ with
corresponding cells. 
More precisely, if $\gamma$ is a virtual representation of $\Pi_G B$, then,
if $X$ is a $\delta$-$G$-CW$(f^*\gamma)$ complex over $A$, then
$f_!X$ is a $\delta$-$G$-CW$(\gamma)$ complex over $B$.

For the algebra, we have
\[
 f_!\colon \stab\Pi_{G,\delta}A \to \stab\Pi_{G,\delta}B,
\]
given by the change-of-base functor $f_!$ on parametrized spectra.
We write $f_! = (f_!)_!$ and $f^* = (f_!)^*$ for the induced functors on
modules.

This leads to the following isomorphism of chain complexes.

\begin{proposition}
If $\gamma$ is a virtual representation of $\Pi_G B$ and
$X$ is a $\delta$-$G$-CW$(f^*\gamma)$ complex over $A$, then
\[
 \Mackey C^{G,\delta}_{\gamma+*}(f_!X) \iso f_! \Mackey C^{G,\delta}_{f^*\gamma+*}(X).
\]
Moreover, this isomorphism respects suspension isomorphisms.
\end{proposition}

\begin{proof}
We define a map $f_! \Mackey C^{G,\delta}_{f^*\gamma+*}(X) \to \Mackey C^{G,\delta}_{\gamma+*}(f_!X)$
as follows.
For simplicity of notation, we write $\gamma$ again for $f^*\gamma$
and, for $a\colon G/H\to A$ an object of $\stab\Pi_{G,\delta}A$, we write
\[
 S^{\gamma+n,a} = G_+\smsh_H S^{\gamma_0(a)-\delta(G/H)+n,a}.
\]
If $b\in \stab\Pi_{G,\delta}B$, we consider the following map:
\begin{align*}
 f_! \Mackey C^{G,\delta}_{f^*\gamma+n}(X)(b)
  &= \int\nolimits^{a\in\stab\Pi_{G,\delta}A} 
     \Mackey C^{G,\delta}_{f^*\gamma+n}(X)(a) \tensor \stab\Pi_{G,\delta}B(b,f_!(a)) \\
  &= [S^{\gamma+n,a}, \susp^\infty_A X^{\gamma+n}/X^{\gamma+n-1}]_{G,A}
     \tensor \stab\Pi_{G,\delta}B(b,f_!(a)) \\
  &\xrightarrow{f_!} [f_!S^{\gamma+n,a},f_! \susp^\infty_A X^{\gamma+n}/X^{\gamma+n-1}]_{G,B}
     \tensor \stab\Pi_{G,\delta}B(b,f_!(a)) \\
  &\iso [S^{\gamma+n,f_!(a)}, \susp^\infty_B (f_!X)^{\gamma+n}/(f_!X)^{\gamma+n-1}]_{G,B}
     \tensor \stab\Pi_{G,\delta}B(b,f_!(a)) \\
  &\to [S^{\gamma+n,b}, \susp^\infty_B (f_!X)^{\gamma+n}/(f_!X)^{\gamma+n-1}]_{G,B} \\
  &= \Mackey C^{G,\delta}_{\gamma+n}(f_!X)(b).
\end{align*}
This description shows that the map is a chain map.
Further, when we look at its effect on one cell in $X^{\gamma+n}/X^{\gamma+n-1}$,
Proposition~\ref{prop:indresadjunction}(5) implies that we have an isomorphism,
so Lemma~\ref{lem:skeletalquotient} shows that the map above is a chain isomorphism.

That this isomorphism respects suspension isomorphisms follows from the definition
of the map above.
\end{proof}

This gives us the following isomorphisms in homology and cohomology.

\begin{theorem}\label{thm:changeofbasespace}
Let $\delta$ be a dimension function for $G$, let $f\colon A\to B$ be a $G$-map,
let $\gamma$ be a virtual representation of $\Pi_G B$,
let $\Mackey T$ be a contravariant $\stab\Pi_{G,\delta}B$-module, and
let $\MackeyOp S$ be a covariant $\stab\Pi_{G,\delta}B$-module.
If $X$ is a $\delta$-$G$-CW$(f^*\gamma)$ complex, then there are natural isomorphisms
\[
 \tCH^{G,\delta}_\gamma(f_!X;\MackeyOp S) \iso \tCH^{G,\delta}_{f^*\gamma}(X;f^*\MackeyOp S)
\]
and
\[
 \tCH_{G,\delta}^\gamma(f_!X;\Mackey T) \iso \tCH_{G,\delta}^{f^*\gamma}(X;f^*\Mackey T).
\]
If $\delta$ is familial, for $X$ a well-based ex-$G$-space there are natural isomorphisms
\[
 \tilde H^{G,\delta}_\gamma(f_!X;\MackeyOp S) \iso \tilde H^{G,\delta}_{f^*\gamma}(X;f^*\MackeyOp S)
\]
and
\[
 \tilde H_{G,\delta}^\gamma(f_!X;\Mackey T) \iso \tilde H_{G,\delta}^{f^*\gamma}(X;f^*\Mackey T).
\]
These isomorphisms all respect suspension isomorphisms.
\end{theorem}

\begin{proof}
If $X$ is a $\delta$-$G$-CW$(f^*\gamma)$ complex, then the preceding proposition gives us
\begin{align*}
 \Mackey C^{G,\delta}_{\gamma+*}(f_!X) \tensor_{\stab\Pi_{G,\delta}B} \MackeyOp S
  &\iso f_!\Mackey C^{G,\delta}_{f^*\gamma+*}(X) \tensor_{\stab\Pi_{G,\delta}B} \MackeyOp S \\
  &\iso \Mackey C^{G,\delta}_{f^*\gamma+*}(X) \tensor_{\stab\Pi_{G,\delta}A} f^*\MackeyOp S.
\end{align*}
On taking homology we get the first isomorphism of the theorem.
Similarly, we have
\begin{align*}
 \Hom_{\stab\Pi_{G,\delta}B}( \Mackey C^{G,\delta}_{\gamma+*}(f_!X), \Mackey T )
  &\iso \Hom_{\stab\Pi_{G,\delta}B}( f_!\Mackey C^{G,\delta}_{f^*\gamma+*}(X), \Mackey T ) \\
  &\iso \Hom_{\stab\Pi_{G,\delta}A}( \Mackey C^{G,\delta}_{f^*\gamma+*}(X), f^*\Mackey T ).
\end{align*}
On taking homology we get the second isomorphism of the theorem.
That these isomorphisms respect suspension follows from the preceding proposition.

Now let $X$ be a well-based ex-$G$-space over $A$.
If $\Gamma X \to X$ is an approximation
by a $\delta$-$G$-CW$(f^*\gamma)$ complex,
then $f_!\Gamma X \to f_!X$ is an approximation by a $\delta$-$G$-CW$(\gamma)$ complex.
This follows from Theorems~\ref{thm:gammaequivalence} and~\ref{thm:weakCharacterization}
and the fact that, nonequivariantly, $f_!$ preserves weak equivalence of
well-based ex-spaces \cite[7.3.4]{MaySig:parametrized}.
From this it follows that $f_!$ takes an approximation
$\Gamma\susp_A^\infty X \to \susp_A^\infty X$ by a $G$-CW$(f^*\gamma)$ spectrum 
to an approximation of $\susp_B^\infty f_!X$, 
which gives a chain isomorphism 
$\Mackey C^{G,\delta}_{\gamma+*}(f_!X) \iso f_!\Mackey C^{G,\delta}_{f^*\gamma+*}(X)$.
The homology and cohomology
isomorphisms now follow as for CW complexes.
\end{proof}

We can use these isomorphisms to get the following ``push-forward'' maps:
Let $\MackeyOp S$ be a covariant $\stab\Pi_{G,\delta}A$-module and
let $\Mackey T$ be a contravariant $\stab\Pi_{G,\delta}A$-module.
Using the unit $\eta\colon \MackeyOp S\to f^*f_!\MackeyOp S$ of the adjunction, we can then define
\[
 f_!\colon \tilde H^{G,\delta}_{f^*\gamma}(X;\MackeyOp S) \to
 	\tilde H^{G,\delta}_{\gamma}(f_!X; f_!\MackeyOp S)
\]
to be the composite
\[
 \tilde H^{G,\delta}_{f^*\gamma}(X;\MackeyOp S)
  \xrightarrow{\eta_*} \tilde H^{G,\delta}_{f^*\gamma}(X;f^*f_!\MackeyOp S)
  \iso \tilde H^{G,\delta}_{\gamma}(f_!X;f_!\MackeyOp S).
\]
We define
\[
 f_!\colon \tilde H_{G,\delta}^{f^*\gamma}(X;\Mackey T) \to
 	\tilde H_{G,\delta}^{\gamma}(f_!X; f_!\Mackey T)
\]
similarly. Note that, if $\MackeyOp U$ is a covariant $\stab\Pi_{G,\delta}B$-module,
then the composite
\[
 \tilde H^{G,\delta}_{f^*\gamma}(X;f^*\MackeyOp U)
  \xrightarrow{f_!} \tilde H^{G,\delta}_{\gamma}(f_!X; f_!f^*\MackeyOp U)
  \xrightarrow{\epsilon_*} \tilde H^{G,\delta}_{\gamma}(f_!X; \MackeyOp U)
\]
agrees with the isomorphism of the Theorem~\ref{thm:changeofbasespace}, 
where $\epsilon$ is the counit of the adjunction.
Of course, the similar statement is true for cohomology.

We now look at how the isomorphisms of Theorem~\ref{thm:changeofbasespace} are represented.
The main result is the following.

\begin{proposition}
Let $f\colon A\to B$ be a $G$-map.
Let $\delta$ be a familial dimension function for $G$,
let $\gamma$ be a virtual representation of $\Pi_G B$,
let $\MackeyOp S$ be a covariant $\stab\Pi_{G,\delta}B$-module,
and let $\Mackey T$ be a contravariant $\stab\Pi_{G,\delta}B$-module.
Then we have stable equivalences
\[
 f^* H^\delta\MackeyOp S_\gamma \hmtpc H^\delta(f^*\MackeyOp S)_{f^*\gamma}
\]
and
\[
 f^* H_\delta\Mackey T^\gamma \hmtpc H_\delta(f^*\Mackey T)^{f^*\gamma}
\]
\end{proposition}

\begin{proof}
If $a\colon G/K\to A$ is a $G$-map, then
\[
 a^*f^* H_\delta\Mackey T^\gamma \hmtpc (fa)^* H_\delta\Mackey T^\gamma
\]
is a nonparametrized Eilenberg-Mac\,Lane spectrum by the characterization
given in the preceding section, with
\[
 \Mackey\pi^{K,\delta}_{f^*\gamma_0(a)}((f^*H_\delta\Mackey T^\gamma)_a)
 \iso \Mackey\pi^{K,\delta}_{\gamma_0(fa)}((H_\delta\Mackey T^\gamma)_{fa})
  \iso (fa)^*\Mackey T,
\]
hence
\[
 \Mackey\pi^\delta_{f^*\gamma} f^*H_\delta\Mackey T^\gamma \iso f^*\Mackey T.
\]
Further,
\[
 F(E\F(\delta)_+, f^*H_\delta\Mackey T^\gamma) 
  \hmtpc f^*F(E\F(\delta)_+, H_\delta\Mackey T^\gamma)
  \hmtpc f^* H_\delta\Mackey T^\gamma.
\]
Hence, $f^* H_\delta\Mackey T^\gamma \hmtpc H_\delta(f^*\Mackey T)^{f^*\gamma}$ as claimed.
The proof for covariant modules is essentially the same.
\end{proof}

The isomorphisms of Theorem~\ref{thm:changeofbasespace} are then represented as follows:
\begin{align*}
 \tilde H^{G,\delta}_\gamma(f_!X; \MackeyOp S)
  &\iso [S, \rho_!(H^\delta\MackeyOp S_\gamma \smsh_B f_!X]_G \\
  &\iso [S, \rho_!f_!(f^*H^\delta\MackeyOp S_\gamma \smsh_A X]_G \\
  &\iso [S, \rho_!(H^\delta(f^*\MackeyOp S)_{f^*\gamma} \smsh_A X]_G \\
  &\iso \tilde H^{G,\delta}_{f^*\gamma}(X; f^*\MackeyOp S)
\end{align*}
and
\begin{align*}
 \tilde H_{G,\delta}^\gamma(f_!X; \Mackey T)
  &\iso [f_!X, H_\delta\Mackey T^\gamma]_{G,B} \\
  &\iso [X, f^*H_\delta\Mackey T^\gamma]_{G,A} \\
  &\iso [X, H_\delta(f^*\Mackey T)^{f^*\gamma}]_{G,A} \\
  &\iso \tilde H_{G,\delta}^{f^*\gamma}(X; f^*\Mackey T).
\end{align*}

As for the push-forward map,
if $\MackeyOp S$ is a covariant $\stab\Pi_{G,\delta}A$-module, the push-forward is
represented as
\begin{align*}
 \tilde H^{G,\delta}_{f^*\gamma}(X;\MackeyOp S)
  &\iso [S, \rho_!(H^\delta\MackeyOp S_{f^*\gamma} \smsh_A X)]_G \\
  &\iso [S, \rho_!f_!(H^\delta\MackeyOp S_{f^*\gamma}\smsh_A X)]_G \\
  &\to [S, \rho_!f_!(H^\delta(f^*f_!\MackeyOp S)_{f^*\gamma}\smsh_A X)]_G \\
  &\iso [S, \rho_!f_!(f^*H^\delta(f_!\MackeyOp S)_{\gamma}\smsh_A X)]_G \\
  &\iso [S, \rho_!(H^\delta(f_!\MackeyOp S)_{\gamma} \smsh_B f_!X)]_G \\
  &\iso \tilde H^{G,\delta}_{\gamma}(f_!X; f_!\MackeyOp S).
\end{align*}
We could write down a similar composite in cohomology, but it's interesting to write
it in a slightly different way.
If $\Mackey T$ is a contravariant $\stab\Pi_{G,\delta}A$-module, we have a map
\[
 f_!H_\delta\Mackey T^{f^*\gamma} \to H_\delta(f_!\Mackey T)^\gamma
\]
adjoint to the map 
\[
 H_\delta\Mackey T^{f^*\gamma} 
   \xrightarrow{\eta} H_\delta(f^*f_!\Mackey T)^{f^*\gamma}
   \xrightarrow{\hmtpc} f^*H_\delta(f_!\Mackey T)^\gamma.
\]
The push-forward map in cohomology is then
\begin{align*}
 \tilde H_{G,\delta}^{f^*\gamma}(X;\Mackey T)
  &\iso [X, H_\delta\Mackey T^{f^*\gamma}]_{G,A} \\
  &\to [f_!X, f_!H_\delta\Mackey T^{f^*\gamma}]_{G,B} \\
  &\to [f_!X, H_\delta(f_!\Mackey T)^\gamma]_{G,B} \\
  &\iso \tilde H_{G,\delta}^\gamma(f_!X; f_!\Mackey T).
\end{align*}
For homology, we have a similar map 
$f_!H^\delta\MackeyOp S_{f^*\gamma} \to H^\delta(f_!\MackeyOp S)_\gamma$.
We also use the map of spectra $f_!(E\smsh_A F)\to f_!E\smsh_B f_!F$
given by the composite
\[
 f_!(E\smsh_A F)\to f_!(f^*f_!E\smsh_A F) \hmtpc f_!E\smsh_B f_!F.
\]
We can then write the push-forward in homology as
\begin{align*}
 \tilde H^{G,\delta}_{f^*\gamma}(X;\MackeyOp S)
  &\iso [S, \rho_!(H^\delta\MackeyOp S_{f^*\gamma} \smsh_A X)]_G \\
  &\iso [S, \rho_!f_!(H^\delta\MackeyOp S_{f^*\gamma} \smsh_A X)]_G \\
  &\to [S, \rho_!(f_!H^\delta\MackeyOp S_{f^*\gamma} \smsh_B f_!X)]_G \\
  &\to [S, \rho_!(H^\delta(f_!\MackeyOp S)_\gamma \smsh_B f_!X)]_G \\
  &\iso \tilde H^{G,\delta}_{\gamma}(f_!X; f_!\MackeyOp S).
\end{align*}  
An interesting fact about this way of writing the push-forward is the following.

\begin{proposition}
Let $f\colon A\to B$ be a $G$-map.
Let $\delta$ be a familial dimension function for $G$,
let $\gamma$ be a virtual representation of $\Pi_G B$,
let $\MackeyOp S$ be a covariant $\stab\Pi_{G,\delta}A$-module,
and let $\Mackey T$ be a contravariant $\stab\Pi_{G,\delta}A$-module.
Then we have
\[
 \Mackey\pi^{\Lie-\delta}_{-\gamma+n} f_!H^\delta\MackeyOp S_{f^*\gamma} \iso 
   \begin{cases}
      f_!\MackeyOp S &\text{if $n=0$} \\
      0 &\text{if $n<0$}.
   \end{cases}
\]
Further, the map
$f_!H^\delta\MackeyOp S_{f^*\gamma} \to H^\delta(f_!\MackeyOp S)_\gamma$
induces an isomorphism in $\Mackey\pi^{\Lie-\delta}_{-\gamma}$.
If $\delta$ is complete, then we have
\[
 \Mackey\pi^\delta_{\gamma+n} f_!H_\delta\Mackey T^{f^*\gamma} \iso 
   \begin{cases}
      f_!\Mackey T &\text{if $n=0$} \\
      0 &\text{if $n<0$.}
   \end{cases}
\]
Further, the map
$f_!H_\delta\Mackey T^{f^*\gamma} \to H_\delta(f_!\Mackey T)^\gamma$
induces an isomorphism in $\Mackey\pi^\delta_{\gamma}$.
\end{proposition}

\begin{proof}
The approach to calculating the homotopy groups is similar
to that used in the proof of Proposition~\ref{prop:EMsmshhomotopy}.
Consider the case of $H^\delta\MackeyOp S_{f^*\gamma}$ and assume
that it has been constructed as in Construction~\ref{con:paramEilenbergMacLane},
so that we have a cofibration
\[
 R^\delta\MackeyOp S_{f^*\gamma} \to F^\delta\MackeyOp S_{f^*\gamma} 
   \to C^\delta\MackeyOp S_{f^*\gamma}
\]
which induces an exact sequence
\[
 \Mackey\pi^{\Lie-\delta}_{-f^*\gamma}R^\delta\MackeyOp S_{f^*\gamma}
  \to \Mackey\pi^{\Lie-\delta}_{-f^*\gamma}F^\delta\MackeyOp S_{f^*\gamma}
  \to \MackeyOp S \to 0.
\]
Now, $F^\delta\MackeyOp S_{f^*\gamma}$ and $R^\delta\MackeyOp S_{f^*\gamma}$ are wedges of spheres
and we have
\[
 f_!(G_+\smsh_K S^{-f^*\gamma_0(a)-(\Lie-\delta)(G/K),a}) 
   = G_+\smsh_K S^{-\gamma_0(fa)-(\Lie-\delta)(G/K),fa},
\]
which implies that
\[
 \Mackey\pi^{\Lie-\delta}_{-\gamma}f_!F^\delta\MackeyOp S_{f^*\gamma}
  \iso f_!\Mackey\pi^{\Lie-\delta}_{-f^*\gamma}F^\delta\MackeyOp S_{f^*\gamma}
\]
and similarly for $R^\delta$.
Because the functor $f_!$ on spectra preserves cofibrations, we have the cofibration sequence
\[
 f_!R^\delta\MackeyOp S_{f^*\gamma} \to f_!F^\delta\MackeyOp S_{f^*\gamma} 
     \to f_!C^\delta\MackeyOp S_{f^*\gamma}
\]
Taking homotopy and using the isomorphism above, we get
exact sequences
\[
 f_!\Mackey\pi^{\Lie-\delta}_{-f^*\gamma+n}R^\delta\MackeyOp S_{f^*\gamma}
  \to f_!\Mackey\pi^{\Lie-\delta}_{-f^*\gamma+n}F^\delta\MackeyOp S_{f^*\gamma}
  \to \Mackey\pi^{\Lie-\delta}_{-\gamma+n}f_!C^\delta\MackeyOp S_{f^*\gamma} \to 0
\]
for $n\leq 0$. Using the fact that the algebraic functor $f_!$ is right exact, this implies that
\[
 \Mackey\pi^{\Lie-\delta}_{-\gamma+n}f_!C^\delta\MackeyOp S_{f^*\gamma} \iso
  \begin{cases}
   f_!\MackeyOp S &\text{if $n=0$} \\
   0 &\text{if $n<0$}.
  \end{cases}
\]
The spectrum $P^\delta\MackeyOp S_{f^*\gamma}$ is obtained from
$C^\delta\MackeyOp S_{f^*\gamma}$ by killing higher homotopy groups, so we
have the same calculation of the lower homotopy groups of $f_!P^\delta\MackeyOp S_{f^*\gamma}$.
This implies the same calculation for
\[
 f_!H^\delta\MackeyOp S_{f^*\gamma} 
  = f_!(P^\delta\MackeyOp S_{f^*\gamma}\smsh E\F(\delta)_+)
  \hmtpc f_!P^\delta\MackeyOp S_{f^*\gamma} \smsh E\F(\delta)_+.
\]

Now consider the map 
\[
 \Mackey\pi^{\Lie-\delta}_{-\gamma}f_!H^\delta\MackeyOp S_{f^*\gamma} 
   \to \Mackey\pi^{\Lie-\delta}_{-\gamma}H^\delta(f_!\MackeyOp S)_\gamma.
\]
By the calculation above, we can build a model of $H^\delta(f_!\MackeyOp S)_\gamma$
by starting with $f_!H^\delta(\MackeyOp S)_{f^*\gamma}$ and killing higher homotopy groups.
This gives us a map 
$f_!H^\delta(\MackeyOp S)_{f^*\gamma}\to H^\delta(f_!\MackeyOp S)_\gamma$
that is an isomorphism in $\Mackey\pi^{\Lie-\delta}_{-\gamma}$, but is it the map
we're interested in? It's adjoint can be written as
\[
 H^\delta(\MackeyOp S)_{f^*\gamma} \to f^*f_!H^\delta(\MackeyOp S)_{f^*\gamma}
  \to f^*H^\delta(f_!\MackeyOp S)_\gamma.
\]
On applying $\Mackey\pi^{G,\delta}_{f^*\gamma}$, the first map is
the unit $\eta\colon\MackeyOp S\to f^*f_!\MackeyOp S$ and the second is the identity, 
so this map of spectra coincides with the map
$H^\delta(\MackeyOp S)_{f^*\gamma} \to H^\delta(f^*f_!\MackeyOp S)_{f^*\gamma}$
induced by $\eta$.
Therefore, we are looking at the correct map and conclude that
$f_!H^\delta(\MackeyOp S)_{f^*\gamma}\to H^\delta(f_!\MackeyOp S)_\gamma$
does give us an isomorphism on applying
$\Mackey\pi^{\Lie-\delta}_{-\gamma}$, as claimed.

When $\delta$ is complete, the results for cohomology follow from the identity
$H_\delta\Mackey T^\gamma = H^{\Lie-\delta}\Mackey T_{-\gamma}$.
\end{proof}

\begin{remarks}
The change-of-base space isomorphisms suggest two interesting alternatives
to considering the homology of spaces over a specified base space.
\begin{enumerate}
\item
If $p\colon X\to B$ is a $G$-space over $B$, we have $p = p_!(\id)$ where
$\id\colon X\to X$ is $X$ considered as a space over itself. Thus, 
if $\gamma$ is a representation of $\Pi_G B$ and
$\MackeyOp S$ is a covariant $\stab\Pi_{G,\delta}B$-module,
we have
\[
 H^{G,\delta}_\gamma(X;\MackeyOp S) \iso H^{G,\delta}_{p^*\gamma}(X;p^*\MackeyOp S),
\]
where the group on the right considers $X$ as a space over itself,
and similarly in cohomology.
In other words, the homology and cohomology of $X$ as a space over $B$ are
intrinsic to $X$ and not dependent on $B$ at all.
$B$ is simply a convenient carrier for $\gamma$ and $\MackeyOp S$.
If $X$ is an ex-$G$-space over $B$, we can say similarly that the reduced
homology and cohomology of $X$ can be described in terms of the pair
$(X,\sigma(B))$ over $X$, whose dependence on $B$ may be reduced further using excision.

This is the way we originally described $RO(\Pi X)$-graded homology and cohomology
in, for example, \cite{CW:duality}. We could have continued in that vein and
described, for example, homology as being defined on a category consisting of triples
$(X,\gamma,\MackeyOp S)$ and compatible maps,
where $\gamma$ is a representation of $\Pi X$ and
$\MackeyOp S$ is a $\stab\Pi_{G,\delta}X$-module.
We've chosen to recast the theories in terms of parametrized spaces for several reasons,
one being that parametrized spectra seem to be the most natural representing objects.

\item
In the other direction, we could ask whether there is some universal
base space we could choose to use rather than define theories 
for each choice of base space.
Ignoring coefficient systems, we have a positive answer
from Theorem~\ref{thm:classification}:
There is a classifying space $Bv\orthosmall Gn$ such that,
if $B$ has the $G$-homotopy type of a $G$-CW complex, then
$[B,Bv\orthosmall Gn]_G$ is in one-to-one correspondence with the
set of $n$-dimensional virtual representations of $\Pi_G B$.
Put another way, the space $Bv\orthosmall Gn$ carries a universal
$n$-dimensional virtual representation. To consider homology graded
on $n$-dimensional representations, it suffices then to consider spaces
over $Bv\orthosmall Gn$.

To do this properly we should consider the disjoint union of these spaces
over all $n$ and the structures on this union implied by suspension and, more
generally, addition of representations.
We should also look at the spectrum over this space representing homology or cohomology,
and investigate what additional structure it has.

This does not address the issue of coefficient systems---the 
classifying space for representations will certainly not carry all the
coefficient systems we may be interested in.
So, we could ask whether there is a space carrying a universal coefficient system.
However, the collection of possible coefficient systems on a space is not a set
but a proper class, so
without some restriction on the possible coefficient systems this obviously can't work.
We have not investigated whether there is a classifying space under such a restriction
or what other approach might work.

\end{enumerate}
\end{remarks}

\section{Change of groups}

\subsection{Subgroups}

Let $\iota\colon K\to G$ be the inclusion of a (closed) subgroup.
in \cite[\S 14.3]{MaySig:parametrized}, May and Sigurdsson define
the forgetful functor 
\[
 \iota^*\colon G\PreSpec{} B \to K\PreSpec{} B
\]
and, for $A$ a $K$-space, an induction functor 
\[
 \iota_!\colon K\PreSpec{} A\to G\PreSpec{}{G\times_K A},
\]
an analogue of the space level functor that takes a $K$-space
$Z$ to the $G$-space $G\times_K Z$.
In fact, $\iota_!$ is an equivalence of categories with inverse
$\nu^*\circ\iota^*$, where $\nu\colon A\to G\times_K A$ is the
$K$-map $\nu(a) = [e,a]$. If $\epsilon\colon G\times_K B\to B$ is the
map $\epsilon[g,b] = gb$, we let
\[
 i^G_K = \epsilon_! \iota_! \colon K\PreSpec{}B\to G\PreSpec{}B.
\]
(Note that this makes sense only when $B$ is a $G$-space.)
Then $i^G_K$ is left adjoint to $\iota^*$: If $X$ is a $K$-spectrum over $B$
and $E$ is a $G$-spectrum over $B$, then
\[
 [ i^G_K X, E ]_{G,B} \iso [X, (\iota_!)^{-1}\epsilon^* E]_{K,B} 
   \iso [X, \iota^* E]_{K,B}.
\]
The second isomorphism follows from the natural isomorphism
$(\iota_!)^{-1}\epsilon^* \iso \nu^*\iota^*\epsilon^* \iso \iota^*$,
the last isomorphism being a special case of
\cite[14.3.3]{MaySig:parametrized}.
We shall also write $G_+\smsh_K X$ for $i^G_K X$
and $E|K$ for $\iota^* E$. With this notation we have
\[
 [G_+\smsh_K X, E]_{G,B} \iso [X, E|K]_{K,B}.
\]

\begin{definition}
Let $B$ be a $G$-space, let $\delta$ be a dimension function for $G$, and
let $K$ be a subgroup of $G$.
Let
\[
 i_K^G = G_+\smsh_K\susp^{-\delta(G/K)}(-) \colon \stab\Pi_{K,\delta} B\to \stab\Pi_{G,\delta} B.
\]
Its effect on an object $b\colon K/L\to B$ is
\[
 i_K^G(b) = G\times_K b\colon G/L = G\times_K (K/L) \to G\times_K B\to B
\]
and it takes a stable map 
$f\colon K_+\smsh_L S^{-\delta(K/L),b} \to K_+\smsh_M S^{-\delta(K/M),c}$ to the stable map
\begin{align*}
 G_+\smsh_L S^{-\delta(G/L),b}
  &= G_+\smsh_L S^{-\delta(K/L)-\delta(G/K),b} \\
  &= G_+\smsh_K \susp^{-\delta(G/K)} (K_+\smsh_L S^{-\delta(K/L),b}) \\
  &\to G_+\smsh_K \susp^{-\delta(G/K)} (K_+\smsh_M S^{-\delta(K/M),c}) \\
  &= G_+\smsh_M S^{-\delta(G/M),c}.
\end{align*}
If $\Mackey T$ is a contravariant $\stab\Pi_{G,\delta}B$-module, let
\[
 \Mackey T|K = (i_K^G)^*\Mackey T = \Mackey T\circ i_K^G,
\]
a $\stab\Pi_{K,\delta}B$-module,
as in Definition~\ref{def:indres}. 
If $\Mackey C$ is a contravariant $\stab\Pi_{K,\delta}B$-module, let 
\[
 G\times_K \Mackey C = (i_K^G)_!\Mackey C, 
\]
a $\stab\Pi_{G,\delta}B$-module.
 \end{definition}
 
It follows from Proposition~\ref{prop:indresadjunction} that
 \[
 \Hom_{\stab\Pi_{G,\delta} B}(G\times_K \Mackey C, \Mackey T)
 \iso \Hom_{\stab\Pi_{K,\delta} B}(\Mackey C, \Mackey T|K)
 \]
for any contravariant $\stab\Pi_{K,\delta}B$-module $\Mackey C$
and contravariant $\stab\Pi_{G,\delta}B$-module $\Mackey T$.
The same adjunction holds for covariant modules.
We also have the isomorphisms
 \[
 (G\times_K \Mackey C)\tensor_{\stab\Pi_{G,\delta} B} \MackeyOp S
 \iso \Mackey C\tensor_{\stab\Pi_{K,\delta} B} (\MackeyOp S|K)
 \]
if $\MackeyOp S$ is a covariant $\stab\Pi_{G,\delta}B$-module, and
 \[
 (G\times_K \MackeyOp D)\tensor_{\stab\Pi_{G,\delta} B} \Mackey T
  \iso \MackeyOp D \tensor_{\stab\Pi_{K,\delta} B} (\Mackey T|K)
 \]
if $\MackeyOp D$ is a covariant $\stab\Pi_{K,\delta}B$-module.
Finally, we have the calculations
 \[
 G\times_K \Mackey A_b \iso \Mackey A_{G\times_K b}
 \]
and
 \[
 G\times_K \MackeyOp A^{b} \iso \MackeyOp A^{G\times_K b}.
 \]
Here, $b\colon K/L\to B$, $\Mackey A_b = \stab\Pi_{K,\delta} B(-,b)$,
and $\MackeyOp A^b = \stab\Pi_{K,\delta} B(b,-)$.

\begin{definition}
Let $\gamma$ be a representation of $\Pi_G B$ and let
$K$ be a subgroup of $G$. We define $\gamma|K$ to be the
representation of $\Pi_K B$ whose value on
$b\colon K/L\to B$ is the restriction to $K/L$
of $\gamma(G\times_K b)$ over $G/L$.
Its values on morphisms is defined by restriction similarly.
\end{definition}

We can now state the Wirthm\"uller isomorphisms on the chain level
and for homology and cohomology. For simplicity, we write
$\gamma$ again for $\gamma|K$.

We first note the following analogue of Proposition~\ref{prop:genInduction},
shown in the same way.

\begin{proposition}
Let $\delta$ be a dimension function for $G$ and let $K\in\F(\delta)$.
If $X$ is a $\delta$-$K$-$(\gamma-\delta(G/K))$-cell complex, then
$G\times_K X$ is a $\delta$-$G$-$\gamma$-cell complex with corresponding cells.
If $X$ is a CW complex, then so is $G\times_K X$ with this structure.
\qed
\end{proposition}

\begin{proposition}
Let $\delta$ be a dimension function for $G$,
let $K\in\F(\delta)$,
and let $\gamma$ be a virtual representation of $\Pi_{G}B$.
Let $X$ be a based $\delta$-$K$-CW$(\gamma-\delta(G/K))$ complex and give
$G_+\smsh_K X$ the $\delta$-$G$-CW$(\gamma)$ structure from the preceding proposition.
Then
\[
 G\times_K\Mackey C^{K,\delta}_{\gamma-\delta(G/K)+*}(X,B) 
   \iso \Mackey C^{G,\delta}_{\gamma+*}(G_+\smsh_K X,B).
\]
This isomorphism respects suspension in the sense that, if $W$ is a representation of $G$,
then the following diagram commutes:
\[
 \xymatrix{
  G\times_K\Mackey C^{K,\delta}_{\gamma-\delta(G/K)+*}(X,B) \ar[r]^-\iso \ar[d]_{\sigma^W}
   & \Mackey C^{G,\delta}_{\gamma+*}(G_+\smsh_K X,B) \ar[d]^{\sigma^W} \\
  G\times_K\Mackey C^{K,\delta}_{\gamma-\delta(G/K)+W+*}(\susp^W X,B) \ar[r]^-\iso
   & \Mackey C^{G,\delta}_{\gamma+W+*}(G_+\smsh_K \susp^W X,B)
 }
\]  
\end{proposition}

The proof is the same as that of Proposition~\ref{prop:wirthmulleronchains}.

\begin{theorem}[Wirthm\"uller Isomorphisms]
Let $\delta$ be a dimension function for $G$,
let $K\in\F(\delta)$,
and let $\gamma$ be a virtual representation of $\Pi_{G}B$.
Then, for $X$ a $\delta$-$K$-CW$(\gamma-\delta(G/K))$ complex, there are natural isomorphisms
 \[
 \tCH^{G,\delta}_{\gamma}(G_+\smsh_K X;\MackeyOp S)
 \iso  \tCH^{K,\delta}_{\gamma - \delta(G/K)}(X;\MackeyOp S|K)
 \]
and
 \[
 \tCH_{G,\delta}^{\gamma}(G_+\smsh_K X;\Mackey T)
 \iso  \tCH_{K,\delta}^{\gamma-\delta(G/K)}(X;\Mackey T|K).
 \]
These isomorphisms respect suspension in the sense that, if $W$ is a representation of $G$,
then the following diagram commutes:
\[
 \xymatrix{
  \tCH^{G,\delta}_{\gamma}(G_+\smsh_K X;\MackeyOp S) \ar[r]^-{\iso} \ar[d]_{\sigma^W}
    & \tCH^{K,\delta}_{\gamma - \delta(G/K)}(X;\MackeyOp S|K) \ar[d]^{\sigma^W} \\
  \tCH^{G,\delta}_{\gamma+W}(G_+\smsh_K \susp^W X;\MackeyOp S) \ar[r]^-{\iso}
    & \tCH^{K,\delta}_{\gamma + W - \delta(G/K)}(\susp^W X;\MackeyOp S|K)
 }
\]
and similarly for cohomology.
If $\delta$ is familial, for $X$ an ex-$K$-space there are natural isomorphisms
 \[
 \tilde H^{G,\delta}_{\gamma}(G_+\smsh_K X;\MackeyOp S)
 \iso  \tilde H^{K,\delta}_{\gamma - \delta(G/K)}(X;\MackeyOp S|K)
 \]
and
 \[
 \tilde H_{G,\delta}^{\gamma}(G_+\smsh_K X;\Mackey T)
 \iso  \tilde H_{K,\delta}^{\gamma-\delta(G/K)}(X;\Mackey T|K).
 \]
These isomorphisms respect suspension in the sense above.
 \end{theorem}

The proof is essentially the same as that of Theorem~\ref{thm:wirthmuller}.

Recall from Definition~\ref{def:paramMackeyvalued} that we can consider 
cellular homology and cohomology to be $\stab\Pi_\delta B$-module-valued, with
\[
 \MackeyOp H^{G,\delta}_\gamma(X;\MackeyOp S)(b\colon G/K\to B)
  = \tilde H^{G,\delta}_\gamma(X\smsh G_+\smsh_K S^{-\delta(G/K),b};\Mackey S)
\]
and
\[
 \Mackey H_{G,\delta}^\gamma(X;\Mackey T)(b\colon G/K\to B)
  = \tilde H_{G,\delta}^\gamma(X\smsh G_+\smsh_K S^{-\delta(G/K),b};\Mackey T).
\]
We then get the following analogue of Corollary~\ref{cor:MackeyStructure},
whose proof needs only the nonparametrized version of the Wirthm\"uller isomorphism.
This result shows that the module-valued theories capture the local behavior of the
theories on the fibers, not the global values.

\begin{corollary}
If $X$ is an ex-$G$-space, $\delta$ is a familial dimension function for $G$,
$\MackeyOp S$ is a covariant $\stab\Pi_\delta B$-module, and
$\Mackey T$ is a contravariant $\stab\Pi_\delta B$-module,
then, for $K\in\F(\delta)$, we have
\[
 \MackeyOp H^{G,\delta}_\gamma(X;\MackeyOp S)(b\colon G/K\to B) \iso
  \tilde H^{K,\delta}_{\gamma_0(b)}(b^*X;b^*\MackeyOp S)
\]
and
\[
 \Mackey H_{G,\delta}^\gamma(X;\Mackey T)(b\colon G/K\to B) \iso
  \tilde H_{K,\delta}^{\gamma_0(b)}(b^*X;b^*\Mackey T|K),
\]
where we implicitly identify $G\PreSpec{}{G/K}$ with $K\PreSpec{}{}$.
Recall that $b^* X$ must be understood as the derived functor, the
homotopy pullback.
\end{corollary}

\begin{proof}
We have
\begin{align*}
 \MackeyOp H^{G,\delta}_\gamma(X;\MackeyOp S)(b\colon G/K\to B)
  &= \tilde H^{G,\delta}_\gamma(X\smsh G_+\smsh_K S^{-\delta(G/K),b};\Mackey S) \\
  &= \tilde H^{G,\delta}_\gamma(X\smsh_B b_!(G\times_K S^{-\delta(G/K)});\Mackey S) \\
  &\iso \tilde H^{G,\delta}_\gamma(b_!(b^* X \smsh_{G/K} (G\times_K S^{-\delta(G/K)});\Mackey S) \\
  &\iso \tilde H^{G,\delta}_{\gamma(b)}(b^* X\smsh_{G/K}(G\times_K S^{-\delta(G/K)});b^*\Mackey S) \\
  &\iso \tilde H^{K,\delta}_{\gamma_0(b)}(b^* X; b^*\Mackey S).
\end{align*}
The proof for cohomology is similar.
\end{proof}

\begin{corollary}
Let $\delta$ be a familial dimension function for $G$ and let $K\in\F(\delta)$.
If $X$ is an ex-$G$-space,
$\MackeyOp S$ is a covariant $\stab\Pi_\delta B$-module, and
$\Mackey T$ is a contravariant $\stab\Pi_\delta B$-module,
then we have
\[
 \MackeyOp H^{G,\delta}_\gamma(X;\MackeyOp S)|K \iso
  \MackeyOp H^{K,\delta}_\gamma(X;\MackeyOp S|K)
\]
and
\[
 \Mackey H_{G,\delta}^\gamma(X;\Mackey T)|K \iso
  \Mackey H_{K,\delta}^\gamma(X;\Mackey T|K).
\]
\qed
\end{corollary}

On the spectrum level, what underlies the Wirthm\"uller isomorphisms is the
following fact.

\begin{proposition}
Let $H_\delta\Mackey T^\gamma$ be a parametrized Eilenberg-Mac\,Lane spectrum
and let $K$ be a subgroup of $G$.
Then 
\[
 (H_\delta\Mackey T^\gamma)|K \hmtpc \susp_K^{-\delta(G/K)}H_\delta(\Mackey T|K)^\gamma.
\]
Similarly, for a covariant $\MackeyOp S$, we have
\[
 (H^\delta\MackeyOp S_\gamma)|K \hmtpc \susp_K^{-(\Lie(G/K)-\delta(G/K))} H^\delta(\MackeyOp S|K)_\gamma.
\]
\end{proposition}

\begin{proof}
We've already mentioned that parametrized Eilenberg-Mac\,Lane spectra are
characterized by their fibers being non-parametrized Eilenberg-Mac\,Lane spectra.
This proposition then follows from Proposition~\ref{prop:emsubgroups}.
\end{proof}

On the represented level, the cohomology Wirthm\"uller isomorphism is then the adjunction
\begin{align*}
 [G_+\smsh_K \susp_B^\infty X, \susp_B^V H_\delta\Mackey T^\gamma]_{G,B}
   &\iso [\susp_B^\infty X, \susp_B^V (H_\delta\Mackey T^\gamma)|K]_{K,B} \\
   &\iso [\susp_B^\infty X, \susp_B^{V-\delta(G/K)} H_\delta(\Mackey T|K)^\gamma]_{K,B}.
\end{align*}
For homology, we have
\begin{align*}
 [S^V,{} &\rho_!(H^\delta\MackeyOp S_\gamma \smsh_B (G_+\smsh_K X))]_{G,B} \\
  &\iso [S^V, G_+\smsh_K \rho_!((H^\delta\MackeyOp S_\gamma)|K \smsh_B X)]_{G,B} \\
  &\iso [S^V, \susp^{\Lie(G/K)} \rho_!((H^\delta\MackeyOp S_\gamma)|K \smsh_B X)]_{K,B} \\
  &\iso [S^V, \rho_!(\susp^{\Lie(G/K)}\susp^{-(\Lie(G/K)-\delta(G/K))}H^\delta(\MackeyOp S|K)_{\gamma} \smsh_B X)]_{K,B} \\
  &= [S^V, \rho_!(\susp^{\delta(G/K)}H^\delta(\MackeyOp S|K)_{\gamma} \smsh_B X)]_{K,B} \\
  &\iso [S^{V-\delta(G/K)}, \rho_!(H^\delta(\MackeyOp S|K)_{\gamma} \smsh_B X)]_{K,B}.
\end{align*}

As in the non-parametrized case, we can use the Wirthm\"uller isomorphism to define
restriction to subgroups.
In the case of cohomology, we use the composite
\[
 \tilde H_{G,\delta}^\gamma(X;\Mackey T)
  \to \tilde H_{G,\delta}^\gamma(G/K_+\smsh X;\Mackey T)
  \iso \tilde H_{K,\delta}^{\gamma-\delta(G/K)}(X;\Mackey T|K).
\]
Here, the first map is induced by the projection $G/K_+\smsh X\to X$ and
the second is the Wirthm\"uller isomorphism.

In homology, we take an embedding of $G/K$ in a representation $V$, with
$c\colon S^V\to G_+\smsh_K S^{V-\Lie(G/K)}$ being the induced collapse map,
and then use the composite
\begin{align*}
 \tilde H_\gamma^{G,\delta}(X;\MackeyOp S)
  &\iso \tilde H_{\gamma+V}^{G,\delta}(\susp^V X; \MackeyOp S) \\
  &\to \tilde H_{\gamma+V}^{G,\delta}(G_+\smsh_K\susp^{V-\Lie(G/K)} X; \MackeyOp S) \\
  &\iso \tilde H_{\gamma+V-\delta(G/K)}^{K,\delta}(\susp^{V-\Lie(G/K)} X; \MackeyOp S|K) \\
  &\iso \tilde H_{\gamma+\Lie(G/K)-\delta(G/K)}^{K,\delta}(X;\MackeyOp S|K).
\end{align*}

As before, we write $a \mapsto a|K$ for this restriction in either cohomology or homology.

On the represented level, these maps are given simply by restriction of stable $G$-maps
to stable $K$-maps. Explicitly, the cohomology restriction is the following:
\begin{align*}
 \tilde H_{G,\delta}^\gamma(X;\Mackey T)
  &\iso [\susp_B^\infty X, H_\delta\Mackey T^\gamma]_{G,B} \\
  &\to [\susp_B^\infty X, (H_\delta\Mackey T^\gamma)|K]_{K,B} \\
  &\iso [\susp_B^\infty X, \susp^{-\delta(G/K)} H_\delta(\Mackey T|K)^\gamma]_{K,B} \\
  &\iso \tilde H_{K,\delta}^{\gamma-\delta(G/K)} (X;\Mackey T|K).
\end{align*}
The homology restriction is the following:
\begin{align*}
 \tilde H^{G,\delta}_\gamma(X;\MackeyOp S)
  &\iso [S, \rho_!( H^\delta\MackeyOp S_\gamma \smsh_B X)]_{G,B} \\
  &\to [S, \rho_!( (H^\delta\MackeyOp S_\gamma)|K \smsh_B X)]_{K,B} \\
  &\iso [S, \rho_!( \susp^{-(\Lie(G/K)-\delta(G/K)} H^{\delta}(\MackeyOp S|K)_{\gamma}
      \smsh_B X)]_{K,B} \\
  &\iso [S^{\Lie(G/K)-\delta(G/K)}, \rho_!(H^{\delta}(\MackeyOp S|K)_{\gamma}
      \smsh_B X)]_{K,B} \\
  &\iso \tilde H^{K,\delta}_{\gamma+\Lie(G/K)-\delta(G/K)}(X; \MackeyOp S|K).
\end{align*}

\begin{remark}
The specializations of these maps to the ordinary and dual theories give us the
following restriction maps:
\begin{align*}
 \tilde H^G_\gamma(X;\MackeyOp S) &{}\to \tilde H^K_{\gamma+\Lie(G/K)}(X;\MackeyOp S|K) 
  & \tilde\H_G^\gamma(X;\MackeyOp S) &{}\to \tilde\H_K^{\gamma-\Lie(G/K)}(X;\MackeyOp S|K) \\
 \tilde H_G^\gamma(X;\Mackey T) &{}\to \tilde H_K^{\gamma}(X;\Mackey T|K)
  & \tilde\H^G_\gamma(X;\Mackey T) &{}\to \tilde\H^K_{\gamma}(X;\Mackey T|K) \\
\end{align*}
\end{remark}

\subsection{Quotient groups}

Let $N$ be a normal subgroup of $G$ and let $\epsilon\colon G\to G/N$ denote the quotient map.
See \cite[\S 14.4]{MaySig:parametrized} for some of the functors on parametrized spectra we use here.

Similar to Proposition~\ref{prop:inducedCWstructure}, if
$\delta$ is a dimension function for $G$,
$A$ is a $G/N$-space,
$\gamma$ is a virtual representation of $\Pi_{G/N}A$,
and $Y$ is a based $\delta$-$G/N$-CW$(\gamma)$ complex
over $A$, then $\epsilon^*Y$ is a based
$\delta$-$G$-CW$(\gamma)$ complex over $\epsilon^* A$.
We will usually not write $\epsilon^*$ when the meaning is clear without it.

Similar to Proposition~\ref{prop:CWfixedsets}, if
$B$ is a $G$-space,
$\gamma$ is a virtual representation of $\Pi_G B$,
and $X$ is a based $\delta$-$G$-CW$(\gamma)$ complex over $B$, then
$X^N$ is a based $\delta$-$G/N$-CW$(\gamma^N)$ complex over $B^N$.

To describe the associated algebra on chains, let
\[
 \theta\colon \stab\Pi_{G/N,\delta}A \to \stab\Pi_{G,\delta}A
\]
be the restriction of $\epsilon^\sharp$, the functor that assigns $G/N$-spectra over $A$
the induced $G$ action and expands from an $N$-trivial universe to a complete $G$-universe.

If $E$ is a $G$-spectrum over a $G$-space $B$, then $E^N$ is a $G/N$-spectrum over $B^N$.
We define the geometric fixed-point construction
\[
 \Phi^N(E) = (\tilde E\F[N]\smsh E)^N
\]
as in Definition~\ref{def:geoFixedPoints}.
In this context, $\Phi^N$ continues to enjoy the properties listed
after that definition.
As in Remark~\ref{rem:zeroOrbit}, consider $\stab\Pi_{G,\delta} B$ as augmented with
a zero object $*$ given by the trivial spectrum.
We write 
\[
 \Phi^N\colon \stab\Pi_{G,\delta} B\to \stab\Pi_{G/N,\delta}B^N
\]
for the restriction
of $\Phi^N$ to the augmented stable fundamental groupoid.

\begin{definition}
If $\MackeyOp S$ is a covariant and $\Mackey T$ a contravariant
$\stab\Pi_{G,\delta}B$-module, let
\begin{align*}
 \MackeyOp S^N &= \Phi^N_!\MackeyOp S \qquad\text{and}\\
 \Mackey T^N &= \Phi^N_!\Mackey T.
\end{align*}
We call this the {\em $N$-fixed point functor.}
If $\MackeyOp U$ is a covariant and $\Mackey V$ is a contravariant
$\stab\Pi_{G/N,\delta}B^N$-module, let
\begin{align*}
 \Inf_{G/N}^G \MackeyOp U &= (\Phi^N)^*\MackeyOp U \qquad\text{and}\\
 \Inf_{G/N}^G \Mackey V &= (\Phi^N)^*\Mackey V.
\end{align*}
We call this the {\em inflation functor.}
 \end{definition}

For $b\colon G/L\to B$, write
\[
 \MackeyOp A^b_G = \stab\Pi_{G,\delta}B(b,-)
\]
and
\[
 \Mackey A_b^G = \stab\Pi_{G,\delta}B(-,b).
\]
We have the following analogue of Proposition~\ref{prop:burnsidefixedsets}.

\begin{proposition}\label{prop:Xburnsidefixedsets}
If $N$ is normal in $G$ and $b\colon G/L\to B$, then
 \[
 (\MackeyOp A^b_G)^N 
  \iso \begin{cases} \MackeyOp A^b_{G/N} & \text{if $N\leq L$} \\
                  0 & \text{if $N\not\leq L$}
    \end{cases}
 \]
and
 \[
 (\Mackey A_b^G)^N 
  \iso \begin{cases} \Mackey A_b^{G/N} & \text{if $N\leq L$} \\
                  0 & \text{if $N\not\leq L$.}
    \end{cases}
 \]
If $B$ is a $G/N$-space, $N\leq L\leq G$, and $b\colon (G/N)/(G/L)\to B$, then
\[
 \theta_! \MackeyOp A^b_{G/N} \iso \MackeyOp A^b_G
\]
and
\[
 \theta_! \Mackey A_b^{G/N} \iso \Mackey A_b^G,
\]
where we write $b$ again for $\epsilon^* b\colon G/L\to B$.
\end{proposition}

\begin{proof}
These are special cases of Proposition~\ref{prop:indresadjunction}(\ref{item:indFree}).
 \end{proof}

We now have the following chain-level calculations, proved in the same way
as Propositions~\ref{prop:inducedchains} and~\ref{prop:fixedSetChainIso}.

\begin{proposition}\label{prop:paraminducedchains}
Let $N$ be a normal subgroup of $G$, let $\delta$ be a dimension function for $G$,
let $B$ be a $G/N$-space,
and let $\gamma$ be a virtual representation of $\Pi_{G/N}B$.
Let $Y$ be a based $\delta$-$G/N$-CW$(\gamma)$ complex over $B$ and give
$\epsilon^* Y$ the corresponding $\delta$-$G$-CW$(\gamma)$ structure. 
Then we have a natural chain isomorphism
\[
 \theta_!\Mackey C^{G/N,\delta}_{\gamma+*}(Y,*) \iso
  \Mackey C^{G,\delta}_{\gamma+*}(Y,*),
\]
writing $Y$ for $\epsilon^*Y$ on the right as usual.
This isomorphism respects suspension in the sense that, if $W$ is a representation
of $G/N$, then the following diagram commutes:
\[
 \xymatrix{
   \theta_!\Mackey C^{G/N,\delta}_{\gamma+*}(Y,*) \ar[r]^\iso \ar[d]_{\sigma^W}
    & \Mackey C^{G,\delta}_{\gamma+*}(Y,*) \ar[d]^{\sigma^W} \\
   \theta_!\Mackey C^{G/N,\delta}_{\gamma+W+*}(\susp^W Y,*) \ar[r]_\iso
    & \Mackey C^{G,\delta}_{\gamma+W+*}(\susp^W Y,*).
 }
\]
\qed
\end{proposition}

\begin{proposition}\label{prop:paramfixedSetChainIso}
Let $N$ be a normal subgroup of $G$, let $\delta$ be a dimension function for $G$,
let $B$ be a $G$-space,
and let $\gamma$ be a virtual representation of $\Pi_{G}B$.
Let $X$ be a based $\delta$-$G$-CW$(\gamma)$ complex over $B$
and give $X^N$ the corresponding $\delta$-$G/N$-CW$(\gamma^N)$
structure.
With this structure, we have a natural chain isomorphism
\[
 \Mackey C^{G,\delta}_{\gamma+*}(X,*)^N \iso \Mackey C^{G/N,\delta}_{\gamma^N+*}(X^N,*).
\]
This isomorphism respects suspension in the sense that, if $W$ is a representation $G$,
then the following diagram commutes:
\[
 \xymatrix{
  \Mackey C^{G,\delta}_{\gamma+*}(X,*)^N \ar[r]^-\iso \ar[d]_{(\sigma^W)^N}
   & \Mackey C^{G/N,\delta}_{\gamma^N+*}(X^N,*) \ar[d]^{\sigma^{W^N}} \\
  \Mackey C^{G,\delta}_{\gamma+W+*}(\susp^W X,*)^N \ar[r]_-\iso
   & \Mackey C^{G/N,\delta}_{\gamma^N+W^N+*}(\susp^{W^N}X^N,*)
 }
\]
\qed
\end{proposition}

The chain isomorphisms lead to the following analogues of
Theorems~\ref{thm:inducedHomology} and~\ref{thm:fixedSetIso}, with similar proofs.

\begin{theorem}\label{thm:paraminducedHomology}
Let $N$ be a normal subgroup of $G$,
let $\delta$ be a dimension function for $G$,
let $B$ be a $G/N$-space,
let $\gamma$ be a virtual representation of $\Pi_{G/N}B$,
let $\MackeyOp S$ be a covariant $\stab\Pi_{G,\delta}B$-module,
and let $\Mackey T$ be a contravariant $\stab\Pi_{G,\delta}B$-module.
Then, for $Y$ a $\delta$-$G/N$-CW$(\gamma)$ complex over $B$, we have natural isomorphisms
\begin{align*}
 \tCH^{G,\delta}_{\gamma}(Y;\MackeyOp S)
   &\iso \tCH^{G/N,\delta}_{\gamma}(Y; \theta^* \MackeyOp S) \qquad\text{and} \\
 \tCH_{G,\delta}^{\gamma}(Y;\Mackey T)
   &\iso \tCH_{G/N,\delta}^{\gamma}(Y; \theta^*\Mackey T).
\end{align*}
These isomorphisms respect suspension in the sense that, if
$W$ is a representation of $G/N$,
then the following diagram commutes:
\[
 \xymatrix{
  \tCH^{G,\delta}_{\gamma}(Y;\MackeyOp S) \ar[r]^-\iso \ar[d]_{\sigma^W}
   & \tCH^{G/N,\delta}_{\gamma}(Y; \theta^*\MackeyOp S) \ar[d]^{\sigma^{W}} \\
  \tCH^{G,\delta}_{\gamma+W}(\susp^W Y;\MackeyOp S) \ar[r]_-\iso
   & \tCH^{G/N,\delta}_{\gamma+W}(\susp^{W}Y; \theta^*\MackeyOp S)
 }
\]
and similarly for cohomology.
If $\delta$ is $N$-closed (as in Definition~\ref{def:Nclosed}) and familial,
then, for $Y$ an ex-$G/N$-space over $B$, we have natural isomorphisms
\begin{align*}
 \tilde H^{G,\delta}_{\gamma}(Y;\MackeyOp S)
   &\iso \tilde H^{G/N,\delta}_{\gamma}(Y; \theta^*\MackeyOp S) \qquad\text{and} \\
 \tilde H_{G,\delta}^{\gamma}(Y;\Mackey T)
   &\iso \tilde H_{G/N,\delta}^{\gamma}(Y; \theta^*\Mackey T).
\end{align*}
These isomorphisms respect suspension in the sense that, if
$Y$ is well-based and $W$ is a representation of $G/N$,
then the following diagram commutes:
\[
 \xymatrix{
  \tilde H^{G,\delta}_{\gamma}(Y;\MackeyOp S) \ar[r]^-\iso \ar[d]_{\sigma^W}
   & \tilde H^{G/N,\delta}_{\gamma}(Y; \theta^*\MackeyOp S) \ar[d]^{\sigma^{W}} \\
  \tilde H^{G,\delta}_{\gamma+W}(\susp^W Y;\MackeyOp S) \ar[r]_-\iso
   & \tilde H^{G/N,\delta}_{\gamma+W}(\susp^{W}Y; \theta^*\MackeyOp S)
 }
\]
and similarly for cohomology.
\qed
\end{theorem}

\begin{theorem}\label{thm:paramFixedSetIso}
Let $\delta$ be a dimension function for $G$,
let $N$ be a normal subgroup of $G$, 
let $B$ be a $G$-space,
let $\MackeyOp S$ be a covariant $\stab\Pi_{G/N,\delta}B^N$-module,
and let $\Mackey T$ be a contravariant $\stab\Pi_{G/N,\delta}B^N$-module.
Then, for $X$ a based $\delta$-$G$-CW$(\gamma)$ complex over $B$, we have natural isomorphisms
\begin{align*}
 \tCH^{G,\delta}_{\gamma}(X;\Inf_{G/N}^G \MackeyOp S)
   &\iso \tCH^{G/N,\delta}_{\gamma^N}(X^N; \MackeyOp S) \qquad\text{and} \\
 \tCH_{G,\delta}^{\gamma}(X;\Inf_{G/N}^G \Mackey T)
   &\iso \tCH_{G/N,\delta}^{\gamma^N}(X^N; \Mackey T).
\end{align*}
These isomorphisms respect suspension in the sense that, if
$W$ is a representation of $G$,
then the following diagram commutes:
\[
 \xymatrix{
  \tCH^{G,\delta}_{\gamma}(X;\Inf_{G/N}^G \MackeyOp S) \ar[r]^-\iso \ar[d]_{\sigma^W}
   & \tCH^{G/N,\delta}_{\gamma^N}(X^N; \MackeyOp S) \ar[d]^{\sigma^{W^N}} \\
  \tCH^{G,\delta}_{\gamma+W}(\susp^W X;\Inf_{G/N}^G \MackeyOp S) \ar[r]_-\iso
   & \tCH^{G/N,\delta}_{\gamma^N+W^N}(\susp^{W^N}X^N; \MackeyOp S)
 }
\]
and similarly for cohomology.
If $\delta$ is familial, then, for $X$ an ex-$G$-space over $B$, we have natural isomorphisms
\begin{align*}
 \tilde H^{G,\delta}_{\gamma}(X;\Inf_{G/N}^G \MackeyOp S)
   &\iso \tilde H^{G/N,\delta}_{\gamma^N}(X^N; \MackeyOp S) \qquad\text{and} \\
 \tilde H_{G,\delta}^{\gamma}(X;\Inf_{G/N}^G \Mackey T)
   &\iso \tilde H_{G/N,\delta}^{\gamma^N}(X^N; \Mackey T).
\end{align*}
These isomorphisms respect suspension in the sense that, if
$X$ is well-based and $W$ is a representation of $G$,
then the following diagram commutes:
\[
 \xymatrix{
  \tilde H^{G,\delta}_{\gamma}(X;\Inf_{G/N}^G \MackeyOp S) \ar[r]^-\iso \ar[d]_{\sigma^W}
   & \tilde H^{G/N,\delta}_{\gamma^N}(X^N; \MackeyOp S) \ar[d]^{\sigma^{W^N}} \\
  \tilde H^{G,\delta}_{\gamma+W}(\susp^W X;\Inf_{G/N}^G \MackeyOp S) \ar[r]_-\iso
   & \tilde H^{G/N,\delta}_{\gamma^N+W^N}(\susp^{W^N}X^N; \MackeyOp S)
 }
\]
and similarly for cohomology.
\qed
\end{theorem}

\begin{remarks}
\begin{enumerate}\item[]
\item
As in the nonparametrized case, these isomorphisms 
are compatible, in the sense that,
if $B$ is a $G/N$-space,
$\gamma$ is a virtual representation of $\Pi_{G/N}B$,
$X$ is an ex-$G/N$-space over $B$, and
$\Mackey T$ is a $\stab\Pi_{G/N,\delta}B$-module, then the following
composite is the identity:
\begin{align*}
 \tilde H_\gamma^{G/N,\delta}(X;\Mackey T)
  &= \tilde H_{\gamma^N}^{G/N,\delta}(X^N;\Mackey T) \\
  &\iso \tilde H_\gamma^{G,\delta}(X;\Inf_{G/N}^G\Mackey T) \\
  &\iso \tilde H_\gamma^{G/N,\delta}(X;\theta^*\Inf_{G/N}^G\Mackey T) \\
  &\iso \tilde H_\gamma^{G/N,\delta}(X;\Mackey T).
\end{align*}
(Here, we use the isomorphism $\theta^*\Inf^G_{G/N}\Mackey T \iso \Mackey T$
that we noted first in the nonparametrized case.)
The similar statement for homology is also true.

\item
The two isomorphisms combine to give a third isomorphism, if
$B$ is a $G$-space,
$\gamma$ is a virtual representation of $\Pi_G B$,
$X$ is an ex-$G$-space over $B$,
and $\Mackey T$ is a $\stab\Pi_{G/N,\delta}B^N$-module:
\begin{align*}
 \tilde H_\gamma^{G,\delta}(X;\Inf_{G/N}^G\Mackey T)
  &\iso \tilde H_{\gamma^N}^{G/N,\delta}(X^N; \Mackey T) \\
  &\iso \tilde H_{\gamma^N}^{G/N,\delta}(X^N; \theta^*\Inf_{G/N}^G\Mackey T) \\
  &\iso \tilde H_{\gamma^N}^{G,\delta}(X^N; \Inf_{G/N}^G\Mackey T).
\end{align*}
In the last group, we consider $X^N$ as an ex-$G$-space over $B^N$.
We get a similar isomorphism in homology.

\end{enumerate}
\end{remarks}

We can now define induction and restriction to fixed sets.

\begin{definition}\label{def:paraminduction}
Let $\delta$ be an $N$-closed familial dimension function for $G$,
let $B$ be a $G/N$-space,
let $\gamma$ be a virtual representation of $\Pi_{G/N}B$,
let $\MackeyOp S$ be a covariant $\stab\Pi_{G/N,\delta}B$-module, and
let $\Mackey T$ be a contravariant $\stab\Pi_{G/N,\delta}B$-module.
If $Y$ is an ex-$G/N$-space over $B$,
we define {\em induction from $G/N$ to $G$} to be the composites
\[
 \epsilon^*\colon \tilde H^{G/N,\delta}_\gamma(Y;\MackeyOp S)
  \to \tilde H^{G/N,\delta}_\gamma(Y;\theta^*\theta_!\MackeyOp S)
  \iso \tilde H^{G,\delta}_\gamma(Y;\theta_!\MackeyOp S)
\]
and
\[
 \epsilon^*\colon \tilde H_{G/N,\delta}^\gamma(Y;\Mackey T)
  \to \tilde H_{G/N,\delta}^\gamma(Y;\theta^*\theta_!\Mackey T)
  \iso \tilde H_{G,\delta}^\gamma(Y;\theta_!\Mackey T).
\]
The first map in each case is induced by the unit of the $\theta_!$-$\theta^*$ adjunction.
\end{definition}

\begin{definition}\label{def:paramRestrictToFixed}
Let $\delta$ be a familial dimension function for $G$,
let $B$ be a $G$-space,
let $\gamma$ be a virtual representation of $\Pi_G B$,
let $\MackeyOp S$ be a covariant $\stab\Pi_{G,\delta}B$-module,
and let $\Mackey T$ be a contravariant $\stab\Pi_{G,\delta}B$-module
If $X$ is an ex-$G$-space over $B$ and
$N$ is a normal subgroup of $G$, define
{\em restriction to fixed sets} to be the composites
\[
 (-)^N\colon \tilde H^{G,\delta}_{\gamma}(X;\MackeyOp S) \to
 \tilde H^{G,\delta}_{\gamma}(X;\Inf_{G/N}^G\MackeyOp S^N) \xrightarrow{\iso}
 \tilde H^{G/N,\delta}_{\gamma^N}(X^N;\MackeyOp S^N)
\]
and
\[
 (-)^N\colon \tilde H_{G,\delta}^{\gamma}(X;\Mackey T) \to
 \tilde H_{G,\delta}^{\gamma}(X;\Inf_{G/N}^G\Mackey T^N) \xrightarrow{\iso}
 \tilde H_{G/N,\delta}^{\gamma^N}(X^N;\Mackey T^N).
\]
The first map in each case is induced by the unit of the $(-)^N$-$\Inf_{G/N}^G$ adjunction.
\end{definition}

As in the nonparametrized case, Theorems~\ref{thm:paraminducedHomology}
and~\ref{thm:paramFixedSetIso} imply that
both induction and restriction to fixed sets
respect suspension.
We also have that induction followed by restriction to fixed sets is the identity.

On the represented level, we get the following results. 

\begin{proposition}
Let $N$ be a normal subgroup of $G$, let $\delta$ be an $N$-closed familial dimension
function for $G$,
let $B$ be a $G/N$-space,
let $\gamma$ be a virtual representation of $\Pi_{G/N} B$,
let $\MackeyOp S$ be a covariant $\stab\Pi_{G,\delta}B$-module
and let $\Mackey T$ be a contravariant $\stab\Pi_{G,\delta}B$-module.
Then
\[
 (H^\delta\MackeyOp S_\gamma)^N \hmtpc H^\delta(\theta^*\MackeyOp S)_\gamma
\]
and
\[
 (H_\delta\Mackey T^\gamma)^N \hmtpc H_\delta(\theta^*\Mackey T)^\gamma.
\]
\end{proposition}

\begin{proof}
This follows from Proposition~\ref{prop:inductEMspectra} applied fiberwise.
For example, if $N\leq L\leq G$ and $b\colon G/L\to B$, we have
\begin{align*}
 b^*((H_\delta\Mackey T^\gamma)^N)
  &\hmtpc (b^*H_\delta\Mackey T^\gamma)^N \\
  &\hmtpc (\susp^{\gamma_0(b)}H_\delta(b^*\Mackey T))^N \\
  &\hmtpc \susp^{\gamma_0(b)}H_\delta(\theta^*b^*\Mackey T) \\
  &\hmtpc \susp^{\gamma_0(b)}H_\delta(b^*\theta^*\Mackey T) \\
  &\hmtpc b^*(H_\delta(\theta^*\Mackey T)^\gamma).
\end{align*}
The rest of the characterization of Eilenberg-Mac\,Lane spectra goes as
in the proof of Proposition~\ref{prop:inductEMspectra}.
\end{proof}

The isomorphisms of Theorem~\ref{thm:paraminducedHomology} are then represented as follows,
for $Y$ an ex-$G/N$-space over the $G/N$-space $B$:
\begin{align*}
 [S, \rho_!(H^\delta\MackeyOp S_\gamma \smsh_B Y)]_G
  &\iso [S, (\rho_!(H^\delta\MackeyOp S_\gamma \smsh_B Y))^N]_{G/N} \\
  &\iso [S, \rho_!((H^\delta\MackeyOp S_\gamma \smsh_B Y)^N)]_{G/N} \\
  &\iso [S, \rho_!((H^\delta\MackeyOp S_\gamma)^N\smsh_B Y)]_{G/N} \\
  &\iso [S, \rho_!(H^\delta(\theta^*\MackeyOp S)_\gamma \smsh_B Y)]_{G/N}
\end{align*}
for homology, where the second isomorphism comes from \cite[14.4.4]{MaySig:parametrized}, and
\begin{align*}
 [\susp^\infty_B Y, H_\delta\Mackey T^\gamma]_{G,B}
  &\iso [\susp^\infty_B Y, (H_\delta\Mackey T^\gamma)^N]_{G/N,B} \\
  &\iso [\susp^\infty_B Y, H_\delta(\theta^*\Mackey T)^\gamma]_{G/N,B}
\end{align*}
for cohomology.

\begin{proposition}\label{prop:paramFixedEMSpectra}
Let $\delta$ be a familial dimension function for $G$,
let $N$ be a normal subgroup of $G$,
let $B$ be a $G$-space,
let $\gamma$ be a virtual representation of $\Pi_G B$,
let $\Mackey T$ be a contravariant
$\stab\Pi_{G/N,\delta}B^N$-module, and let $\MackeyOp S$ be a covariant
$\stab\Pi_{G/N,\delta}B^N$-module.
If $i\colon B^N\to B$ is the inclusion, then
\[
 i_!(\epsilon^\sharp H^\delta\MackeyOp S_{\gamma^N}) \smsh \tE\F[N] \hmtpc
   H^\delta(\Inf_{G/N}^G\MackeyOp S)_\gamma
\]
and
\[
 i_!(\epsilon^\sharp H_\delta\Mackey T^{\gamma^N})\smsh \tE\F[N] \hmtpc
   H_\delta(\Inf_{G/N}^G\Mackey T)^\gamma.
\]
Therefore,
\[
 (H^\delta(\Inf_{G/N}^G\MackeyOp S)_\gamma)^N 
   \hmtpc \Phi^N H^\delta(\Inf_{G/N}^G\MackeyOp S)_\gamma 
   \hmtpc H^\delta\MackeyOp S_{\gamma^N}
\]
and
\[
 (H_\delta(\Inf_{G/N}^G\Mackey T)^\gamma)^N 
  \hmtpc \Phi^N H_\delta(\Inf_{G/N}^G\Mackey T)^\gamma 
  \hmtpc H_\delta\Mackey T^{\gamma^N}.
\]
\end{proposition}

\begin{proof}
If $b\colon G/K\to B$ with $N\not\leq K$, then
$(\Inf_{G/N}^G\MackeyOp S)(b) = 0$.
Combined with the equivalence
\[
 S^{\gamma(b)}\smsh \tE\F[N] \hmtpc S^{\gamma(b)^N}\smsh\tE\F[N],
\]
the first equivalence follows from Proposition~\ref{prop:fixedEMSpectra} applied fiberwise.
A similar argument applies to the second. These imply the last two equivalences
as in the proof of Proposition~\ref{prop:fixedEMSpectra}.
\end{proof}

The isomorphisms of Theorem~\ref{thm:paramFixedSetIso} are then represented as
\begin{align*}
 [S, \rho_!(H^{\delta}(\Inf_{G/N}^G\MackeyOp S)_\gamma \smsh_B X)]_{G}
  &\iso [S, \rho_!(H^{\delta}(\Inf_{G/N}^G\MackeyOp S)_\gamma \smsh_B X)^N]_{G/N} \\
  &\iso [S, \rho_!\Phi^N((H^{\delta}\Inf_{G/N}^G\MackeyOp S)_\gamma \smsh_B X)]_{G/N} \\
  &\iso [S, \rho_!(\Phi^N H^{\delta}(\Inf_{G/N}^G\MackeyOp S)_\gamma \smsh_{B^N} X^N)]_{G/N} \\
  &\iso [S, \rho_!(H^{\delta}\MackeyOp S_{\gamma^N} \smsh_{B^N} X^N)]_{G/N}
\end{align*}
and
\begin{align*}
 [\susp_G^\infty X, (H_\delta\Inf_{G/N}^G\Mackey T)^\gamma]_{G,B}
  &\iso [\Phi^N(\susp_G^\infty X), \Phi^N H_\delta(\Inf_{G/N}^G\Mackey T)^\gamma]_{G/N,B^N} \\
  &\iso [\susp_{G/N}^\infty X^N, H_\delta\Mackey T^{\gamma^N}]_{G/N,B^N}.
\end{align*}

Finally, we can define restriction to $K$-fixed sets
when $K$ is not normal by first restricting to the normalizer $NK$ and then taking
$K$-fixed sets. As in the nonparametrized case, the possible dimension shift goes away and we get
the following fixed set maps:
\begin{align*}
 (-)^K\colon \tilde H^{G,\delta}_\gamma(X;\MackeyOp S)
    &\to \tilde H^{WK,\delta}_{\gamma^K}(X^K;\MackeyOp S^K)\qquad\text{and}\\
 (-)^K\colon \tilde H_{G,\delta}^\gamma(X;\Mackey T)
    &\to \tilde H_{WK,\delta}^{\gamma^K}(X^K;\Mackey T^K)
\end{align*}
where $WK = NK/K$.

\subsection{Subgroups of quotient groups}

As in the nonparametrized case, we now look at how induction and restriction
to fixed sets interact with restriction to subgroups. Once again, the main result
needed is that the coefficient systems involved agree.

Suppose that $N$ is a normal subgroup of $G$, $N\leq L \leq G$, and
$\delta$ is familial with $L\in\F(\delta)$.
We then have the following commutative diagrams:
\[
 \xymatrix{
  \stab\Pi_{L/N,\delta}B^N \ar[r]^-\theta \ar[d]_{i_{L/N}^{G/N}}
     & \stab\Pi_{L,\delta}B \ar[d]^{i_L^G} \\
  \stab\Pi_{G/N,\delta}B^N \ar[r]_-\theta 
     & \stab\Pi_{G,\delta}B
 }
\]
and
\[
 \xymatrix{
  \stab\Pi_{L,\delta}B \ar[d]_{i_L^G} \ar[r]^-{\Phi^N} 
    & \stab\Pi_{L/N,\delta}B^N \ar[d]^{i_{L/N}^{G/N}} \\
  \stab\Pi_{G,\delta}B \ar[r]_-{\Phi^N} & \stab\Pi_{G/N,\delta}B^N.
 }
\]
By Proposition~\ref{prop:commutingInducedMaps}, we have a natural map
$\xi\colon \theta_! i^*\to i^*\theta_!$.

\begin{lemma}
The natural transformation $\xi\colon \theta_! i^*\to i^*\theta_!$
is an isomorphism.
\end{lemma}

The impatient reader may skip the proof, which is simply an elaboration
of the proof of Lemma~\ref{lem:inductFixedIso}.

\begin{proof}
We give the argument for contravariant Mackey functors; the proof for
covariant functors is similar or we can appeal to duality (replacing $\delta$
with $\Lie-\delta$).

By Proposition~\ref{prop:commutingInducedMaps}, the result will follow if we show that
\[
 \xi\colon \int\nolimits^{z\in\stab\Pi_{L/N,\delta}B^N}
  \stab\Pi_{G/N,\delta}B^N(G\times_L z, b) \tensor \stab\Pi_{L,\delta}B(x,\theta z) 
   \to \stab\Pi_{G,\delta}B(G\times_L x, \theta b),
\]given by $\xi(f\tensor g) = \theta f \circ (G_+\smsh_L g)$,
is an isomorphism for all $x\in \stab\Pi_{L,\delta}B$ and $b \in \stab\Pi_{G/N,\delta}B^N$.
Fix $x\colon L/H\to B$ and $b\colon G/K\to B^N$, where $N\leq K\leq G$.

Define a map $\zeta$ inverse to $\xi$ as follows:
Using Theorem~\ref{thm:StableMapsOrbitsOverB}, let 
\[
 \xymatrix{
   [G\times_L x & w \ar[l]_-{p} \ar@{=>}[r]^-{q} & \theta b]
 }
\]
be a generator of $\stab\Pi_{G,\delta}B(G\times_L x, \theta b)$,
with $w\colon G/M\to B$,
so $M\leq H\leq L$ and $p$ is the projection. 
Then $w = G\times_L w'$ with $w\colon L/M\to B$ and
$p = G\times_L p'$ where $p'\colon w'\to x$.
By assumption, $N\leq K$, so $q$ factors
as $w \laxto \theta v \to \theta b$ with $v\colon G/MN\to B^N$; because $MN\leq L$ we can write
$v = G\times_L v'$ with $v'\colon L/MN\to B^N$.
We then let
\[
 \zeta[G\times_L x \from w \laxto \theta b] 
   = [G\times_L v' \to b] \tensor [x \xleftarrow{p'} w' \laxto \theta v'].
\]
Clearly, $\xi\circ\zeta$ is the identity. On the other hand, a typical element in the coend
is a sum of elements of the form
\[
 [G\times_L z \xleftarrow{p} w \laxto b] \tensor g,
\]
where $z\colon L/J\to B$ and $w\colon G/M\to B$
with $N\leq M\leq J\leq L$ and $p$ the projection.
Write $w = G\times_L w'$ so $p = G\times_L p'$ with $p'\colon w'\to z$. 
Finally, factor $w\laxto b$ as $G\times_L w'\laxto G\times_L v'\to b$ as above
and let $v = G\times_L v'$, so
\[
 [G\times_L z \xleftarrow{p} w \laxto b]
  = [G\times_L z \xleftarrow{G\times_L p'} G\times_L w' \laxto G\times_L v' \to b].
\]
We then have
\begin{align*}
 [G\times_L z \xleftarrow{p} w \laxto b] \tensor g
  &\sim [v \to b] \tensor [z \xleftarrow{p'} w' \laxto v']\circ g \\
  &= \sum_i [v \to b] \tensor [x \from u_i \laxto v'] \\
  &= \sum_i [v \to b] \tensor [x \from u_i \laxto t_i \to v'] \\
  &\sim \sum_i [G\times_L t_i \to v \to b] \tensor [x \from u_i \laxto t_i]
\end{align*}
which is in the image of $\zeta$. 
(Here, $u_i\colon L/P_i\to B$ for some $P_i$ and $t_i\colon L/P_iN\to B$.)
So, $\zeta$ is an epimorphism, hence an isomorphism
and the inverse of $\xi$.
\end{proof}

\begin{proposition}\label{prop:paraminductionfixedsets}
Let $N$ be a normal subgroup of $G$,
let $\delta$ be an $N$-closed familial dimension function for $G$,
let $\gamma$ be a virtual representation of $\Pi_{G/N}B^N$,
let $\MackeyOp S$ be a covariant $\stab\Pi_{G/N,\delta}B^N$-module, and
let $\Mackey T$ be a contravariant $\stab\Pi_{G/N,\delta}B^N$-module.
Let $N\leq L\leq G$ with $L\in\F(\delta)$.
If $Y$ is an ex-$G/N$-space over $B^N$
and $y\in\tilde H_\gamma^{G/N,\delta}(Y;\MackeyOp S)$, then
\[
 \epsilon^*(y|L/N) = (\epsilon^* y)|L
  \in \tilde H_\gamma^{L,\delta}(Y; \theta_!(\MackeyOp S|L/N))
  \iso \tilde H_\gamma^{L,\delta}(Y; (\theta_!\MackeyOp S)|L).
\]
Similarly, if $y\in \tilde H^\gamma_{G/N,\delta}(Y;\Mackey T)$, then
\[
 \epsilon^*(y|L/N) = (\epsilon^* y)|L
  \in \tilde H^\gamma_{L,\delta}(Y; \theta_!(\Mackey T|L/N))
  \iso \tilde H^\gamma_{L,\delta}(Y; (\theta_!\Mackey T)|L).
\]
\end{proposition}

\begin{proof}
The proof is formally the same as the proof of Proposition~\ref{prop:inductionfixedsets}.
\end{proof}

We now turn to the relationship between restriction to fixed sets and restriction to subgroups.
The development is very similar to the nonparametrized case.
Referring back to a diagram at the start of this section,
by Proposition~\ref{prop:commutingInducedMaps}, we have a natural map
\[
 \xi\colon \Phi^N_!i^* \to i^*\Phi^N_!.
\]
As in the nonparametrized case, 
we first use Theorem~\ref{thm:StableMapsOrbitsOverB} to get an explicit description of $\Phi^N$.

\begin{proposition}\label{prop:paramStableFixedPoints}
Let $x\colon G/H\to B$ and $y\colon G/K\to B$ be objects of $\stab\Pi_{G,\delta}B$.
Consider the map
\[
 \Phi^N\colon \stab\Pi_{G,\delta}B(x,y) \to \stab\Pi_{G/N,\delta}B^N(x^N,y^N).
\]
If either $N\not\leq H$ or $N\not\leq K$, the target is the trivial group.
If $N\leq H$ and $N\leq K$, then,
on generators, $\Phi^N$ is given by
\begin{align*}
 \Phi^N[x\from z\laxto y]
  &= [x^N \from z^N \laxto y^N] \\
  &= \begin{cases}
       [x\from z\laxto y] & \text{if $N\leq J$} \\
       0 & \text{otherwise,}
     \end{cases}
\end{align*}
where $z\colon G/J\to B$.
Therefore, it is a split epimorphism with kernel generated by those
diagrams $[x\from z\laxto y]$ with $z\colon G/J\to B$ and $N\not\leq J$.
\end{proposition}

\begin{proof}
The proof is essentially the same as the proof of Proposition~\ref{prop:stableFixedPoints}.
\end{proof}

Now we can prove that $\xi$ is an isomorphism.

\begin{lemma}\label{lem:paramfixedAndSubgroupEquality}
If $N$ is a normal subgroup of $G$, $N\leq L\leq G$, and $\Mackey T$ and $\MackeyOp S$ are, respectively, contravariant and covariant $\stab\Pi_{\delta}B$-modules, then
$\xi\colon (\Mackey T|L)^N \to \Mackey T^N|(L/N)$ and
$\xi\colon (\MackeyOp S|L)^N \to \MackeyOp S^N|(L/N)$
are isomorphisms.
\end{lemma}

Again, the proof is an elaboration
of the proof of Lemma~\ref{lem:fixedAndSubgroupEquality}.

\begin{proof}
By Proposition~\ref{prop:commutingInducedMaps}, the result will follow if we show that
\[
 \xi\colon \int\nolimits^{z\in\stab\Pi_{L,\delta} B} \stab\Pi_{G,\delta}B(G\times_L z,b)\tensor 
         \stab\Pi_{L/N,\delta}B^N(x,z^N) 
   \to \stab\Pi_{G/N,\delta}B^N(G\times_L x,b^N),
\]
given by $\xi(f\tensor g) = f^N\circ(G_+\smsh_L g)$,
is an isomorphism for all $x\in\stab\Pi_{L/N,\delta}B^N$ and
$b\in\stab\Pi_{G,\delta}B$.
Fix $x\colon L/H\to B^N$ and $b\colon G/K\to B$.

If $N\not\leq K$, then $b^N = *$ and the target of $\xi$ is 0. For the source,
consider a typical generator 
\[
 f = [G\times_L z \xleftarrow{p} w\xrightarrow{q} b]
\]
of $\stab\Pi_{G,\delta}B(G\times_L z,b)$.
Because $p$ is a projection, we can write $p = G\times_L p'$ where
$p'\colon L/M\to L/J$ is a projection.
Because $M$ is subconjugate to $K$, $N\not\leq M$, hence $(L/M)^N = 0$.  
Therefore, for any $g$,
\[
 f\tensor g = t(p')^*q\tensor g = q\tensor t((p')^N)g = 0
\]
in the coend.

So, assume that $N\leq K$ so $(G/K)^N = G/K$.
Define a map $\zeta$ inverse to $\xi$ as follows: If 
\[
 [G\times_L x\xleftarrow{p} z \xrightarrow{q} b]
\]
is a generator of $\stab\Pi_{G/N,\delta}B^N(G\times_L x,b^N)$ 
with $z\colon G/J\to B^N$ and $N\leq J\leq H$, let
\[
 \zeta[G\times_L x\xleftarrow{p} z \xrightarrow{q} b] = q\tensor t(p')
\]
where $p = G\times_L p'$ and $p'\colon L/J\to L/H$ is the projection.
Clearly, $\xi\circ\zeta$ is the identity.
On the other hand, a typical element in the coend is a sum of elements of the form
\[
 [z\xleftarrow{p}w\laxto b]\tensor g,
\]
where we may assume $z\colon G/J\to B$ with $N\leq J$, because otherwise
the element would live in a 0 group. If $w\colon G/M\to B$, this  is equivalent to
$[w\laxto b]\tensor t(p')g$ where $p'\colon L/M\to L/J$ (and we may assume $N\leq M$,
otherwise this element is 0). 
In turn, such an element can be
written as a sum of elements of the form
$[w\laxto b]\tensor [x\from v\laxto w']$ (where $w = G\times_L w'$), which is equivalent to
\[
 [G\times_L v\laxto w\laxto b]\tensor [x\from v],
\]
which is in the image of $\zeta$.
Thus, $\zeta$ is an isomorphism inverse to $\xi$.

The argument for covariant modules is similar, or we can appeal to duality
(replacing $\delta$ with $\Lie-\delta$).
\end{proof}

The following now follows as it did in the nonparametrized case.

\begin{proposition}\label{prop:fixedandsubgroupParam}
Let $N$ be a normal subgroup of $G$, let $\delta$ be a familial dimension function for $G$,
let $\gamma$ be a virtual representation of $\Pi_G B$,
let $\MackeyOp S$ be a covariant $\stab\Pi_{G,\delta}B$-module and
let $\Mackey T$ be a contravariant $\stab\Pi_{G,\delta}B$-module.
Let $N\leq L\leq G$ with $L\in\F(\delta)$.
If $X$ is an ex-$G$-space over $B$ and $x\in \tilde H_\gamma^{G,\delta}(X;\MackeyOp S)$, then
\[
 (x|L)^N = x^N|L/N \in 
  \tilde H_{\gamma^N-\delta(G/L)}^{L/N,\delta}(X^N; (\MackeyOp S|L)^N) \iso
  \tilde H_{\gamma^N-\delta(G/L)}^{L/N,\delta}(X^N; \MackeyOp S^N|L/N).
\]
Similarly, if $x\in \tilde H^\gamma_{G,\delta}(X;\Mackey T)$, then
\[
 (x|L)^N = x^N|L/N \in 
  \tilde H^{\gamma^N-\delta(G/L)}_{L/N,\delta}(X^N; (\Mackey T|L)^N) \iso
  \tilde H^{\gamma^N-\delta(G/L)}_{L/N,\delta}(X^N; \Mackey T^N|L/N).
\]
\qed
\end{proposition}

\section{Products}

We now turn to various pairings that we saw in the nonparametrized case.

\subsection{Cup products}

As in the nonparametrized case, we start with the external cup product.
The appropriate tensor product of modules is defined as follows.

\begin{definition}\label{def:parammackeyproducts}
Let $G$ and $K$ be two compact Lie groups.
Let $\delta$ be a dimension function
for $G$, let $\epsilon$ be a dimension function for $K$,
and let $\delta\times\epsilon$ denote their product
as in Definition~\ref{def:productDimFcn}.
Let $\zeta$ be a dimension function for $G\times K$ with
$\zeta\dimpred \delta\times\epsilon$,
as in Definition~\ref{def:dimensionpo}.
Let $A$ be a $G$-space and let $B$ be a $K$-space.
\begin{enumerate}
\item
Let $\stab\Pi_{G,\delta}A\tensor\stab\Pi_{K,\epsilon}B$ denote the preadditive category whose objects
are pairs of objects $(a,b)$, as in the product category, with
\[
 (\stab\Pi_{G,\delta}A\tensor\stab\Pi_{K,\epsilon}B)((a,b),(c,d)) 
  = \stab\Pi_{G,\delta}A(a,c)\tensor \stab\Pi_{K,\epsilon}B(b,d).
\]
Let 
\[
 p\colon \stab\Pi_{G,\delta}A\tensor\stab\Pi_{K,\epsilon}B \to h(G\times K)\Spec{}{A\times B}
\]
be the restriction of the external smash product.

Let $i\colon \stab\Pi_{G\times K,\zeta}(A\times B)\to h(G\times K)\Spec{}{A\times B}$ denote the inclusion of the subcategory.

\item
If $\Mackey T$ is a contravariant $\stab\Pi_{G,\delta}A$-module and $\Mackey U$ is a
contravariant $\stab\Pi_{K,\epsilon}B$-module, let $\Mackey T\tensor \Mackey U$ denote the 
$(\stab\Pi_{G,\delta}A\tensor\stab\Pi_{K,\epsilon}B)$-module defined by
 \[
  (\Mackey T\tensor \Mackey U)(a,b) = \Mackey T(a)\tensor \Mackey U(b).
 \]
Define $\MackeyOp S\tensor\MackeyOp V$ similarly for covariant modules.

\item
Let $\Mackey T \exboxprod \Mackey U= \Mackey T \exboxprod_\zeta \Mackey U$ 
(the {\em external box product}) be
the contravariant $\stab\Pi_{G\times K,\zeta}(A\times B)$-module defined by
 \[
 \Mackey T \exboxprod \Mackey U= \Mackey T \exboxprod_\zeta \Mackey U 
   = i^*p_!(\Mackey T\tensor \Mackey U),
 \]
where $i^*$ and $p_!$ are defined in Definition~\ref{def:indres}.
Define the external box product of covariant modules similarly.

\item
If $G=K$ and $A=B$, so $\Mackey T$ and $\Mackey U$ are both $\stab\Pi_{G,\delta}B$-modules,
and $\Delta\in\F(\zeta)$, where $\Delta\leq G\times G$ is the diagonal subgroup, define the
{\em (internal) box product} to be
 \[
  \Mackey T \boxprod \Mackey U = \Mackey T \boxprod_\zeta \Mackey U
   = (\Mackey T\exboxprod \Mackey U) | \Delta(B)
 \]
where $\Delta(B)\subset B\times B$ is the diagonal subspace and we also restrict to
the $\Delta\leq G\times G$, so that
$\Mackey T \boxprod \Mackey U$ is a $\stab\Pi_{G,\delta}B$-module.
Define the box product of covariant modules similarly.
\end{enumerate}
 \end{definition}

\begin{proposition}
Let $\delta$ be a dimension function for $G$, 
let $\epsilon$ be a dimension function for $K$,
and let $\zeta \dimpred \delta\times\epsilon$. 
Let $a\colon G/H\to A$ be an object of $\stab\Pi_{G,\delta}A$ and let
$b\colon K/L\to B$ be an object of $\stab\Pi_{K,\epsilon}B$.
Then
\[
 \Mackey A_{a} \exboxprod_\zeta \Mackey A_{b} 
   \iso h(G\times K)\Spec{}{A\times B}(-,a\times b)
\]
and
\[
 \MackeyOp A^{a} \exboxprod_\zeta \MackeyOp A^{b} 
   \iso h(G\times K)\Spec{}{A\times B}(a\times b,-).
\]
\end{proposition}

\begin{proof}
It's clear from the definitions that
\[
 \Mackey A_{a} \tensor \Mackey A_{b}
  \iso (\stab\Pi_{G,\delta}A\tensor\stab\Pi_{K,\epsilon}B)(-, (a,b)).
\]
Applying $p_!$ and using Proposition~\ref{prop:indresadjunction} gives the result.
The proof for covariant functors is similar.
\end{proof}

Let $\Mackey A_{G/G}$ denote the free $G$-Mackey functor over a point as
in Proposition~\ref{prop:burnsideIdentity} and also denote its pullback to
a $\stab\Pi_{G,\delta}B$-module along the projection $B\to *$.
The following proposition is proved in the same way
as Proposition~\ref{prop:burnsideIdentity}.

\begin{proposition}\label{prop:paramburnsideIdentity}
Let $\delta$ be a complete dimension function for $G$ and let $\Mackey T$ be
a contravariant $\stab\Pi_{G,\delta}B$-module. Then
\[
 \Mackey A_{G/G} \boxprod_{\delta_\Delta} \Mackey T \iso \Mackey T.
\]
\qed
\end{proposition}

The following result identifies the chain complex of a product of CW complexes.
Its proof is the same as that of Proposition~\ref{prop:productChains}.

\begin{proposition}\label{prop:paramproductChains}
Let $\delta$ be a dimension function for $G$ and
let $\epsilon$ be a dimension function for $K$.
Let $A$ be a $G$-space and let $B$ be a $K$-space.
If $X$ is a based $\delta$-$G$-CW($\beta$) complex over $A$ and
$Y$ is a based $\epsilon$-$K$-CW($\gamma$) complex over $B$, then
$X\smsh_{A\times B} Y$ is a based $(\delta\times\epsilon)$-$(G\times K)$-CW($\beta+\gamma$) complex
over $A\times B$ with
\[
 \Mackey C^{G,\delta}_{\beta+*}(X,*) \exboxprod_{\delta\times\epsilon} \Mackey C^{K,\epsilon}_{\gamma+*}(Y,*)
  \iso \Mackey C^{G\times K,\delta\times\epsilon}_{\beta+\gamma + *} (X\smsh_{A\times B} Y,*).
\]
Moreover, this isomorphism respects suspension in each of $X$ and $Y$.
\qed
\end{proposition}

As we noted in the nonparametrized case, 
because of the limitations of $\delta\times\epsilon$, this result by itself is not
as useful as we would like.
To get more useful results, we first note that level-wise smash product induces a map
\[
 \smsh\colon G\PreSpec{\V}{A}\tensor K\PreSpec{\W}{B}
  \to (G\times K)\PreSpec{(\V\dirsum\W)}{A\times B}
\]
where $\V\dirsum\W = \{ V_i\dirsum W_i \}$.
This pairing extends to semistable maps as well.

We also need the following.

\begin{definition}
Suppose that $X$ is a $G$-space over $B$.
We say that $X$ is {\em compactly supported} if there exists a compact subspace
$C$ of $B$ and a $G$-space $X'$ over $C$ such that $X = i_! X'$,
where $i\colon C\to B$ is the inclusion.
Similarly, suppose that $E \in G\PreSpec{\V}{B}$ is a prespectrum over $B$.
We say that $E$ is {\em compactly supported} if there exists a compact subspace
$i\colon C\to B$ and a prespectrum $E'$ over $C$ such that $E = i_! E'$.
\end{definition}

\begin{proposition}\label{prop:pullbackapprox}
Let $\delta$ be a familial dimension function for $G$ and let
$\epsilon$ be a familial dimension function for $K$.
Let $X$ be a based $G$-space over $A$ and let
$Y$ be a based $K$-space over $B$.
Take a $\delta$-$G$-CW($\alpha$) approximation
$\Gamma^\delta\susp_{G,A}^\infty X\to \susp_{G,A}^\infty X$ and
an $\epsilon$-$K$-CW($\beta$) approximation
$\Gamma^\epsilon\susp_{K,B}^\infty Y\to \susp_{K,B}^\infty Y$.
If $j\colon C\to A\times B$ is the inclusion of a compact subspace, then
there exists an integer $N$ (depending on $C$) such that
\[
 j^*(\Gamma^\delta\susp_{G,A}^\infty X\smsh \Gamma^\epsilon\susp_{K,B}^\infty Y)(V_i\dirsum W_i)
  \to j^*(\susp_{G,A}^{V_i} X \smsh \susp_{K,B}^{W_i} Y)
\]
is an $\overline{\F(\delta)\times\F(\epsilon)}$-equivalence for all $i\geq N$
(where we first make the prespectra fibrant before taking $j^*$; equivalently,
we take the homotopy pullback).
\end{proposition}

\begin{proof}
Weak equivalence is determined (homotopy) fiberwise.
Because $C$ is compact, its fundamental groupoid has a skeleton with only
finitely many objects.
Therefore, there are only finitely many fibers to consider, and a sufficiently large suspension
allows us to apply the analogue of Proposition~\ref{prop:productApproximation}.
\end{proof}

\begin{lemma}
Let $E$ be a $\delta$-$G$-CW$(\alpha)$ prespectrum over a $\sigma$-compact base\-space $B$.
Then there is a sequence of compactly supported subcomplexes $\{E_i\}$ such that
$E = \colim_i E_i$ and 
$\Mackey C_{\alpha+*}^{G,\delta}(E) \iso \colim_i \Mackey C_{\alpha+*}^{G,\delta}(E_i)$.
\end{lemma}

\begin{proof}
Write $B = \union_i B_i$ where $B_1\subset B_2\subset\cdots$ and each $B_i$ is compact.
Construct $E'_k$ over $B_k$ as follows: Take as its vertices all 0-dimensional cells
in $E$ that lie over $B_k$. Inductively, take as its $n$-cells all $n$-cells of
$E$ that lie over $B_k$ and attach to lower-dimensional cells already chosen.
Let $E_k = i_! E'_k$ where $i\colon B_1\to B$ is the inclusion.

Because $B_{k-1}\subset B_k$, we get that $E_{k-1}\subset E_k$.
Because every cell in $E$ is contained in a finite subcomplex, it will appear in some $E_k$.
Therefore, $E$ is the sequential colimit of the $E_k$ and the statement
about chain complexes is clear.
\end{proof}

\begin{proposition}
Let $\delta$ be a familial dimension function for $G$, let
$\epsilon$ be a familial dimension function for $K$,
and let $\zeta$ be a familial dimension function for $G\times K$ with $\zeta\dimpred\delta\times\epsilon$.
Let $A$ be a $\sigma$-compact $G$-space and let $B$ be a $\sigma$-compact $K$-space.
Let $X$ be a based $G$-space over $A$ and let $Y$ be a based $K$-space over $B$.
Let $\Gamma^\delta\susp_{G,A}^\infty X\to \susp_{G,A}^\infty X$ and
$\Gamma^\epsilon\susp_{K,B}^\infty Y\to \susp_{K,B}^\infty Y$
be, respectively, $\delta$-$G$-CW($\alpha$) and
$\epsilon$-$K$-CW($\beta$) approximations
and let $\Gamma^\zeta\susp_{G\times K,A\times B}^\infty(X\smsh Y) 
\to \susp_{G\times K,A\times B}^\infty(X\smsh Y)$
be a $\zeta$-$(G\times K)$-CW($\alpha+\beta$) approximation.
Write
\[
 \Gamma^\zeta\susp_{G\times K,A\times B}^\infty(X\smsh Y)
  = \colim_i \Gamma_i
\]
where each $\Gamma_i$ is compactly supported, by the preceding lemma.
Then there exist compatible semistable cellular lax maps
\[
 \mu_i\colon \Gamma_i \laxto
  \Gamma^\delta\susp_{G,A}^\infty X \smsh \Gamma^\epsilon\susp_{K,B}^\infty Y
\]
over $\susp_{G\times K,A\times B}^\infty(X\smsh Y)$, 
the collection of maps being unique up to semistable cellular homotopy.
\end{proposition}

\begin{proof}
Write $A\times B = \union C_i$, where each $C_i$ is compact, and
$\Gamma_i = j_! \Gamma'_i$ where $j\colon C_i\to A\times B$ is the inclusion.
Consider the following diagram, in which $j^*$ indicates homotopy pullback.
\[
 \xymatrix{
  \Gamma'_{i-1} \ar[r] \ar[d] 
    & j^*(\Gamma^\delta\susp_{G,A}^\infty X \smsh \Gamma^\epsilon\susp_{K,B}^\infty Y) \ar[d] \\
  \Gamma'_i \ar[r] \ar@{-->}[ur] &
    j^*(\susp_{G,A}^\infty X \smsh \susp_{K,B}^\infty Y)
 }
\]
A semistable lift exists by the relative Whitehead theorem, using
Proposition~\ref{prop:pullbackapprox}.
The adjoint is the map
\[
 \mu_i\colon \Gamma_i = j_!\Gamma_i \laxto \Gamma^\delta\susp_{G,A}^\infty X \smsh \Gamma^\epsilon\susp_{K,B}^\infty Y.
\]
Uniqueness up to semistable cellular homotopy follows from cellular approximation.
\end{proof}

\begin{definition}\label{def:paramchainProduct}
Let $\delta$ be a familial dimension function for $G$, let
$\epsilon$ be a familial dimension function for $K$,
and let $\zeta$ be a familial dimension function for $G\times K$ with $\zeta\dimpred\delta\times\epsilon$.
If $X$ is an ex-$G$-space over $A$ and $Y$ is an ex-$K$-space over $B$,
where $A$ and $B$ are $\sigma$-compact, let
\[
 \mu_*\colon \Mackey C^{G\times K,\zeta}_{\alpha+\beta+*}(X\smsh Y)
  \to \Mackey C^{G,\delta}_{\alpha+*}(X)\exboxprod_{\zeta} \Mackey C^{K,\epsilon}_{\beta+*}(Y)
\]
be the chain map induced by the family $\{\mu_i\}$ from the preceding corollary
(using Proposition~\ref{prop:paramproductChains} to identify the chain complex on the right).
It is well-defined up to chain homotopy.
\end{definition}

We can now define pairings in cohomology.
Let $\delta$, $\epsilon$, and $\zeta$ be familial with $\zeta\dimpred \delta\times\epsilon$, as above,
let $X$ be an ex-$G$-space over $A$, let $Y$ be an ex-$K$-space over $B$,
with $A$ and $B$ $\sigma$-compact,
let $\Mackey T$ be a $\stab\Pi_{G,\delta}A$-module, and
let $\Mackey U$ be a $\stab\Pi_{K,\epsilon}B$-module.
The external box product $\exboxprod = \exboxprod_\zeta$ and the
map $\mu_*$ induce a natural chain map
\begin{align*}
 \Hom_{\stab\Pi_{G,\delta}A}(\Mackey C^{G,\delta}_{\alpha+*}(X), \Mackey T) \tensor {}&
       \Hom_{\stab\Pi_{K,\epsilon}B}(\Mackey C^{K,\epsilon}_{\beta+*}(Y), \Mackey U) \\
 &\to
   \Hom_{\stab\Pi_{G\times K,\zeta}(A\times B)}(
    \Mackey C^{G,\delta}_{\alpha+*}(X)\exboxprod \Mackey C^{K,\epsilon}_{\beta+*}(Y), 
    \Mackey T \exboxprod \Mackey U) \\
 &\to \Hom_{\stab\Pi_{G\times K,\zeta}(A\times B)}(
    \Mackey C^{G\times K,\zeta}_{\alpha+\beta+*}(X\smsh Y), 
    \Mackey T \exboxprod \Mackey U).
\end{align*}
This induces the (external) cup product
 \[
 -\cup - \colon 
 \tilde H^{\alpha}_{G,\delta}(X;\Mackey T) \tensor
 \tilde H^{\beta}_{K,\epsilon}(Y;\Mackey U) 
  \to \tilde H^{\alpha+\beta}_{G\times K,\zeta}(X\smsh_{A\times B} Y;\Mackey T\exboxprod \Mackey U).
 \]
When $G = K$, $A=B$, and $\F(\zeta)$ contains the diagonal subgroup $\Delta\leq G\times G$, we can follow the external cup product with the restriction to $\Delta(B)\subset B\times B$ and
$\Delta\leq G\times G$. This gives the internal cup product
\[
 -\cup - \colon 
 \tilde H^{\alpha}_{G,\delta}(X;\Mackey T) \tensor
 \tilde H^{\beta}_{G,\epsilon}(Y;\Mackey U) 
 \to \tilde H^{\alpha+\beta-\zeta(G\times G/\Delta)}_{G,\zeta|\Delta}(X\smsh_B Y;\Mackey T\boxprod \Mackey U).
\]
Of course, when $X = Y$ we can apply restriction along the diagonal $X\to X\smsh_B X$ to
get
\[
 -\cup - \colon 
 \tilde H^{\alpha}_{G,\delta}(X;\Mackey T) \tensor
 \tilde H^{\beta}_{G,\epsilon}(X;\Mackey U) 
 \to \tilde H^{\alpha+\beta-\zeta(G\times G/\Delta)}_{G,\zeta|\Delta}(X;\Mackey T\boxprod \Mackey U).
\]
A useful special case is $\delta_\Delta \dimpred 0\times \delta$ with $\delta$ complete, which gives us pairings
\[
 -\cup - \colon 
 \tilde H^{\alpha}_{G}(X;\Mackey T) \tensor
 \tilde H^{\beta}_{G,\delta}(Y;\Mackey U) 
 \to \tilde H^{\alpha+\beta}_{G,\delta}(X\smsh_B Y;\Mackey T\boxprod \Mackey U)
\]
and
\[
 -\cup - \colon 
 \tilde H^{\alpha}_{G}(X;\Mackey T) \tensor
 \tilde H^{\beta}_{G,\delta}(X;\Mackey U) 
 \to \tilde H^{\alpha+\beta}_{G,\delta}(X;\Mackey T\boxprod \Mackey U).
\]

Say that a contravariant $\stab\Pi_G B$-module $\Mackey T$ is a {\em ring} 
if there is an associative multiplication 
$\Mackey T\boxprod \Mackey T\to \Mackey T$. 
This then makes $\tilde H_G^*(X;\Mackey T)$ a ring.
For example, recall that we let $\Mackey A_{G/G} = \Mackey A_{G/G,0}$, pulled back
to $B$ along the projection $B\to *$.
$\Mackey A_{G/G}$ is a ring
by Proposition~\ref{prop:paramburnsideIdentity}.
By that same proposition,
$\Mackey A_{G/G} \boxprod_{\delta_\Delta} \Mackey T \iso \Mackey T$ for any contravariant
$\delta$-$G$-Mackey functor $\Mackey T$. 
This makes $\Mackey T$ a {\em module} over
$\Mackey A_{G/G}$, in the sense that there is an associative pairing
$\Mackey A_{G/G}\boxprod_{\delta_\Delta} \Mackey T\to \Mackey T$ (namely, the isomorphism).
This makes every cellular cohomology theory a module over ordinary cohomology
with coefficients in $\Mackey A_{G/G}$, using the cup product
\[
 -\cup - \colon 
 \tilde H^{\alpha}_{G}(X;\Mackey A_{G/G}) \tensor
 \tilde H^{\beta}_{G,\delta}(X;\Mackey T) 
 \to \tilde H^{\alpha+\beta}_{G,\delta}(X;\Mackey T).
\]

The following theorem records the main properties of the external cup product, 
from which similar properties
of the other products follow by naturality.

\begin{theorem}\label{thm:paramcupproduct}
Let $\delta$ be a familial dimension function for $G$, 
let $\epsilon$ be a familial dimension function for $K$,
and let $\zeta$ be a familial dimension function for $G\times K$
with $\zeta\dimpred \delta\times\epsilon$.
The external cup product
 \[
 -\cup - \colon 
 \tilde H^{\alpha}_{G,\delta}(X;\Mackey T) \tensor
 \tilde H^{\beta}_{K,\epsilon}(Y;\Mackey U) 
  \to \tilde H^{\alpha+\beta}_{G\times K,\zeta}(X\smsh_{A\times B} Y;\Mackey T\exboxprod \Mackey U)
 \]
satisfies the following.
 \begin{enumerate}
 \item It is natural: $f^*(x)\cup g^*(y) = (f\smsh g)^*(x\cup y)$.
 \item It respects suspension: For any representation $V$ of $G$, the following diagram commutes:
\[
 \xymatrix{
 \tilde H^{\alpha}_{G,\delta}(X;\Mackey T) \tensor \tilde H^{\beta}_{K,\epsilon}(Y;\Mackey U)
  \ar[r]^\cup \ar[d]_{\sigma^V\tensor\id}
  & \tilde H^{\alpha+\beta}_{G\times K,\zeta}(X\smsh_{A\times B} Y;\Mackey T\exboxprod \Mackey U)
     \ar[dd]^{\sigma^V} \\
 \tilde H^{\alpha+V}_{G,\delta}(\susp^V X;\Mackey T) \tensor \tilde H^{\beta}_{K,\epsilon}(Y;\Mackey U)
  \ar[d]_\cup \\
 \tilde H^{\alpha+V+\beta}_{G\times K,\zeta}(\susp^V X\smsh_{A\times B} Y;\Mackey T\exboxprod \Mackey U)
     \ar[r]_-\iso
  & \tilde H^{\alpha+\beta+V}_{G\times K,\zeta}(\susp^V(X\smsh_{A\times B} Y);\Mackey T\exboxprod \Mackey U)
 }
\]
The horizontal isomorphism at the bottom of the diagram comes from the identification
$\alpha+V+\beta \iso \alpha+\beta+V$.
The similar diagram for suspension of $Y$ also commutes.
 \item It is associative: $(x\cup y)\cup z = x\cup (y\cup z)$ when we identify
gradings using the obvious identification $(\alpha+\beta)+\gamma \iso \alpha + (\beta+\gamma)$.
 \item\label{item:commutativityParam} It is commutative: If $x\in \tilde H^\alpha_{G,\delta}(X;\Mackey T)$ and
$y\in \tilde H^\beta_{K,\epsilon}(Y;\Mackey U)$ then $x\cup y = \iota(y\cup x)$ where $\iota$ is the
evident isomorphism
\[
 \iota\colon \tilde H^{\alpha+\beta}_{G\times K,\zeta}(X\smsh_{A\times B} Y;\Mackey T\exboxprod \Mackey U)
  \iso \tilde H^{\beta+\alpha}_{K\times G,\tilde\zeta}(Y\smsh_{B\times A} X;\Mackey U\exboxprod \Mackey T);
\]
$\tilde\zeta$ is the dimension function on $K\times G$ induced by $\zeta$
and $\iota$ uses the isomorphism of $\alpha+\beta$ and $\beta+\alpha$ that
switches the direct summands.
 \item It is unital: The map
\[
 \tilde H^0_G(S^0;\Mackey A_{G/G}) \tensor \tilde H^\alpha_{G,\delta}(X;\Mackey T)
   \to \tilde H^\alpha_{G,\delta}(X;\Mackey T)
\]
takes $1\tensor x \mapsto x$, where $1 \in \tilde H^0_G(S^0;\Mackey A_{G/G})\iso A(G)$
is the unit.
 \item\label{item:paramcupsubgroup}
 It respects restriction to subgroups: 
$(x|J)\cup (y|L) = (x\cup y)|(J\times L)$. (But see Remark~\ref{rem:paramprodCoeffs}.)
 \item\label{item:paramcupfixedsets}
 It respects restriction to fixed sets:
$x^J \cup y^L = (x\cup y)^{J\times L}$. (But, again, see Remark~\ref{rem:paramprodCoeffs}.)
 \end{enumerate}
 \end{theorem}

The proofs are all standard except for the last two points,
but the proofs of these are just the obvious generalizations of the proofs
of the last two parts of Theorem~\ref{thm:vcupproduct}.

\begin{remark}\label{rem:paramprodCoeffs}
Part (\ref{item:paramcupsubgroup}) of Theorem~\ref{thm:paramcupproduct} is stated
a bit loosely. In fact,
 \[
 (x|J)\cup(y|L) \in H_{J\times L}^*(X\smsh Y; 
			(\Mackey T|J)\exboxprod (\Mackey U|L))
 \]
while
 \[
 (x\cup y)|(J\times L) \in H_{J\times L}^*(X\smsh Y;
  (\Mackey T\exboxprod \Mackey U)|(J\times L)).
 \]
The theorem should say that, when we apply the map induced by the homomorphism
$(\Mackey T|J)\exboxprod (\Mackey U|L) \to 
(\Mackey T\exboxprod \Mackey U)|(J\times L)$,
the element $(x|J)\cup(y|L)$ maps to $(x\cup y)|(J\times L)$.

A similar comment applies to part (\ref{item:paramcupfixedsets}) in general,
when we have to first restrict to normalizers.
However, if $J$ is normal in $G$ and $L$ is normal in $K$, we can use the isomorphism
$\Mackey T^J \exboxprod \Mackey U^L \iso (\Mackey T \exboxprod \Mackey U)^{J\times L}$
to identify the two cohomology groups.
\end{remark}

Now we look at how the cup product is represented.
Let $A$ be a $G$-space and let $B$ be a $K$-space.
Let $\alpha$ be a virtual representation of $\Pi_G A$ and let $\beta$ be
a virtual representation of $\Pi_K B$.
Given the $G$-spectrum $H_\delta\Mackey T^\alpha$ over $A$
and the $K$-spectrum $H_\epsilon\Mackey U^\beta$ over $B$,
we can form the $(G\times K)$-spectrum
$H_\delta\Mackey T \smsh H_\epsilon\Mackey U$ over $A\times B$.
The external cup product should then be represented by a $(G\times K)$-map
\[
 H_\delta\Mackey T^\alpha \smsh H_\epsilon\Mackey U^\beta
  \to H_\zeta(\Mackey T\exboxprod \Mackey U)^{\alpha+\beta}
\]
that is an isomorphism in $\Mackey\pi_{\alpha+\beta}^{G\times K,\zeta}$. 
An explicit construction
is based on the following calculation.

\begin{proposition}\label{prop:paramEMsmshhomotopy}
Let $\delta$ be a familial dimension function for $G$, let $\epsilon$ be
a familial dimension function for $K$, and let $\zeta\dimpred \delta\times\epsilon$.
Let $\Mackey T$ be a contravariant $\stab\Pi_{G,\delta}A$-module and 
let $\Mackey U$ be a contravariant $\stab\Pi_{K,\epsilon}B$-module.
Then we have
\[
 \Mackey\pi_{\alpha+\beta+n}^{G\times K,\zeta}(H_\delta\Mackey T^\alpha \smsh H_\epsilon\Mackey U^\beta)
  = \begin{cases}
      0 & \text{if $n < 0$} \\
      \Mackey T \exboxprod_\zeta \Mackey U & \text{if $n = 0$.}
    \end{cases}
\]
\end{proposition}

\begin{proof}
Let $c\colon (G\times K)/J \to A\times B$ be a $(G\times K)$-map with
$J\in \F(\zeta)$.
Let $H = p_1(J)\leq G$ and $L = p_2(J)\leq K$;
because $\zeta\dimpred \delta\times\epsilon$, we must have
$H\in\F(\delta)$ and $L\in\F(\epsilon)$.
If $c(eJ) = (x,y)$, then we have maps $a\colon G/H\to A$ with $a(eH) = x$ and
$b\colon K/L\to B$ with $b(eL) = y$, and $c$ factors as $(a\times b)\circ p$
where $p\colon (G\times K)/J\to G/H\times K/L$ is the projection.
Further, $a\times b$ is initial with the properties that it is a product and $c$ factors
through it.
We need to compute
\begin{align*}
 \Mackey\pi_{\alpha+\beta+n}^{G\times K,\zeta}
   (H_\delta\Mackey T^\alpha \smsh H_\epsilon\Mackey U^\beta)(c)
  &\iso \pi_{\alpha_0(a)+\beta_0(b)+n}^{J,\zeta}
       (c^* (H_\delta\Mackey T^\alpha \smsh H_\epsilon\Mackey U^\beta)) \\
  &\iso \pi_{\alpha_0(a)+\beta_0(b)+n}^{H\times L,\zeta}
       ((a\times b)^* (H_\delta\Mackey T^\alpha \smsh H_\epsilon\Mackey U^\beta)),
\end{align*}
where we implicitly consider a $(G\times K)$-spectrum over $(G\times K)/J$ to
be a $J$-spectrum, and similarly for spectra over $G/H\times K/L$.
Using the calculation of the fibers given at the beginning of
Section~\ref{sec:paramrepspectra} we get
\[
 (a\times b)^* (H_\delta\Mackey T^\alpha \smsh H_\epsilon\Mackey U^\beta)
  \hmtpc \susp^{\alpha_0(a)}H_\delta(a^*\Mackey T)\smsh 
          \susp^{\beta_0(b)}H_\epsilon(b^*\Mackey U).
\]
The calculation then follows from Proposition~\ref{prop:EMsmshhomotopy}.
\end{proof}

As usual, for example, by killing higher homotopy groups, it follows that there
is a map of $(G\times K)$-spectra over $A\times B$
\[
 H_\delta\Mackey T^\alpha \smsh H_\epsilon\Mackey U^\beta
  \to H_\zeta(\Mackey T\exboxprod_\zeta \Mackey U)^{\alpha+\beta}
\]
that is an isomorphism in $\Mackey\pi_{\alpha+\beta}^{G\times K,\zeta}$.
That this represents the cup product in cohomology that we constructed on the chain level
follows by considering its effect on products of cells,
which reduces to the statement that the maps of fibers represent
nonparametrized cup products.

When $G = K$ and $A=B$ we can restrict to diagonals to
get a map of $G$-spectra
\[
 H_\delta\Mackey T^\alpha \smsh H_\epsilon\Mackey U^\beta
   \to \susp^{-\zeta(G\times G/\Delta)}H_\zeta(\Mackey T\boxprod_\zeta \Mackey U)^{\alpha+\beta}
\]
that is an isomorphism in $\Mackey\pi_{\alpha+\beta}^{G,\zeta}$.

Finally, let's point out the specializations to several interesting choices of $\delta$ and $\epsilon$,
using the internal cup products.

\begin{remark}
The general cup product gives us the following special cases.
\begin{enumerate}
\item
Taking $\delta = \epsilon = 0$ on $G$ and $\zeta=0$ on $G\times G$, we have the cup product
\[
 \tilde H_G^\alpha(X;\Mackey T)\tensor \tilde H_G^\beta(Y;\Mackey U)
   \to \tilde H_G^{\alpha+\beta}(X\smsh Y; \Mackey T\boxprod \Mackey U).
\]
This product is represented by a $G$-map
\[
 H\Mackey T^\alpha \smsh H\Mackey U^\beta \to H(\Mackey T \boxprod \Mackey U)^{\alpha+\beta}.
\]

\item
Taking $\delta=\epsilon=\Lie$ on $G$ and $\zeta=\Lie$ on $G\times G$, we have the cup product
\[
 \tilde\H_G^\alpha(X;\MackeyOp R)\tensor \tilde\H_G^\beta(Y;\MackeyOp S)
   \to \tilde\H_G^{\alpha+\beta-\Lie(G)}(X\smsh Y; \MackeyOp R\boxprod_\Lie \MackeyOp S).
\]
This product is represented by a $G$-map
\[
 H_\Lie\MackeyOp R^\alpha \smsh H_\Lie\MackeyOp S^\beta 
 \to \susp^{-\Lie(G)}H_\Lie(\MackeyOp R \boxprod_\Lie \MackeyOp S)^{\alpha+\beta}.
\]

\item
Taking $\delta=0$, $\epsilon=\Lie$, and $\zeta= \Lie_\Delta$, we have the cup product
\[
 \tilde H_G^\alpha(X;\Mackey T)\tensor \tilde\H_G^\beta(Y;\MackeyOp S)
   \to \tilde\H_G^{\alpha+\beta}(X\smsh Y; \Mackey T \boxprod_{\Lie_\Delta} \MackeyOp S).
\]
This product is represented by a $G$-map
\[
 H\Mackey T^\alpha \smsh H_\Lie\MackeyOp S^\beta 
 \to H_\Lie(\Mackey T \boxprod_{\Lie_\Delta} \MackeyOp S)^{\alpha+\beta}.
\]
Note that $\Lie_\Delta$ restricts to $\Lie$ on the diagonal, so $\Mackey T \boxprod_{\Lie_\Delta} \MackeyOp S$ is a contravariant
$\stab\Pi_{G,\Lie}B$-module, 
which we may also consider as a covariant $\stab\Pi_{G,0}B$-module.

\item
Taking $\delta = 0$, $\epsilon=\Lie$, and $\zeta= \Lie - \Lie_\Delta$, we have the cup product
\[
 \tilde H_G^\alpha(X;\Mackey T)\tensor \tilde\H_G^\beta(Y;\MackeyOp S)
   \to \tilde H_G^{\alpha+\beta-\Lie(G)}(X\smsh Y; \Mackey T \boxprod_{\Lie - \Lie_\Delta} \MackeyOp S).
\]
This product is represented by a $G$-map
\[
 H\Mackey T^\alpha \smsh H_\Lie\MackeyOp S^\beta
  \to \susp^{-\Lie(G)}H(\Mackey T \boxprod_{\Lie - \Lie_\Delta} \MackeyOp S)^{\alpha+\beta}.
\]
Note that $\Mackey T \boxprod_{\Lie - \Lie_\Delta} \MackeyOp S$ is a
contravariant $\stab\Pi_{G,0}B$-module because $\Lie - \Lie_\Delta$ restricts to 0 on the diagonal.
\end{enumerate}
\end{remark}

\subsection{Slant products, evaluations, and cap products}

We now consider evaluation maps and cap products.
See the comments at the beginning of Section~\ref{sec:capproducts}.

\begin{definition}
Let $\delta$ be a dimension function for $G$, let $\epsilon$ be a
dimension function for $K$, and let $\zeta$ be a dimension function for $G\times K$
with $\zeta \dimpred \delta\times\epsilon$.
Let $A$ be a $G$-space and let $B$ be a $K$-space.
Suppose that $\Mackey T$ is a contravariant $\stab\Pi_{K,\epsilon}B$-module and
$\MackeyOp U$ is a covariant $\stab\Pi_{G\times K,\zeta}(A\times B)$-module.
Then we define a covariant $\stab\Pi_{G,\delta}A$-module $\Mackey T\mixprod \MackeyOp U$ by
\[
 (\Mackey T\mixprod \MackeyOp U)(a) 
  = \Mackey T \tensor_{\stab\Pi_{K,\epsilon}B} (p^*i_!\MackeyOp U)(a\tensor -)
\]
where 
\[
 p\colon \stab\Pi_{G,\delta}A\tensor\stab\Pi_{K,\epsilon}B\to h(G\times K)\Spec{}{A\times B}
\]
and 
\[
 i\colon\stab\Pi_{G\times K,\zeta}(A\times B) \to h(G\times K)\Spec{}{A\times B}
\]
are the functors given in Definition~\ref{def:parammackeyproducts}.
\end{definition}

\begin{example}\label{prop:paramEMmixedproduct}
As an example and as a calculation we'll need later, we show that
\[
 \Mackey A_{b}\mixprod \MackeyOp A^{c}
  \iso h(G\times K)\Spec{}{A\times B}(c, -\smsh b)
\]
when $b\in \stab\Pi_{K,\epsilon}B$ and $c\in \stab\Pi_{G\times K,\zeta}(A\times B)$.
For, if $a$ is an object in $\stab\Pi_{G,\delta}A$, we have
\begin{align*}
 (\Mackey A_{b}\mixprod \MackeyOp A^{c})(a)
  &= \Mackey A_{b}\tensor_{\stab\Pi_{K,\epsilon}B} 
           (p^*i_!\MackeyOp A^{c})(a\tensor -) \\
  &= \Mackey A_{b}\tensor_{\stab\Pi_{K,\epsilon}B} 
           h(G\times K)\Spec{}{A\times B}(c, a \smsh -) \\
  &\iso h(G\times K)\Spec{}{A\times B}(c, a \smsh b).
\end{align*}
\end{example}

For $X$ an ex-$G$-space over $A$ and $Y$ an ex-$K$-space over $B$, 
the external slant product will be a map
\[
 -\slant - \colon 
  \tilde H_{K,\epsilon}^\beta(Y;\Mackey T) \tensor
    \tilde H^{G\times K,\zeta}_{\alpha+\beta}(X\smsh Y; \MackeyOp U)
  \to \tilde H^{G,\delta}_\alpha(X; \Mackey T \mixprod \MackeyOp U).
\]
On the chain level, we take the following map:
\begin{align*}
 \Hom&{}_{\stab\Pi_{K,\epsilon}B}(\Mackey C^{K,\epsilon}_{\beta+*}(Y),\Mackey T) \tensor
  (\Mackey C^{G\times K,\zeta}_{\alpha+\beta+*}(X\smsh Y)
     \tensor_{\stab\Pi_{G\times K,\zeta}(A\times B)} \MackeyOp U) \\
  &\to \Hom_{\stab\Pi_{K,\epsilon}B}(\Mackey C^{K,\epsilon}_{\beta+*}(Y),\Mackey T) \tensor
    ((\Mackey C^{G,\delta}_{\alpha+*}(X)\exboxprod_\zeta \Mackey C^{K,\epsilon}_{\beta+*}(Y))
    \tensor_{\stab\Pi_{G\times K,\zeta}(A\times B)} \MackeyOp U) \\
  &\xrightarrow{\nu} (\Mackey C^{G,\delta}_{\alpha+*}(X)\exboxprod_\zeta\Mackey T )
    \tensor_{\stab\Pi_{G\times K,\zeta}(A\times B)} \MackeyOp U \\
  &\iso (\Mackey C^{G,\delta}_{\alpha+*}(X)\tensor\Mackey T)
    \tensor_{\stab\Pi_{G,\delta}A\tensor\stab\Pi_{K,\epsilon}B} p^*i_!\MackeyOp U \\
  &\iso \Mackey C^{G,\delta}_{\alpha+*}(X) \tensor_{\stab\Pi_{G,\delta}A}
    (\Mackey T \tensor_{\stab\Pi_{K,\epsilon}B} p^*i_!\MackeyOp U) \\
  &= \Mackey C^{G,\delta}_{\alpha+*}(X) \tensor_{\stab\Pi_{G,\delta}A}
    (\Mackey T\mixprod \MackeyOp U).
\end{align*}
The first arrow is induced by the map $\mu_*$ from Definition~\ref{def:paramchainProduct}.
The map $\nu$ is evaluation, with the sign
\[
 \nu(y\tensor a\tensor b\tensor u) = (-1)^{pq}a\tensor y(b)\tensor u
\]
if $y\in \Hom_{\stab\Pi_{K,\epsilon}B}(\Mackey C^{K,\epsilon}_{\beta+p}(Y),\Mackey T)$ and
$a\in \Mackey C^{G,\delta}_{\alpha+q}(X)$.
Our various sign conventions imply that the composite above is a chain map.
Taking homology defines our slant product.
The following properties follow easily from the definition.

\begin{theorem}\label{thm:paramslantproduct}
Let $\delta$ be a familial dimension function for $G$, let $\epsilon$ be a familial
dimension function for $K$, and let $\zeta$ be a familial dimension function for $G\times K$
with $\zeta \dimpred \delta\times\epsilon$.
Let $\alpha$ be a virtual representation of $\Pi_G A$ and let $\beta$ be a virtual representation
of $\Pi_K B$.
The slant product
\[
 -\slant - \colon 
  \tilde H_{K,\epsilon}^\beta(Y;\Mackey T) \tensor
    \tilde H^{G\times K,\zeta}_{\alpha+\beta}(X\smsh Y; \MackeyOp U)
  \to \tilde H^{G,\delta}_\alpha(X; \Mackey T \mixprod \MackeyOp U).
\]
has the following properties.
\begin{enumerate}
\item
 It is natural in the following sense: Given a $G$-map
$f\colon X\laxto X'$ over $A$, a $K$-map $g\colon Y\laxto Y'$ over $B$, and elements
$y'\in \tilde H_{K,\epsilon}^\beta(Y';\Mackey T)$ and
$z\in \tilde H^{G\times K,\zeta}_{\alpha+\beta}(X\smsh Y; \MackeyOp U)$,
we have
\[
 y' \slant (f \smsh g)_*(z) = f_*(g^*(y') \slant z).
\]
Put another way, the slant product is a natural transformation in its adjoint form
\[
 \tilde H^{G\times K,\zeta}_{\alpha+\beta}(X\smsh Y; \MackeyOp U)
  \to \Hom(\tilde H_{K,\epsilon}^\beta(Y;\Mackey T), 
  		\tilde H^{G,\delta}_\alpha(X; \Mackey T \mixprod \MackeyOp U)).
\]

\item
It respects suspension in the sense that
\[
 (\sigma^W y) \slant (\sigma^{V+W} z) = \sigma^V(y \slant z).
\]

\item
It is associative in the following sense. Given
$y\in \tilde H_{K,\epsilon}^\beta(Y;\Mackey R)$,
$z\in \tilde H_{L,\zeta}^\gamma(Z;\Mackey S)$, and
$w\in \tilde H^{G\times K\times L,\eta}_{\alpha+\beta+\gamma}(X\smsh Y\smsh Z;\MackeyOp U)$,
where $\eta\dimpred \delta\times\epsilon\times\zeta$,
we have
\[
 (y\cup z)\slant w = y\slant(z\slant w).
\]

\qed
 \end{enumerate}
 \end{theorem}

We're most interested in the internalization of the slant product to the diagonal
$\delta\leq G\times G$.
So, we also call the following a slant product.

\begin{definition}\label{def:paraminternalslant}
Let $\delta$ and $\epsilon$ be familial dimension functions for $G$ and let
$\zeta\dimpred \delta\times\epsilon$ be a familial dimension function for $G\times G$;
assume that $\Delta\in\F(\zeta)$ and
write $\zeta$ again for $\zeta|\Delta$.
Let $A$ and $B$ be $G$-spaces, 
let $X$ be an ex-$G$-space over $A$ and let $Y$ be an ex-$G$-space over $B$.
If $\Mackey T$ is a contravariant $\stab\Pi_{G,\epsilon}B$-module and
$\MackeyOp U$ is a covariant $\stab\Pi_{G,\zeta}(A\times B)$-module, write
\[
 \Mackey T\mixprod \MackeyOp U = \Mackey T \mixprod [(G\times G)\times_\Delta \MackeyOp U],
\]
a covariant $\stab\Pi_{G,\delta}A$-module defined using the external version of $\mixprod$ on the right.
The {\em internal slant product}
\[
 -\slant- \colon 
   \tilde H^\beta_{G,\epsilon}(Y;\Mackey T) \tensor
   \tilde H^{G,\zeta}_{\alpha+\beta-\zeta(G\times G/\Delta)}(X\smsh Y;\MackeyOp U)
  \to \tilde H^{G,\delta}_\alpha(X;\Mackey T\mixprod \MackeyOp U)
\]
is then the composite
\begin{align*}
 \tilde H^\beta_{G,\epsilon}(Y;\Mackey T) &\tensor
		\tilde H^{G,\zeta}_{\alpha+\beta-\zeta(G\times G/\Delta)}(X\smsh Y;\MackeyOp U) \\
  &\to \tilde H^\beta_{G,\epsilon}(Y;\Mackey T) \tensor
       \tilde H^{G,\zeta}_{\alpha+\beta-\zeta(G\times G/\Delta)}(X\smsh Y;((G\times G)\times_\Delta\MackeyOp U)|\Delta) \\
  &\iso \tilde H^\beta_{G,\epsilon}(Y;\Mackey T) \tensor
       \tilde H^{G\times G,\zeta}_{\alpha+\beta}((G\times G)_+\smsh_\Delta(X\smsh Y); (G\times G)\times_\Delta\MackeyOp U) \\
  &\to \tilde H^\beta_{G,\epsilon}(Y;\Mackey T) \tensor
       \tilde H^{G\times G,\zeta}_{\alpha+\beta}(X\smsh Y; (G\times G)\times_\Delta\MackeyOp U) \\
  &\to \tilde H^{G,\delta}_\alpha(X;\Mackey T\mixprod \MackeyOp U).
\end{align*}
The first map is the unit $\MackeyOp U\to ((G\times G)\times_\Delta\MackeyOp U)|\Delta$,
the second is the Wirthm\"uller isomorphism, the third is induced by the $(G\times G)$-map
$(G\times G)_+\smsh_\Delta(X\smsh Y) \to X\smsh Y$
over $A\times B$, and the last map is the external slant product.
\end{definition}

We can now use the internal slant product to define evaluation and the cap product.

\begin{definition}
Let $\delta$ be a familial dimension function for $G$,
let $\alpha$ be a virtual representation of $G$, and let
$\beta$ be a virtual representation of $\Pi_G B$.
Let $\Mackey T$ be a contravariant $\stab\Pi_{G,\delta}B$-module and
let $\Mackey U$ be a covariant $\stab\Pi_{G,\delta}B$-module.
The {\em evaluation map}
\[
 \langle -,- \rangle \colon
  \tilde H_{G,\delta}^\beta(X;\Mackey T) 
     \tensor \tilde H^{G,\delta}_{\alpha+\beta}(X;\MackeyOp U)
  \to \tilde H^{G,0}_\alpha(S^0;\Mackey T\mixprod \MackeyOp U)
\]
is the slant product
\[
 -\slant- \colon
  \tilde H_{G,\delta}^\beta(X;\Mackey T) 
     \tensor \tilde H^{G,\delta_\Delta}_{\alpha+\beta}(S^0\smsh X;\MackeyOp U)
  \to \tilde H^{G,0}_\alpha(S^0;\Mackey T\mixprod \MackeyOp U)
\]
where $X$ is an ex-$G$-space over $B$ and $S^0$ is considered as an ex-$G$-space over $*$.
Here we are using $\delta_\Delta\dimpred 0\times\delta$ and the fact that
$\delta_\Delta(G\times G/\Delta) = 0$ and $\delta_\Delta|\Delta = \delta$.
\end{definition}

Note that we can express the naturality of evaluation by saying that the adjoint map
\[
 \tilde H_{G,\delta}^\beta(X;\Mackey T) \to
  \Hom(\tilde H^{G,\delta}_{\alpha+\beta}(X;\MackeyOp U), \tilde H^{G,0}_\alpha(S^0;\Mackey T\mixprod \MackeyOp U))
\]
is contravariant in $X$. When $\alpha = n\in\Z$, 
$\tilde H^{G,0}_n(S^0;\Mackey T\mixprod\MackeyOp U)$ is nonzero only when $n = 0$ and is then
$(\Mackey T\mixprod\MackeyOp U)(G/G) \iso \Mackey T\tensor_{\stab\Pi_{G,\delta}B}\MackeyOp U$, giving the evaluation
\[
 \tilde H_{G,\delta}^\beta(X;\Mackey T) \to
  \Hom(\tilde H^{G,\delta}_{\beta}(X;\MackeyOp U), \Mackey T\tensor_{\stab\Pi_{G,\delta}B}\MackeyOp U).
\]
There are many other interesting variations available, which we leave to the imagination
of the reader.

\begin{definition}\label{def:paramcapproduct}
Let $\delta$ and $\epsilon$ be familial dimension functions for $G$ and let
$\zeta\dimpred \delta\times\epsilon$ be a familial dimension function for $G\times G$;
assume that $\Delta\in\F(\zeta)$ and
write $\zeta$ again for $\zeta|\Delta$.
Let $B$ be a $G$-space and let $X$ be an ex-$G$-space over $B$.
Let $\Mackey T$ be a contravariant $\stab\Pi_{G,\epsilon}B$-module and
let $\MackeyOp U$ be a covariant $\stab\Pi_{G,\zeta}B$-module.
The {\em cap product}
\[
 -\cap- \colon
  \tilde H_{G,\epsilon}^{\beta}(X;\Mackey T)
    \tensor \tilde H^{G,\zeta}_{\alpha+\beta-\zeta(G\times G/\Delta)}(X;\MackeyOp U)
  \to \tilde H^{G,\delta}_{\alpha}(X;\Mackey T\mixprod\Delta^B_!\MackeyOp U)
\]
is the composite
\begin{align*}
 \tilde H_{G,\epsilon}^{\beta}(X;\Mackey T)
    &\tensor \tilde H^{G,\zeta}_{\alpha+\beta-\zeta(G\times G/\Delta)}(X;\MackeyOp U) \\
 &\to \tilde H_{G,\epsilon}^{\beta}(X;\Mackey T)
    \tensor \tilde H^{G,\zeta}_{\alpha+\beta-\zeta(G\times G/\Delta)}(X\smsh X;\Delta^B_!\MackeyOp U) \\
 &\to \tilde H^{G,\delta}_{\alpha}(X;\Mackey T\mixprod\Delta^B_!\MackeyOp U)
\end{align*}
where the first map is induced by the diagonal $X\to X\smsh X$ and the
unit $\MackeyOp U \to (\Delta^B)^*\Delta^B_!\MackeyOp U$, 
$\Delta^B\colon B\to B\times B$ being the diagonal.
\end{definition}

We get interesting special cases by considering particular choices of $\delta$, $\epsilon$, and $\zeta$.
One we'll use later is the case $\delta = \Lie - \epsilon$ and $\zeta = \Lie_\Delta$.
Using the facts that $\Lie_\Delta(G\times G/\Delta) = 0$ and
$\Lie_\Delta|\Delta = \Lie$, the cap product in this case takes the form
\[
 -\cap- \colon
  \tilde H_{G,\epsilon}^{\beta}(X;\Mackey T)
    \tensor \tilde \H^G_{\alpha+\beta}(X;\MackeyOp U)
  \to \tilde H^{G,\Lie-\epsilon}_{\alpha}(X;\Mackey T\mixprod\Delta^B_!\MackeyOp U).
\]
Specializing further, we can use $\MackeyOp U = \Mackey A_{G/G} = \rho^*\Mackey A_{G/G}$
(considered as a covariant $\stab\Pi_{G,\Lie}B$-module) and the following map.

\begin{proposition}
If $\Mackey T$ is a contravariant $\stab\Pi_{G,\epsilon}B$-module, then
there is a natural map
$\Mackey T\mixprod \Delta^B_!\Mackey A_{G/G} \to \Mackey T$
of covariant $\stab\Pi_{G,\Lie-\epsilon}B$-modules,
which agrees with the isomorphism of Proposition~\ref{prop:mixprodidentity}
when $B=*$.
\end{proposition}

\begin{proof}
Write $\delta = \Lie-\epsilon$. Let
\[
 \Delta\colon \stab\Pi_{G,\Lie}B \to \stab\Pi_{G\times G,\Lie_\Delta}(B\times B)
\]
be the combination of the inclusion of the diagonal of $B\times B$
and extension from the diagonal copy of $G$ to $G\times G$.
For $a\in \stab\Pi_{G,\delta}B$, let $D(a)$ denote the homological dual of $a$,
considered as an element of $\stab\Pi_{G,\epsilon}B$.
We first define a map
\[
 (p^*i_! \Delta_! \Mackey A_{G/G})(a\tensor -)
   \to \stab\Pi_{G,\epsilon}B(D(a),-)
\]
as follows, for $b\in\stab\Pi_{G,\epsilon}B$:
\begin{align*}
 (p^*i_! \Delta_! \Mackey A_{G/G})&(a\tensor b) \\
  &\iso \int\nolimits^{c\in\stab\Pi_{G,\Lie}B}
    h(G\times G)\Spec{}{B\times B}(\Delta_!(c),a\smsh_{B\times B} b)
    \tensor hG\Spec{}{}(S, \rho_!(c)) \\
  &\iso \int\nolimits^{c\in\stab\Pi_{G,\Lie}B}
    hG\Spec{}{B}(c,\Delta^*(a\smsh_{B\times B} b))
    \tensor hG\Spec{}{}(S, \rho_!(c)) \\
  &\iso \int\nolimits^{c\in\stab\Pi_{G,\Lie}B}
    hG\Spec{}{B}(c,a\smsh_{B} b)
    \tensor hG\Spec{}{}(S, \rho_!(c)) \\
  &\to hG\Spec{}{}(S, \rho_!(a\smsh_{B} b)) \\
  &\iso hG\Spec{}{B}(D(a), b)\\
  &= \stab\Pi_{G,\epsilon}B(D(a),b),
\end{align*}
where the map is given by composition.
The map of the proposition is then
\begin{align*}
 (\Mackey T\mixprod \Delta^B_! \Mackey A_{G/G})(a)
  &= \Mackey T\tensor_{\stab\Pi_{G,\epsilon}B} 
  		(p^*i_!\Delta_! \Mackey A_{G/G})(a\tensor -) \\
  &\to \Mackey T\tensor_{\stab\Pi_{G,\epsilon}B} \stab\Pi_{G,\epsilon}B(D(a),-) \\
  &\iso \Mackey T(D(a)).
\end{align*} 
\end{proof}

Using this map, we have a cap product
\[
 -\cap- \colon
  \tilde H_{G,\epsilon}^{\beta}(X;\Mackey T)
    \tensor \tilde \H^G_{\alpha+\beta}(X;\Mackey A_{G/G})
  \to \tilde H^{G,\Lie-\epsilon}_{\alpha}(X;\Mackey T).
\]

Both the evaluation map and the cap product inherit properties
from Theorem~\ref{thm:paramslantproduct}. In particular, the associativity property
gives us the following:
\[
 (x\cup y)\cap z = x\cap (y\cap z)
\]
and
\[
 \langle x\cup y, z \rangle = \langle x, y\cap z \rangle
\]
when $x$, $y$, and $z$ lie in appropriate groups.

We now look at these pairings on the spectrum level. The main calculation
is the following.

\begin{proposition}
Let $\delta$ be a familial dimension function for $G$,
let $\epsilon$ be a familial dimension function for $K$, and
let $\zeta$ be a familial dimension function for $G\times K$ with $\zeta\dimpred \delta\times\epsilon$.
Let $A$ be a $G$-space, let $\alpha$ be a virtual representation of $\Pi_G A$,
let $B$ be a $K$-space, and let $\beta$ be a virtual representation of $\Pi_K B$.
Let $\Mackey T$ be a contravariant $\stab\Pi_{K,\epsilon}B$-module and
let $\MackeyOp U$ be a covariant $\stab\Pi_{G\times K,\zeta}(A\times B)$-module.
Then, if $q\colon A\times B\to A$ is the projection, we have
\[
 \Mackey\pi_{-\alpha+n}^{G,\Lie-\delta}
   (q_!(H_\epsilon\Mackey T^\beta\smsh_B H^\zeta\MackeyOp U_{\alpha+\beta})^K)
   \iso
   \begin{cases}
     0 & \text{if $n < 0$} \\
     \Mackey T\mixprod \MackeyOp U & \text{if $n = 0$.}
   \end{cases}
\]
\end{proposition}

\begin{proof}
The proof is similar to that of Proposition~\ref{prop:EMmixedproduct}.
If $a\colon G/H\to A$, write as shorthand
\[
 S^{\alpha,\delta,a} = G_+\smsh_H S^{\alpha_0(a)-\delta(G/H),a}.
\]
The necessary calculational input is that,
for $a\colon G/H\to A$, $b\colon K/L\to B$, and
$c\colon G\times K/M \to A\times B$,
\begin{align*}
 hG\Spec{}{A}( S^{-\alpha,\Lie-\delta,a}, 
         &q_!( S^{\beta,\epsilon,b}\smsh_B S^{-\alpha-\beta,\Lie-\zeta,c} )^K ) \\
  &\iso h(G\times K)\Spec{}{A}( S^{-\alpha,\Lie-\delta,a}, 
         q_!( S^{\beta,\epsilon,b}\smsh_B S^{-\alpha-\beta,\Lie-\zeta,c} ) ) \\
  &\iso h(G\times K)\Spec{}{A\times B}( S^{\alpha+\beta,\zeta,c},
         S^{\alpha,\delta,a}\smsh S^{\beta,\epsilon,b} ) 
\end{align*}
using homological duality, so that
\[
 hG\Spec{}{A}( S^{-\alpha,\Lie-\delta,a}, q_!( S^{\beta,\epsilon,b}\smsh_B S^{-\alpha-\beta,\Lie-\zeta,c} )^K )
  \iso (\Mackey A_b \mixprod \MackeyOp A^c)(a)
\]
as functors in $a\in\stab\Pi_{G,\delta}A$, using Example~\ref{prop:paramEMmixedproduct}
and Theorem~\ref{thm:StableMapsOrbitsOverB}.
The result now follows by analyzing the structure provided by
Construction~\ref{con:paramEilenbergMacLane}.
\end{proof}

It follows that there is a map
\[
 q_!(H_\epsilon\Mackey T^\beta\smsh_B H^\zeta\MackeyOp U_{\alpha+\beta})^K 
     \to P^\delta(\Mackey T\mixprod\MackeyOp U)_\alpha
\]
of $G$-spectra over $A$ that is an isomorphism on $\Mackey\pi_{-\alpha}^{G,\Lie-\delta}$,
where $P^\delta(\Mackey T\mixprod\MackeyOp U)_\alpha$ is
a spectrum with $\Mackey\pi^{G,\Lie-\delta}_{-\alpha+n}$ concentrated at $n=0$, as 
in Construction~\ref{con:paramEilenbergMacLane}. 
We then use that $E\F(\zeta)\times p_1^*E\F(\delta) \hmtpc E\F(\zeta)$ because
$\F(\zeta)\subset \F(\delta)\times\A(K)$, where $\A(K)$ is the collection of all subgroups of $K$,
so that
\[
 H^\zeta\MackeyOp U_{\alpha+\beta}\smsh p_1^*E\F(\delta)_+ \hmtpc H^\zeta\MackeyOp U_{\alpha+\beta},
\]
to get a map
\begin{align*}
  q_!(H_\epsilon\Mackey T^\beta\smsh_B H^\zeta\MackeyOp U_{\alpha+\beta})^K
  &\hmtpc q_!(H_\epsilon\Mackey T^\beta\smsh_B H^\zeta\MackeyOp U_{\alpha+\beta}\smsh p_1^* E\F(\delta)_+)^K \\
  &\hmtpc q_!(H_\epsilon\Mackey T^\beta\smsh_B H^\zeta\MackeyOp U_{\alpha+\beta})^K\smsh E\F(\delta)_+ \\
  &\to P^\delta(\Mackey T\mixprod\MackeyOp U)_\alpha \smsh E\F(\delta)_+ \\
  &\hmtpc H^\delta(\Mackey T\mixprod\MackeyOp U)_\alpha.
\end{align*}
The slant product
\[
 -\slant - \colon 
  \tilde H_{K,\epsilon}^\beta(Y;\Mackey T) \tensor
    \tilde H^{G\times K,\zeta}_{\alpha+\beta}(X\smsh Y; \MackeyOp U)
  \to \tilde H^{G,\delta}_\alpha(X; \Mackey T \mixprod \MackeyOp U).
\]
is then represented as follows:
\begin{align*}
 [Y, H_\epsilon\Mackey T^\beta]_{K,B}
   \tensor{}& [S, \rho_!(H^\zeta\MackeyOp U_{\alpha+\beta} \smsh_{A\times B} (X\smsh Y))]_{G\times K} \\
  &\to [S, \rho_!(H^\zeta\MackeyOp U_{\alpha+\beta} 
  			\smsh_{A\times B} (X\smsh H_\epsilon\Mackey T^\beta))]_{G\times K} \\
  &\iso [S, \rho_!((H_\epsilon\Mackey T^\beta 
  			\smsh_B H^\zeta\MackeyOp U_{\alpha+\beta})\smsh_A X)]_{G\times K} \\
  &\iso [S, \rho_!((H_\epsilon\Mackey T^\beta\smsh_B H^\zeta\MackeyOp U_{\alpha+\beta})\smsh_A X)^K]_G \\
  &\iso [S, \rho_!(q_!(H_\epsilon\Mackey T^\beta\smsh_B H^\zeta\MackeyOp U_{\alpha+\beta})^K \smsh_A X)]_G \\
  &\to [S, \rho_!(H^\delta(\Mackey T\mixprod\MackeyOp U)_\alpha\smsh_A X)]_G.
\end{align*}

\section{The Thom Isomorphism and Poincar\'e Duality}\label{sec:genThom}

\subsection{The Thom Isomorphism}

\begin{definition}
Let $p\colon E\to B$ be a $G$-vector bundle and let
$\gamma$ be the associated representation of $\Pi_G B$.
Let $T(p) = D(p)/_B S(p)$ be the {\em fiberwise Thom space} of $p$.
A {\em Thom class} for $p$ is an element
$t \in \tilde H_G^\gamma(T(p);\Mackey A_{G/G})$ such that, for each
$G$-map $b\colon G/K\to B$,
 \begin{align*}
 b^*(t) \in{} &\tilde H_G^{\gamma(b)}(b^*T(p);\Mackey A_{G/G}) \\
 &\iso \tilde H_G^{\gamma(b)}(G\times_K S^V;\Mackey A_{G/G}) \\
 &\iso \tilde H_K^{V}(S^V;\Mackey A_{K/K}) \\
 &\iso A(K)
 \end{align*}
is a generator.
Here, $V$ is a representation of $K$ such that
$\gamma(b) \hmtpc G\times_K S^V$.
 \end{definition}

As in Remark~\ref{rem:Thomclasshome}, the Thom class must live in ordinary cohomology
$\tilde H_{G,\delta}^*$ with $\delta = 0$.

Since a global Thom class is characterized by its local behavior, and
locally a $G$-vector bundle is a $V$-bundle, the next proposition
follows from the $RO(G)$-graded analogue, Proposition~\ref{prop:vThomclass}.

\begin{proposition}\label{prop:gThomclass}
The following are equivalent for a cohomology class $t\in \tilde H_G^\gamma(T(p);\Mackey A_{G/G})$.
 \begin{enumerate}
 \item
$t$ is a Thom class for $p$.
 \item
For every subgroup $K\subset G$,
$t|K\in \tilde H_K^{\gamma|K}(T(p);\Mackey A_{K/K})$
is a Thom class for $p$ as a $K$-bundle.
 \item
For every subgroup $K\subset G$, 
$t^K\in \tilde H_{WK}^{\gamma^K}(T(p)^K; \Mackey A_{WK/WK})$
is a Thom class for $p^K\colon E^K\to B^K$ as a $WK$-bundle.
 \item
For every subgroup $K\subset G$,
$t^K|e \in \tilde H^{|\gamma^K|}(T(p)^K;\Z)$
is a Thom class for $p^K$ as a nonequivariant bundle.
 \end{enumerate}
 \qed
 \end{proposition}

The following generalizes Theorem~C of \cite{CW:duality}.

\begin{theorem}[Thom Isomorphism]
If $p\colon E\to B$ is any $G$-vector bundle
and $\gamma$ is the associated representation of $\Pi_G B$, 
then there exists a Thom class
$t\in \tilde H_G^\gamma(T(p);\Mackey A_{G/G})$. For any Thom class $t$,
the maps
 \[
 t\cup - \colon H_{G,\delta}^\alpha(B;\Mackey T)
 \to \tilde H_{G,\delta}^{\alpha+\gamma}(T(p);\Mackey T)
 \]
is an isomorphism for any familial $\delta$.
 \end{theorem}

\begin{proof}
As in the proof of Theorem~\ref{thm:vthomiso},
the theorem is clear in the special case that $p$ is a bundle over an orbit.
The general case follows, as it does in the nonequivariant case for
twisted coefficients, by a Mayer-Vietoris patching argument.
The key point is that we can choose a compatible collection of local
classes because the action of $\Pi_G B$ on $\gamma$ is the same
as the action on the fibers of $p$.
 \end{proof}

\subsection{Poincar\'e Duality}

We are now in a position to describe Poincar\'e duality for 
arbitrary compact smooth $G$-manifolds.
(Again, the noncompact case can be handled using cohomology with
compact supports.)

\begin{definition}
Let $M$ be a closed smooth $G$-manifold and let $\tau$ be the tangent representation
of $\Pi_G M$, i.e., the representation associated with the tangent bundle. 
Think of $M$ as a $G$-space over itself in the following.
A {\em fundamental class} of $M$ is
a class $[M]\in \H_\tau^G(M;\Mackey A_{G/G})$ such that, for each 
point $m\in M$, thought of as the map
$m\colon G/G_m\to M$
 with image $Gm$, and tangent plane $\tau(m) = G\times_{G_m}V$,
the image of $[M]$ in
 \begin{align*}
 \H_\tau^G(M,M-Gm;\Mackey A_{G/G}) 
 &\iso \tilde\H_V^G(G_+\smsh_{G_m}S^{V-\Lie(G/G_m)};
    \Mackey A_{G/G}) \\
 &\iso \tilde\H_V^{G_m}(S^V; \Mackey A_{G/G}) \\
 &\iso A(G_m)
 \end{align*}
is a generator.
 \end{definition}

The fundamental class $[M]$ is related to fundamental classes of the fixed
submanifolds $M^K$ as in
Proposition~\ref{prop:vFundclass}.

\begin{proposition}\label{prop:gFundclass}
The following are equivalent for a dual homology class 
$\mu\in \H_\tau^G(M;\Mackey A_{G/G})$.
 \begin{enumerate}
 \item
$\mu$ is a fundamental class for $M$ as a $G$-manifold.
 \item
For every subgroup $K\subset G$,
$\mu|K \in \H_{\tau|K}^K(M;\Mackey A_{K/K})$ is a fundamental class for
$M$ as a $K$-manifold.
 \item
For every subgroup $K\subset G$,
$\mu^K \in \H_{\tau^K}^{WK} (M^K; \Mackey A_{WK/WK})$ is a fundamental class for
$M^K$ as a $WK$-manifold.
 \item
For every subgroup $K\subset G$,
$\mu^K|e \in H_{|\tau^K|}(M^K;\Z)$ is a fundamental class for $M^K$
as a nonequivariant manifold. \qed
 \end{enumerate}
\end{proposition}

We now have sufficient machinery in place to prove Poincar\'e duality for
arbitrary $G$-manifolds either geometrically, along the lines of
\cite{MS:characteristicclasses} or \cite{CW:duality}, or
homotopically as in Section~\ref{sec:vThomIso}.
The homotopical approach uses, of course, the duality
of Theorem~\ref{thm:manifoldduality}.

\begin{theorem}[Poincar\'e Duality]
Every closed smooth $G$-manifold $M$ 
has a fundamental class $[M]\in\H_\tau^G(M;\Mackey A_{G/G})$, and
 \[
 -\cap [M] \colon H_{G,\delta}^\gamma(M;\Mackey T) \to H_{\tau-\gamma}^{G,\Lie-\delta}(M;\Mackey T)
 \]
is an isomorphism for every familial $\delta$ such that $\F(\delta)$ contains every isotropy
subgroup of $M$.
\qed
 \end{theorem}

In particular, when $\delta = 0$ we get the isomorphism
 \[
 -\cap [M] \colon H_{G}^\gamma(M;\Mackey T) \to \H_{\tau-\gamma}^{G}(M;\Mackey T)
 \]
and when $\delta = \Lie$ we get the isomorphism
 \[
 -\cap [M] \colon \H_{G}^\gamma(M;\Mackey T) \to H_{\tau-\gamma}^{G}(M;\Mackey T).
 \]

If $M$ is a compact $G$-manifold with boundary, then we get relative,
or Lefschetz, duality.

\begin{definition}
Let $M$ be a compact $G$-manifold with boundary, with tangent representation $\tau$. 
A {\em fundamental class} of $M$ is a dual homology
class $[M,\bndry M]\in \H_\tau^G(M,\bndry M;\Mackey A_{G/G})$ such that, for each 
point $m\in M-\bndry M$,
thought of as the map
$m\colon G/G_m\to M$
 with image $Gm$, and tangent plane $\tau(m) = G\times_{G_m}V$,
the image of $[M,\bndry M]$ in
 \begin{align*}
 \H_\tau^G(M,M-Gm;\Mackey A_{G/G}) 
 &\iso \tilde\H_V^G(G_+\smsh_{G_m}S^{V-\Lie(G/G_m)};
    \Mackey A_{G/G}) \\
 &\iso \tilde\H_V^{G_m}(S^V; \Mackey A_{G/G}) \\
 &\iso A(G_m)
 \end{align*}
is a generator.
 \end{definition}

There is an obvious relative version of Proposition~\ref{prop:gFundclass}.

\begin{theorem}[Lefschetz Duality]
Every compact smooth $G$-manifold $M$ has a fundamental class 
$[M,\bndry M]\in\H_V^G(M,\bndry M;\Mackey A_{G/G})$,
and the following are isomorphisms:
 \[
 -\cap [M,\bndry M] \colon H_{G,\delta}^\gamma(M;\Mackey T) \to
   H_{\tau-\gamma}^{G,\Lie-\delta}(M,\bndry M;\Mackey T)
 \]
and
 \[
 -\cap [M,\bndry M] \colon H_{G,\delta}^\alpha(M,\bndry M;\Mackey T) \to
   H_{\tau-\alpha}^{G,\Lie-\delta}(M;\Mackey T).
 \]
 \qed
 \end{theorem}

As in the nonparametrized version, we get the following special cases when
$\delta=0$ or $\delta=\Lie$:
\begin{align*}
 -\cap [M,\bndry M] &\colon H_{G}^\gamma(M;\Mackey T) \to
   \H_{\tau-\gamma}^{G}(M,\bndry M;\Mackey T),
\\
 -\cap [M,\bndry M] &\colon H_{G}^\gamma(M,\bndry M;\Mackey T) \to
   \H_{\tau-\gamma}^{G}(M;\Mackey T),
\\
 -\cap [M,\bndry M] &\colon \H_{G}^\gamma(M;\Mackey T) \to
   H_{\tau-\gamma}^{G}(M,\bndry M;\Mackey T), \qquad\text{and}
\\
 -\cap [M,\bndry M] &\colon \H_{G}^\gamma(M,\bndry M;\Mackey T) \to
   H_{\tau-\gamma}^{G}(M;\Mackey T).
\end{align*}

\backmatter
\bibliographystyle{amsplain}
\bibliography{Topology}

\end{document}